МИНИСТЕРСТВО ОБРАЗОВАНИЯ И НАУКИ УКРАИНЫ
ЗАПОРОЖСКИЙ НАЦИОНАЛЬНЫЙ УНИВЕРСИТЕТ

**С. В. Курапов**
**М. В. Давидовский**

# АЛГОРИТМИЧЕСКИЕ МЕТОДЫ КОНЕЧНЫХ ДИСКРЕТНЫХ СТРУКТУР

# ИЗОМОРФИЗМ НЕСЕПАРАБЕЛЬНЫХ ГРАФОВ

**(монография)**

**Курапов Сергей Всеволодович**

**Давидовский Максим Владимирович**

Для решения задачи изоморфизма графов, в качестве математической структуры, предлагается использовать понятие спектров реберных разрезов и реберных циклов графа. Реберный разрез определяется ребром и инцидентными к нему вершинами. В отличие от порождения итерированных реберных графов, рассматривается итерированная цепочка квалиразрезов исходного графа, порождаемая реберными разрезами и определяемая рекуррентным соотношением. Реберный цикл определяется множеством изометрических циклов графа. В работе рассматриваются вопросы построения спектров реберных разрезов $W_s$ и спектра реберных циклов $T_c$ графа G. Показана, что в основе формирования спектров находится матрицой инциденций графа. Показана независимость построения структуры графа от нумерации вершин и ребер. Показана необходимость и достаточность спектров реберных разрезов и спектра реберных циклов для определения изоморфизма графовых структур. Рассмотрена связь внутренних структур графа с теоремой Уитни.

Для научных работников, студентов и аспирантов высших учебных заведений использующих методы прикладной математики.

## Содержание

# ВВЕДЕНИЕ

Одним из основных классов задач прикладной теории графов является класс задач различения графов и различения расположения фрагментов графов. История создания и развития методов решения задач различения графов, насчитывает более 50 лет. В настоящее время в решении задач распознавания изоморфизма графов, распознавания изоморфного вложения графов и смежных с ними задач достигнут большой теоретический и практический прогресс.

Теоретически алгоритм проверки пары графов G и H на изоморфизм существует. Будем переставлять строки и соответствующие столбцы матрицы смежностей графа G до тех пор, пока она не превратится в другую, равную матрице смежностей графа H, или остановимся после n! перестановок, если графы не изоморфны. Однако такое решение подразумевает полный перебор вариантов [2,5,11].

Поэтому решение задачи распознавания изоморфизма графов заключается в нахождении таких графовых структур (если это возможно), которые позволяют произвести распознавание графов G и H не переборными алгоритмами.

Как правило, почти для любой алгоритмической задачи теории графов удается построить полиномиальный алгоритм или доказать ее принадлежность к классу NP-полных задач [6,8]. Задача определения изоморфизма графов – это одна из немногих классических задач теории графов, для которой не удалось осуществить ни то, ни другое (хотя для некоторых специальных классов графов удалось построить полиномиальные алгоритмы [8,17]).

Как удачно подметили Погребной Ан.В. и Погребной В.К. «Следует признать, что в теории графов сложилась парадоксальная ситуация, которую можно обозначить как отсутствие возможностей для однозначной идентификации структуры графа. Действительно, принимаясь за исследование очередного графа, мы не можем установить – является ли он одним из ранее рассмотренных или нет. Для однородных графов такая ситуация возникает уже, например, при приближении значений $n$ к 15–20 вершинам. Что касается неоднородных графов, содержащих десятки и сотни вершин, то для таких размерностей даже применение алгоритмов, учитывающих специфику отдельных графов, объёмы вычислений часто становятся также нереалистичными.

Получается, что в теории графов нет инструмента для описания основного объекта исследований – структуры графа. Произвольная нумерация вершин, которая традиционно используется для представления графа в компьютере, не является инвариантным описанием его структуры, т. е. независимым от нумерации вершин. Один и тот же граф с произвольными нумерациями вершин воспринимается как два разных графа.»[30].

В настоящее время известны полиномиальные алгоритмы для следующих классов графов: графы ограниченного рода, графы ограниченной степени, графы с ограниченной кратностью собственных значений, *k*-разделимые графы, *k*-стягиваемые графы, графы не стягиваемые на $K_{3,q}$. Особое внимание заслуживают сильнорегулярные графы. Из [25] известно, что сильнорегулярные графы образуют класс наиболее трудных для распознавания изоморфизма графов и большинство контрпримеров для эвристических алгоритмов распознавания изоморфизма графов входит в указанный класс.

Для решения задачи распознавания изоморфизма графов наметилось несколько подходов. Один из подходов связан с построением непереборных алгоритмов, рекурсивно улучшающих свою эффективность в смысле полноты или чувствительности используемых характеристик графа, неизменных (инвариантных) относительно изоморфизма графов. Такие характеристики графа называются инвариантами. Поскольку инвариант графа не меняет своих значений на изоморфных графах, то равенство инвариантов является необходимым условием изоморфизма графов.

При анализе структур графов, в частности в компьютерной химии, часто используются числовые характеристики, которые отражают отдельные свойства структур и не зависят от нумерации вершин. Такие характеристики названы инвариантами. Примером инварианта является *n* -мерный вектор локальных степеней с упорядоченными по возрастанию значениями степеней вершин графа G. Можно привести примеры и других инвариантов. Однако применяемые в настоящее время инварианты не гарантируют изоморфизма графов, т.е. не являются полными [31].

Примерами полного инварианта графа является максимальный и минимальный двоичные коды матриц смежностей этого графа. Для получения двоичного кода по матрице смежности графа пронумеруем ее элементы подряд слева направо и сверху вниз. Тогда полученная матрица может рассматриваться как $n^2$-разрядное двоичное число для графа G с *n* вершинами. Наименьший из них называется минимальным кодом графа G, а наибольший – максимальным кодом. По любому из этих кодов и количеству вершин графа можно восстановить одну из его матриц смежности, а значит и сам граф с точностью до изоморфизма. Однако процесс вычисления минимального или максимального двоичного кода матрицы смежностей для заданного *n*-вершинного графа столь же труден, как и лобовой перебор *n*! соответствий вершин двух графов.

Другой подход к разработке алгоритмов распознавания изоморфизма графов, как раз отличается тем, что обязательно включает в себя процедуру направленного перебора на одном из этапов поиска изоморфной подстановки. Большинство известных переборных алгоритмов распознавания изоморфизма графов основано на разбиении множества вершин каждого из графов на подмножества, которые в случае изоморфизма также должны соответствовать друг другу. Подход Ласло Бабаи использует информацию о специальном строении, которым должно обладать группа автоморфизма графа и глубое применение алгоритмических методоы теории перестановок. Подробно метод рассмотрен в работах [8,40-41].

Для нахождения изоморфизма используют методы спектральной теории графов. Характеристический многочлен $|\lambda I — A|$ матрицы смежности A графа G называется характеристическим многочленом графа G и обозначается как $Pg(\lambda)$. Собственные значения матрицы A (т. е. нули многочлена $|\lambda I — A|$) и спектр матрицы A, состоящий из собственных значений, называются собственными значениями в спектре графа G. Если $\lambda_1,\ldots,\lambda_n$ — собственные значения графа G, то весь спектр обозначается как $Sp(G) = [\lambda_1,\ldots,\lambda_n]$. Очевидно, что изоморфные графы имеют один и тот же спектр.

Два графа называются изоспектральными или коспектральными, если матрицы смежности графов имеют одинаковые мультимножества собственных значений. Изоспектральные графы не обязательно изоморфны, но изоморфные графы всегда изоспектральны [35].

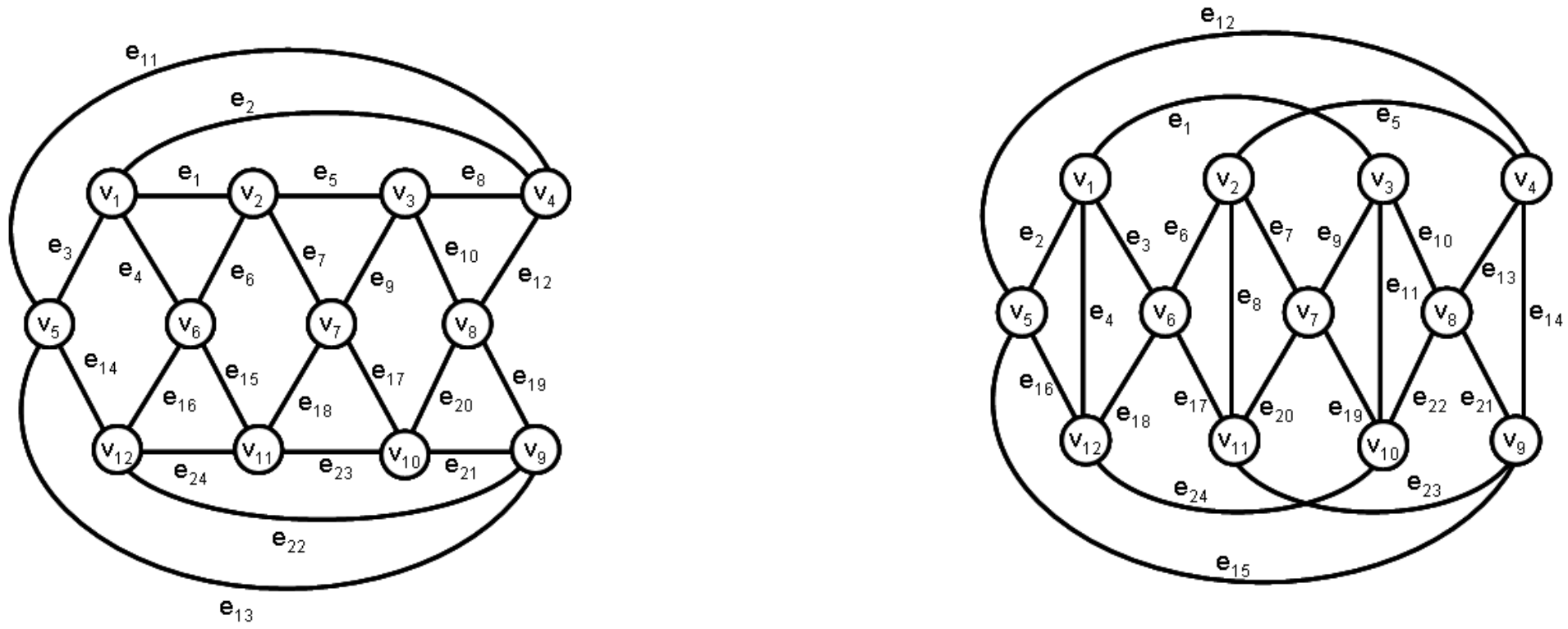

Рис. 1. Не изоморфные изоспектральные графы.

Общим для указанных на рисунке графов спектром является (4,2,2,2,0,0,0,-2,-2,-2,-2,-2).

По мнению ряда специалистов, в настоящее время самым быстродействующим считается алгоритм Ласло Бабаи [40-41]. Этот алгоритм основан на методе виртуального «окрашивания» вершин графа. Сначала случайным образом выбираются несколько вершин, они «окрашиваются» в разные цвета. Затем выбираются несколько вершин во втором графе, предположительно соответствующих вершинам из первого графа, им присваиваются те же цвета. В конце концов, перебираются все варианты. После первоначального выбора варианта раскраски вершин в графах G и H, алгоритм окрашивает предположительно изоморфные вершины, соседствующие с первоначально окрашенными в другие цвета. Процесс продолжается до тех пор, пока не закончатся связи между вершинами. Существуют и другие подходы к решению задачи распознования изоморфизма графов, но как правило, все они основаны на рассмотрении матрицы смежностей графа [2,5].

Длительные попытки решения задачи распознования графов, на основе использования только разрезов графа, не привели к желаемому результату. Такая ситуация вызвана тем фактом, что однозначно определить множество центральных разрезов, не представляет большого труда, это можно проделать с помощью карандаша и бумаги, даже для графов большого размера. Для определения множества простых циклов, применяется фундаментальная система циклов, основанная на случайном построении дерева графа, носящая неоднозначный характер относительно длин циклов. В отличии от традиционного подхода, разработанная в последние годы, методика построения изометрических циклов, представляющих собой определенную теоретико-множественную структуру графа, определена однозначно. Однако для выделения множества изометрических циклов графа, уже нужно применить компьютер.

Важнейшим понятием теории графов является понятие  - суграфа графа. Независимая система однореберных суграфов и операция кольцевого суммирования, индуцирует все множество суграфов графа, и формирует линейное пространство суграфов £(G). Элементами пространства суграфов £(G) являются, подпространство циклов C(G) и подпространство разрезов S(G) [34].

Алгебраические методы теории графов связывают в единое целое и суграфы, и разрезы, и циклы. Однобокое применение только части пространства суграфов, порождает не точное представление о природе явления и некорректное решение. К примеру, можно привести спектральную теорию графов [35]. Спектральная теория графов, на основе матрицы смежностей графа (прообраз множества центральных разрезов графа), определяет собственные значения матрицы смежностей графа A (т. е. нули характеристического многочлена $|\lambda I — A|$) и строит спектр собственных значений. Очевидно, что изоморфные графы имеют один и тот же спектр. Но, к сожалению, известнв случаи равенства собственных значений для неизоморфных графов. К недостаткам подхода, следует отнести невозможность постановки в соответствие вершине графа – элемент спектра собственнух значений.

Примерами построения математических моделей основанных на связи двух подпространств C(G) и S(G) графа, является взаимодействие матрицы фундаментальных циклов и матрицы фундаментальных разрезов графа [1,9,10,26,33,34,36], построение топологического рисунка плоского суграфа [16,21,22]. Вращение вершин в топологическом рисунке (циклическая запись центральных разрезов), порождает систему независимых простых циклов с нулевым значением функционала Маклейна и наоборот. Другим классическим примером, является теорема Уитни [24]. Теорема Уитни об изоморфизме графов, сформулированная Хасслером Уитни в 1932 году, гласит, что два связных графа с одинаковым количеством вершин и ребер изоморфны тогда и только тогда, когда изоморфны их рёберные графы [36,55,56]. Изометрические циклы реберного графа L(G), являются прообразами центральных разрезов и изометрических циклов графа G. И это определяет процесс различения двух графов G и H.

Если в качестве математической структуры использовать основные реберные разрезы и изометрические циклы графа, то применение преобразования $\gamma$ порождает элементы спектров реберных разрезов и реберных циклов графа. Определение количества участия каждого ребра во множестве суграфов спектра реберных разрезов и реберных циклов порождает качественную картину и определяет его меру (вес ребра) [23,24]. Совокупность весов ребер порождает аналог полного инварианта графа служащий его характеристикой. С помощью данной характеристики, графы можно сопоставлять и сравнивать между собой.

Применение в прикладной теории графов, таких структур, как множество изометрических циклов и реберных разрезов графа с определенной мерой, позволило определить решение многих задач:

- создать методы вычисления векторных инвариантов для распознавания изоморфизма графов;
- на базе векторных инвариантов и свойства устойчивости подмножества вершин, построить методы для выделения образующих и орбит группы автоморфизма графа;
- создать вычислительные методы определения биективного соответствия вершин графов G и H;
- произвести удобную классификацию видов графов, с целью определения применения вычислительных алгоритмов по назначению, в зависимости от вида графа.

Выделение множество несепарабельных графов в отдельный вид, по законам диакоптики, [19] позволило многие задачи прикладной теории графов, свести к методам дискретной оптимизации

В данной работе представлены необходимые и первоначальные сведения о структурах и математических моделях применяемых для решения задачи различения графов.

В целях более удобного представления структур, иногда обозначение ребер, вершин и циклов заменяется цифрами и тогда принадлежность элементов определяется контекстом.

Работа предназначена для прикладных математиков, занимающихся разработкой и внедрением методов и алгоритмов дискретной математики в различных областях промышленности. С целью детальной проработки алгоритмов представлено множество примеров. Это способствует более углубленному процессу изучения и описания математических структур дискретной математики студентами и преподавателями.

## Основные обозначения

G(V,E;P) – граф;

V – множество вершин графа;

E – множество ребер графа;

P – трехместный предикат графа, ставящий в соответствие одному ребру графа две вершины;

$n$ = card V – количество вершин в графе;

$m$ = card E – количество ребер в графе;

$\rho(G)$ – вектор локальных степеней вершин;

A(G) – матрица смежностей графа;

B(G) – матрица инциденций графа;

A(L(G)) – матрица смежностей реберного графа G;

$C_\tau$ – множество изометрических циклов графа G;

$W_s(G)$ – спектр реберных разрезов графа G;

$T_c(G)$ – спектр реберных циклов графа G;

$c_i$ – суграф цикла с номером $i$;

S(G) – подпространство разрезов графа G;

C(G) – подпространство циклов графа G;

$£_G$ – пространство суграфов графа G;

$s(v_j)$ – суграф центрального разреза вершины $v_j$;

$w_l(e_i)$ – суграф реберного разреза уровня $l$ в спектре реберных разрезов текущего ребра $e_i$;

$\tau_l(e_i)$ – суграф реберного цикла уровня $l$ в спектре реберных циклов текущего ребра $e_i$;

$\beta(G)$ – вектор – столбец ребер графа G;

$l_j \in L_s$ – номер уровня в спектре реберных разрезов графа G;

$L_s$ – множество номеров уровней в спектре реберных разрезов графа G;

$card\ L_s$ – количество уровней в спектре реберных разрезов графа G;

$\oplus$ – операция кольцевой суммы суграфов;

$\varepsilon(w(e_i))$ – кортеж весов ребер для строки с аргументом $e_i$ спектра реберных разрезов;

$\xi(w(l_j))$ – кортеж весов ребер для столбца уровня $l_j$ спектра реберных разрезов;

$\xi_w(G)$ – кортеж весов ребер в спектре реберных разрезов;

$\varepsilon(\tau(e_i))$ – кортеж весов ребер для строки с аргументом $e_i$ спектра реберных циклов;

$\xi(\tau(l_j))$ – кортеж весов ребер для столбца уровня $l_j$ матрицы $T_c$ спектра реберных циклов;

$\zeta_w(G)$ – кортеж весов вершин в спектре реберных разрезов;

$\zeta_\tau(G)$ – кортеж весов вершин в спектре реберных циклов;

$\xi_\tau(e_i)$ – кортеж весов ребра в спектре реберных циклов;

$\xi_L(e_i)$ – кортеж весов ребра в цифровом инварианте реберного графа;

$\zeta_L(e_i)$ – кортеж весов вершин в цифровом инварианте реберного графа;

$F_w(\xi(G))\,\&\,F_w(\zeta(G))$ или $F(\xi(w(e_i)))\,\&\,F(\zeta(w(e_i)))$ – векторный инвариант спектра реберных разрезов графа G;

$F_\tau(\xi(G))\,\&\,F_\tau(\zeta(G))$ или $F(\xi(\tau(e_i)))\,\&\,F(\zeta(\tau(e_i)))$ – векторный инвариант спектра реберных циклов графа G;

$F_w(\xi(G))\,\&\,F_w(\zeta(G))\,\&\,F_\tau(\xi(G))\,\&\,F_\tau(\zeta(G))$ $F_{es}(G)$ или

$F(\xi(w(e_i)))\,\&\,F(\zeta(w(e_i)))\,\&\,F(\xi(\tau(e_i)))\,\&\,F(\zeta(\tau(e_i)))$ – интегральный инвариант графа G;

$\varnothing$ – пустое множество;

$K_n$ – полный граф с $n$ вершинами;

$K_{n1,n2}$ – двудольный граф с $n_1,n_2$ вершинами.

# Глава 1. ЛИНЕЙНОЕ ПРОСТРАНСТВО СУГРАФОВ ГРАФА

## 1.1 Классификация графов

Важным аспектом описания графов является их четкая классификация, так как для различных классов графов применяются различные методы и алгоритмы преобразования. В общем, все графы можно разделить на ориентированные и неориентированные. В свою очередь, неориентированные графы можно разбить на три крупных класса – класс неориентированных сепарабельных графов, класс неориентированных несепарабельных графов и класс неориентированных ациклических графов.

**Определение 1.1.** Трехсвязный неориентированный граф G, не имеющий мостов и точек сочленения, без петель и кратных ребер и без вершин с локальной степенью меньшей или равной двум, называется *несепарабельным графом* G.

Характеристикой несепарабельного графа является то, что локальная степень любой вершины этого графа больше либо равна трем.

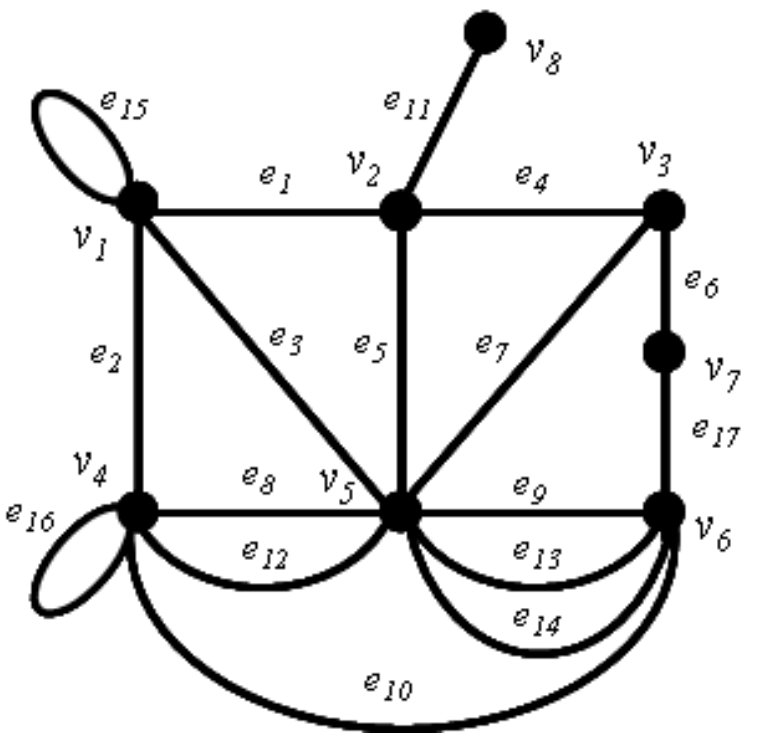

Рис. 1.1. Сепарабельный граф.

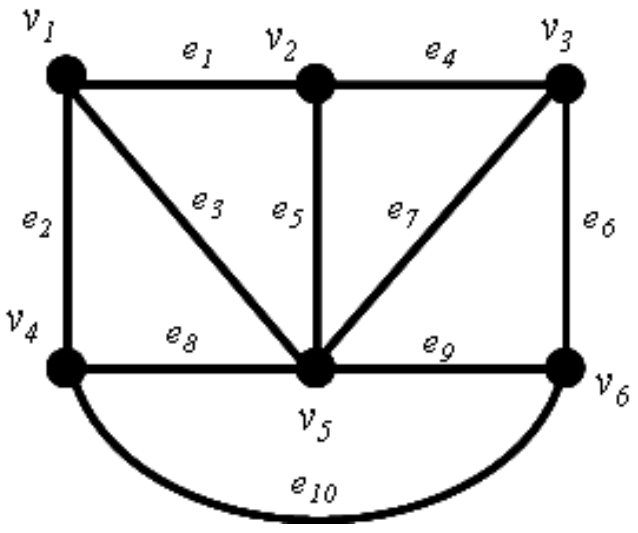

Рис. 1.2. Несепарабельная часть сепарабельного графа

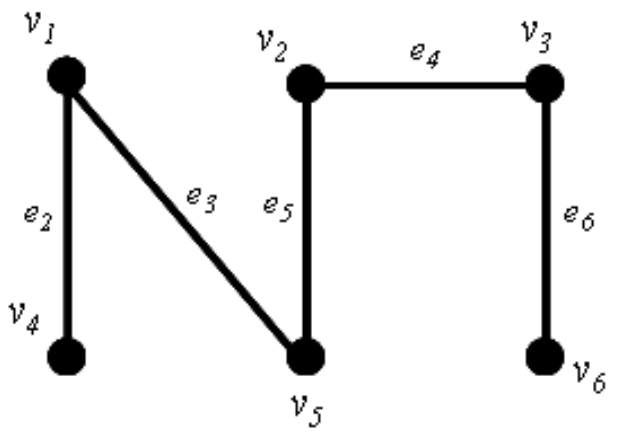

Рис. 1.3. Ациклический граф.

В любом сепарабельном графе (рис. 1.1) всегда можно выделить несепарабельную часть (рис. 1.2), осуществив его последовательное преобразование. В этом случае сепарабельные графы можно представить как кольцевое сложение несепарабельного подграфа $\overline{G}$ и суграфов отдельных частей $G^R$ [34].

Таким образом, для решения прикладных задач теории графов возникает подход, *основанный на цикломатических свойствах* графа. Этот подход состоит из двух последовательных этапов:

1. выделение несепарабельной части графа;

2. добавление удалённых частей до сепарабельного графа (пусть даже с нулевым количеством удалённых частей).

В этом случае, из сепарабельного графа

• удаляются все петли;

• кратные ребра заменяются одним ребром;

• разрезы по точкам сочленения разбивают граф на связанные блоки;

• удаление мостов также разбивает граф на связные блоки;

• цепочка ребер, состоящая из ребер с двухвалентными вершинами, заменяется одним ребром;

• удаляются все ребра, имеющие в своем составе вершины с валентностью равной единице.

В зависимости от способа представления и вида графы можно разделить на следующие подклассы:

• планарных;

• непланарных;

• регулярных;

• полных;

• сильно регулярных;

• изоспектральных;

• двудольных;

• графов Муна-Мозера [37].

Кроме этого, имеются графы с характерными отличительными чертами: пустые графы ($\varnothing$-графы); граф Петерсена; граф тетраэдра и другие графы.

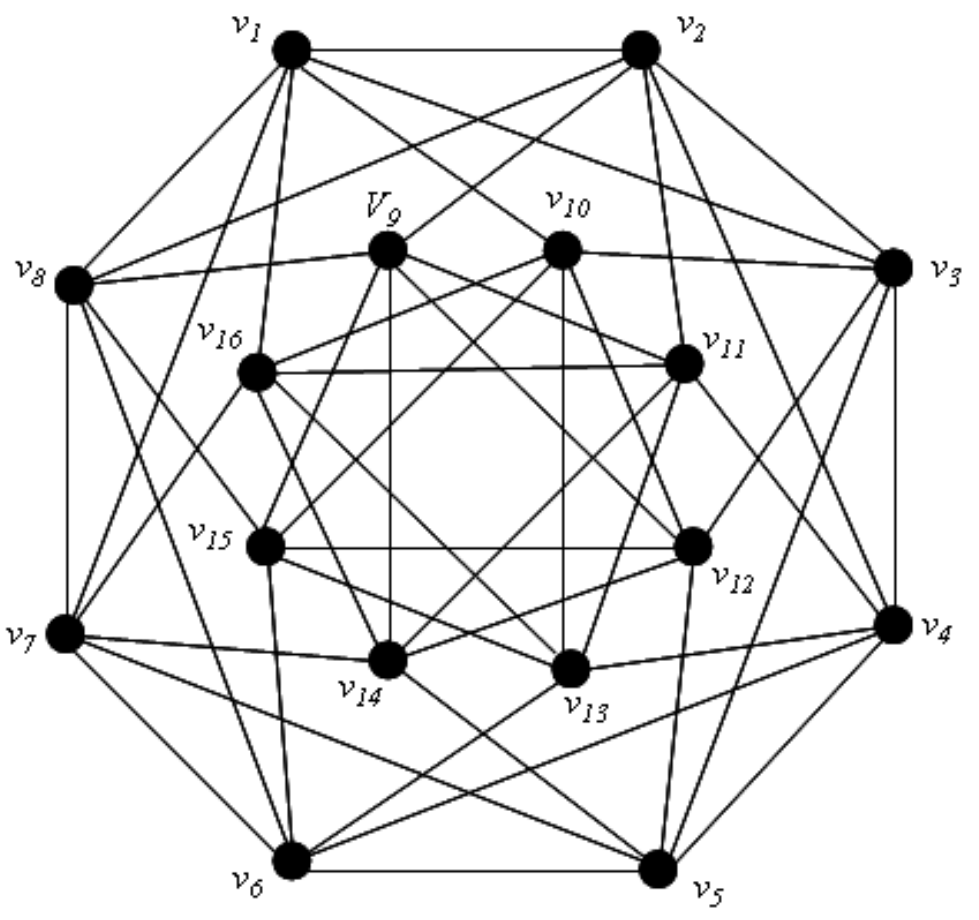

Рис. 1.4. Сильно регулярный граф.

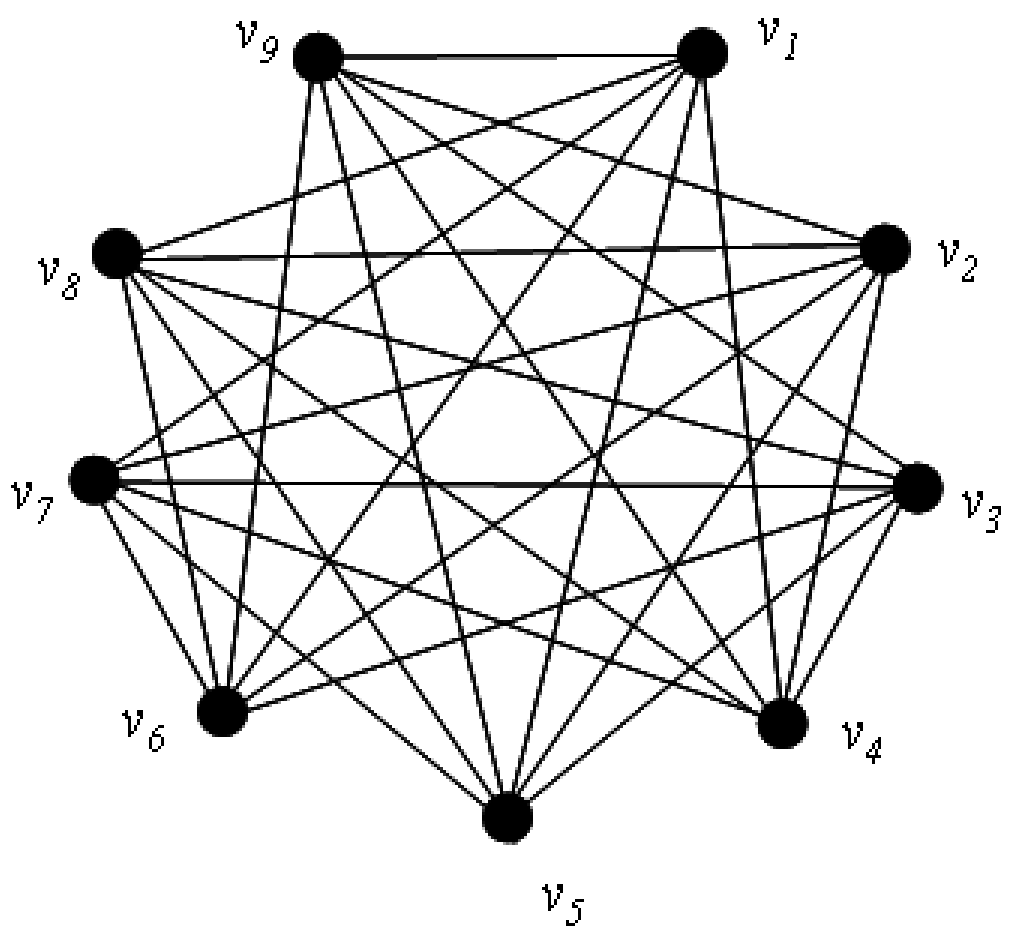

Рис. 1.5. Граф Муна-Мозера.

Впредь будем рассматривать только несепарабельные неориентированные графы за исключением некоторых простых примеров.

## 1.2. Линейное пространство

**Определение 1.2.** Множество £ называется *линейным* или *векторным пространством,* если для всех элементов (векторов) этого множества определены операции сложения и умножения на число и справедливы следующие аксиомы [3,4,18]:

**Аксиома 1**. Каждой паре элементов **x** и **y** из £ отвечает элемент **x** + **y** из £, называемый *суммой* **x** и **y**, причём:

**x** + **y** = **y** + **x** – сложение коммутативно; (1.1)

**x** + **(y** + **z)** = **(x + y) + z** – сложение ассоциативно; (1.2)

**x** + **0** = **x** – существует единственный *нулевой* элемент **0**, такой, что сумма **0 + x** = **x** для любого **x** из £; (1.3)

**x** + (– **x**) = **0** – для каждого элемента **x** из £ существует единственный *противоположный* элемент **–x**, такой что **x + (–x) = 0** для любого **x** из £. (1.4)

**Аксиома 2**. Каждой паре **x** и $\alpha$, где $\alpha$ – число, а **x** элемент из £, отвечает элемент $\alpha\cdot$**x**, называемый *произведением* $\alpha$ и **x**, причём:

$\alpha\cdot(\beta\,\mathbf{x})$ = $(\alpha\beta)\,\mathbf{x}$ – умножение на число ассоциативно; (1.5)

**x** = **x** – для любого элемента **x** из £. (1.6)

**Аксиома 3**. Операции сложения и умножения на число связаны соотношениями:

$\alpha(\mathbf{x} + \mathbf{y}) = \alpha\mathbf{x} + \alpha\mathbf{y}$ – умножение на число дистрибутивно относительно сложения элементов; (1.7)

$(\alpha+\beta)\mathbf{x} = \alpha\mathbf{x} + \beta\mathbf{x}$ – умножение на вектор дистрибутивно относительно сложения чисел. (1.8)

Примером линейного пространства является пространство геометрических радиусов-векторов на плоскости **L** = $\mathbf{R}^2$:

a = $x_1$·i + $y_1$·j, b = $x_2$·i + $y_2$·j,

a + b = ($x_1$ + $x_2$)·i + ($y_1$ + $y_2$)·j, α·*a* = ($\alpha x_1$)·i + ($\alpha y_1$)·j,

0 = 0·i + 0 j, –a = (–$x_1$)·i + (–$y_1$)·j.

Первоначально понятие линейного пространства относилось к геометрическим векторам евклидового пространства, где справедливость аксиом линейного пространства рассматривалась исходя из свойств операций сложения и умножения в области действительных чисел.

В 1888 году Пеано на базе исчисления Грассмана впервые в явном виде сформулировал аксиомы линейного пространства (векторных пространств над полем действительных чисел, в том числе бесконечномерных) и применил обозначения, сохранившиеся в употреблении по сей день. Тёплиц в начале 1910-х годов обнаружил, что при помощи аксиоматизации линейного пространства для доказательства основных теорем линейной алгебры не требуется прибегать к понятию определителя, что позволяет распространить их результаты на случай бесконечного числа измерений. Аксиоматическое определение векторного и евклидова пространства было впервые чётко сформулировано в начале XX века практически одновременно Вейлем и фон Нейманом, исходя из запросов квантовой механики.

Дальнейшее развитие математики показало, что понятие линейного пространства можно применить и для описания множества ребер графа G и их свойств.

## 1.3. Основные положения теории графов

Точное определение графа состоит в том, что задаются два множества (первое из них обязательно непустое) и предикат, указывающий, какую пару элементов первого множества соединяет тот или иной элемент второго. Именно, дан граф G = (V,E;P), если даны два множества V, E и трехместный предикат $P \Leftrightarrow P_G \Leftrightarrow P(x,y,u)$, удовлетворяющий следующим двум условиям:

а) предикат P определен на всех таких упорядоченных тройках элементов $x, y, u$ для которых $x, y \in V$ и $u \in E$;

$$\text{б) } \forall u \exists x, y\{P(x,u,y) \ \& \ \forall x^*, y^*[P(x^*,u,y^*) \Rightarrow (x = x^* \ \& \ y = y^*) \vee (x = y^* \& \ y = x^*)]\} \qquad (1.9)$$

Элементы множества V называются вершинами, элементы множества E называются ребрами, а предикат P – инцидентор графа G; высказывание $P(x,u,y)$ читается так: ребро $u$ соединяет вершину $x$ с вершиной $y$, или $u$ соединяет упорядоченную пару вершин $\overline{xy}$. Условие (б) говорит о том, что каждое ребро графа соединяет какую-либо пару $\overline{xy}$ его вершин, но кроме этой пары может соединять еще только обратную пару $\overline{yx}$ [9,10].

Впредь будем пользоваться общепринятым обозначением графа G = (V,E) без указания трехместного предиката, так как обычно роль предиката выполняет матрица инциденций графа [7,9].

Любые две вершины $x,y \in V$ графа G = (V,E) называются *смежными*, если существует ребро $u$, соединяющее эти вершины, т.е. $u = (x,y)$.

Если ребро $u \in E$ графа G = (V,E) соединяет вершины $x, y \in V$ т. е. $u = (x,y)$, то говорят, что ребро $u$ инцидентно вершинам $x, y$ и наоборот, вершины $x, y$ инцидентны ребру $u$.

Все ребра с одинаковыми концевыми вершинами называются *кратными* или *параллельными*. Кроме того, концевые вершины ребра не обязательно различны. Если $u = (x,x)$, то ребро $u$ называется *петлей*. Граф называется *простым*, если он не содержит петель и кратных ребер [32].

Граф, не имеющий ребер, называется *пустым*. Граф, не имеющий вершин (и, следовательно, ребер) называется *нуль-графом* [37].

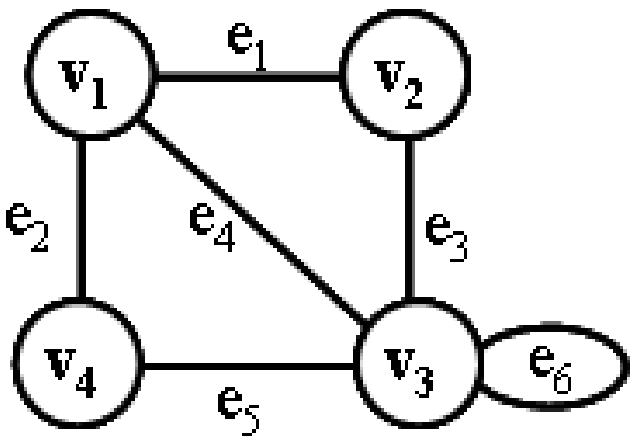

Рис. 1.6. Граф **G.**

*Локальной степенью* $d(v)$ вершины $v \in V$ называют число ребер инцидентных этой вершине. Иногда степень вершины называется также её *валентностью* [9,17].

Количество ребер графа обозначается буквой $m$, а количество вершин в графе буквой $n$.

Графически граф может быть представлен диаграммой, в которой вершина изображена точкой или кружком, а ребро – отрезком линии, соединяющим точки или кружки, соответствующие концевым вершинам графа. Например, если $V = \{v_1,v_2,v_3,v_4\}$ и $E = \{e_1,e_2,e_3,e_4,e_5,e_6\}$ такие, что $e_1 = (v_1,v_2)$, $e_2 = (v_1,v_4)$, $e_3 = (v_2,v_3)$, $e_4 = (v_1,v_3)$, $e_5 = (v_3,v_4)$, $e_6 =(v_3,v_3)$ и тогда этот граф G = (V,E) может быть представлен диаграммой (рис. 1.6).

## 1.4. Часть графа

Граф $G^* = (V^*,E^*)$ называется частью графа G = (V,E), если $V^* \subseteq V$ и $U^* \subseteq U$ т.е. часть графа образуется из исходного графа удалением некоторых вершин и ребер.

Особо важную роль играют два типа частей графа: подграф и суграф [10,34].

**Определение 1.3.** Часть $G^* = (V^*,E^*)$ называют *подграфом* графа G = (V,E), если $E^* = \{xy \in E / x,y \in V^*\}$. Подграф образуется из исходного графа некоторым количеством выделенных вершин и некоторым количеством инцидентных им ребер (рис. 1.7, б).

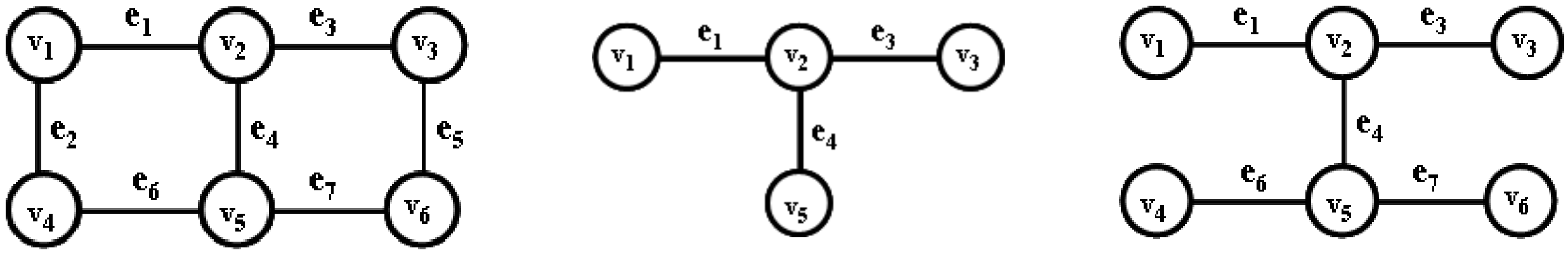

а) исходный граф б) подграф в) суграф

Рис 1.7. Получение из исходного графа частей графа.

**Определение 1.4.** Часть графа $G^* = (V^*,E^*)$ называют *суграфом* графа G = (X,U) если $X^* = X$, т.е. суграф образуется из исходного графа удалением только ребер, без удаления вершин (рис. 1.7,в).

## 1.5. Операции над суграфами

Введем несколько операций над суграфами графа G (рис. 1.8,а). Рассмотрим суграфы $G_1 = (V,E_1)$ и $G_2 = (V, E_2)$ (рис. 1.8,б и рис. 1.8,в).

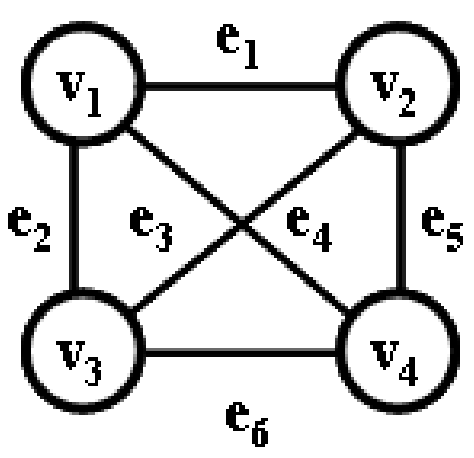

а) граф G.

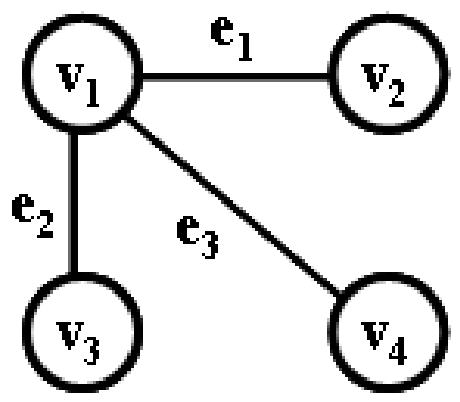

б) суграф $G_1$.

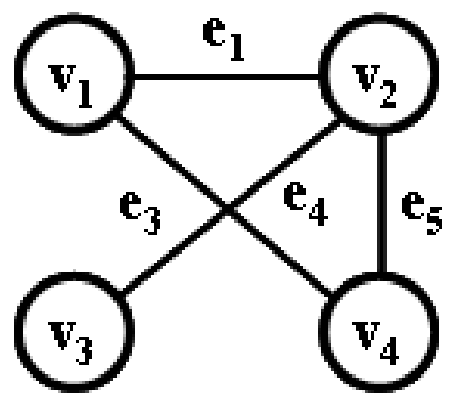

в) суграф $G_2$.

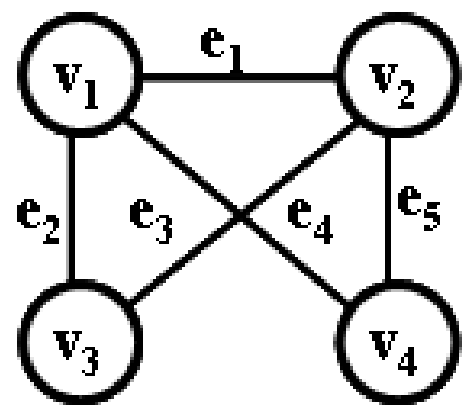

г) объединение суграфов $G_1$ и $G_2$.

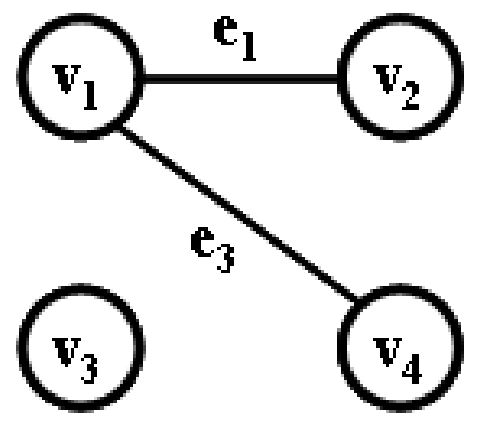

д) пересечение суграфов $G_1$ и $G_2$.

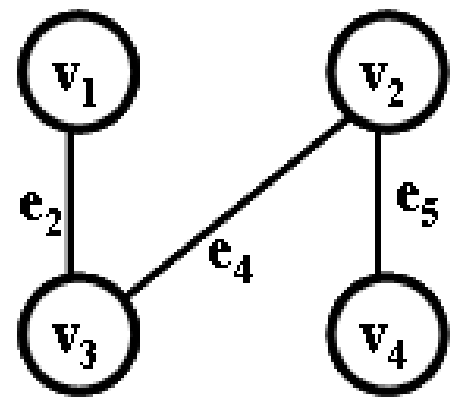

е) кольцевая сумма суграфов $G_1$ и $G_2$.

Рис. 1.8. Суграфы.

*Объединение* суграфов $G_1$ и $G_2$, обозначаемое как $G_1 \cup G_2$, представляет собой такой суграф $G_3 = (V, E_1 \cup E_2)$, что множество его ребер является объединением $E_1$ и $E_2$. Например, суграфы $G_1$ и $G_2$ и их объединение представлено на рис. 1.8,г.

*Пересечение* суграфов $G_1$ и $G_2$, обозначаемое как $G_1 \cap G_2$, представляет собой граф $G_3 = (V,E_1 \cap E_2)$. Таким образом, множество ребер $G_3$ состоит только из ребер, присутствующих одновременно в $G_1$ и $G_2$. Пересечение суграфов $G_1$ и $G_2$ показано на рис. 1.8,д.

*Кольцевая сумма* двух суграфов $G_1$ и $G_2$, обозначаемая как $G_1 \oplus G_2$, представляет собой суграф $G_3$, порожденный на множестве ребер $E_1 \oplus E_2$. Другими словами, суграф $G_3$ состоит только из ребер, присутствующих либо в $G_1$, либо в $G_2$, но не в обоих суграфах одновременно [34]. Кольцевая сумма суграфов показана на рис. 1.8,е.

Легко убедиться в том, что три рассмотренные операции коммутативны, т.е. $G_1 \cup G_2 = G_2 \cup G_1$, $G_1 \cap G_2 = G_2 \cap G_1$, $G_1 \oplus G_2 = G_2 \oplus G_1$.

Заметим также, что эти операции бинарны, т. е. определены по отношению к двум суграфам.

Таким образом, появляется операции сложения суграфов (будем называть ее кольцевой суммой [34]), операция отлична от известной арифметической операции сложения

$$(V,E_1;P) \oplus (V,E_2;P) = (V,(E_1 \cup E_2)\backslash(E_1 \cap E_2);P). \quad (1.10)$$

## 1.6. GF(2) – поле чисел по модулю 2

Рассматривая множество элементов векторного пространства, мы убеждаемся, что вещественные числа (а также комплексные или одни только рациональные числа) можно складывать и перемножать по известным правилам арифметики, получая при этом такие же числа. Это можно выразить словами: вещественные числа (а также комплексные, рациональные числа) образуют *поле* [3,34].

**Определение 1.5.** *Полем* называется множество F элементов, для которых определены две алгебраические операции – сложение и умножение (так что сумма **a** + **b** и произведение **ab** любых двух элементов **a** и **b** из F принадлежат F), причем выполняются следующие условия (аксиомы поля):

1. **a** + **b** = **b** + **a** для всех **a**, **b** из F (*сложение коммутативно*);

2. (**a** + **b**) + **c** = **a** + (**b** + **c**) для всех **a**, **b**, **c** из F (*сложение ассоциативно*);

3. в множестве F имеется нуль, т.е. такой элемент 0, что для каждого а из F сумма **a** + 0 = **a**;

4. для каждого **a** из F существует такой элемент **-a**, что **a** + (**-a**) = 0; такой элемент называется противоположный **a**;

5. **ab** = **ba** для всех **a**, **b** из F (*умножение коммутативно*);

6. (**ab**)**c** = **a**(**bc**) для всех **a**, **b** и **c** из F (*умножение ассоциативно*);

7. в множестве F имеется единица – такой элемент 1, что для всякого **a** из F имеем **a** • 1 = **a**;

8. для каждого отличного от нуля элемента **a** из F имеется такой (обратный **a**) элемент $\mathbf{a}^{-1}$, что $\mathbf{a} \bullet \mathbf{a}^{-1} = 1$;

9. (**a** + **b**)**c** = **ac** + **bc** для всех **a**, **b**, **c** из F (умножение дистрибутивно относительно сложения).

Ясно, что если коэффициенты системы принадлежат полю F, то и результат операций ее (если он существует) следует искать среди наборов из элементов поля F.

Особый интерес представляет поле GF(2) – множество целых чисел с операциями сложения и умножения по модулю 2. В этом поле: 0+0 = 0, 1+0 = 0+1 = 1, 1+1 = 0, $0\times 0 = 0$, $1\times 0 = 0\times 1 = 0$, $1\times 1 = 1$.

## 1.7. Алгебра суграфов

Пусть G = (V,E;P) – граф с пронумерованным множеством ребер $E = \{e_1, e_2, ..., e_m\}$, тогда суграф графа можно представлять в виде характеристического вектора с

коэффициентами

$$\alpha_{\text{i}} = \begin{cases} 1, \text{если } e_i \in E; \\ 0, \text{если } \text{e}_i \notin E. \end{cases}$$

Следующее множество элементов (1,0,....0), (0,1,.....,0),....,(0,0,....,1) представлявляет собой характеристические вектора для однореберных суграфов. Размерность характеристического вектора определяется количеством ребер графа.

В качестве примера рассмотрим граф G представленный на рис. 1.8(а). Характеристический вектор графа G = (1,1,1,1,1,1) можно записать в виде множества элементов $\{e_1,e_2,e_3,e_4,e_5,e_6\}$, рассматривая каждый элемент множества как результат произведения на соответствующий коэффициент из характеристического вектора. Последовательность записи коэффициентов в характеристическом векторе осуществляется справа налево.

Суграф $G_1$ представленный на рис. 1.8,б запишется в виде характеристического вектора (0,0,0,1,1,1) или в виде подмножества элементов $\{e_1,e_2,e_3\}$. Суграф $G_2$ представленный на рис. 1.8(в) запишется в виде характеристического вектора (0,1,1,1,0,1) или в виде подмножества элементов $\{e_1,e_3,e_4,e_5\}$. Суграф $G_1 \cup G_2$ представленный на рис. 1.8,г запишется в виде характеристического вектора (0,1,1,1,1,1) или в виде подмножества элементов $\{e_1,e_2,e_3,e_4,e_5\}$. Суграф $G_1 \cap G_2$ представленный на рис. 1.8(д) запишется в виде характеристического вектора (0,0,0,1,0,1) или в виде подмножества элементов $\{e_1,e_3\}$. Суграф $G_1 \oplus G_2$ представленный на рис. 1.8(е) запишется в виде характеристического вектора (0,1,1,0,1,0) или в виде подмножества элементов $\{e_2,e_4,e_5\}$.

$$G_1 \oplus G_2 = (0,0,0,1,1,1) + (0,1,1,1,0,1) = (0,1,1,0,1,0)$$

или

$$\{e_1,e_2,e_3) \oplus \{e_1,e_3,e_4,e_5\} = \{e_2,e_4,e_5\};$$

$$G_1 \oplus G_2 = (G_1 \cup G_2) \backslash (G_1 \cap G_2) = (0,1,1,1,1,1) + (0,0,0,1,0,1) = (0,1,1,0,1,0)$$

или

$$\{e_1,e_2,e_3,e_4,e_5) \oplus \{e_1,e_3\} = \{e_2,e_4,e_5\}.$$

## 1.8. Пространство суграфов графа

В дальнейшем множество суграфов графа G удобно рассматривать как линейное пространство $£_G$ над полем коэффициентов GF(2) = {0,1}, называемое пространством суграфов данного графа G. Размерность этого пространства dim $£_G$ = $m$ ибо множество элементов (1,0,...,0), (0,1,...,0),..., (0,0,...,1), представляющих однореберные суграфы, образует базис [4,34].

Рассмотрим множество $£_G$ всех $m$-векторов (суграфов) над полем GF(2). Если $\omega_1 = (\alpha_1, \alpha_2, ..., \alpha_m)$ и $\omega_2 = (\beta_1, \beta_2, ..., \beta_m)$ – элементы $£_G$, то $\omega_1 \oplus \omega_2 = (\alpha_1 \oplus \beta_2, \alpha_2 \oplus \beta_2, ..., \alpha_m \oplus \beta_m)$.

Если $\lambda$ принадлежит GF(2), то $\lambda\ \omega = (\lambda\ \alpha_1,\ \lambda\ \alpha_2,\ ...,\ \lambda\ \alpha_m)$.

Нетрудно установить, что множество $£_G$ удовлетворяет первой аксиоме в определении векторного пространства. Легко убедиться в том, что элементы множеств $£_G$ и GF(2) удовлетворяют и другим аксиомам этого определения.

Таким образом, $£_G$ векторное пространство над полем GF(2). Если E = {$e_1$,$e_2$,…,$e_m$}, то подмножества {$e_1$},{$e_2$},…,{$e_m$} образуют базис для $£_G$. Из того, что каждый реберно-порожденный суграф графа G соответствует единственному подмножеству множества E и того, что кольцевой сумме любых 2-реберно-порожденных суграфов можно поставить в соответствие кольцевую сумму двух соответствующих множеств ребер, следует, что множество всех реберно-порожденных суграфов графа G является векторным пространством. Заметим, что $£_G$ включает в себя нуль-граф $\varnothing$.

Приведем следующую теорему без доказательства.

**Теорема 1.1** [34]. Для графа G пространство $£_G$ является $m$-мерным векторным пространством над полем GF(2).

## 1.9. Подпространства циклов и разрезов

Следующие подмножества $£_G$ пространства суграфов являются подпространствами:

1. C(G) – множество всех циклов (включая и нуль-граф $\varnothing$) и объединений реберно-непересекающихся циклов (квазициклов) графа G.

2. S(G) – множество всех разрезающих множеств (включая и нуль-граф $\varnothing$) и объединений реберно-непересекающихся разрезающих множеств (квалиразрезов) графа G.

Если мы покажем, что C(G) и S**(**G**)** являются замкнутыми по отношению к операции сложения, то из этого и будет следовать, что C(G) и S(G) – подпространства.

**Теорема 1.2.** [34] Множество C(G) всех циклов и объединений реберно-непересекающихся циклов графа G является подпространством векторного пространства $£_G$ графа G.

*Доказательство.* Суграф можно представить в виде объединения реберно-непересекающихся циклов тогда и только тогда, когда каждая вершина в графе имеет четную степень. Следовательно, C(G) можно рассматривать как множество всех реберно-порожденных суграфов графа G, все вершины которого имеют четную степень.

Рассмотрим два любых элемента $c_1$ и $c_2$ множества C(G). Они являются реберно-порожденными суграфами, вершины которых имеют четную степень. Пусть $c_3$ – кольцевая сумма для элементов $c_1$ и $c_2$. Для доказательства необходимо показать, что сумма $c_3$ принадлежит множеству C(G). Иными словами, нужно показать, что каждая вершина в $c_3$ имеет четную степень.

Рассмотрим любую вершину *x*, принадлежащую $c_3$. Очевидно, эта вершина должна присутствовать, по крайней мере, в одном из суграфов: $c_1$ или $c_2$. Обозначим как $E_i$ множество ребер инцидентных вершине v в суграфе $c_i$, а $|E_i|$ – число ребер во множестве $E_i$. Таким образом, $|E_i|$ – степень вершины *x* в суграфе $c_i$. Заметим, что числа $|E_1|$ и $|E_2|$ – четные и могут даже равняться нулю. Поскольку $c_3 = c_1 \oplus c_2$, получаем, что $|E_3| = |E_1| + |E_2|$ - $2|E_1 \cap E_2|$. Из этого равенства видно, что $|E_3|$ является четной степенью, потому что обе степени $|E_1|$ и $|E_2|$ – четные. Другими словами, вершина *v* в $c_3$ имеет четную степень. А так как вершину *v* мы выбрали произвольно, то $c_3$ принадлежит C(G). *Теорема доказана.*

Покажем теперь, что множество S(G) всех разрезов (разрезающих множеств) и объединений реберно-непересекающихся разрезающих множеств графа G является подпространством пространства суграфов графа G.

Разрез является разрезающим множеством или объединением нескольких реберно-непересекающихся разрезающих множеств. Таким образом, каждый разрез графа G принадлежит множеству S(G). Докажем, что каждый элемент последнего является разрезом и одновременно покажем, что S(G) – подпространство пространства $L_G$.

**Теорема 1.3** [34]. Сумма двух разрезов графа G также является разрезом графа G.

*Доказательство*. Рассмотрим два любых разреза: $S_1 = \langle V_1, V_2\rangle$ и $S_2 = \langle V_3, V_4\rangle$ графа G = (V, E). Заметим, что $V_1 \cup V_2 = V_3 \cup V_4 = V$ и $V_1 \cap V_2 = V_3 \cap V_4 = \varnothing$. Пусть $A = V_1 \cap V_3$, $B = V_1 \cap V_4$, $C = V_2 \cap V_3$, $D = V_2 \cap V_4$.

Нетрудно видеть, что множества A, B, C, D взаимно не пересекаются. Тогда $\mathbf{S}_1 = \langle A \cup B, C \cup D\rangle = \langle A, C\rangle \cup \langle A, D\rangle \cup \langle B, C\rangle \cup \langle B, D\rangle$ и $\mathbf{S}_2 = \langle A \cup C, B \cup D\rangle = \langle A, B\rangle \cup \langle A, D\rangle \cup \langle C, B\rangle \cup \langle C, D\rangle$.

Следовательно, получим $\mathbf{S}_1 \oplus \mathbf{S}_2 = \langle A, C\rangle \cup \langle B, D\rangle \cup \langle A, B\rangle \cup \langle C, D\rangle$.

Поскольку $\langle A \cup D, B \cup C\rangle = \langle A, C\rangle \cup \langle B, D\rangle \cup \langle A, B\rangle \cup \langle C, D\rangle$, можно записать $\mathbf{S}_1 \oplus \mathbf{S}_2 = \langle A \cup D, B \cup C\rangle$.

Из того, что $A \cup D$ и $B \cup C$ взаимно не пересекаются и вместе включают в себя все вершины X, следует, что $S_1 \oplus S_2$ является разрезом в графе G. *Теорема доказана.*

Следующую теорему приведем без доказательства.

**Теорема 1.4.** [34]. Множество $S_G$ всех разрезов и объединений реберно-непересекающихся разрезов графа G является подпространством векторного пространства $£_G$ графа G.

## 1.10. Линейные операторы

*Оператором* называется правило [3], по которому каждому элементу $\boldsymbol{x}$ некоторого непустого множества **X** ставится в соответствие единственный элемент $\boldsymbol{y}$ некоторого непустого множества **Y**. Говорят, что оператор действует из **X** в **Y**.

Действие оператора обозначают $\boldsymbol{y} = \mathbf{A}(\boldsymbol{x})$, **y** – образ $\boldsymbol{x}$, **x** – прообраз $\boldsymbol{y}$.

Если каждый элемент $\boldsymbol{y}$ из **Y** имеет единственный прообраз $\boldsymbol{x}$ из $\boldsymbol{X}$, $\boldsymbol{y} = \mathbf{A}(\boldsymbol{x})$, оператор называют *взаимно однозначным отображением* **X** в **Y** или преобразованием **X**, **X** – область определения оператора.

Пусть **X** и **Y** два линейных пространства. Оператор **A**, действующий из **X** в **Y**, называется *линейным оператором*, если для любых двух элементов $\boldsymbol{u}$ и $\boldsymbol{v}$ из $\boldsymbol{X}$ и любого числа α справедливо: $\mathbf{A}(\boldsymbol{u} + \boldsymbol{v}) = \mathbf{A}(\boldsymbol{u}) + \mathbf{A}(\boldsymbol{v})$, $\mathbf{A}(\alpha\cdot\boldsymbol{u}) = \alpha\cdot\mathbf{A}(\boldsymbol{u})$.

Будем говорить, что задано отображение $\phi$ пространства X в пространство Y $\phi: X \to Y$, если каждому вектору $x$ из X поставлен однозначно в соответствие вектор $\phi(x)$ из Y.

## 1.11. Операторы в пространстве суграфов графа

Рассмотрим векторное пространство суграфов £(G) графа G. Известно, что множество суграфов графа G образует векторное пространство, называемое пространством суграфов графа £(G), где каждое ребро может быть представлено элементом базиса пространства £(G). В пространстве суграфов сложение двух суграфов $a$ и $b$ определяется операцией кольцевого суммирования, т.е. суммированием векторов по модулю 2 [34].

$$a \oplus b = (a \cup b)/(a \cap b) \qquad (1.11)$$

Таким образом, множество всех суграфов графа G может быть образовано кольцевым сложением из $m$-векторов векторного пространства £(G) над полем GF(2).

Рассмотрим сказанное на примере графа $G_1$ представленного на рис. 1.9.

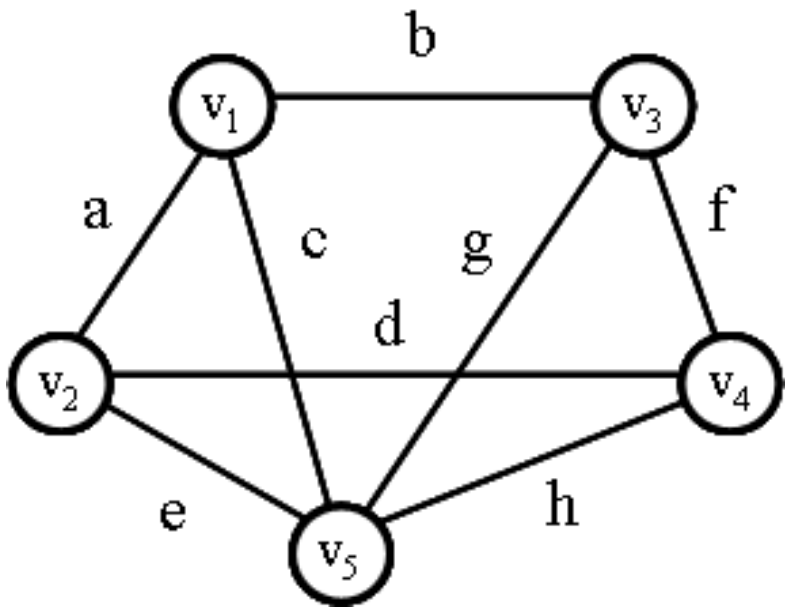

Рис. 1.9. Граф $G_1$.

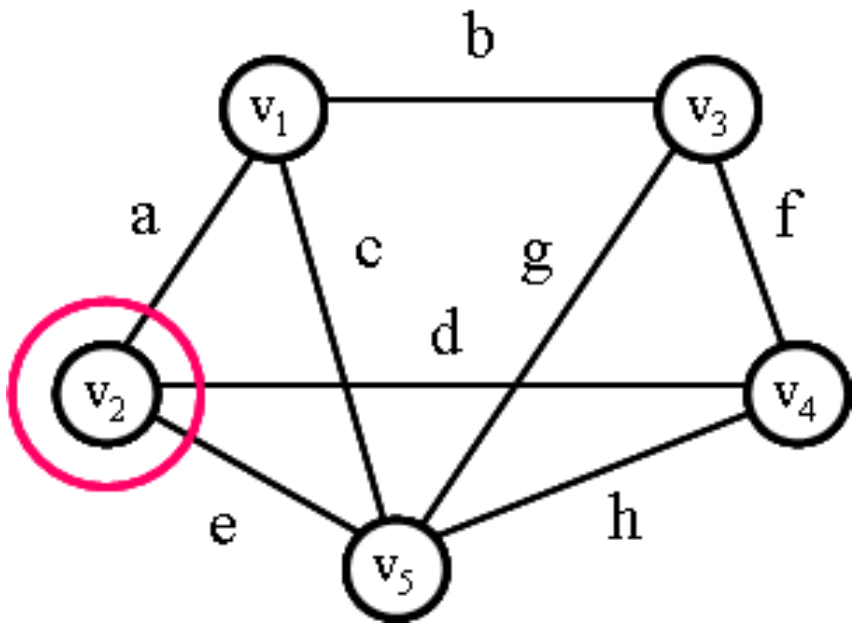

Рис. 1.10. Центральный разрез графа $G_1$.

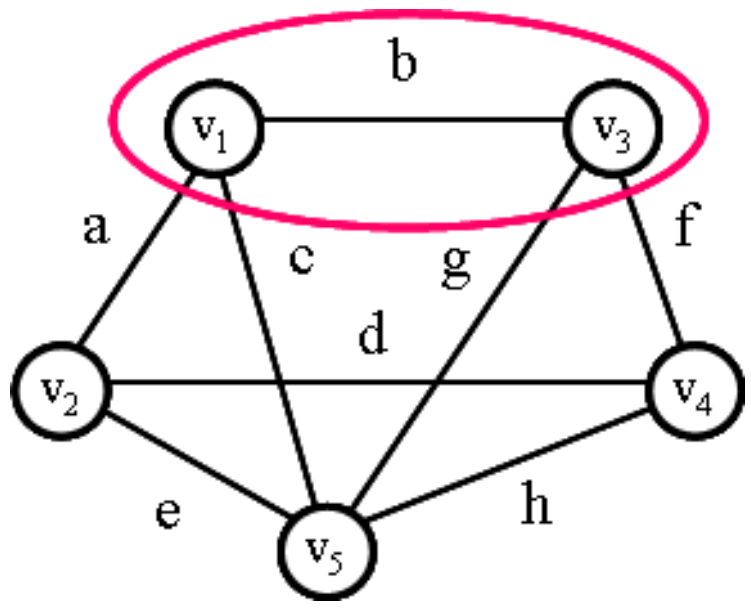

Рис. 1.11. Базовый реберный разрез графа $G_1$.

Представим ребра графа $G_1$ как восьмимерные вектора:

$a = (0,0,0,0,0,0,0,1)$; $b = (0,0,0,0,0,0,1,0)$; $c = (0,0,0,0,0,1,0,0)$;
$d = (0,0,0,0,1,0,0,0)$; $e = (0,0,0,1,0,0,0,0)$; $f = (0,0,1,0,0,0,0,0)$;
$g = (0,1,0,0,0,0,0,0)$; $h = (1,0,0,0,0,0,0,0)$.

Тогда суграф состоящий из ребра *a*, ребра *d* и ребра *g* имеет вид $\{a,d,g\}$= (1,0,0,1,0,0,1,0).

Матрица инциденций графа $G_1$ имеет вид:

B =

| | a | b | c | d | e | f | g | h |
|---|---|---|---|---|---|---|---|---|
| $v_1$ | 1 | 1 | 1 | | | | | |
| $v_2$ | 1 | | | 1 | 1 | | | |
| $v_3$ | | 1 | | | | 1 | 1 | |
| $v_4$ | | | | 1 | | 1 | | 1 |
| $v_5$ | | | 1 | | 1 | | 1 | 1 |

Рассмотрим матрицу B [7,33] как линейный оператор в пространстве суграфов £(G) [20]. Оператор B осуществляет сюръективное отображение ребер графа из пространств суграфов £($G_1$) во множество центральных разрезов вершин из подпространства разрезов. Обозначим вектор-столбец ребер как $\beta$($G_1$) = (a,b,c,d,e,f,g,h), тогда вектор-столбец S($G_1$) центральных разрезов можно определить с помощью выражения

$$S^T(G) = B \times \beta^T(G). \qquad (1.12)$$

S(G) =

| |
|---|
| $s_1 = \{a,b,c\}$ |
| $s_2 = \{a,d,e\}$ |
| $s_3 = \{b,f,g\}$ |
| $s_4 = \{d,f,h\}$ |
| $s_5 = \{c,t,g,h\}$ |

=

| a | b | c | d | e | f | g | h |
|---|---|---|---|---|---|---|---|
| 1 | 1 | 1 | | | | | |
| 1 | | | 1 | 1 | | | |
| | 1 | | | | 1 | 1 | |
| | | | 1 | | 1 | | 1 |
| | | 1 | | 1 | | 1 | 1 |

×

| |
|---|
| *a* |
| *b* |
| *c* |
| *d* |
| *e* |
| *f* |
| *g* |
| *h* |

Для графа $G_1$ множество центральных разрезов графа S состоит из следующих элементов:

$S^T(G_1) = \{s_1 = \{a,b,c\}, s_2 = \{a,d,e\}, s_3 = \{b,f,g\}, s_4 = \{d,f,h\}, s_5 = \{c,e,g,h\}\}$.

Транспонированная матрица $B^T$ также является оператором в пространстве суграфов £(G). Оператор $B^T$ осуществляет сюръективное отображение центральных разрезов вершин из подпространства разрезов S во множество реберных разрезов $W_0$ пространства суграфов £(G). Вектор-столбец $W_0(G)$ базовых реберных разрезов графа можно определить с помощью выражения

$$W_0^T(G) = B^T \times S^T(G). \quad (1.13)$$

Для графа $G_1$ множество базовых реберных разрезов графа $W_0(G)$ состоит из следующих элементов образованных кольцевым суммированием соответствующих центральных разрезов:

$W_0(G) =$

| |
|---|
| $w_0(a) = \{b,c,e,d\}$ |
| $w_0(b) = \{a,c,f,g\}$ |
| $w_0(c) = \{a,b,e,g,h\}$ |
| $w_0(d) = \{a,e,f,h\}$ |
| $w_0(e) = \{a,d,c,g,h\}$ |
| $w_0(f) = \{b,d,g,h\}$ |
| $w_0(g) = \{b,c,e,f,h\}$ |
| $w_0(h) = \{d,e,f,g,h\}$ |

=

| | $s_1$ | $s_2$ | $s_3$ | $s_4$ | $s_5$ |
|---|---|---|---|---|---|
| a | 1 | 1 | | | |
| b | 1 | | 1 | | |
| c | 1 | | | | 1 |
| d | | 1 | | 1 | |
| e | | 1 | | | 1 |
| f | | | 1 | 1 | |
| g | | | 1 | | 1 |
| h | | | | 1 | 1 |

×

| |
|---|
| $s_1 = \{a,b,c\}$ |
| $s_2 = \{a,d,e\}$ |
| $s_3 = \{b,f,g\}$ |
| $s_4 = \{d,f,h\}$ |
| $s_5 = \{c,t,g,h\}$ |

$w_0(G_1) = \{w_0(a) = s_1 \oplus s_2 = \{b,c,e,d\}, w_0(b) = s_1 \oplus s_3 = \{a,c,f,g), w_0(c) = s_1 \oplus s_5 = \{a,b,e,g,h\},$
$w_0(d) = s_2 \oplus s_4 = \{a,e,f,h\}, w_0(e) = s_2 \oplus s_5 = \{a,d,c,g,h\}, w_0(f) = s_3 \oplus s_4 = \{b,d,g,h\},$
$w_0(g) = s_3 \oplus s_5 = \{b,c,e,f,h\}, w_0(h) = s_4 \oplus s_5 = \{d,e,f,g,h\}\}$.

Базовый реберный разрез ребра $w_0(e_i)$ можно представлять как функцию, определив в качестве аргумента ребро $e_i$ графа G.

На рис. 1.10 показан центральный разрез вершины $v_2$. На рис. 1.11 показан базовый реберный разрез ребра *b* графа $G_1$. На рис. 1.12 представлено суръективное отображение множества ребер графа $\beta(G_1)$ во множество центральных разрезов вершин $S(G_1)$. На рис. 1.13 представлено отображение множества центральных разрезов вершин $S(G_1)$ графа $G_1$ во множество базовых реберных разрезов $W_0(G_1)$.

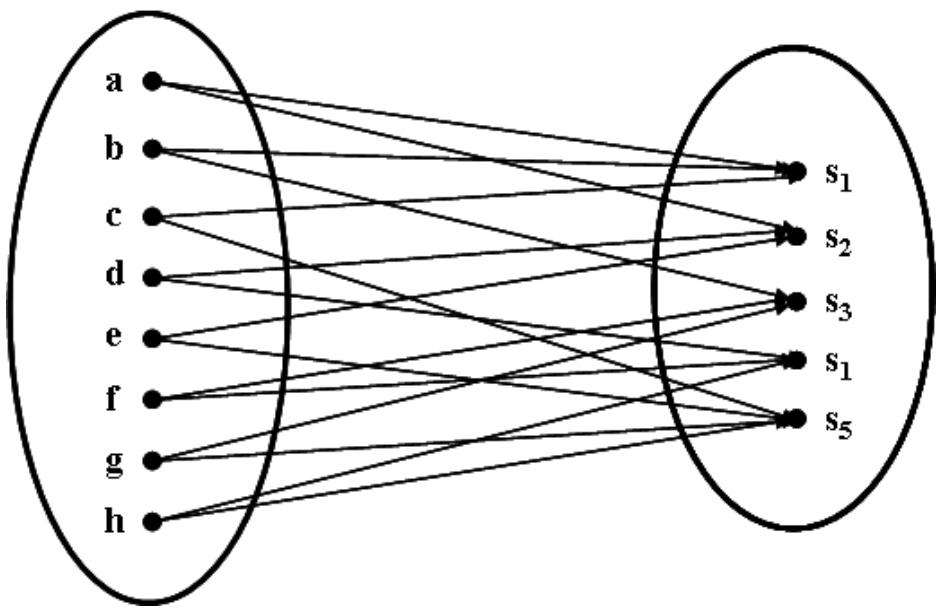

Рис. 1.12. Отображение множества ребер во множество центральных разрезов вершин.

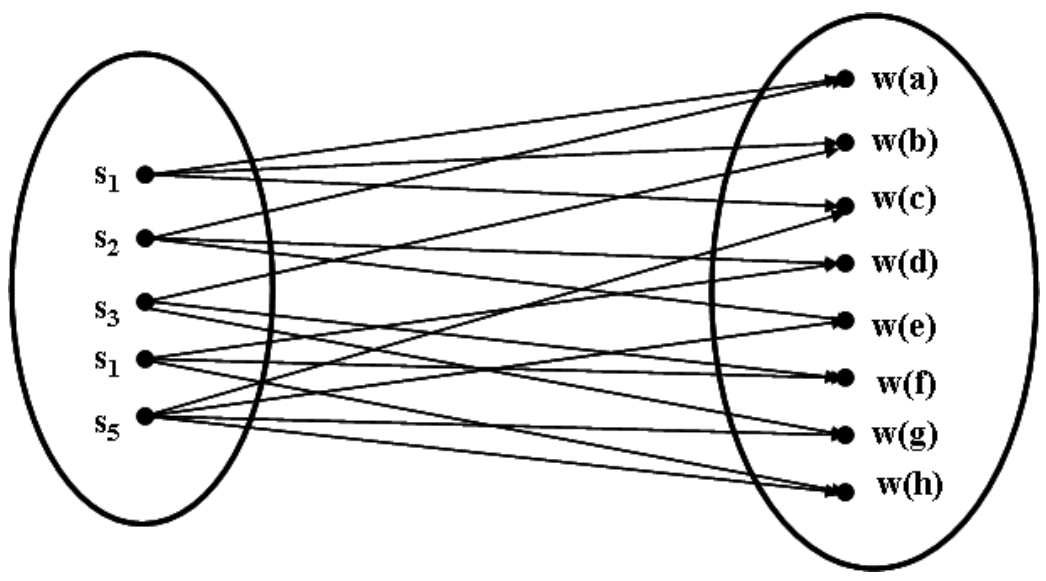

Рис. 1.13. Отображение множества центральных разрезов во множество реберных разрезов.

Известно, что *матрица смежностей реберного графа* A(L(G)) графа G в поле вещественных чисел может быть вычислена по формуле [36]:

$$A(L(G)) = B(G)^T \times B(G) - 2 \times I \tag{1.14}$$

Очевидно, что в поле чисел GF(2) *матрица смежностей реберного графа* A(L(G)) графа G равна

$$A(L(G)) = B(G)^T \times B(G)\,. \tag{1.15}$$

Подставляя последнее в формулу (1.13) с учетом (1.12), получим

$$W^T(G) = B^T \times S^T(G) = B^T \times B \times \beta^T(G) = A(L(G)) \times \beta^T(G)\,. \tag{1.16}$$

Отсюда следует, что матрицу смежностей реберного графа A(L(G)) можно рассматривать как оператор в пространстве суграфов £(G). Обозначим этот оператор как $W$ .

$W(G_1) =$

| | a | b | c | d | e | f | g | h |
|---|---|---|---|---|---|---|---|---|
| *a* | | 1 | 1 | 1 | 1 | | | |
| *b* | 1 | | 1 | | | 1 | 1 | |
| *c* | 1 | 1 | | | 1 | | 1 | 1 |
| *d* | 1 | | | | 1 | 1 | | 1 |
| *e* | 1 | | 1 | 1 | | | 1 | 1 |
| *f* | | 1 | | 1 | | | 1 | 1 |
| *g* | | 1 | 1 | | 1 | 1 | | 1 |
| *h* | | | 1 | 1 | 1 | 1 | 1 | |

Наличие такого оператора в пространстве суграфов графа G позволяет осуществлять преобразование векторного пространства суграфов само в себя. То есть имеет смысл рассматривать операторы $W^2 = W \times W$ , $W^3 = W \times W \times W$ ,…, $W^p = W \times W \times ... \times W$ ,...

Введем понятие реберного разреза графа.

**Определение 1.6.** *Базовым реберным разрезом* будем называть суграф (квалиразрез), образованный кольцевым суммированием соответствующих концевых центральных разрезов и состоящий из инцидентных ребер принадлежащих двум концевым вершинам данного ребра, за исключением самого ребра. Обозначим его как $w(e_i)$ для ребра $e_i$ соединяющего вершины $v_{j1}$ и $v_{j2}$ (рис. 1.11) [20].

Все множество ребер графа запишем в виде вектора–столбца

$$\beta(G)^T = \{e_1, e_2, ..., e_m\}^T\,. \tag{1.17}$$

$W^2(G_1) =$

| | a | b | c | d | e | f | g | h |
|---|---|---|---|---|---|---|---|---|
| *a* | | 1 | 1 | 1 | 1 | | | |
| *b* | 1 | | 1 | | | 1 | 1 | |
| *c* | 1 | 1 | | | 1 | | 1 | 1 |
| *d* | 1 | | | | 1 | 1 | | 1 |
| *e* | 1 | | 1 | 1 | | | 1 | 1 |
| *f* | | 1 | | 1 | | | 1 | 1 |
| *g* | | 1 | 1 | | 1 | 1 | | 1 |
| *h* | | | 1 | 1 | 1 | 1 | 1 | |

Вычислим операторы для нашего примера. Определим оператор $W^2$ как результат умножения двух матриц A(L(G)), применяя закон сложения и умножения элементов по mod 2.

$W^3(G_1) =$

| | *a* | *b* | *c* | *d* | *e* | *f* | *g* | *h* |
|---|---|---|---|---|---|---|---|---|
| *a* | | 1 | 1 | 1 | 1 | | | |
| *b* | 1 | | 1 | | | 1 | 1 | |
| *c* | 1 | 1 | | | 1 | | 1 | 1 |
| *d* | 1 | | | | 1 | 1 | | 1 |
| *e* | 1 | | 1 | 1 | | | 1 | 1 |
| *f* | | 1 | | 1 | | | 1 | 1 |
| *g* | | 1 | 1 | | 1 | 1 | | 1 |
| *h* | | | 1 | 1 | 1 | 1 | 1 | |

$= 0$

Закончим построение цепочки операторов в случае, когда имеется такой показатель степени оператора со значением $q$, при котором все элементы линейного оператора $W^q$ либо равны нулю, либо равны элементам базового оператора $W$. То есть, по условию остановки процесса порождения реберных разрезов все элементы оператора $W^q$ со значением степени $q$ и выше – пустое множество. В нашем примере оператор третьей степени совпадает с базовым оператором $W(G_1)$ и поэтому должен быть равен $\varnothing$.

Исходя из цепочки операторов формируются множества реберных разрезов различного уровня, которые определяются по формулам:

$$W_0(G_1)^T = W\,(G_1)\times\beta(G_1)^T; \quad (1.18)$$

$$W_1(G_1)^T = W(G_1)\times W_0(G_1)^T = W^2(G_1)\times\beta(G_1)^T = W(G_1)\times W(G_1)\times\beta(G_1)^T; \quad (1.19)$$

$$W_2(G_1)^T = W(G_1)\times W_1(G_1)^T = W^2(G_1)\times W_1(G_1)^T = W(G_1)\times W^2(G_1)\times\beta(G_1)^T; \quad (1.20)$$

и т.д.

## Комментарии

Матрицу инциденций графа можно рассматривать как линейный оператор, действующий в пространстве суграфов графа G. Данный оператор переводит вектор-столбец ребер во множество суграфов характеризующих центральные разрезы вершин графа. Транспонированную матрицу инциденций также можно рассматривать как линейный оператор, отображающий множество центральных разрезов графа во множество базовых реберных разрезов графа W(G). В свою очередь, матрица базовых реберных разрезов W(G) по сути, является матрицей смежностей A(L(G)) реберного графа и также может рассматриваться как линейный оператор, действующий в пространстве суграфов. Линейный оператор A(L(G)) является оператором W в поле GF(2). Как мы увидим в дальнейшем, линейный оператор A(L(G)) порождает конечную цепочку реберных разрезов графа.

Так как мв будем уметь дело только с несепарабельными графами, то удобным способом представления графа является множество суграфов. Множечтво суграфов и операция кольцевого суммирования суграфов порождает пространство суграфов. В свою очередь, пространство суграфов состоит из подпространства разрезов и подпространства циклов графа.

# Глава 2. СПЕКТР РЕБЕРНЫХ РАЗРЕЗОВ ГРАФА

## 2.1. Реберный разрез графа и его свойства

Рассмотрим свойства разрезов графа G.

**Определение 2.1.** *Центральным разрезом* называется суграф графа G (квалиразрез), состоящий из инцидентных ребер принадлежащих данной вершине. Обозначим его как $s(v_i)$ для вершины *i* (рис. 2.1) [10].

**Определение 2.2.** *Базовым реберным разрезом* называется суграф графа G (квалиразрез), образованный из инцидентных ребер принадлежащих двум концевым вершинам данного ребра, за исключением самого ребра. Обозначим его как $w(e_i)$ для ребра $e_i$ соединяющего вершины $v_{j1}$ и $v_{j2}$ (рис. 2.2).

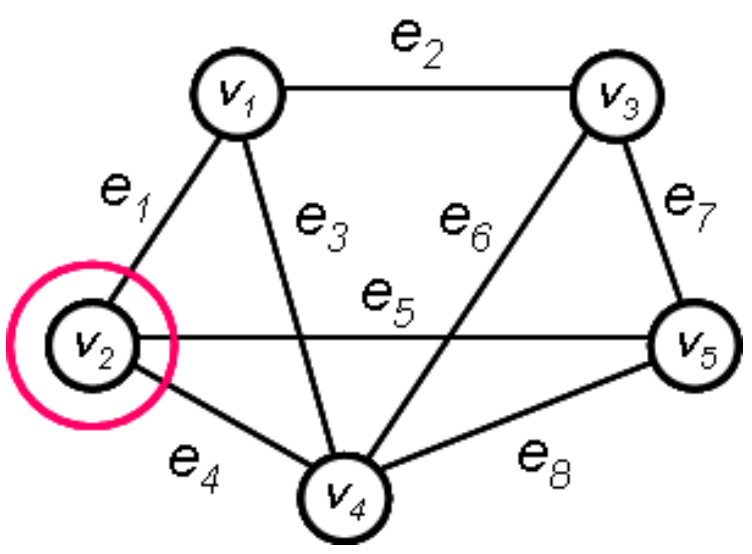

Рис. 2.1 Центральный разрез $\{e_1,e_6,e_7\}$.

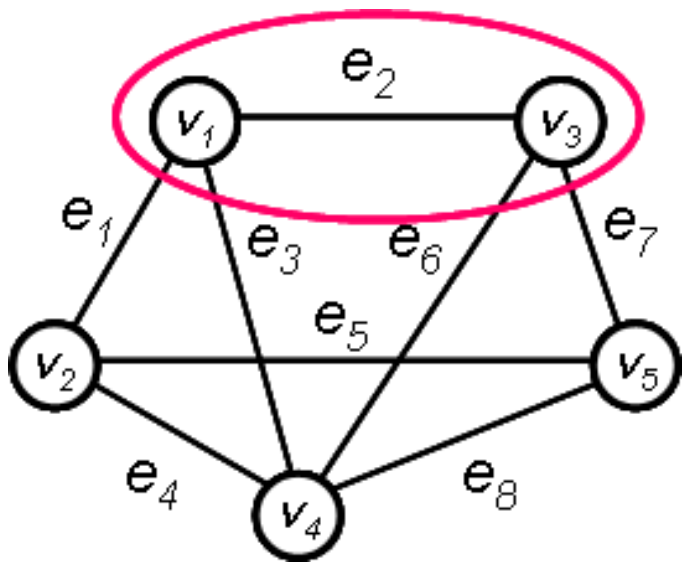

Рис. 2.2. Реберный разрез $\{e_1,e_3,e_4,e_5\}$.

**Определение 2.3.** Количество ребер в квалиразрезе (суграфе) назовём *длиной разреза*. Множество всех базовых реберных разрезов графа G обозначим как $W_0(G)$.

Между центральными разрезами и квазициклами графа существует связь, устанавливаемая следующими леммами:

**Лемма 2.1.** Кольцевая сумма множества центральных разрезов графа есть пустое множество.

*Доказательство.* Этот факт следует из линейной независимости любых *n-1* вершин [4,34].

**Лемма 2.2.** Кольцевая сумма реберных разрезов для любого квазицикла графа есть пустое множество.

*Доказательство.* Пусть длина квазицикла равна *p* и реберный разрез равен кольцевой сумме центральных разрезов для концевых вершин ребра $w(e_i) = s(v_{j1}) \oplus s(v_{j2})$. Тогда

$$w(e_1) \oplus w(e_2) \oplus ... \oplus w(e_p) = (s(v_1) \oplus s(v_2)) \oplus (s(v_2) \oplus s(v_3)) \oplus ... \oplus (s(v_p) \oplus s(v_1)) = \varnothing . \quad (2.1)$$

Следовательно, кольцевая сумма есть пустое множество, так как центральные разрезы присутствуют в правой части выражения (2.1) дважды (рис. 2.3).

**Лемма 2.3.** Кольцевая сумма пересечения центральных разрезов для ребер любого квазицикла графа есть квазицикл.

*Доказательство.* Любое ребро графа можно представить в виде $e_i = \mathrm{s}(v_{j1}) \cap \mathrm{s}(v_{j2})$, а любой квазицикл длиной *p* можно представить как кольцевую сумму ребер. Тогда

$$c = e_1 \oplus e_2 \oplus ... \oplus e_p = (\mathrm{s}(v_1) \cap \mathrm{s}(v_2)) \oplus (\mathrm{s}(v_2) \cap \mathrm{s}(v_3) \oplus ... \oplus (\mathrm{s}(v_p) \cap \mathrm{s}(v_1)) . \quad (2.2)$$

Отсюда следует: $$\mathrm{c} = \sum_{k=1}^{p} \mathrm{s}(v_i) \cap \mathrm{s}(v_{i+1}) . \quad (2.3)$$

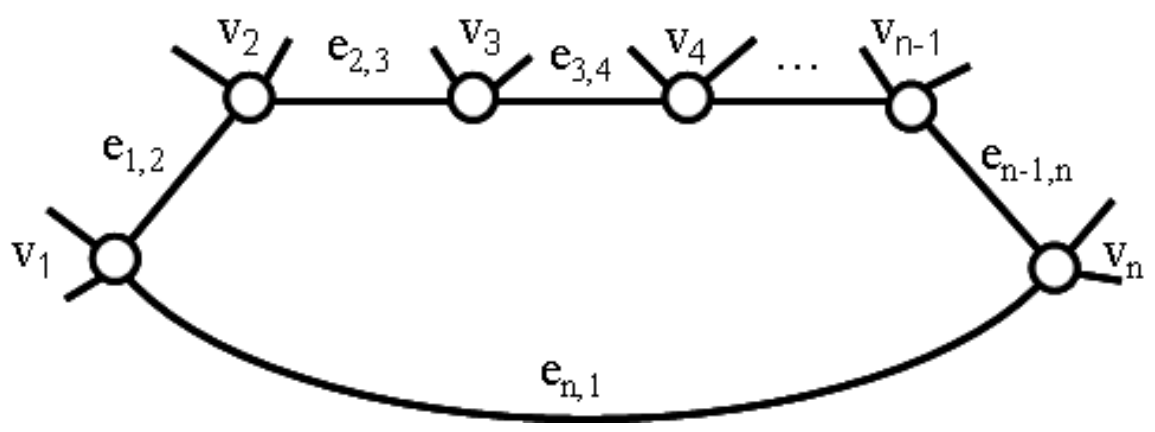

Рис. 2.3. Цикл графа

Для иллюстрации лемм рассмотрим граф $\mathbf{G_2}$ представленный на рис. 2.4.

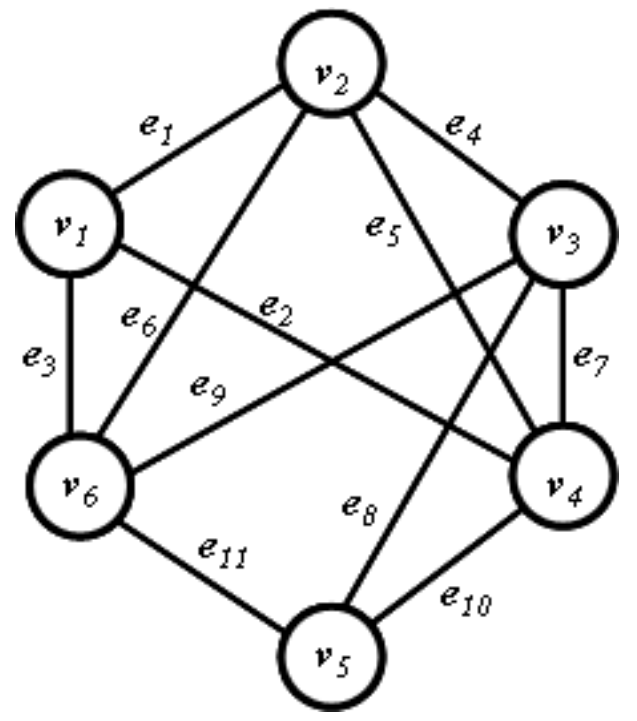

Рис. 2.4. Граф $G_2$.

Множество центральных разрезов графа $\mathbf{G_2}$ определяется выражением

$$S = \{s_1,s_2,s_3,s_4,s_5,s_6\},$$

где:

$s_1 = \{e_1,e_2,e_3\}$; $s_2 = \{e_1,e_4,e_5,e_6\}$; $s_3 = \{e_4,e_7,e_8,e_9\}$;
$s_4 = \{e_2,e_5,e_7,e_{10}\}$; $s_5 = \{e_8,e_{10},e_{11}\}$; $s_6 = \{e_3,e_6,e_9,e_{11}\}$.

Кольцевая сумма центральных разрезов графа $\mathbf{G_2}$ равна:

$s_1 \oplus s_2 \oplus s_3 \oplus s_4 \oplus s_5 \oplus s_6 =$

$=\{e_1,e_2,e_3\} \oplus \{e_1,e_4,e_5,e_6\} \oplus \{e_4,e_7,e_8,e_9\} \oplus \{e_2,e_5,e_7,e_{10}\} \oplus \{e_8,e_{10},e_{11}\} \oplus \{e_3,e_6,e_9,e_{11}\} = \varnothing$.

Выберем цикл $\{e_4,e_6,e_9\}$. Кольцевая сумма реберных разрезов равна:

$w(e_4) \oplus w(e_6) \oplus w(e_9) = \{e_1,e_5,e_6,e_7,e_8,e_9\} \oplus \{e_1,e_3,e_4,e_5,e_9,e_{11}\} \oplus \{e_3,e_4,e_6,e_7,e_8,e_{11}\} = \varnothing$.

Цикл представляется как сумма пересечений центральных разрезов:

$c = (\{e_1,e_4,e_5,e_6\} \cap \{e_4,e_7,e_8,e_9\}) \oplus (\{e_1,e_4,e_5,e_6\} \cap \{e_3,e_6,e_9,e_{11}\}) \oplus (\{e_3,e_6,e_9,e_{11}\} \cap$
$\cap \{e_4,e_7,e_8,e_9\}) = \{e_4\} \oplus \{e_6\} \oplus \{e_9\} = \{e_4,e_6,e_9\}$.

**Лемма 2.4.** Кольцевая сумма двух реберных разрезов порождается четным числом центральных разрезов.

*Доказательство.* Так как реберный разрез определяется двумя центральными разрезами, то количество центральных разрезов $q$ порожденного реберного разреза можно определить по формуле

$$q = 2 \times q_1 - 2 \times q_2, \qquad (2.4)$$

где $q_1$ – количество центральных разрезов в первом реберном разрезе, $q_2$ – количество центральных разрезов во втором реберном разрезе. Так как первоначально реберный разрез состоит из двух центральных разрезов и сложение производится по законам кольцевого суммирования, то количество центральных разрезов в сумме всегда четно. *Лемма доказана.*

В качестве примера сформируем множество базовых реберных разрезов графа $G_2$:

$w_0(e_1) = s_1 \oplus s_2 = \{e_1,e_2,e_3\} \oplus \{e_1,e_4,e_5,e_6\} = \{e_2,e_3,e_4,e_5,e_6\}$;
$w_0(e_2) = s_1 \oplus s_4 = \{e_1,e_2,e_3\} \oplus \{e_2,e_5,e_7,e_{10}\} = \{e_1,e_3,e_5,e_7,e_{10}\}$;
$w_0(e_3) = s_1 \oplus s_6 = \{e_1,e_2,e_3\} \oplus \{e_3,e_6,e_9,e_{11}\} = \{e_1,e_2,e_6,e_9,e_{11}\}$;
$w_0(e_4) = s_2 \oplus s_3 = \{e_1,e_4,e_5,e_6\} \oplus \{e_4,e_7,e_8,e_9\} = \{e_1,e_5,e_6,e_7,e_8,e_9\}$;
$w_0(e_5) = s_2 \oplus s_4 = \{e_1,e_4,e_5,e_6\} \oplus \{e_2,e_5,e_7,e_{10}\} = \{e_1,e_2,e_4,e_6,e_7,e_{10}\}$;
$w_0(e_6) = s_2 \oplus s_6 = \{e_1,e_4,e_5,e_6\} \oplus \{e_3,e_6,e_9,e_{11}\} = \{e_1,e_3,e_4,e_5,e_9,e_{11}\}$;
$w_0(e_7) = s_3 \oplus s_4 = \{e_4,e_7,e_8,e_9\} \oplus \{e_2,e_5,e_7,e_{10}\} = \{e_2,e_4,e_5,e_8,e_9,e_{10}\}$;
$w_0(e_8) = s_3 \oplus s_5 = \{e_4,e_7,e_8,e_9\} \oplus \{e_8,e_{10},e_{11}\} = \{e_4,e_7,e_9,e_{10},e_{11}\}$;
$w_0(e_9) = s_3 \oplus s_6 = \{e_4,e_7,e_8,e_9\} \oplus \{e_3,e_6,e_9,e_{11}\} = \{e_3,e_4,e_6,e_7,e_8,e_{11}\}$;
$w_0(e_{10}) = s_4 \oplus s_5 = \{e_2,e_5,e_7,e_{10}\} \oplus \{e_8,e_{10},e_{11}\} = \{e_2,e_5,e_7,e_8,e_{11}\}$;
$w_0(e_{11}) = s_5 \oplus s_6 = \{e_8,e_{10},e_{11}\} \oplus \{e_3,e_6,e_9,e_{11}\} = \{e_3,e_6,e_8,e_9,e_{10}\}$.

Рассмотрим основные свойства реберного разреза [20,24].

**Свойство 2.1.** Если реберный разрез представляет собой квазицикл, то согласно лемме 2.2 такой реберный разрез есть пустое множество.

Например, для графа $G_2$ реберный разрез $w(e_1) \oplus w(e_8) = s_1 \oplus s_2 \oplus s_3 \oplus s_5 = \{e_2,e_3,e_5,e_6,e_7,e_9,e_{10},e_{11}\}$ есть пустое множество.

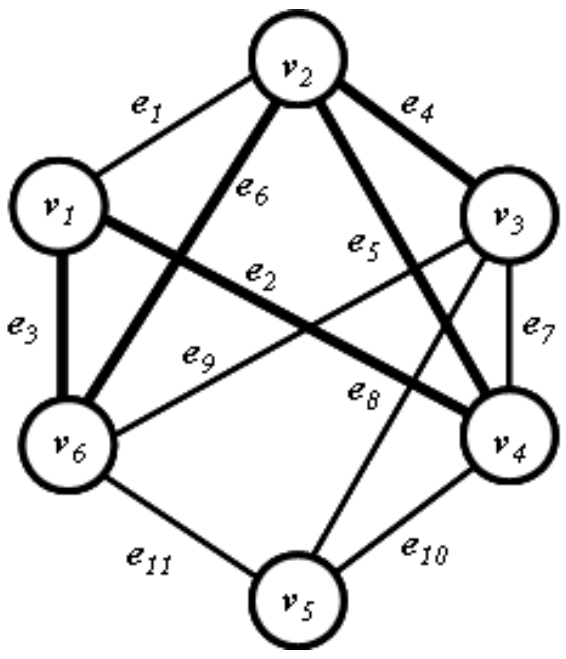

Рис. 2.5. Суграф $\{e_2,e_3,e_4,e_5,e_6\}$

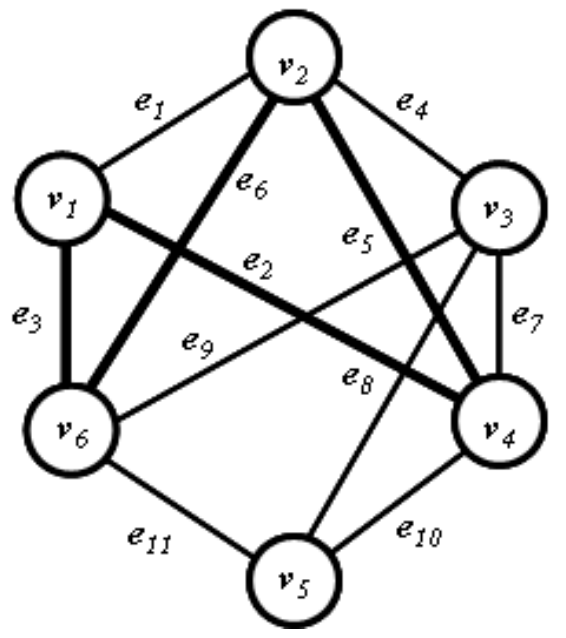

Рис. 2.6. Суграф с четной валентностью вершин.

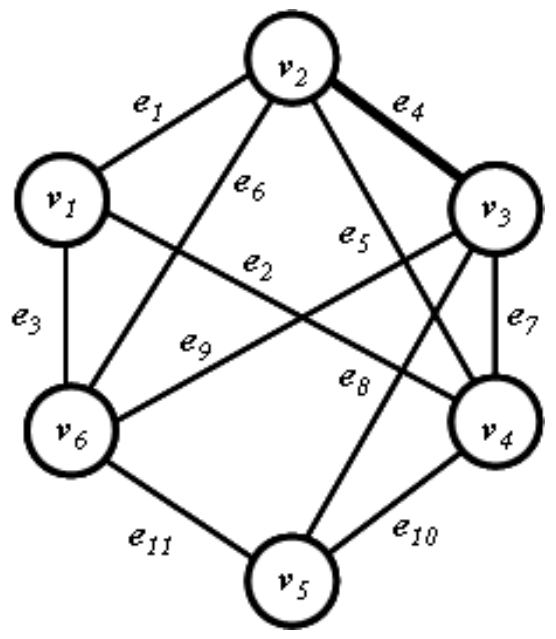

Рис. 2.7. Суграф с нечетной валентностью вершин.

**Свойство 2.2.** Реберный разрез может быть представлен двумя суграфами, один из которых имеет четную (рис. 2.6), а другой – нечетную валентность вершин (рис. 2.7).

## 2.2. Базовый реберный разрез графа и матрица совместимостей графа

Как правило, для представления графа используется матрица инциденций B(G). Матрица смежностей графа A(G) в поле действительных чисел может быть вычислена по формуле:

$$A(G) = B(G) \times B(G)^{T} - \rho(G) \times I, \qquad (2.5)$$

где $\rho$ (G) – вектор локальных степеней графа, I – единичная матрица [30,32].

Если вместо единиц матрицы A(G) в соответствующих ячейках записать идентификаторы ребер согласно трехместному предикату графа, то получим матрицу совместимостей графа AB(G). Сказанное рассмотрим на примере графа представленного на рис. 2.8.

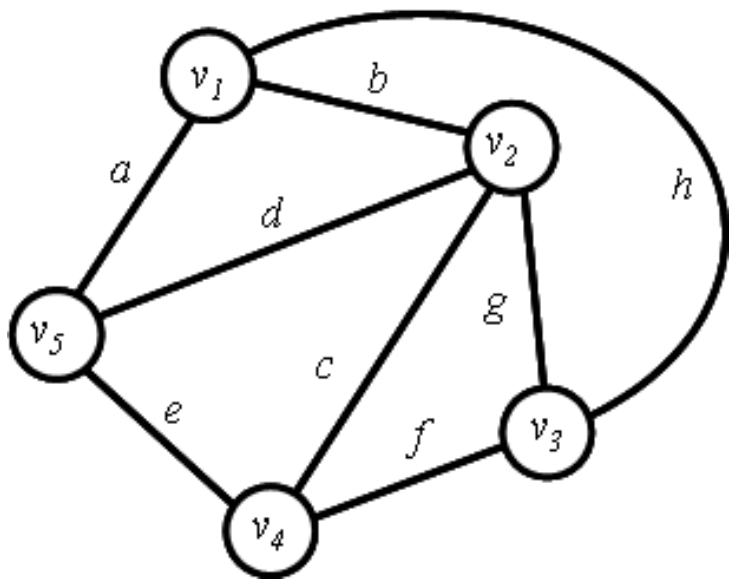

Рис. 2.8. Граф $G_3$.

Матрицы A(G) и AB(G) для графа $G_3$:

A(G) =

| | $v_1$ | $v_2$ | $v_3$ | $v_4$ | $v_5$ |
|---|---|---|---|---|---|
| $v_1$ | | 1 | 1 | | 1 |
| $v_2$ | 1 | | 1 | 1 | 1 |
| $v_3$ | 1 | 1 | | 1 | |
| $v_4$ | | 1 | 1 | | 1 |
| $v_5$ | 1 | 1 | | 1 | |

AB(G)=

| | $v_1$ | $v_2$ | $v_3$ | $v_4$ | $v_5$ |
|---|---|---|---|---|---|
| $v_1$ | | *b* | *h* | | *a* |
| $v_2$ | *b* | | *g* | *c* | *d* |
| $v_3$ | *h* | *g* | | *f* | |
| $v_4$ | | *c* | *f* | | *e* |
| $v_5$ | *a* | *d* | | *e* | |

Из матрицы AB(G) можно сформировать и матрицу смежностей графа A(G), и матрицу инциденций графа B(G). Кроме того, рассматривая матрицу AB(G) можно записать предикат P, определяющий местоположение ребра относительно вершин:

$P(G_2) = a \rightarrow [(v_1,v_5) \vee (v_5,v_1)];$
$b \rightarrow [(v_1,v_2) \vee (v_2,v_1)];$
$c \rightarrow [(v_2,v_4) \vee (v_5,v_2)];$
$d \rightarrow [(v_2,v_5) \vee (v_5,v_2)];$
$e \rightarrow [(v_1,v_5) \vee (v_5,v_1)];$
$f \rightarrow [(v_3,v_4) \vee (v_4,v_3)];$
$g \rightarrow [(v_2,v_3) \vee (v_3,v_2)];$
$h \rightarrow [(v_1,v_3) \vee (v_3,v_1)].$

Рассматривая строки (столбцы) матрицы AB(G), можно определить центральные разрезы графа G путем перечисления элементов строк (столбцов):

$s(v_1) = s_1 = \{b,h,a\}$;
$s(v_2) = s_2 = \{b,g,c,d\}$;
$s(v_3) = s_3 = \{h,g,f\}$;
$s(v_4) = s_4 = \{c,f,e\}$;
$s(v_5) = s_5 = \{a,d,e\}$.

Рассматривая объединение строки и столбца матрицы AB(G) для выбранного текущего ребра $e_i \in E$, можно определить базовые реберные разрезы графа G (реберные разрезы 0-го уровня) [20,21]:

$w_0(a) = s_1 \oplus s_5 = \{b,h,d,e\}$;
$w_0(b) = s_1 \oplus s_2 = \{a,h,g,c,d\}$;
$w_0(c) = s_2 \oplus s_4 = \{b,g,d,f,e\}$;
$w_0(d) = s_2 \oplus s_5 = \{a,e,b,g,c\}$;
$w_0(e) = s_4 \oplus s_5 = \{c,f,a,d\}$;
$w_0(f) = s_3 \oplus s_4 = \{h,g,c,e\}$;
$w_0(g) = s_2 \oplus s_3 = \{h,f,b,c,d\}$;
$w_0(h) = s_1 \oplus s_3 = \{b,a,g,f\}$.

Будем менять местами вершины графа $G_3$ (рис. 2.9). Построим граф $G_4$ осуществляя перестановку (1 4)(2 3 5) вершин в графе $G_3$ (рис. 2.10). Построим граф $G_5$ осуществляя перестановку (1 2 5 3)(4) вершин в графе $G_3$ (рис. 2.11) .

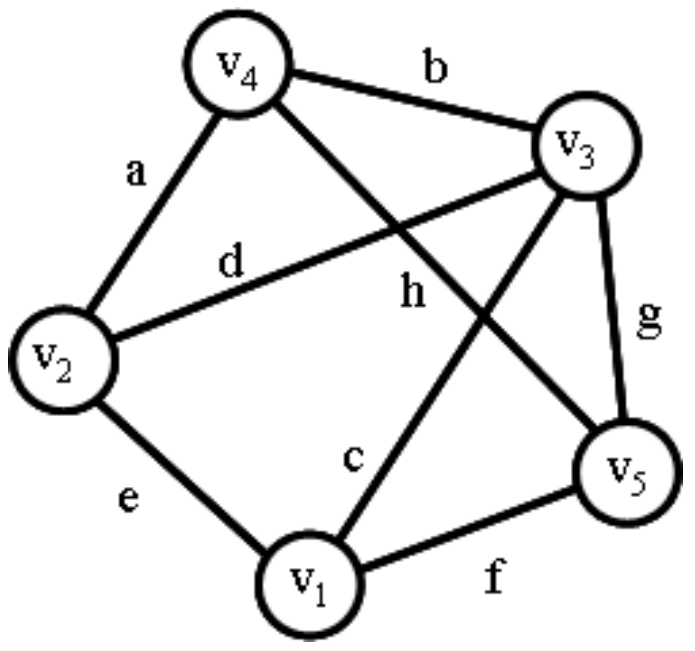

Рис. 2.9. Граф $G_3$.

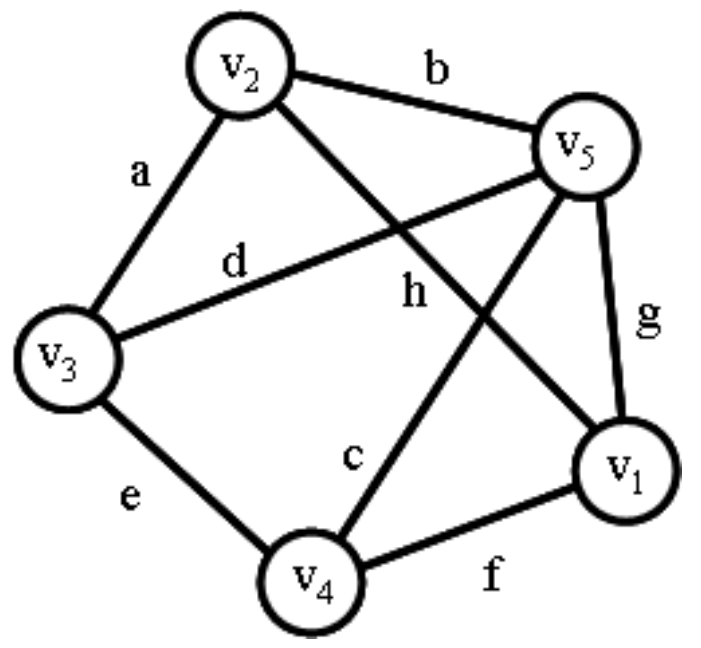

Рис. 2.10. Граф $G_4$.

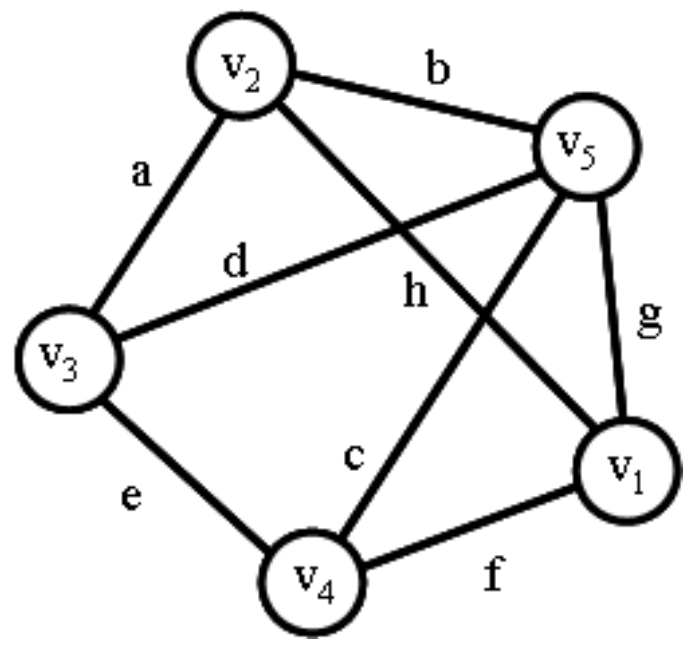

Рис. 2.11. Граф $G_5$.

Матрицы совместимостей графов $G_4$ и $G_5$ представлены ниже:

$AB(G_4)$ =

| | $v_1$ | $v_2$ | $v_3$ | $v_4$ | $v_5$ |
|---|---|---|---|---|---|
| $v_1$ | | e | c | | f |
| $v_2$ | e | | d | a | |
| $v_3$ | c | d | | b | g |
| $v_4$ | | a | b | | h |
| $v_5$ | f | | g | h | |

$AB(G_5)$=

| | $v_1$ | $v_2$ | $v_3$ | $v_4$ | $v_5$ |
|---|---|---|---|---|---|
| $v_1$ | | *h* | | *f* | *g* |
| $v_2$ | *h* | | *a* | | *b* |
| $v_3$ | | *a* | | *e* | *d* |
| $v_4$ | *f* | | *e* | | *c* |
| $v_5$ | *g* | *b* | *d* | *c* | |

Обратим внимание на тот факт, что базовый реберный разрез $w_0(a)$ для ребра *a* (серые цвета ячеек) остается неизменным при изменении нумерации вершин. Поэтому базовый реберный разрез можно рассматривать как дискретную функцию, где в качестве аргумента рассматривается ребро графа.

Рассмотрим изоморфные графы, представленные на рис. 2.9 – 2.11. Оставим без изменения идентификацию ребер, но изменим нумерацию вершин. Тогда базовые реберные разрезы графа:

$w_0(a) = [s_1 \oplus s_5]_{G_3} = [s_2 \oplus s_4]_{G_4} = [s_2 \oplus s_4]_{G_5} = \{b,h,d,e\}$;

$w_0(b) = [s_1 \oplus s_2]_{G_3} = [s_3 \oplus s_4]_{G_4} = [s_2 \oplus s_5]_{G_5} = \{a,h,g,c,d\}$;

$w_0(c) = [s_2 \oplus s_4]_{G_3} = [s_1 \oplus s_3]_{G_4} = [s_4 \oplus s_5]_{G_5} = \{b,g,d,f,e\}$;

$w_0(d) = [s_2 \oplus s_5]_{G_3} = [s_2 \oplus s_3]_{G_4} = [s_3 \oplus s_5]_{G_5} = \{a,e,b,g,c\}$;

$w_0(e) = [s_4 \oplus s_5]_{G_3} = [s_1 \oplus s_2]_{G_4} = [s_3 \oplus s_4]_{G_5} = \{c,f,a,d\}$;

$w_0(f) = [s_3 \oplus s_4]_{G_3} = [s_1 \oplus s_5]_{G_4} = [s_1 \oplus s_4]_{G_5} = \{h,g,c,e\}$;

$w_0(g) = [s_2 \oplus s_3]_{G_3} = [s_3 \oplus s_5]_{G_4} = [s_1 \oplus s_5]_{G_5} = \{h,f,b,c,d\}$;

$w_0(h) = [s_1 \oplus s_3]_{G_3} = [s_4 \oplus s_5]_{G_4} = [s_1 \oplus s_2]_{G_5} = \{b,a,g,f\}$.

Следует заметить, что перестановка вершин не изменяет значение и состав множества базовых реберных разрезов графа (рис. 2.9 – 2.11). Таким образом, множество базовых реберных разрезов графа является постоянной структурой графа и не зависит от нумерации вершин графа. Такие постоянные структуры в математике принято называть инвариантами.

Реберный разрез графа можно описывать в виде множества ребер, то есть в виде суграфа, но иногда выгоднее применять запись реберного разреза в виде характеристического кортежа.

**Определение 2.4**. *Характеристическим кортежем реберного разреза* $\mu(e_i)$ называется кортеж, состоящий из значений указывающих, является ли $e_i$, элементом множества $e_i \in$ E, $i$ = (1,2,…,$m$):

$$\mu(e_i) = \begin{cases} 1,\ если\, e_i \in E, \\ 0,\ если\, e_i \notin E \end{cases}. \tag{2.6}$$

Например, характеристический кортеж реберного разреза для суграфа w($b$) = {$a,h,g,c,d$} имеет вид: <1,0,1,1,0,0,1,1>.

$W_0(G_3)$ =

| | a | b | c | d | e | f | g | h |
|---|---|---|---|---|---|---|---|---|
| $w_0$(a) | | 1 | | 1 | 1 | | | 1 |
| $w_0$(b) | 1 | | 1 | 1 | | | 1 | 1 |
| $w_0$(c) | | 1 | | 1 | 1 | 1 | 1 | |
| $w_0$(d) | 1 | 1 | 1 | | 1 | | 1 | |
| $w_0$(e) | 1 | | 1 | 1 | | 1 | | |
| $w_0$(f) | | | 1 | | 1 | | 1 | 1 |
| $w_0$(g) | | 1 | 1 | 1 | | 1 | | 1 |
| $w_0$(h) | 1 | 1 | | | | 1 | 1 | |

Характеристические кортежи реберных разрезов позволяют легко строить матрицу базовых реберных разрезов $W_0(G)$ для графа G. Например, для графа $G_3$ матрица базовых реберных разрезов представлена выше.

Как видно, матрица базовых реберных разрезов $W_0(G)$ симметрична относительно главной диагонали.

## 2.3. Нильпотентный оператор реберных разрезов графа

Для установления свойств реберных разрезов графа, рассмотрим матрицу смежностей реберного графа $A(L(G_2))$ графа $G_2$ как линейный оператор $W_\lambda(G_2)$.

$W_\lambda(G_2) =$

| | $e_1$ | $e_2$ | $e_3$ | $e_4$ | $e_5$ | $e_6$ | $e_7$ | $e_8$ | $e_9$ | $e_{10}$ | $e_{11}$ |
|---|---|---|---|---|---|---|---|---|---|---|---|
| $e_1$ | | 1 | 1 | 1 | 1 | 1 | | | | | |
| $e_2$ | 1 | | 1 | | 1 | | 1 | | | 1 | |
| $e_3$ | 1 | 1 | | | | 1 | | | 1 | | 1 |
| $e_4$ | 1 | | | | 1 | 1 | 1 | 1 | 1 | | |
| $e_5$ | 1 | 1 | | 1 | | 1 | 1 | | | 1 | |
| $e_6$ | 1 | | 1 | 1 | 1 | | | | 1 | | 1 |
| $e_7$ | | 1 | | 1 | 1 | | | 1 | 1 | 1 | |
| $e_8$ | | | | 1 | | | 1 | | 1 | 1 | 1 |
| $e_9$ | | | 1 | 1 | | 1 | 1 | 1 | | | 1 |
| $e_{10}$ | | 1 | | | 1 | | 1 | 1 | | | 1 |
| $e_{11}$ | | | 1 | | | 1 | | 1 | 1 | 1 | |

Наличие такого оператора в пространстве суграфов графа $G_2$ позволяет осуществлять преобразование векторного пространства суграфов само в себя. То есть имеет смысл рассматривать операторы $W_\lambda^2 = W_\lambda \times W_\lambda, W_\lambda^3 = W_\lambda \times W_\lambda \times W_\lambda, \ldots, W_\lambda^q = W_\lambda \times W_\lambda \times \ldots \times W_\lambda$

**Определение 2.5**. Оператор $W_\lambda$ называется *нильпотентным оператором в пространстве суграфов графа G*, если существует такое натуральное число *q*, что $W_\lambda^q = \varnothing$ или $W_\lambda^q = W_\lambda = \varnothing$ [18,24].

Данной матрице оператора $W_\lambda$ соответствует следующий реберный граф (рис. 2.12).

Исходя из рассмотрения свойств суграфов графа G, оператор A(L(G)) является нильпотентным оператором в пространстве суграфов £(G).

Рассчитаем цепочку операторов для нашего примера. Определим оператор $W_\lambda^2$, как результат умножения двух матриц A(L(G)), применяя закон сложения элементов по mod 2.

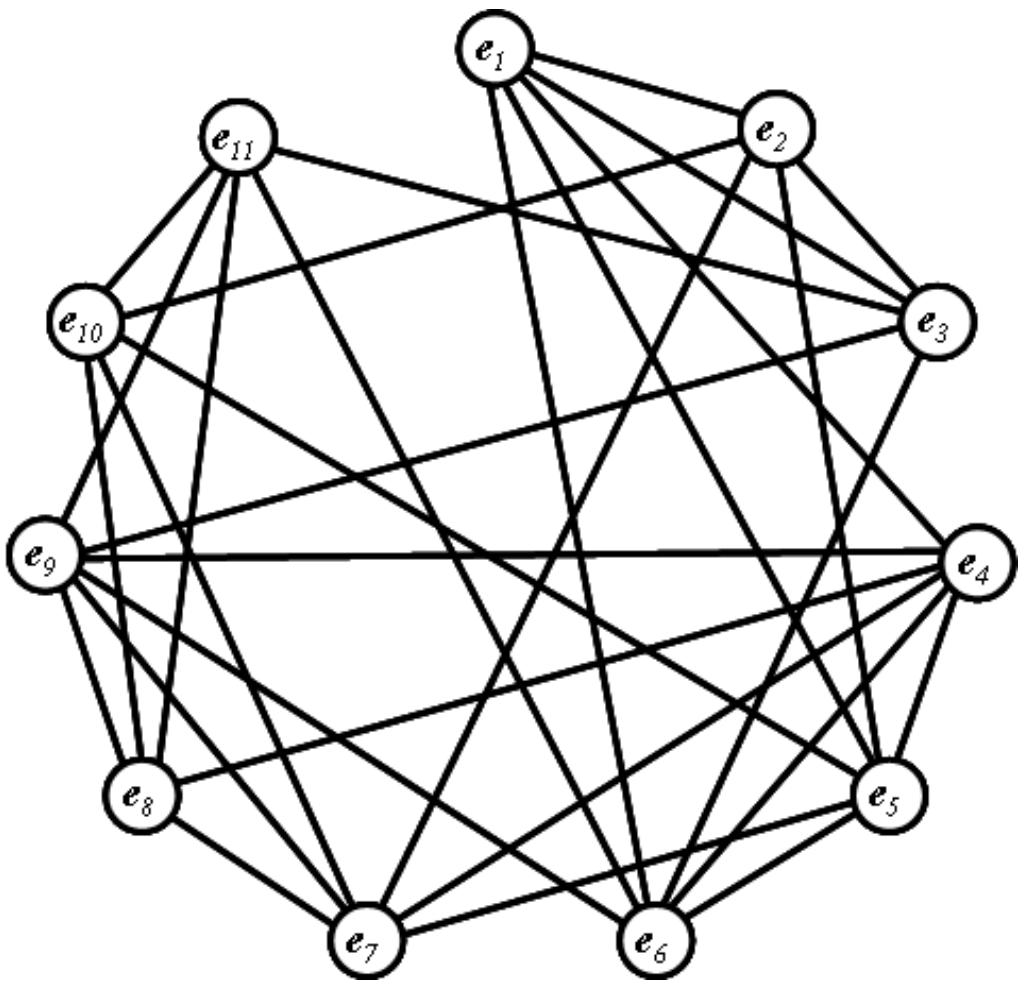

Рис. 2.12. Нильпотентный реберный граф L(G2).

Рассмотрим оператор $W_\lambda^2$.

$W_\lambda^2$ =

| | $e_1$ | $e_2$ | $e_3$ | $e_4$ | $e_5$ | $e_6$ | $e_7$ | $e_8$ | $e_9$ | $e_{10}$ | $e_{11}$ |
|---|---|---|---|---|---|---|---|---|---|---|---|
| $e_1$ | 1 | | | | 1 | 1 | 1 | 1 | 1 | | |
| $e_2$ | | 1 | 1 | 1 | 1 | 1 | | | | | |
| $e_3$ | | 1 | 1 | 1 | 1 | 1 | | | | | |
| $e_4$ | | 1 | 1 | | 1 | 1 | 1 | | 1 | 1 | 1 |
| $e_5$ | 1 | 1 | 1 | 1 | | | 1 | 1 | 1 | | |
| $e_6$ | 1 | 1 | 1 | 1 | | | 1 | 1 | 1 | | |
| $e_7$ | 1 | | | 1 | 1 | 1 | | 1 | | 1 | 1 |
| $e_8$ | 1 | | | | 1 | 1 | 1 | 1 | 1 | | |
| $e_9$ | 1 | | | 1 | 1 | 1 | | 1 | | 1 | 1 |
| $e_{10}$ | | | | 1 | | | 1 | | 1 | 1 | 1 |
| $e_{11}$ | | | | 1 | | | 1 | | 1 | 1 | 1 |

$W_\lambda^3$ =

| | $e_1$ | $e_2$ | $e_3$ | $e_4$ | $e_5$ | $e_6$ | $e_7$ | $e_8$ | $e_9$ | $e_{10}$ | $e_{11}$ |
|---|---|---|---|---|---|---|---|---|---|---|---|
| $e_1$ | | 1 | 1 | | 1 | 1 | 1 | | 1 | 1 | 1 |
| $e_2$ | 1 | | | | 1 | 1 | 1 | 1 | 1 | | |
| $e_3$ | 1 | | | | 1 | 1 | 1 | 1 | 1 | | |
| $e_4$ | | | | | | | | | | | |
| $e_5$ | 1 | 1 | 1 | | | | | 1 | | 1 | 1 |
| $e_6$ | 1 | 1 | 1 | | | | | 1 | | 1 | 1 |
| $e_7$ | 1 | 1 | 1 | | | | | 1 | | 1 | 1 |
| $e_8$ | | 1 | 1 | | 1 | 1 | 1 | | 1 | 1 | 1 |
| $e_9$ | 1 | 1 | 1 | | | | | 1 | | 1 | 1 |
| $e_{10}$ | 1 | | | | 1 | 1 | 1 | 1 | 1 | | |
| $e_{11}$ | 1 | | | | 1 | 1 | 1 | 1 | 1 | | |

$W_\lambda^4 =$

| | e1 | e2 | e3 | e4 | e5 | e6 | e7 | e8 | e9 | e10 | e11 |
|---|---|---|---|---|---|---|---|---|---|---|---|
| e1 | | | | | | | | | | | |
| e2 | | 1 | 1 | | 1 | 1 | 1 | | 1 | 1 | 1 |
| e3 | | 1 | 1 | | 1 | 1 | 1 | | 1 | 1 | 1 |
| e4 | | | | | | | | | | | |
| e5 | | 1 | 1 | | 1 | 1 | 1 | | 1 | 1 | 1 |
| e6 | | 1 | 1 | | 1 | 1 | 1 | | 1 | 1 | 1 |
| e7 | | 1 | 1 | | 1 | 1 | 1 | | 1 | 1 | 1 |
| e8 | | | | | | | | | | | |
| e9 | | 1 | 1 | | 1 | 1 | 1 | | 1 | 1 | 1 |
| e10 | | 1 | 1 | | 1 | 1 | 1 | | 1 | 1 | 1 |
| e11 | | 1 | 1 | | 1 | 1 | 1 | | 1 | 1 | 1 |

Прекратим построение цепочки операторов, когда имеется показатель степени оператора со значением $q$, при котором все элементы линейного оператора $W_\lambda^q$ либо равны нулю, либо все элементы линейного оператора равны элементам базового оператора $W_\lambda$. То есть, по условию остановки процесса порождения реберных разрезов, все элементы оператора $W_\lambda^q$ со значением степени $q$ и выше – пустое множество. В нашем примере оператор пятой степени равен нулю

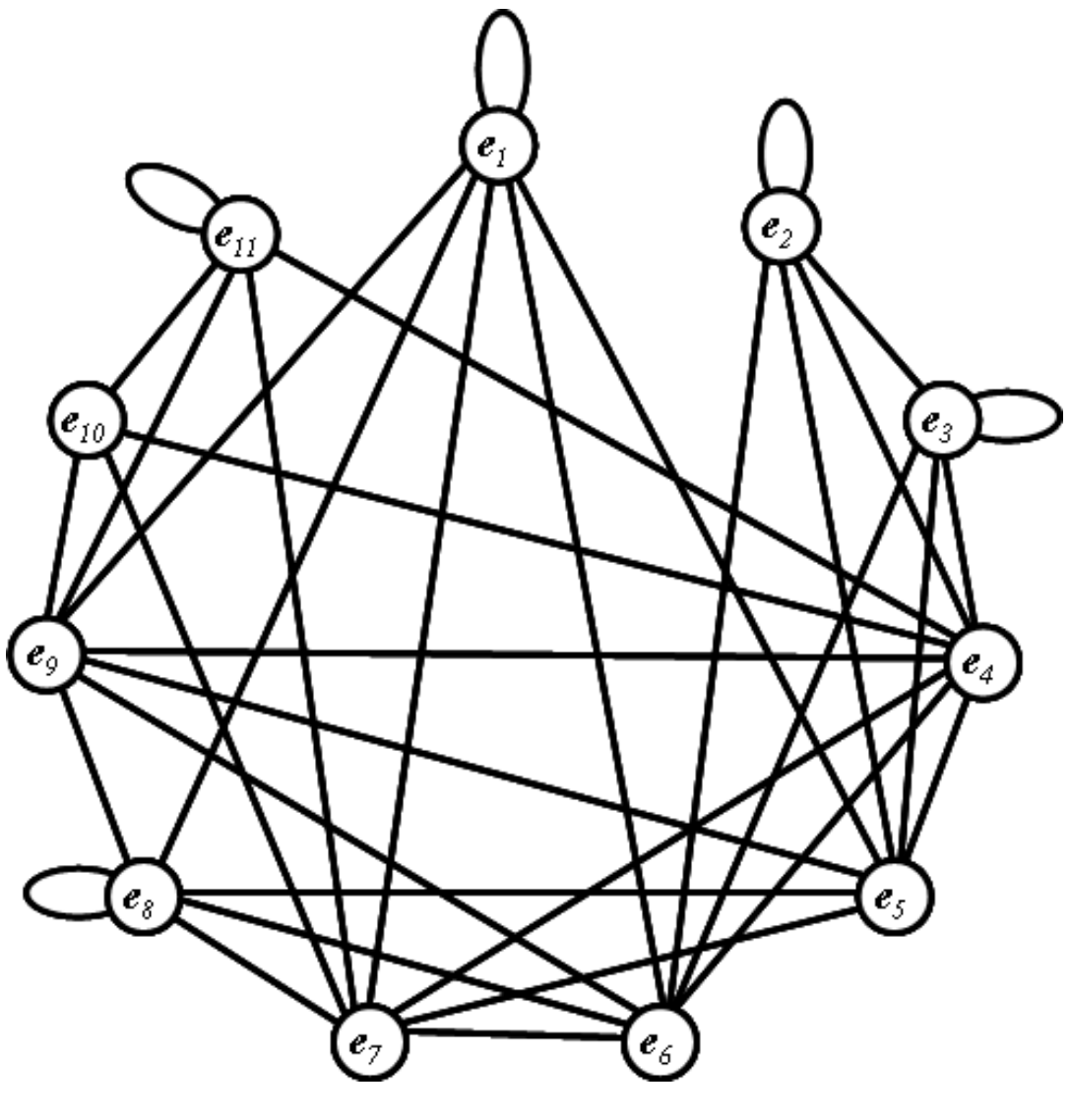

Рис. 2.13. Реберный граф $W_\lambda^2$.

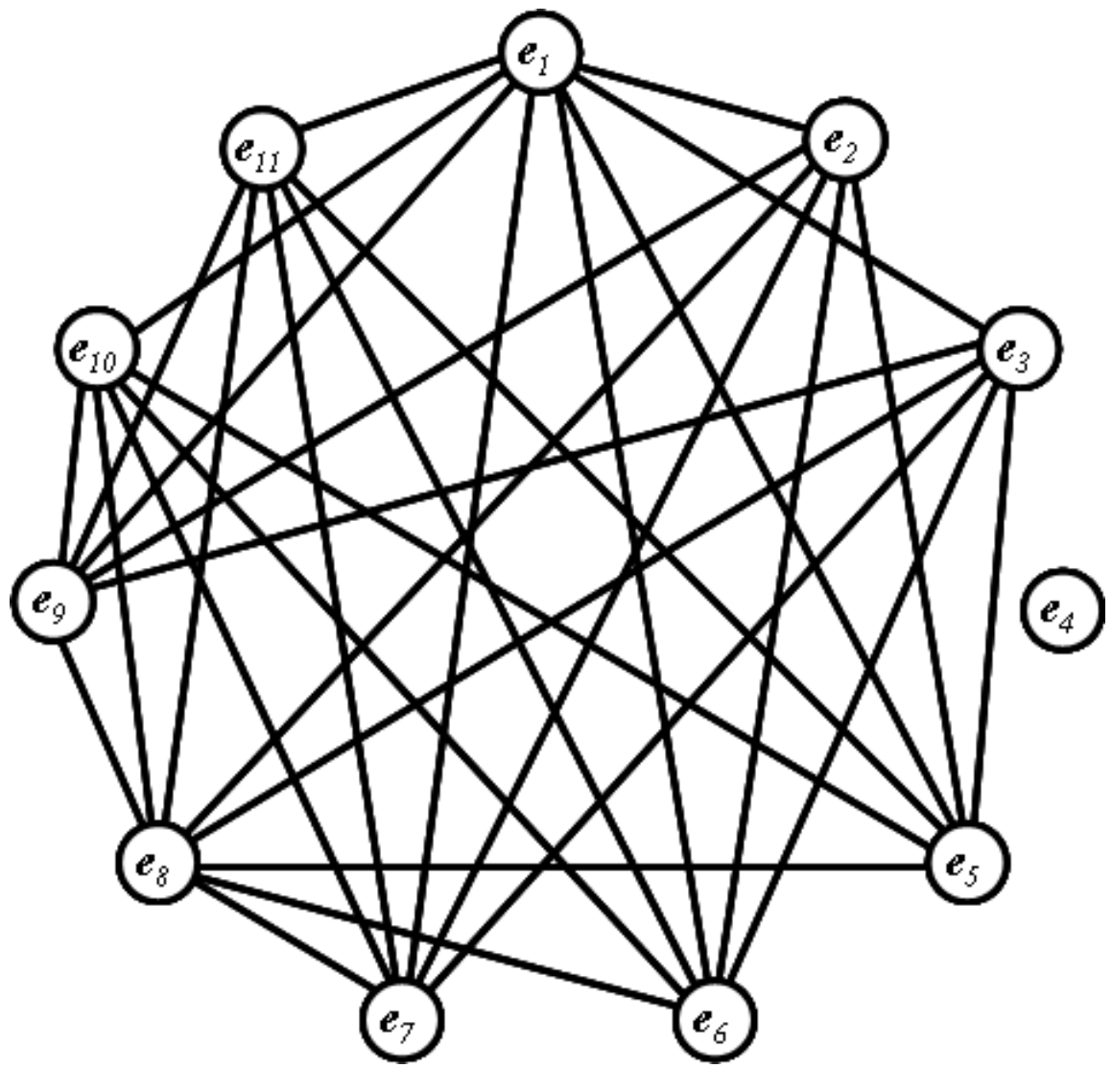

Рис. 2.14. Реберный граф $W_\lambda^3$.

Исходя из цепочки операторов формируются множества реберных разрезов различного уровня, которые определяются по формулам:

$$W_0(G_2)^T = W_\lambda(G_2) \times R(G_2)^T; \quad (2.7)$$

$$W_1(G_2)^T = W_\lambda(G_2) \times W_0(G_2)^T = W_\lambda^2(G_2) \times R(G_1)^T = W_\lambda(G_2) \times W_\lambda(G_2) \times R(G_2)^T; \quad (2.8)$$

$$W_2(G_2)^T = W_\lambda(G_2) \times W_1(G_2)^T = W_\lambda^2(G_2) \times W_0(G_2)^T = W_\lambda^3(G_2) \times R(G_2)^T. \quad (2.9)$$

и т.д.

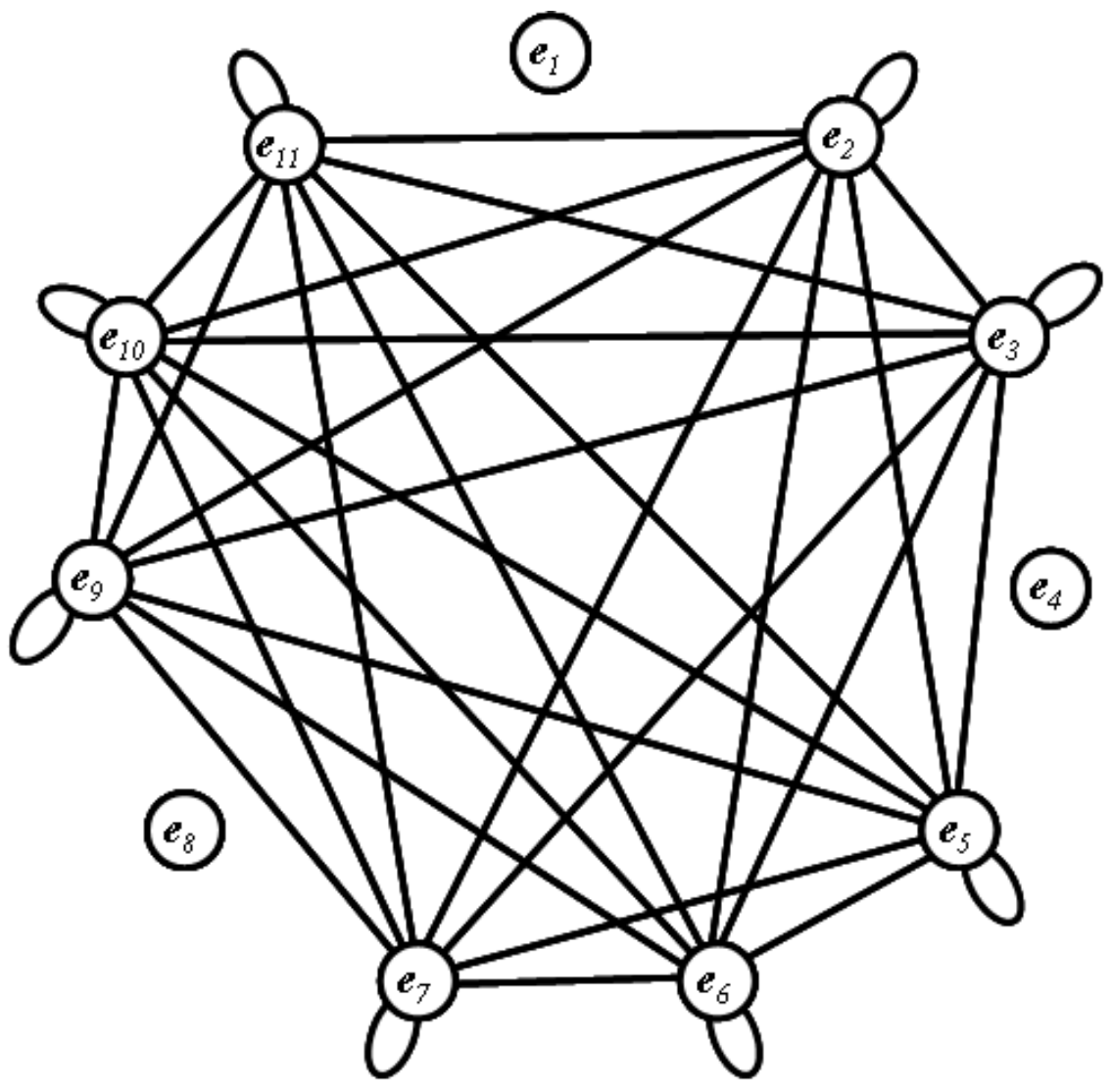

Рис. 2.15. Реберный граф $W_\lambda^4$.

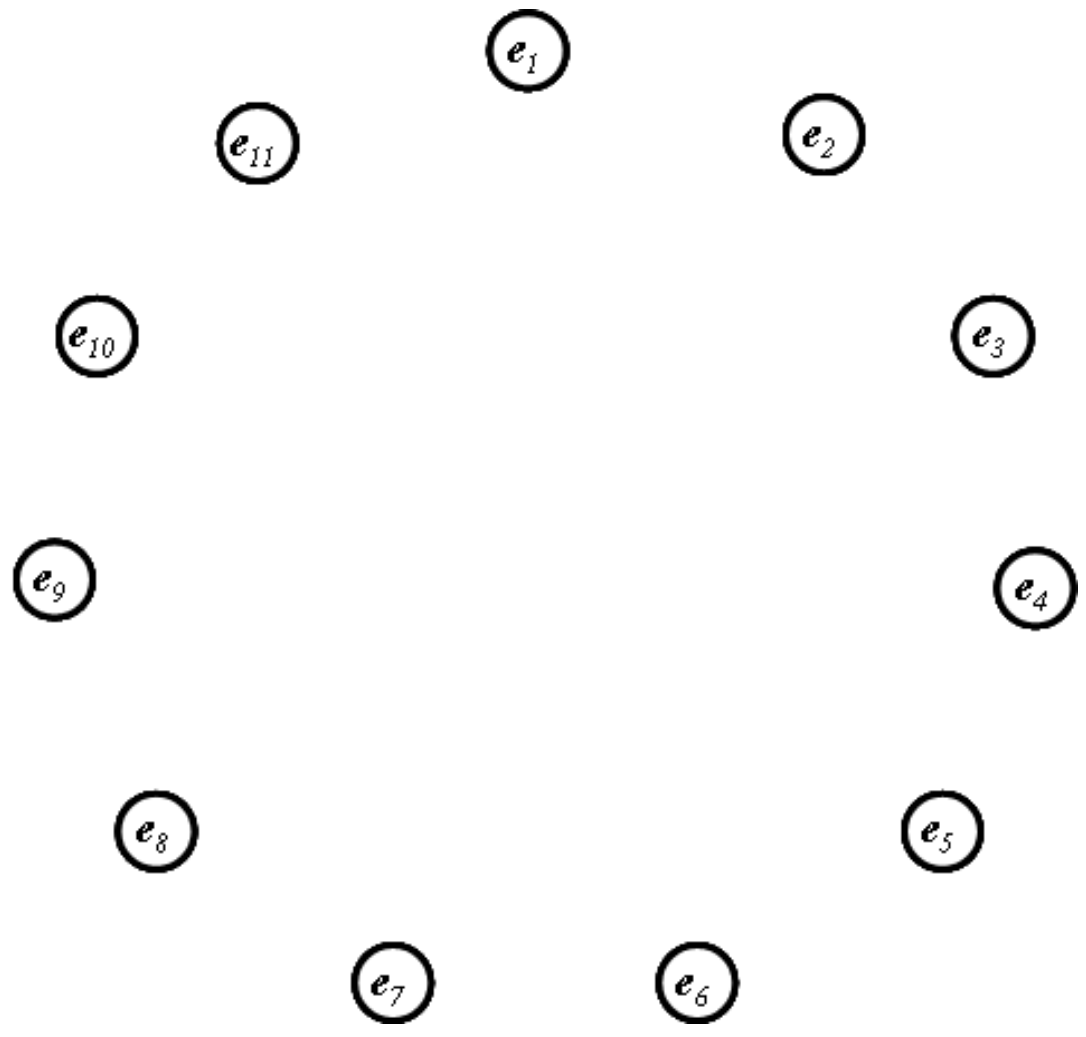

Рис. 2.16. Реберный граф $W_\lambda^5$.

В процессе построения цепочки операторов, на определенном уровне появляются нулевые строчки. Появление нулевых строчек определяет общий совпадающий элемент из подпространства разрезов графа S и подпространства циклов C. Для графа $G_2$ существует суграф, у которого кольцевая сумма центральных разрезов равна квазициклу, на основании Леммы 2.3: $s_1 \oplus s_2 \oplus s_3 \oplus s_5 = \{e_2, e_3, e_5, e_6, e_7, e_9, e_{10}, e_{11}\}$.

## 2.4. Спектр реберных разрезов графа

Реберный разрез $w(e_i)$ можно рассматривать как элемент множества разрезов графа S(G). В классической теории графов ребра реберного графа L(G) порождают новый реберный граф L(L(G)) и т. д. [32]. В общем случае итерированный реберный граф можно определить рекуррентным соотношением $L_n(G) = L(L_{n-1}(G)), n \geq 2$. Но стоит заметить, что в таком представлении количество ребер при переходе от реберного графа $L_{n-1}(G)$ к графу $L_n(G)$ увеличивается.

В отличие от порождения итерированных реберных графов, рассмотрим цепочку порождения итерированных квалиразрезов исходного графа, определяемых рекуррентным соотношением $W_k(G) = W \times W_{k-1}(G), k \geq 2$ по принципу: квалиразрез *k*-го уровня (яруса) $W_k(G)$ порождается квалиразрезом *k*-1 уровня $W_{k-1}(G)$. В отличие от цепочки реберных графов, такая цепочка порождает реберные разрезы исходного графа всегда с количеством не более величины $m$ – количества ребер, и такое порождение всегда конечно. Окончание построения цепочки реберных разрезов определяется значением показателя степени оператора W при ограничении $\mathrm{W}^{\mathrm{q}} = \mathrm{W}_0 \to \mathrm{W}^{\mathrm{q}} = \varnothing$.

С другой стороны, так как

$$\mathrm{W}_k(\mathrm{G}) = \mathrm{W} \times \mathrm{W}_{k-1}(\mathrm{G}) = \mathrm{W}^{k-1} \times \mathrm{W} \times \beta(\mathrm{G}) = \mathrm{W}^{k-1} \times \mathrm{W}_1(\mathrm{G}), \tag{2.10}$$

то порождение цепочек квалиразрезов *k*-го уровня (яруса) графа G для каждого ребра осуществляется путем сложения по модулю 2 суграфов базовых реберных разрезов определенных реберным разрезом (*k-1*)-го уровня этого ребра. Каждое новое подмножество ребер формируется рекурсивно с помощью предыдущего множества с помощью преобразования $\gamma$ (будем обозначать преобразование $\gamma$ греческой букыой «гамма»):

$$\mathrm{w}_k(e_\mathrm{i}) = \begin{cases} \mathrm{w}_k(e_\mathrm{i}) = \gamma(\mathrm{w}_{\mathrm{k\text{-}1}}(e_\mathrm{i})); \\ \textit{если } \mathrm{w}_\mathrm{k}(e_\mathrm{i}) = \mathrm{w}_\mathrm{h}(e_\mathrm{i}), \textit{ то } \mathrm{w}_\mathrm{k}(e_\mathrm{i}) = \varnothing. \end{cases} \tag{2.11}$$

где $i$ – номер ребра, $k$ – номер уровня, $h$ – номер произвольного предыдущего уровня.

Преобразование $\gamma(w_{k-1}(e_i))$ определим следующим образом:

$$\gamma(\{a,b,c,...,n\}) = \mathrm{w}_0(a) \oplus \mathrm{w}_0(\mathrm{b}) \oplus \mathrm{w}_0(c) \oplus ... \oplus \mathrm{w}_0(n). \tag{2.12}$$

$\gamma(w_{k-1}(e_i))$ – *порождение реберного разреза последующего уровня* определяет строку из предыдущего реберного разреза суграфа$\{a,b,c,...,n\}$ ребра $e_i$ состоящего из множества ребер базовых реберных разрезов $\mathrm{w}_0(a) \oplus \mathrm{w}_0(\mathrm{b}) \oplus \mathrm{w}_0(c) \oplus ... \oplus \mathrm{w}_0(n)$ [20].

Все уровневые подмножества (реберные разрезы) также будут представлять собой квалиразрезы исходного графа G. Так, первый уровень будет состоять из множества базовых реберных разрезов графа G, а все последующие реберные разрезы – порождаться преобразованием $\gamma(w_{k-1}(e_i))$. После получения циклического повтора подмножеств, данному подмножеству ставится в соответствие пустое множество.

На каждом уровне образуются $m$ подмножеств (суграфов), зависящих от выбранного ребра, что дает нам возможность построить прямоугольную матрицу спектра реберных разрезов $\mathrm{W_s}$ размером $m \times k$, где $k$ – количество уровней. Очевидно, что элементы этой матрицы могут быть записаны в виде суграфов исходного графа. Таким образом, реберные разрезы графа состоят из подмножества ребер и порождают ярусные квалиразрезы графа. Любой реберный разрез можно представить в виде функции $\mathrm{w_l}(e_\mathrm{i})$ зависящей от ребра $e_\mathrm{i}$ в качестве аргумента. Учитывая, что окончательное количество уровней заранее не известно, предварительно введем определение спектра реберных разрезов.

**Определение 2.6.** Совокупность всех реберных разрезов порожденных множеством базовых реберных разрезов при ограничении $W^q = W_0 \rightarrow W^q = \varnothing$ называется *спектром реберных разрезов графа*.

$W_s(G)$ =

| | уровни $l_0$ | $l_1$ | .. | $l_k$ |
|---|---|---|---|---|
| $e_1$ | $w_0(e_1)$ | $w_1(e_1)$ | … | $w_k(e_1$ |
| $e_2$ | $w_0(e_2)$ | $w_1(e_2)$ | … | $w_k(e_2$ |
| $e_3$ | $w_0(e_3)$ | $w_1(e_3)$ | … | $w_k(e_3$ |
| … | … | … | … | … |
| $e_m$ | $w_0(e_m)$ | $w_1(e_m)$ | … | $w_k(e_m$ |

В качестве примера выделения спектра реберных разрезов и построения матрицы $W_s(G)$, рассмотрим следующий граф $G_2$ (рис. 2.4).

Множество центральных разрезов графа:

$s(v_1) = s_1 = \{e_1,e_2,e_3\}$;
$s(v_2) = s_2 = \{e_1,e_4,e_5,e_6\}$;
$s(v_3) = s_3 = \{e_4,e_7,e_8,e_9\}$;
$s(v_4) = s_4 = \{e_2,e_5,e_7,e_{10}\}$;
$s(v_5) = s_5 = \{e_8,e_{10},e_{11}\}$;
$s(v_6) = s_6 = \{e_3,e_6,e_9,e_{11}\}$.

Множество базовых реберных разрезов графа (0-уровень квалиразрезов):

$w_0(e_1) = s_1 \oplus s_2 = \{e_1,e_2,e_3\} \oplus \{e_1,e_4,e_5,e_6\} = \{e_2,e_3,e_4,e_5,e_6\}$;
$w_0(e_2) = s_1 \oplus s_4 = \{e_1,e_2,e_3\} \oplus \{e_2,e_5,e_7,e_{10}\} = \{e_1,e_3,e_5,e_7,e_{10}\}$;
$w_0(e_3) = s_1 \oplus s_6 = \{e_1,e_2,e_3\} \oplus \{e_3,e_6,e_9,e_{11}\} = \{e_1,e_2,e_6,e_9,e_{11}\}$;
$w_0(e_4) = s_2 \oplus s_3 = \{e_1,e_4,e_5,e_6\} \oplus \{e_4,e_7,e_8,e_9\} = \{e_1,e_5,e_6,e_7,e_8,e_9\}$;
$w_0(e_5) = s_2 \oplus s_4 = \{e_1,e_4,e_5,e_6\} \oplus \{e_2,e_5,e_7,e_{10}\} = \{e_1,e_2,e_4,e_6,e_7,e_{10}\}$;
$w_0(e_6) = s_2 \oplus s_6 = \{e_1,e_4,e_5,e_6\} \oplus \{e_3,e_6,e_9,e_{11}\} = \{e_1,e_3,e_4,e_5,e_9,e_{11}\}$;
$w_0(e_7) = s_3 \oplus s_4 = \{e_4,e_7,e_8,e_9\} \oplus \{e_2,e_5,e_7,e_{10}\} = \{e_2,e_4,e_5,e_8,e_9,e_{10}\}$;
$w_0(e_8) = s_3 \oplus s_5 = \{e_4,e_7,e_8,e_9\} \oplus \{e_8,e_{10},e_{11}\} = \{e_4,e_7,e_9,e_{10},e_{11}\}$;
$w_0(e_9) = s_3 \oplus s_6 = \{e_4,e_7,e_8,e_9\} \oplus \{e_3,e_6,e_9,e_{11}\} = \{e_3,e_4,e_6,e_7,e_8,e_{11}\}$;
$w_0(e_{10}) = s_4 \oplus s_5 = \{e_2,e_5,e_7,e_{10}\} \oplus \{e_8,e_{10},e_{11}\} = \{e_2,e_5,e_7,e_8,e_{11}\}$;
$w_0(e_{11}) = s_5 \oplus s_6 = \{e_8,e_{10},e_{11}\} \oplus \{e_3,e_6,e_9,e_{11}\} = \{e_3,e_6,e_8,e_9,e_{10}\}$.

Множество квалиразрезов 1-го уровня (яруса):

$w_1(e_1) = \gamma(w_0(e_1)) = \gamma(\{e_2,e_3,e_4,e_5,e_6\}) = s_1 \oplus s_4 \oplus s_1 \oplus s_6 \oplus s_2 \oplus s_3 \oplus s_2 \oplus s_4 \oplus s_2 \oplus s_6 =$
$= s_2 \oplus s_3 = \{e_1,e_4,e_5,e_6\} \oplus \{e_4,e_7,e_8,e_9\} = \{e_1,e_5,e_6,e_7,e_8,e_9\}$;
$w_1(e_2) = \gamma(w_0(e_2)) = \gamma(\{e_1,e_3,e_5,e_7,e_{10}\}) = s_1 \oplus s_2 \oplus s_1 \oplus s_6 \oplus s_2 \oplus s_4 \oplus s_3 \oplus s_4 \oplus s_4 \oplus s_5 =$
$= s_3 \oplus s_4 \oplus s_5 \oplus s_6 = \{e_4,e_7,e_8,e_9\} \oplus \{e_2,e_5,e_7,e_{10}\} \oplus \{e_8,e_{10},e_{11}\} \oplus \{e_3,e_6,e_9,e_{11}\} =$
$= \{e_2,e_3,e_4,e_5,e_6\}$;
$w_1(e_3) = \gamma(w_0(e_3)) = \gamma(\{e_1,e_2,e_6,e_9,e_{11}\}) = s_1 \oplus s_2 \oplus s_1 \oplus s_4 \oplus s_2 \oplus s_6 \oplus s_3 \oplus s_6 \oplus s_5 \oplus s_6 =$
$= s_3 \oplus s_4 \oplus s_5 \oplus s_6 = \{e_4,e_7,e_8,e_9\} \oplus \{e_2,e_5,e_7,e_{10}\} \oplus \{e_8,e_{10},e_{11}\} \oplus \{e_3,e_6,e_9,e_{11}\} =$
$= \{e_2,e_3,e_4,e_5,e_6\}$;
$w_1(e_4) = \gamma(w_0(e_4)) = \gamma(\{e_1,e_5,e_6,e_7,e_8,e_9\}) = s_1 \oplus s_2 \oplus s_2 \oplus s_4 \oplus s_2 \oplus s_6 \oplus s_3 \oplus s_4 \oplus s_3 \oplus s_5 \oplus$
$\oplus s_3 \oplus s_6 = s_1 \oplus s_2 \oplus s_3 \oplus s_5 = \{e_1,e_2,e_3\} \oplus \{e_1,e_4,e_5,e_6\} \oplus \{e_4,e_7,e_8,e_9\} \oplus \{e_8,e_{10},e_{11}\} =$
$= \{e_2,e_3,e_5,e_6,e_7,e_9,e_{10},e_{11}\}$;
$w_1(e_5) = \gamma(w_0(e_5)) = \gamma(\{e_1,e_2,e_4,e_6,e_7,e_{10}\}) = s_1 \oplus s_2 \oplus s_1 \oplus s_4 \oplus s_2 \oplus s_3 \oplus s_2 \oplus s_6 \oplus s_3 \oplus s_4 \oplus$
$\oplus s_4 \oplus s_5 = s_2 \oplus s_4 \oplus s_5 \oplus s_6 = \{e_1,e_4,e_5,e_6\} \oplus \{e_2,e_5,e_7,e_{10}\} \oplus \{e_8,e_{10},e_{11}\} \oplus$
$\oplus \{e_3,e_6,e_9,e_{11}\} = \{e_1,e_2,e_3,e_4,e_7,e_8,e_9\}$;

$w_1(e_6) = \gamma(w_0(e_6)) = \gamma(\{e_1,e_3,e_4,e_5,e_9,e_{11}\}) =$
$= s_1 \oplus s_2 \oplus s_1 \oplus s_6 \oplus s_2 \oplus s_3 \oplus s_2 \oplus s_4 \oplus s_3 \oplus s_6 \oplus s_5 \oplus s_6 =$
$= s_2 \oplus s_4 \oplus s_5 \oplus s_6 = \{e_1,e_4,e_5,e_6\} \oplus \{e_2,e_5,e_7,e_{10}\} \oplus \{e_8,e_{10},e_{11}\} \oplus \{e_3,e_6,e_9,e_{11}\} =$
$= \{e_1,e_2,e_3,e_4,e_7,e_8,e_9\}$;
$w_1(e_7) = \gamma(w_0(e_7)) = \gamma(\{e_2,e_4,e_5,e_8,e_9,e_{10}\}) =$
$= s_1 \oplus s_4 \oplus s_2 \oplus s_3 \oplus s_2 \oplus s_4 \oplus s_3 \oplus s_5 \oplus s_3 \oplus s_6 \oplus s_5 \oplus s_6 =$
$= s_1 \oplus s_3 \oplus s_4 \oplus s_6 = \{e_1,e_2,e_3\} \oplus \{e_4,e_7,e_8,e_9\} \oplus \{e_2,e_5,e_7,e_{10}\} \oplus \{e_3,e_6,e_9,e_{11}\} =$
$= \{e_1,e_4,e_5,e_6,e_8,e_{10},e_{11}\}$;
$w_1(e_8) = \gamma(w_0(e_8)) = \gamma(\{e_4,e_7,e_9,e_{10},e_{11}\}) = s_2 \oplus s_3 \oplus s_3 \oplus s_4 \oplus s_3 \oplus s_6 \oplus s_4 \oplus s_5 \oplus s_5 \oplus s_6 =$
$= s_2 \oplus s_3 = \{e_1,e_4,e_5,e_6\} \oplus \{e_4,e_7,e_8,e_9\} = \{e_1,e_5,e_6,e_7,e_8,e_9\}$;
$w_1(e_9) = \gamma(w_0(e_9)) = \gamma(\{e_3,e_4,e_6,e_7,e_8,e_{11}\}) =$
$= s_1 \oplus s_6 \oplus s_2 \oplus s_3 \oplus s_2 \oplus s_6 \oplus s_3 \oplus s_4 \oplus s_3 \oplus s_5 \oplus s_5 \oplus s_6 =$
$= s_1 \oplus s_3 \oplus s_4 \oplus s_6 = \{e_1,e_2,e_3\} \oplus \{e_4,e_7,e_8,e_9\} \oplus \{e_2,e_5,e_7,e_{10}\} \oplus \{e_3,e_6,e_9,e_{11}\} =$
$= \{e_1,e_4,e_5,e_6,e_8,e_{10},e_{11}\}$;
$w_1(e_{10}) = \gamma(w_0(e_{10})) = \gamma(\{e_2,e_5,e_7,e_8,e_{11}\}) = s_1 \oplus s_4 \oplus s_2 \oplus s_4 \oplus s_3 \oplus s_4 \oplus s_3 \oplus s_5 \oplus s_5 \oplus s_6 =$
$= s_1 \oplus s_2 \oplus s_4 \oplus s_6 = \{e_1,e_2,e_3\} \oplus \{e_1,e_4,e_5,e_6\} \oplus \{e_2,e_5,e_7,e_{10}\} \oplus \{e_3,e_6,e_9,e_{11}\} =$
$= \{e_4,e_7,e_9,e_{10},e_{11}\}$;
$w_1(e_{11}) = \gamma(w_0(e_{11})) = \gamma(\{e_3,e_6,e_8,e_9,e_{10}\}) = s_1 \oplus s_6 \oplus s_2 \oplus s_6 \oplus s_3 \oplus s_5 \oplus s_3 \oplus s_6 \oplus s_5 \oplus s_6 =$
$= s_1 \oplus s_2 \oplus s_4 \oplus s_6 = \{e_1,e_2,e_3\} \oplus \{e_1,e_4,e_5,e_6\} \oplus \{e_2,e_5,e_7,e_{10}\} \oplus \{e_3,e_6,e_9,e_{11}\} =$
$= \{e_4,e_7,e_9,e_{10},e_{11}\}$.

Множество квалиразрезов 2-го уровня (яруса):

$w_2(e_1) = \gamma(\{e_1,e_5,e_6,e_7,e_8,e_9\}) = s_1 \oplus s_2 \oplus s_2 \oplus s_4 \oplus s_2 \oplus s_6 \oplus s_3 \oplus s_4 \oplus s_3 \oplus s_5 \oplus s_3 \oplus s_6 =$
$= s_1 \oplus s_2 \oplus s_3 \oplus s_5 = \{e_1,e_2,e_3\} \oplus \{e_1,e_4,e_5,e_6\} \oplus \{e_4,e_7,e_8,e_9\} \oplus \{e_8,e_{10},e_{11}\} =$
$= \{e_2,e_3,e_5,e_6,e_7,e_9,e_{10},e_{11}\}$;
$w_2(e_2) = \gamma(\{e_2,e_3,e_4,e_5,e_6\}) = s_1 \oplus s_4 \oplus s_1 \oplus s_6 \oplus s_2 \oplus s_3 \oplus s_2 \oplus s_4 \oplus s_2 \oplus s_6 = s_2 \oplus s_3 =$
$= \{e_1,e_4,e_5,e_6\} \oplus \{e_4,e_7,e_8,e_9\} = \{e_1,e_5,e_6,e_7,e_8,e_9\}$;
$w_2(e_3) = \gamma(\{e_2,e_3,e_4,e_5,e_6\}) = s_1 \oplus s_4 \oplus s_1 \oplus s_6 \oplus s_2 \oplus s_3 \oplus s_2 \oplus s_4 \oplus s_2 \oplus s_6 = s_2 \oplus s_3 =$
$= \{e_1,e_4,e_5,e_6\} \oplus \{e_4,e_7,e_8,e_9\} = \{e_1,e_5,e_6,e_7,e_8,e_9\}$;
$w_2(e_4) = \gamma(\{e_2,e_3,e_5,e_6,e_7,e_9,e_{10},e_{11}\}) = s_1 \oplus s_4 \oplus s_1 \oplus s_6 \oplus s_2 \oplus s_4 \oplus s_2 \oplus s_6 \oplus s_3 \oplus s_4 \oplus$
$\oplus s_3 \oplus s_6 \oplus s_4 \oplus s_5 \oplus s_5 \oplus s_6 = \varnothing$;
$w_2(e_5) = \gamma(\{e_1,e_2,e_3,e_4,e_7,e_8,e_9\}) = s_1 \oplus s_2 \oplus s_1 \oplus s_4 \oplus s_1 \oplus s_6 \oplus s_2 \oplus s_3 \oplus s_3 \oplus s_4 \oplus s_3 \oplus$
$\oplus s_5 \oplus s_3 \oplus s_6 = s_1 \oplus s_5 = \{e_1,e_2,e_3\} \oplus \{e_8,e_{10},e_{11}\} = \{e_1,e_2,e_3,e_8,e_{10},e_{11}\}$;
$w_2(e_6) = \gamma(\{e_1,e_2,e_3,e_4,e_7,e_8,e_9\}) = s_1 \oplus s_2 \oplus s_1 \oplus s_4 \oplus s_1 \oplus s_6 \oplus s_2 \oplus s_3 \oplus s_3 \oplus s_4 \oplus s_3 \oplus$
$\oplus s_5 \oplus s_3 \oplus s_6 = s_1 \oplus s_5 = \{e_1,e_2,e_3\} \oplus \{e_8,e_{10},e_{11}\} = \{e_1,e_2,e_3,e_8,e_{10},e_{11}\}$;
$w_2(e_7) = \gamma(\{e_1,e_4,e_5,e_6,e_8,e_{10},e_{11}\}) = s_1 \oplus s_2 \oplus s_2 \oplus s_3 \oplus s_2 \oplus s_4 \oplus s_2 \oplus s_6 \oplus s_3 \oplus s_5 \oplus s_4 \oplus$
$\oplus s_5 \oplus s_5 \oplus s_6 = s_1 \oplus s_5 = \{e_1,e_2,e_3\} \oplus \{e_8,e_{10},e_{11}\} = \{e_1,e_2,e_3,e_8,e_{10},e_{11}\}$;
$w_2(e_8) = \gamma(\{e_1,e_5,e_6,e_7,e_8,e_9\}) = s_1 \oplus s_2 \oplus s_2 \oplus s_4 \oplus s_2 \oplus s_6 \oplus s_3 \oplus s_4 \oplus s_3 \oplus s_5 \oplus s_3 \oplus s_6 =$
$= s_1 \oplus s_2 \oplus s_3 \oplus s_5 = \{e_1,e_2,e_3\} \oplus \{e_1,e_4,e_5,e_6\} \oplus \{e_4,e_7,e_8,e_9\} \oplus \{e_8,e_{10},e_{11}\} =$
$= \{e_2,e_3,e_5,e_6,e_7,e_9,e_{10},e_{11}\}$;
$w_2(e_9) = \gamma(\{e_1,e_4,e_5,e_6,e_8,e_{10},e_{11}\}) = s_1 \oplus s_2 \oplus s_2 \oplus s_3 \oplus s_2 \oplus s_4 \oplus s_2 \oplus s_6 \oplus s_3 \oplus s_5 \oplus s_4 \oplus$
$\oplus s_5 \oplus s_5 \oplus s_6 = s_1 \oplus s_5 = \{e_1,e_2,e_3\} \oplus \{e_8,e_{10},e_{11}\} = \{e_1,e_2,e_3,e_8,e_{10},e_{11}\}$;
$w_2(e_{10}) = \gamma(\{e_4,e_7,e_9,e_{10},e_{11}\}) = s_2 \oplus s_3 \oplus s_3 \oplus s_4 \oplus s_3 \oplus s_6 \oplus s_4 \oplus s_5 \oplus s_5 \oplus s_6 =$
$= s_2 \oplus s_3 = \{e_1,e_4,e_5,e_6\} \oplus \{e_4,e_7,e_8,e_9\} = \{e_1,e_5,e_6,e_7,e_8,e_9\}$;
$w_2(e_{11}) = \gamma(\{e_4,e_7,e_9,e_{10},e_{11}\}) = s_2 \oplus s_3 \oplus s_3 \oplus s_4 \oplus s_3 \oplus s_6 \oplus s_4 \oplus s_5 \oplus s_5 \oplus s_6 =$
$= s_2 \oplus s_3 = \{e_1,e_4,e_5,e_6\} \oplus \{e_4,e_7,e_8,e_9\} = \{e_1,e_5,e_6,e_7,e_8,e_9\}$.

Множество квалиразрезов 3-го уровня (яруса):

$w_3(e_1) = \gamma(\{e_2,e_3,e_5,e_6,e_7,e_9,e_{10},e_{11}\}) = s_1 \oplus s_4 \oplus s_1 \oplus s_6 \oplus s_2 \oplus s_4 \oplus s_2 \oplus s_6 \oplus s_3 \oplus s_4 \oplus$

$\oplus\ s_3 \oplus s_6 \oplus s_4 \oplus s_5 \oplus s_5 \oplus s_6 = \varnothing$;
$w_3(e_2) = \gamma(\{e_1,e_5,e_6,e_7,e_8,e_9\}) = s_1 \oplus s_2 \oplus s_2 \oplus s_4 \oplus s_2 \oplus s_6 \oplus s_3 \oplus s_4 \oplus s_3 \oplus s_5 \oplus s_3 \oplus s_6 =$
$= s_1 \oplus s_2 \oplus s_3 \oplus s_5 = \{e_1,e_2,e_3\} \oplus \{e_1,e_4,e_5,e_6\} \oplus \{e_4,e_7,e_8,e_9\} \oplus \{e_8,e_{10},e_{11}\} =$
$= \{e_2,e_3,e_5,e_6,e_7,e_9,e_{10},e_{11}\}$;
$w_3(e_3) = \gamma(\{e_1,e_5,e_6,e_7,e_8,e_9\}) = s_1 \oplus s_2 \oplus s_2 \oplus s_4 \oplus s_2 \oplus s_6 \oplus s_3 \oplus s_4 \oplus s_3 \oplus s_5 \oplus s_3 \oplus s_6 =$
$= s_1 \oplus s_2 \oplus s_3 \oplus s_5 = \{e_1,e_2,e_3\} \oplus \{e_1,e_4,e_5,e_6\} \oplus \{e_4,e_7,e_8,e_9\} \oplus \{e_8,e_{10},e_{11}\} =$
$= \{e_2,e_3,e_5,e_6,e_7,e_9,e_{10},e_{11}\}$;
$w_3(e_4) = \varnothing$;
$w_3(e_5) = \gamma(\{e_1,e_2,e_3,e_8,e_{10},e_{11}\}) = s_1 \oplus s_2 \oplus s_1 \oplus s_4 \oplus s_1 \oplus s_6 \oplus s_3 \oplus s_5 \oplus s_4 \oplus s_5 \oplus s_5 \oplus$
$\oplus s_6 = s_1 \oplus s_2 \oplus s_3 \oplus s_5 = \{e_1,e_2,e_3\} \oplus \{e_1,e_4,e_5,e_6\} \oplus \{e_4,e_7,e_8,e_9\} \oplus \{e_8,e_{10},e_{11}\} =$
$= \{e_2,e_3,e_5,e_6,e_7,e_9,e_{10},e_{11}\}$;
$w_3(e_6) = \gamma(\{e_1,e_2,e_3,e_8,e_{10},e_{11}\}) = s_1 \oplus s_2 \oplus s_1 \oplus s_4 \oplus s_1 \oplus s_6 \oplus s_3 \oplus s_5 \oplus s_4 \oplus s_5 \oplus s_5 \oplus$
$\oplus s_6 = s_1 \oplus s_2 \oplus s_3 \oplus s_5 = \{e_1,e_2,e_3\} \oplus \{e_1,e_4,e_5,e_6\} \oplus \{e_4,e_7,e_8,e_9\} \oplus \{e_8,e_{10},e_{11}\} =$
$= \{e_2,e_3,e_5,e_6,e_7,e_9,e_{10},e_{11}\}$;
$w_3(e_7) = \gamma(\{e_1,e_2,e_3,e_8,e_{10},e_{11}\}) = s_1 \oplus s_2 \oplus s_1 \oplus s_4 \oplus s_1 \oplus s_6 \oplus s_3 \oplus s_5 \oplus s_4 \oplus s_5 \oplus s_5 \oplus$
$\oplus s_6 = s_1 \oplus s_2 \oplus s_3 \oplus s_5 = \{e_1,e_2,e_3\} \oplus \{e_1,e_4,e_5,e_6\} \oplus \{e_4,e_7,e_8,e_9\} \oplus \{e_8,e_{10},e_{11}\} =$
$= \{e_2,e_3,e_5,e_6,e_7,e_9,e_{10},e_{11}\}$;
$w_3(e_8) = \gamma(\{e_2,e_3,e_5,e_6,e_7,e_9,e_{10},e_{11}\}) = s_1 \oplus s_4 \oplus s_1 \oplus s_6 \oplus s_2 \oplus s_4 \oplus s_2 \oplus s_6 \oplus s_3 \oplus s_4 \oplus$
$\oplus s_3 \oplus s_6 \oplus s_4 \oplus s_5 \oplus s_5 \oplus s_6 = \varnothing$;
$w_3(e_9) = \gamma(\{e_1,e_2,e_3,e_8,e_{10},e_{11}\}) = s_1 \oplus s_2 \oplus s_1 \oplus s_4 \oplus s_1 \oplus s_6 \oplus s_3 \oplus s_5 \oplus s_4 \oplus s_5 \oplus s_5 \oplus$
$\oplus s_6 = s_1 \oplus s_2 \oplus s_3 \oplus s_5 = \{e_1,e_2,e_3\} \oplus \{e_1,e_4,e_5,e_6\} \oplus \{e_4,e_7,e_8,e_9\} \oplus \{e_8,e_{10},e_{11}\} =$
$= \{e_2,e_3,e_5,e_6,e_7,e_9,e_{10},e_{11}\}$;
$w_3(e_{10}) = \gamma(\{e_1,e_5,e_6,e_7,e_8,e_9\}) = s_1 \oplus s_2 \oplus s_2 \oplus s_4 \oplus s_2 \oplus s_6 \oplus s_3 \oplus s_4 \oplus s_3 \oplus s_5 \oplus s_3 \oplus s_6 =$
$= s_1 \oplus s_2 \oplus s_3 \oplus s_5 = \{e_1,e_2,e_3\} \oplus \{e_1,e_4,e_5,e_6\} \oplus \{e_4,e_7,e_8,e_9\} \oplus \{e_8,e_{10},e_{11}\} =$
$= \{e_2,e_3,e_5,e_6,e_7,e_9,e_{10},e_{11}\}$;
$w_3(e_{11}) = \gamma(\{e_1,e_5,e_6,e_7,e_8,e_9\}) = s_1 \oplus s_2 \oplus s_2 \oplus s_4 \oplus s_2 \oplus s_6 \oplus s_3 \oplus s_4 \oplus s_3 \oplus s_5 \oplus s_3 \oplus s_6 =$
$= s_1 \oplus s_2 \oplus s_3 \oplus s_5 = \{e_1,e_2,e_3\} \oplus \{e_1,e_4,e_5,e_6\} \oplus \{e_4,e_7,e_8,e_9\} \oplus \{e_8,e_{10},e_{11}\} =$
$= \{e_2,e_3,e_5,e_6,e_7,e_9,e_{10},e_{11}\}$.

Множество квалиразрезов 4-го уровня (яруса):

$w_4(e_1) = \varnothing$;
$w_4(e_2) = \gamma(\{e_2,e_3,e_5,e_6,e_7,e_9,e_{10},e_{11}\}) = \varnothing$;
$w_4(e_3) = \gamma(\{e_2,e_3,e_5,e_6,e_7,e_9,e_{10},e_{11}\}) = \varnothing$;
$w_4(e_4) = \varnothing$;
$w_4(e_5) = \gamma(\{e_2,e_3,e_5,e_6,e_7,e_9,e_{10},e_{11}\}) = \varnothing$
$w_4(e_6) = \gamma(\{e_2,e_3,e_5,e_6,e_7,e_9,e_{10},e_{11}\}) = \varnothing$;
$w_4(e_7) = \gamma(\{e_2,e_3,e_5,e_6,e_7,e_9,e_{10},e_{11}\}) = \varnothing$;
$w_4(e_8) = \varnothing$;
$w_4(e_9) = \gamma(\{e_2,e_3,e_5,e_6,e_7,e_9,e_{10},e_{11}\}) = \varnothing$;
$w_4(e_{10}) = \gamma(\{e_2,e_3,e_5,e_6,e_7,e_9,e_{10},e_{11}\}) = \varnothing$;
$w_4(e_{11}) = \gamma(\{e_2,e_3,e_5,e_6,e_7,e_9,e_{10},e_{11}\}) = \varnothing$.

Множество квалиразрезов 4-го уровня (яруса) пусто:

Цепочка порожденных реберных разрезов для ребра $e_5$ имеет вид:

$w_0(e_5) = s_2 \oplus s_4 = \{e_1,e_2,e_4,e_6,e_7,e_{10}\}$;
$w_1(e_5) = \gamma(\{e_1,e_2,e_4,e_6,e_7,e_{10}\}) = s_1 \oplus s_2 \oplus s_1 \oplus s_4 \oplus s_2 \oplus s_3 \oplus s_2 \oplus s_6 \oplus s_3 \oplus s_4 \oplus s_4 \oplus s_5 =$
$= s_2 \oplus s_4 \oplus s_5 \oplus s_6 = \{e_1,e_4,e_5,e_6\} \oplus \{e_2,e_5,e_7,e_{10}\} \oplus \{e_8,e_{10},e_{11}\} \oplus \{e_3,e_6,e_9,e_{11}\} =$
$= \{e_1,e_2,e_3,e_4,e_7,e_8,e_9\}$;

$w_2(e_5) = \gamma(\{e_1,e_2,e_3,e_4,e_7,e_8,e_9\}) = s_1 \oplus s_2 \oplus s_1 \oplus s_4 \oplus s_1 \oplus s_6 \oplus s_2 \oplus s_3 \oplus s_3 \oplus s_4 \oplus s_3 \oplus$
$\oplus s_5 \oplus s_3 \oplus s_6 = s_1 \oplus s_5 = \{e_1,e_2,e_3\} \oplus \{e_8,e_{10},e_{11}\} = \{e_1,e_2,e_3,e_8,e_{10},e_{11}\};$
$w_3(e_5) = \gamma(\{e_1,e_2,e_3,e_8,e_{10},e_{11}\}) = s_1 \oplus s_2 \oplus s_1 \oplus s_4 \oplus s_1 \oplus s_6 \oplus s_3 \oplus s_5 \oplus s_4 \oplus s_5 \oplus s_5 \oplus$
$\oplus s_6 = s_1 \oplus s_2 \oplus s_3 \oplus s_5 = \{e_1,e_2,e_3\} \oplus \{e_1,e_4,e_5,e_6\} \oplus \{e_4,e_7,e_8,e_9\} \oplus \{e_8,e_{10},e_{11}\} =$
$= \{e_2,e_3,e_5,e_6,e_7,e_9,e_{10},e_{11}\};$
$w_4(e_5) = \gamma(\{e_2,e_3,e_5,e_6,e_7,e_9,e_{10},e_{11}\}) = s_1 \oplus s_4 \oplus s_1 \oplus s_6 \oplus s_2 \oplus s_4 \oplus s_2 \oplus s_6 \oplus s_3 \oplus s_4 \oplus s_3 \oplus$
$\oplus s_6 \oplus s_4 \oplus s_5 \oplus s_5 \oplus s_6 = \varnothing.$

Цепочку порожденных реберных разрезов можно представить в виде:

$w_0(e_5) = s_2 \oplus s_4 \rightarrow w_1(e_5) = s_2 \oplus s_4 \oplus s_5 \oplus s_6 \rightarrow w_2(e_5) = s_1 \oplus s_5 \rightarrow$
$\rightarrow w_3(e_5) = s_1 \oplus s_2 \oplus s_3 \oplus s_5 \rightarrow w_4(e_5) = \varnothing.$

И так для каждого ребра. Уровень (ярус) $W_l(G)$ будем обозначать буквой $l$. Матрица спектра реберных разрезов $W_s$, записанная как перечисление суграфов подпространства S(G) графа $G_2$, имеет вид:

| ребро | $l_0$ | $l_1$ | $l_2$ | $l_3$ |
|---|---|---|---|---|
| $e_1$ | $\{e_2,e_3,e_4,e_5,e_6\}$ | $\{e_1,e_5,e_6,e_7,e_8,e_9\}$ | $\{e_2,e_3,e_5,e_6,e_7,e_9,e_{10},e_{11}\}$ | $\varnothing$ |
| $e_2$ | $\{e_1,e_3,e_5,e_7,e_{10}\}$ | $\{e_2,e_3,e_4,e_5,e_6\}$ | $\{e_1,e_5,e_6,e_7,e_8,e_9\}$ | $\{e_2,e_3,e_5,e_6,e_7,e_9,e_{10},e_{11}\}$ |
| $e_3$ | $\{e_1,e_2,e_6,e_9,e_{11}\}$ | $\{e_2,e_3,e_4,e_5,e_6\}$ | $\{e_1,e_5,e_6,e_7,e_8,e_9\}$ | $\{e_2,e_3,e_5,e_6,e_7,e_9,e_{10},e_{11}\}$ |
| $e_4$ | $\{e_1,e_5,e_6,e_7,e_8,e_9\}$ | $\{e_2,e_3,e_5,e_6,e_7,e_9,e_{10},e_{11}\}$ | $\varnothing$ | $\varnothing$ |
| $e_5$ | $\{e_1,e_2,e_4,e_6,e_7,e_{10}\}$ | $\{e_1,e_2,e_3,e_4,e_7,e_8,e_9\}$ | $\{e_1,e_2,e_3,e_8,e_{10},e_{11}\}$ | $\{e_2,e_3,e_5,e_6,e_7,e_9,e_{10},e_{11}\}$ |
| $e_6$ | $\{e_1,e_3,e_4,e_5,e_9,e_{11}\}$ | $\{e_1,e_2,e_3,e_4,e_7,e_8,e_9\}$ | $\{e_1,e_2,e_3,e_8,e_{10},e_{11}\}$ | $\{e_2,e_3,e_5,e_6,e_7,e_9,e_{10},e_{11}\}$ |
| $e_7$ | $\{e_2,e_4,e_5,e_8,e_9,e_{10}\}$ | $\{e_1,e_4,e_5,e_6,e_8,e_{10},e_{11}\}$ | $\{e_1,e_2,e_3,e_8,e_{10},e_{11}\}$ | $\{e_2,e_3,e_5,e_6,e_7,e_9,e_{10},e_{11}\}$ |
| $e_8$ | $\{e_4,e_7,e_9,e_{10},e_{11}\}$ | $\{e_1,e_5,e_6,e_7,e_8,e_9\}$ | $\{e_2,e_3,e_5,e_6,e_7,e_9,e_{10},e_{11}\}$ | $\varnothing$ |
| $e_9$ | $\{e_3,e_4,e_6,e_7,e_8,e_{11}\}$ | $\{e_1,e_4,e_5,e_6,e_8,e_{10},e_{11}\}$ | $\{e_1,e_2,e_3,e_8,e_{10},e_{11}\}$ | $\{e_2,e_3,e_5,e_6,e_7,e_9,e_{10},e_{11}\}$ |
| $e_{10}$ | $\{e_2,e_5,e_7,e_8,e_{11}\}$ | $\{e_4,e_7,e_9,e_{10},e_{11}\}$ | $\{e_1,e_5,e_6,e_7,e_8,e_9\}$ | $\{e_2,e_3,e_5,e_6,e_7,e_9,e_{10},e_{11}\}$ |
| $e_{11}$ | $\{e_3,e_6,e_8,e_9,e_{10}\}$ | $\{e_4,e_7,e_9,e_{10},e_{11}\}$ | $\{e_1,e_5,e_6,e_7,e_8,e_9\}$ | $\{e_2,e_3,e_5,e_6,e_7,e_9,e_{10},e_{11}\}$ |

Таким образом, каждый элемент матрицы $w_l(e_i) \in W_s$, характеризуется двумя параметрами: строкой с номером $i$ и столбцом с номером $l$. Каждая строка характеризуется подмножеством реберных разрезов для выбранного ребра $e_i$, где текущее ребро $e_i$ рассматривается как аргумент реберного разреза $w_l(e_i)$.

Заметим, что при распозновании изоморфизма с использованием матрицы смежностей графа приходится одновременно переставлять строки и столбцы матрицы A(G). Структура матрицы $W_s(G)$, в отличие от матрицы смежностей графа, позволяет рассматривать, переставлять, оценивать, и сравнивать между собой отдельно только строки, сохраняя оценку столбцов при перестановке строк.

## 2.5. Вес ребра

Построение спектра реберных разрезов графа позволяет определить количественную характеристику участия ребра во множестве суграфов $W_s$ – вес ребра. Влияние такой количественной характеристики ребра на внутреннюю структуру графа усиливается, если рассматривать не только множество базовых реберных разрезов графа, но и множество всех порожденных преобразованием $\gamma(w_{k-1}(e_i))$ реберных разрезов графа. Введение такой количественной характеристики позволяет в дальнейшем построить интегральный инвариант графа, используя множество реберных разрезов как инвариант.

**Определение 2.7.** Весом ребра $\xi_j(e_i)$ (греческая буква «кси») называется количество суграфов в уровне $l_j$ или строке спектра реберных разрезов $W_s$ с участием ребра $e_i$, где $i$ изменяется от 1 до $m$ ($m$ - число ребер в графе).

Рассмотрим граф $G_2$. Каждой строке $w(e_i)$ матрицы $W_s$ можно поставить в соответствие кортеж весов $\varepsilon(w(e_i))$ (греческая буква «эпсилон»), элементы которого называются *строчным весом ребра* $e_i$ в графе и определяются количеством суграфов строки с участием ребра $e_i$:

$\varepsilon(\mathrm{w}(e_1)) = <1,2,2,1,3,3,2,1,2,1,1>$;
$\varepsilon(\mathrm{w}(e_2)) = <2,2,3,1,4,3,3,1,2,2,1>$;
$\varepsilon(\mathrm{w}(e_3)) = <2,3,2,1,3,4,2,1,3,1,2>$;
$\varepsilon(\mathrm{w}(e_4)) = <1,1,1,0,2,2,2,1,2,1,1>$;
$\varepsilon(\mathrm{w}(e_5)) = <3,4,3,2,1,2,3,2,2,3,2>$;
$\varepsilon(\mathrm{w}(e_6)) = <3,3,4,2,2,1,2,2,3,2,3>$;
$\varepsilon(\mathrm{w}(e_7)) = <2,3,2,2,3,2,1,3,2,4,3>$;
$\varepsilon(\mathrm{w}(e_8)) = <1,1,1,1,2,2,3,1,3,2,2>$;
$\varepsilon(\mathrm{w}(e_9)) = <2,2,3,2,2,3,2,3,1,3,4>$;
$\varepsilon(\mathrm{w}(e_{10})) = <1,2,1,1,3,2,4,2,3,2,3>$;
$\varepsilon(\mathrm{w}(e_{11})) = <1,1,2,1,2,3,3,2,4,3,2>$.

Каждому столбцу матрицы $w(l_j)$ можно поставить в соответствие кортеж $\xi(w(l_j))$ элементы которого (веса ребер) характеризуют количество суграфов уровня $l_j$ с участием ребра $e_i$.

$\xi(\mathrm{w}(l_0)) = <5,5,5,6,6,6,6,5,6,5,5>$;
$\xi(\mathrm{w}(l_1)) = <6,5,5,8,7,7,7,6,7,5,5>$;
$\xi(\mathrm{w}(l_2)) = <8,6,6,0,6,6,6,8,6,6,6>$;
$\xi(\mathrm{w}(l_3)) = <0,8,8,0,8,8,8,0,8,8,8>$.

Построим суммарный кортеж графа G как сумму элементов кортежей $\varepsilon_{\mathrm{w}}(G) = \sum_{i=1}^{m} \varepsilon(\mathrm{w}(e_i))$ или $\xi_{\mathrm{w}}(G) = \sum_{j=1}^{k} \xi(\mathrm{w}(l_j))$:

$\varepsilon_{\mathrm{w}}(G) = <\varepsilon(\mathrm{w}(e_1)), \varepsilon(\mathrm{w}(e_2)), \varepsilon(\mathrm{w}(e_3)), \ldots, \varepsilon(\mathrm{w}(e_m))>$;
$\xi_{\mathrm{w}}(G) = <\xi(\mathrm{w}(l_0)), \xi(\mathrm{w}(l_1)), \xi(\mathrm{w}(l_2)), \ldots, \xi(\mathrm{w}(l_k))>$.

С учетом постронгтя матрицы $W_s$, можно записать

$$\varepsilon_{\mathrm{w}}(G) = \xi_{\mathrm{w}}(G) \qquad (2.13)$$

Для графа $G_2$, кортеж $\varepsilon_{\mathrm{w}}(G_2) = <19,24,24,14,27,27,27,19,27,24,24>$.

<1,2,2,1,3,3,2,1,2,1,1>
\+ <2,2,3,1,4,3,3,1,2,2,1>
\+ <2,3,2,1,3,4,2,1,3,1,2>
\+ <1,1,1,0,2,2,2,1,2,1,1>
\+ <3,4,3,2,1,2,3,2,2,3,2>
\+ <3,3,4,2,2,1,2,2,3,2,3>
\+ <2,3,2,2,3,2,1,3,2,4,3>
\+ <1,1,1,1,2,2,3,1,3,2,2>
\+ <2,2,3,2,2,3,2,3,1,3,4>
\+ <1,2,1,1,3,2,4,2,3,2,3>
\+ <1,1,2,1,2,3,3,2,4,3,2>
= <19,24,24,14,27,27,27,19,27,24,24

или как результат сложения кортежей $\xi$ (w($l_j$)) для столбцов:

  <5,5,5,6,6,6,6,5,6,5,5>
\+ <6,5,5,8,7,7,7,6,7,5,5>
\+ <8,6,6,0,6,6,6,8,6,6,6>
\+ <0,8,8,0,8,8,8,0,8,8,8>
= <19,24,24,14,27,27,27,19,27,24,24>

Для установления изоморфизма графов G и H необходим инвариант, на основе которого мы можем установить изоморфизм. Пусть $f$ – функция, относящая каждому элементу графа G некоторый элемент $f$(G) из множества M произвольной природы. Эту функцию называют инвариантом, если на изоморфных графах G и H (G $\equiv$ H) ее значения совпадают, т.е. для любых G и H

$$G \equiv H \rightarrow f(G) = f(H). \qquad (2.14)$$

Инвариант $f$ называется полным, если для любых G и H

$$f(G) = f(H) \rightarrow G \equiv H. \qquad (2.15)$$

Объединяя оба определения, назовем полным инвариантом графа такую функцию $f$(G) (со значениями в произвольном множестве), для которой f(G) = $f$(H) тогда и только тогда, когда G $\equiv$ H.

Понятие веса ребра позволяет нам поставить в соответчтвие каждому столбцу - вектор размера $m$, с координатами веса для каждого ребра

$\xi(\mathrm{w}(l_0)) = <\xi_0(e_1), \xi_0(e_2), ..., \xi_0(e_m)>$.

Например, для нашего случая:

$\xi(\mathrm{w}(l_0))$ = <5,5,5,6,6,6,6,5,6,5,5> - кортеж ребер базовых реберных разрезов графа $G_2$;
$\xi(\mathrm{w}(l_1))$ = <6,5,5,8,7,7,7,6,7,5,5> - кортеж ребер реберных разрезов графа $G_2$ уровня 1;
$\xi(\mathrm{w}(l_2))$ = <8,6,6,0,6,6,6,8,6,6,6> - кортеж ребер реберных разрезов графа $G_2$ уровня 2;
$\xi(\mathrm{w}(l_3))$ = <0,8,8,0,8,8,8,0,8,8,8> - кортеж ребер реберных разрезов графа $G_2$ уровня 3.

Естественно, что такой вектор удобно записывать в виде кортежа, где вес соответствующего ребра определяется его местоположением.

В нашем случае, ребро $e_5$, в базовом реберном разрезе имеет вес 6, а ребро $e_8$ в реберном разрезе уровня 3, имеет вес равный 0.

Для удобства сравнения векторов, лучше всего расположить веса ребер из кортежа $\varsigma_{lw}$ (G) по неубыванию. Поэтому можно записать вектор в виде упорядоченного значения координат. Например, для нашего примера

$F_w(\xi_0(G_2))$ = (5,5,5,5,5,5,6,6,6,6,6) или (6×5,5×6);
$F_w(\xi_1(G_2))$ = (5,5,5,5,6,6,7,7,7,7,8) или (4×5,2×6,4×7,1×8);
$F_w(\xi_2(G_2))$ = (0,6,6,6,6,6,6,6,6,8,8) или (1×0,8×6,2×8)>;
$F_w(\xi_3(G_2))$ = <0,0,0,8,8,8,8,8,8,8,8) или (3×0,8×8).

Вектор весов ребер яруса является первой частью упорядоченного векторного инварианта. Теперь можно задать вес каждой вершины яруса как сумму весов инцидентных ребер.

**Определение 2.8.** Весом вершины $\zeta_l$($v_i$) в ярусе $l$ (греческая буква «дзета»), будем называть сумму весов инцидентных вершине ребер, где $i$ номер инцидентного ребра в диапазоне от 1 до $n$.

Например, для графа $G_2$:

$\zeta_0(v_1) = \xi_0(e_1)+\xi_0(e_2)+\xi_0(e_3) \rightarrow 5 + 5 + 5 = 15;$
$\zeta_0(v_2) = \xi_0(e_1)+\xi_0(e_4)+\xi_0(e_5)+\xi_0(e_6) \rightarrow 5 + 6 + 6 +6 = 23;$
$\zeta_0(v_3) = \xi_0(e_4)+\xi_0(e_7)+\xi_0(e_8)+\xi_0(e_9) \rightarrow 6 + 6 + 5 + 6 = 23;$
$\zeta_0(v_4) = \xi_0(e_2)+\xi_0(e_5)+\xi_0(e_7)+\xi_0(e_{10}) \rightarrow 5 + 6 + 6 +5 = 22;$
$\zeta_0(v_5) = \xi_0(e_8)+\xi_0(e_{10})+\xi_0(e_{11}) \rightarrow 5 + 5 + 5 = 15;$
$\zeta_0(v_6) = \xi_0(e_3)+\xi_0(e_6)+\xi_0(e_9)+\xi_0(e_{11}) \rightarrow 5 + 6 + 6 +5 = 22.$

Построим кортеж $\zeta_l$(G) для вершин графа яруса $l$:

$$\zeta_l(G) = <\zeta_l(v_1), \zeta_l(v_2), \dots, \zeta_l(v_n)>. \quad (2.16)$$

где $\zeta_l(v_j)$ – вес соответствующей вершины яруса $l$, где $j$ = (1,2,…,$n$).

$\zeta_0$($G_2$) = <15,23,23,22,15,22>.

Взяв за основу кортеж $\zeta_l$(G), построим вектор весов вершин $F_w(\zeta_0(G_2))$ для базавого уровня, распологая веса вершин по неубыванию $F_w(\zeta_0(G_2))$ = (15,15,22,22,23,23).

Вектор весов вершин является второй частью векторного инварианта яруса спектра реберных разрезов. Векторный инвариант спектра реберных разрезов графа G яруса $l$ запишем в виде:

$$F_w(\xi_l(G_2)) \ \& \ F_w(\zeta_l(G_2)) \quad (2.17)$$

Для графа $G_2$ инвариант базовых реберных разрезов имеет вид:

$F(\xi_0(G_2))$ & $F(\zeta_0(G_2))$ = (6×5,5×6) & (2×15,2×22,2×23).

**Определение 2.9.** ***Векторным инвариантом*** будем называть вектор с координатами равными весу его элеиентов (ребер, вершин). Векторный инвариант представляется в виде двух записей. Одна запись в виде кортежа весов элементов, другая запись в виде упорядоченной записи по невозрастанию весов элементов.

Следует различать векторные инварианты для спектра элементов, и векторные инварианты для уровней спектра. В записи уровневых векторных инвариантах обязательно должет присутствовать номер яруса.

Количество элементарных операций для построения суграфа как базового элемента реберных разрезов, можно определить как сумму по модулю 2 базисных характеристических векторов размером $m$. Количество элементарных операций для получения суграфа спектра реберных разрезов можно определить как результат сложения всех базисных суграфов $m\times m = m^2$. Тогда количество суграфов в ярусе определяется как $m\times m^2 = m^3$. В свою очередь, размер матрицы $W_S$ можно определить как количество ярусов, умноженное на количество суграфов в ярусе $q\times m^3$. Итого, для построения матрицы $W_S(G)$ нужно применить $q\times m^3$ элементарных операций. Для определения веса ребра нужно перебрать все элементы спектра реберных разрезов графа, числом операций равным $q\times m^3$. Общее количество элементарных операций можно определить как сумму операций вычисления элементов матрицы $W_s(G)$ и вычисление веса ребра $q\times m^3+q\times m^3 = 2q\times m^3$.

Определим состав ярусного векторного инварианта реберного разреза графа:

$$IS_l(G) = F_w(\xi_l(G_2)) \ \& \ F_w(\zeta_l(G_2)) \qquad (2.18)$$

Опишем алгоритм для определения векторного инварианта реберных разрезов графа, основанный на понятии спектра реберных разрезов.

**Инициализация.** Задан граф G.

**шаг 1. [Построение матрицы смежностей реберного графа**]. Строим линейный оператор в виде матрицы смежностей реберного графа A(L(G)).

**шаг 2. [Определяем реберные разрезы для уровней спектра реберных разрезов]**. Применяем операцию преобразования $\gamma(w_{k-1}(e_i))$ для построения линейного оператора $W_\lambda^k$ графа, получаем реберные разрезы уровня $k$. Если множество реберных разрезов уровня $k$ не пусто, то идем на шаг 2. Иначе идем на шаг 3.

**шаг 3. [Определение количество уровней]**. Определяем количество уровней $l_S$.

**шаг 4. [Определение кортежа весов уровня]**. Создаем кортеж весов ребер для всех уровней матрицы $W_S$ графа.

*Конец работы алгоритма.*

В результате работы алгоритма, будут построены уровневые векторные инварианты реберных разрезов для ребер и уровневые векторные инварианты реберных разрезов для вершин:

$\zeta(\mathrm{w}(l_0)) =$ <15,23,23,22,15,22>;

$\zeta(\mathrm{w}(l_1)) =$ <16,28,28,24,16,24>;

$\zeta(\mathrm{w}(l_2)) =$ <20,20,20,24,20,24>;

$\zeta(\mathrm{w}(l_2)) =$ <16,16,16,32,16,32>.

Вектор весов вершин для спектра реберных разрезов графа $G_2$ имеет вид:

$\zeta(v_1) = \xi(e_1)+\xi(e_2)+\xi(e_3) \to 19 + 24 + 24 = 67;$

$\zeta(v_2) = \xi(e_1)+\xi(e_4)+\xi(e_5)+\xi(e_6) \to 19 + 14 + 27 +27 = 87;$

$\zeta(v_3) = \xi(e_4)+\xi(e_7)+\xi(e_8)+\xi(e_9) \to 14 + 27 + 19 +27 = 87;$

$\zeta(v_4) = \xi(e_2)+\xi(e_5)+\xi(e_7)+\xi(e_{10}) \to 24 + 27 + 27 +24 = 102;$

$\zeta(v_5) = \xi(e_8)+\xi(e_{10})+\xi(e_{11}) \to 19 + 24 + 24 = 67;$

$\zeta(v_6) = \xi(e_3)+\xi(e_6)+\xi(e_9)+\xi(e_{11}) \to 24 + 27 + 27 +24 = 102.$

В результате построения, получим систему векторных инвариантов для каждого яруса:

$$\begin{cases} \mathrm{IS}_0(\mathrm{G}) = \mathrm{F_w}(\xi_0(\mathrm{G}_2))\ \&\ \mathrm{F_w}(\zeta_0(\mathrm{G}_2)) = (6\times 5{,}5\times 6)\ \&\ (2\times 15{,}2\times 22{,}2\times 23); \\ \mathrm{IS}_1(\mathrm{G}) = \mathrm{F_w}(\xi_1(\mathrm{G}_2))\ \&\ \mathrm{F_w}(\zeta_1(\mathrm{G}_2)) = (4\times 5{,}2\times 6{,}4\times 7{,}8)\ \&\ (2\times 16{,}2\times 24{,}2\times 28); \\ \mathrm{IS}_2(\mathrm{G}) = \mathrm{F_w}(\xi_2(\mathrm{G}_2))\ \&\ \mathrm{F_w}(\zeta_2(\mathrm{G}_2)) = (0{,}8\times 6{,}2\times 8)\ \&\ (4\times 20{,}2\times 24); \\ \mathrm{IS}_3(\mathrm{G}) = \mathrm{F_w}(\xi_3(\mathrm{G}_2))\ \&\ \mathrm{F_w}(\zeta_3(\mathrm{G}_2)) = (3\times 0{,}8\times 8)\ \&\ (4\times 16{,}2\times 32). \end{cases}$$

И векторный инвариант спектра реберных разрезов графа $G_2$ представленный в виде:

$\mathrm{IS}(G_2) = \mathrm{F_w}(\xi(\mathrm{G}_2))\ \&\ \mathrm{F_w}(\zeta(\mathrm{G}_2)) = (14{,}2\times 19{,}4\times 24{,}4\times 27)\ \&\ (2\times 67{,}2\times 87{,}2\times 102).$

На данном этапе построения структур, мы не можем говорить о полном инварианте графа. Мы можем говорить только об изоморфизме спектра реберных разрезов графов G и H в результате сравнения.

## 2.6. Гиперграфы реберных разрезов

Удобным представлением элементов спектра реберных разрезов являются гиперграфы. Гиперграфы позволяют представить не только отдельные элементы спектра реберных разрезов, но также строки и столбцы матрицы реберных разрезов.

На рис. 2.17 показан гиперграф отдельного элемента $w_3(e_1)$ спектра реберных разрезов графа $G_2$. Множество вершин гиперграфа характеризует ребра графа $G_2$, а ребро гиперграфа характеризует отдельный элемент спектра реберных разрезов. На рис. 2.18 представлен гиперграф строки для ребра $e_1$ спектра реберных разрезов графа $G_2$. Множество вершин гиперграфа характеризует ребра графа $G_2$, а ребра гиперграфа характеризуют реберные разрезы порожденные ребром $e_1$. На рис. 2.19 представлен гиперграф столбца $l_1$ спектра реберных разрезов графа $G_2$. Здесь множество вершин гиперграфа характеризует ребра графа $G_2$, а ребра гиперграфа характеризуют суграфы спектра реберных разрезов уровня 1.

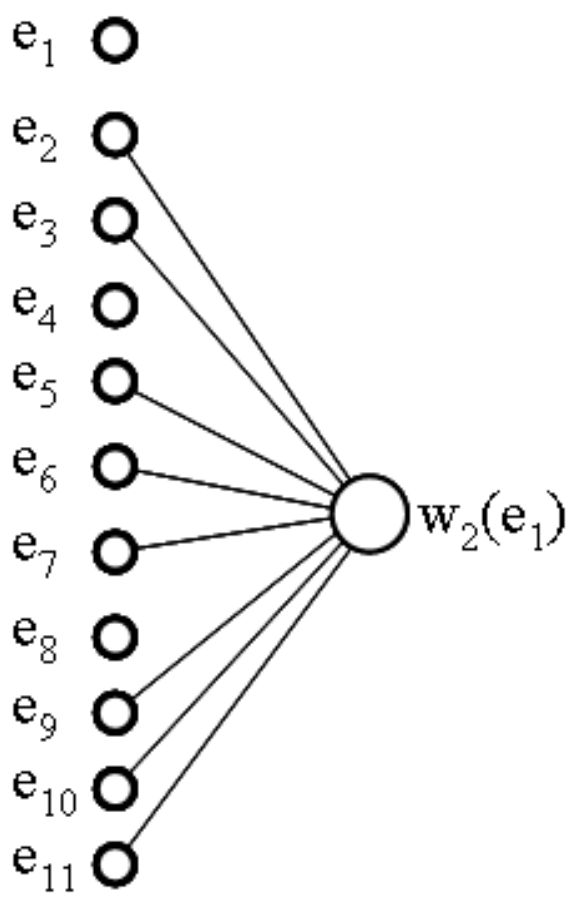

Рис. 2.17. Гиперграф элемента $w_2(e_1)$.

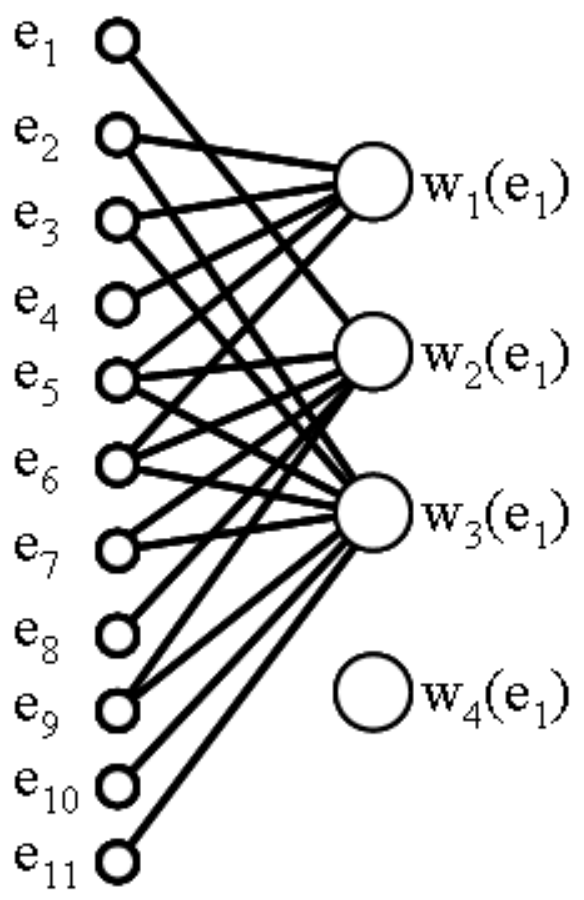

Рис. 2.18. Гиперграф строки ребра $e_1$.

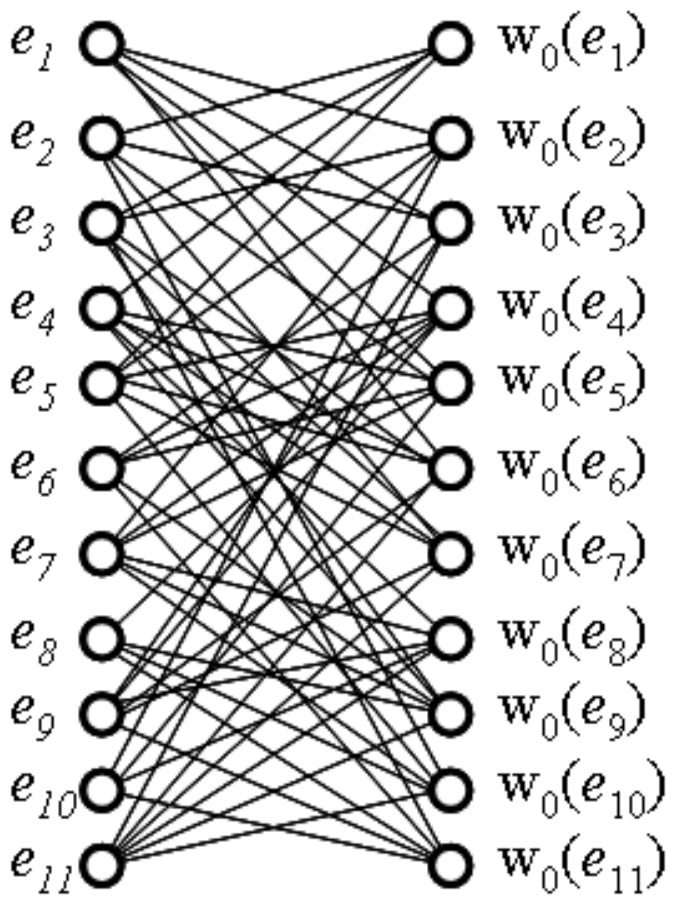

Рис. 2.19. Гиперграф столбца $l_0$.

## 2.7. Программа определения спектра реберных разрезов

```
program Raschet54;

type
        TMasy = array[1..100000] of integer;
        TMass = array[1..400000] of integer;
var
        F1,F2 : text;
        i,ii,j,jj,iii,jjj,p11,K11,K12,Nv,M,Np,My,Prikaz : integer;
        KKK,KKK0 : integer;
        Masy1: TMasy;
        Mass1: TMass;
        Massi: TMass;
        Masy2: TMasy;
        Mass2 : TMass;
        Masy3: TMasy;
        Mass3 : TMass;
        Masy4: TMasy;
        Mass4 : TMass;
        Masy5: TMasy;
        Mass5 : TMass;
        Masy6 : TMasy;
        Mass6 : TMass;
        MasMdop : TMasy;
        MasKol : TMasy;
        MasWesReb : TMasy;
        MasWesVer : TMasy;
{*****************************************************************}
 procedure  Shell(var N : integer;
                 var A : TMasy);
{       процедура Шелла для упорядочивания элементов                }
{                                                                   }
{     N - количество элементов в массиве;                           }
{     A - сортируемый массив;                                       }
var D,Nd,I,J,L,X : integer;
label 1,2,3,4,5;
```

```
begin
  D:=1;
1:D:=2*D;
  if D<=N then goto 1;
2:D:=D-1;
  D:=D div 2;
  if D=0 then goto 5;
  Nd:=N-D;
  for I:=1 to Nd do
  begin
    J:=I;
3:  L:=J+D;
    if A[L]>=A[J] then goto 4;
    X:=A[J];
    A[J]:=A[L];
    A[L]:=X;
    J:=J-D;
    if J>0 then goto 3;
4:end;
  goto 2;
5:end;{Shell}
{*******************************************************************}
procedure FormIncide(var Nv : integer;
                     var My : TMasy;
                     var Ms : TMass;
                     var Ms3 : TMass);
 { Nv - количество вершин в графе;                                   }
 { My  - массив указателей для матрицы смежностей;                   }
 { Ms  - массив элементов матрицы смежностей;                        }
 { Ms3 - массив элементов матрицы инциденций.                        }
 {   Формируется матрица инциденций графа в массиве                  }
 {   Ms3.                                                            }
 {                                                                   }
        var I,J,K,NNN,P,M,L : integer;
        begin
 {      инициализация                                                }
        NNN:=My[Nv+1]-1;
        K:=0;
        for J:= 1 to NNN do Ms3[J]:= 0;
 {   определение номера элемента                      }
        for I:= 1 to Nv do
         for M:= My[I] to My[I+1]-1 do
          if Ms3[M]=0 then
          begin
            P:=Ms[M];
            K:=K+1;
            Ms3[M]:=K;
            for L:=My[P] to My[P+1]-1 do
             if Ms[L]=I then Ms3[L]:=K;
          end;
        end; {FormIncide}
{************************************************************}
 procedure FormRebRas1(var Vt,Vs,Dl : integer;
                       var Masy1 : TMasy;
                       var Mass1 : TMass;
                       var Massi  : TMass;
                       var MasMdop : TMasy);
{ Vt -  номер текущей вершины;                                   }
{ Vs -  номер смежной вершины;                                   }
{ Masy1 -  массив указателей для маттрицы смежностей графа;      }
{ Mass1 -  массив элементов матрицы смежностей;                  }
{ Massi -  массив элементов матрицы инциденций;                  }
{ MasMdop -  вспомогательный массив.                             }
```

```
{                                                                    }
{   Данная процедура формирует множество базовых разрезов            }
{   для ребра с вершинами (Vt,Vs)                                    }
 var i,j,K11 : integer;
 begin
       K11:=0;
       for i:=Masy1[Vt] to Masy1[Vt+1]-1 do
       begin {1}
          K11:=K11+1;
          MasMdop[K11]:=Massi[i];
          {writeln(F2,'MasMdop[',K11,'] = ',MasMdop[k11]); }
       end; {1}
       for i:=Masy1[Vs] to Masy1[Vs+1]-1 do
       begin {2}
          K11:=K11+1;
          MasMdop[K11]:=Massi[i];
          {writeln(F2,'MasMdop[',K11,'] = ',MasMdop[k11]);}
       end; {2}
{В массиве MasMdop находится базовый реберный разрез ребра (Vt,Vs)    }
       for i:=1 to K11-1 do
       begin {3}
         for j:=i+1 to K11 do
         if MasMdop[j]=MasMdop[i] then
         begin {4}
            MasMdop[j]:=0;
            MasMdop[i]:=0;
        end; {4}
       end; {3}
       Dl:=0;
       for i:=1 to K11 do
       if MasMdop[i]<>0 then
       begin  {5}
        Dl:=Dl+1;
        MasMdop[Dl]:=MasMdop[i];
        {writeln(F2,'MasMdop[',Dl,'] = ',MasMdop[Dl]); }
       end;  {5}
       {writeln(F2,'Dl = ',Dl);}
       Shell(Dl,MasMdop);
 end; {FormRebRas1}
{*************************************************************}
 procedure FormRasReb2(var Nv,K11 : integer;
                  var Masy1 : TMasy;
                  var Mass1 : TMass;
                  var Masy2 : TMasy;
                  var Mass2 : TMass;
                  var Massi : TMass;
                  var MasMdop : TMasy);
{Процедура создания множества базовых реберных разрезов графа        }
{                                                                    }
{ Nv - количество вершин в графе;                                    }
{ Masy1  - массив указателей для матрицы смежностей;                 }
{ Mass1  - массив элементов матрицы смежностей;                      }
{ Masy2 - массив указателей для множества базовых разрезов;          }
{ Mass2 - массив элементов множеств базовых реберных разрезов;       }
{ MasMdop -  вспомогательный массив для хранения суграфа.            }
var i,j,ii,Dl,Vt,Vs,K20 : integer;
begin
       K11:=0;
       Masy2[1]:=1;
       for i:= 1 to Nv do
       begin {1}
          Vt:=i;
          {writeln(F2,'Vt = ',Vt); }
```

```
          for j:= Masy1[i] to Masy1[i+1]-1 do
          begin {2}
             Vs:=Mass1[j];
             {writeln(F2,'Vs = ',Vs); }
             if Vs>Vt then
             begin {3}
               FormRebRas1(Vt,Vs,Dl,Masy1,Mass1,Massi,MasMdop);
               K11:=K11+1;
               Masy2[K11+1]:=Masy2[K11]+Dl;
               for ii:=1 to Dl do
               begin  {4}
                  K20:=Masy2[K11]-1+ii;
                  Mass2[K20]:=MasMdop[ii];
               end;   {4}
             end;   {3}
         end;   {2}
       end;   {1}
end;{ FormRasReb2}
{************************************************************}
 procedure Perezapis1(var M : integer;
                var Masy2 : TMasy;
                var Mass2 : TMass;
                var Masy3 : TMasy;
                var Mass3 : TMass);
{Процедура перезаписи массива                                    }
{                                                                }
{ MY - номер уровня;                                             }
{                                                                }
{ Masy2 - массив указателей для множества реберных разрезов 1-уровня;}
{ Mass2 - массив элементов множеств беберных разрезов 1-уровня;       }
{ Masy3 - массив указателей для множества реберных разрезов т-уровня; }
{ Mass3 - массив элементов множеств беберных разрезов n-уровня;        }
var i,j,ii,jj,iii,jjj,pop : integer;
begin
       for i:=1 to M+1 do
       begin
         Masy3[i]:= Masy2[i];
         {writeln(F2,'Masy3[',i,'] = ',Masy3[i]); }
       end;
       pop:= Masy2[M+1]-1;
       {writeln(F2,'pop = ',pop); }
       for i:= 1 to pop do
       begin
         Mass3[i]:= Mass2[i];
         {writeln(F2,'Mass3[',i,'] = ',Mass3[i]);   }
       end;
end;{Perezapis1}
{************************************************************}
 procedure Perezapis2(var My,M : integer;
                var Masy3 : TMasy;
                var Mass3 : TMass;
                var Masy4 : TMasy;
                var Mass4 : TMass);
{Процедура перезаписи массива                                    }
{                                                                }
{ MY - номер уровня;                                             }
{                                                                }
{ Masy3 - массив указателей для множества реберных разрезов n-уровня;}
{ Mass3 - массив элементов множеств беберных разрезов n-уровня;       }
{ Masy4 - массив указателей для множества реберных разрезов т-уровня; }
{ Mass4 - массив элементов множеств беберных разрезов n-уровня;        }
var i,j,ii,jj,iii,jjj,pop,pop1 : integer;
begin
```

```
          pop:=Masy3[(MY-1)*M+1]-1;
          {writeln(F2,'pop = ',pop); }
          for i:=1 to M+1 do
          begin
            Masy4[i]:= Masy3[(MY-1)*M+i]-pop;
            {writeln(F2,'Masy4[',i,'] = ',Masy4[i]); }
          end;
          pop1:= Masy3[(MY-1)*M+M+1]-Masy3[(My-1)*M+1];
          {writeln(F2,'pop1 = ',pop1); }
          for i:= 1 to pop1 do
          begin
            Mass4[i]:= Mass3[pop+i];
            {writeln(F2,'Mass4[',i,'] = ',Mass4[i]); }
          end;
end;{Perezapis2}
{*************************************************************}
 procedure FormRasRebN(var MY,M : integer;
                        var Masy2 : TMasy;
                        var Mass2 : TMass;
                        var Masy4 : TMasy;
                        var Mass4 : TMass;
                        var Masy6 : TMasy;
                        var Mass6 : TMass;
                        var MasMdop : TMasy);
{Процедура создания множества базовых реберных разрезов графа       }
{                                                                   }
{ MY - номер уровня;                                                }
{ Masy1  - массив указателей для базовых реберных разрезов;         }
{ Mass2  - массив элементов базовых реберных разрезов               }
{ Masy2 - массив указателей для множества реберных разрезов n-уровня;}
{ Mass2 - массив элементов множеств беберных разрезов n-уровня;     }
{ MasMdop -  вспомогательный массив для хранения суграфа.           }
var i,j,ii,jj,iii,jjj,DD,K23,K24,K21,pop,pop1 : integer;
begin
          Masy5[1]:=1;
          for i:=1 to M do
          begin  {1}
            MasMdop[1]:=0;
            K23:=1;
            for j:= Masy4[i] to Masy4[i+1]-1 do
            begin {2}
              K24:=Mass4[j]; {Номер ребра}
              {writeln(F2,'K24 = ',K24); }
              if K24<>0 then
              Begin
              for jj:= Masy2[K24] to Masy2[K24+1]-1 do
              begin {3}
                K23:=K23+1;
                MasMdop[K23]:=Mass2[jj];
                {writeln(F2,'MasMdop[',K23,'] = ',MasMdop[K23]);}
              end;  {3}
              end;
              {for ii:=1 to K23 do
              begin
                if ii<>K23 then write(F2,MasMdop[ii],' ');
                if ii=K23 then writeln(F2,MasMdop[ii]);
              end; }
              for ii:=1 to K23-1 do
              if MasMdop[ii]<>0 then
               begin {4}
                for jj:=ii+1 to K23 do
                if MasMdop[jj]=MasMdop[ii] then
                begin {5}
```

```
              MasMdop[jj]:=0;
              MasMdop[ii]:=0;
            end; {5}
          end; {4}
          {for ii:=1 to K23 do
          begin
            if ii<>K23 then write(F2,MasMdop[ii],' ');
            if ii=K23 then writeln(F2,MasMdop[ii]);
          end; }
          DD:=0;
          for iii:=1 to K23 do
          if MasMdop[iii]<>0 then
          begin  {6}
            DD:=DD+1;
            {writeln(F2,'MasMdop[',iii,'] = ',MasMdop[iii]); }
            MasMdop[DD]:=MasMdop[iii];
            {writeln(F2,'MasMdop[',DD,'] = ',MasMdop[DD]); }
          end;  {6}
          {for iii:=1 to DD do
          begin
            if iii<>DD then write(F2,MasMdop[iii],' ');
            if iii=DD then writeln(F2,MasMdop[iii]);
          end;}
        K23:=DD;
        end;  {2}
        {for iii:=1 to DD do
          begin
            if iii<>DD then write(F2,MasMdop[iii],' ');
            if iii=DD then writeln(F2,MasMdop[iii]);
          end; }
        if DD=0 then
        begin
          Masy5[i+1]:=Masy5[i]+1;
          pop:=Masy5[i]-1;
          Mass5[pop+1]:=0;
        end;
        if DD<>0 then
        begin
          Shell(DD,MasMdop);
          Masy5[i+1]:=Masy5[i]+DD;
          pop:=Masy5[i]-1;
          {writeln(F2,pop); }
          K21:=0;
          for ii:=1 to DD do
          begin
            K21:=K21+1;
            Mass5[pop+ii]:=MasMdop[K21];
            {writeln(F2,'Mass3[',pop+ii,'] = ',Mass3[pop+ii]); }
          end;
        end;
      end;   {1}
end;{FormRasRebN}
{*************************************************************}
procedure Ravenstvo(var IR : integer;
                    var Masy4 : TMasy;
                    var Mass4 : TMass;
                    var Masy5 : TMasy;
                    var Mass5 : TMass;
                    var MasKol : TMasy);
{Процедура проверки суграфов на совпадение                   }
{                                                            }
{ IR - текущее ребро;                                        }
{ Masy4  - массив указателей для предыдущих реберных разрезов;   }
```

```
{ Mass4  - массив элементов предыдущих реберных разрезов               }
{ Masy5 - массив указателей для текущего реберных разрезов;            }
{ Mass5 - массив элементов множеств реберных разрезов n-уровня;      }
{ MasKol -  вспомогательный массив для признаков.                        }
var i,j,ii,jj,iii,jjj,pop,pop1,priz : integer;
var a1,a2,b1,b2,c1,c2 : integer;
label 1;
        begin
        {writeln(F2,'IR = ',IR); }
        priz:=0;
        {writeln(F2,'  MasKol[',IR,'] = ',MasKol[IR]); }
        if MasKol[IR]=0 then goto 1;
        if Masy5[IR]=0 then goto 1;
        c1:=Masy5[IR+1]-Masy5[IR];
        c2:=Masy4[IR+1]-Masy4[IR];
        {writeln(F2,'  c1 = ',c1,'  c2 = ',c2); }
        if c1=c2 then
        begin  {1}
           pop:=Masy5[IR]-1;
           pop1:=Masy4[IR]-1;
           for ii:=1 to c1 do
           begin
           {writeln(F2,' Mass5[',pop+ii,'] = ',Mass5[pop+ii],
           '  Mass4[',pop1+ii,'] = ',Mass4[pop1+ii]); }
           if Mass5[pop+ii]<>Mass4[pop1+ii] then
             begin  {2}
               priz:=1;
               goto 1;
             end;  {2}
           end;
        end;   {1}
        if c1<>c2 then priz:=1;
1:       MasKol[IR]:=priz;
        {writeln(F2,'MasKol[',IR,'] = ',MasKol[IR]); }
end;{Ravenstvo}
{*****************************************************************}
procedure Perezapis3(var M : integer;
                     var Masy5 : TMasy;
                     var Mass5 : TMass;
                     var Masy6 : TMasy;
                     var Mass6 : TMass;
                     var MasKol : TMasy);
{Процедура перезаписи массива                                            }
{  Mass5  в  Mass6                                                       }
{ M - количество ребер в графе;                                          }
{                                                                        }
{ Masy5 - массив указателей для множества реберных разрезов n-уровня;}
{ Mass5 - массив элементов множеств беберных разрезов n-уровня;         }
{ Masy6 - массив указателей для множества реберных разрезов т-уровня; }
{ Mass6 - массив элементов множеств беберных разрезов n-уровня;          }
var i,j,ii,jj,iii,jjj,pop,pop1,a1,a2,c1,b1 : integer;
begin
        Masy6[1]:=1;
        pop:=0;
        for i:=1 to M do
        begin  {1}
           if MasKol[i]=0 then
           begin   {2}
             Masy6[i+1]:=Masy6[i]+1;
             pop:=pop+1;
             Mass6[pop]:=0;
           end;   {2}
           if MasKol[i]=1 then
```

```
           begin   {3}
             a1:=Masy5[i];
             a2:=Masy5[i+1]-1;
             c1:=Masy5[i+1]-Masy5[i];
             Masy6[i+1]:=Masy6[i]+c1;
             b1:=a1-1;
             for ii:= 1 to c1 do
             begin  {4}
               pop:=pop+1;
               Mass6[pop]:=Mass5[b1+ii];
             end;   {4}
           end;   {3}
         end;   {1}
end;{Perezapis3}
{*************************************************************}
procedure FormPrikaz(var Prikaz,M : integer;
                     var MasKol : TMasy);
{Процедура формирования признака на последующее выполнение     }
{   если Prikaz=1 то идти, если Prikaz=0 конец работы программы       }
{                                                                              }
{ Prikaz - признак останова расчета;                                 }
{ MasKol -  вспомогательный массив для признаков.                     }
var i,j,ii,jj,iii,jjj: integer;
     begin
       Prikaz:=0;
       for i:=1 to M do if MasKol[i]=1 then Prikaz:=1;
end;{ FormPrikaz}
{*************************************************************}
procedure Perezapis4(var M,My : integer;
                     var Masy3 : TMasy;
                     var Mass3 : TMass;
                     var Masy6 : TMasy;
                     var Mass6 : TMass);
{Процедура перезаписи массива                                               }
{                                                                              }
{ M - количество ребер в графе;                                             }
{ MY - номер уровня;                                                        }
{ Masy3 - массив указателей для множества реберных разрезов n-уровня;    }
{ Mass3 - массив элементов множеств беберных разрезов n-уровня;            }
{ Masy6 - массив указателей для множества реберных разрезов т-уровня;     }
{ Mass6 - массив элементов множеств беберных разрезов n-уровня;             }
var i,j,ii,jj,iii,jjj,pop,pop1 : integer;
begin
     {writeln(F2,' Начало работы Perezapis4 '); }
     for i:=1 to M do
     begin  {1}
       Masy3[(My-1)*M+1+i]:= Masy3[(MY-1)*M+i]+Masy6[i+1]-Masy6[i];
       {writeln(F2,'Masy3[',(My-1)*M+1+i,'] = ',Masy3[(My-1)*M+1+i]); }
     end;
     pop:= Masy6[M+1]-1;
     pop1:=Masy3[(My-1)*M+1]-1;
     {writeln(F2,'pop = ',pop,'pop1 = ',pop1);}
     for ii:= 1 to pop do
     begin
       Mass3[pop1+ii]:= Mass6[ii];
       {writeln(F2,'Mass3[',pop1+ii,'] = ',Mass3[pop1+ii]);}
     end;  {1}
     {writeln(F2,' Конец работы Perezapis4 '); }
end;{Perezapis4}
{*************************************************************}
procedure Perezapis64(var M : integer;
                     var Masy4 : TMasy;
                     var Mass4 : TMass;
```

```
                   var Masy6 : TMasy;
                   var Mass6 : TMass);
{Процедура перезаписи массива                                        }
{                                                                    }
{ M - количество ребер в графе;                                      }
{ MY - номер уровня;                                                 }
{ Masy4 - массив указателей для множества реберных разрезов n-уровня;  }
{ Mass4 - массив элементов множеств беберных разрезов n-уровня;        }
{ Masy6 - массив указателей для множества реберных разрезов т-уровня;  }
{ Mass6 - массив элементов множеств беберных разрезов n-уровня;        }
var i,j,ii,jj,iii,jjj,pop,pop1 : integer;
begin
       for i:=1 to M+1 do
            Masy4[i]:= Masy6[i];
       pop:= Masy6[M+1]-1;
       for j:= 1 to pop do
          Mass4[j]:= Mass6[j];
end;{Perezapis64}
{***********************************************************}
 procedure Printura1(var M : integer;
                   var May : TMasy;
                   var Mas : TMass);
{Процедура перезаписи печати массива                    }
var i,j,ii,jj,iii,jjj : integer;
begin
       for i:=1 to M+1 do
       begin
         if i<> M+1 then write(F2, May[i],' ');
         if i= M+1 then writeln(F2, May[i]);
       end;
       for i:=1 to M do
       begin
        for j:=May[i] to May[i+1]-1 do
        begin
          if j<>May[i+1]-1 then write(F2,Mas[j],' ');
          if j=May[i+1]-1 then writeln(F2,Mas[j]);
        end;
       end;
end;{Printura1}
{***********************************************************}
procedure Printura2(var M,My : integer;
                   var May : TMasy;
                   var Mas : TMass);
{Процедура перезаписи печати массива                    }
var i,j,ii,jj,iii,jjj : integer;
begin
       writeln(F2,' Начало работы Printura2 ');
       for i:=1 to MY*M+1 do
       begin
         if i<> MY*M+1 then write(F2, May[i],' ');
         if i= MY*M+1 then writeln(F2, May[i]);
       end;
       for i:=1 to MY*M do
       begin
        for j:=May[i] to May[i+1]-1 do
        begin
          if j<>May[i+1]-1 then write(F2,Mas[j],' ');
          if j=May[i+1]-1 then writeln(F2,Mas[j]);
        end;
       end;
       writeln(F2,' Окончание работы Printura2 ');
end;{Printura2}
{***********************************************************}
```

```
procedure Printura3(var M : integer;
                    var Mas : TMasy);
{Процедура перезаписи печати массива                    }
var i,j,ii,jj,iii,jjj : integer;
begin
       for i:=1 to M do
       begin
         if i<> M then write(F2, Mas[i],' ');
         if i= M then writeln(F2, Mas[i]);
       end;
end;{Printura3}
{***********************************************************}
 procedure FormMasKol(var M : integer;
                    var Masy5 : TMasy;
                    var Mass5 : TMass;
                    var MasKol : TMasy);
{Процедура формирования  массива  MasKol                    }
var i,j,ii,jj,iii,jjj : integer;
begin
       for i:=1 to M do
       begin
         if Mass5[Masy5[i]]=0 then
         MasKol[i]:=0;
         if Mass5[Masy5[i]]<>0 then
         MasKol[i]:=1;
       end;
end;{FormMasKol}
{***********************************************************}
 procedure PriNasReb(var M,My,Nv : integer;
                    var Masy3 : TMasy;
                    var Mass3 : TMass;
                    var MasWesReb : TMasy;
                    var Masy1 : TMasy;
                    var Massi : TMass;
                    var MasWesVer : TMasy);
{Процедура перезаписи печати массива                    }
var i,j,ii,jj,iii,jjj,pop,pop1,a1 : integer;
begin
       for j:=1 to MY do
       begin
       writeln(F2,'  Номер уровня (ярус):  ',j);
        for i:=1 to M do
        begin
          pop:=(j-1)*M+i;
          pop1:=(j-1)*M+i+1;
          for ii:= Masy3[pop] to Masy3[pop1]-1 do
          begin
           if ii<>Masy3[pop1]-1 then write(F2,Mass3[ii],' ');
           if ii=Masy3[pop1]-1 then writeln(F2,Mass3[ii]);
          end;
        end;
       end;
       for i:=1 to M do MasWesReb[i]:=0;
       pop:=Masy3[M*MY+1]-1;
       for i:=1 to pop do
       begin
        pop1:=Mass3[i];
        if pop1<>0 then MasWesReb[pop1]:= MasWesReb[pop1]+1;
       end;
       writeln(F2,'  Кортеж весов ребер :  ');
       for j:=1 to M do
       begin
           if j<>M then write(F2,MasWesReb[j],' ');
```

```
      if j=M then writeln(F2,MasWesReb[j]);
    end;
    for j:=1 to Nv do
    begin
     MasWesVer[j]:=0;
     {writeln(F2,'  Номер уровня (ярус):  ',j); }
     pop:=Masy1[j+1]-Masy1[j];
     pop1:=Masy1[j]-1;
     for i:=1 to pop do
     begin
      a1:=Massi[pop1+i];
      MasWesVer[j]:= MasWesVer[j]+MasWesReb[a1];
     end;
    end;
    writeln(F2,'  Кортеж весов вершин :  ');
    for j:=1 to Nv do
    begin
      if j<>Nv then write(F2,MasWesVer[j],' ');
      if j=Nv then writeln(F2,MasWesVer[j]);
    end;
end;{PriNasReb}
{**************************************************************}
{**************************************************************}
  label 1,2;
    begin
    assign(F1,'D:\Isomorf\GRF\6a10v04.grf');
    reset(F1);
    readln(F1,Nv);
    for I:=1 to Nv+1 do
    begin
     if I<>Nv+1 then read(F1,Masy1[I]);
     if I=Nv+1 then readln(F1,Masy1[I]);
    end;
    Np:=Masy1[Nv+1]-1;
    for I:=1 to Np do
    begin
     if I<>Np then read(F1,Mass1[I]);
     if I=Np then read(F1,Mass1[I]);
    end;
    close (F1);
    { Создаём новый файл и открываем его в режиме "для чтения и записи"}
    Assign(F2,'D:\Isomorf\DUB\6a10v04.dub');
    Rewrite(F2);
    FormIncide(Nv,Masy1,Mass1,Massi);
    writeln(F2,'   Матрица смежностей графа');
    for I:=1 to Nv+1 do
    begin
     if i<>Nv+1 then write(F2,Masy1[i],' ');
     if i=Nv+1 then writeln(F2,Masy1[i]);
    end;
    for I:=1 to Nv do
    for j:=Masy1[i] to Masy1[i+1]-1 do
    begin
     if j<>Masy1[i+1]-1 then write(F2,Mass1[j],' ');
     if j=Masy1[i+1]-1 then writeln(F2,Mass1[j]);
    end;
    writeln(F2,'   Матрица инциденций графа');
    for I:=1 to Nv do
    for j:=Masy1[i] to Masy1[i+1]-1 do
    begin
     if j<>Masy1[i+1]-1 then write(F2,Massi[j],' ');
     if j=Masy1[i+1]-1 then writeln(F2,Massi[j]);
    end;
```

```
FormRasReb2(Nv,K11,Masy1,Mass1,Masy2,Mass2,Massi,MasMdop);
writeln(F2,'  Количество вершин в графе = ',Nv);
M:=Np div 2;
writeln(F2,'  Количество ребер в графе = ',M);
writeln(F2,'  Вычислен массив базовых реберных разрезов графа.');
Printura1(M,Masy2,Mass2);
MY:=1;
Perezapis1(M,Masy2,Mass2,Masy3,Mass3);
{writeln(F2,'  Помещен массив базовых реберных разрезов Mass3 в базу');}
{Printura1(M,Masy3,Mass3); }
Perezapis2(My,M,Masy3,Mass3,Masy4,Mass4);
{writeln(F2,'  Сформирован массив Mass4'); }
{Printura1(M,Masy4,Mass4); }
prikaz:=1;
while prikaz=1 do
begin
  {writeln(F2,'  Входим на следующий уровень'); }
  {writeln(F2,'   Номер уровня My = ',My+1);}
  {Printura1(M,Masy2,Mass2); }
  {Printura1(M,Masy4,Mass4); }
  FormRasRebN(MY,M,Masy2,Mass2,Masy4,Mass4,Masy5,Mass5,MasMdop);
  {writeln(F2,'   Вычислен массив следующего уровня Mass5'); }
  {Printura1(M,Masy5,Mass5); }
  FormMasKol(M,Masy5,Mass5,MasKol);
  for iii:=1 to M do
  begin
    KKK:=iii;
    Ravenstvo(KKK,Masy4,Mass4,Masy5,Mass5,MasKol);
  end;
  {Printura3(M,MasKol);  }
  FormPrikaz(Prikaz,M,MasKol);
  if Prikaz=0 then goto 1;
  if My=1 then
  begin
    K12:=MY+1;
    Perezapis3(M,Masy5,Mass5,Masy6,Mass6,MasKol);
    {writeln(F2,'  Вычислен массив реберных разрезов Mass5'); }
    {Printura1(M,Masy5,Mass5); }
    {writeln(F2,'  Помещен массив реберных разрезов Mass6 '); }
    {Printura1(M,Masy6,Mass6);  }
    Perezapis4(M,K12,Masy3,Mass3,Masy6,Mass6);
    {writeln(F2,'  Помещен массив базовых реберных разрезов Mass3 в базу'); }
    {Printura2(M,K12,Masy3,Mass3); }
    Perezapis64(M,Masy4,Mass4,Masy6,Mass6);
    {writeln(F2,'  Вычислен массив реберных разрезов Mass4'); }
    {Printura1(M,Masy4,Mass4); }
    for j:=1 to M do MasWesReb[j]:=0;
    for j:=Masy3[(My-1)*M+1] to Masy3[(My-1)*M+M+1]-1 do
    begin
      p11:=Mass3[j];
      {writeln(F2,'  j =  ',j,'   p11 =  ',p11);}
      if p11<>0 then MasWesReb[p11]:= MasWesReb[p11]+1;
    end;
    writeln(F2,'  Кортеж весов ребер уровня  ',My,':  ');
    for j:=1 to M do
    begin
        if j<>M then write(F2,MasWesReb[j],' ');
        if j=M then writeln(F2,MasWesReb[j]);
    end;
  end;
  if My>1 then
  begin
    for jjj:=1 to My do
```

```
      begin
       {writeln(F2,' jjj = ',jjj); }
       K11:=MY+1-jjj;
       {writeln(F2,' K11 = ',K11); }
       Perezapis2(K11,M,Masy3,Mass3,Masy4,Mass4);
       {writeln(F2,'   Вычислен массив следующего уровня Mass4'); }
       {Printura1(M,Masy4,Mass4);}
       FormMasKol(M,Masy5,Mass5,MasKol);
       for iii:=1 to M do
       begin
       KKK0:=iii;
       Ravenstvo(KKK0,Masy4,Mass4,Masy5,Mass5,MasKol);
       end;
       {Printura3(M,MasKol); }
       FormPrikaz(Prikaz,M,MasKol);
       if Prikaz=0 then goto 1;
       K12:=MY+1;
       Perezapis3(M,Masy5,Mass5,Masy6,Mass6,MasKol);
       {writeln(F2,'  Вычислен массив реберных разрезов Mass5'); }
       {Printura1(M,Masy5,Mass5); }
       {writeln(F2,'  Помещен массив реберных разрезов Mass6 '); }
       {Printura1(M,Masy6,Mass6);}
       {Printura3(M,MasKol); }
       FormPrikaz(Prikaz,M,MasKol);
       Perezapis4(M,K12,Masy3,Mass3,Masy6,Mass6);
       {writeln(F2,'  Помещен массив базовых реберных разрезов Mass3 в базу'); }
       {Printura2(M,K12,Masy3,Mass3); }
       Perezapis64(M,Masy4,Mass4,Masy6,Mass6);
       {writeln(F2,'  Вычислен массив реберных разрезов Mass4'); }
       {Printura1(M,Masy4,Mass4); }
       if Prikaz=0 then goto 1;
      end;
    end;
    My:=MY+1;
    {writeln(F2,'  Переходим на следующий уровень'); }
    for j:=1 to M do MasWesReb[j]:=0;
    for j:=Masy3[(My-1)*M+1] to Masy3[(My-1)*M+M+1]-1 do
    begin
     p11:=Mass3[j];
     {writeln(F2,'  j =  ',j,'   p11 =  ',p11);}
     if p11<>0 then MasWesReb[p11]:= MasWesReb[p11]+1;
    end;
    writeln(F2,'  Кортеж весов ребер уровня  ',My,':  ');
    for j:=1 to M do
    begin
       if j<>M then write(F2,MasWesReb[j],' ');
       if j=M then writeln(F2,MasWesReb[j]);
    end;
  end;
1:  writeln(F2,'  Последний уровень = ',MY);
  {Printura2(M,K12,Masy3,Mass3); }
  PriNasReb(M,My,Nv,Masy3,Mass3,MasWesReb,Masy1,Massi,MasWesVer);
  Shell(M,MasWesReb);
  writeln(F2,'  Вектор весов ребер :  ');
  for j:=1 to M do
  begin
   if j<>M then write(F2,MasWesReb[j],' ');
   if j=M then writeln(F2,MasWesReb[j]);
  end;
  Shell(Nv,MasWesVer);
  writeln(F2,'  Вектор весов вершин :  ');
  for j:=1 to Nv do
  begin
```

```
    if j<>Nv then write(F2,MasWesVer[j],' ');
    if j=Nv then writeln(F2,MasWesVer[j]);
  end;
  {  определение веса ребер        }
  close (F2);
  writeln('  Конец расчета  ');

end.
```

## 2.8. Файлы программы Raschet54

### Входной файл 6a10v04.grf

```
6                              {количество вершин графа}
1 4 7 10 14 17 21              {массив указателей}
2  4  6                        {смежные вершины с вершиной 1}
1  4  5                        {смежные вершины с вершиной 2}
4  5  6                        {смежные вершины с вершиной 3}
1  2  3  6                     {смежные вершины с вершиной 4}
2  3  6                        {смежные вершины с вершиной 5}
1  3  4  5                     {смежные вершины с вершиной 6}
```

### Выходной файл 6a10v04.dub

```
  Матрица смежностей графа
1 4 7 10 14 17 21              {массив указателей}
2 4 6                          {смежные вершины с вершиной 1}
1 4 5                          {смежные вершины с вершиной 2
4 5 6                          {смежные вершины с вершиной 3}
1 2 3 6                        {смежные вершины с вершиной 4}
2 3 6                          {смежные вершины с вершиной 5}
1 3 4 5                        {смежные вершины с вершиной 6}
  Матрица инциденций графа
1 2 3                          {инцидентные ребра для вершины 1}
1 4 5                          {инцидентные ребра для вершины 2}
6 7 8                          {инцидентные ребра для вершины 3}
2 4 6 9                        {инцидентные ребра для вершины 4}
5 7 10                         {инцидентные ребра для вершины 5}
3 8 9 10                       {инцидентные ребра для вершины 6}
 Количество вершин в графе = 6
 Количество ребер в графе = 10
 Вычислен массив базовых реберных разрезов графа.
1 5 10 15 20 24 29 33 38 44 49 {массив указателей}
2 3 4 5                        {базовый реберный разрез ребра 1}
1 3 4 6 9                      {базовый реберный разрез ребра 2}
1 2 8 9 10                     {базовый реберный разрез ребра 3}
1 2 5 6 9                      {базовый реберный разрез ребра 4}
1 4 7 10                       {базовый реберный разрез ребра 5}
2 4 7 8 9                      {базовый реберный разрез ребра 6}
5 6 8 10                       {базовый реберный разрез ребра 7}
3 6 7 9 10                     {базовый реберный разрез ребра 8}
2 3 4 6 8 10                   {базовый реберный разрез ребра 9}
3 5 7 8 9                      {базовый реберный разрез ребра 10}
```

Кортеж весов ребер уровня 1:
4 5 5 5 4 5 4 5 6 5
Кортеж весов ребер уровня 2:
4 5 5 3 4 5 4 5 6 3
Кортеж весов ребер уровня 3:
3 4 4 3 4 4 3 4 6 3
Кортеж весов ребер уровня 4:
3 4 4 3 4 4 3 4 6 3
Последний уровень = 4
Номер уровня (ярус): 1
2 3 4 5 {реберный разрез ребра 1}
1 3 4 6 9 {реберный разрез ребра 2}
1 2 8 9 10 {реберный разрез ребра 3}
1 2 5 6 9 {реберный разрез ребра 4}
1 4 7 10 {реберный разрез ребра 5}
2 4 7 8 9 {реберный разрез ребра 6}
5 6 8 10 {реберный разрез ребра 7}
3 6 7 9 10 {реберный разрез ребра 8}
2 3 4 6 8 10 {реберный разрез ребра 9}
3 5 7 8 9 {реберный разрез ребра 10}
Номер уровня (ярус): 2
3 5 7 8 9 {реберный разрез ребра 1}
2 4 7 8 9 {реберный разрез ребра 2}
1 3 4 6 9 {реберный разрез ребра 3}
2 3 4 5 {реберный разрез ребра 4}
0 {реберный разрез ребра 5}
3 6 7 9 10 {реберный разрез ребра 6}
1 2 5 6 9 {реберный разрез ребра 7}
1 2 8 9 10 {реберный разрез ребра 8}
1 2 3 6 7 8 {реберный разрез ребра 9}
5 6 8 10 {реберный разрез ребра 10}
Номер уровня (ярус): 3
5 6 8 10 {реберный разрез ребра 1}
3 6 7 9 10 {реберный разрез ребра 2}
2 4 7 8 9 {реберный разрез ребра 3}
3 5 7 8 9 {реберный разрез ребра 4}
0 {реберный разрез ребра 5}
1 2 8 9 10 {реберный разрез ребра 6}
2 3 4 5 {реберный разрез ребра 7}
1 3 4 6 9 {реберный разрез ребра 8}
0 {реберный разрез ребра 9}
1 2 5 6 9 {реберный разрез ребра 10}
Номер уровня (ярус): 4
1 2 5 6 9 {реберный разрез ребра 1}
1 2 8 9 10 {реберный разрез ребра 2}
3 6 7 9 10 {реберный разрез ребра 3}
5 6 8 10 {реберный разрез ребра 4}
0 {реберный разрез ребра 5}
1 3 4 6 9 {реберный разрез ребра 6}
3 5 7 8 9 {реберный разрез ребра 7}
2 4 7 8 9 {реберный разрез ребра 8}
0 {реберный разрез ребра 9}

2 3 4 5 {реберный разрез ребра 10}
Кортеж весов ребер :
14 18 18 14 16 18 14 18 24 14
Кортеж весов вершин :
50 44 50 74 44 74
Вектор весов ребер :
14 14 14 14 16 18 18 18 18 24
Вектор весов вершин :
44 44 50 50 74 74

## Комментарии

В данной главе рассмотрены основные свойства множеств реберных разрезов графа G. Показано, что линейный нильпотентный оператор A(L(G)) порождает конечную цепочку реберных разрезов, представленных в виде суграфов подпространства S(G) графа, принадлежащих спектру реберных разрезов. Это, в свою очередь, позволяет оценить численный вклад каждого ребра в спектр реберных разрезов, называемый весом ребра. Таким образом, появляется возможность создать ряд числовых характеристик, описывающих структуру графа G в виде векторных инвариантов спектра реберных разрезов графа G. Числовые характеристики позволяют проводить сравнительный анализ и синтез структур графов.

Вычислительная сложность алгоритма построения векторных инвариантов спектра реберных разрезов, равна $O(m^3)$. Следовательно, задача определения векторных инвариантов спектра реберных разрезов графа, относится к классу P – полиномиальных алгоритмов.

# Глава 3. ИЗОМЕТРИЧЕСКИЕ ЦИКЛЫ ГРАФА

## 3.1. Фундаментальные циклы и разрезы

Как известно, в пространстве суграфов можно выделить два подпространства называемых подпространством разрезов S и подпространством циклов C. Обычно в теории графов для поиска базисов применяется фундаментальная система циклов и разрезов. Данная система образуется в результате выделения случайного дерева графа (ациклического суграфа), тем самым разделяя ребра графа на ветви дерева и хорды. Ребра, принадлежащие дереву, называются ветвями, а не принадлежащие дереву – хордами. Каждый фундаментальный цикл образуется как объединение одной хорды и ветвей дерева. Рассматривая все хорды для выделенного дерева графа, строим матрицу фундаментальных циклов. Если представить (0,1)-матрицу фундаментальных циклов в виде единичной матрицы хорд и блочной матрицы π состоящей из ветвей дерева, то можно получить матрицу фундаментальных разрезов графа в виде единичной матрицы ветвей дерева и блочной матрицы $\rho = \pi^t$ состоящей из хорд.

$C_ф$=

| | $e_2$ | $e_4$ | $e_5$ | $e_6$ | $e_8$ | $e_9$ | $e_1$ | $e_3$ | $e_7$ | $e_{10}$ | $e_{11}$ |
|---|---|---|---|---|---|---|---|---|---|---|---|
| $e_2$ | 1 | | | | | | | 1 | | 1 | 1 |
| $e_4$ | | 1 | | | | | 1 | 1 | 1 | 1 | 1 |
| $e_5$ | | | 1 | | | | 1 | 1 | | 1 | 1 |
| $e_6$ | | | | 1 | | | 1 | 1 | | | |
| $e_8$ | | | | | 1 | | | | 1 | 1 | |
| $e_9$ | | | | | | 1 | | | 1 | 1 | 1 |

Единичная блочная подматрица | Блочная подматрица π

Матрица фундаментальных разрезов имеет вид

$S_ф$=

| | $e_1$ | $e_3$ | $e_7$ | $e_{10}$ | $e_{11}$ | $e_2$ | $e_4$ | $e_5$ | $e_6$ | $e_8$ | $e_9$ |
|---|---|---|---|---|---|---|---|---|---|---|---|
| $e_1$ | 1 | | | | | | 1 | 1 | 1 | | |
| $e_3$ | | 1 | | | | 1 | 1 | 1 | 1 | | |
| $e_7$ | | | 1 | | | | 1 | | | 1 | 1 |
| $e_{10}$ | | | | 1 | | 1 | 1 | 1 | | 1 | 1 |
| $e_{11}$ | | | | | 1 | 1 | 1 | | | | 1 |

Единичная блочная подматрица | Транспонированная блочная подматрица $\rho = \pi^t$

Например, для графа $G_2$ (рис. 3.1) относительно дерева T = {$e_1$,$e_3$,$e_7$,$e_{10}$,$e_{11}$} матрица фундаментальных циклов представлена выше.

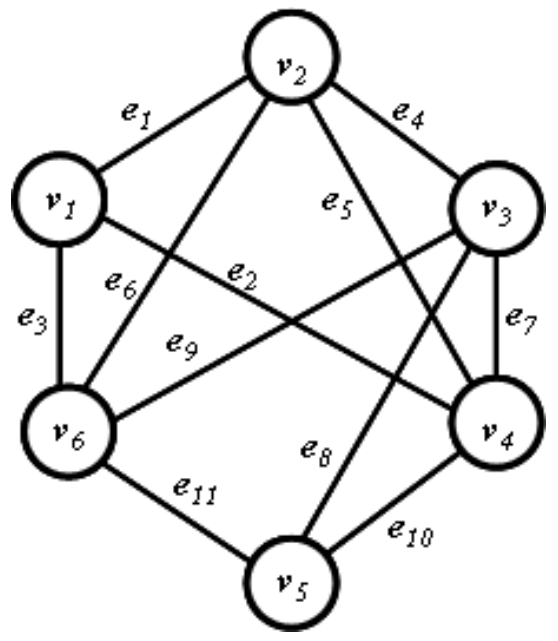

Рис. 3.1. Граф $G_2$.

Таким образом, относительно ациклического суграфа (дерева) строится система фундаментальных циклов и фундаментальных разрезов графа. Фундаментальная система циклов и фундаментальная система разрезов, в свою очередь, служат для формирования базиса подпространства циклов и базиса подпространства разрезов.

Для определения базиса подпространства циклов и базиса подпространства разрезов в теории графов применяется фундаментальная система циклов и разрезов. Количество фундаментальных циклов определяется цикломатическим числом графа $\nu(\mathbf{G}) = m-n+1$, а количество фундаментальных разрезов определяется рангом графа $\rho(\mathbf{G}) = n - 1$ [7,9,10,34].

Любой суграф, принадлежащий подпространству разрезов S графа G, в общем случае, является квалиразрезом. Количество ребер составляющих квалиразрез называется *длиной квалиразреза.* Например, фундаментальные разрезы графа G относительно дерева T = $\{e_1,e_3,e_7,e_{10},e_{11}\}$ (рис. 3.2):

$s_1 = \{\mathbf{e}_1,e_4,e_5,e_6\}$; $s_2 = \{e_2,\mathbf{e}_3,e_4,e_5,e_6\}$; $s_3 = \{e_4,\mathbf{e}_7,e_8,e_9\}$;

$s_4 = \{e_2,e_4,e_5,e_8,e_9,\mathbf{e}_{10}\}$; $s_5 = \{e_2,e_4,e_5,e_9,\mathbf{e}_{11}\}$

имеют суммарную длину $l = l_1 + l_2 + l_3 + l_4 + l_5 = 4 + 5 + 4 + 6 + 5 = 24$.

А следующие базисные разрезы графа G полученные путем линейной комбинации фундаментальных разрезов

$s_6 = s_1 \oplus s_2 = \{e_1,e_4,e_5,e_6\} \oplus \{e_2,e_3,e_4,e_5,e_6\} = \{\mathbf{e}_1,e_2,\mathbf{e}_3\}$;
$s_7 = s_1 = \{\mathbf{e}_1,e_4,e_5,e_6\}$;
$s_8 = s_3 = \{e_4,\mathbf{e}_7,e_8,e_9\}$;
$s_9 = s_3 \oplus s_4 = \{e_4,e_7,e_8,e_9\} \oplus \{e_2,e_4,e_5,e_8,e_9,e_{10}\} = \{e_2,e_5,\mathbf{e}_7,\mathbf{e}_{10}\}$;
$s_{10} = s_4 \oplus s_5 = \{e_2,e_4,e_5,e_8,e_9,e_{10}\} \oplus \{e_2,e_4,e_5,e_9,e_{11}\} = \{e_8,\mathbf{e}_{10},\mathbf{e}_{11}\}$,

имеют меньшую суммарную длину $l = l_6 + l_7 + l_8 + l_9 + l_{10} = 3+4+4+4+3 = 18$. Такие разрезы, имеющие минимально возможную длину, будем называть центральными разрезами графа G [10].

С другой стороны, длина центрального разреза ставит в соответствие локальную степень вершины и наоборот.

Рассмотрим фундаментальные циклы для выбранного дерева графа $G_2$ (рис. 3.1):

$c_1 = \{e_2, \mathbf{e}_3, \mathbf{e}_{10}, \mathbf{e}_{11}\}$; $c_2 = \{\mathbf{e}_1, \mathbf{e}_3, e_4, \mathbf{e}_7, e_{10}, \mathbf{e}_{11}\}$; $c_3 = \{\mathbf{e}_1, \mathbf{e}_3, e_5, \mathbf{e}_{10}, \mathbf{e}_{11}\}$;
$c_4 = \{\mathbf{e}_1, \mathbf{e}_3, e_6\}$; $c_5 = \{\mathbf{e}_7, e_8, \mathbf{e}_{10}\}$; $c_5 = \{\mathbf{e}_7, e_9, \mathbf{e}_{10}, \mathbf{e}_{11}\}$

имеют суммарную длину $l = l_1 + l_2 + l_3 + l_4 + l_5 + l_6 = 4 + 6 + 5 + 3 + 3 + 4 = 25$.

А следующие базисные циклы графа $G_2$, полученные путем линейной комбинации фундаментальных циклов

$c_7 = c_1 = \{e_2, \mathbf{e}_3, \mathbf{e}_{10}, \mathbf{e}_{11}\}$;
$c_8 = c_4 = \{\mathbf{e}_1, \mathbf{e}_3, e_6\}$;
$c_9 = c_5 = \{\mathbf{e}_7, e_8, \mathbf{e}_{10}\}$;
$c_{10} = c_5 \oplus c_6 = \{e_7, e_8, e_{10}\} \oplus \{e_7, e_9, e_{10}, e_{11}\} = \{e_8, e_9, \mathbf{e}_{11}\}$;
$c_{11} = c_2 \oplus c_4 \oplus c_6 = \{e_1, e_3, e_4, e_7, e_{10}, e_{11}\} \oplus \{e_7, e_9, e_{10}, e_{11}\} \oplus \{e_1, e_3, e_6\} = \{e_4, e_6, e_9\}$;
$c_{12} = c_2 \oplus c_3 = \{e_1, e_3, e_4, e_7, e_{10}, e_{11}\} \oplus \{e_1, e_3, e_5, e_{10}, e_{11}\} = \{e_4, e_5, \mathbf{e}_7\}$

имеют меньшую суммарную длину $l = l_7 + l_8 + l_9 + l_{10} + l_{11} + l_{12} = 4 + 3 + 3 + 3 + 3 + 3 = 19$.

Будем рассматривать такие базисные циклы, суммарная длина которых минимальна.

## 3.2. Метрика графов. Расстояние в графе

Введем понятие расстояния между двумя вершинами графа.

*Расстоянием* $\rho$(x,y) графе **G** между вершинами x и y графа G = (V,E) называется длина кратчайшего из маршрутов (а значит – кратчайшей из простых цепей), соединяющих эти вершины; если x и y отделены в G, то $\rho$ (x,y) = $+\infty$. Функция $\rho = \rho$(x,y) определенная на множестве всех пар вершин графа G и принимающая целые неотрицательные значения (к числу которых относится и бесконечное), заслуживает названия метрики графа, поскольку она удовлетворяет трем аксиомам Фреше [9]:

$$\forall x, y \in X[\rho(x, y) = 0 \Leftrightarrow x = y], \tag{3.1}$$

$$\forall x, y \in X[\rho(x, y) = \rho(y, x_1)], \tag{3.2}$$

$$\forall x, y, z \in X[\rho(x, y) + \rho(y, z) \geq \rho(x, z)]. \tag{3.3}$$

Выполнение первых двух аксиом тривиально, проверим третью (неравенство треугольника).

Если вершины x, y или вершины y, z отделены, то по крайней мере одна из двух величин $\rho$ (x, y) и $\rho$ (y, z) есть $\infty$. Если же ни x и y, ни z и y не отделены, то пусть

$x\ e_1\ x_1\ e_2\ x_2 \ldots x_{\rho(x, y)-1}\ e_{\rho(x, y)}\ y$;

и

$y\ u_1\ y_1\ u_2\ y_2 \ldots y_{\rho(y, z)-1}\ u_{\rho(y, z)}\ z$

– какие-либо из кратчайших цепей, соединяющих эти пары вершин.

Маршрут

$x\ e_1\ x_1\ e_2\ x_2 \ldots x_{\rho(x,y)-1}\ e_{\rho(x,y)}\ y\ u_1\ y_1\ u_2\ y_2 \ldots y_{\rho(y,z)-1}\ u_{\rho(y,z)}\ z$

обладает длиной $\rho$(x,y) + $\rho$(y,z), значит длина $\rho$(x,z) кратчайшей цепи между x и z не превышает $\rho$(x,y) + $\rho$(y,z). Таким образом, в обоих случаях неравенство треугольника выполнено.

Введем следующее понятие, связанное с метрикой графа

**Определение 3.1**[51]. *Изометрический подграф* – подграф $G^*$ графа G, у которого все расстояния внутри $G^*$ те же самые, что и в G.

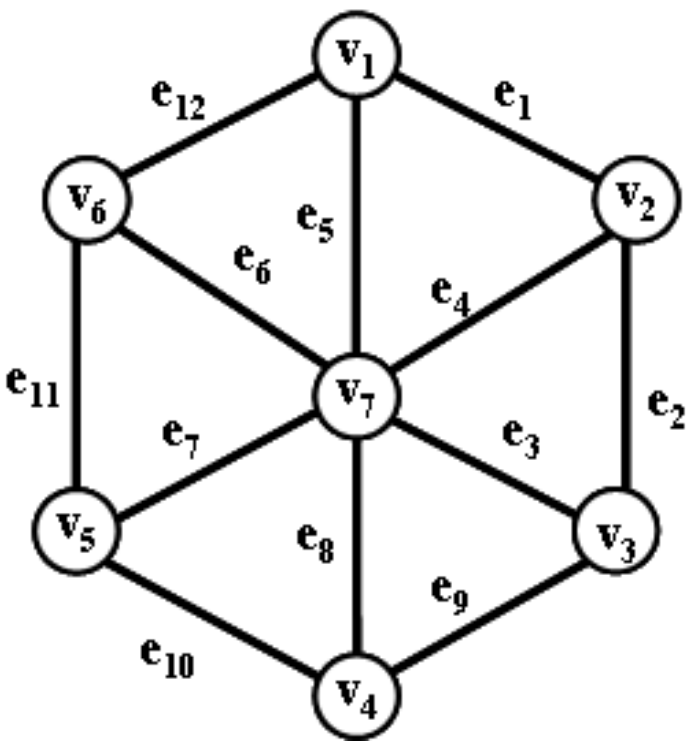

Рис. 3.2. Граф G

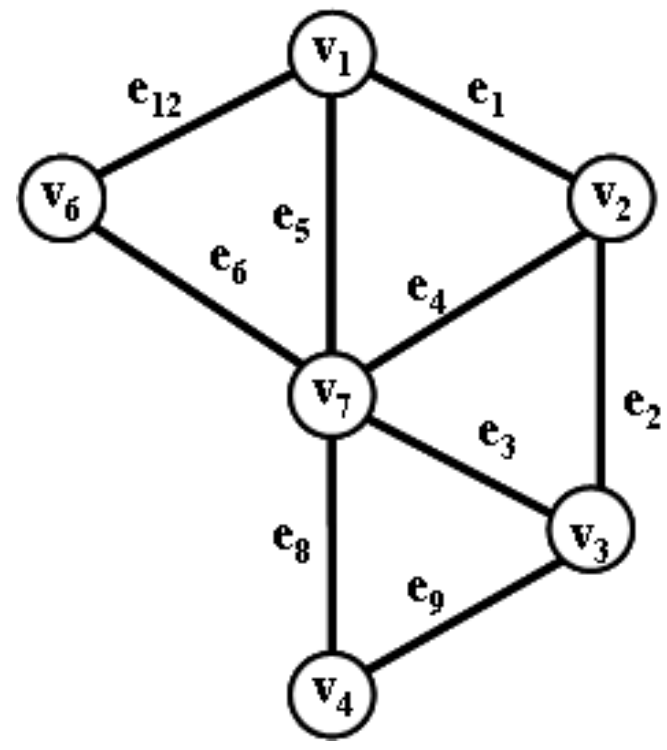

Рис. 3.3. Изометрический подграф $G'$

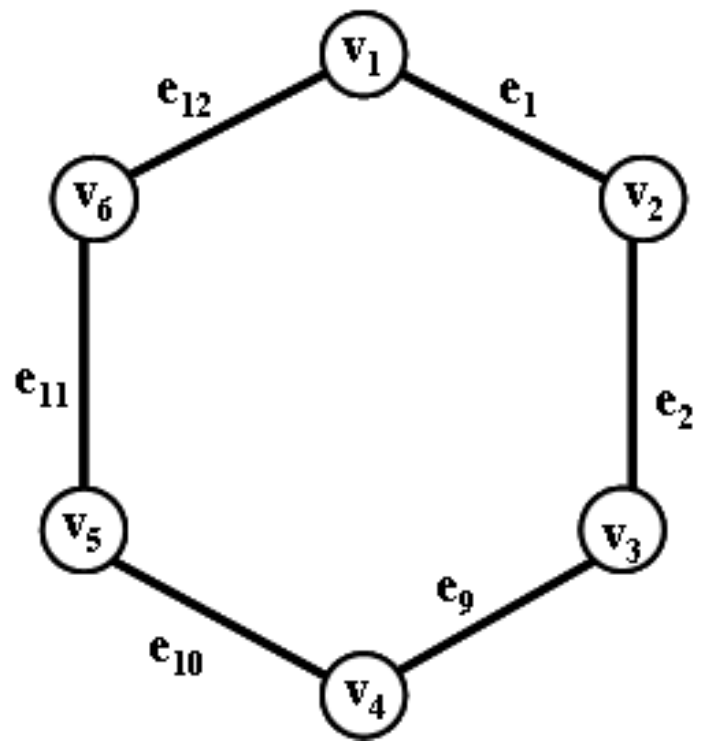

Рис. 3.4. Неизометрический подграф $G'$

На рис. 3.2. представлен граф G, имеющий следующую матрицу расстояний:

$\rho =$

| | $v_1$ | $v_2$ | $v_3$ | $v_4$ | $v_5$ | $v_6$ | $v_7$ |
|---|---|---|---|---|---|---|---|
| $v_1$ | | 1 | 2 | 2 | 2 | 1 | 1 |
| $v_2$ | 1 | | 1 | 2 | 2 | 2 | 1 |
| $v_3$ | 2 | 1 | | 1 | 2 | 2 | 1 |
| $v_4$ | 2 | 2 | 1 | | 1 | 2 | 1 |
| $v_5$ | 2 | 2 | 2 | 1 | | 1 | 1 |
| $v_6$ | 1 | 2 | 2 | 2 | 1 | | 1 |
| $v_7$ | 1 | 1 | 1 | 1 | 1 | 1 | |

Подматрица расстояний для изометрического подграфа представленного на рис. 3.3.

$\rho =$

| | $v_1$ | $v_2$ | $v_3$ | $v_4$ | $v_6$ | $v_7$ |
|---|---|---|---|---|---|---|
| $v_1$ | | 1 | 2 | 2 | 1 | 1 |
| $v_2$ | 1 | | 1 | 2 | 2 | 1 |
| $v_3$ | 2 | 1 | | 1 | 2 | 1 |
| $v_4$ | 2 | 2 | 1 | | 2 | 1 |
| $v_6$ | 1 | 2 | 2 | 2 | | 1 |
| $v_7$ | 1 | 1 | 1 | 1 | 1 | |

Подматрица расстояний для неизометрического подграфа представленного на рис. 3.4.

$\rho =$

| | $v_1$ | $v_2$ | $v_3$ | $v_4$ | $v_5$ | $v_6$ |
|---|---|---|---|---|---|---|
| $v_1$ | | 1 | 2 | 3 | 2 | 1 |
| $v_2$ | 1 | | 1 | 2 | 3 | 2 |
| $v_3$ | 2 | 1 | | 1 | 2 | 3 |
| $v_4$ | 3 | 2 | 1 | | 1 | 2 |
| $v_5$ | 2 | 3 | 2 | 1 | | 1 |
| $v_6$ | 1 | 2 | 3 | 2 | 1 | |

Как видим, матрица расстояний для подграфа на рис. 3.3 является подматрицей расстояний графа G, а матрица расстояний для подграфа на рис. 3.4 не является подматрицей расстояний графа G.

## 3.3. Множества изометрических циклов и центральных разрезов

Что касается базисной изометрической системы циклов и центральных разрезов, принцип их выбора отличается от выбора фундаментальной системы циклов и разрезов, так как понятие выделенного дерева здесь не несет полезной информации.

Центральные разрезы образуют множество размерностью равное количеству вершин графа G. Обозначим данное множество символом $S_e$.

**Свойство 3.1.** Множество центральных разрезов графа обладает следующим свойством: кольцевая сумма всех центральных разрезов для графа G с $n$ вершинами есть пустое множество:

$$\sum_{i=1}^{n} s_i = \varnothing. \qquad (3.4)$$

Количество центральных разрезов равное рангу графа **G** определяет базис подпространства разрезов. Алгоритм выделения множества центральных разрезов довольно прост, он может быть построен перечислением инцидентных ребер для *n-1* вершин графа **G** за линейное время.

Для графа, представленного на рис. 3.5, множество центральных разрезов $\mathbf{S_e}$ = $\{s_1,s_2,s_3,s_4,s_5,s_6\}$, где:

$s_1 = \{e_1,e_2,e_3\}$; $s_2 = \{e_1,e_4,e_5,e_6\}$; $s_3 = \{e_4,e_7,e_8,e_9\}$; $s_4 = \{e_2,e_5,e_7,e_9\}$;
$s_5 = \{e_8,e_{10},e_{11}\}$; $s_6 = \{e_3,e_6,e_{10},e_{11}\}$.

Любой суграф, принадлежащий подпространству циклов C, в общем случае является квазициклом. *Простые циклы* – это *квазициклы*, у которых локальная степень вершин в точности равна двум [9,10].

Мощность подмножества простых циклов в графе меньше мощности множества квазициклов. Подмножество простых циклов обозначим $C_R$:

$$\text{card } C_R \le \text{card } C \qquad (3.5)$$

Однако существует подмножество с мощностью еще меньшей, чем подмножество простых циклов, обладающее определенными характерными свойствами.

**Определение 3.2.** *Изометрическим циклом* в графе называется простой цикл, для которого кратчайший путь между любыми двумя его вершинами состоит из рёбер этого цикла. Изометрический цикл – частный случай изометрического подграфа [16,47].

Или, другими словами, изометрическим циклом в графе называется подграф G′ в виде простого цикла, если между двумя любыми несмежными вершинами данного подграфа в соответствуюшем графе G не существует маршрутов меньшей длины, чем маршруты, принадлежащие данному циклу.

Подмножество, состоящее из изометрических циклов, назовем подмножеством изометрических циклов и обозначим $C_\tau$. Сказанное поясним на примерах. Рассмотрим суграф, состоящий из ребер $\{e_1,e_3,e_{13},e_{15}\}$ графа $G_a$, представленного на рис. 3.5,а. Как видно, это простой цикл. Но в то же время это не изометрический цикл, так как между вершинами $x_7$ и $x_8$ в графе существуют маршрут меньшей длины, проходящий по ребру $e_{14}$.

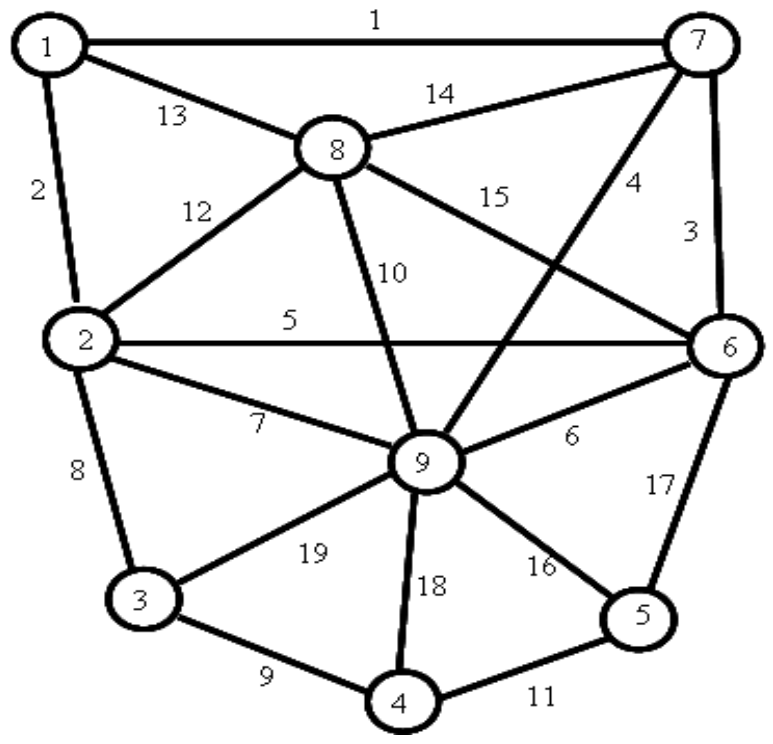

а) Граф $G_a$

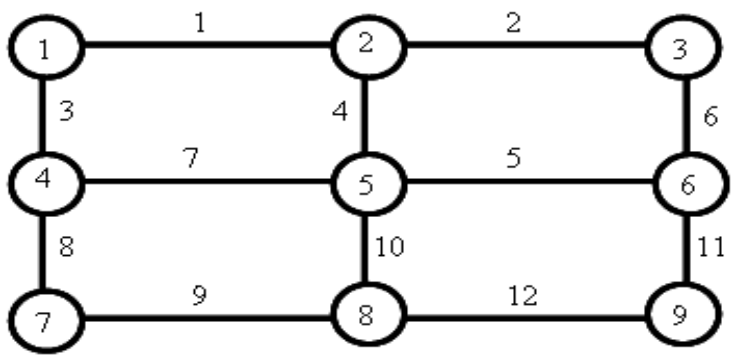

б) Граф $G_b$

Рис. 3.5. Графы $G_a$ и $G_b$.

Рассмотрим граф $G_b$, представленный на рис. 3.5,б. Пусть цикл состоит из ребер $e_1,e_2,e_3,e_6,e_8,e_9,e_{11},e_{12}$. Данный суграф есть простой цикл. Однако этот суграф не может быть изометрическим циклом, так как в соответствующем графе между вершинами $x_2$ и $x_8$ имеется маршрут меньшей длины (а именно, маршрут, проходящий по ребрам $e_4$ и $e_{10}$), чем маршруты, принадлежащие этому суграфу (например, маршрут, проходящий по ребрам $e_1,e_3,e_8,e_9$ или $e_2,e_6,e_{11},e_{12}$).

Следует заметить, что в полных графах множество изометрических циклов совпадает с множеством циклов минимальной длины. В целях сокращения записи иногда будем обозначать вершины и ребра целыми числами.

Для изучения свойств изометрических циклов нам понадобится следующая теорема.

**Теорема 3.1.** Линейное подпространство квазициклов несепарабельного графа имеет базис, состоящий из независимого подмножество изометрических циклов с мощностью равной цикломатическому числу графа.

*Доказательство.* Рассмотрим множество деревьев M(T) графа G.

Рассмотрим систему фундаментальных циклов порожденных деревом **T**. Пусть вершины $A_1,A_2,...,A_p$ образуют фундаментальный цикл содержащий хорду ($A_p$,$A_1$). Если между несмежными вершинами этого цикла в графе не существует путей меньшей или равной длины, чем пути, принадлежащие циклу – то это изометрический цикл. Если в цикле существуют две несмежные вершины графа $A_i$ и $A_j$ ($i<j$), а в графе существует путь меньшей длины для выбранных вершин, чем путь по циклу $A_i,B_1,B_2,...,B_r,A_j$, то образуются циклы, кольцевая сумма которых есть исходный цикл. Оставляем цикл, содержащий хорду. Это и есть изометрический цикл. Так как количество фундаментальных циклов независимо и определяется цикломатическим числом, то количество изометрических циклов, полученных описанным выше способом, также независимо и равно цикломатическому числу графа включая все хорды. *Теорема доказана.*

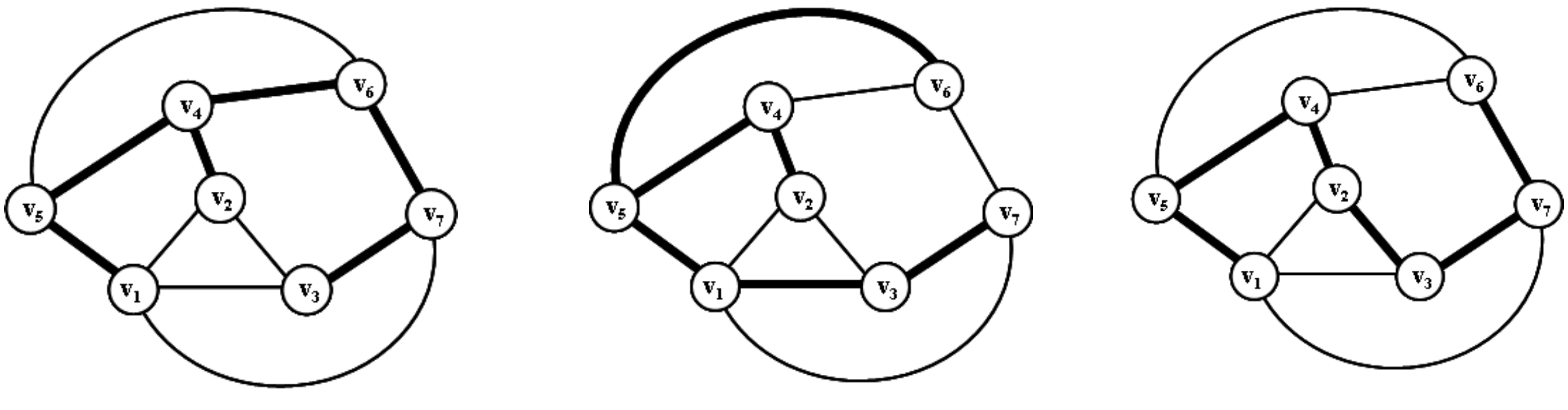

Рис. 3.6. К теореме 3.1.

### 3.4. Методы выделения множества изометрических циклов графа

Понятие изометрического цикла графа G тесно связано с минимальными (s-t) маршрутами графа [36].

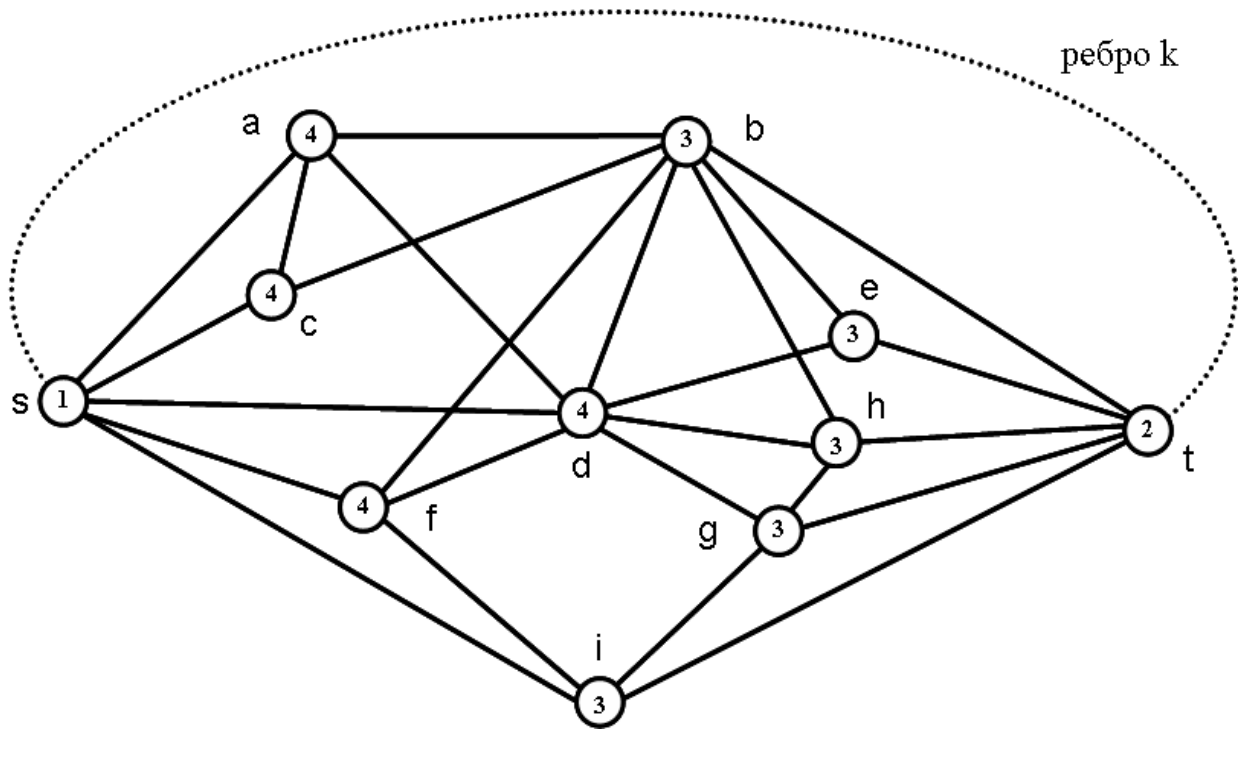

Рис. 3.7. Прямая разметка вершин для ребра k.

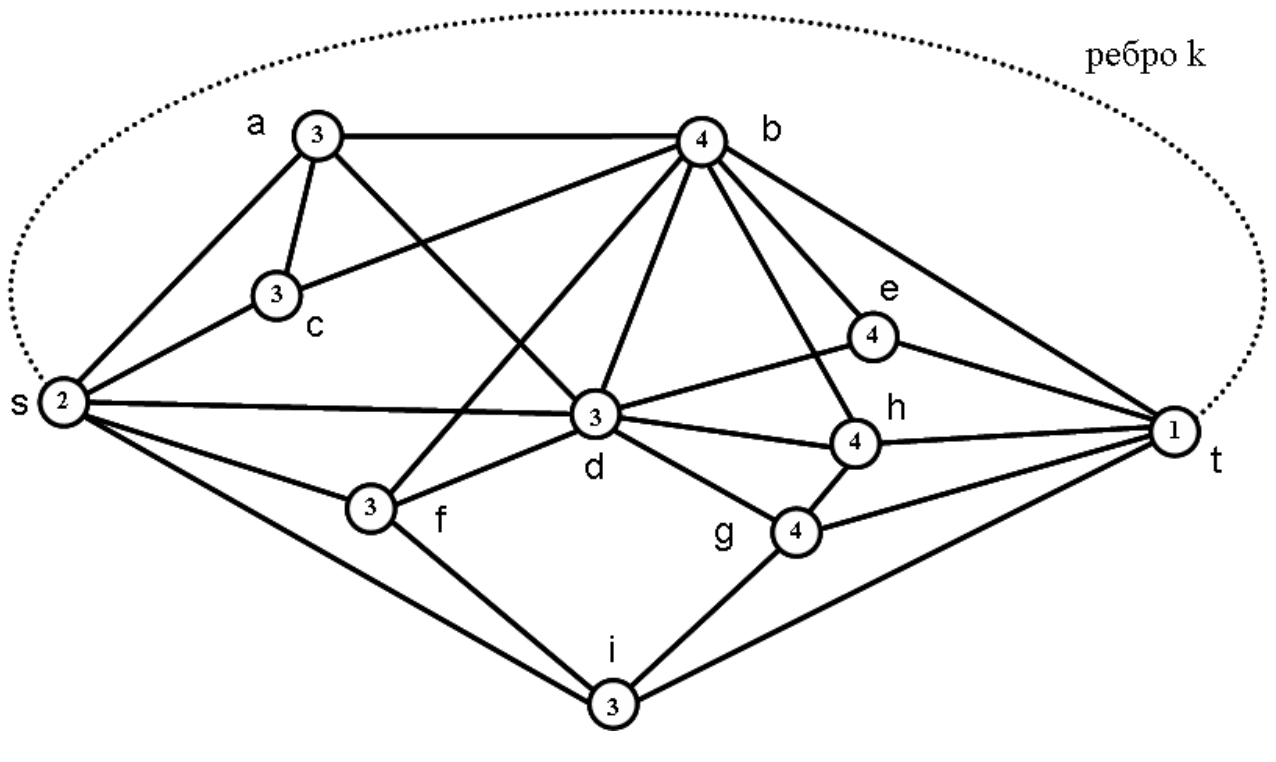

Рис. 3.8. Обратная разметка вершин для ребра k.

С этой целью рассмотрим изометрические циклы, проходящие по k-му ребру соединяющему вершины s и t графа G (рис. 3.7). Удалим из графа ребро k. Получим граф

**G**-k, где вершины s и t теперь несмежны. Применим алгоритм поиска в ширину. Вершине s поставим в соответствие фронт волны 1, а вершине t – фронт волны 2. Тогда вершины, смежные с вершиной 2 и еще не помеченные пометим цифрой 3 и так далее. Другими словами, применим алгоритм поиска в ширину относительно вершины t. Выделим все простые цепи минимальной длины, образованные алгоритмом поиска в ширину, осуществляя проход от вершины с большим номером к вершине с меньшим номером (рис. 3.7). Сформируем множество $C_{st}$ циклов, где элементами множеств являются вершины:

$C_{st}$ = {{s,a,b,t},{s,c,b,t},{s,f,b,t},{s,d,b,t},{s,d,e,t},{s,d,h,t},{s,d,g,t},{s,f,b,t},
{s,f,i,t},{s,i,t}}.

Теперь вершине t поставим в соответствие фронт волны 1, а вершине s – фронт волны 2. Тогда вершины, смежные с вершиной 2 и еще не помеченные пометим цифрой 3, и так далее Другими словами, применим алгоритм поиска в ширину относительно вершины s. Выделим все простые цепи минимальной длины, образованные алгоритмом поиска в ширину, осуществляя проход от вершины с большим номером к вершине с меньшим номером (рис. 3.8). Сформируем множество $C_{ts}$ циклов:

$C_{ts}$ = {{s,a,b,t},{s,c,b,t},{s,f,b,t},{s,d,b,t},{s,d,e,t},{s,d,h,t},{s,d,g,t},{s,f,b,t},
{s,i,g,t},{s,i,t}}.

Циклы {s,a,b,t},{s,c,b,t},{s,f,b,t},{s,d,b,t},{s,d,e,t},{s,d,h,t},{s,d,g,t},{s,f,b,t},{s,i,t} принадлежащие одновременно и множеству $M_{st}$ и множеству $M_{ts}$ суть изометрические циклы. Циклы {s,f,i,t} и {s,i,g,t} принадлежат только одному из множеств и поэтому не являются изометрическими циклами.

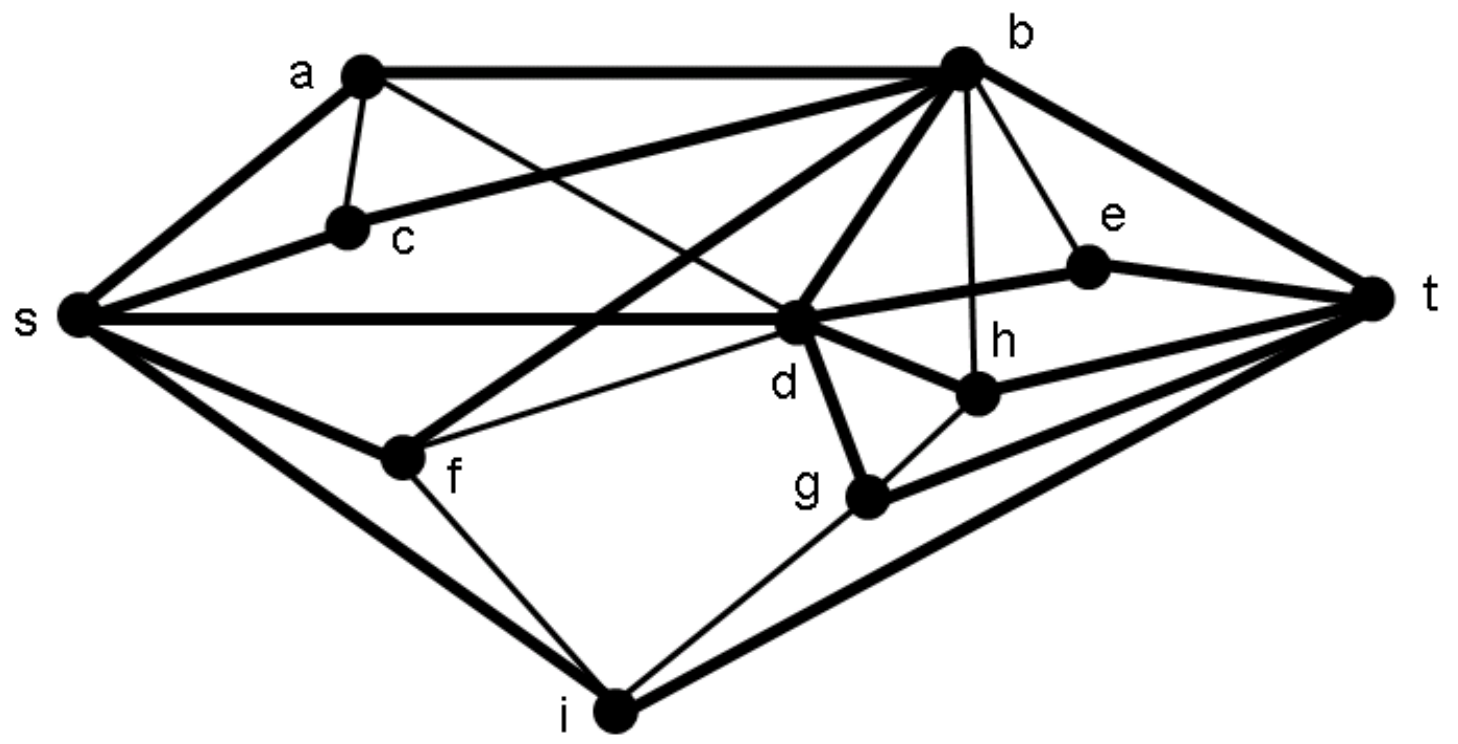

Рис. 3.9. Минимальные s-t цепи

Данные рассуждения можно применить ко всем ребрам графа G и сформировать множество изометрических циклов графа.

Вычислительная сложность такого алгоритма определится громоздкостью операции сравнения циклов выделенных при прямой и обратной разметке вершин. Здесь максимальное количество циклов определяется как (n-2) для каждой вершины, а количество сравнений тогда будет равно $(n-2)^2$. Полученное выражение нужно умножить на количество ребер = n(n-1)/2. Окончательно получим, что сложность вычисления можно определить по ыормуле:

$$f(n) = n(n-1)(n-2)^2/2. \quad (3.6)$$

Откуда вычислительную сложность можно определить как O(n4)/

Количество изометрических циклов в полном графе определяется по формуле

$$card(C_e) = n(n-1)(n-2)/6. \quad (3.7)$$

Рассмотрим другой способ построения изометрических циклов. Выделим в полном графе множество циклов длиной три. Это, очевидно, будет множество изометрических циклов для полного графа.

Будем последовательно удалять ребра из полного графа $K_n$. Естественно, что тогда будут удалены и изометрические циклы, содержащие данное ребро, или два удаляемых цикла образуют новый изометрический цикл равный их кольцевой сумме в случае, если вновь образованный цикл не образован оставшимися изометрическими циклами. Таким образом, процесс удаления ребер из полного графа приводит к уменьшению количества изометрических циклов в графе, и, естественно, что их количество в произвольном графе не может превышать величины n(n-1)(n-2)/6.

Сказанное рассмотрим на примере графа $K_5$. Множество изометрических циклов $\mathbf{C}_\tau$ для графа $K_5$ (рис. 3.10):

$c_1 = \{e_1,e_2,e_5\}$; $c_2 = \{e_1,e_3,e_6\}$; $c_3 = \{e_1,e_4,e_7\}$; $c_4 = \{e_2,e_3,e_8\}$; $c_5 = \{e_2,e_4,e_9\}$;
$c_6 = \{e_3,e_4,e_{10}\}$; $c_7 = \{e_5,e_6,e_8\}$; $c_8 = \{e_5,e_7,e_9\}$; $c_9 = \{e_6,e_7,e_{10}\}$;
$c_{10} = \{e_8,e_9,e_{10}\}$.

Удалим из графа ребро $e_{10}$. Тогда из множества изометрических циклов $C_\tau$ удаляются все циклы включающие 10-ое ребро:

$c_6 = \{e_3,e_4,e_{10}\}$; $c_9 = \{e_6,e_7,e_{10}\}$; $c_{10} = \{e_8,e_9,e_{10}\}$.

Остаются изометрические циклы:

$c_1 = \{e_1,e_2,e_5\}$; $c_2 = \{e_1,e_3,e_6\}$; $c_3 = \{e_1,e_4,e_7\}$; $c_4 = \{e_2,e_3,e_8\}$; $c_5 = \{e_2,e_4,e_9\}$;
$c_7 = \{e_5,e_6,e_8\}$; $c_8 = \{e_5,e_7,e_9\}$.

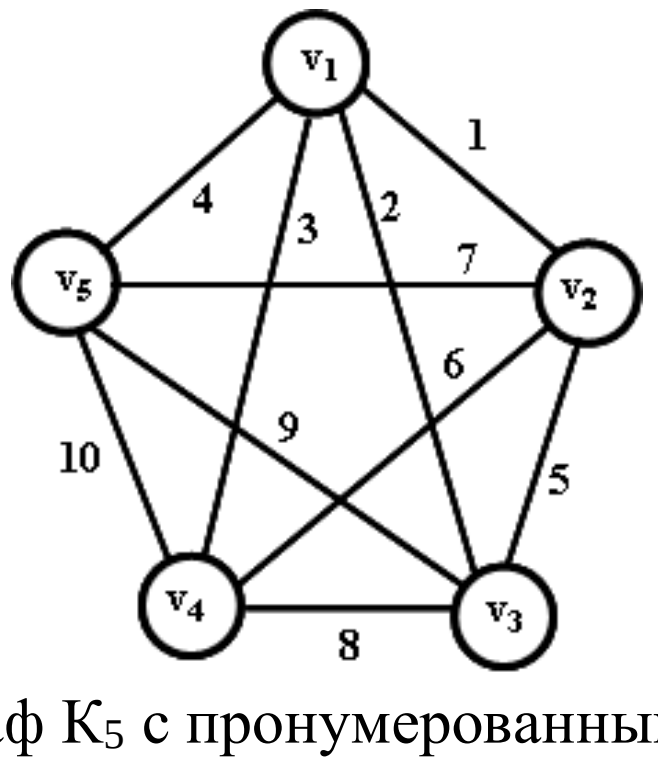

Рис. 3.10. Граф $K_5$ с пронумерованными ребрами.

В перспективе должны образоваться новые изометрические циклы длиной четыре,

образованные из удаленных циклов:

$c_6 \oplus c_9 = \{e_3,e_4,e_{10}\} \oplus \{e_6,e_7,e_{10}\} = \{e_3,e_4,e_6,e_7\}$;
$c_6 \oplus c_{10} = \{e_3,e_4,e_{10}\} \oplus \{e_8,e_9,e_{10}\} = \{e_3,e_4,e_8,e_9\}$;
$c_9 \oplus c_{10} = \{e_6,e_7,e_{10}\} \oplus \{e_8,e_9,e_{10}\} = \{e_6,e_7,e_8,e_9\}$.

Однако их включение во множество оставшихся изометрических циклов невозможно, так как они могут быть образованы как результат кольцевого суммирования из оставшихся изометрических циклов:

$c_2 \oplus c_3 = \{e_1,e_3,e_6\} \oplus \{e_1,e_4,e_7\} = \{e_3,e_4,e_6,e_7\}$;
$c_4 \oplus c_5 = \{e_2,e_3,e_8\} \oplus \{e_2,e_4,e_9\} = \{e_3,e_4,e_8,e_9\}$;
$c_7 \oplus c_8 = \{e_5,e_6,e_8\} \oplus \{e_5,e_7,e_9\} = \{e_6,e_7,e_8,e_9\}$.

Если мы продолжим удаление ребра $e_2$ из графа, то из множества изометрических циклов $C_\tau = \{c_1,c_2,c_3,c_4,c_5,c_7,c_8\}$ удаляются все циклы включающие 2-ое ребро:

$c_1 = \{e_1,e_2,e_5\}$; $c_4 = \{e_2,e_3,e_8\}$; $c_5 = \{e_2,e_4,e_9\}$.

Остаются изометрические циклы:

$c_2 = \{e_1,e_3,e_6\}$; $c_3 = \{e_1,e_4,e_7\}$; $c_7 = \{e_5,e_6,e_8\}$; $c_8 = \{e_5,e_7,e_9\}$.

В перспективе должны образоваться новые изометрические циклы длиной четыре:

$c_1 \oplus c_4 = \{e_1,e_2,e_5\} \oplus \{e_2,e_3,e_8\} = \{e_1,e_3,e_5,e_8\}$;
$c_1 \oplus c_5 = \{e_1,e_2,e_5\} \oplus \{e_2,e_4,e_9\} = \{e_1,e_4,e_5,e_9\}$;
$c_4 \oplus c_5 = \{e_2,e_3,e_8\} \oplus \{e_2,e_4,e_9\} = \{e_3,e_4,e_8,e_9\}$.

Однако включение во множество оставшихся изометрических циклов двух первых невозможно, так как они могут быть образованы как результат кольцевого суммирования из оставшихся изометрических циклов:

$c_1 \oplus c_4 = \{e_1,e_3,e_6\} \oplus \{e_5,e_6,e_8\} = \{e_1,e_3,e_5,e_8\}$;
$c_1 \oplus c_5 = \{e_1,e_4,e_7\} \oplus \{e_5,e_7,e_9\} = \{e_1,e_4,e_5,e_9\}$.

Вновь образованный изометрический цикл $c_{4,5}$ включается во множество оставшихся изометрических циклов: $c_4 \oplus c_5 = \{e_3,e_4,e_8,e_9\}$.

Таким образом, множество изометрических циклов $C_\tau$ для графа полученного путем удаления 10-го и 2-го ребер из графа $K_5$ состоит из следующих изометрических циклов:

$c_2 = \{e_1,e_3,e_6\}$; $c_3 = \{e_1,e_4,e_7\}$; $c_7 = \{e_5,e_6,e_8\}$; $c_8 = \{e_5,e_7,e_9\}$;
$c_{4,5} = \{e_3,e_4,e_8,e_9\}$.

Характерная и особая роль изометрических циклов в теории графов определяется тем, что в плоских графах они являются границами граней. В свою очередь, характерная особенность центральных разрезов проявляется в том, что их длина определяет локальную степень вершин.

Ввиду важности вопроса выделения конечного множества изометрических циклов из множества квазициклов, предлагается алгоритм выделения множества изометрических циклов в графе.

Построение алгоритма начинается с выделения всех ребер в графе G. Выберем очередное ребро графа. Одну из вершин такого выбранного ребра пометим индексом 1, другую – индексом 2. Вершины графа смежные с вершиной, имеющей индекс 2, и ещё не помеченные, пометим индексом 3. Вершины графа смежные с вершиной, имеющей индекс 3, и ещё не помеченные, пометим индексом 4 и т. д. Число, выражающее индекс последней помеченной вершины (вершин) графа, называется глубиной проникновения волны, относительно выбранного ребра. Данный процесс представляет собой разметку вершин графа, относительно выбранного ребра волновым алгоритмом (алгоритмом поиска в ширину).

Построим простые циклы, проходящие по выбранному ребру, относительно первоначальной ориентации. С этой целью выберем все вершины графа G смежные с вершиной, помеченной индексом 1. Будем идти от любой выбранной вершины, имеющей глубину проникновения d, к вершинам, имеющим глубину проникновения (d-1), проходя при этом по ребрам графа, затем от вершины (d-1) к вершинам (d-2) и т.д. Остановим этот процесс тогда, когда подойдем к вершине, имеющей индекс 2. Пройдя по всем таким образом построенным маршрутам, построим систему циклов, проходящих по выбранному ребру j. Обозначим такое множество циклов через $S_j^1$. Переориентируем направление разметки, т.е. вершина, имеющая индекс 1, будет иметь индекс 2, а вершина, имеющая индекс 2, будет иметь индекс 1. И вновь построим разметку вершин. Описанным выше методом выделим систему циклов. Изометрические циклы, проходящие по выбранному ребру j, будут образованы как:

$$C_j = C_j^1 \cap C_j^2 . \qquad (3.8)$$

Множество изометрических циклов графа G будет образовано как объединение всех циклов, проходящих по всем ребрам графа:

$$C_\tau = \bigcup_{i=1}^{m} C_i \ (i = 1,2,...,m). \qquad (3.9)$$

## 3.5. Алгоритм выделения множества изометрических циклов графа

**Алгоритм 3.1. [Выделение множества изометрических циклов методом поиска в ширину]**

**Шаг 1.** [**Выбор ребра**]. Выбираем ребро, идем на шаг 2. Если количество ребер исчерпано, то конец работы алгоритма.

**Шаг 2**. [**Прямая разметка вершин относительно вершины s ребра**]. Алгоритмом поиска в ширину производим прямую разметку вершин относительно вершины s выбранного ребра. Идем на шаг 3.

**Шаг 3**. [**Формирование множества циклов $C_s$ при прямой разметке вершин**]. Производим формирование множества циклов при прямой разметке вершин. Идем на шаг 4.

**Шаг 4**. [**Обратная разметка вершин относительно вершины t ребра**]. Алгоритмом поиска в ширину производим обратную разметку вершин относительно вершины t выбранного ребра. Идем на шаг 5.

**Шаг 5**. [**Формирование множества циклов $C_t$ при обратной разметке вершин**]. Производим формирование множества циклов при обратной разметке вершин. Идем на шаг 6.

**Шаг 6.** [**Проверка циклов $C_s$ и $C_t$ на совпадение**]. Приверяем циклы $C_s$ и $C_t$ для выбранного ребра на совпадение. Несовпадающие циклы исключаем из рассмотрения. Идем на шаг 7.

**Шаг 7.** [**Запись циклов во множество изометрических циклов**]. Проверяем сформированные на предыдущем шаге циклы с ранее записанными циклами во множестве изометрических циклов и, в случае их отсутствия, добавляем их во множество изометрических циклов. Идем на шаг 1.

***Пример 3.1.*** В качестве примера рассмотрим граф G (рис. 3.11).

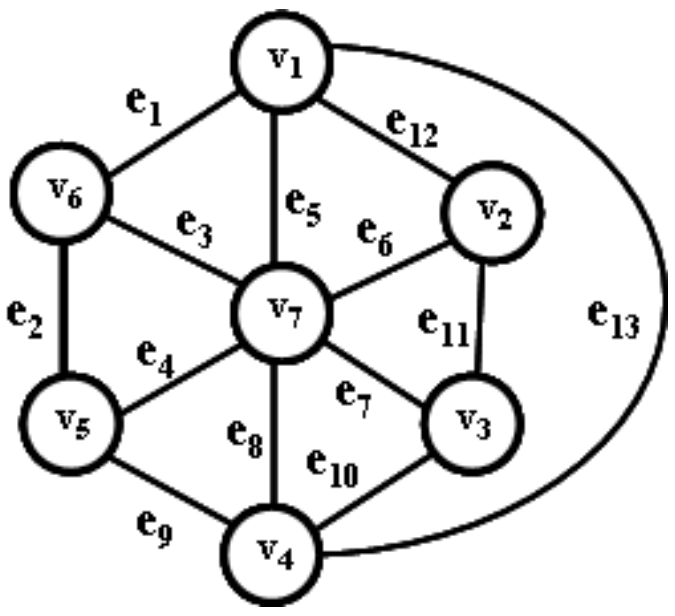

Рис. 3.11. Граф $G_6$.

Если выбрать ребро $e_{13}$, то процесс разметки вершин имеет вид, представленный на рис. 3.12. Система циклов проходящих по ребру $e_{13}$ для разметки, показанной на рис. 3.12,*а* (номера разметки представлены внутри вершины):

$$C_{13}^{1} = \{\{e_5,e_8,e_{13}\}, \{e_1,e_3,e_8,e_{13}\}, \{e_1,e_2,e_9,e_{13}\}, \{e_6,e_8,e_{12},e_{13}\},\{e_{10},e_{11},e_{12},e_{13}\}\}.$$

Система циклов проходящих по ребру $e_{13}$ для разметки, представленной на рис. 3.12,б (номера разметки представлены внутри вершины):

$$C_{13}^{2} = \{\{e_5,e_8,e_{13}\},\{e_4,e_5,e_9,e_{13}\},\{e_1,e_2,e_9,e_{13}\},\{e_{10},e_{11},e_{12},e_{13}\},\{e_5,e_7,e_{10},e_{13}\}\}.$$

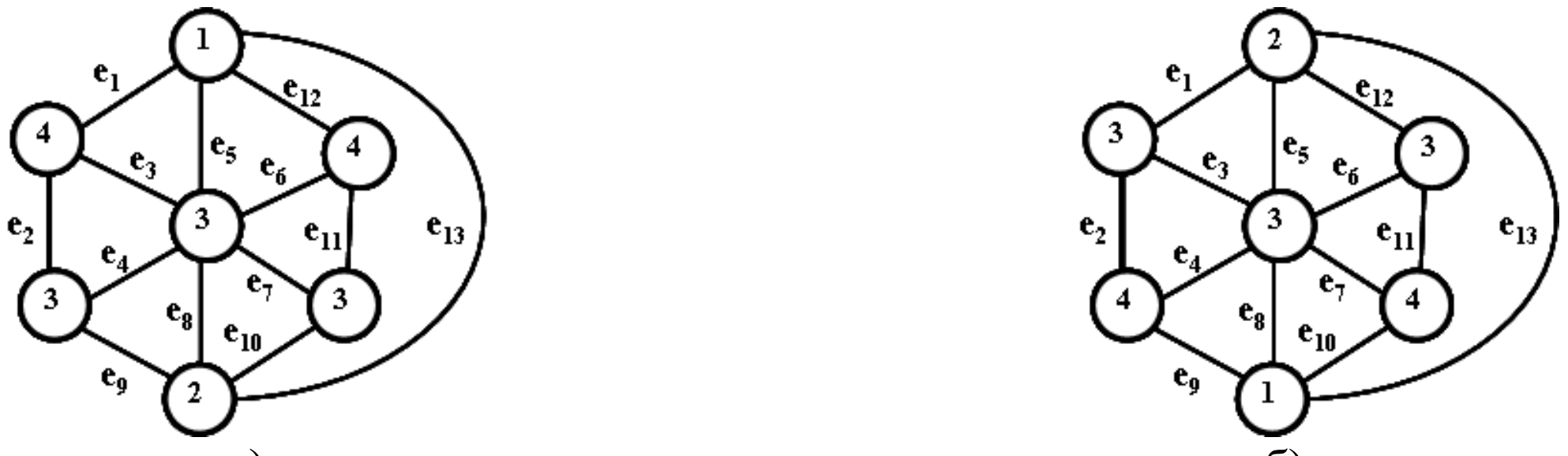

а) б)

Рис. 3.12. Прямой и обратный процесс разметки вершин (числа в вершинах фронт волны)

Пересечение множеств $C^1_{13}$ и $C^2_{13}$:

$C_{13} = C^1_{13} \cap C^2_{13} = \{\{e_5,e_8,e_{13}\},\{e_1,e_2,e_9,e_{13}\},\{e_{10},e_{11},e_{12},e_{13}\}\}$.

Каждому ребру принадлежат следующие изометрические циклы:

$C_1 = \{\{e_1,e_3,e_5\},\{e_1,e_2,e_9,e_{13}\}\}$;
$C_2 = \{\{e_2,e_3,e_4\},\{e_1,e_2,e_9,e_{13}\}\}$;
$C_3 = \{\{e_2,e_3,e_4\},\{e_1,e_3,e_5\}\}$;
$C_4 = \{\{e_4,e_8,e_9\},\{e_2,e_3,e_4\}\}$;
$C_5 = \{\{e_1,e_3,e_5\},\{e_5,e_6,e_{12}\},\{e_5,e_8,e_{13}\}\}$;
$C_6 = \{\{e_5,e_6,e_{12}\},\{e_6,e_7,e_{11}\}\}$;
$C_7 = \{\{e_6,e_7,e_{11}\},\{e_7,e_8,e_{10}\}\}$;
$C_8 = \{\{e_5,e_8,e_{13}\},\{e_7,e_8,e_{10}\},\{e_3,e_8,e_9\}\}$;
$C_9 = \{\{e_4,e_8,e_9\},\{e_1,e_2,e_9,e_{13}\}\}$;
$C_{10} = \{\{e_7,e_8,e_{10}\},\{e_{10},e_{11},e_{12},e_{13}\}\}$;
$C_{11} = \{\{e_6,e_7,e_{11}\},\{e_{10},e_{11},e_{12},e_{13}\}\}$;
$C_{12} = \{\{e_5,e_6,e_{12}\},\{e_{10},e_{11},e_{12},e_{13}\}\}$;
$C_{13} = \{\{e_5,e_8,e_{13}\},\{e_1,e_2,e_9,e_{13}\},\{e_{10},e_{11},e_{12},e_{13}\}\}$.

Множество изометрических циклов получим как объединение:

$C_\tau = C_1 \cup C_2 \cup C_3 \cup C_4 \cup C_5 \cup C_6 \cup C_7 \cup C_8 \cup C_9 \cup C_{10} \cup C_{11} \cup C_{12} \cup C_{13} =$
$= \{\{e_1,e_3,e_5\}, \{e_2,e_3,e_4\}, \{e_4,e_8,e_9\}, \{e_5,e_6,e_{12}\}, \{e_5,e_8,e_{13}\}, \{e_6,e_7,e_{11}\},$
$\{e_7,e_8,e_{10}\},\{e_1,e_2,e_9,e_{13}\}, \{e_{10},e_{11},e_{12},e_{13}\}\}$.

Таким образом, множество изометрических циклов состоит из 9-ти элементов. Цикломатическое число графа $G_6$ равно 7. Следовательно, для построения базиса нужно удалить два изометрических цикла. Очевидно, что для любого трехсвязного и более графа G множество изометрических циклов имеет мощность меньшую, чем мощность множества простых циклов, но большую или равную цикломатическому числу графа:

$$\nu(G) \le \mathrm{card}\ C_\tau \le \mathrm{card}\ C_R \le \mathrm{card}\ C. \qquad (3.10)$$

Теперь покажем, что построение множества изометрических циклов должно производиться относительно всего множества ребер графа.

***Пример 3.2.*** Следующий пример демонстрирует невозможность получения полного множества изометрических циклов, если построение производится только относительно хорд для выбранного дерева графа. Рассмотрим граф представленный на рис. 3.13.

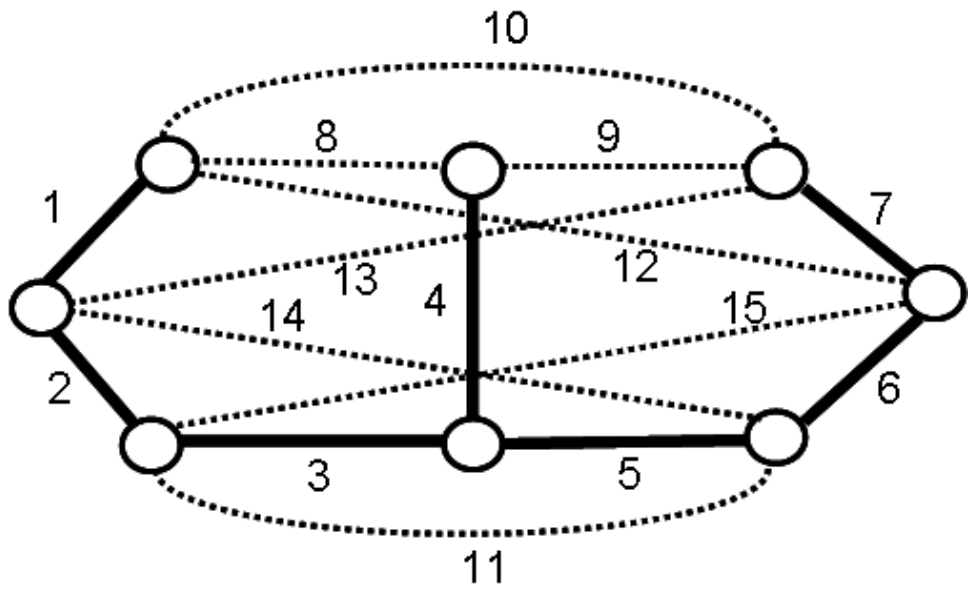

Рис. 3.13. Граф $G_7$ и его дерево.

Рассмотрим изометрические циклы относительно 8-й хорды:

| | |
|---|---|
| $c_1 = \{e_8,e_9,e_{10}\}$; | $c_1' = \{e_8,e_9,e_{10}\}$; |
| $c_2 = \{e_1,e_8,e_9,e_{13}\}$; | $c_2' = \{e_4,e_5,e_6,e_8,e_{12}\}$; |
| $c_3 = \{e_7,e_8,e_9,e_{12}\}$; | $c_3' = \{e_1,e_2,e_3,e_4,e_8\}$; |
| | $c_4' = \{e_4,e_5,e_8,e_{12},e_{15}\}$; |
| | $c_5' = \{e_4,e_5,e_9,e_{13},e_{14}\}$. |

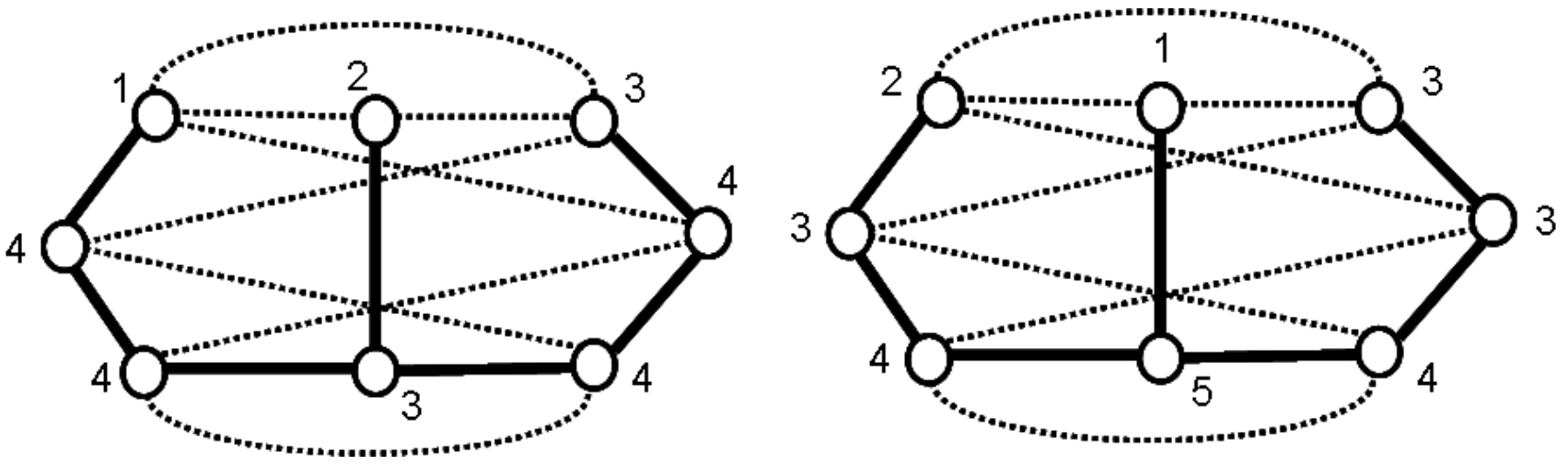

Рис. 3.14. Изометрические циклы относительно 8-й хорды.

Изометрический цикл относительно 8-й хорды: $\{e_8,e_9,e_{10}\}$.

Рассмотрим изометрические циклы относительно 9-й хорды:

| | |
|---|---|
| $c_1 = \{e_8,e_9,e_{10}\}$; | $c_1' = \{e_8,e_9,e_{10}\}$; |
| $c_2 = \{e_4,e_5,e_6,e_8,e_{12}\}$; | $c_2' = \{e_1,e_8,e_9,e_{13}\}$; |
| $c_3 = \{e_2,e_3,e_4,e_9,e_{13}\}$; | $c_3' = \{e_7,e_8,e_9,e_{12}\}$. |
| $c_4 = \{e_3,e_4,e_7,e_9,e_{15}\}$; | |
| $c_5 =$ | |

Изометрический цикл относительно 9-й хорды: $\{e_8,e_9,e_{10}\}$.

Рассмотрим изометрические циклы относительно 10-й хорды:

| | |
|---|---|
| $c_1 = \{e_8,e_9,e_{10}\}$; | $c_1' = \{e_8,e_9,e_{10}\}$; |
| $c_2 = \{e_7,e_{10},e_{12}\}$; | $c_2' = \{e_7,e_{10},e_{12}\}$; |
| $c_3 = \{e_1,e_{10},e_{13}\}$; | $c_3' = \{e_1,e_{10},e_{13}\}$. |

Изометрические циклы относительно 10-й хорды: $\{e_8,e_9,e_{10}\}$, $\{e_1,e_{10},e_{13}\}$, $\{e_7,e_{10},e_{12}\}$.

Рассмотрим изометрические циклы относительно 11-й хорды:

| | |
|---|---|
| $c_1 = \{e_3,e_5,e_{11}\}$; | $c_1' = \{e_3,e_5,e_{11}\}$; |
| $c_2 = \{e_2,e_{11},e_{14}\}$; | $c_2' = \{e_2,e_{11},e_{14}\}$; |
| $c_3 = \{e_6,e_{11},e_{15}\}$; | $c_3' = \{e_6,e_{11},e_{15}\}$. |

Изометрические циклы относительно 11-й хорды: $\{e_3,e_5,e_{11}\}$, $\{e_2,e_{11},e_{14}\}$, $\{e_6,e_{11},e_{15}\}$.

Рассмотрим изометрические циклы относительно 12-й хорды:

| | |
|---|---|
| $c_1 = \{e_7,e_{10},e_{12}\}$; | $c_1' = \{e_7,e_{10},e_{12}\}$; |
| $c_2 = \{e_7,e_8,e_9,e_{12}\}$; | $c_2' = \{e_1,e_2,e_{12},e_{15}\}$; |
| $c_3 = \{e_1,e_2,e_{12},e_{15}\}$; | $c_3' = \{e_1,e_6,e_{12},e_{14}\}$. |
| $c_4 = \{e_1,e_7,e_{12},e_{13}\}$; | |
| $c_5 = \{e_1,e_6,e_{12},e_{14}\}$; | |

Изометрические циклы относительно 12-й хорды: $\{e_7,e_{10},e_{12}\}$, $\{e_1,e_2,e_{12},e_{15}\}$, $\{e_1,e_6,e_{12},e_{14}\}$.

Рассмотрим изометрические циклы относительно 13-й хорды:

| | |
|---|---|
| $c_1 = \{e_1,e_{10},e_{13}\}$; | $c_1' = \{e_1,e_{10},e_{13}\}$; |
| $c_2 = \{e_6,e_7,e_{13},e_{14}\}$; | $c_2' = \{e_1,e_7,e_{12},e_{13}\}$; |
| $c_3 = \{e_2,e_7,e_{13},e_{15}\}$; | $c_3' = \{e_6,e_7,e_{13},e_{14}\}$; |
| | $c_4' = \{e_2,e_7,e_{13},e_{15}\}$; |
| | $c_5' = \{e_1,e_8,e_9,e_{13}\}$. |

Изометрические циклы относительно 13-й хорды: $\{e_1,e_{10},e_{13}\}$, $\{e_6,e_7,e_{13},e_{14}\}$, $\{e_2,e_7,e_{13},e_{15}\}$.

Рассмотрим изометрические циклы относительно 14-й хорды:

| | |
|---|---|
| $c_1 = \{e_2,e_{11},e_{14}\}$; | $c_1' = \{e_2,e_{11},e_{14}\}$; |
| $c_2 = \{e_6,e_7,e_{13},e_{14}\}$; | $c_2' = \{e_2,e_6,e_{14},e_{15}\}$; |
| $c_3 = \{e_1,e_6,e_{12},e_{14}\}$; | $c_3' = \{e_2,e_3,e_5,e_{14}\}$; |
| | $c_4' = \{e_6,e_7,e_{13},e_{14}\}$; |
| | $c_5' = \{e_1,e_6,e_{12},e_{14}\}$. |

Изометрические циклы относительно 14-й хорды: $\{e_2,e_{11},e_{14}\}$, $\{e_6,e_7,e_{13},e_{14}\}$, $\{e_1,e_6,e_{12},e_{14}\}$.

Рассмотрим изометрические циклы относительно 15-й хорды:

| | |
|---|---|
| $c_1 = \{e_6,e_{11},e_{15}\}$; | $c_1' = \{e_6,e_{11},e_{15}\}$; |
| $c_2 = \{e_3,e_5,e_6,e_{15}\}$; | $c_2' = \{e_2,e_7,e_{13},e_{15}\}$; |
| $c_3 = \{e_2,e_6,e_{14},e_{15}\}$; | $c_3' = \{e_1,e_2,e_{12},e_{15}\}$. |
| $c_4 = \{e_2,e_7,e_{13},e_{15}\}$; | |
| $c_5 = \{e_1,e_2,e_{12},e_{15}\}$; | |

Изометрические циклы относительно 15-й хорды: $\{e_6,e_{11},e_{15}\}$, $\{e_2,e_7,e_{13},e_{15}\}$, $\{e_1,e_2,e_{12},e_{15}\}$.

Рассмотрим изометрические циклы относительно 4-го ребра:

| | |
|---|---|
| $c_1 = \{e_1,e_2,e_3,e_4,e_8\}$; | $c_1' = \{e_1,e_2,e_3,e_4,e_8\}$; |
| $c_2 =$ | $c_2' =$ |
| $c_3 = \{e_4,e_5,e_6,e_8,e_{12}\}$; | $c_3' = \{e_4,e_5,e_6,e_8,e_{12}\}$; |
| $c_4 = \{e_4,e_5,e_6,e_7,e_9\}$; | $c_4' = \{e_4,e_5,e_6,e_7,e_9\}$; |
| $c_5 = \{e_2,e_3,e_4,e_9,e_{13}\}$; | $c_5' = \{e_2,e_3,e_4,e_9,e_{13}\}$; |
| $c_6 =$ | $c_6' =$ |

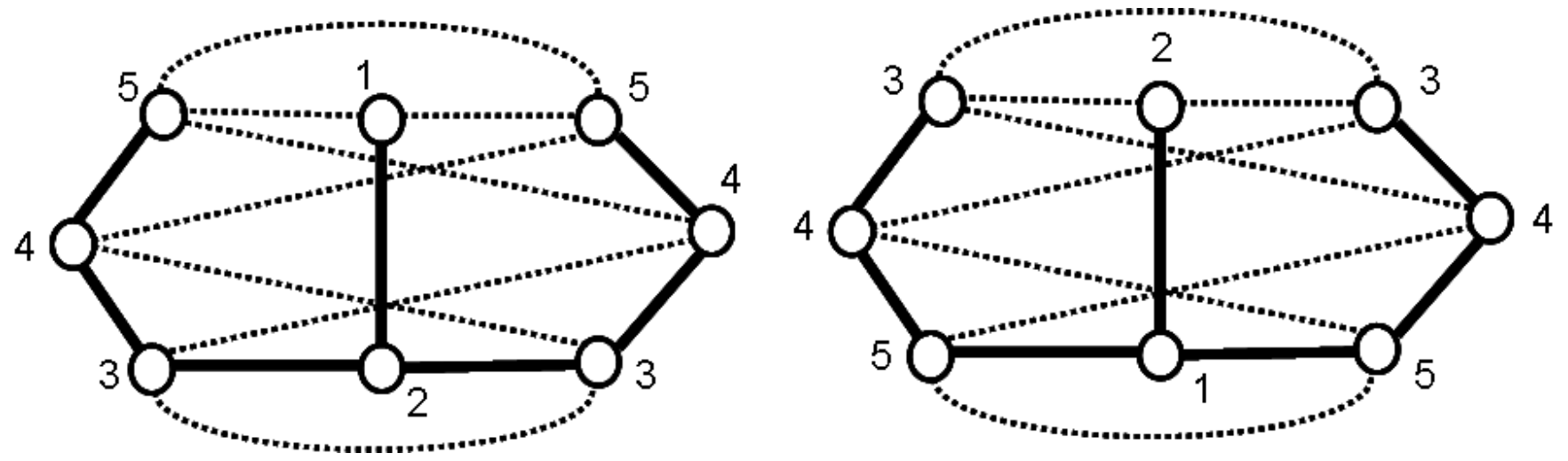

Рис. 3.15. Изометрические циклы относительно 4-го ребра дерева.

Изометрические циклы относительно 4-го ребра: $\{e_1,e_2,e_3,e_4,e_8\}$, $\{e_3,e_4,e_8,e_{12},e_{15}\}$, $\{e_4,e_5,e_6,e_8,e_{12}\}$, $\{e_4,e_5,e_6,e_7,e_9\}$, $\{e_2,e_3,e_4,e_9,e_{13}\}$, $\{e_4,e_5,e_9,e_{13},e_{14}\}$.

Как видно из данного примера, если построение производится только относительно хорд для выбранного дерева графа $G_7$, то множество изометрических циклов будет не полно. В данном примере во множество изометрических циклов не вошли изометрические циклы проходящие по четвертому ребру.

## 3.6. Свойства изометрических циклов и центральных разрезов

Рассмотрим более подробно свойства изометрических циклов и центральных разрезов на примере графа $K_5$.

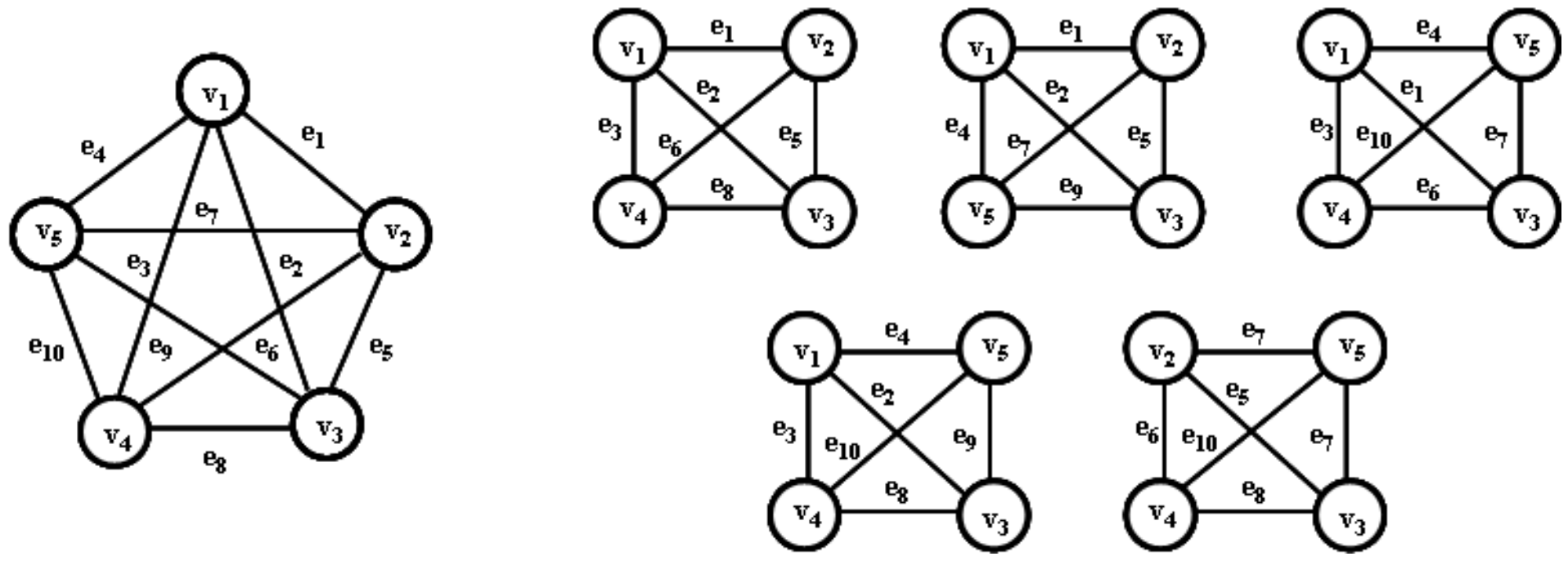

Рис. 3.16. Граф $K_5$ и его пять четырехвершинных подграфов.

Множество центральных разрезов **S** для графа $K_5$:

$s_1 = \{e_1,e_2,e_3,e_4\}$; $s_2 = \{e_1,e_5,e_6,e_7\}$; $s_3 = \{e_2,e_5,e_8,e_9\}$; $s_4 = \{e_3,e_6,e_8,e_{10}\}$;

$s_5 = \{e_4,e_7,e_9,e_{10}\}$.

Множество изометрических циклов C для графа $K_5$:

$c_1 = \{e_1,e_2,e_5\}$; $c_2 = \{e_1,e_3,e_6\}$; $c_3 = \{e_1,e_4,e_7\}$; $c_4 = \{e_2,e_3,e_8\}$; $c_5 = \{e_2,e_4,e_9\}$;

$c_6 = \{e_3,e_4,e_{10}\}$; $c_7 = \{e_5,e_6,e_8\}$; $c_8 = \{e_5,e_7,e_9\}$; $c_9 = \{e_6,e_7,e_{10}\}$; $c_{10} = \{e_8,e_9,e_{10}\}$.

Как видно, все изометрические циклы принадлежат четырехвершинным подграфам (рис. 3.16).

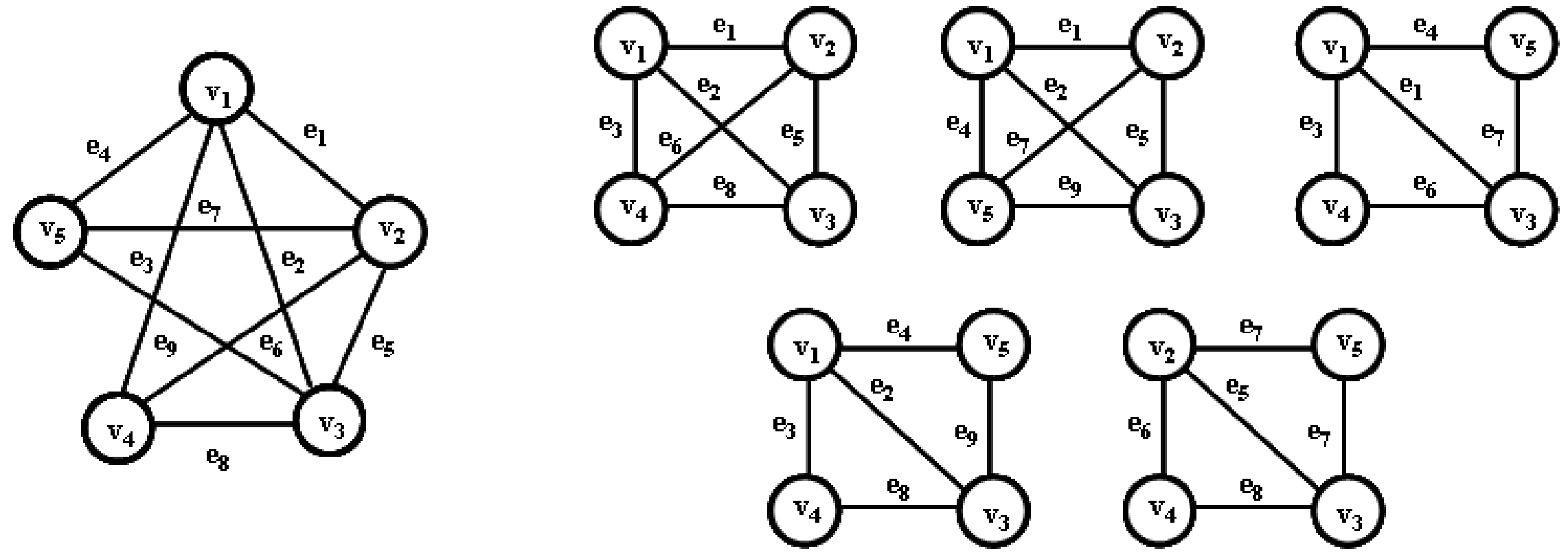

Рис. 3.17. Граф $K_5$ с удаленным ребром $e_{10}$ и его пять четырехвершинных подграфов.

Любой граф G можно получить из полного графа путем удаления соответствующих ребер.

Проиллюстрируем данный процесс на примере графа $K_5$. Удалим из графа ребро $e_{10}$ (рис. 3.17).

Из множества изометрических циклов C удаляются все циклы включающие 10-е ребро:

$c_6 = \{e_3,e_4,e_{10}\}$; $c_9 = \{e_6,e_7,e_{10}\}$; $c_{10} = \{e_8,e_9,e_{10}\}$.

Остаются изометрические циклы:

$c_1 = \{e_1,e_2,e_5\}$; $c_2 = \{e_1,e_3,e_6\}$; $c_3 = \{e_1,e_4,e_7\}$; $c_4 = \{e_2,e_3,e_8\}$; $c_5 = \{e_2,e_4,e_9\}$;
$c_7 = \{e_5,e_6,e_8\}$; $c_8 = \{e_5,e_7,e_9\}$.

В перспективе должны образоваться новые изометрические циклы длиной четыре:

$c_6 \oplus c_9 = \{e_3,e_4,e_{10}\} \oplus \{e_6,e_7,e_{10}\} = \{e_3,e_4,e_6,e_7\}$;
$c_6 \oplus c_{10} = \{e_3,e_4,e_{10}\} \oplus \{e_8,e_9,e_{10}\} = \{e_3,e_4,e_8,e_9\}$;
$c_9 \oplus c_{10} = \{e_6,e_7,e_{10}\} \oplus \{e_8,e_9,e_{10}\} = \{e_6,e_7,e_8,e_9\}$.

Однако их включение во множество оставшихся изометрических циклов невозможно, потому что они образованы как результат попарного кольцевого суммирования оставшихся изометрических циклов (рис. 3.17):

$c_2 \oplus c_3 = \{e_1,e_3,e_6\} \oplus \{e_1,e_4,e_7\} = \{e_3,e_4,e_6,e_7\}$;
$c_4 \oplus c_5 = \{e_2,e_3,e_8\} \oplus \{e_2,e_4,e_9\} = \{e_3,e_4,e_8,e_9\}$;
$c_7 \oplus c_8 = \{e_5,e_6,e_8\} \oplus \{e_5,e_7,e_9\} = \{e_6,e_7,e_8,e_9\}$.

Таким образом, множество изометрических циклов $C_1$ для графа $G_1$ (рис. 3.17), полученного путем удаления 10-го ребра из графа $K_5$, состоит из следующих изометрических циклов:

$c_1 = \{e_1,e_2,e_5\}$; $c_2 = \{e_1,e_3,e_6\}$; $c_3 = \{e_1,e_4,e_7\}$; $c_4 = \{e_2,e_3,e_8\}$; $c_5 = \{e_2,e_4,e_9\}$;
$c_7 = \{e_5,e_6,e_8\}$; $c_8 = \{e_5,e_7,e_9\}$.

Продолжаем удалять ребра из графа $K_5$. Удалим ребро $e_2$ (рис. 3.18).

Из множества изометрических циклов $C_1 = \{c_1,c_2,c_3,c_4,c_5,c_7,c_8\}$ удаляются все циклы включающие 2-е ребро:

$c_1 = \{e_1,e_2,e_5\}$; $c_4 = \{e_2,e_3,e_8\}$; $c_5 = \{e_2,e_4,e_9\}$.

Остаются изометрические циклы:

$c_2 = \{e_1,e_3,e_6\}$; $c_3 = \{e_1,e_4,e_7\}$; $c_7 = \{e_5,e_6,e_8\}$; $c_8 = \{e_5,e_7,e_9\}$.

В перспективе должны образоваться новые изометрические циклы длиной четыре:

$c_1 \oplus c_4 = \{e_1,e_2,e_5\} \oplus \{e_2,e_3,e_8\} = \{e_1,e_3,e_5,e_8\}$;
$c_1 \oplus c_5 = \{e_1,e_2,e_5\} \oplus \{e_2,e_4,e_9\} = \{e_1,e_4,e_5,e_9\}$;
$c_4 \oplus c_5 = \{e_2,e_3,e_8\} \oplus \{e_2,e_4,e_9\} = \{e_3,e_4,e_8,e_9\}$.

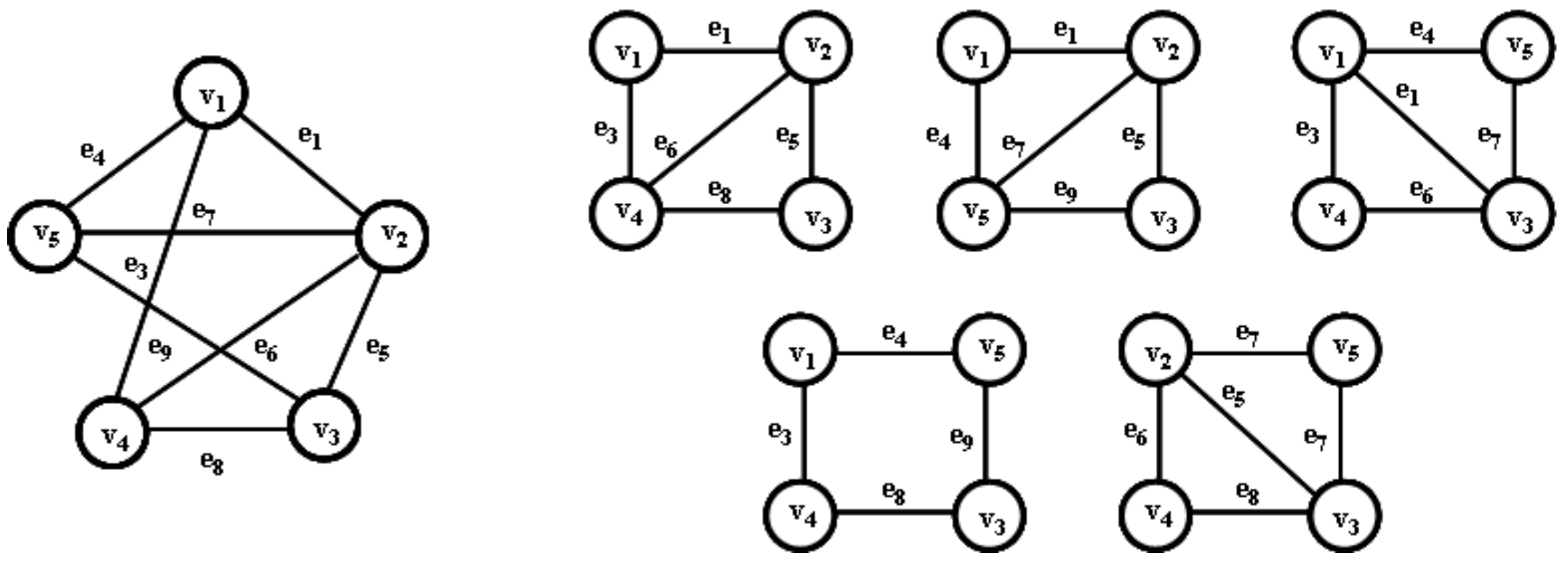

Рис. 3.18. Граф $K_5$ с удаленными ребрами $e_{10}$ и $e_2$ и его пять четырехвершинных подграфов.

Однако включение во множество оставшихся изометрических циклов двух первых невозможно, потому что они образованы как результат попарного кольцевого суммирования оставшихся изометрических циклов (рис. 3.18):

$c_1 \oplus c_4 = \{e_1,e_3,e_6\} \oplus \{e_5,e_6,e_8\} = \{e_1,e_3,e_5,e_8\}$;
$c_1 \oplus c_5 = \{e_1,e_4,e_7\} \oplus \{e_5,e_7,e_9\} = \{e_1,e_4,e_5,e_9\}$.

В то же время, вновь образованный изометрический цикл $c_{4,5}$ может быть включен во множество оставшихся изометрических циклов $c_4 \oplus c_5 = \{e_3,e_4,e_8,e_9\}$, так как для него не существует попарного кольцевого суммирования оставшихся изометрических циклов.

Таким образом, множество изометрических циклов $C_2$ для графа $G_2$ (рис. 3.19), полученного путем удаления 10-го и 2-го ребер из графа $K_5$, состоит из следующих изометрических циклов:

$c_2 = \{e_1,e_3,e_6\}$; $c_3 = \{e_1,e_4,e_7\}$; $c_7 = \{e_5,e_6,e_8\}$; $c_8 = \{e_5,e_7,e_9\}$; $c_{4,5} = \{e_3,e_4,e_8,e_9\}$.

Продолжаем удалять ребра из графа $K_5$. Удалим ребро $e_7$ (рис. 3.19).

Из множества изометрических циклов $C_2 = \{c_2,c_3,c_7,c_8\ c_{4,5}\}$ удаляются все циклы включающие 7-ое ребро:

$c_3 = \{e_1,e_4,e_7\}$; $c_8 = \{e_5,e_7,e_9\}$.

Остаются изометрические циклы:

$c_2 = \{e_1,e_3,e_6\}$; $c_7 = \{e_5,e_6,e_8\}$, $c_{4,5} = \{e_3,e_4,e_8,e_9\}$.

И должны, в перспективе, образоваться новые изометрические циклы:

$c_3 \oplus c_8 = \{e_1,e_4,e_7\} \oplus \{e_5,e_7,e_9\} = \{e_1,e_4,e_5,e_9\}$.

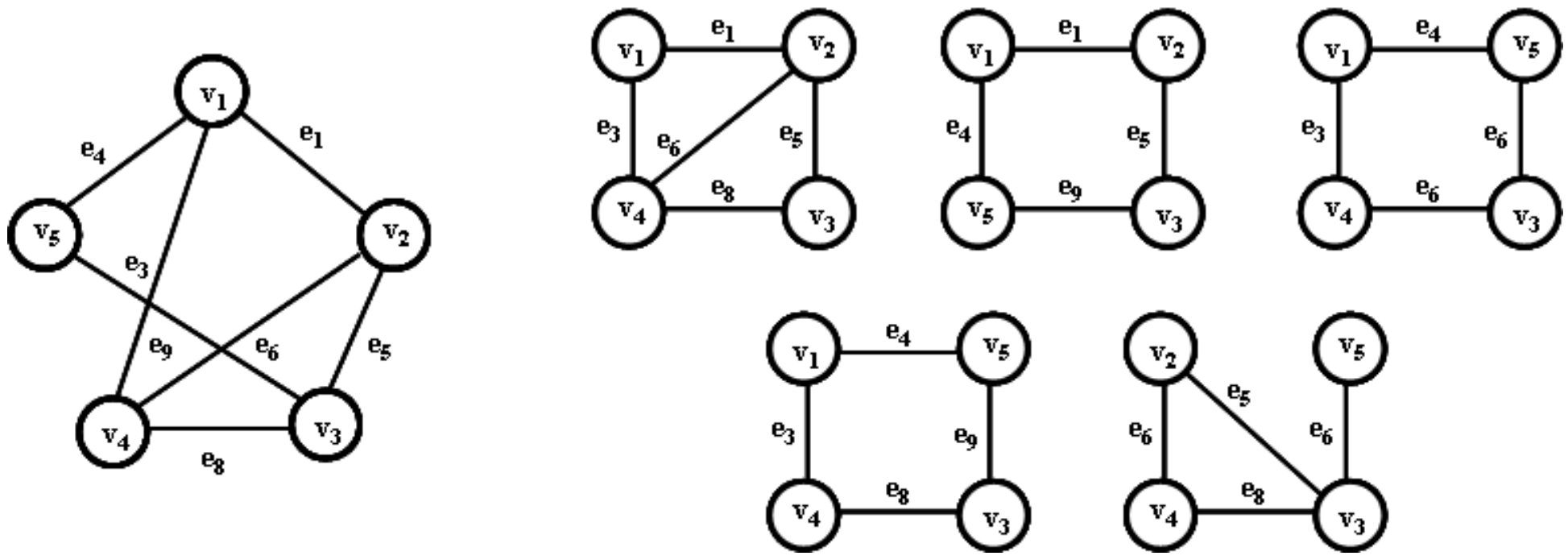

Рис. 3.19. Граф $K_5$ с удаленными $e_{10},e_2,e_7$ ребрами и его пять четырехвершинных подграфов.

Данный цикл $c_{3,8}$ может быть включен во множество оставшихся изометрических циклов, потому что он не может быть образован как результат попарного кольцевого суммирования из оставшихся изометрических циклов (рис. 3.19):

Таким образом, множество изометрических циклов $C_3$ для графа $G_3$, полученного путем удаления 10-го, 2-го и 7-го ребер из графа $K_5$, состоит из следующих изометрических циклов:

$c_2 = \{e_1,e_3,e_6\}$; $c_7 = \{e_5,e_6,e_8\}$, $c_{4,5} = \{e_3,e_4,e_8,e_9\}$; $c_{3,8} = \{e_1,e_4,e_5,e_9\}$.

Продолжаем удалять ребра из графа $K_5$. Удалим ребро $e_1$ (рис. 3.20).

Из множества изометрических циклов $C_3 = \{c_2,c_7,c_{4,5},c_{3,8}\}$ удаляются все циклы, включающие 1-е ребро:

$c_2 = \{e_1,e_3,e_6\}$; $c_{3,8} = \{e_1,e_4,e_5,e_9\}$.

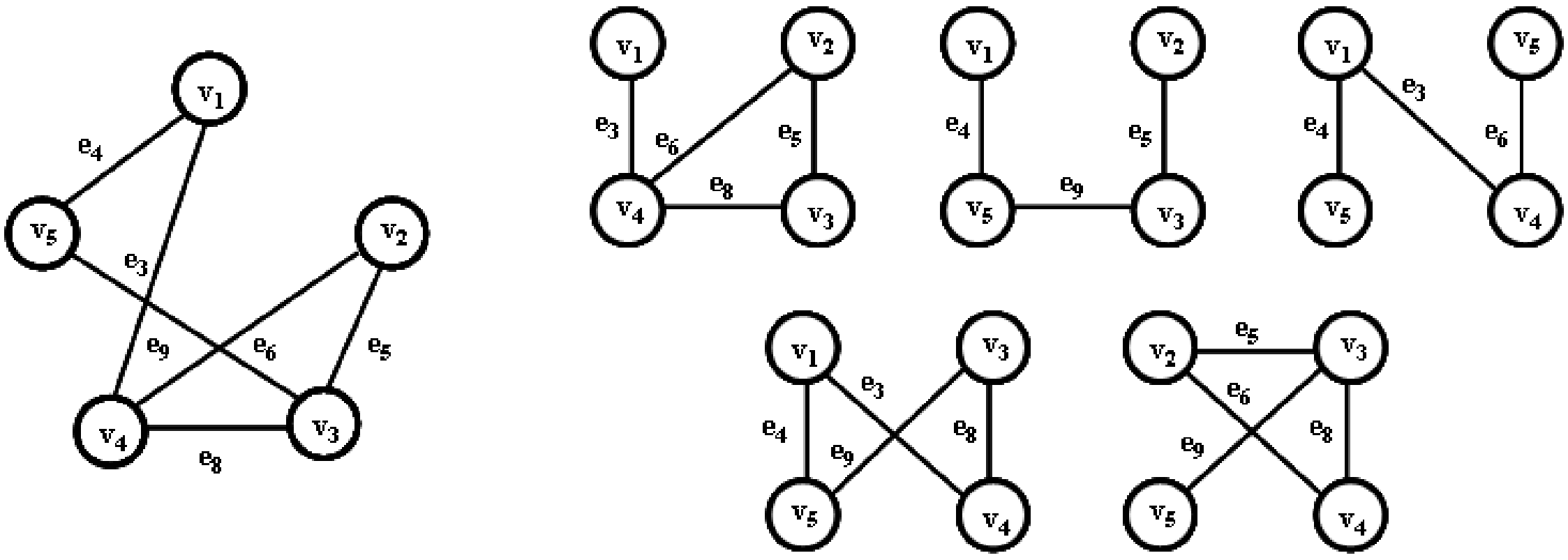

Рис. 3.20. Граф $K_5$ с удаленными ребрами $e_{10},e_2,e_7,e_1$ и его пять четырехвершинных подграфов.

Остаются изометрические циклы:

$c_7 = \{e_5,e_6,e_8\}$, $c_{4,5} = \{e_3,e_4,e_8,e_9\}$.

В перспективе должны образоваться новые изометрические циклы:

$c_2 \oplus c_{3,8} = \{e_1,e_3,e_6\} \oplus \{e_1,e_4,e_5,e_9\} = \{e_3,e_4,e_5,e_6,e_9\}$.

Дальнейшее образование путем удаления ребер можно приостановить, так как полученный цикл равен кольцевой сумме оставшихся циклов:

$c_7 \oplus c_{4,5} = \{e_5,e_6,e_8\} \oplus \{e_3,e_4,e_8,e_9\} = \{e_3,e_4,e_5,e_6,e_9\}$.

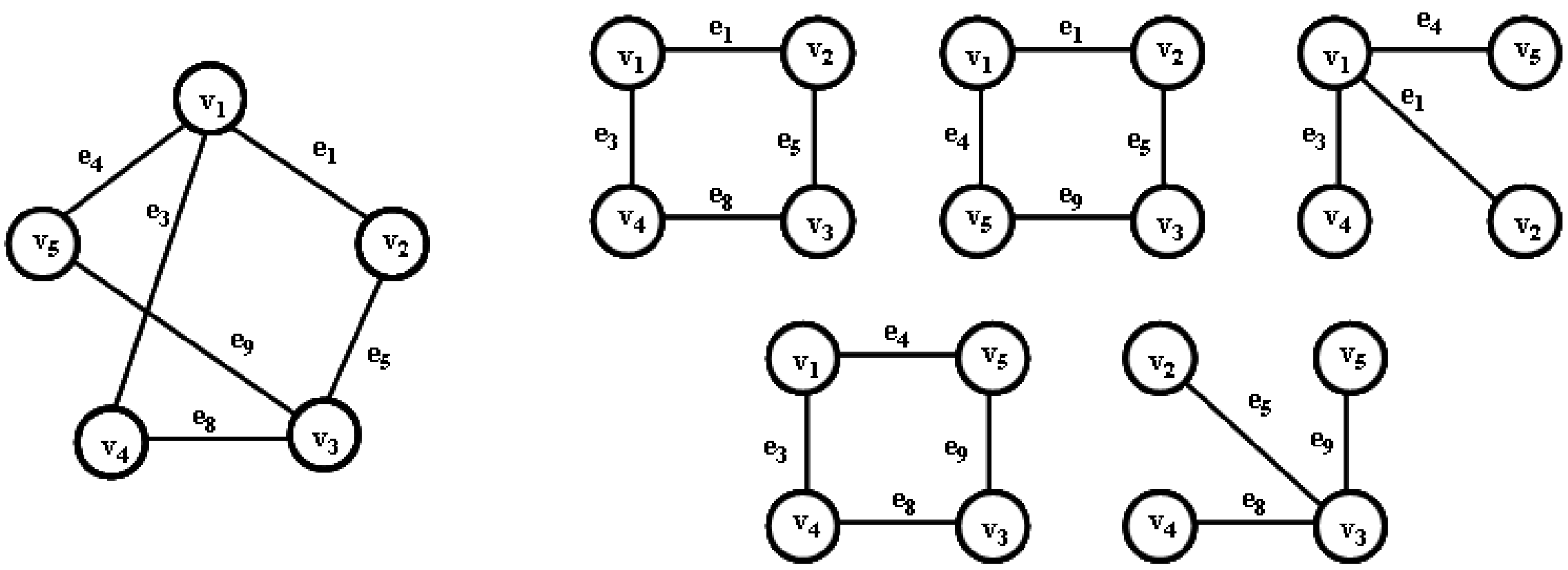

Рис. 3.21. Граф $K_5$ с удаленными ребрами $e_{10},e_2,e_7,e_6$ и его пять четырехвершинных подграфов.

Таким образом, множество изометрических циклов для графа $G_5$ (рис. 3.20), полученного путем удаления 10-го, 2-го, 7-го и 1-го ребер из графа $K_5$ состоит из следующих изометрических циклов: $c_7 = \{e_5,e_6,e_8\}$ и $c_{4,5} = \{e_3,e_4,e_8,e_9\}$.

Если удалить 6-ое ребро вместо 1-го ребра (рис. 3.21), то из множества изометрических циклов $C_3 = \{c_2,c_7, c_{4,5},c_{3,8}\}$ удаляются все циклы включающие 6-е ребро:

$c_2 = \{e_1,e_3,e_6\}$; $c_7 = \{e_5,e_6,e_8\}$.

Остаются изометрические циклы:

$c_{4,5} = \{e_3,e_4,e_8,e_9\}$; $c_{3,8} = \{e_1,e_4,e_5,e_9\}$.

И должны образоваться новые изометрические циклы:

$c_2 \oplus c_7 = \{e_1,e_4,e_7\} \oplus \{e_5,e_7,e_9\} = \{e_1,e_3,e_5,e_8\}$.

Дальнейшее образование путем удаления ребер можно приостановить, так как полученный цикл $c_2,c_7$ равен кольцевой сумме оставшихся циклов:

$c_{4,5} \oplus c_{3,8} = \{e_3,e_4,e_8,e_9\} \oplus \{e_1,e_4,e_5,e_9\} = \{e_1,e_3,e_5,e_8\}$.

Таким образом, множество изометрических циклов для графа $G_5$, полученного путем удаления 10-го, 2-го, 7-го и 6-го ребер из графа $K_5$, состоит из следующих изометрических циклов:

$c_{4,5} = \{e_3,e_4,e_8,e_9\}$ и $c_{3,8} = \{e_1,e_4,e_5,e_9\}$.

Аналогичные рассуждения можно провести для любого графа G с *n* вершинами, удаляя соответствующие ребра из полного графа $K_n$.

Опишем алгоритм формирования множества изометрических циклов из полного графа методом удаления ребер

**[Инициализация].** Существует множество изометрических циклов полного графа $K_n$. Изометрические циклы записанны в виде подмножества ребер $M_e$ и в виде подмножества вершин $M_v$, хранящихся в виде связного списка. Имеется массив номеров исключаемых ребер из полного графа.

**Шаг 1.** [**Перебор исключаемых ребер**]. Последовательно выбираем текущее ребро для исключения. Если список ребер не исчерпан, идем на шаг 2. Иначе конец работы алгоритма.

**Шаг 2.** [**Выбор исключенного ребра графа**]. Выбираем очередное исключаемое ребро графа. Формируем множество новых простых циклов как кольцевую сумму всех попарно пересекающихся по данному ребру изометрических циклов в записи по ребрам. Одновременно формируем такое же связное множество циклов в записи по вершинам как объединение выбранных циклов по вершинам. Идем на шаг 3.

**Шаг 3.** [**Удаление изометрических циклов**]. Связно удаляем изометрические циклы, имеющие исключаемое ребро. Идем на шаг 4.

**Шаг 4.** [**Проверка на включение**]. Проверяем множество новых простых циклов на включение с множеством изометрических циклов по записям в виде вершин. Если включение имеется, то такие простые циклы не включаются во множество изометрических циклов. Если включения нет, то производится запись в массив изометрических циклов. Идем на шаг 1.

Сказанное рассмотрим на примере графа $K_5$ (рис. 3.16).

Пусть задано множество ребер $\{e_2,e_6,e_7\}$ для удаления из полного графа.

Связное множество изометрических циклов запишем в виде связного списка:

| **Цикл** | **Множество $M_e$** | | **Множество $M_v$** |
|---|---|---|---|
| $c_1$ | $\{e_1,e_2,e_5\}$ | $\to$ | $\{v_1,v_2,v_3\}$ |
| $c_2$ | $\{e_1,e_3,e_6\}$ | $\to$ | $\{v_1,v_2,v_4\}$ |
| $c_3$ | $\{e_1,e_4,e_7\}$ | $\to$ | $\{v_1,v_2,v_5\}$ |
| $c_4$ | $\{e_2,e_3,e_8\}$ | $\to$ | $\{v_1,v_3,v_4\}$ |
| $c_5$ | $\{e_2,e_4,e_9\}$ | $\to$ | $\{v_1,v_3,v_5\}$ |
| $c_6$ | $\{e_3,e_4,e_{10}\}$ | $\to$ | $\{v_1,v_4,v_5\}$ |
| $c_7$ | $\{e_5,e_6,e_8\}$ | $\to$ | $\{v_2,v_3,v_4\}$ |
| $c_8$ | $\{e_5,e_7,e_9\}$ | $\to$ | $\{v_2,v_3,v_5\}$ |
| $c_9$ | $\{e_6,e_7,e_{10}\}$; | $\to$ | $\{v_2,v_4,v_5\}$ |
| $c_{10}$ | $\{e_8,e_9,e_{10}\}$ | $\to$ | $\{v_3,v_4,v_5\}$ |

Исключаем второе ребро $e_2$. Формируем новое множество простых циклов:

| **Цикл** | | **Новое $M_e$** | | | **Новое $M_v$** |
|---|---|---|---|---|---|
| $c_1 \oplus c_4 =$ | $\{e_1,e_2,e_5\} \oplus \{e_2,e_3,e_8\}$ | $=\{e_1,e_3,e_5,e_8\}$ | $\to$ | $\{v_1,v_2,v_3\} \cup \{v_1,v_3,v_4\}$ | $=\{v_1,v_2,v_3,v_4\}$ |
| $c_1 \oplus c_5 =$ | $\{e_1,e_2,e_5\} \oplus \{e_2,e_4,e_9\}$ | $=\{e_1,e_4,e_5,e_9\}$ | $\to$ | $\{v_1,v_2,v_3\} \cup \{v_1,v_3,v_5\}$ | $=\{v_1,v_2,v_3,v_5\}$ |
| $c_4 \oplus c_5 =$ | $\{e_2,e_3,e_8\} \oplus \{e_2,e_4,e_9\}$ | $=\{e_1,e_3,e_4,e_5\}$ | $\to$ | $\{v_1,v_3,v_4\} \cup \{v_1,v_3,v_5\}$ | $=\{v_1,v_3,v_4,v_5\}$ |

Как видно, после исключения изометрических циклов $c_1,c_4,c_5$, имеющих ребро $e_2$, новые простые циклы включают следующие изометрические циклы: $c_2,c_7,c_9,c_6,c_3,c_8,c_{10}$. Например, новый цикл $(c_1 \oplus c_4) \to \{v_1,v_2,v_3,v_4\}$ включает $\{v_1,v_2,v_4\}$ и $\{v_2,v_3,v_4\}$. Новый цикл $(c_1 \oplus c_5) \to$ $\{v_1,v_2,v_3,v_5\}$ включает $\{v_1,v_2,v_5\}$ и $\{v_2,v_3,v_5\}$. Новый цикл ($c_4 \oplus c_5$) $\rightarrow$ $\{v_1,v_3,v_4,v_5\}$ включает $\{v_1,v_4,v_5\}$ и $\{v_3,v_4,v_5\}$. Поэтому, во множестве новых циклов отсутствуют новые изометрические циклы для включения во множество изометрических циклов.

| **Цикл** | **Множество $M_e$** | | **Множество $M_v$** |
|---|---|---|---|
| $c_2$ | $\{e_1,e_3,e_6\}$ | $\rightarrow$ | $\{v_1,v_2,v_4\}$ |
| $c_3$ | $\{e_1,e_4,e_7\}$ | $\rightarrow$ | $\{v_1,v_2,v_5\}$ |
| $c_6$ | $\{e_3,e_4,e_{10}\}$ | $\rightarrow$ | $\{v_1,v_4,v_5\}$ |
| $c_7$ | $\{e_5,e_6,e_8\}$ | $\rightarrow$ | $\{v_2,v_3,v_4\}$ |
| $c_8$ | $\{e_5,e_7,e_9\}$ | $\rightarrow$ | $\{v_2,v_3,v_5\}$ |
| $c_9$ | $\{e_6,e_7,e_{10}\}$; | $\rightarrow$ | $\{v_2,v_4,v_5\}$ |
| $c_{10}$ | $\{e_8,e_9,e_{10}\}$ | $\rightarrow$ | $\{v_3,v_4,v_5\}$ |

Исключаем ребро $e_6$. Формируем новое множество простых циклов:

| **Цикл** | | **Новое $M_e$** | | | **Новое $M_v$** |
|---|---|---|---|---|---|
| $c_2 \oplus c_7=$ | $\{e_1,e_3,e_6\}\oplus\{e_5,e_6,e_8\}$ | $=\{e_1,e_3,e_5,e_8\}$ | $\rightarrow$ | $\{v_1,v_2,v_4\}\cup\{v_2,v_3,v_4\}$ | $=\{v_1,v_2,v_3,v_4\}$ |
| $c_2 \oplus c_9=$ | $\{e_1,e_3,e_6\}\oplus\{e_6,e_7,e_{10}\}$ | $=\{e_1,e_3,e_7,e_{10}\}$ | $\rightarrow$ | $\{v_1,v_2,v_4\}\cup\{v_2,v_4,v_5\}$ | $=\{v_1,v_2,v_4,v_5\}$ |
| $c_7 \oplus c_9=$ | $\{e_5,e_6,e_8\}\oplus\{e_6,e_7,e_{10}\}$ | $=\{e_5,e_7,e_8,e_{10}\}$ | $\rightarrow$ | $\{v_2,v_3,v_4\}\cup\{v_2,v_4,v_5\}$ | $=\{v_2,v_3,v_4,v_5\}$ |

После исключения изометрических циклов $c_2,c_7,c_7$, имеющих ребро $e_6$, некоторые новые простые циклы включают в себя оставшиеся изометрические циклы, а некоторые – нет.

Например, новый цикл ($c_2 \oplus c_7$) $\rightarrow$ $\{v_1,v_2,v_3,v_4\}$ не включает в себя оставшиеся изометрические циклы и поэтому может быть включен во множество изометрических циклов. Новый цикл ($c_2 \oplus c_9$) $\rightarrow$ $\{v_1,v_2,v_4,v_5\}$ включает $\{v_1,v_2,v_5\}$ и $\{v_1,v_4,v_5\}$. Новый цикл ($c_7 \oplus c_9$) $\rightarrow$ $\{v_2,v_3,v_4,v_5\}$ включает $\{v_2,v_3,v_5\}$ и $\{v_3,v_4,v_5\}$. Поэтому циклы ($c_2 \oplus c_9$) и ($c_7 \oplus c_9$) не могут быть включены во множество изометрических циклов.

| **Цикл** | **Множество $M_e$** | | **Множество $M_v$** |
|---|---|---|---|
| $c_3$ | $\{e_1,e_4,e_7\}$ | $\rightarrow$ | $\{v_1,v_2,v_5\}$ |
| $c_6$ | $\{e_3,e_4,e_{10}\}$ | $\rightarrow$ | $\{v_1,v_4,v_5\}$ |
| $c_8$ | $\{e_5,e_7,e_9\}$ | $\rightarrow$ | $\{v_2,v_3,v_5\}$ |
| $c_{10}$ | $\{e_8,e_9,e_{10}\}$ | $\rightarrow$ | $\{v_3,v_4,v_5\}$ |

Множество изометрических циклов изменится:

| **Цикл** | **Множество $M_e$** | | **Множество $M_v$** |
|---|---|---|---|
| $c_3$ | $\{e_1,e_4,e_7\}$ | $\rightarrow$ | $\{v_1,v_2,v_5\}$ |
| $c_6$ | $\{e_3,e_4,e_{10}\}$ | $\rightarrow$ | $\{v_1,v_4,v_5\}$ |
| $c_8$ | $\{e_5,e_7,e_9\}$ | $\rightarrow$ | $\{v_2,v_3,v_5\}$ |
| $c_{10}$ | $\{e_8,e_9,e_{10}\}$ | $\rightarrow$ | $\{v_3,v_4,v_5\}$ |
| $c_2 \oplus c_7$ | $\{e_1,e_3,e_5,e_8\}$ | $\rightarrow$ | $\{v_1,v_2,v_3,v_4\}$ |

Исключаем ребро $e_7$. Формируем новое множество простых циклов:

| **Цикл** | | **Новое $M_e$** | | | **Новое $M_v$** |
|---|---|---|---|---|---|
| $c_3 \oplus c_8$ | $\{e_1,e_4,e_7\} \oplus \{e_5,e_7,e_9\}$ | $=\{e_1,e_4,e_5,e_9\}$ | $\rightarrow$ | $\{v_1,v_2,v_5\} \cup \{v_2,v_3,v_5\}$ | $=\{v_1,v_2,v_3,v_5\}$ |

После исключения изометрических циклов $c_3,c_8$, имеющих ребро $e_7$, новый простой цикл не включает в себя оставшиеся изометрические циклы. Поэтому, он может быть включен во множество изометрических циклов.

| Цикл | Множество $M_e$ | | Множество $M_v$ |
|---|---|---|---|
| $c_6$ | $\{e_3,e_4,e_{10}\}$ | $\rightarrow$ | $\{v_1,v_4,v_5\}$ |
| $c_{10}$ | $\{e_8,e_9,e_{10}\}$ | $\rightarrow$ | $\{v_3,v_4,v_5\}$ |
| $c_2 \oplus c_7$ | $\{e_1,e_3,e_5,e_8\}$ | $\rightarrow$ | $\{v_1,v_2,v_3,v_4\}$ |

Множество изометрических циклов изменится.

| Цикл | Множество $M_e$ | | Множество $M_v$ |
|---|---|---|---|
| $c_6$ | $\{e_3,e_4,e_{10}\}$ | $\rightarrow$ | $\{v_1,v_4,v_5\}$ |
| $c_{10}$ | $\{e_8,e_9,e_{10}\}$ | $\rightarrow$ | $\{v_3,v_4,v_5\}$ |
| $c_2 \oplus c_7$ | $\{e_1,e_3,e_5,e_8\}$ | $\rightarrow$ | $\{v_1,v_2,v_3,v_4\}$ |
| $c_3 \oplus c_8$ | $\{e_1,e_4,e_5,e_9\}$ | $\rightarrow$ | $\{v_1,v_2,v_3,v_5\}$ |

Таким образом, для графа с удаленными ребрами $e_2,e_6,e_7$ из полного графа, оставшиеся изометрические циклы образуют множество изометрических циклов усеченного графа.

Можно сказать, что формирование множества изометрических циклов можно осуществить двумя способами:

- первый способ основан на выделении и сравнении циклов, проходящих по каждому ребру алгоритмом поиска в ширину с учетом минимальных s,t-маршрутов графа;
- второй способ основан на методе удаления ребер из полного графа с соответствующим удалением изометрических циклов полного графа и, в случае необходимости, включения в оставшееся множество изометрических циклов, которые образуются как кольцевая сумма удаленных изометрических циклов.

Произведя сравнительный анализ, можно сказать следующее:

- количество изометрических циклов в графе является постоянной величиной равной или большей цикломатического числа графа и не зависит от способа их выделения, в то время как множество фундаментальных циклов в точности равно цикломатическому числу графа и зависит от выбора дерева;
- длина центральных разрезов графа определяет локальные степени вершин.

Рассмотрим основные свойства множества изометрических циклов графа. Введем фундаментальное понятие 0-подмножества изометрических циклов и рассмотрим их основные свойства.

По поводу изометрических циклов следует сказать следующее. Так как каждый суграф графа G представляет собой вектор из пространства суграфов $L_G$ размерностью $m$, то система точек $x_0,x_1,x_2,\ldots,x_k$ $m$-мерного линейного пространства $L_G$ называется независимой, если система векторов:

$$(x_1 - x_0), (x_2 - x_0),\ldots, (x_k - x_0) \tag{3.11}$$

линейно независима. Очевидно, что независимость возможна и при $k < m$. Система (3.11) линейно независима тогда и только тогда, когда из соотношений

$$\lambda_0 x_0 + \lambda_1 x_1 + \lambda_2 x_2 + \ldots + \lambda_k x_k = 0 \tag{3.12}$$

$$\lambda_0 + \lambda_1 + \lambda_2 + \ldots + \lambda_k = 0 \tag{3.13}$$

следует:

$$\lambda_0 = \lambda_1 = \lambda_2 = \ldots = \lambda_k = 0 \tag{3.14}$$

Здесь $\lambda_0,\lambda_1,\lambda_2,\ldots\lambda_k$ – действительные числа. Таким образом, свойство системы $x_0,x_1,x_2,\ldots,x_k$ быть независимой не зависит от порядка нумерации точек. Более того, ясно, что если система точек независима, то всякая её подсистема также независима.

Покажем, что если система векторов (3.11) линейно независима, то из соотношений (3.12) и (3.10) вытекает (3.14). В силу (3.13) соотношение (3.12) переписывается в форме: $(\lambda_1 + \lambda_2 + \ldots + \lambda_k)\, x_0 + \lambda_1 x_1 + \lambda_2 x_2 + \ldots + \lambda_k x_k = 0$ или иначе: $\lambda_1(x_1 - x_0) + \lambda_2(x_2 - x_0) + \ldots + \lambda_k(x_k - x_0) = 0$.

Но, так как система (3.11) линейно независима, то из последнего вытекает $\lambda_0 = \lambda_1 = \lambda_2 = \ldots = \lambda_k = 0$, а отсюда ввиду (3.13) следует и $\lambda_0 = 0$. Покажем теперь, что если из соотношения (3.12) и (3.13) следует (3.14), то система (3.11) линейно независима.

Пусть:

$$\lambda_1(x_1 - x_0) + \lambda_2(x_2 - x_0) + \ldots + \lambda_k(x_k - x_0) = 0. \tag{3.15}$$

Полагая $\lambda_0 = -(\lambda_1 + \lambda_2 + \ldots + \lambda_k)$, мы можем переписать соотношение (3.16) в виде:

$$\lambda_0 x_0 + \lambda_1 x_1 + \lambda_2 x_2 + \ldots + \lambda_k x_k = 0,$$

причем для чисел $\lambda_0,\lambda_1,\lambda_2,\ldots\lambda_k$ выполнено условие (3.13). Таким образом, в силу предположения имеем $\lambda_0 = \lambda_1 = \lambda_2 = \ldots = \lambda_k = 0$, т.е. из (3.16) вытекает $\lambda_1 = \lambda_2 = \ldots = \lambda_k = 0$, а это и означает линейную независимость системы (3.12).

**Определение 3.5.** Кольцевую сумму всех изометрических циклов будем называть ***ободом графа.***

Введем понятие 0-подмножества графа как зависимой системы изометрических циклов и рассмотрим основные свойства.

**Определение 3.6.** ***0-подмножество изометрических циклов*** – это множество изометрических циклов, кольцевая сумма которых есть пустое множество.

Например, пусть заданы изометрические циклы {a,b,c}, {c,d,e}, {b,g,e}, {a,d,g}. Здесь a,b,c,d,e,g – ребра графа (рис. 3.20). Их кольцевая сумма есть пустое множество.

Кольцевая сумма изометрических циклов полного двудольного графа и обода есть пустое множество. Кольцевая сумма полного графа с четным количеством вершин есть пустое множество. Кольцевая сумма полного графа с нечетным количеством вершин и его обода есть пустое иножество.

**Определение 3.7.** ***Дубль-цикл*** – это простой цикл, который допускает, по крайней мере, два различных нетривиальных разложения в сумму изометрических циклов.

В нашем случае может быть образован следующий дубль-цикл:

$\{a,b,c\} \oplus \{c,d,e\} = \{b,g,e\} \oplus \{a,d,g\} = \{a,b,d,e\}$,

также можно организовать следующие дубль-циклы:

$\{a,b,c\} \oplus \{b,g,e\} = \{c,d,e\} \oplus \{a,d,g\} = \{a,c,e,g\}$,
$\{a,b,c\} \oplus \{c,d,e\} \oplus \{b,g,e\} = \{a,d,g\}$.

По сути, множество дубль-циклов – это подмножество простых циклов. Поэтому, дубль-циклы обладают всеми свойствами простых циклов. Изометрические циклы, простые циклы и дубль-циклы являются частными случаями квазициклов.

***Пример 3.3***. Рассмотрим следующий граф $K_4$ (рис. 3.22). Выделим изометрические циклы в этом графе:

$C^{\tau} = \{\{e_1,e_4,e_5\},\{e_2,e_3,e_5\},\{e_1,e_2,e_6\},\{e_3,e_4,e_6\}\}$.

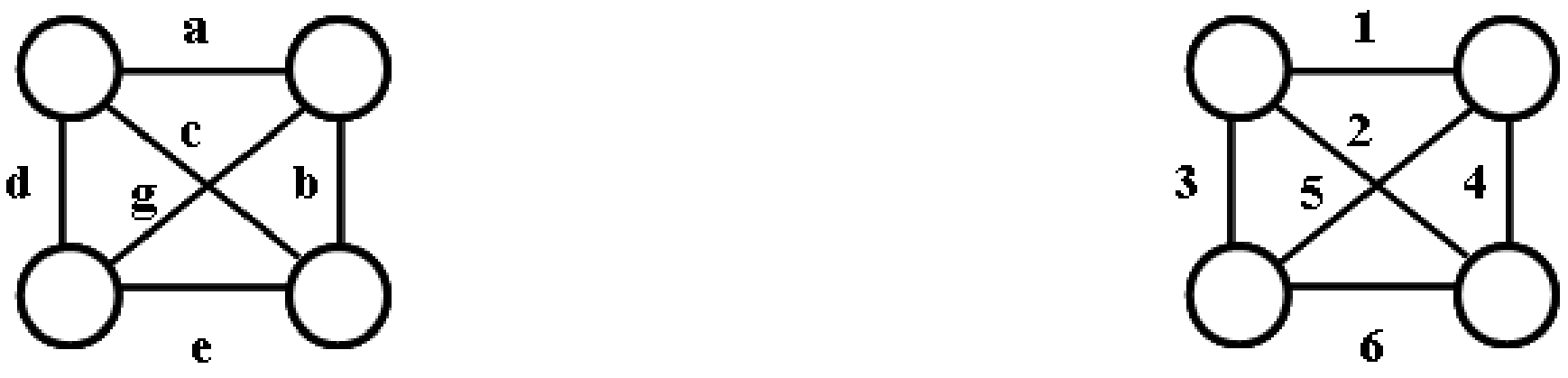

Рис. 3.22. Граф $K_4$.

Кольцевая сумма всех циклов есть пустое множество $\{e_1,e_2,e_4\} \oplus \{e_2,e_3,e_6\} \oplus \{e_1,e_3,e_5\} \oplus \{e_4,e_5,e_6\} = \varnothing$.

Теперь выделим дубль-циклы длины четыре:

$C_1^d = \{e_1,e_2,e_4\} \oplus \{e_2,e_3,e_6\} = \{e_1,e_3,e_5\} \oplus \{e_4,e_5,e_6\} = \{e_1,e_3,e_4,e_6\}$;
$C_2^d = \{e_1,e_2,e_4\} \oplus \{e_1,e_3,e_5\} = \{e_2,e_3,e_6\} \oplus \{e_4,e_5,e_6\} = \{e_2,e_3,e_4,e_5\}$;
$C_3^d = \{e_1,e_2,e_4\} \oplus \{e_4,e_5,e_6\} = \{e_2,e_3,e_6\} \oplus \{e_1,e_3,e_5\} = \{e_1,e_2,e_5,e_6\}$.

## 3.7. Инварианты, построенные на множестве изометрических циклов и центральных разрезов графа

Как мы успели убедиться, запись замкнутого маршрута может быть осуществлена через подмножество ребер или через подмножество вершин графа. Запись цикла через подмножество ребер будем называть реберной записью цикла. Соответственно, вершинная запись цикла состоит из подмножества вершин, принадлежащих рассматриваемому циклу. В основном для записи суграфов и операций над ними применяется реберная запись. Вершинная запись применяется реже и характеризует несколько иные свойства циклов.

Имея множество изометрических циклов графа, можно построить вектор количества изометрических циклов, проходящих по ребру (впредь будем называть его вектором циклов по ребрам). Например, для графа $G_6$ представленного на рис. 3.11 множество изометрических циклов в реберной записи имеет вид:

$c_1 = \{e_1,e_3,e_5\}$; $c_2 = \{e_2,e_3,e_4\}$; $c_3 = \{e_4,e_8,e_9\}$; $c_4 = \{e_5,e_6,e_{12}\}$;
$c_5 = \{e_5,e_8,e_{13}\}$; $c_6 = \{e_6,e_7,e_{11}\}$; $c_7 = \{e_7,e_8,e_{10}\}$; $c_8 = \{e_1,e_2,e_9,e_{13}\}$;
$c_9 = \{e_{10},e_{11},e_{12},e_{13}\}$.

И тогда вектор циклов по ребрам можно записать в виде кортежа

$P_e$ =

| $e_1$ | $e_2$ | $e_3$ | $e_4$ | $e_5$ | $e_6$ | $e_7$ | $e_8$ | $e_9$ | $e_{10}$ | $e_{11}$ | $e_{12}$ | $e_{13}$ |
|---|---|---|---|---|---|---|---|---|---|---|---|---|
| 2 | 2 | 2 | 2 | 3 | 2 | 2 | 3 | 2 | 2 | 2 | 2 | 3 |

или в виде $P_e$ = <2,2,2,2,3,2,2,3,2,2,2,2,3>.

Если записать изометрические циклы через вершины

$c_1 = \{v_1,v_6,v_7\}$; $c_2 = \{v_5,v_6,v_7\}$; $c_3 = \{v_4,v_5,v_7\}$; $c_4 = \{v_1,v_2,v_7\}$;
$c_5 = \{v_1,v_4,v_7\}$; $c_6 = \{v_2,v_3,v_7\}$; $c_7 = \{v_3,v_4,v_7\}$; $c_8 = \{v_1,v_4,v_5,v_6\}$;
$c_9 = \{v_1,v_2,v_3,v_4\}$,

то вектор количества изометрических циклов, проходящих по вершинам графа (впредь будем называть его вектором циклов по вершинам) запишется в виде кортежа:

$P_v$ =

| $v_1$ | $v_2$ | $v_3$ | $v_4$ | $v_5$ | $v_6$ | $v_7$ |
|---|---|---|---|---|---|---|
| 5 | 3 | 3 | 5 | 3 | 3 | 7 |

или в виде $P_v$ = <5,3,3,5,3,3,7>.

Однако существует запись циклов в виде замкнутого ориентированного маршрута, так как любое неориентированное ребро может быть представлено двумя разнонаправленными ориентированными ребрами. Такая запись циклов (замкнутых маршрутов) характерна только для плоских графов с учетом заданного направления обхода. Например, для графа представленного на рис. 3.11 базисная система изометрических циклов и обод, характеризующие плоский граф, могут быть записаны в векторном виде:

$c_1 = (v_1,v_7) + (v_7,v_6) + (v_6,v_1)$;
$c_2 = (v_6,v_7) + (v_7,v_5) + (v_5,v_6)$;
$c_3 = (v_5,v_7) + (v_7,v_4) + (v_4,v_5)$;
$c_4 = (v_1,v_2) + (v_2,v_7) + (v_7,v_1)$;
$c_6 = (v_2,v_3) + (v_3,v_7) + (v_7,v_2)$;
$c_7 = (v_7,v_3) + (v_3,v_4) + (v_4,v_7)$;
$c_9 = (v_1,v_4) + (v_4,v_3) + (v_3,v_2) + (v_2,v_1)$;
$c_0 = c_8 = (v_1,v_6) + (v_6,v_5) + (v_5,v_5) + (v_4,v_1)$.

Инвариант графа – это число (функция) графа G, которое принимает одно и то же значение на любом графе изоморфном G. Пусть f – функция, соотносящая каждому графу G некоторый элемент f(G) из множества M произвольной природы (элементами множества M чаще всего служат числа и системы чисел, векторы, многочлены, матрицы). Эту функцию будем называть инвариантом, если на изоморфных графах её значения совпадают, т.е. для любых G и G':

$$G \cong G' \Rightarrow f(G) = f(G'). \qquad (3.16)$$

Подпространство разрезов и подпространство циклов являются нормированными пространствами, так как любому элементу подпространства можно поставить в соответствие неотрицательное вещественное число $\|l\|$, называемое *нормой*. В данном случае для подпространства разрезов это длина разреза, а для подпространства циклов – длина цикла. Причем, введенное понятие удовлетворяет следующим условиям:

- $\|l\| > 0$ при $l \neq 0$, $\|0\| = 0$,
- $\|l_1 + l_2\| \leq \|l_1\| + \|l_2\|$ для любых $l_1 \in R$, $l_2 \in R$,
- $\|\alpha l\| = |\alpha| \|l\|$ для любого $l \in R$ и вещественного числа $\alpha$.

Множеству центральных разрезов можно поставить в соответствие так называемый вектор локальных степеней вида, который также будет инвариантом графа:

$$P_s = (p_1 \times l_1, p_2 \times l_2, \ldots), \qquad (3.17)$$

где $p_1$ – количество центральных разрезов во множестве **S** длиной $l_1$; $p_2$ – количество центральных разрезов длиной $l_2$ во множестве **S** и т.д. Причем $l_1<l_2<l_3<\ldots$, то есть длина циклов расставлена в порядке неубывания.

Множеству изометрических циклов можно также поставить в соответствие вектор вида (3.16), который также будет инвариантом графа:

$$P_c = (p_1 \times l_1, p_2 \times l_2, \ldots), \qquad (3.18)$$

где $p_1$ – количество изометрических циклов во множестве $C_\tau$ длиной $l_1$; $p_2$ – количество изометрических циклов длиной $l_2$ во множестве $C_\tau$ и т.д. Причем $l_1<l_2<l_3<\ldots$, то есть длина циклов, расставленая в порядке неубывания.

### 3.8. Текст программы Raschet1 определения изометрических циклов как множество ребер

```
program Raschet1;

type
    TMasy = array[1..1000] of integer;
    TMass = array[1..4000] of integer;
var
    F1,F2 : text;
    i,ii,j,jj,K,K1,Np,Nv,Kzikl,M,MakLin : integer;
    Ziklo,Nr,KKK,AB,K9 : integer;
    Masy: TMasy;
    Mass: TMass;
```

```
          Masi: TMass;
          MasyT: TMasy;
          MassT : TMass;
          MasMdop : TMasy;
          MasMy1 : TMasy;
          MasMs1 : TMass;
          MasMy2 : TMasy;
          MasMs2 : TMass;
          MasMy3 : TMasy;
          MasMs3 : TMass;
          MasMcg : TMasy;
          MasMcg1 : TMasy;
          MasKol : TMasy;
          Mass1 : TMass;
{*****************************************************************}
procedure FormVolna(var Nv,Nv1,Nv2 : integer;
                    var My : TMasy;
                    var Ms : TMass;
                    var Mdop : TMasy);
{ Nv -  количество вершин в графе;                                  }
{ Nv1 -  номер первой вершины;                                      }
{ Nv2 -  номер второй вершины;                                      }
{ My -  массив указателей для маттрицы смежностей графа;            }
{ Ms -  массив элементов матрицы смежностей;                        }
{ Mdop -  массив глубины распространения волны.                     }
{                                                                   }
{   Данная процедура формирует массив распространения               }
{   волны, здесь номер уровень волны.                               }
 var I,Im,J,Kum,Istart,Istop : integer;
 label 1,2,3,4;
 begin
       for I:=1 to Nv do Mdop[I]:=0;
       Mdop[Nv1]:=1;
       Mdop[Nv2]:=2;
       Im:=2;
      1:  Im:=Im+1;
       for I:=1 to Nv do
        if Mdop[I]=0 then goto 2;
       goto 3;
      2:  for J:=1 to Nv do
       begin
         if Mdop[J]<>Im-1 then goto 4;
         Istart:=My[J];
         Istop:=My[J+1]-1;
         for I:=Istart to Istop do
         begin
           Kum:=Ms[I];
           if Mdop[Kum]=0 then Mdop[Kum]:=Im;
         end;
      4:  end;
       goto 1;
      3:;
      end; {FormVolna}
{**********************************************************}
 procedure FormKpris(var Kzikl : integer;
                    var Myy : TMasy;
                    var Mss : TMass;
                    var My : TMasy;
                    var Ms : TMass;
                    var Ms2 : TMass);
{ Kzikl  - количество t-циклов в графе;                 }
{ My  - массив указателей для матрицы смежностей;       }
{ Ms  - массив элементов матрицы смежностей;            }
```

```
{ Ms2 -  массив элементов матрицы инциденций;                     }
{ Myy - массив указателей для матрицы t-циклов;                   }
{ Mss - массив элементов матрицы t-циклов;                        }
{                                                                 }
{   Процедура переводит запись циклов в виде вершин               }
{   в запись в виде ребер.                                        }
{                                                                 }
    var I,J,JJ,Ip,Ip1,JJJ,Npn,KK : integer;
    label 1,2;
    begin
        for J:= 1 to Kzikl do
        begin
         for JJ:=Myy[J] to Myy[J+1]-2 do
         begin
           Ip:= Mss[JJ];
           Ip1:=Mss[JJ+1];
           for JJJ:=My[Ip] to My[IP+1]-1 do
            begin
              if Ms[JJJ]<>Ip1 then goto 1;
              Mss[JJ]:=Ms2[JJJ];
1:           end;
        end;
        Mss[Myy[J+1]-1]:=0;
        end;
        Npn:=Myy[Kzikl+1]-1;
        KK:=0;
        for I:=1 to Npn do
        begin
           if Mss[I]=0 then goto 2;
           KK:=KK+1;
           Mss[KK]:=Mss[I];
2:       end;
        for I:=1 to Kzikl do Myy[I+1]:=Myy[I+1]-I;
  end; {FormKpris}
{********************************************************}
  procedure FormDiz(var Mm1,Mm2,Mm4 : integer;
                    var M1 : TMasy;
                    var M2 : TMasy;
                    var M4 : TMasy);
 { Mm1 - количество элементов в первом цикле;                     }
 { Mm2 - количество элементов во втором цикле;                    }
 { Mm4 - количество элементов в их пересечении;                   }
 { M1 - массив элементов первого цикла;                           }
 { M2 - массив элементов второго цикла;                           }
 { M4 - массив элементов пересечения.                             }
 {                                                                }
 {   Процедура определения пересечения двух циклов                }
 {                                                                }
  var J,I : integer;
  label 1,2;
        begin
         Mm4:=0;
         for I:=1 to Mm1 do
         begin
            for J:=1 to Mm2 do
            begin
               if M1[I]<>M2[J] then goto 1;
               Mm4:=Mm4+1;
               M4[Mm4]:=M2[J];
               goto 2;
        1:     end;
        2: end;
        end; {FormDiz}
```

```
{**********************************************************}
procedure FormSwigug(var Kzikl1,Kzikl2 : integer;
                      var My1 : TMasy;
                      var Ms1 : TMass;
                      var My2 : TMasy;
                      var Ms2 : TMass;
                      var Mcg : TMasy;
                      var Mcg1 : TMasy;
                      var Kol : Tmasy);
{ Kzikl1 - количество первых t-циклов в графе;                 }
{ Kzikl2 - количество вторых t-циклов в графе;                 }
{ My1 - массив указателей для первых t-циклов в графе;         }
{ Ms1 - массив для первых t-циклов в графе;                    }
{ My2 - массив указателей для вторых t-циклов в графе;         }
{ Ms2 - массив для вторых t-циклов в графе;                    }
{ Mcg -  вспомогательный массив;                               }
{ Mcg1 -  вспомогательный массив;                              }
{ Kol -  вспомогательный массив.                               }
{                                                              }
{    Процедура формирования и сравнения t-циклов               }
{                                                              }
      label 1;
      var I,J,JJ,II,NN1,NN2,NN3 : integer;
      begin
       for I:=1 to Kzikl1 do
       begin
         NN1:=My1[I+1]-My1[I];
         for II:=1 to NN1 do Mcg[II]:=Ms1[My1[I]-1+II];
         for J:=1 to Kzikl2 do
         begin
           NN2:=My2[J+1]-My2[J];
           for JJ:=1 to NN2 do Mcg1[JJ]:=Ms2[My2[J]-1+JJ];
           FormDiz(NN1,NN2,NN3,Mcg,Mcg1,Kol);
           if NN3=NN1 then goto 1;
         end;
         for II:=1 to NN1 do Ms1[My1[I]-1+II]:=0;
      1:  end;
      end; {FormSwigug}
{************************************************************}
procedure FormDozas(var Kzikl,Kzikl1 : integer;
                      var My1 : TMasy;
                      var Ms1 : TMass;
                      var MyT : TMasy;
                      var MsT : TMass;
                      var Mcg : TMasy;
                      var Mcg1 : TMasy;
                      var Kol : Tmasy);
{ Kzikl  - количество t-циклов в графе;                          }
{ Ms1  - массив указателей для матрицы первых t-циклов;          }
{ Ms1  - массив элементов матрицы первых t-циклов;               }
{ MyT - массив указателей для матрицы t-циклов;                  }
{ MsT - массив элементов матрицы t-циклов;                       }
{ Mcg -  вспомогательный массив;                                 }
{ Mcg1 -  вспомогательный массив;                                }
{ Kol -  вспомогательный массив.                                 }
{                                                                }
{    Процедура формирования t-циклов                             }
{                                                                }
label 1,2,3,4;
var I,J,JJ,II,III,Kp,Ks,NN1,NN2,NN3 : integer;
      begin
       if Kzikl=0  then goto 1;
       Kp:=Kzikl;
```

```
        Ks:=MyT[Kzikl+1]-1;
        for I:=1 to Kzikl1 do
        begin
          NN1:=My1[I+1]-My1[I];
          for II:=1 to NN1 do
          begin
            if Ms1[My1[I]-1+II]=0 then goto 2;
            Mcg[II]:=Ms1[My1[I]-1+II];
          end;
          for J:=1 to Kzikl do
          begin
            NN2:=MyT[J+1]-MyT[J];
            for JJ:=1 to NN2 do Mcg1[JJ]:=MsT[MyT[J]-1+JJ];
            FormDiz(NN1,NN2,NN3,Mcg,Mcg1,Kol);
            if NN3=NN1 then goto 2;
          end;
          Kp:=Kp+1;
          MyT[Kp+1]:=MyT[Kp]+NN1;
          for III:=1 to NN1 do
          begin
            Ks:=Ks+1;
            MsT[Ks]:=Mcg[III];
          end;
       2:  end;
        Kzikl:=Kp;
        goto 3;
       1:  Ks:=0;
        MyT[1]:=1;
        for I:=1 to Kzikl1 do
        begin
          NN1:=My1[I+1]-My1[I];
          for II:=1 to NN1 do
          begin
            if Ms1[My1[I]-1+II]=0 then goto 4;
            Ks:=Ks+1;
            MsT[Ks]:=Ms1[My1[I]-1+II];
          end;
          Kzikl:=Kzikl+1;
          MyT[Kzikl+1]:=MyT[Kzikl]+NN1;
       4:  end;
       3:;
       end; {FormDozas}
{**********************************************************}
 procedure FormSoasda(var Nv1,L,I2 : integer;
                      var Key1 : Boolean;
                      var My1 : TMasy;
                      var Ms1 : TMass;
                      var My : TMasy;
                      var Ms : TMass;
                      var Mcg : TMasy;
                      var Mcg1 : TMasy;
                      var Kol : TMasy);
{ Nv1 - номер вершины в графе;                                 }
{ Key1  - признак;                                             }
{ I2 - признак;                                                }
{ L - длина цикла;                                             }
{ My  - массив указателей для матрицы смежностей;              }
{ Ms - массив для элементов мматрицы смежностей;               }
{ My1 - массив указателей для матрицы t-циклов;                }
{ Ms1 - массив элементов матрицы t-циклов;                     }
{ Mcg -  вспомогательный массив;                               }
{ Mcg1 -  вспомогательный массив;                              }
{ Kol -  вспомогательный массив.                               }
```

```
{                                                                   }
{    Процедура формирования цикла заданной длины                    }
{                                                                   }
     label 2,3,14,11,5,6,12;
     var I,J,Ip1,Isu,Ip2 : integer;
     begin
      Key1:=false;
      if I2>0 then goto 2;
      for I:=1 to L do
      begin
        Mcg[I]:=1;
        Kol[I]:=My1[I+1]-My1[I];
      end;
      Mcg[L]:=0;
     2:  Mcg[L]:=Mcg[L]+1;
     14:;
      for I:=1 to L do
      begin
        if Mcg[I]<=Kol[I] then goto 3;
        if Mcg[1]>Kol[1] then goto 11;
        Mcg[I-1]:=Mcg[I-1]+1;
        for J:=I to L do Mcg[J]:=1;
        goto 14;
     3:  end;
      Mcg1[1]:=Nv1;
      Mcg1[L+1]:=Nv1;
      for I:=2 to L do
      begin
        Ip1:=My1[I]+Mcg[I]-1;
        Isu:=Ms1[Ip1];
        Mcg1[I]:=Isu;
        Ip2:=Mcg1[I-1];
        for J:=My[Ip2] to My[Ip2+1]-1 do
        begin
          if Ms[J]<>Isu then goto 6;
          goto 5;
     6:     end;
        Mcg[I]:=Mcg[I]+1;
        goto 14;
     5:  end;
      Key1:=true;
      I2:=1;
      goto 12;
     11: Key1:=false;
     12:;
     end; {FormSoasda}
{************************************************************}
 procedure FormWegin(var Nv,Nv1,Nv2,Kzikl : integer;
                     var My : TMasy;
                     var Ms : TMass;
                     var My1 : TMasy;
                     var Ms1 : TMass;
                     var MyT : TMasy;
                     var MsT : TMass;
                     var Mdop : TMasy;
                     var Mcg : TMasy;
                     var Mcg1 : TMasy;
                     var Kol : Tmasy);
{ Nv - количество вершин в графе;                                   }
{ Nv1 -  номер первой вершины;                                      }
{ Nv2 -  номер второй вершины;                                      }
{ Kzikl  - количество t-циклов в графе;                             }
{ My  - массив указателей для матрицы смежностей;                   }
```

```
{ Ms  - массив элементов матрицы смежностей;                          }
{ My1 - массив для формирования указателей t-циклов;                  }
{ Ms1 -  массив формирования t-циклов;                                }
{ MyT - массив указателей для матрицы t-циклов;                       }
{ MsT - массив элементов матрицы t-циклов;                            }
{ Mdop -  вспомогательный массив;                                     }
{ Mcg -  вспомогательный массив;                                      }
{ Mcg1 -  вспомогательный массив;                                     }
{ Kol -  вспомогательный массив.                                      }
{ Ms2 -  текущая начальная строка элементов;                          }
{                                                                     }
{    Процедура построения единичных циклов                            }
{                                                                     }
label 1,2,3,22,5,7,8,10,11,13,9;
var I,J,JJ,II,Kot,Istart,Istop,JJJ,KKK,JI,Imum : integer;
var Ji1,Kot1,Kot2,Isu,Irr,I2 : integer;
var Key1 : Boolean;
  begin
      Kot:=0;
      MyT[1]:=1;
      Kzikl:=0;
      for I:=1 to Nv do
      begin
        if Mdop[I]<=2 then goto 1;
        for J:=My[I] to My[I+1]-1 do
        begin
          if Ms[J]<>Nv1 then goto 2;
          My1[1]:=1;
          Isu:=Mdop[I];
          My1[2]:=2;
          Ms1[1]:=Nv1;
          My1[3]:=3;
          Ms1[2]:=I;
          Irr:=Isu;
          Imum:=2;
      3:     Irr:=Irr-1;
          if Irr=1 then goto 22;
          Istart:=My1[Isu-Irr+1];
          Istop:=My1[Isu-Irr+2]-1;
          for II:=1 to Nv do
          begin
            if Mdop[II]<>Irr then goto 5;
            for JJ:=My[II] to My[II+1]-1 do
            begin
              for JJJ:=Istart to Istop do
              begin
                if Ms[JJ]<>Ms1[JJJ] then goto 7;
                Imum:=Imum+1;
                Ms1[Imum]:=II;
      7:         end;
            end;
      5:    end;
         Kot1:=Istop+1-Imum;
         if Kot1=0 then goto 9;
         for II:=Istop+1 to Imum-1 do
         begin
           if Ms1[II]=0 then goto 8;
           for JJ:=II+1 to Imum do
            if Ms1[JJ]=Ms1[II] then Ms1[JJ]:=0;
      8:    end;
      9:    KKK:=Istop;
         for II:=Istop+1 to Imum do
         begin
```

```
           if Ms1[II]=0 then goto 10;
           KKK:=KKK+1;
           Ms1[KKK]:=Ms1[II];
      10:    end;
          Imum:=KKK;
          My1[Isu-Irr+3]:=Imum+1;
          Kot2:=Isu-Irr+2;
          goto 3;
      2:   end;
        goto 1;
      22:  I2:=0;
      11:  Kzikl:=Kzikl+1;
        FormSoasda(Nv1,Isu,I2,Key1,My1,Ms1,My,Ms,Mcg,Mcg1,Kol);
        if Key1=false then goto 13;
        MyT[Kzikl+1]:=MyT[Kzikl]+Isu+1;
        for Ji:=1 to Isu+1 do
        begin
          Kot:=Kot+1;
          MsT[Kot]:=Mcg1[JI];
        end;
        goto 11;
      13:  Kzikl:=Kzikl-1;
      1:  end;
   end; {FormWegin}
{************************************************************}
procedure FormIncide(var Nv : integer;
                     var My : TMasy;
                     var Ms : TMass;
                     var Ms3 : TMass);
{ Nv - количество вершин в графе;                              }
{ My  - массив указателей для матрицы смежностей;              }
{ Ms  - массив элементов матрицы смежностей;                   }
{ Ms3 - массив элементов матрицы инциденций.                   }
{   Формируется матрица инциденций графа в массиве             }
{   Ms3.                                                       }
{                                                              }
     var I,J,K,NNN,P,M,L : integer;
     begin
{     инициализация                        }
     NNN:=My[Nv+1]-1;
     K:=0;
     for J:= 1 to NNN do Ms3[J]:= 0;
{   определение номера элемента                 }
     for I:= 1 to Nv do
      for M:= My[I] to My[I+1]-1 do
       if Ms3[M]=0 then
       begin
         P:=Ms[M];
         K:=K+1;
         Ms3[M]:=K;
         for L:=My[P] to My[P+1]-1 do
          if Ms[L]=I then Ms3[L]:=K;
       end;
     end; {FormIncide}
{**********************************************************}
 procedure EinZikle(var Nv,Kzikl,M : integer;
                    var Masy : TMasy;
                    var Mass : TMass;
                    var Mass1 : TMass;
                    var Masi : TMass;
                    var MasyT : TMasy;
                    var MassT : TMass;
                    var MasMy1 : TMasy;
```

```
                    var MasMs1 : TMass;
                    var MasMy2 : TMasy;
                    var MasMs2 : TMass;
                    var MasMy3 : TMasy;
                    var MasMs3 : TMass;
                    var MasMdop : Tmasy;
                    var MasMcg : Tmasy;
                    var MasMcg1 : Tmasy;
                    var MasKol : Tmasy);
{Процедура создания множества единичных циклов графа                  }
{                                                                     }
{ Nv - количество вершин в графе;                                     }
{ Kzikl  - количество единичных циклов в графе;                       }
{ Masy  - массив указателей для матрицы смежностей;                   }
{ Mass  - массив элементов матрицы смежностей;                        }
{ Mass1 - массив для несмежных элементов строки.                      }
{ Masi -  массив элементов матрицы инциденций;                        }
{ MasyT - массив указателей для матрицы единичных циклов;             }
{ MassT - массив элементов матрицы единичных циклов;                  }
{ MasMy1 -  вспомогательный массив;                                   }
{ MasMs1 -  вспомогательный массив;                                   }
{ MasMy2 -  вспомогательный массив;                                   }
{ MasMs2 -  вспомогательный массив;                                   }
{ MasMy3 -  вспомогательный массив;                                   }
{ MasMs3 -  вспомогательный массив;                                   }
{ MasMdop -  вспомогательный массив для хранения уровней;             }
{ MasMcg -  вспомогательный массив;                                   }
{ MasMcg1 -  вспомогательный массив;                                  }
{ MasKol -  вспомогательный массив.                                   }
label 4;
var I,J,JJ,Nv1,Nv2,Pr,Kzikl1,Kzikl2,Ip1: integer;
begin
    Pr:= Masy[Nv+1]-1;
    M:= Pr div 2;
    FormIncide(Nv,Masy,Mass,Masi);
    Kzikl:= 0;
    MasyT[1]:= 1;
    for I:= 1 to M do
    begin {1}
    Ip1:= I;
      for J:= 1 to Nv do
      begin {2}
        for JJ:=Masy[J] to Masy[J+1]-1 do
        begin {3}
          if Masi[JJ] = Ip1 then
          begin {4}
            Nv1:= J; {определение первой вершины ребра}
            Nv2:= Mass[JJ]; {определение второй вершины ребра}
            goto 4; {концевые вершины ребра определены}
          end;   {4}
        end;   {3}
      end;   {2}
4:      FormVolna(Nv,Nv1,Nv2,Masy,Mass,MasMdop); {алгоритм поиска в ширину}
         {для ориентированного ребра (Nv1,Nv2)}
         FormWegin(Nv,Nv1,Nv2,Kzikl1,Masy,Mass,MasMy3,MasMs3,
         MasMy1,MasMs1,MasMdop,MasMcg,MasMcg1,MasKol);
{построение кандидатов в единичные циклы проходящих по ребру (Nv1,Nv2)}
         FormKpris(Kzikl1,MasMy1,MasMs1,Masy,Mass,Masi);
         FormVolna(Nv,Nv2,Nv1,Masy,Mass,MasMdop);
{алгоритм поиска в ширину для ориентированного ребра (Nv2,Nv1)}
         FormWegin(Nv,Nv2,Nv1,Kzikl2,Masy,Mass,MasMy3,MasMs3,
         MasMy2,MasMs2,MasMdop,MasMcg,MasMcg1,MasKol);
{построение кандидатов в единичные циклы проходящих по ребру (Nv2,Nv1)}
```

```
            FormKpris(Kzikl2,MasMy2,MasMs2,Masy,Mass,Masi);
            FormSwigug(Kzikl1,Kzikl2,MasMy1,MasMs1,MasMy2,MasMs2,
            MasMcg,MasMcg1,MasKol);
            FormDozas(Kzikl,Kzikl1,MasMy1,MasMs1,MasyT,MassT,
            MasMcg,MasMcg1,MasKol);
          end;   {1}
end;{EinZikle}
{************************************************************}
procedure  Shell(var N : integer;
                 var A : TMasy);
{       процедура Шелла для упорядочивания элементов                }
{                                                                   }
{     N - количество элементов в массиве;                           }
{     A - сортируемый массив;                                       }
var D,Nd,I,J,L,X : integer;
label 1,2,3,4,5;
begin
  D:=1;
1:D:=2*D;
  if D<=N then goto 1;
2:D:=D-1;
  D:=D div 2;
  if D=0 then goto 5;
  Nd:=N-D;
  for I:=1 to Nd do
  begin
   J:=I;
3:  L:=J+D;
   if A[L]>=A[J] then goto 4;
   X:=A[J];
   A[J]:=A[L];
   A[L]:=X;
   J:=J-D;
   if J>0 then goto 3;
4:end;
  goto 2;
5:end;{Shell}
{************************************************************}
procedure  ProzYpor(var N : integer;
           var Masy: TMasy;
           var Mass: TMass;
           var A : TMasy);
{Расположение элементов массива в порядке возрастания               }
var i,j,K: integer;
begin
   for i:= 1 to N do
   begin
     K:=0;
     for j:= MasY[i] to MasY[i+1]-1 do
     begin
       K:=K+1;
       A[K]:= Mass[j];
     end;
     Shell(K,A);
     for j:= 1 to K do Mass[MasY[i]-1+j]:=A[j];
   end;
end;{ProzYpor}
{******************************************************************}
begin
          assign(F1,'D:\Isomorf\GRF\Петерсен.grf');
          reset(F1);
          readln(F1,Nv);
          for I:=1 to Nv+1 do
```

```
    begin
      if I<>Nv+1 then read(F1,Masy[I]);
      if I=Nv+1 then readln(F1,Masy[I]);
    end;
    Np:=Masy[Nv+1]-1;
    for I:=1 to Np do
    begin
      if I<>Np then read(F1,Mass[I]);
      if I=Np then read(F1,Mass[I]);
    end;
    close (F1);
    { Создаём новый файл и открываем его в режиме "для чтения и записи"}
    Assign(F2,'D:\Isomorf\EZI\Петерсен.ezi');
    Rewrite(F2);
    EinZikle(Nv,Kzikl,M,Masy,Mass,Mass1,Masi,MasyT,MassT,
            MasMy1,MasMs1,MasMy2,MasMs2,MasMy3,MasMs3,MasMdop,
            MasMcg,MasMcg1,MasKol);
    ProzYpor(Kzikl,MasyT,MassT,MasKol);
    writeln(F2,Nv);
    writeln(F2,M);
    writeln(F2,Kzikl);
    for I:=1 to Nv+1 do
    begin
      if i<>Nv+1 then write(F2,Masy[i],' ');
      if i=Nv+1 then writeln(F2,Masy[i]);
    end;
    for I:=1 to Nv do
    for j:=Masy[i] to Masy[i+1]-1 do
    begin
      if j<>Masy[i+1]-1 then write(F2,Mass[j],' ');
      if j=Masy[i+1]-1 then writeln(F2,Mass[j]);
    end;
    for I:=1 to Nv do
    for j:=Masy[i] to Masy[i+1]-1 do
    begin
      if j<>Masy[i+1]-1 then write(F2,Masi[j],' ');
      if j=Masy[i+1]-1 then writeln(F2,Masi[j]);
    end;
    for i:=1 to Kzikl+1 do
    begin
       if i<> Kzikl+1 then write(F2, MasyT [i],' ');
       if i= Kzikl+1 then writeln(F2, MasyT [i]);
    end;
    for I:=1 to Kzikl do
    begin
      for j:=MasyT[i] to MasyT[i+1]-1 do
      begin
        if j<>MasyT[i+1]-1 then write(F2,MassT[j],' ');
        if j=MasyT[i+1]-1 then writeln(F2,MassT[j]);
      end;
    end;
    close (F2);
    writeln('Конец расчета!');
end.
```

## 3.9. Входные и выходные файлы процедуры Raschet1

В качестве примера рассмотрим граф Петерсена (см. рис. 2.23).

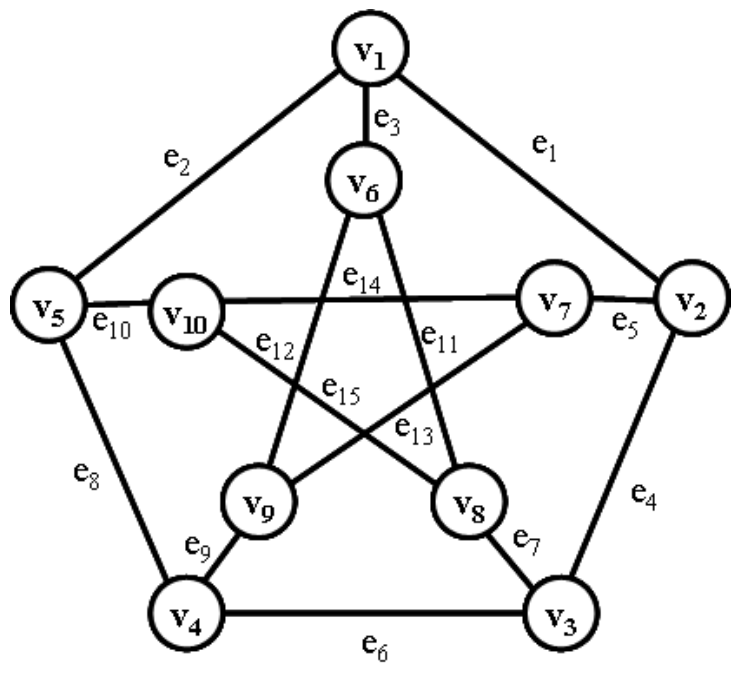

Рис .2.23. Граф Петерсена.

Входной файл Петерсен.grf

| | |
|---|---|
| 10 | {количество вершин} |
| 1 4 7 10 13 16 19 22 25 28 31 | {массив указателей} |
| 2 5 6 | {смежность вершины 1} |
| 1 3 7 | {смежность вершины 2} |
| 2 4 8 | {смежность вершины 3} |
| 3 5 9 | {смежность вершины 4} |
| 1 4 10 | {смежность вершины 5} |
| 1 8 9 | {смежность вершины 6} |
| 2 9 10 | {смежность вершины 7} |
| 3 6 10 | {смежность вершины 8} |
| 4 6 7 | {смежность вершины 9} |
| 5 7 8 | {смежность вершины 10} |

Выходной файл Петерсен.ezi

| | |
|---|---|
| 10 | {количество вершин} |
| 15 | {количество ребер} |
| 12 | {количество изометрических циклов} |
| 1 4 7 10 13 16 19 22 25 28 31 | {массив указателей для матрицы смежностей} |
| 2 5 6 | {смежность вершины 1} |
| 1 3 7 | {смежность вершины 2} |
| 2 4 8 | {смежность вершины 3} |
| 3 5 9 | {смежность вершины 4} |
| 1 4 10 | {смежность вершины 5} |
| 1 8 9 | {смежность вершины 6} |
| 2 9 10 | {смежность вершины 7} |
| 3 6 10 | {смежность вершины 8} |
| 4 6 7 | {смежность вершины 9} |
| 5 7 8 | {смежность вершины 10} |
| 1 2 3 | {совместимость ребер для вершины 1} |
| 1 4 5 | {совместимость ребер для вершины 2} |
| 4 6 7 | {совместимость ребер для вершины 3} |
| 6 8 9 | {совместимость ребер для вершины 4} |
| 2 8 10 | {совместимость ребер для вершины 5} |
| 3 11 12 | {совместимость ребер для вершины 6} |
| 5 13 14 | {совместимость ребер для вершины 7} |
| 7 11 15 | {совместимость ребер для вершины 8} |
| 9 12 13 | {совместимость ребер для вершины 9} |
| 10 14 15 | {совместимость ребер для вершины 10} |
| 1 6 11 16 21 26 31 36 41 46 51 56 61 | {массив указателей для циклов} |
| 1 2 4 6 8 | {изометрический цикл 1 в ребрах} |
| 1 2 5 10 14 | {изометрический цикл 2 в ребрах } |
| 1 3 4 7 11 | {изометрический цикл 3 в ребрах } |
| 1 3 5 12 13 | {изометрический цикл 4 в ребрах } |
| 2 3 10 11 15 | {изометрический цикл 5 в ребрах } |
| 2 3 8 9 12 | {изометрический цикл 6 в ребрах } |
| 4 5 6 9 13 | {изометрический цикл 7 в ребрах } |
| 4 5 7 14 15 | {изометрический цикл 8 в ребрах } |
| 6 7 9 11 12 | {изометрический цикл 9 в ребрах } |
| 6 7 8 10 15 | {изометрический цикл 10 в ребрах } |
| 8 9 10 13 14 | {изометрический цикл 11 в ребрах } |
| 11 12 13 14 15 | {изометрический цикл 12 в ребрах } |

**Примечания**

**примечание 1**. Массив указателей формируется как сложение локальных степеней вершин с указанием местоположения в списке смежных вершин

Например: 1 4 7 10 13 16 19 22 25 28 31 {MY - массив указателей}

Нужно определить смежность для вершины 4.

Элементы для 4-ой вершины начинаются с MY(4)=10 и заканчивается MY(5)-1=13-1=12.

Следовательно, вершины смежные к вершине 4 находятся с 10 по 12 местах в списке вершин

2 5 6 1 3 7 2 4 8 **3 5 9** 1 4 10 1 8 9 2 9 10 3 6 10 4 6 7 5 7 8

**примечание 2.** Цифрами, в машинной информации, обозначаются либо вершины, либо ребра графа в зависимости от их принадлежности.

### 3.10. Текст программы PrintZikNabc для вычисления изометрических циклов как множество вершин

```
type
      TMasy = array[1..1000] of integer;
      TMass = array[1..4000] of integer;
var
      F1,F2 : text;
      i,ii,j,jj,K,K1,M,Nv,Kzikl,Nr : integer;
      Siz,Np,Nk,M1,KCF,KCF1,Ln : integer;
      Masy: TMasy;
      Mass: TMass;
      Masi: TMass;
      MasMy1: TMasy;
      MasMs1: TMass;
      MasMy2: TMasy;
      MasMs2: TMass;
      MasMy5: TMasy;
      MasMs5: TMass;
      MasMy4: TMasy;
```

```
          MasMs4: TMass;
          MasMdop: TMasy;
{*************************************************************}
procedure FormVerBasis1(var Nv,Ziklo : integer;
                        var My : TMasy;
                        var Ms : TMass;
                        var My1 : TMasy;
                        var Ms1 : TMass;
                        var Mass1 : TMass;
                        var Mdop : TMasy;
                        var Masi : TMass);
{Процедура записи базиса циклов через вершины                   }
{                                                               }
{ Nv - количество вершин в графе;                               }
{ Ziklo - цикломатическое число графа;                          }
{ My - массив указателей для матрицы смежностей;                }
{ Ms - массив элементов матрицы смежностей;                     }
{ My1 - массив указателей для матрицы базисных циклов           }
{ Ms1 : массив элементов матрицы базисных циклов;               }
{ Masi - массив элементов матрицы инциденций;                   }
{ Mass1 - элементов матрицы базисных циклов через вершины;      }
{ Mdop - вспомогательный массив;                                }
{***************************************************************}
label 1,2,3;
var i,j,ii,jj,iii,jjj,Nr,KKK,K1 :integer;
begin
  for i:= 1 to Ziklo do
  begin{1}
   {writeln(F2,'FormVerBasis: i = ',i);}
    KKK:=0;
    j:=My1[i]-1;
    2:
    j:=j+1;
    if j = My1[i+1] then goto 3 else
    begin {2}
     Nr:= Ms1[j];
     {writeln(F2,'FormVerBasis: Nr = ',Nr);}
     for ii:=1 to Nv do
     begin {3}
      for jj:= My[ii] to My[ii+1]-1 do
      begin {4}
       if Masi[jj]= Nr  then
       begin {5}
        KKK:=KKK+1;
        Mdop[KKK]:=ii;
        KKK:=KKK+1;
        Mdop[KKK]:=Ms[jj];
        {for jjj:= 1 to KKK do
        begin
         if jjj <> KKK then write(F2, Mdop[jjj],' ');
         if jjj = KKK then writeln(F2, Mdop[jjj],' ');
        end;}
        goto 2;
       end; {5}
      end;{4}
     end;{3}
    end; {2}
    3: {writeln(F2,'Запись [FormVerBasis]');}
    for jjj:=1 to KKK-1 do if Mdop[jjj]<>0 then
    begin {7}
     for iii:=jjj+1 to KKK do if Mdop[iii]<>0 then
     begin {8}
      if Mdop[iii]=Mdop[jjj] then Mdop[iii]:=0;
```

```
          begin {9}
           K1:=0;
           for jj:=1 to KKK do if Mdop[jj]<>0 then
           begin {6}
             K1:=K1+1;
             Mdop[K1]:=Mdop[jj];
           end; {6}
          end; {9}
         end; {8}
        end; {7}
        for jjj:=1 to K1 do Mass1[My1[i]-1 +jjj]:=Mdop[jjj];
   end; {1}
   {for i:= 1 to Ziklo do
     begin
       for j:= My1[i] to My1[i+1]-1 do
       begin
        if j<> My1[i+1]-1 then write(F2, Mass1[j],' ');
        if j= My1[i+1]-1 then writeln(F2, Mass1[j],' ');
       end;
     end;}
end;{FormVerBasis1}
{*******************************************************************}
procedure  Shell(var N : integer;
                 var A : TMasy);
{       процедура Шелла для упорядочивания элементов        }
{                                                    }
{     N - количество элементов в массиве;                  }
{     A - сортируемый массив;                          }
var D,Nd,I,J,L,X : integer;
label 1,2,3,4,5;
begin
  D:=1;
1:D:=2*D;
  if D<=N then goto 1;
2:D:=D-1;
  D:=D div 2;
  if D=0 then goto 5;
  Nd:=N-D;
  for I:=1 to Nd do
  begin
    J:=I;
3:  L:=J+D;
    if A[L]>=A[J] then goto 4;
    X:=A[J];
    A[J]:=A[L];
    A[L]:=X;
    J:=J-D;
    if J>0 then goto 3;
4:end;
  goto 2;
5:end;{Shell}
{************************************************************}
label 1;
begin

          assign(F1,'D:\Isomorf\EZI\Петерсен.ezi');
          reset(F1);
          readln(F1,Nv);
          readln(F1,Nr);
          readln(F1,Kzikl);
          {Вводим матрицу смежностей неориентированного графа }
          for I:=1 to Nv+1 do
          begin {1}
```

```
  if i<>Nv+1 then read(F1,Masy[i]);
  if i=Nv+1 then readln(F1,Masy[i]);
end;  {1}
for I:=1 to Nv do
for j:=Masy[i] to Masy[i+1]-1 do
begin {2}
  if j<>Masy[i+1]-1 then read(F1,Mass[j]);
  if j=Masy[i+1]-1 then readln(F1,Mass[j]);
end;  {2}
{Вводим матрицу инциденций неориентированного графа }
for I:=1 to Nv do
for j:=Masy[i] to Masy[i+1]-1 do
begin
  if j<>Masy[i+1]-1 then read(F1,Masi[j]);
  if j=Masy[i+1]-1 then readln(F1,Masi[j]);
end;
{Вводим единичные циклы }
for I:=1 to Kzikl+1 do
begin
  if i<>Kzikl+1 then read(F1,MasMy4[i]);
  if i=Kzikl+1 then readln(F1,MasMy4[i]);
end;
for I:=1 to Kzikl do
begin
  for j:=MasMy4[i] to MasMy4[i+1]-1 do
  begin
     if j<>MasMy4[i+1]-1 then read(F1,MasMs4[j]);
     if j=MasMy4[i+1]-1 then readln(F1,MasMs4[j]);
  end;
end;
close (F1);
FormVerBasis1(Nv,Kzikl,Masy,Mass,MasMy4,MasMs4,MasMs2,MasMy5,Masi);
{ Создаём новый файл и открываем его в режиме "для чтения и записи"}
Assign(F2,'D:\Isomorf\MY4\Петерсен.my4');
Rewrite(F2);
writeln(F2,'  Количество вершин графа = ',Nv);
writeln(F2,'  Количество ребер графа = ',Nr);
writeln(F2,'  Количество единичных циклов = ',Kzikl);
writeln(F2,' ');
writeln(F2,'  Матрица смежностей графа: ');
writeln(F2,' ');
for I:=1 to Nv do
begin
    write(F2,' вершина ',I:3,':');
    Siz:=0;
    for K:=Masy[I] to Masy[I+1]-1 do
    begin
        Siz:=Siz+1;
        Np:=Siz mod 10;
        Nk:=Siz div 10;
        if K<>Masy[I+1]-1 then
        begin
         if (Np=1) and (Nk>0) then
         write(F2,' ',Mass[K]:3);
         if (Np=1) and (Nk=0) then
         write(F2,' ',Mass[K]:3);
         if (Np>1) and (Np<=9) then
         write(F2,' ',Mass[K]:3);
         if Np=0 then
          writeln(F2,' ',Mass[K]:3);
       end;
       if K=Masy[I+1]-1 then
       begin
```

```
        if (Np=1) and (Nk>0) then
         writeln(F2,' ',Mass[K]:3)
        else
       writeln(F2,' ',Mass[K]:3);
      end;
    end;
end;
writeln(F2,' ');
writeln(F2,'  Элементы матрицы инциденций: ');
writeln(F2,' ');
M1:=Masy[Nv+1]-1;
M:=M1 div 2;
for I:=1 to M do
begin
  for  K:=1 to Nv do
  begin
    for  M1:=Masy[K] to Masy[K+1]-1 do
    begin
      if Masi[M1]=I then
      begin
        write(F2,' ребро ',I:3,':');
        Np:=Mass[M1];
        write(F2,'  ( ',K:3);
        write(F2,' ',Np:3,' )  или  ');
        write(F2,' ( ',Np:3);
        writeln(F2,'  ',K:3,' )');
        goto 1;
      end;
    end;
  end;
1: end;
writeln(F2,' ');
writeln(F2,'  Множество единичных циклов графа: ');
writeln(F2,' ');
for  I:=1 to Kzikl do
begin
    write(F2,' цикл ',I:3,':');
    Siz:=0;
    for K:=MasMy4[I] to MasMy4[I+1]-1 do
    begin
       Siz:=Siz+1;
       Np:=Siz mod 10;
       Nk:=Siz div 10;
       if K<>MasMy4[I+1]-1 then
       begin
        if (Np=1) and (Nk>0) then
          write(F2,' ',MasMs4[K]:3);
        if (Np=1) and (Nk=0) then
         write(F2,' ',MasMs4[K]:3);
        if (Np>1) and (Np<=9) then
         write(F2,' ',MasMs4[K]:3);
        if Np=0 then
         writeln(F2,' ',MasMs4[K]:3);
      end;
      if K=MasMy4[I+1]-1 then
      begin
       if (Np=1) and (Nk>0) then
          writeln(F2,' ',MasMs4[K]:3)
       else
         writeln(F2,' ',MasMs4[K]:3);
      end;
    end;
end;
```

```
  writeln(F2,' ');
  for I:=1 to Kzikl do
  begin
   Ln:= MasMy4[i+1] -MasMy4[I];
   ii:=0;
   for j:=MasMy4[i] to MasMy4[I+1]-1 do
   begin
    ii:=ii+1;
    MasMdop[ii]:=MasMs2[j];
   end;
   Shell(Ln,MasMdop);
   ii:=0;
   for j:=MasMy4[i] to MasMy4[I+1]-1 do
   begin
    ii := ii+1;
    MasMs2[j]:=MasMdop[ii];
   end;
  end;
  writeln(F2,'  Множество единичных циклов графа в записи вершин: ');
  writeln(F2,' ');
  for I:=1 to Kzikl do
  begin
    write(F2,' цикл ',I:3,':');
    Siz:=0;
    for K:=MasMy4[I] to MasMy4[I+1]-1 do
    begin
      Siz:=Siz+1;
      Np:=Siz mod 10;
      Nk:=Siz div 10;
      if K<>MasMy4[I+1]-1 then
      begin
       if (Np=1) and (Nk>0) then
         write(F2,' ',MasMs2[K]:3);
       if (Np=1) and (Nk=0) then
        write(F2,' ',MasMs2[K]:3);
       if (Np>1) and (Np<=9) then
        write(F2,' ',MasMs2[K]:3);
       if Np=0 then
        writeln(F2,' ',MasMs2[K]:3);
     end;
     if K=MasMy4[I+1]-1 then
     begin
      if (Np=1) and (Nk>0) then
         writeln(F2,' ',MasMs2[K]:3)
      else
        writeln(F2,' ',MasMs2[K]:3);
     end;
    end;
  end;
  close (F2);
  writeln('Конец расчета!');
end.
```

## 3.11. Входные и выходные файлы программы PrintZikNabc

Выходной файл Петерсен.my4

Количество вершин графа = 10
Количество ребер графа = 15
Количество единичных циклов = 12

Матрица смежностей графа:

вершина 1: 2 5 6
вершина 2: 1 3 7
вершина 3: 2 4 8
вершина 4: 3 5 9
вершина 5: 1 4 10
вершина 6: 1 8 9
вершина 7: 2 9 10
вершина 8: 3 6 10
вершина 9: 4 6 7
вершина 10: 5 7 8

Элементы матрицы инциденций:

ребро 1: ( 1 2 ) или ( 2 1 )
ребро 2: ( 1 5 ) или ( 5 1 )
ребро 3: ( 1 6 ) или ( 6 1 )
ребро 4: ( 2 3 ) или ( 3 2 )
ребро 5: ( 2 7 ) или ( 7 2 )
ребро 6: ( 3 4 ) или ( 4 3 )
ребро 7: ( 3 8 ) или ( 8 3 )
ребро 8: ( 4 5 ) или ( 5 4 )
ребро 9: ( 4 9 ) или ( 9 4 )
ребро 10: ( 5 10 ) или ( 10 5 )
ребро 11: ( 6 8 ) или ( 8 6 )
ребро 12: ( 6 9 ) или ( 9 6 )
ребро 13: ( 7 9 ) или ( 9 7 )
ребро 14: ( 7 10 ) или ( 10 7 )
ребро 15: ( 8 10 ) или ( 10 8 )

Множество изометрических циклов графа в записи ребер:

цикл 1: 1 2 4 6 8
цикл 2: 1 2 5 10 14
цикл 3: 1 3 4 7 11
цикл 4: 1 3 5 12 13
цикл 5: 2 3 10 11 15
цикл 6: 2 3 8 9 12
цикл 7: 4 5 6 9 13
цикл 8: 4 5 7 14 15
цикл 9: 6 7 9 11 12
цикл 10: 6 7 8 10 15
цикл 11: 8 9 10 13 14
цикл 12: 11 12 13 14 15

Множество изометрических циклов графа в записи вершин:

цикл 1: 1 2 3 4 5
цикл 2: 1 2 5 7 10
цикл 3: 1 2 3 6 8

цикл 4: 1 2 6 7 9
цикл 5: 1 5 6 8 10
цикл 6: 1 4 5 6 9
цикл 7: 2 3 4 7 9
цикл 8: 2 3 7 8 10
цикл 9: 3 4 6 8 9
цикл 10: 3 4 5 8 10
цикл 11: 4 5 7 9 10
цикл 12: 6 7 8 9 10

## Примеры к главе 3

***Пример 3.4.*** Выделить множество изометрических циклов в графе **$G_1$**.

Количество вершин графа = 10.
Количество ребер графа = 23.
Количество изометрических циклов графа = 32.

Матрица смежностей графа, представленная в виде списка смежных вершин:

вершина $v_1$: $\{v_2,v_6,v_7,v_{10}\}$;
вершина $v_2$: $\{v_1,v_3,v_5,v_7\}$;
вершина $v_3$: $\{v_2,v_4,v_5,v_9\}$;
вершина $v_4$: $\{v_3,v_5,v_6,v_7,v_9\}$;
вершина $v_5$: $\{v_2,v_3,v_4,v_6,v_8,v_{10}\}$;
вершина $v_6$: $\{v_1,v_4,v_5,v_7,v_9\}$;
вершина $v_7$: $\{v_1,v_2,v_4,v_6,v_8\}$;
вершина $v_8$: $\{v_5,v_7,v_9,v_{10}\}$;
вершина $v_9$: $\{v_3,v_4,v_6,v_8,v_{10}\}$;
вершина $v_{10}$: $\{v_1,v_5,v_8,v_9\}$.

Элементы матрицы инциденций:

ребро $e_1$: $(v_1,v_2)$ или $(v_2,v_1)$;
ребро $e_2$: $(v_1,v_6)$ или $(v_6,v_1)$;
ребро $e_3$: $(v_1,v_7)$ или $(v_7,v_1)$;
ребро $e_4$: $(v_1,v_{10})$ или $(v_{10},v_1)$;
ребро $e_5$: $(v_2,v_3)$ или $(v_3,v_2)$;
ребро $e_6$: $(v_2,v_5)$ или $(v_5,v_2)$;
ребро $e_7$: $(v_2,v_7)$ или $(v_7,v_2)$;
ребро $e_8$: $(v_3,v_4)$ или $(v_4,v_3)$;
ребро $e_9$: $(v_3,v_5)$ или $(v_5,v_3)$;
ребро $e_{10}$: $(v_3,v_9)$ или $(v_9,v_3)$;
ребро $e_{11}$: $(v_4,v_5)$ или $(v_5,v_4)$;
ребро $e_{12}$: $(v_4,v_6)$ или $(v_6,v_4)$;
ребро $e_{13}$: $(v_4,v_7)$ или $(v_7,v_4)$;
ребро $e_{14}$: $(v_4,v_9)$ или $(v_9,v_4)$;
ребро $e_{15}$: $(v_5,v_6)$ или $(v_6,v_5)$;
ребро $e_{16}$: $(v_5,v_8)$ или $(v_8,v_5)$;
ребро $e_{17}$: $(v_5,v_{10})$ или $(v_{10},v_5)$;
ребро $e_{18}$: $(v_6,v_7)$ или $(v_7,v_6)$;
ребро $e_{19}$: $(v_6,v_9)$ или $(v_9,v_6)$;
ребро $e_{20}$: $(v_7,v_8)$ или $(v_8,v_7)$;
ребро $e_{21}$: $(v_8,v_9)$ или $(v_9,v_8)$;
ребро $e_{22}$: $(v_8,v_{10})$ или $(v_{10},v_8)$;
ребро $e_{23}$: $(v_9,v_{10})$ или $(v_{10},v_9)$.

Множество изометрических циклов графа:

| **Циклы** | **Множество изометрических циклов графа в виде ребер:** | **Множество изометрических циклов графа в виде вершин:** |
|---|---|---|
| цикл $c_1$ | $\{e_1,e_2,e_6,e_{15}\}$ | $\{v_1,v_2,v_5,v_6\}$; |
| цикл $c_2$ | $\{e_1,e_3,e_7\}$ | $\{v_1,v_2,v_7\}$; |
| цикл $c_3$ | $\{e_1,e_4,e_6,e_{17}\}$ | $\{v_1,v_2,v_5,v_{10}\}$; |
| цикл $c_4$ | $\{e_2,e_3,e_{18}\}$ | $\{v_1,v_6,v_7\}$; |
| цикл $c_5$ | $\{e_2,e_4,e_{15},e_{17}\}$ | $\{v_1,v_5,v_6,v_{10}\}$; |
| цикл $c_6$ | $\{e_2,e_4,e_{19},e_{23}\}$ | $\{v_1,v_6,v_9,v_{10}\}$; |
| цикл $c_7$ | $\{e_3,e_4,e_{20},e_{22}\}$ | $\{v_1,v_7,v_8,v_{10}\}$; |
| цикл $c_8$ | $\{e_1,e_2,e_5,e_{10},e_{19}\}$ | $\{v_1,v_2,v_3,v_6,v_9\}$; |
| цикл $c_9$ | $\{e_1,e_4,e_5,e_{10},e_{23}\}$ | $\{v_1,v_2,v_3,v_9,v_{10}\}$; |
| цикл $c_{10}$ | $\{e_5,e_6,e_9\}$ | $\{v_2,v_3,v_5\}$; |
| цикл $c_{11}$ | $\{e_5,e_7,e_8,e_{13}\}$ | $\{v_2,v_3,v_4,v_7\}$; |
| цикл $c_{12}$ | $\{e_6,e_7,e_{11},e_{13}\}$ | $\{v_2,v_4,v_5,v_7\}$; |
| цикл $c_{13}$ | $\{e_6,e_7,e_{15},e_{18}\}$ | $\{v_2,v_5,v_6,v_7\}$; |
| цикл $c_{14}$ | $\{e_6,e_7,e_{16},e_{20}\}$ | $\{v_2,v_5,v_7,v_8\}$; |
| цикл $c_{15}$ | $\{e_8,e_9,e_{11}\}$ | $\{v_3,v_4,v_5\}$; |
| цикл $c_{16}$ | $\{e_8,e_{10},e_{14}\}$ | $\{v_3,v_4,v_9\}$; |
| цикл $c_{17}$ | $\{e_9,e_{10},e_{15},e_{19}\}$ | $\{v_3,v_5,v_6,v_9\}$; |
| цикл $c_{18}$ | $\{e_9,e_{10},e_{16},e_{21}\}$ | $\{v_3,v_5,v_8,v_9\}$; |
| цикл $c_{19}$ | $\{e_9,e_{10},e_{17},e_{23}\}$ | $\{v_3,v_5,v_9,v_{10}\}$; |
| цикл $c_{20}$ | $\{e_{11},e_{12},e_{15}\}$ | $\{v_4,v_5,v_6\}$; |
| цикл $c_{21}$ | $\{e_{11},e_{13},e_{16},e_{20}\}$ | $\{v_4,v_5,v_7,v_8\}$; |
| цикл $c_{22}$ | $\{e_{11},e_{14},e_{16},e_{21}\}$ | $\{v_4,v_5,v_8,v_9\}$; |
| цикл $c_{23}$ | $\{e_{11},e_{14},e_{17},e_{23}\}$ | $\{v_4,v_5,v_9,v_{10}\}$; |
| цикл $c_{24}$ | $\{e_{12},e_{13},e_{18}\}$ | $\{v_4,v_6,v_7\}$; |
| цикл $c_{25}$ | $\{e_{12},e_{14},e_{19}\}$ | $\{v_4,v_6,v_9\}$; |
| цикл $c_{26}$ | $\{e_{13},e_{14},e_{20},e_{21}\}$ | $\{v_4,v_7,v_8,v_9\}$; |
| цикл $c_{27}$ | $\{e_{15},e_{16},e_{18},e_{20}\}$ | $\{v_5,v_6,v_7,v_8\}$; |
| цикл $c_{28}$ | $\{e_{15},e_{16},e_{19},e_{21}\}$ | $\{v_5,v_6,v_8,v_9\}$; |
| цикл $c_{29}$ | $\{e_{15},e_{17},e_{19},e_{23}\}$ | $\{v_5,v_6,v_9,v_{10}\}$; |
| цикл $c_{30}$ | $\{e_{16},e_{17},e_{22}\}$ | $\{v_5,v_8,v_{10}\}$; |
| цикл $c_{31}$ | $\{e_{18},e_{19},e_{20},e_{21}\}$ | $\{v_6,v_7,v_8,v_9\}$; |
| цикл $c_{32}$ | $\{e_{21},e_{22},e_{23}\}$ | $\{v_8,v_9,v_{10}\}$. |

***Пример 3.5.*** Выделить множество изометрических циклов в графе $G_2$.

Количество вершин графа = 12
Количество ребер графа = 35
Количество изометрических циклов графа = 56

Матрица смежностей графа, представленная в виде списка смежных вершин:

вершина $v_1$: $\{v_2,v_3,v_4,v_5,v_6,v_7,v_8,v_9,v_{10},v_{12}\}$;
вершина $v_2$: $\{v_1,v_3,v_4,v_7,v_8,v_{12}\}$;
вершина $v_3$: $\{v_1,v_2,v_4,v_6,v_7,v_8,v_{10}\}$;
вершина $v_4$: $\{v_1,v_2,v_3,v_6,v_{11}\}$;
вершина $v_5$: $\{v_1,v_8,v_9\}$;
вершина $v_6$: $\{v_1,v_3,v_4,v_7,v_9,v_{12}\}$;
вершина $v_7$: $\{v_1,v_2,v_3,v_6,v_8,v_9,v_{11},v_{12}\}$;
вершина $v_8$: $\{v_1,v_2,v_3,v_5,v_7,v_{10},v_{11}\}$;
вершина $v_9$: $\{v_1,v_5,v_6,v_7,v_{10}\}$;
вершина $v_{10}$: $\{v_1,v_3,v_8,v_9,v_{12}\}$;
вершина $v_{11}$: $\{v_4,v_7,v_8\}$;
вершина $v_{12}$: $\{v_1,v_2,v_6,v_7,v_{10}\}$.

Элементы матрицы инциденций:

ребро $e_1$: $(v_1,v_2)$ или $(v_2,v_1)$;
ребро $e_2$: $(v_1,v_3)$ или $(v_3,v_1)$;
ребро $e_3$: $(v_1,v_4)$ или $(v_4,v_1)$;
ребро $e_4$: $(v_1,v_5)$ или $(v_5,v_1)$;
ребро $e_5$: $(v_1,v_6)$ или $(v_6,v_1)$;
ребро $e_6$: $(v_1,v_7)$ или $(v_7,v_1)$;
ребро $e_7$: $(v_1,v_8)$ или $(v_8,v_1)$;
ребро $e_8$: $(v_1,v_9)$ или $(v_9,v_1)$;
ребро $e_9$: $(v_1,v_{10})$ или $(v_{10},v_1)$;
ребро $e_{10}$: $(v_1,v_{12})$ или $(v_{12},v_1)$;
ребро $e_{11}$: $(v_2,v_3)$ или $(v_3,v_2)$;
ребро $e_{12}$: $(v_2,v_4)$ или $(v_4,v_2)$;
ребро $e_{13}$: $(v_2,v_7)$ или $(v_7,v_2)$;
ребро $e_{14}$: $(v_2,v_8)$ или $(v_8,v_2)$;
ребро $e_{15}$: $(v_2,v_{12})$ или $(v_{12},v_2)$;
ребро $e_{16}$: $(v_3,v_4)$ или $(v_4,v_3)$;
ребро $e_{17}$: $(v_3,v_6)$ или $(v_6,v_3)$;
ребро $e_{18}$: $(v_3,v_7)$ или $(v_7,v_3)$;
ребро $e_{19}$: $(v_3,v_8)$ или $(v_8,v_3)$;
ребро $e_{20}$: $(v_3,v_{10})$ или $(v_{10},v_3)$;
ребро $e_{21}$: $(v_4,v_6)$ или $(v_6,v_4)$;
ребро $e_{22}$: $(v_4,v_{11})$ или $(v_{11},v_4)$;
ребро $e_{23}$: $(v_5,v_8)$ или $(v_8,v_5)$;
ребро $e_{24}$: $(v_5,v_9)$ или $(v_9,v_5)$;
ребро $e_{25}$: $(v_6,v_7)$ или $(v_7,v_6)$;
ребро $e_{26}$: $(v_6,v_9)$ или $(v_9,v_6)$;
ребро $e_{27}$: $(v_6,v_{12})$ или $(v_{12},v_6)$;
ребро $e_{28}$: $(v_7,v_8)$ или $(v_8,v_7)$;
ребро $e_{29}$: $(v_7,v_9)$ или $(v_9,v_7)$;
ребро $e_{30}$: $(v_7,v_{11})$ или $(v_{11},v_7)$;
ребро $e_{31}$: $(v_7,v_{12})$ или $(v_{12},v_7)$;
ребро $e_{32}$: $(v_8,v_{10})$ или $(v_{10},v_8)$;
ребро $e_{33}$: $(v_8,v_{11})$ или $(v_{11},v_8)$;
ребро $e_{34}$: $(v_9,v_{10})$ или $(v_{10},v_9)$;

ребро $e_{35}$: $(v_{10},v_{12})$ или $(v_{12},v_{10})$.

## Множество изометрических циклов графа:

| Циклы | Множество изометрических циклов графа в виде ребер: | Множество изометрических циклов графа в виде вершин: |
|---|---|---|
| цикл $c_1$ | $\{e_1,e_2,e_{11}\}$; | $\{v_1,v_2,v_3\}$; |
| цикл $c_2$ | $\{e_1,e_3,e_{12}\}$; | $\{v_1,v_2,v_4\}$; |
| цикл $c_3$ | $\{e_1,e_6,e_{13}\}$; | $\{v_1,v_2,v_7\}$; |
| цикл $c_4$ | $\{e_1,e_7,e_{14}\}$; | $\{v_1,v_2,v_8\}$; |
| цикл $c_5$ | $\{e_1,e_{10},e_{15}\}$; | $\{v_1,v_2,v_{12}\}$; |
| цикл $c_6$ | $\{e_2,e_3,e_{16}\}$; | $\{v_1,v_3,v_4\}$; |
| цикл $c_7$ | $\{e_2,e_5,e_{17}\}$; | $\{v_1,v_3,v_6\}$; |
| цикл $c_8$ | $\{e_2,e_6,e_{18}\}$; | $\{v_1,v_3,v_7\}$; |
| цикл $c_9$ | $\{e_2,e_7,e_{19}\}$; | $\{v_1,v_3,v_8\}$; |
| цикл $c_{10}$ | $\{e_2,e_9,e_{20}\}$; | $\{v_1,v_3,v_{10}\}$; |
| цикл $c_{11}$ | $\{e_3,e_5,e_{21}\}$; | $\{v_1,v_4,v_6\}$; |
| цикл $c_{12}$ | $\{e_3,e_6,e_{22},e_{30}\}$; | $\{v_1,v_4,v_7,v_{11}\}$; |
| цикл $c_{13}$ | $\{e_3,e_7,e_{22},e_{33}\}$; | $\{v_1,v_4,v_8,v_{11}\}$; |
| цикл $c_{14}$ | $\{e_4,e_7,e_{23}\}$; | $\{v_1,v_5,v_8\}$; |
| цикл $c_{15}$ | $\{e_4,e_8,e_{24}\}$; | $\{v_1,v_5,v_9\}$; |
| цикл $c_{16}$ | $\{e_5,e_6,e_{25}\}$; | $\{v_1,v_6,v_7\}$; |
| цикл $c_{17}$ | $\{e_5,e_8,e_{26}\}$; | $\{v_1,v_6,v_9\}$; |
| цикл $c_{18}$ | $\{e_5,e_{10},e_{27}\}$; | $\{v_1,v_6,v_{12}\}$; |
| цикл $c_{19}$ | $\{e_6,e_7,e_{28}\}$; | $\{v_1,v_7,v_8\}$; |
| цикл $c_{20}$ | $\{e_6,e_8,e_{29}\}$; | $\{v_1,v_7,v_9\}$; |
| цикл $c_{21}$ | $\{e_6,e_{10},e_{31}\}$; | $\{v_1,v_7,v_{12}\}$; |
| цикл $c_{22}$ | $\{e_7,e_9,e_{32}\}$; | $\{v_1,v_8,v_{10}\}$; |
| цикл $c_{23}$ | $\{e_8,e_9,e_{34}\}$; | $\{v_1,v_9,v_{10}\}$; |
| цикл $c_{24}$ | $\{e_9,e_{10},e_{35}\}$; | $\{v_1,v_{10},v_{12}\}$; |
| цикл $c_{25}$ | $\{e_{11},e_{12},e_{16}\}$; | $\{v_2,v_3,v_4\}$; |
| цикл $c_{26}$ | $\{e_{11},e_{13},e_{18}\}$; | $\{v_2,v_3,v_7\}$; |
| цикл $c_{27}$ | $\{e_{11},e_{14},e_{19}\}$; | $\{v_2,v_3,v_8\}$; |
| цикл $c_{28}$ | $\{e_{11},e_{15},e_{17},e_{27}\}$; | $\{v_2,v_3,v_6,v_{12}\}$; |
| цикл $c_{29}$ | $\{e_{11},e_{15},e_{20},e_{35}\}$; | $\{v_2,v_3,v_{10},v_{12}\}$; |
| цикл $c_{30}$ | $\{e_{12},e_{13},e_{21},e_{25}\}$; | $\{v_2,v_4,v_6,v_7\}$; |
| цикл $c_{31}$ | $\{e_{12},e_{13},e_{22},e_{30}\}$; | $\{v_2,v_4,v_7,v_{11}\}$; |
| цикл $c_{32}$ | $\{e_{12},e_{14},e_{22},e_{33}\}$; | $\{v_2,v_4,v_8,v_{11}\}$; |
| цикл $c_{33}$ | $\{e_{12},e_{15},e_{21},e_{27}\}$; | $\{v_2,v_4,v_6,v_{12}\}$; |
| цикл $c_{34}$ | $\{e_{13},e_{14},e_{28}\}$; | $\{v_2,v_7,v_8\}$; |
| цикл $c_{35}$ | $\{e_{13},e_{15},e_{31}\}$; | $\{v_2,v_7,v_{12}\}$; |
| цикл $c_{36}$ | $\{e_{14},e_{15},e_{32},e_{35}\}$; | $\{v_2,v_8,v_{10},v_{12}\}$; |
| цикл $c_{37}$ | $\{e_{16},e_{17},e_{21}\}$; | $\{v_3,v_4,v_6\}$; |
| цикл $c_{38}$ | $\{e_{16},e_{18},e_{22},e_{30}\}$; | $\{v_3,v_4,v_7,v_{11}\}$; |
| цикл $c_{39}$ | $\{e_{16},e_{19},e_{22},e_{33}\}$; | $\{v_3,v_4,v_8,v_{11}\}$; |
| цикл $c_{40}$ | $\{e_{17},e_{18},e_{25}\}$; | $\{v_3,v_6,v_7\}$; |
| цикл $c_{41}$ | $\{e_{17},e_{20},e_{26},e_{34}\}$; | $\{v_3,v_6,v_9,v_{10}\}$; |
| цикл $c_{42}$ | $\{e_{17},e_{20},e_{27},e_{35}\}$; | $\{v_3,v_6,v_{10},v_{12}\}$; |
| цикл $c_{43}$ | $\{e_{18},e_{19},e_{28}\}$; | $\{v_3,v_7,v_8\}$; |
| цикл $c_{44}$ | $\{e_{18},e_{20},e_{29},e_{34}\}$; | $\{v_3,v_7,v_9,v_{10}\}$; |
| цикл $c_{45}$ | $\{e_{18},e_{20},e_{31},e_{35}\}$; | $\{v_3,v_7,v_{10},v_{12}\}$; |
| цикл $c_{46}$ | $\{e_{19},e_{20},e_{32}\}$; | $\{v_3,v_8,v_{10}\}$; |
| цикл $c_{47}$ | $\{e_{21},e_{22},e_{25},e_{30}\}$; | $\{v_4,v_6,v_7,v_{11}\}$; |
| цикл $c_{48}$ | $\{e_{23},e_{24},e_{28},e_{29}\}$; | $\{v_5,v_7,v_8,v_9\}$; |
| цикл $c_{49}$ | $\{e_{23},e_{24},e_{32},e_{34}\}$; | $\{v_5,v_8,v_9,v_{10}\}$; |
| цикл $c_{50}$ | $\{e_{25},e_{26},e_{29}\}$; | $\{v_6,v_7,v_9\}$; |
| цикл $c_{51}$ | $\{e_{25},e_{27},e_{31}\}$; | $\{v_6,v_7,v_{12}\}$; |
| цикл $c_{52}$ | $\{e_{26},e_{27},e_{34},e_{35}\}$; | $\{v_6,v_9,v_{10},v_{12}\}$; |
| цикл $c_{53}$ | $\{e_{28},e_{29},e_{32},e_{34}\}$; | $\{v_7,v_8,v_9,v_{10}\}$; |
| цикл $c_{54}$ | $\{e_{28},e_{30},e_{33}\}$; | $\{v_7,v_8,v_{11}\}$; |
| цикл $c_{55}$ | $\{e_{28},e_{31},e_{32},e_{35}\}$; | $\{v_7,v_8,v_{10},v_{12}\}$; |
| цикл $c_{56}$ | $\{e_{29},e_{31},e_{34},e_{35}\}$. | $\{v_7,v_9,v_{10},v_{12}\}$. |

## Комментарии

Используя метрические свойства теории графов, введено понятие изометрического цикла графа. Рассмотрены основные свойства изометрических циклов графа. Представлен алгоритм выделения множества изометрических циклов графа $C_\tau(G)$. Показано, что выделение всего множества изометрических циклов должно осуществляться на всем множестве ребер графа G. Доказана теорема о существовании базиса в подпространстве циклов графа C(G), состоящего только из изометрических циклов.

Следует заметить, что множество изометрических циклов графа, является уникальной характеристикой графа, и может раммматриваться как некий инвариант.

Вычислительная сложность алгоритма выделения множества изометрических циклов рарна $O(n^4)$.

# Глава 4. СПЕКТР РЕБЕРНЫХ ЦИКЛОВ ГРАФА

## 4.1. Реберные циклы графа и их свойства

Ранее мы рассмотрели свойства суграфов принадлежащих подпространству разрезов графа. Рассмотрим свойства суграфов принадлежащих подпространству циклов графа.

В качестве примера рассмотрим следующий граф $G_8$ (рис. 4.1).

Выделим множество изометрических циклов графа $C_\tau$.

**Определение 4.1.** *Обод графа* определим как сумму по модулю 2 (кольцевая сумма) всех изометрических циклов графа.

$$c_0 = \sum_{i=0}^{k} c_i \, ; \tag{4.1}$$

где $k$ – мощность множества изометрических циклов графа $k = card\, C_\tau$. Ободом графа может быть пустое множество.

Определим суммарный цикл, проходящий по $i$-тому ребру графа. Для этого произведем кольцевое суммирование всех циклов имеющих в своем составе ребру $e_i$, где $i$ = 1,2,…,$m$. Обозначим такой цикл $\tau_0$ ($e_i$), и будем называть его базовым реберным циклом ребра $e_i$.

**Определение 4.2**. *Базовым реберным циклом* ребра $e_i$ называется кольцевая сумма изометрических циклов и обода графа, имеющих в своем составе ребро $e_i$.

Пусть задан следующий граф (рис. 4.1).

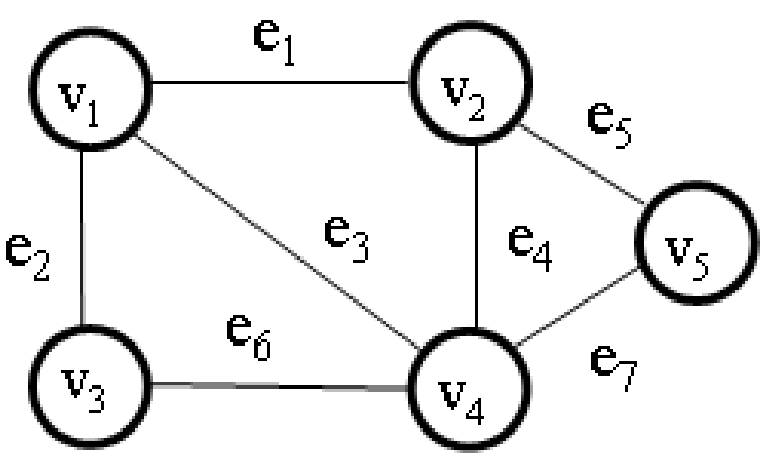

Рис. 4.1. Граф $G_8$.

Множество изометрических циклов графа $G_8$:

$c_1=\{e_2,e_3,e_6\}$;
$c_2=\{e_1,e_3,e_4\}$;
$c_3=\{e_4,e_5,e_7\}$;

$c_0=\{e_1,e_2,e_5,e_6,e_7\}$ – обод графа (кольцевая сумма множества изометрических циклов графа).

Для нашего примера множество базовых реберных циклов имеет вид:

$\tau_0$ ($e_1$) = $c_2 \oplus c_0$ = $\{e_1,e_3,e_4\} \oplus \{e_1,e_2,e_5,e_6,e_7\}$ = $\{e_2,e_3,e_4,e_5,e_6,e_7\}$;

$\tau_0$ ($e_2$) = $c_1 \oplus c_0$ = $\{e_2,e_3,e_6\} \oplus \{e_1,e_2,e_5,e_6,e_7\}$ = $\{e_1,e_3,e_5,e_7\}$;

$\tau_0$ ($e_3$) = $c_1 \oplus c_2$ = {$e_2,e_3,e_6$} $\oplus$ {$e_1,e_3,e_4$} = {$e_1,e_2,e_4,e_6$};
$\tau_0$ ($e_4$) = $c_2 \oplus c_3$ = {$e_1,e_3,e_4$} $\oplus$ {$e_4,e_5,e_7$} = {$e_1,e_3,e_5,e_7$};
$\tau_0$ ($e_5$) = $c_3 \oplus c_0$ = {$e_4,e_5,e_7$} $\oplus$ {$e_1,e_2,e_5,e_6,e_7$} = {$e_1,e_2,e_4,e_6$};
$\tau_0$ ($e_6$) = $c_1 \oplus c_0$ = {$e_2,e_3,e_6$} $\oplus$ {$e_1,e_2,e_5,e_6,e_7$} = {$e_1,e_3,e_5,e_7$};
$\tau_0$ ($e_7$) = $c_3 \oplus c_0$ = {$e_4,e_5,e_7$} $\oplus$ {$e_1,e_2,e_5,e_6,e_7$} = {$e_1,e_2,e_4,e_6$}.

Определим суграфы следующих уровней с помощью преобразования:

$$\tau_k(e_\mathrm{i}) = \begin{cases} \tau_k(e_\mathrm{i}) = \gamma(\tau_{\mathrm{k}\text{-}1}(e_\mathrm{i})); \\ \textit{если } \tau_\mathrm{k}(e_\mathrm{i}) = \tau_\mathrm{h}(e_\mathrm{i}), \textit{ то } \tau_\mathrm{k}(e_\mathrm{i}) = \varnothing. \end{cases} \quad (4.2)$$

где $i$ – номер ребра, $k$ – номер уровня, $h$ – номер произвольного предыдущего уровня.

Преобразование $\gamma(\tau_{k-1}(e_i))$ определим следующим образом:

$$\gamma(\{a,b,c,...,m\}) = \tau_0(a) \oplus \tau_0(\mathrm{b}) \oplus \tau_0(c) \oplus ... \oplus \tau_0(m)\,. \quad (4.3)$$

$\gamma(\tau_{k-1}(e_i))$ – *порождение реберного цикла последующего уровня* определяет строку из предыдущего реберного разреза суграфа {$a,b,c,...,m$} ребра $e_i$ состоящего из множества ребер базовых реберных циклов $\tau_0(a) \oplus \tau_0(\mathrm{b}) \oplus \tau_0(c) \oplus ... \oplus \tau_0(m)$.

По аналогии с суграфами реберных разрезов, следующее порожденное подмножество суграфов характеризует уровень реберных циклов (4.2):

$\tau_1$($e_1$) = $\gamma$ {$e_2,e_3,e_4,e_5,e_6,e_7$} = $\varnothing$;
$\tau_1$($e_2$) = $\gamma$ {$e_1,e_3,e_5,e_7$} = {$e_1,e_3,e_5,e_7$} = $\varnothing$;
$\tau_1$($e_3$) = $\gamma$ {$e_1,e_2,e_4,e_6$} = {$e_1,e_2,e_4,e_6$} = $\varnothing$;
$\tau_1$($e_4$) = $\gamma$ {$e_1,e_3,e_5,e_7$} = {$e_1,e_3,e_5,e_7$} = $\varnothing$;
$\tau_1$($e_5$) = $\gamma$ {$e_1,e_2,e_4,e_6$} = {$e_1,e_2,e_4,e_6$} = $\varnothing$;
$\tau_1$($e_6$) = $\gamma$ {$e_1,e_3,e_5,e_7$} = {$e_1,e_3,e_5,e_7$} = $\varnothing$;
$\tau_1$($e_7$) = $\gamma$ {$e_1,e_2,e_4,e_6$} = {$e_1,e_2,e_4,e_6$} = $\varnothing$.

Отсюда, по аналогии, вычисляем кортеж весов ребер $\xi_\tau$($G_2$) и кортеж весов вершин $\zeta_\tau$($G_2$), так как имеется всего один базовый ярус:

$\xi_\tau$($G_2$) = <6,4,4,4,4,4,4>;
$\zeta_\tau$($G_2$) = <14,14,8,16,8>.

## 4.2. Спектр реберных циклов графа

Все уровневые подмножества (реберные циклы) также будут представлять собой квазициклы исходного графа G. Так нулевой уровень будет состоять из множества базовых реберных циклов графа G, а все последующие реберные циклы порождаться преобразованием $\gamma(\tau_{k-1}(e_i))$. После получения циклического повтора подмножеств последующему повторяющемуся подмножеству ставится в соответствие пустое множество.

На каждом уровне образуются $m$ подмножеств (суграфов), зависящих от выбранного ребра, что дает нам возможность построить прямоугольную матрицу спектра реберных циклов Т размером $m \times k$, где $k$ – количество уровней. Очевидно, что элементы этой матрицы могут быть записаны в виде суграфов исходного графа. Таким образом, реберные циклы графа состоят из подмножества ребер и порождают ярусные квазициклы графа. Любой реберный цикл можно представить в виде функции $w_l(e_i)$ зависящей от ребра $e_i$ применяемого в качестве аргумента. Введем предварительно понятие спектра реберных циклов.

**Определение 4.3.** Совокупность всех реберных циклов порожденных множеством базовых реберных циклов будем называть *спектром реберных циклов графа.*

уровни

T(G) =

| | $l_0$ | $l_1$ | .. | $l_k$ |
|---|---|---|---|---|
| $e_1$ | $\tau_0(e_1)$ | $\tau_1(e_1)$ | … | $\tau_k(e_1)$ |
| $e_2$ | $\tau_0(e_2)$ | $\tau_1(e_2)$ | … | $\tau_k(e_2)$ |
| $e_3$ | $\tau_0(e_3)$ | $\tau_1(e_3)$ | … | $\tau_k(e_3)$ |
| … | … | … | … | … |
| $e_m$ | $\tau_0(e_m)$ | $\tau_1(e_m)$ | … | $\tau_k(e_m$ |

Еще раз отметим, что ободом графа может быть пустое множество.

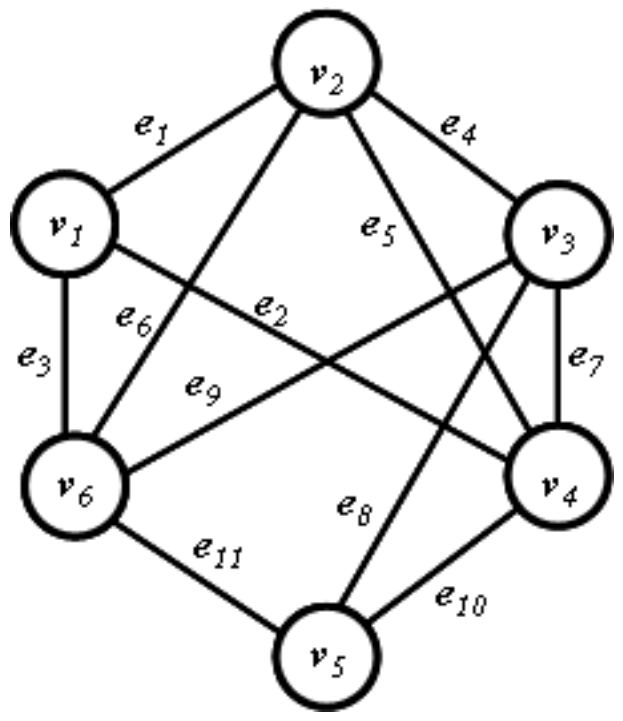

Рис. 4.2. Граф $G_2$.

Рассмотрим тот же граф $G_2$, что и спектра реберных разрезов (рис. 4.2). Множество изометрических циклов графа $G_2$:

$c_1 = \{e_1,e_2,e_5\}$;
$c_2 = \{e_1,e_3,e_6\}$;
$c_3 = \{e_2,e_3,e_7,e_9\}$;
$c_4 = \{e_2,e_3,e_{10},e_{11}\}$;
$c_5 = \{e_4,e_5,e_7\}$;
$c_6 = \{e_4,e_6,e_9\}$;
$c_7 = \{e_5,e_6,e_{10},e_{11}\}$;
$c_8 = \{e_7,e_8,e_{10}\}$;

$c_9 = \{e_8,e_9,e_{11}\}$.

Обод графа состоит из множества изометрических циклов графа

$$c_0 = \sum_{i=1}^{9} c_i = \{e_2,e_3,e_5,e_6,e_7,e_9,e_{10},e_{11}\}.$$

Определим базовые реберные циклы для каждого ребра:

$\tau_0(e_1) = c_1 \oplus c_2 = \{e_1,e_2,e_5\} \oplus \{e_1,e_3,e_6\} = \{e_2,e_3,e_5,e_6\}$;

$\tau_0(e_2) = c_1 \oplus c_3 \oplus c_4 \oplus c_0 = \{e_1,e_2,e_5\} \oplus \{e_2,e_3,e_7,e_9\} \oplus \{e_2,e_3,e_{10},e_{11}\} \oplus$
$\oplus \{e_2,e_3,e_5,e_6,e_7,e_9,e_{10},e_{11}\} = \{e_1,e_3,e_6\}$;

$\tau_0(e_3) = c_2 \oplus c_3 \oplus c_4 \oplus c_0 = \{e_1,e_3,e_6\} \oplus \{e_2,e_3,e_7,e_9\} \oplus \{e_2,e_3,e_{10},e_{11}\} \oplus$
$\oplus \{e_2,e_3,e_5,e_6,e_7,e_9,e_{10},e_{11}\} = \{e_1,e_2,e_5\}$;

$\tau_0(e_4) = c_5 \oplus c_6 = \{e_4,e_5,e_7\} \oplus \{e_4,e_6,e_9\} = \{e_5,e_6,e_7,e_9\}$;

$\tau_0(e_5) = c_1 \oplus c_5 \oplus c_7 \oplus c_0 = \{e_1,e_2,e_5\} \oplus \{e_4,e_5,e_7\} \oplus \{e_5,e_6,e_{10},e_{11}\} \oplus \{e_2,e_3,e_5,e_6,e_7,e_9,e_{10},e_{11}\} =$
$= \{e_1,e_3,e_4,e_9\}$;

$\tau_0(e_6) = c_2 \oplus c_6 \oplus c_7 \oplus c_0 = \{e_1,e_3,e_6\} \oplus \{e_4,e_6,e_9\} \oplus \{e_5,e_6,e_{10},e_{11}\} \oplus \{e_2,e_3,e_5,e_6,e_7,e_9,e_{10},e_{11}\} =$
$= \{e_1,e_2,e_4,e_7\}$;

$\tau_0(e_7) = c_3 \oplus c_5 \oplus c_8 \oplus c_0 = \{e_2,e_3,e_7,e_9\} \oplus \{e_4,e_5,e_7\} \oplus \{e_7,e_8,e_{10}\} \oplus \{e_2,e_3,e_5,e_6,e_7,e_9,e_{10},e_{11}\} =$
$= \{e_4,e_6,e_8,e_{11}\}$;

$\tau_0(e_8) = c_8 \oplus c_9 = \{e_7,e_8,e_{10}\} \oplus \{e_8,e_9,e_{11}\} = \{e_7,e_9,e_{10},e_{11}\}$;

$\tau_0(e_9) = c_3 \oplus c_6 \oplus c_9 \oplus c_0 = \{e_2,e_3,e_7,e_9\} \oplus \{e_4,e_6,e_9\} \oplus \{e_8,e_9,e_{11}\} \oplus \{e_2,e_3,e_5,e_6,e_7,e_9,e_{10},e_{11}\} =$
$= \{e_4,e_5,e_8,e_{10}\}$;

$\tau_0(e_{10}) = c_4 \oplus c_7 \oplus c_8 \oplus c_0 = \{e_2,e_3,e_{10},e_{11}\} \oplus \{e_5,e_6,e_{10},e_{11}\} \oplus \{e_7,e_8,e_{10}\} \oplus$
$\oplus \{e_2,e_3,e_5,e_6,e_7,e_9,e_{10},e_{11}\} = \{e_8,e_9,e_{11}\}$;

$\tau_0(e_{11}) = c_4 \oplus c_7 \oplus c_9 \oplus c_0 = \{e_2,e_3,e_{10},e_{11}\} \oplus \{e_5,e_6,e_{10},e_{11}\} \oplus \{e_8,e_9,e_{11}\} \oplus$
$\oplus \{e_2,e_3,e_5,e_6,e_7,e_9,e_{10},e_{11}\} = \{e_7,e_8,e_{10}\}$.

Построим 1-й уровень реберных циклов:

$\tau_1(e_1) = \gamma\,\{e_2,e_3,e_5,e_6\} = \{e_1,e_3,e_6\} \oplus \{e_1,e_2,e_5\} \oplus \{e_1,e_3,e_4,e_9\} \oplus \{e_1,e_2,e_4,e_7\} = \{e_5,e_6,e_7,e_9\}$;

$\tau_1(e_2) = \gamma\,\{e_1,e_3,e_6\} = \{e_2,e_3,e_5,e_6\} \oplus \{e_1,e_2,e_5\} \oplus \{e_1,e_2,e_4,e_7\} = \{e_2,e_3,e_4,e_6,e_7\}$;

$\tau_1(e_3) = \gamma\,\{e_1,e_2,e_5\} = \{e_2,e_3,e_5,e_6\} \oplus \{e_1,e_3,e_6\} \oplus \{e_1,e_3,e_4,e_9\} = \{e_2,e_3,e_4,e_5,e_9\}$;

$\tau_1(e_4) = \gamma\,\{e_5,e_6,e_7,e_9\} = \{e_1,e_3,e_4,e_9\} \oplus \{e_1,e_2,e_4,e_7\} \oplus \{e_4,e_6,e_8,e_{11}\} \oplus \{e_4,e_5,e_8,e_{10}\} =$
$= \{e_2,e_3,e_5,e_6,e_7,e_9,e_{10},e_{11}\}$;

$\tau_1(e_5) = \gamma\,\{e_1,e_3,e_4,e_9\} = \{e_2,e_3,e_5,e_6\} \oplus \{e_1,e_2,e_5\} \oplus \{e_5,e_6,e_7,e_9\} \oplus \{e_4,e_5,e_8,e_{10}\} =$
$= \{e_1,e_3,e_4,e_7,e_8,e_9,e_{10}\}$;

$\tau_1(e_6) = \gamma\,\{e_1,e_2,e_4,e_7\} = \{e_2,e_3,e_5,e_6\} \oplus \{e_1,e_3,e_6\} \oplus \{e_5,e_6,e_7,e_9\} \oplus \{e_4,e_6,e_8,e_{11}\} =$
$= \{e_1,e_2,e_4,e_7,e_8,e_9,e_{11}\}$;

$\tau_1(e_7) = \gamma\,\{e_4,e_6,e_8,e_{11}\} = \{e_5,e_6,e_7,e_9\} \oplus \{e_1,e_2,e_4,e_7\} \oplus \{e_7,e_9,e_{10},e_{11}\} \oplus \{e_7,e_8,e_{10}\} =$
$= \{e_1,e_2,e_4,e_5,e_6,e_8,e_{11}\}$;

$\tau_1(e_8) = \gamma\,\{e_7,e_9,e_{10},e_{11}\} = \{e_4,e_6,e_8,e_{11}\} \oplus \{e_4,e_5,e_8,e_{10}\} \oplus \{e_8,e_9,e_{11}\} \oplus \{e_7,e_8,e_{10}\} =$
$= \{e_5,e_6,e_7,e_9\}$;

$\tau_1(e_9) = \gamma\,\{e_4,e_5,e_8,e_{10}\} = \{e_5,e_6,e_7,e_9\} \oplus \{e_1,e_3,e_4,e_9\} \oplus \{e_7,e_9,e_{10},e_{11}\} \oplus \{e_8,e_9,e_{11}\} =$
$= \{e_1,e_3,e_4,e_5,e_6,e_8,e_{10}\}$;

$\tau_1(e_{10}) = \gamma\,\{e_8,e_9,e_{11}\} = \{e_7,e_9,e_{10},e_{11}\} \oplus \{e_4,e_5,e_8,e_{10}\} \oplus \{e_7,e_8,e_{10}\} = \{e_4,e_5,e_9,e_{10},e_{11}\}$;

$\tau_1(e_{11}) = \gamma\,\{e_7,e_8,e_{10}\} = \{e_4,e_6,e_8,e_{11}\} \oplus \{e_7,e_9,e_{10},e_{11}\} \oplus \{e_8,e_9,e_{11}\} = \{e_4,e_6,e_7,e_{10},e_{11}\}$.

Построим 2-й уровень реберных циклов:

$\tau_2(e_1) = \gamma\{e_5,e_6,e_7,e_9\} = \{e_1,e_3,e_4,e_9\} \oplus \{e_1,e_2,e_4,e_7\} \oplus \{e_4,e_6,e_8,e_{11}\} \oplus \{e_4,e_5,e_8,e_{10}\} =$
$= \{e_2,e_3,e_5,e_6,e_7,e_9,e_{10},e_{11}\}$;
$\tau_2(e_2) = \gamma\{e_2,e_3,e_4,e_6,e_7\} = \{e_1,e_3,e_6\} \oplus \{e_1,e_2,e_5\} \oplus \{e_5,e_6,e_7,e_9\} \oplus \{e_1,e_2,e_4,e_7\} \oplus$
$\oplus \{e_4,e_6,e_8,e_{11}\} = \{e_1,e_3,e_6,e_8,e_9,e_{11}\}$;
$\tau_2(e_3) = \gamma\{e_2,e_3,e_4,e_5,e_9\} = \{e_1,e_3,e_6\} \oplus \{e_1,e_2,e_5\} \oplus \{e_5,e_6,e_7,e_9\} \oplus \{e_1,e_3,e_4,e_9\} \oplus$
$\oplus \{e_4,e_5,e_8,e_{10}\} = \{e_1,e_2,e_5,e_7,e_8,e_{10}\}$;
$\tau_2(e_4) = \gamma\{e_2,e_3,e_5,e_6,e_7,e_9,e_{10},e_{11}\} = \{e_1,e_3,e_6\} \oplus \{e_1,e_2,e_5\} \oplus \{e_1,e_3,e_4,e_9\} \oplus \{e_1,e_2,e_4,e_7\} \oplus$
$\oplus \{e_4,e_6,e_8,e_{11}\} \oplus \{e_4,e_5,e_8,e_{10}\} \oplus \{e_8,e_9,e_{11}\} \oplus \{e_7,e_8,e_{10}\} = \varnothing$;
$\tau_2(e_5) = \gamma\{e_1,e_3,e_4,e_7,e_8,e_9,e_{10}\} = \{e_2,e_3,e_5,e_6\} \oplus \{e_1,e_2,e_5\} \oplus \{e_5,e_6,e_7,e_9\} \oplus \{e_4,e_6,e_8,e_{11}\} \oplus$
$\oplus \{e_7,e_9,e_{10},e_{11}\} \oplus \{e_4,e_5,e_8,e_{10}\} \oplus \{e_8,e_9,e_{11}\} = \{e_1,e_3,e_6,e_8,e_9,e_{11}\}$;
$\tau_2(e_6) = \gamma\{e_1,e_2,e_4,e_7,e_8,e_9,e_{11}\} = \{e_2,e_3,e_5,e_6\} \oplus \{e_1,e_3,e_6\} \oplus \{e_5,e_6,e_7,e_9\} \oplus \{e_4,e_6,e_8,e_{11}\} \oplus$
$\oplus \{e_7,e_9,e_{10},e_{11}\} \oplus \{e_4,e_5,e_8,e_{10}\} \oplus \{e_7,e_8,e_{10}\} = \{e_1,e_2,e_5,e_7,e_8,e_{10}\}$;
$\tau_2(e_7) = \gamma\{e_1,e_2,e_4,e_5,e_6,e_8,e_{11}\} = \{e_2,e_3,e_5,e_6\} \oplus \{e_1,e_3,e_6\} \oplus \{e_5,e_6,e_7,e_9\} \oplus \{e_1,e_3,e_4,e_9\} \oplus$
$\oplus \{e_1,e_2,e_4,e_7\} \oplus \{e_7,e_9,e_{10},e_{11}\} \oplus \{e_7,e_8,e_{10}\} = \{e_1,e_3,e_6,e_8,e_9,e_{11}\}$;
$\tau_2(e_8) = \gamma\{e_5,e_6,e_7,e_9\} = \{e_1,e_3,e_4,e_9\} \oplus \{e_1,e_2,e_4,e_7\} \oplus \{e_4,e_6,e_8,e_{11}\} \oplus \{e_4,e_5,e_8,e_{10}\} =$
$= \{e_2,e_3,e_5,e_6,e_7,e_9,e_{10},e_{11}\}$;
$\tau_2(e_9) = \gamma\{e_1,e_3,e_4,e_5,e_6,e_8,e_{10}\} = \{e_2,e_3,e_5,e_6\} \oplus \{e_1,e_2,e_5\} \oplus \{e_5,e_6,e_7,e_9\} \oplus \{e_1,e_3,e_4,e_9\} \oplus$
$\oplus \{e_1,e_2,e_4,e_7\} \oplus \{e_7,e_9,e_{10},e_{11}\} \oplus \{e_8,e_9,e_{11}\} = \{e_1,e_2,e_5,e_7,e_8,e_{10}\}$;
$\tau_2(e_{10}) = \gamma\{e_4,e_5,e_9,e_{10},e_{11}\} = \{e_5,e_6,e_7,e_9\} \oplus \{e_1,e_3,e_4,e_9\} \oplus \{e_4,e_5,e_8,e_{10}\} \oplus$
$\oplus \{e_8,e_9,e_{11}\} \oplus \{e_7,e_8,e_{10}\} = \{e_1,e_3,e_6,e_8,e_9,e_{11}\}$;
$\tau_2(e_{11}) = \gamma\{e_4,e_6,e_7,e_{10},e_{11}\} = \{e_5,e_6,e_7,e_9\} \oplus \{e_1,e_2,e_4,e_7\} \oplus \{e_4,e_6,e_8,e_{11}\} \oplus$
$\oplus \{e_8,e_9,e_{11}\} \oplus \{e_7,e_8,e_{10}\} = \{e_1,e_2,e_5,e_7,e_8,e_{10}\}$.

Строим 3-й уровень реберных циклов:

$\tau_3(e_1) = \gamma\{e_2,e_3,e_5,e_6,e_7,e_9,e_{10},e_{11}\} = \varnothing$;
$\tau_3(e_2) = \gamma\{e_1,e_3,e_6,e_8,e_9,e_{11}\} = \{e_2,e_3,e_5,e_6\} \oplus \{e_1,e_2,e_5\} \oplus \{e_1,e_2,e_4,e_7\} \oplus \{e_7,e_9,e_{10},e_{11}\} \oplus$
$\oplus \{e_4,e_5,e_8,e_{10}\} \oplus \{e_7,e_8,e_{10}\} = \{e_2,e_3,e_5,e_6,e_7,e_9,e_{10},e_{11}\}$;
$\tau_3(e_3) = \gamma\{e_1,e_2,e_5,e_7,e_8,e_{10}\} = \{e_2,e_3,e_5,e_6\} \oplus \{e_1,e_3,e_6\} \oplus \{e_1,e_3,e_4,e_9\} \oplus \{e_4,e_6,e_8,e_{11}\} \oplus$
$\oplus \{e_7,e_9,e_{10},e_{11}\} \oplus \{e_8,e_9,e_{11}\} = \{e_2,e_3,e_5,e_6,e_7,e_9,e_{10},e_{11}\}$;
$\tau_3(e_4) = \varnothing$;
$\tau_3(e_5) = \gamma\{e_1,e_3,e_6,e_8,e_9,e_{11}\} = \{e_2,e_3,e_5,e_6\} \oplus \{e_1,e_2,e_5\} \oplus \{e_1,e_2,e_4,e_7\} \oplus \{e_7,e_9,e_{10},e_{11}\} \oplus$
$\oplus \{e_4,e_5,e_8,e_{10}\} \oplus \{e_7,e_8,e_{10}\} = \{e_2,e_3,e_5,e_6,e_7,e_9,e_{10},e_{11}\}$;
$\tau_3(e_6) = \gamma\{e_1,e_2,e_5,e_7,e_8,e_{10}\} = \{e_2,e_3,e_5,e_6\} \oplus \{e_1,e_3,e_6\} \oplus \{e_1,e_3,e_4,e_9\} \oplus \{e_4,e_6,e_8,e_{11}\} \oplus$
$\oplus \{e_7,e_9,e_{10},e_{11}\} \oplus \{e_8,e_9,e_{11}\} = \{e_2,e_3,e_5,e_6,e_7,e_9,e_{10},e_{11}\}$;
$\tau_3(e_7) = \gamma\{e_1,e_3,e_6,e_8,e_9,e_{11}\} = \{e_2,e_3,e_5,e_6\} \oplus \{e_1,e_2,e_5\} \oplus \{e_1,e_2,e_4,e_7\} \oplus \{e_7,e_9,e_{10},e_{11}\} \oplus$
$\oplus \{e_4,e_5,e_8,e_{10}\} \oplus \{e_7,e_8,e_{10}\} = \{e_2,e_3,e_5,e_6,e_7,e_9,e_{10},e_{11}\}$;
$\tau_3(e_8) = \gamma\{e_2,e_3,e_5,e_6,e_7,e_9,e_{10},e_{11}\} = \varnothing$;
$\tau_3(e_9) = \gamma\{e_1,e_2,e_5,e_7,e_8,e_{10}\} = \{e_2,e_3,e_5,e_6\} \oplus \{e_1,e_3,e_6\} \oplus \{e_1,e_3,e_4,e_9\} \oplus \{e_4,e_6,e_8,e_{11}\} \oplus$
$\oplus \{e_7,e_9,e_{10},e_{11}\} \oplus \{e_8,e_9,e_{11}\} = \{e_2,e_3,e_5,e_6,e_7,e_9,e_{10},e_{11}\}$;
$\tau_3(e_{10}) = \gamma\{e_1,e_3,e_6,e_8,e_9,e_{11}\} = \{e_2,e_3,e_5,e_6\} \oplus \{e_1,e_2,e_5\} \oplus \{e_1,e_2,e_4,e_7\} \oplus \{e_7,e_9,e_{10},e_{11}\} \oplus$
$\oplus \{e_4,e_5,e_8,e_{10}\} \oplus \{e_7,e_8,e_{10}\} = \{e_2,e_3,e_5,e_6,e_7,e_9,e_{10},e_{11}\}$;
$\tau_3(e_{11}) = \gamma\{e_1,e_2,e_5,e_7,e_8,e_{10}\} = \{e_2,e_3,e_5,e_6\} \oplus \{e_1,e_3,e_6\} \oplus \{e_1,e_3,e_4,e_9\} \oplus \{e_4,e_6,e_8,e_{11}\} \oplus$
$\oplus \{e_7,e_9,e_{10},e_{11}\} \oplus \{e_8,e_9,e_{11}\} = \{e_2,e_3,e_5,e_6,e_7,e_9,e_{10},e_{11}\}$;

Все суграфы 4-го уровня пусты.

С целью определения веса ребра построим кортеж количества суграфов для строки матрицы T с участием ребра $e_i$

$\varepsilon(\tau(e_1)) = <0,2,2,0,3,3,2,0,2,1,1>$;

$\varepsilon(\tau(e_2)) = <2,2,4,1,1,4,2,1,2,1,2>$;

$\varepsilon(\tau(e_3)) = <2,4,2,1,4,1,2,1,2,2,1>$;

$\varepsilon(\tau(e_4)) = <0,1,1,0,2,2,2,0,2,1,1>$;

$\varepsilon(\tau(e_5)) = <3,1,4,2,1,2,2,2,4,2,2>$;

$\varepsilon(\tau(e_6)) = <3,4,1,2,2,1,4,2,2,2,2>$;

$\varepsilon(\tau(e_7)) = <2,2,2,2,2,4,1,3,2,1,4>$;

$\varepsilon(\tau(e_8)) = <0,1,1,0,2,2,3,0,3,2,2>$;

$\varepsilon(\tau(e_9)) = <2,2,2,2,4,2,2,3,1,4,1>$;

$\varepsilon(\tau(e_{10})) = <1,1,2,1,2,2,1,2,4,2,4>$;

$\varepsilon(\tau(e_{11})) = <1,2,1,1,2,2,4,2,1,4,2>$.

Каждому столбцу матрицы $w(l_j)$ поставим в соответствие кортеж $\xi(\tau(l_j))$, элементы которого (веса ребер) характеризуют количество подмножеств уровня с участием ребра $e_i$.

$\xi(\tau(l_1)) = <4,3,3,4,4,4,4,4,4,3,3>$;

$\xi(\tau(l_2)) = <4,5,5,8,7,7,7,4,7,5,5>$;

$\xi(\tau(l_3)) = <8,6,6,0,6,6,6,8,6,6,6>$;

$\xi(\tau(l_4)) = <0,8,8,0,8,8,8,0,8,8,8>$.

Построим суммарные кортежи графа $G_2$ как сумму элементов кортежей $\varepsilon_\tau(\mathrm{G}) = \sum_{i=1}^{m} \varepsilon(w(e_i))$ или $\xi_\tau(\mathrm{G}) = \sum_{j=1}^{k} \xi(\tau(l_j))$:

$\varepsilon_\tau(\mathrm{G}) = <\varepsilon(e_1), \varepsilon(e_2), \varepsilon(e_3), \ldots, \varepsilon(e_m)>$;

$\xi_\tau(\mathrm{G}) = <\xi(l_1), \xi(l_2), \xi(l_3), \ldots, \xi(l_k)>$.

Из построения следует:

$$\varepsilon_\tau(\mathrm{G}) = \xi_\tau(\mathrm{G}) \qquad (4.4)$$

Для графа $G_2$ кортеж $\xi_\tau(G_2) = <16,22,22,12,25,25,25,16,25,22,22>$.

Теперь можно задать вес каждой вершины в спектре реберных циклов как сумму весов инцидентных ребер.

Например, для графа $G_2$:

$\zeta(v_1) = \xi(e_1) + \xi(e_2) + \xi(e_3) \rightarrow 16 + 22 + 22 = 60$;

$\zeta(v_2) = \xi(e_1) + \xi(e_4) + \xi(e_5) + \xi(e_6) \rightarrow 16 + 12 + 25 + 25 = 78$;

$\zeta(v_3) = \xi(e_4) + \xi(e_7) + \xi(e_8) + \xi(e_9) \rightarrow 16 + 12 + 25 + 25 = 78$;

$\zeta(v_4) = \xi(e_2) + \xi(e_5) + \xi(e_7) + \xi(e_{10}) \rightarrow 22 + 25 + 25 + 22 = 94$;

$\zeta(v_5) = \xi(e_8) + \xi(e_{10}) + \xi(e_{11}) \rightarrow 16 + 22 + 22 = 60$;

$\zeta(v_6) = \xi(e_3) + \xi(e_6) + \xi(e_9) + \xi(e_{11}) \rightarrow 22 + 25 + 25 + 22 = 94$.

Построим кортеж $\zeta_\tau(G_2)$ для вершин графа

$$\zeta_\tau(\mathrm{G}) = <\zeta_\tau(v_1), \zeta_\tau(v_2), \ldots, \zeta_\tau(v_n)>. \qquad (4.5)$$

где $\zeta_{\tau}(v_j)$ – общий вес соответствующей вершины, где $j = (1,2,\ldots,n)$.

$\zeta_{\tau}(G_1)$ = <60,78,78,94,60,94>.

Построим матрицу базовых реберных циклов графа $G_2$.

T =

| | $e_1$ | $e_2$ | $e_3$ | $e_4$ | $e_5$ | $e_6$ | $e_7$ | $e_8$ | $e_9$ | $e_{10}$ | $e_{11}$ |
|---|---|---|---|---|---|---|---|---|---|---|---|
| $e_1$ | | 1 | 1 | | 1 | 1 | | | | | |
| $e_2$ | 1 | | 1 | | | 1 | | | | | |
| $e_3$ | 1 | 1 | | | 1 | | | | | | |
| $e_4$ | | | | | 1 | 1 | 1 | | 1 | | |
| $e_5$ | 1 | | 1 | 1 | | | | | 1 | | |
| $e_6$ | 1 | 1 | | 1 | | | 1 | | | | |
| $e_7$ | | | | 1 | | 1 | | 1 | | | 1 |
| $e_8$ | | | | | | | 1 | | 1 | 1 | 1 |
| $e_9$ | | | | 1 | 1 | | | 1 | | 1 | |
| $e_{10}$ | | | | | | | | 1 | 1 | | 1 |
| $e_{11}$ | | | | | | | 1 | 1 | | 1 | |

Как видно, матрица базовых реберных циклов является симметричной относительно главной диагонали. С другой стороны, она является нильпотентным оператором в пространстве суграфов графа. Точно так же является линейным оператором и матрица базовых разрезов ребер в пространстве суграфов графа.

$T^2(e_i)$ =

| | $e_1$ | $e_2$ | $e_3$ | $e_4$ | $e_5$ | $e_6$ | $e_7$ | $e_8$ | $e_9$ | $e_{10}$ | $e_{11}$ |
|---|---|---|---|---|---|---|---|---|---|---|---|
| $e_1$ | | | | | 1 | 1 | 1 | | 1 | | |
| $e_2$ | | 1 | 1 | 1 | | 1 | 1 | | | | |
| $e_3$ | | 1 | 1 | 1 | 1 | | | | 1 | | |
| $e_4$ | | 1 | 1 | | 1 | 1 | 1 | | 1 | 1 | 1 |
| $e_5$ | 1 | | 1 | 1 | | | 1 | 1 | 1 | 1 | |
| $e_6$ | 1 | 1 | | 1 | | | 1 | 1 | 1 | | 1 |
| $e_7$ | 1 | 1 | | 1 | 1 | 1 | | 1 | | | 1 |
| $e_8$ | | | | | 1 | 1 | 1 | | 1 | | |
| $e_9$ | 1 | | 1 | 1 | 1 | 1 | | 1 | | 1 | |
| $e_{10}$ | | | | 1 | 1 | | | | 1 | 1 | 1 |
| $e_{11}$ | | | | 1 | | 1 | 1 | | | 1 | 1 |

$T^3(e_i)$ =

| | $e_1$ | $e_2$ | $e_3$ | $e_4$ | $e_5$ | $e_6$ | $e_7$ | $e_8$ | $e_9$ | $e_{10}$ | $e_{11}$ |
|---|---|---|---|---|---|---|---|---|---|---|---|
| $e_1$ | | 1 | 1 | | 1 | 1 | 1 | | 1 | 1 | 1 |
| $e_2$ | 1 | | 1 | | | 1 | | 1 | 1 | | 1 |
| $e_3$ | 1 | 1 | | | 1 | | 1 | 1 | | 1 | |
| $e_4$ | | | | | | | | | | | |
| $e_5$ | 1 | | 1 | | | 1 | | 1 | 1 | | 1 |
| $e_6$ | 1 | 1 | | | 1 | | 1 | 1 | | 1 | |
| $e_7$ | 1 | | 1 | | | 1 | | 1 | 1 | | 1 |
| $e_8$ | | 1 | 1 | | 1 | 1 | 1 | | 1 | 1 | 1 |
| $e_9$ | 1 | 1 | | | 1 | | 1 | 1 | | 1 | |
| $e_{10}$ | 1 | | 1 | | | 1 | | 1 | 1 | | 1 |
| $e_{11}$ | 1 | 1 | | | 1 | | 1 | 1 | | 1 | |

| | | $e_1$ | $e_2$ | $e_3$ | $e_4$ | $e_5$ | $e_6$ | $e_7$ | $e_8$ | $e_9$ | $e_{10}$ | $e_{11}$ |
|---|---|---|---|---|---|---|---|---|---|---|---|---|
| | $e_1$ | | | | | | | | | | | |
| | $e_2$ | | 1 | 1 | | 1 | 1 | 1 | | 1 | 1 | 1 |
| | $e_3$ | | 1 | 1 | | 1 | 1 | 1 | | 1 | 1 | 1 |
| | $e_4$ | | | | | | | | | | | |
| | $e_5$ | | 1 | 1 | | 1 | 1 | 1 | | 1 | 1 | 1 |
| $T^4(e_i)$ = | $e_6$ | | 1 | 1 | | 1 | 1 | 1 | | 1 | 1 | 1 |
| | $e_7$ | | 1 | 1 | | 1 | 1 | 1 | | 1 | 1 | 1 |
| | $e_8$ | | | | | | | 1 | | 1 | 1 | 1 |
| | $e_9$ | | 1 | 1 | | 1 | 1 | 1 | | 1 | 1 | 1 |
| | $e_{10}$ | | 1 | 1 | | 1 | 1 | 1 | | 1 | 1 | 1 |
| | $e_{11}$ | | 1 | 1 | | 1 | 1 | 1 | | 1 | 1 | 1 |

Все элементы матрицы $T^5(e_i) = \varnothing$ .

## 4.3. Спектр реберных циклов и реберных разрезов графа

Проведем сравнительный анализ спектра реберных циклов и реберных разрезов графа.

Таблица 4.1. Сравнительный анализ реберных разрезов и реберных циклов

| Реберные разрезы графа $G_2$ | Реберные циклы графа $G_2$ |
|---|---|
| 1 | 2 |
| Множество базовых реберных разрезов: | Множество базовых реберных циклов: |
| $w_0(e_1) =\{e_2,e_3,e_4,e_5,e_6\}$; | $\tau_0(e_1) = \{e_2,e_3,e_5,e_6\}$; |
| $w_0(e_2) = \{e_1,e_3,e_5,e_7,e_{10}\}$; | $\tau_0(e_2) = \{e_1,e_3,e_6\}$; |
| $w_0(e_3) = \{e_1,e_2,e_6,e_9,e_{11}\}$; | $\tau_0(e_3) = \{e_1,e_2,e_5\}$; |
| $w_0(e_4) = \{e_1,e_5,e_6,e_7,e_8,e_9\}$; | $\tau_0(e_4) = \{e_5,e_6,e_7,e_9\}$; |
| $w_0(e_5) = \{e_1,e_2,e_4,e_6,e_7,e_{10}\}$; | $\tau_0(e_5) = \{e_1,e_3,e_4,e_9\}$; |
| $w_0(e_6) = \{e_1,e_3,e_4,e_5,e_9,e_{11}\}$; | $\tau_0(e_6) = \{e_1,e_2,e_4,e_7\}$; |
| $w_0(e_7) = \{e_2,e_4,e_5,e_8,e_9,e_{10}\}$; | $\tau_0(e_7) = \{e_4,e_6,e_8,e_{11}\}$; |
| $w_0(e_8) = \{e_4,e_7,e_9,e_{10},e_{11}\}$; | $\tau_0(e_8) = \{e_7,e_9,e_{10},e_{11}\}$; |
| $w_0(e_9) = \{e_3,e_4,e_6,e_7,e_8,e_{11}\}$; | $\tau_0(e_9) = \{e_4,e_5,e_8,e_{10}\}$; |
| $w_0(e_{10}) = \{e_2,e_5,e_7,e_8,e_{11}\}$; | $\tau_0(e_{10}) = \{e_8,e_9,e_{11}\}$; |
| Количество суграфов для строки спектра разрезов $W_s$ с участием ребра $e_i$: | Количество суграфов для строки спектра циклов $T_c$ с участием ребра $e_i$: |
| $\varepsilon(w(e_1)) = \langle 1,2,2,1,3,3,2,1,2,1,1\rangle$; | $\varepsilon(\tau(e_1)) = \langle 0,2,2,0,3,3,2,0,2,1,1\rangle$; |
| $\varepsilon(w(e_2)) = \langle 2,2,3,1,4,3,3,1,2,2,1\rangle$; | $\varepsilon(\tau(e_2)) = \langle 2,2,4,1,1,4,2,1,2,1,2\rangle$; |
| $\varepsilon(w(e_3)) = \langle 2,3,2,1,3,4,2,1,3,1,2\rangle$; | $\varepsilon(\tau(e_3)) = \langle 2,4,2,1,4,1,2,1,2,2,1\rangle$; |
| $\varepsilon(w(e_4)) = \langle 1,1,1,0,2,2,2,1,2,1,1\rangle$; | $\varepsilon(\tau(e_4)) = \langle 0,1,1,0,2,2,2,0,2,1,1\rangle$; |
| $\varepsilon(w(e_5)) = \langle 3,4,3,2,1,2,3,2,2,3,2\rangle$; | $\varepsilon(\tau(e_5)) = \langle 3,1,4,2,1,2,2,2,4,2,2\rangle$; |
| $\varepsilon(w(e_6)) = \langle 3,3,4,2,2,1,2,2,3,2,3\rangle$; | $\varepsilon(\tau(e_6)) = \langle 3,4,1,2,2,1,4,2,2,2,2\rangle$; |
| $\varepsilon(w(e_7)) = \langle 2,3,2,2,3,2,1,3,2,4,3\rangle$; | $\varepsilon(\tau(e_7)) = \langle 2,2,2,2,2,4,1,3,2,1,4\rangle$; |
| $\varepsilon(w(e_8)) = \langle 1,1,1,1,2,2,3,1,3,2,2\rangle$; | $\varepsilon(\tau(e_8)) = \langle 0,1,1,0,2,2,3,0,3,2,2\rangle$; |
| $\varepsilon(w(e_9)) = \langle 2,2,3,2,2,3,2,3,1,3,4\rangle$; | $\varepsilon(\tau(e_9)) = \langle 2,2,2,2,4,2,2,3,1,4,1\rangle$; |
| $\varepsilon(w(e_{10})) = \langle 1,2,1,1,3,2,4,2,3,2,3\rangle$; | $\varepsilon(\tau(e_{10})) = \langle 1,1,2,1,2,2,1,2,4,2,4\rangle$; |
| $\varepsilon(w(e_{11})) = \langle 1,1,2,1,2,3,3,2,4,3,2\rangle$. | $\varepsilon(\tau(e_{11})) = \langle 1,2,1,1,2,2,4,2,1,4,2\rangle$. |
| Количество суграфов уровней с участием ребра $e_i$.<br>$\zeta(\mathrm{w}(l_1))$ = <5,5,5,6,6,6,6,5,6,5,5>;<br>$\zeta(\mathrm{w}(l_2))$ = <6,5,5,8,7,7,7,6,7,5,5>;<br>$\zeta(\mathrm{w}(l_3))$ = <8,6,6,0,6,6,6,8,6,6,6>;<br>$\zeta(\mathrm{w}(l_4))$ = <0,8,8,0,8,8,8,0,8,8,8>. | Количество суграфов уровней с участием ребра $e_i$.<br>$\zeta(\tau(l_1))$ = <4,3,3,4,4,4,4,4,4,3,3>;<br>$\zeta(\tau(l_2))$ = <4,5,5,8,7,7,7,4,7,5,5>;<br>$\zeta(\tau(l_3))$ = <8,6,6,0,6,6,6,8,6,6,6>;<br>$\zeta(\tau(l_4))$ = <0,8,8,0,8,8,8,0,8,8,8>. |
| Суммарный реберный кортеж весов:<br>$\xi_{\mathrm{w}}$(G) = <19,24,24,14,27,27,27,19,27,24,24><br>Вершинный кортеж весов:<br>$\zeta_{\mathrm{w}}$(G) = <67,87,87,102,67,102>. | Суммарный реберный кортеж весов:<br>$\xi_\tau$(G) = <16,22,22,12,25,25,25,16,25,22,22>.<br>Вершинный кортеж весов:<br>$\zeta_\tau$(G) = <60,78,78,94,60,94>. |

Как видим, есть определенная синхронность в весах ребер и вершин. Кроме того, процессы формирования спектра реберных циклов аналогичны процессу формирования спектра реберных разрезов графа со всеми ограничениями.

Таким образом, процесс порождения реберных циклов также можно ограничить следующими факторами: $\mathrm{T}_\lambda^{\mathrm{q}} = \varnothing$ если $\mathrm{T}_\lambda^{\mathrm{q}} = \mathrm{T}_\lambda$, или если $\mathrm{T}_\lambda^{\mathrm{q}}$ состоит из пустых строк.

Как будет показано дальше, для задачи определения изоморфизма графов можно ограничиться участием только базовых реберных циклов. Тогда вычислительную сложность алгоритма построения спектра реберных циклов можно определить, построив множество изометрических циклов графа. Вычислительная сложность алгоритма построения множества изометрических циклов равна O($n^4$). Следующий шаг – формирование суграфа, представляющегося собой кольцевую сумму циклов, проходящих по ребру $e_i$. Вычислительную сложность данного процесса можно определить как вычислительную сложность построения множества изометрических циклов, то есть O($n^4$). Количество реберных циклов можно определить как количество ребер умноженное на количество циклов проходящих по ребру (в предположении, что такие циклы длиной три,), то их количество равно n-2. Получаем, что количество реберных циклов равно $m \times (n\text{-}2)$. Отсюда, вычислительная сложность для определения реберных циклов равна количеству уровней в спектре реберных циклов умноженному $q \times m \times (n\text{-}2)$. Тогда вычислительная сложность построения матрицы $T_C$(G), в предположении, что количество уровней равно единице, можно определить как O($n^4$)+O($n^3$). Окончательно вычислительную сложность алгоритма построения инварианта спектра реберных циклов графа, можно определить относительно вершин графа как O($n^4$).

Таким образом, вычислительная сложность алгоритма построения инварианта спектра реберных циклов, равна O($n^4$). Задача определения инварианта спектра реберных циклов графа, относится к классу P - полиномиальных алгоритмов.

Для графа G инвариант спектра реберных циклов имеет вид аналогичный базовому инварианту спектра реберных разрезов: $\mathrm{F}_{\tau}(\xi(\mathrm{G}))\,\&\,\mathrm{F}_{\tau}(\zeta(\mathrm{G}))$ где $\mathrm{F}_{\tau}(\xi(\mathrm{G}))$ – упорядоченный вектор весов ребер для спектра реберных циклов графа, а $\mathrm{F}_{\tau}(\zeta(\mathrm{G}))$ – упорядоченный вектор весов вершин для спектра реберных циклов графа.

Будем считать достаточным для различения графов, применение только множества базовых реберных циклов.

В качестве примера, рассмотрим цифровой состав инварианта реберных циклов некоторых графов с 12 вершинами и 30 ребрами.

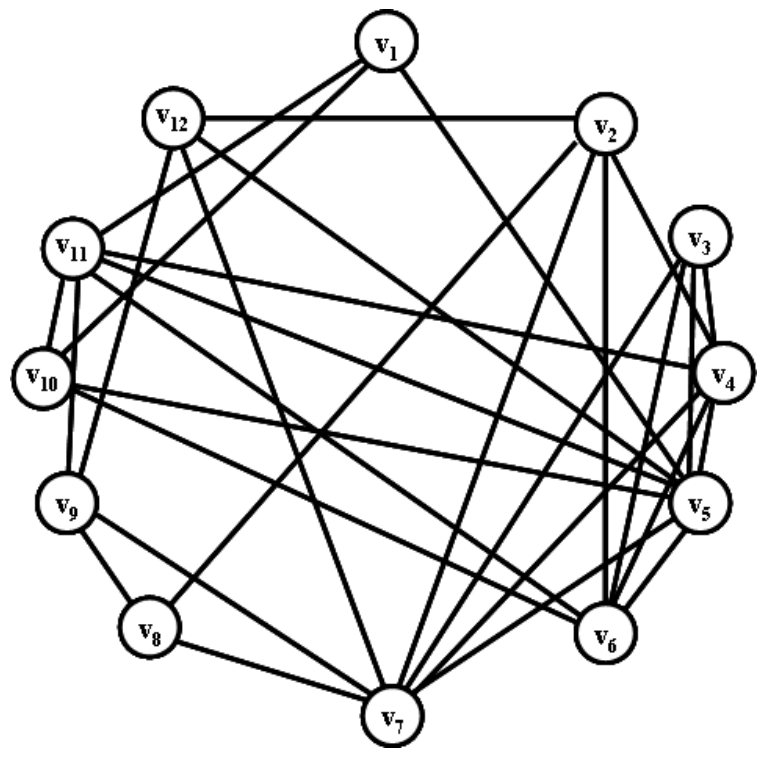

Рис. 4.3. Граф $G_9$.

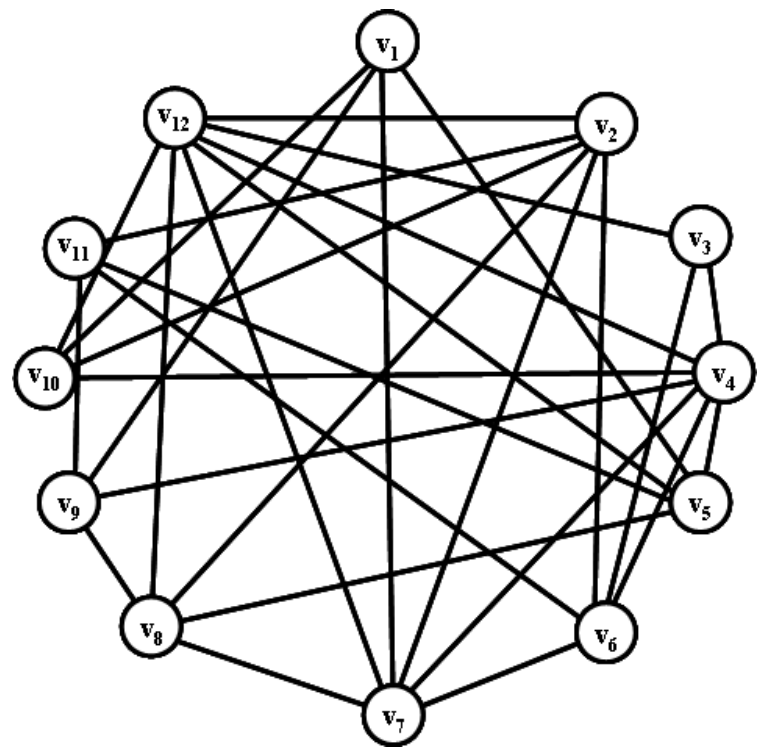

Рис. 4.4.. Граф $G_{10}$.

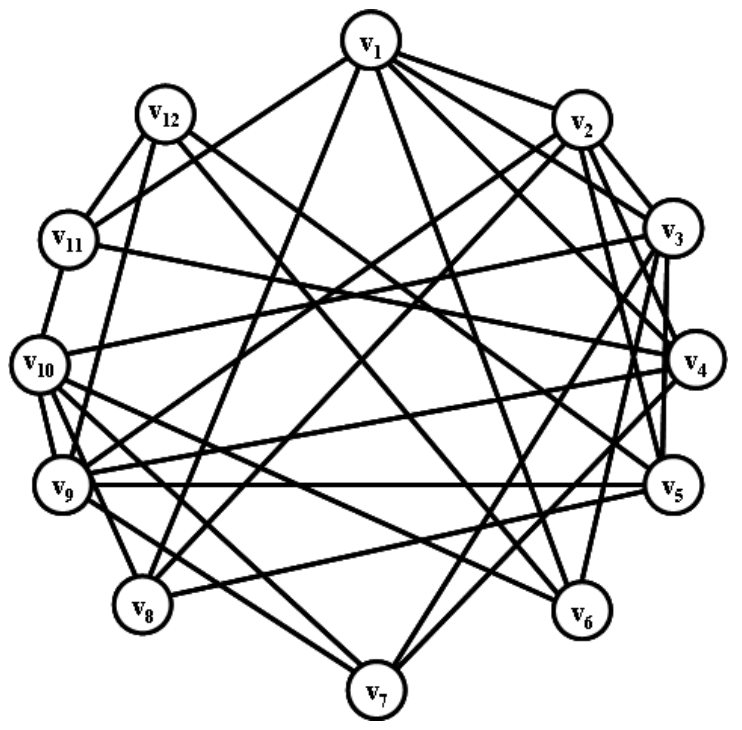

Рис. 4.5. Граф $G_{11}$.

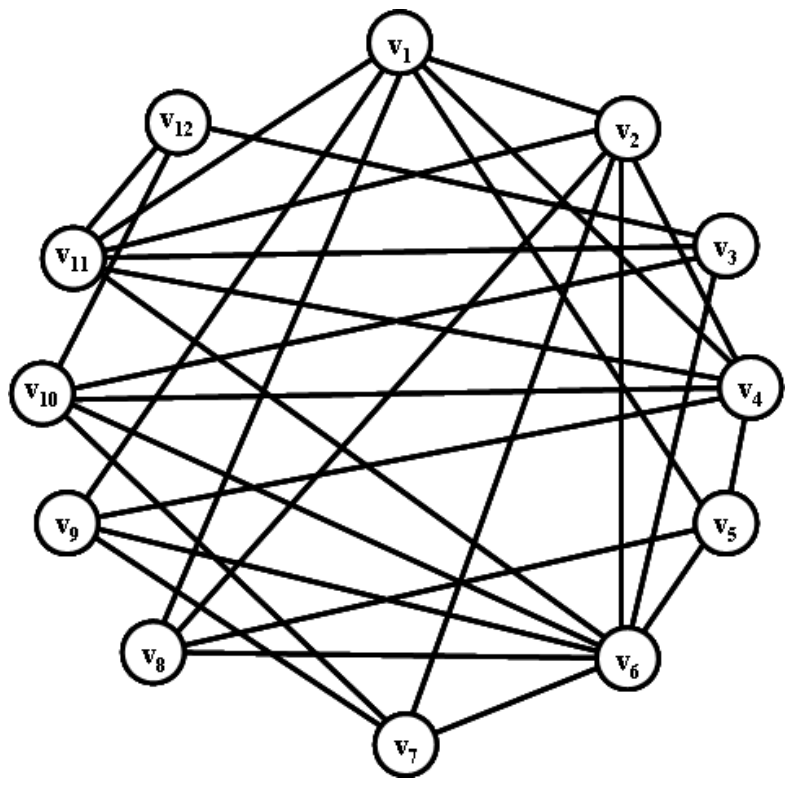

Рис. 4.6. Граф $G_{12}$.

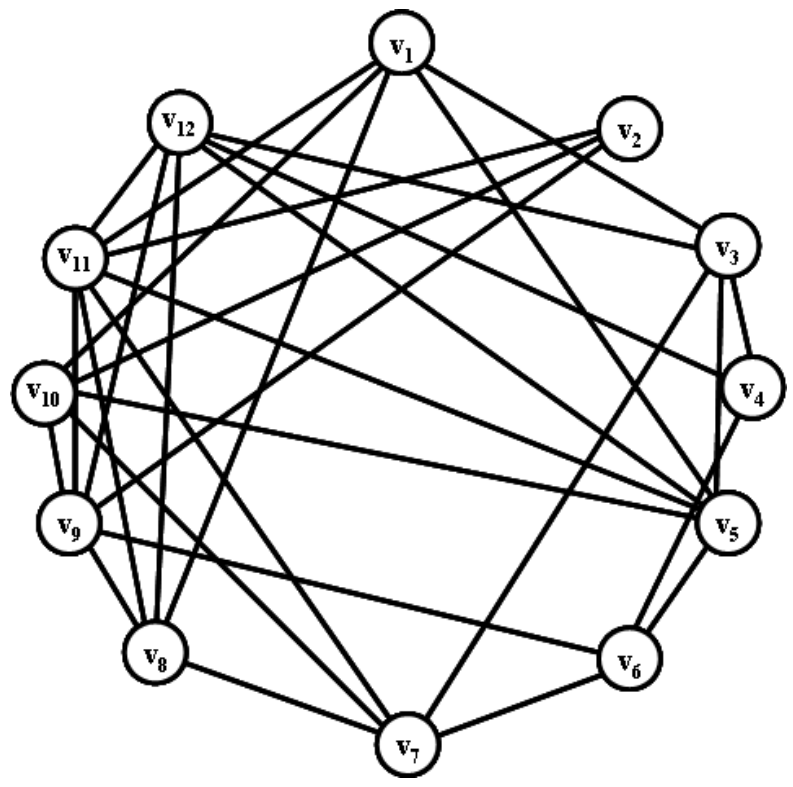

Рис. 4.7. Граф $G_{13}$.

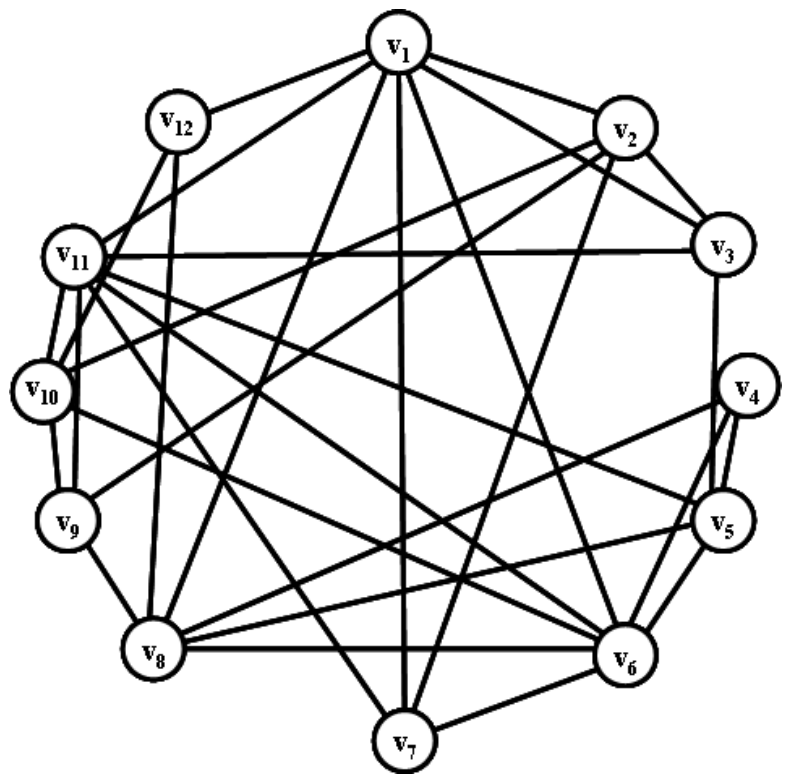

Рис. 4.8. Граф $G_{14}$.

Количество изометрических циклов в графе $G_9$ = 32. Векторный инвариант спектра реберных циклов графа $G_9$:

$\mathrm{IC}(\mathrm{G}_9) = \mathrm{F}_{\tau}(\xi(\mathrm{G}_9))\,\&\,\mathrm{F}_{\tau}(\zeta(\mathrm{G}_9))$ =
= (4,4,4,4,4,8,8,9,9,9,9,10,10,11,12,12,12,13,13,14,15,15,16,16,16,17,18,19,19,20) &
& (12,38,46,46,50,52,60,68,74,74,82,98).

Количество изометрических циклов в графе $G_{10}$ = 44. Векторный инвариант спектра реберных циклов графа $G_{10}$:

$\mathrm{IC}(\mathrm{G}_{10}) = \mathrm{F}_{\tau}(\xi(\mathrm{G}_{10}))\,\&\,\mathrm{F}_{\tau}(\zeta(\mathrm{G}_{10}))$ =
= (4,7,7,8,8,8,8,8,9,10,10,10,11,11,11,11,12,12,12,12,13,13,13,13,14,14,14,15,18,18) &

& (26,36,42,44,50,56,56,58,66,68,80,86).

Количество изометрических циклов в графе $G_{11}$ = 40. Векторный инвариант спектра реберных циклов графа $G_{11}$:

$$IC(G_{11}) = F_{\tau}(\xi(G_{11})) \,\&\, F_{\tau}(\zeta(G_{11})) =$$
$$= (7,7,7,7,8,8,8,8,8,8,9,9,9,10,10,10,10,10,11,11,11,11,12,12,12,13,13,14,14,15)\ \&$$
$$\&\ (32,38,38,40,42,44,52,56,58,60,66,78).$$

Количество изометрических циклов в графе $G_{12}$ = 42. Векторный инвариант спектра реберных циклов графа $G_{12}$:

$$IC(G_{12}) = F_{\tau}(\xi(G_{12})) \,\&\, F_{\tau}(\zeta(G_{12})) =$$
$$(4,4,7,7,8,8,8,10,11,11,11,11,12,12,12,13,13,13,14,14,14,15,15,15,16,17,17,18,18,20)\ \&$$
$$\&\ (30,36,44,52,54,58,66,70,72,78,82,94).$$

Количество изометрических циклов в графе $G_{13}$ = 49. Векторный инвариант спектра реберных циклов графа $G_{13}$:

$$IC(G_{13}) = F_{\tau}(\xi(G_{13})) \,\&\, F_{\tau}(\zeta(G_{13})) =$$
$$= (4,6,8,8,9,9,9,9,9,10,10,10,11,11,12,12,12,13,13,14,14,14,14,15,15,16,16,17,19)\ \&$$
$$\&\ (22,28,46,48,54,60,66,66,72,72,74,88).$$

Количество изометрических циклов в графе $G_{14}$ = 40. Векторный инвариант спектра реберных циклов графа $G_{14}$:

$$IC(G_{14}) = F_{\tau}(\xi(G_{14})) \,\&\, F_{\tau}(\zeta(G_{14})) =$$
$$= (4,4,4,8,8,8,8,9,9,9,9,9,10,10,10,11,11,12,12,13,13,13,13,13,15,16,18,18,18,19)\ \&$$
$$\&\ (12,38,44,44,46,54,60,60,72,74,80,84).$$

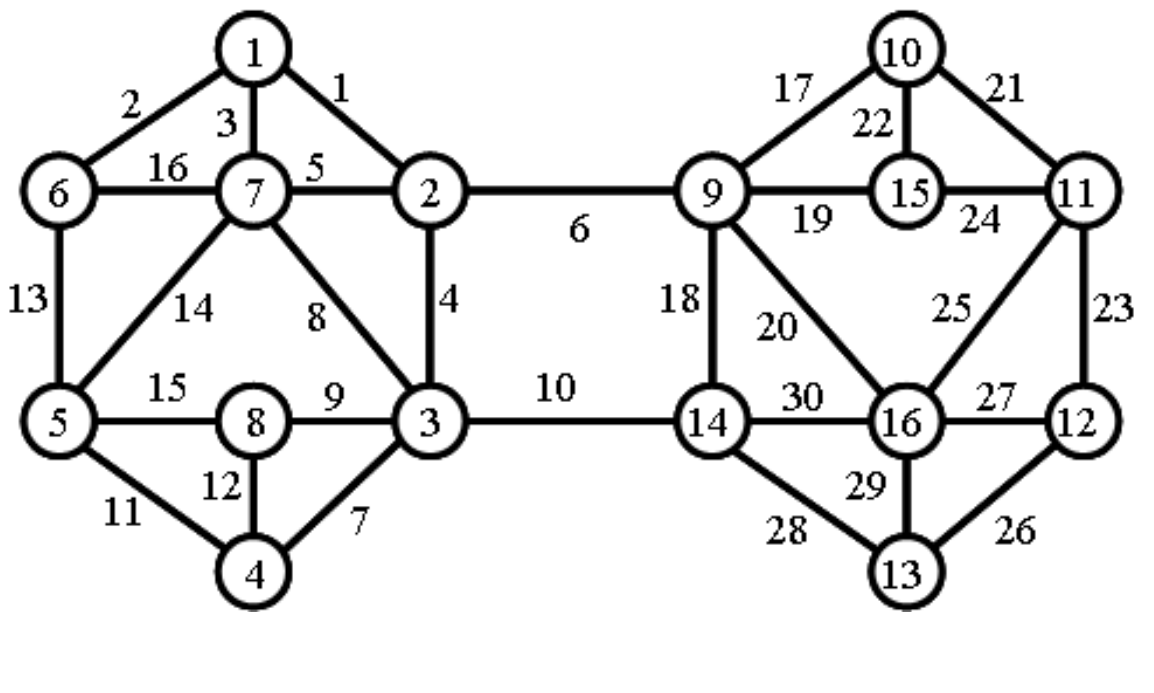

Рис. 4.9. Граф $G_{15}$.

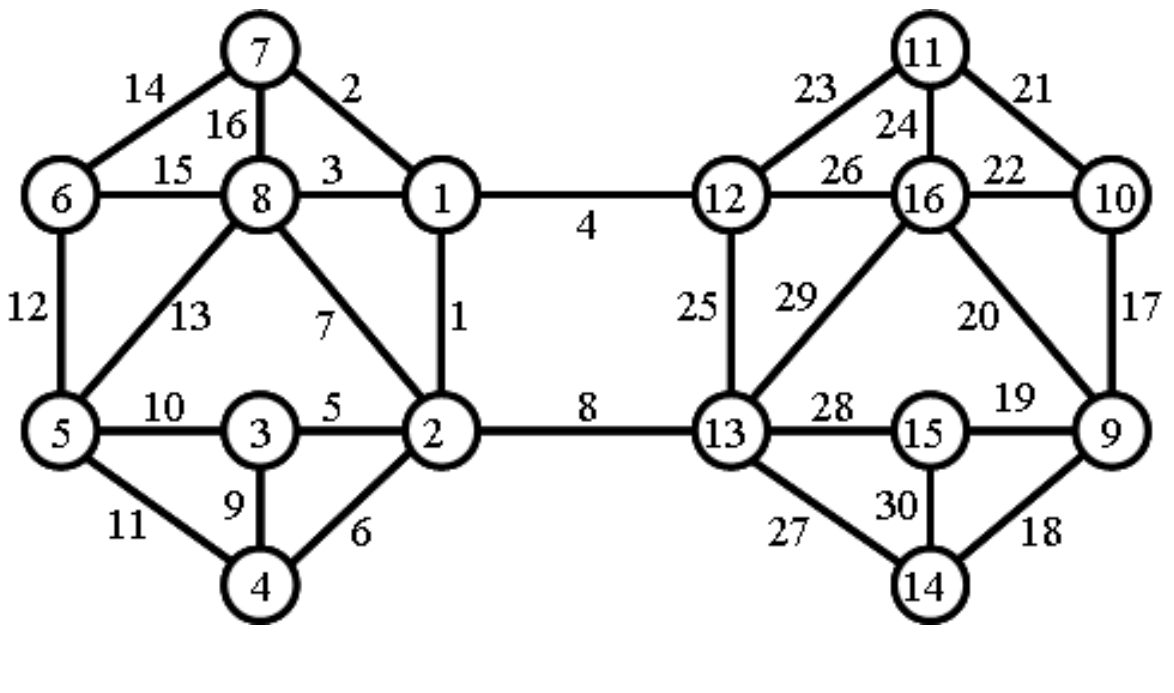

Рис. 4.10. Граф $G_{16}$.

Инвариант реберных циклов не может служить полным инвариантом для распознавания изоморфизма графов. В качестве примера можно привести два неизоморфных графа $G_{15}$ и $G_{16}$, имеющих один и тот же инвариант реберных циклов.

Количество изометрических циклов в графе $G_{15}$ = 17. Векторный инвариант спектра реберных циклов графа $G_{15}$:

$$IC(G_{15}) = F_{\tau}(\xi(G_{15})) \,\&\, F_{\tau}(\zeta(G_{15})) =$$
$$= (4,4,4,4,4,4,4,4,5,5,5,5,5,5,5,5,5,5,11,11,12,12,13,13,13,13,15,15,17,17)\ \&$$

& (14,14,14,14,28,28,30,30,34,34,36,36,44,44,44,44).

Количество изометрических циклов в графе $G_{16}$ = 17. Векторный инвариант спектра реберных циклов графа $G_{16}$:

$IC(G_{16}) = F_{\tau}(\xi(G_{16})) \& F_{\tau}(\zeta(G_{16}))$ =

= (4,4,4,4,4,4,4,4,5,5,5,5,5,5,5,5,5,5,11,11,12,12,13,13,13,13,15,15,17,17) &

& (14,14,14,14,28,28,30,30,34,34,36,36,44,44,44,44).

## Комментарии

На основе множества изометрических циклов графа представлен метод формирования реберного цикла относительно ребра графа. Базовые реберные циклы являются основой для построения спектра реберных циклов графа, что позволяет вводить числовые характеристики элементов графа. Создание векторных инвариантов спектра реберных циклов определяет возможность проводить сравнительный анализ структур графа на основе спектров реберных циклов и спектра реберных разрезов графов.

Вычислительная сложность построения инварианта для спектра реберных циклов графа определяется как $O(n^4)$.

# Глава 5. РЕБЕРНЫЕ ГРАФЫ

## 5.1. Структура реберного графа

Для распознавания изоморфизма графов существует теорема Уитни [32,51,52].

**Теорема Уитни**. Два графа G и H, с одинаковым количеством вершин и ребер, изоморфны тогда и только тогда, когда изоморфны их реберные графы L(G) и L(H) [].

Определим структуру реберного графа как отображение элементов графа G. Построение такой структуры, для облегчения восприятия, будем проводить на примере несепарабельного неориентированного графа G и его реберного графа L(G), параллельно используя более общие определения и обозначения.

Рассмотрим следующий граф $G_{17}$.

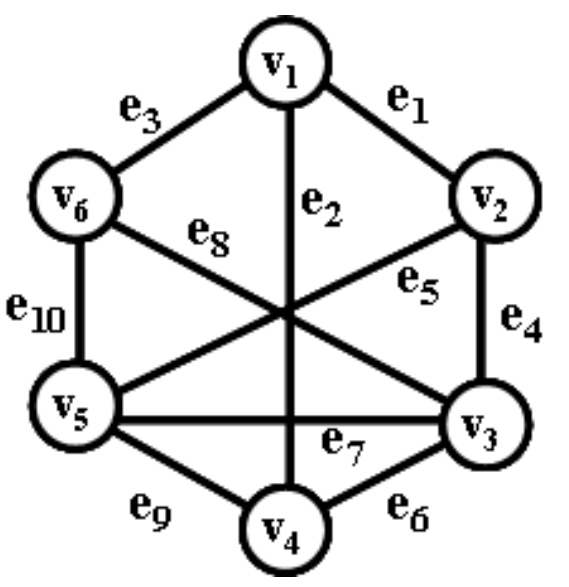

Рис. 5.1. Граф $G_{17}$.

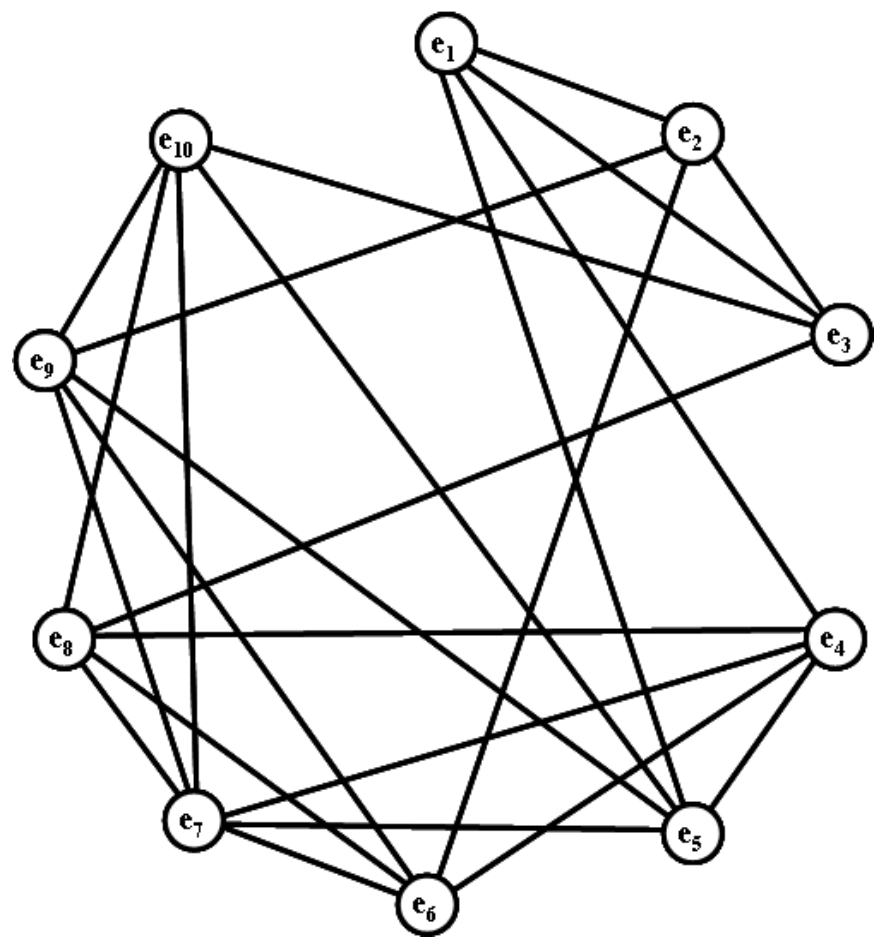

Рис. 5.2. Реберный граф $L(G_{17})$.

Выделим множество центральных разрезов $S(G_{17})$ графа $G_{17}$[20]:

$s_1 = \{e_1,e_2,e_3\}$;
$s_2 = \{e_1,e_4,e_5\}$;
$s_3 = \{e_4,e_6,e_7,e_8\}$;
$s_4 = \{e_2,e_6,e_9\}$;
$s_5 = \{e_5,e_7,e_9,e_{10}\}$;
$s_6 = \{e_3,e_8,e_{10}\}$.

Определим множество базовых реберных разрезов $W_0(G_{17})$ графа $G_{17}$:

$w_0(e_1) = s_1 \oplus s_2 = \{e_2,e_3,e_4,e_5\}$;
$w_0(e_2) = s_1 \oplus s_4 = \{e_1,e_3,e_6,e_9\}$;
$w_0(e_3) = s_1 \oplus s_6 = \{e_1,e_2,e_8,e_{10}\}$;
$w_0(e_4) = s_2 \oplus s_3 = \{e_1,e_5,e_6,e_7,e_8\}$;
$w_0(e_5) = s_2 \oplus s_5 = \{e_1,e_4,e_7,e_9,e_{10}\}$;
$w_0(e_6) = s_3 \oplus s_4= \{e_2,e_4,e_7,e_8,e_9\}$;
$w_0(e_7) = s_3 \oplus s_5 = \{e_4,e_5,e_6,e_8,e_9,e_{10}\}$;
$w_0(e_8) = s_3 \oplus s_6 = \{e_3,e_4,e_6,e_7,e_{10}\}$;
$w_0(e_9) = s_4 \oplus s_5 = \{e_2,e_5,e_6,e_7,e_{10}\}$;
$w_0(e_{10}) = s_5 \oplus s_6 = \{e_3,e_5,e_7,e_8,e_9\}$.

Множество изометрических циклов $C_\tau$ графа $G_{17}$:

$c_1 = \{e_1,e_2,e_4,e_6\}$;
$c_2 = \{e_1,e_2,e_5,e_9\}$;
$c_3 = \{e_1,e_3,e_4,e_8\}$;
$c_4 = \{e_1,e_3,e_5,e_{10}\}$;
$c_5 = \{e_2,e_3,e_6,e_8\}$;
$c_6 = \{e_2,e_3,e_9,e_{10}\}$;
$c_7 = \{e_4,e_5,e_7\}$;
$c_8 = \{e_6,e_7,e_9\}$;
$c_9 = \{e_7,e_9,e_{10}\}$.

Запишем множество изометрических циклов $C_\tau^L$ графа $L(G_{17})$ в виде вершин:

$c_1^L = \{v_1,v_2,v_3\}$;
$c_2^L = \{v_1,v_2,v_4,v_6\}$;
$c_3^L = \{v_1,v_2,v_5,v_9\}$;
$c_4^L = \{v_1,v_3,v_4,v_8\}$;
$c_5^L = \{v_1,v_3,v_5,v_{10}\}$;
$c_6^L = \{v_1,v_4,v_5\}$;
$c_7^L = \{v_2,v_3,v_6,v_8\}$;
$c_8^L = \{v_2,v_3,v_9,v_{10}\}$;
$c_9^L = \{v_2,v_6,v_9\}$;
$c_{10}^L = \{v_3,v_8,v_{10}\}$;
$c_{11}^L = \{v_4,v_5,v_6,v_9\}$;
$c_{12}^L = \{v_4,v_5,v_7\}$;
$c_{13}^L = \{v_4,v_5,v_8,v_{10}\}$;
$c_{14}^L = \{v_4,v_6,v_7\}$;
$c_{15}^L = \{v_4,v_6,v_8\}$;
$c_{16}^L = \{v_4,v_7,v_8\}$;
$c_{17}^L = \{v_5,v_7,v_9\}$;
$c_{18}^L = \{v_5,v_7,v_{10}\}$;
$c_{19}^L = \{v_5,v_9,v_{10}\}$;
$c_{20}^L = \{v_6,v_7,v_8\}$;
$c_{21}^L = \{v_6,v_7,v_9\}$;
$c_{22}^L = \{v_6,v_8,v_9,v_{10}\}$;
$c_{23}^L = \{v_7,v_8,v_{10}\}$;
$c_{24}^L = \{v_7,v_9,v_{10}\}$.

Поставим в соответствие вершинам графа $L(G_{17})$ ребра графа $G_{17}$:

$$\varphi : v_i \to e_i,\ i = (1,2,\ldots,m). \qquad (5.1)$$

Тогда центральным разрезам $S_0^L$ графа $L(G_{17})$ будут соответствовать базовые реберные разрезы графа $G_{17}$:

$\varphi : s_1^L \to w_0(e_1) = \{e_2,e_3,e_4,e_5\};$
$\varphi : s_2^L \to w_0(e_2) = \{e_1,e_3,e_6,e_9\};$
$\varphi : s_3^L \to w_0(e_3) = \{e_1,e_2,e_8,e_{10}\};$
$\varphi : s_4^L \to w_0(e_4) = \{e_1,e_5,e_6,e_7,e_8\};$
$\varphi : s_5^L \to w_0(e_5) = \{e_1,e_4,e_7,e_9,e_{10}\};$
$\varphi : s_6^L \to w_0(e_6) = \{e_2,e_4,e_7,e_8,e_9\};$
$\varphi : s_7^L \to w_0(e_7) = \{e_4,e_5,e_6,e_8,e_9,e_{10}\};$
$\varphi : s_8^L \to w_0(e_8) = \{e_3,e_4,e_6,e_7,e_{10}\};$
$\varphi : s_9^L \to w_0(e_9) = \{e_2,e_5,e_6,e_7,e_{10}\};$
$\varphi : s_{10}^L \to w_0(e_{10}) = \{e_3,e_5,e_7,e_8,e_9\}.$

Множество изометрических циклов $C_\tau^L$ графа $L(G_{17})$ можно построить как объединение следующих трех соответствующих подмножеств суграфов графа $G_{17}$.

1. Множество изометрических циклов соответствующих множеству центральных разрезов графа $G_{17}$, будем обозначать как $S_R^L$:

$\varphi : c_1^L \to s_1 = \{e_1,e_2,e_3\};$
$\varphi : c_6^L \to s_2 = \{e_1,e_4,e_5\};$
$\varphi : (c_{14}^L \cup c_{15}^L \cup c_{16}^L \cup c_{20}^L) \to s_3 = \{e_4,e_6,e_7,e_8\};$
$\varphi : (c_{17}^L \cup c_{18}^L \cup c_{19}^L \cup c_{24}^L) \to s_5 = \{e_5,e_7,e_9,e_{10}\};$
$\varphi : c_9^L \to s_4 = \{e_2,e_6,e_9\};$
$\varphi : c_{10}^L \to s_6 = \{e_3,e_8,e_{10}\}.$

2. Множество изометрических циклов графа $L(G_{17})$ соответствующих множеству изометрических циклов графа $G_{17}$, будем обозначать как $C_\tau^L$:

$\varphi : c_2^L \to c_1 = \{e_1,e_2,e_4,e_6\};$
$\varphi : c_3^L \to c_2 = \{e_1,e_2,e_5,e_9\};$
$\varphi : c_4^L \to c_3 = \{e_1,e_3,e_4,e_8\};$
$\varphi : c_5^L \to c_4 = \{e_1,e_3,e_5,e_{10}\};$
$\varphi : c_7^L \to c_5 = \{e_2,e_3,e_6,e_8\};$
$\varphi : c_8^L \to c_6 = \{e_2,e_3,e_9,e_{10}\};$
$\varphi : c_{12}^L \to c_7 = \{e_4,e_5,e_7\};$
$\varphi : c_{21}^L \to c_8 = \{e_6,e_7,e_9\};$
$\varphi : c_{23}^L \to c_9 = \{e_7,e_8,e_{10}\}.$

3. Множество изометрических циклов соответствующих множеству простых циклов графа $G_{17}$, будем обозначать как $C_d^L$:

$\varphi : c_{11}^L \to c_1 \oplus c_2 = c_7 \oplus c_8 = \{e_4,e_5,e_6,e_9\};$
$\varphi : c_{13}^L \to c_3 \oplus c_4 = c_7 \oplus c_9 = \{e_4,e_5,e_8,e_{10}\};$

$\varphi : c_{22}^{L} \rightarrow c_5 \oplus c_6 = c_8 \oplus c_9 = c_1 \oplus c_3 \oplus c_6 = \{e_6, e_8, e_9, e_{10}\}$.

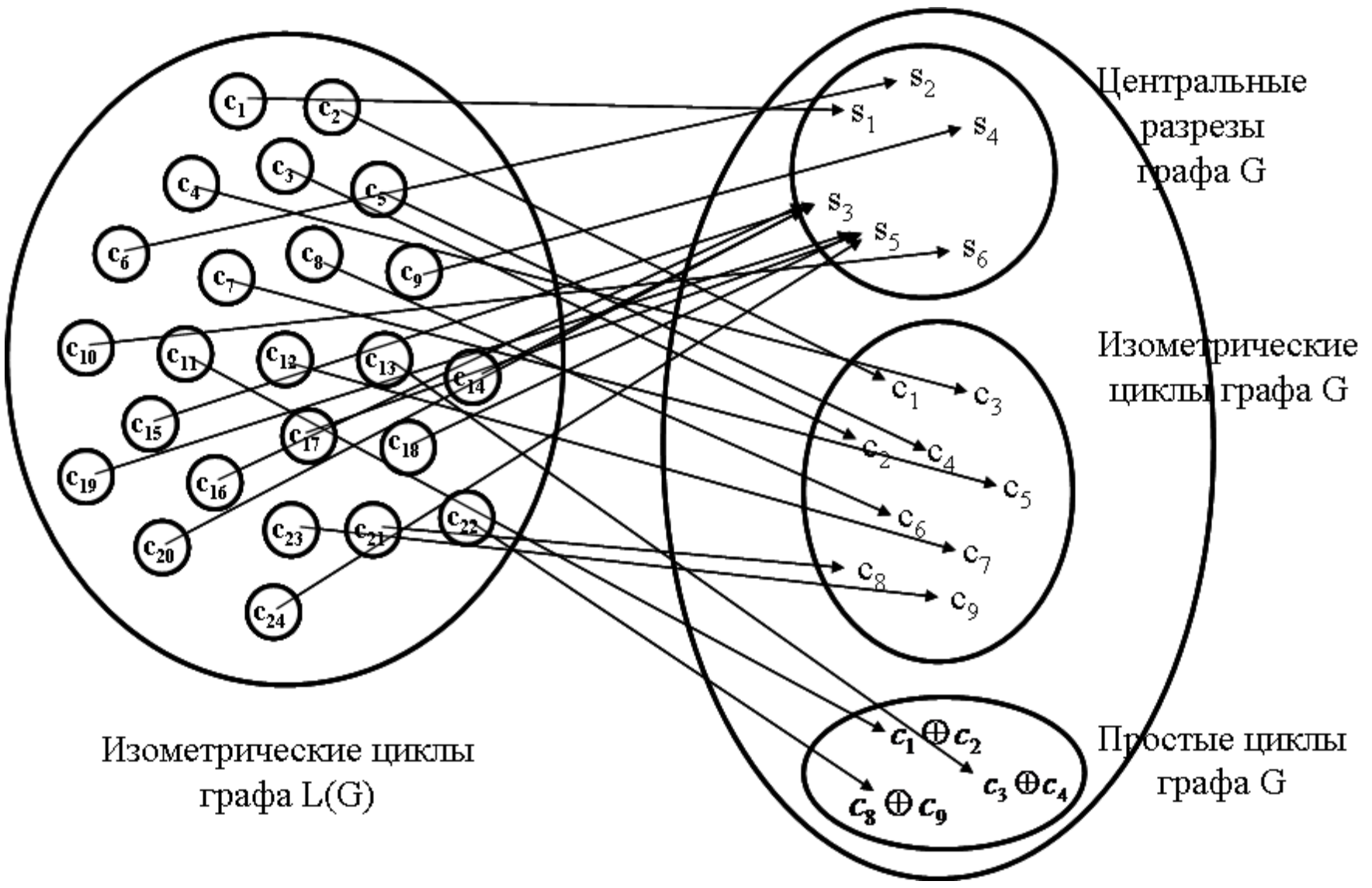

Рис. 5.3. Соответствия между изометрическими циклами графа L($G_{17}$) и суграфами графа $G_{17}$.

Каждый простой цикл является дубль-циклом. Каждый дубль-цикл характеризует систему зависимых изометрических циклов в $G_{17}$:

$c_1 \oplus c_2 \oplus c_7 \oplus c_8 = \varnothing$;
$c_3 \oplus c_4 \oplus c_7 \oplus c_9 = \varnothing$;
$c_8 \oplus c_9 \oplus c_1 \oplus c_3 \oplus c_6 = \varnothing$.

В свою очередь зависимая система циклов порождает топологический рисунок плоского суграфа. Например, система зависимых циклов $c_8 \oplus c_9 \oplus c_1 \oplus c_3 \oplus c_6$ порождает топологический рисунок плоского суграфа с двумя удаленными ребрами, представленный на рис. 5.4.

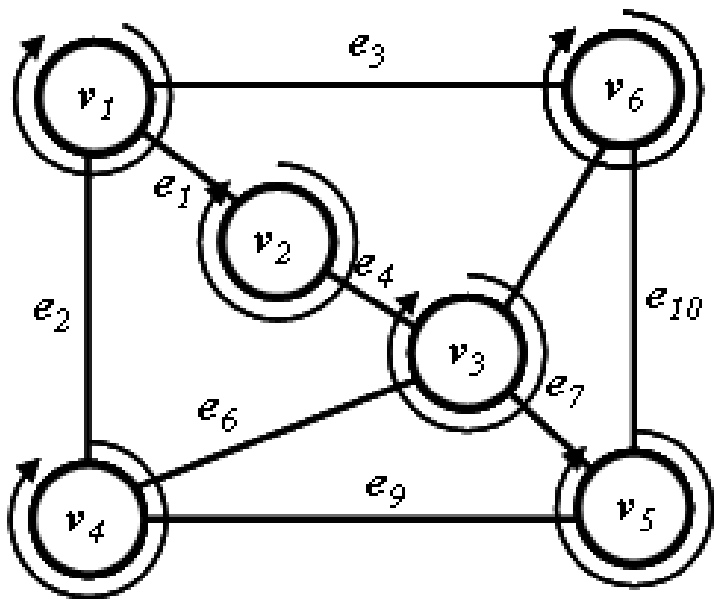

Рис. 5.4. Топологический рисунок плоского суграфа графа $G_{17}$.

**Определение 5.1.** Будем называть количество суграфов с участием ребра $e_i$ во всем множестве изометрических циклов реберного графа L(G) *весом ребра,* и обозначать $\xi_L(e_i)$. Построим следующий кортеж весов ребер i = 1,2,..,m.

Кортеж весов ребер : $\xi_L(e_i)$ = <6,6,6,9,9,9,9,9,9,9>

Это значит, что вес ребра $e_1$ = 6, вес ребра $e_5$ = 9 и т.д.

Определим вес вершины графа G как сумму весов инцидентных ребер. Построим следующий кортеж $\varsigma_L$ ($G_{17}$) для весов вершин j = 1,2,…,n.

Кортеж весов вершин: $\varsigma_L$ ($G_{17}$) = <18,24,36,24,36,24>.

Если расположить элементы кортежей по неубыванию, то можно получить инвариант графа, который является его индивидуальной характеристикой. Данный инвариант будем называть цифровым инвариантом реберного графа.

**Определение 5.2.** *Цифровым инвариантом реберного графа* называется числовая функция расположения весов ребер и вершин по не убыванию, согласно их кортежам.

IL(G) = F($\xi_L$(G)) & F($\zeta_L$(G))

В нашем примере:

IL($G_1$) = F($\xi_L$($G_{17}$)) & F($\zeta_L$($G_{17}$)) = (6,6,6,9,9,9,9,9,9,9) & (18,24,36,24,36,24) = $(3\times 6,7\times 9)$ & $(1\times 18,3\times 24,2\times 36)$.

Следует отметить, что цифровой инвариант реберного графа не зависят от нумерации вершин и ребер графа.

Таким образом, на основании теоремы Уитни можно утверждать следующее: если графы G и H с одинаковым количеством вершин и ребер имеют одинаковые цифровые инварианты реберных графов L(G) и L(H), то они изоморфны. Вычислительная сложность алгоритма определения цифрового инварианта реберного графа в основном определяется построением множества изометрических циклов реберного графа и равна $O(m^4)$ [23].

Структура реберного графа L($G_{17}$):

- количество вершин в графе L($G_{17}$) $n$ =10;
- количество ребер в графе L($G_{17}$) $m$ = 24;
- количество изометрических циклов реберного графа L($G_{17}$) = $card(C_\tau(L(G_{17})))$ = 24;
- количество изометрических циклов графа $G_{17}$ = $card(C_\tau)$ = 9;
- количество троек во множестве изоморфных циклов L($G_{17}$) $k_3$ = 10;
- количество четверок во множестве изоморфных циклов L($G_{17}$) $k_4$ = 3;
- цикломатическое число L($G_{17}$) = m-n+1=24-10+1 = 15;
- $card(C_\tau(L(G_{17})))$ = $card(C_\tau)$ +$k_3$+$k_4$ = 24 = 12+9+3.

Определим состав цифрового инварианта реберного графа L(G):

$$IL(G) = F(\xi_L(G)) \,\&\, F(\zeta_L(G)) \qquad (5.2)$$

- матрица инциденций графа B(G);
- $n$ – количество вершин графа G;
- $m$ – количество ребер графа G;

- $card(C_\tau)$ - количество изометрических циклов графа G;
- $W_0(G)$ – базовые реберные разрезы графа;
- $C_\tau^L$ – множество изометрических циклов реберного графа L(G);
- $\xi_L(G)$ – кортеж весов ребер;
- $\zeta_L(G)$ – кортеж весов вершин;
- IL(G) = F($\xi_L$(G)) & F($\zeta_L$(G)) – цифровой инвариант реберного графа L(G);

Опишем алгоритм для определения цифрового инварианта реберного графа.

**Инициализация.** Задан граф G.

**шаг 1. [Построение базовых реберных разрезов**]. В графе G выделяем множество базовых реберных разрезов графа.

**шаг 2. [Построение реберного графа]**. Строим реберный граф L(G)

**шаг 3. [Выделение множества изометрических циклов реберного графа]**. Формируем множество изометрических циклов реберного графа L(G).

**шаг 4. [Построение цифрового инварианта реберного графа]**. На основании выделенного множества изометрических циклов реберного графа строим цифровой инвариант реберного графа IL(G).

Конец работы алгоритма.

## 5.2. Суграфы реберного графа

Изометрические циклы реберного графа порождают три типа суграфов графа:

- суграфы, состоящие из трех ребер графа характеризующие центральные разрезы графа G (будем называть их тройками);
- суграфы, отображающие изометрические циклы графа G;
- суграфы, состоящие из четырех ребер графа G (будем называть их четверками).

Рассмотрим следующий граф $G_{18}$.

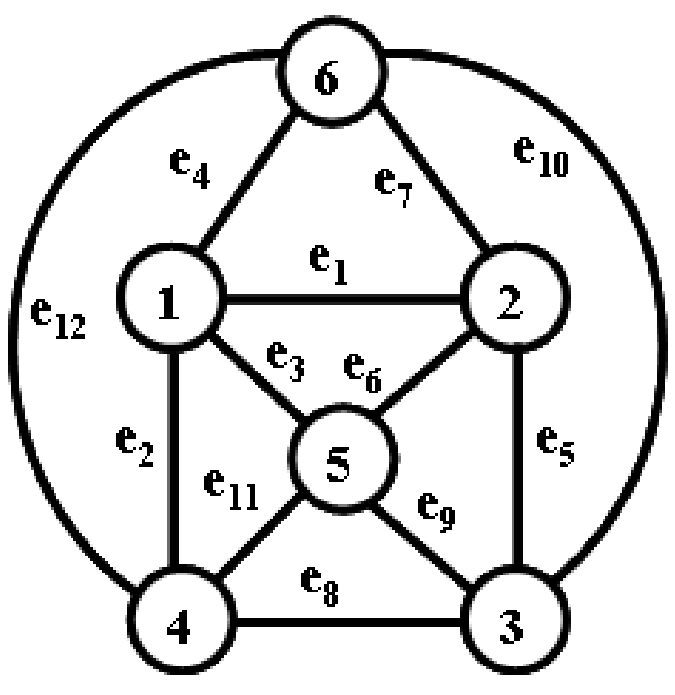

Рис. 5.5. Граф $G_{18}$.

Смежность графа $G_{18}$:

$v_1$: {$v_2,v_3,v_5,v_6$};
$v_2$: {$v_1,v_4,v_5,v_6$};
$v_3$: {$v_1,v_4,v_5,v_6$};
$v_4$; {$v_2,v_3,v_5,v_6$};
$v_5$; {$v_1,v_2,v_3,v_4$};
$v_6$; {$v_1,v_2,v_3,v_4$}.

Инцидентность графа $G_{18}$:

$v_1$: {$e_1,e_2,e_3,e_4$};
$v_2$: {$e_1,e_5,e_6,e_7$};
$v_3$: {$e_2,e_8,e_9,e_{10}$};
$v_4$: {$e_5,e_8,e_{11},e_{12}$};
$v_5$: {$e_3,e_6,e_9,e_{11}$};
$v_6$: {$e_4,e_7,e_{10},e_{12}$}.

Количество базовых реберных разрезов графа $G_{18}$= 12

Базовые реберные разрезы графа $G_{18}$:

$w_0(e_1)$ = {$e_2,e_3,e_4,e_5,e_6,e_7$};
$w_0(e_2)$ = {$e_1,e_3,e_4,e_8,e_9,e_{10}$};
$w_0(e_3)$ = {$e_1,e_2,e_4,e_6,e_9,e_{11}$};
$w_0(e_4)$ = {$e_1,e_2,e_3,e_7,e_{10},e_{12}$};
$w_0(e_5)$ = {$e_1,e_6,e_7,e_8,e_{11},e_{12}$};
$w_0(e_6)$ = {$e_1,e_3,e_5,e_7,e_9,e_{11}$};
$w_0(e_7)$ = {$e_1,e_4,e_5,e_6,e_{10},e_{12}$};
$w_0(e_8)$ = {$e_2,e_5,e_9,e_{10},e_{11},e_{12}$};
$w_0(e_9)$ = {$e_2,e_3,e_6,e_8,e_{10},e_{11}$};
$w_0(e_{10})$ = {$e_2,e_4,e_7,e_8,e_9,e_{12}$};
$w_0(e_{11})$ = {$e_3,e_5,e_6,e_8,e_9,e_{12}$};
$w_0(e_{12})$ = {$e_4,e_5,e_7,e_8,e_{10},e_{11}$}.

Кортеж весов ребер $\xi_w(G_{18})$ = <12×14>, кортеж весов вершин $\zeta_w(G_{18})$ = <6×56>.

Количество уровней = 2. Множество изометрических циклов:

$c_1$ = {$e_1,e_3,e_6$}; $c_2$ = {$e_2,e_3,e_{11}$}; $c_3$ = {$e_8,e_9,e_{11}$};
$c_4$ = {$e_5,e_6,e_9$}; $c_5$ = {$e_2,e_4,e_{12}$}; $c_6$ = {$e_5,e_7,e_{10}$};
$c_7$ = {$e_1,e_4,e_7$}; $c_8$ = {$e_8,e_{10},e_{12}$}; $c_9$ = {$e_1,e_2,e_5,e_8$};
$c_{10}$ = {$e_3,e_4,e_9,e_{10}$}; $c_{11}$ = {$e_6,e_7,e_{11},e_{12}$}.

Образы изометрических циклов реберного графа $L(G_{18})$ в виде ребер графа $G_{18}$:

$c_1^L$ = {$e_1,e_2,e_3$} тройка;
$c_2^L$ = {$e_1,e_2,e_4$} тройка;
$c_3^L$ = {$e_1,e_2,e_5,e_8$} изометрический цикл;
$c_4^L$ = {$e_1,e_2,e_6,e_{11}$} дубль-цикл;
$c_5^L$ = {$e_1,e_2,e_7,e_{12}$} дубль-цикл;
$c_6^L$ = {$e_1,e_3,e_4$} тройка;
$c_7^L$ = {$e_1,e_3,e_5,e_9$} дубль-цикл;
$c_8^L$ = {$e_1,e_3,e_6$} изометрический цикл;

$c_9^L = \{e_1,e_4,e_5,e_{10}\}$ дубль-цикл;

$c_{10}^L = \{e_1,e_4,e_7\}$ изометрический цикл;

$c_{11}^L = \{e_1,e_5,e_6\}$ тройка;

$c_{12}^L = \{e_1,e_5,e_7\}$ тройка;

$c_{13}^L = \{e_1,e_6,e_7\}$ тройка;

$c_{14}^L = \{e_2,e_3,e_4\}$ тройка;

$c_{15}^L = \{e_2,e_3,e_8,e_9\}$ дубль-цикл;

$c_{16}^L = \{e_2,e_3,e_{11}\}$ изометрический цикл;

$c_{17}^L = \{e_2,e_4,e_8,e_{10}\}$ дубль-цикл;

$c_{18}^L = \{e_2,e_4,e_{12}\}$ изометрический цикл;

$c_{19}^L = \{e_2,e_8,e_{11}\}$ тройка;

$c_{20}^L = \{e_2,e_8,e_{12}\}$ тройка;

$c_{21}^L = \{e_2,e_{11},e_{12}\}$ тройка;

$c_{22}^L = \{e_3,e_4,e_6,e_7\}$ дубль - цикл;

$c_{23}^L = \{e_3,e_4,e_9,e_{10}\}$ дубль - цикл;

$c_{24}^L = \{e_3,e_4,e_{11},e_{12}\}$ изометрический цикл;

$c_{25}^L = \{e_3,e_6,e_9\}$ тройка;

$c_{26}^L = \{e_3,e_6,e_{11}\}$ тройка;

$c_{27}^L = \{e_3,e_9,e_{11}\}$ тройка;

$c_{28}^L = \{e_4,e_7,e_{10}\}$ тройка;

$c_{29}^L = \{e_4,e_7,e_{12}\}$ тройка;

$c_{30}^L = \{e_4,e_{10},e_{12}\}$ тройка;

$c_{31}^L = \{e_5,e_6,e_7\}$ тройка;

$c_{32}^L = \{e_5,e_6,e_8,e_{11}\}$ дубль - цикл;

$c_{33}^L = \{e_5,e_6,e_9\}$ изометрический цикл;

$c_{34}^L = \{e_5,e_7,e_8,e_{12}\}$ дубль - цикл;

$c_{35}^L = \{e_5,e_7,e_{10}\}$ изометрический цикл;

$c_{36}^L = \{e_5,e_8,e_9\}$ тройка;

$c_{37}^L = \{e_5,e_8,e_{10}\}$ тройка;

$c_{38}^L = \{e_5,e_9,e_{10}\}$ тройка;

$c_{39}^L = \{e_6,e_7,e_9,e_{10}\}$ изометрический цикл;

$c_{40}^L = \{e_6,e_7,e_{11},e_{12}\}$ дубль - цикл

$c_{41}^L = \{e_6,e_9,e_{11}\}$ тройка;

$c_{42}^L = \{e_7,e_{10},e_{12}\}$ тройка;

$c_{43}^L = \{e_8,e_9,e_{10}\}$ тройка;

$c_{44}^L = \{e_8,e_9,e_{11}\}$ изометрический цикл;

$c_{45}^L = \{e_8,e_{10},e_{12}\}$ изометрический цикл;

$c_{46}^L = \{e_8,e_{11},e_{12}\}$ тройка;

$c_{47}^{L}$ = {$e_9$,$e_{10}$,$e_{11}$,$e_{12}$}  дубль – цикл.

Кортеж весов ребер: $\xi_L$ ($e_i$) = <13,13,13,13,13,13,13,13,13,13,13,13>

Кортеж весов вершин: $\varsigma_L$ ($e_i$) = <52,52,52,52,52,52>.

Вектор весов ребер: F($\xi_L$ ($G_{18}$)) = (13,13,13,13,13,13,13,13,13,13,13,13>.

Вектор весов вершин: F($\zeta_L$ ($G_{18}$)) = (52,52,52,52,52,52).

Структура реберного графа L($G_{18}$):

- количество вершин в графе L($G_{18}$) $n$ =12;
- количество ребер в графе L($G_{18}$) $m$ = 36;
- количество изометрических циклов реберного графа L($G_{18}$) = $card(C_\tau(L(G_{18}))$ = 47;
- количество изометрических циклов графа $G_{18}$ = $card(C_\tau)$ = 11;
- количество троек во множестве изоморфных циклов L($G_{18}$), $k_3$ = 24;
- количество дубль-циклов во множестве изоморфных циклов L($G_{18}$), $k_4$ = 12;
- цикломатическое число L($G_{18}$) = m-n+1=36-12+1 = 25;
- $card(C_\tau(L(G_{18})))$ = $card(C_\tau)$ +$k_3$+$k_4$ = 47 = 11+24+12.

Рассмотрим непланарный граф Петерсена $G_{19}$ (см. рис. 5.6).

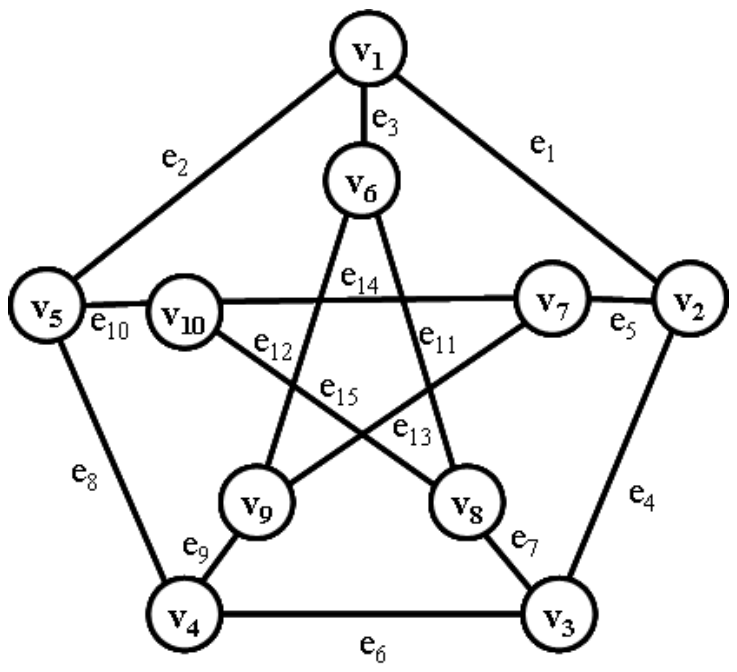

Рис. 5.6. Граф Петерсена $G_{19}$.

Смежность графа $G_{19}$:

$v_1$: {$v_2$,$v_5$,$v_6$};
$v_2$: {$v_1$,$v_3$,$v_7$};
$v_3$: {$v_2$,$v_4$,$v_8$};
$v_4$: {$v_3$,$v_5$,$v_9$};
$v_5$: {$v_1$,$v_4$,$v_{10}$};
$v_6$: {$v_1$,$v_8$,$v_9$};
$v_7$: {$v_2$,$v_9$,$v_{10}$};
$v_8$: {$v_3$,$v_6$,$v_{10}$};
$v_9$: {$v_4$,$v_6$,$v_7$};
$v_{10}$: {$v_5$,$v_7$,$v_8$}.

Инцидентность графа:

$v_1$: {$e_1$,$e_2$,$e_3$};

$v_2$: {$e_1$,$e_4$,$e_5$};
$v_3$: {$e_4$,$e_6$,$e_7$};
$v_4$: {$e_6$,$e_8$,$e_9$};
$v_5$: {$e_2$,$e_8$,$e_{10}$};
$v_6$: {$e_3$,$e_{11}$,$e_{12}$};
$v_7$: {$e_5$,$e_{13}$,$e_{14}$};
$v_8$: {$e_7$,$e_{11}$,$e_{15}$};
$v_9$: {$e_9$,$e_{12}$,$e_{13}$};
$v_{10}$: {$e_{10}$,$e_{14}$,$e_{15}$};

Количество базовых реберных разрезов графа = 15

Базовые реберные разрезы графа:

$w_0(e_1)$ = {$e_2$,$e_3$,$e_4$,$e_5$};
$w_0(e_2)$ = {$e_1$,$e_3$,$e_8$,$e_{10}$};
$w_0(e_3)$ = {$e_1$,$e_2$,$e_{11}$,$e_{12}$};
$w_0(e_4)$ = {$e_1$,$e_5$,$e_6$,$e_7$};
$w_0(e_5)$ = {$e_1$,$e_4$,$e_{13}$,$e_{14}$};
$w_0(e_6)$ = {$e_4$,$e_7$,$e_8$,$e_9$};
$w_0(e_7)$ = {$e_4$,$e_6$,$e_{11}$,$e_{15}$};
$w_0(e_8)$ = {$e_2$,$e_6$,$e_9$,$e_{10}$};
$w_0(e_9)$ = {$e_6$,$e_8$,$e_{12}$,$e_{13}$};
$w_0(e_{10})$ = {$e_2$,$e_8$,$e_{14}$,$e_{15}$};
$w_0(e_{11})$ = {$e_3$,$e_7$,$e_{12}$,$e_{15}$};
$w_0(e_{12})$ = {$e_3$,$e_9$,$e_{11}$,$e_{13}$};
$w_0(e_{13})$ = {$e_5$,$e_9$,$e_{12}$,$e_{14}$};
$w_0(e_{14})$ = {$e_5$,$e_{10}$,$e_{13}$,$e_{15}$};
$w_0(e_{15})$ = {$e_7$,$e_{10}$,$e_{11}$,$e_{14}$}.

Количество вершин в реберном графе = 15.

Количество ребер в реберном графе = 30.

Количество изометрических циклов в реберном графе = 22.

Образы изометрических циклов реберного графа в виде ребер графа $G_{19}$:

$c_1^L$ = {$e_1$,$e_2$,$e_3$} тройка;
$c_2^L$ = {$e_1$,$e_2$,$e_4$,$e_6$,$e_8$} изометрический цикл;
$c_3^L$ = {$e_1$,$e_2$,$e_5$,$e_{10}$,$e_{14}$} изометрический цикл;
$c_4^L$ = {$e_1$,$e_3$,$e_4$,$e_7$,$e_{11}$} изометрический цикл;
$c_5^L$ = {$e_1$,$e_3$,$e_5$,$e_{12}$,$e_{13}$} изометрический цикл;
$c_6^L$ = {$e_1$,$e_4$,$e_5$} тройка;
$c_7^L$ = {$e_2$,$e_3$,$e_8$,$e_9$,$e_{12}$} изометрический цикл;
$c_8^L$ = {$e_2$,,$e_3$,$e_{10}$,$e_{11}$,$e_{15}$} изометрический цикл;
$c_9^L$ = {$e_2$,$e_8$,$e_{10}$} тройка;
$c_{10}^L$ = {$e_3$,$e_{11}$,$e_{12}$} тройка;
$c_{11}^L$ = {$e_4$,$e_5$,$e_6$,$e_9$,$e_{13}$} изометрический цикл;
$c_{12}^L$ = {$e_4$,$e_5$,$e_7$,$e_{14}$,$e_{15}$} изометрический цикл;
$c_{13}^L$ = {$e_4$,$e_6$,$e_7$} тройка;
$c_{14}^L$ = {$e_5$,$e_{13}$,$e_{14}$} тройка;

$c_{15}^{L}$ = {$e_6,e_7,e_8,e_{10},e_{15}$} изометрический цикл;
$c_{16}^{L}$ = {$e_6,e_7,e_9,e_{11},e_{12}$} изометрический цикл;
$c_{17}^{L}$ = {$e_6,e_8,e_9$} тройка;
$c_{18}^{L}$ = {$e_7,e_{11},e_{15}$} тройка;
$c_{19}^{L}$ = {$e_8,e_9,e_{10},e_{13},e_{14}$} изометрический цикл;
$c_{20}^{L}$ = {$e_9,e_{12},e_{13}$} тройка;
$c_{21}^{L}$ = {$e_{10},e_{14},e_{15}$} тройка;
$c_{22}^{L}$ = {$e_{11},e_{12},e_{13},e_{14},e_{15}$} изометрический цикл;

Кортеж весов ребер : $\xi_L$ ($G_{19}$) = <6,6,6,6,6,6,6,6,6,6,6,6,6,6,6>.

Кортеж весов вершин : $\varsigma_L$ ( $G_{19}$) = <18,18,18,18,18,18,18,18,18,18>.

Вектор весов ребер : F($\xi_L$ ($G_{19}$)) =(6,6,6,6,6,6,6,6,6,6,6,6,6,6,6>.

Вектор весов вершин : F($\zeta_L$ ($G_{19}$)) =(18,18,18,18,18,18,18,18,18,18).

Структуру реберного графа L($G_{19}$) не имеет дубль-циклов.

Структура реберного графа L($G_{19}$):

- количество вершин в графе L($G_{19}$) $n$ =15;
- количество ребер в графе L($G_{19}$) $m$ = 30;
- количество изометрических циклов реберного графа L($G_{19}$) = $card(C_\tau(L(G_{19}))$ = 22;
- количество изометрических циклов графа $G_{19}$ = $card(C_\tau)$ = 12;
- количество троек во множестве изоморфных циклов L($G_{19}$) $k_3$ = 10;
- количество четверок во множестве изоморфных циклов L($G_{19}$) $k_4$ = 0;
- цикломатическое число L($G_{19}$) = m-n+1=30-15+1 = 16;
- $card(C_\tau(L(G_{19}))$ = $card(C_\tau)$ +$k_3$+$k_4$ = 22 = 10+12+0.

Рассмотрим двудольный граф $G_{20}$ (см. рис. 5.7).

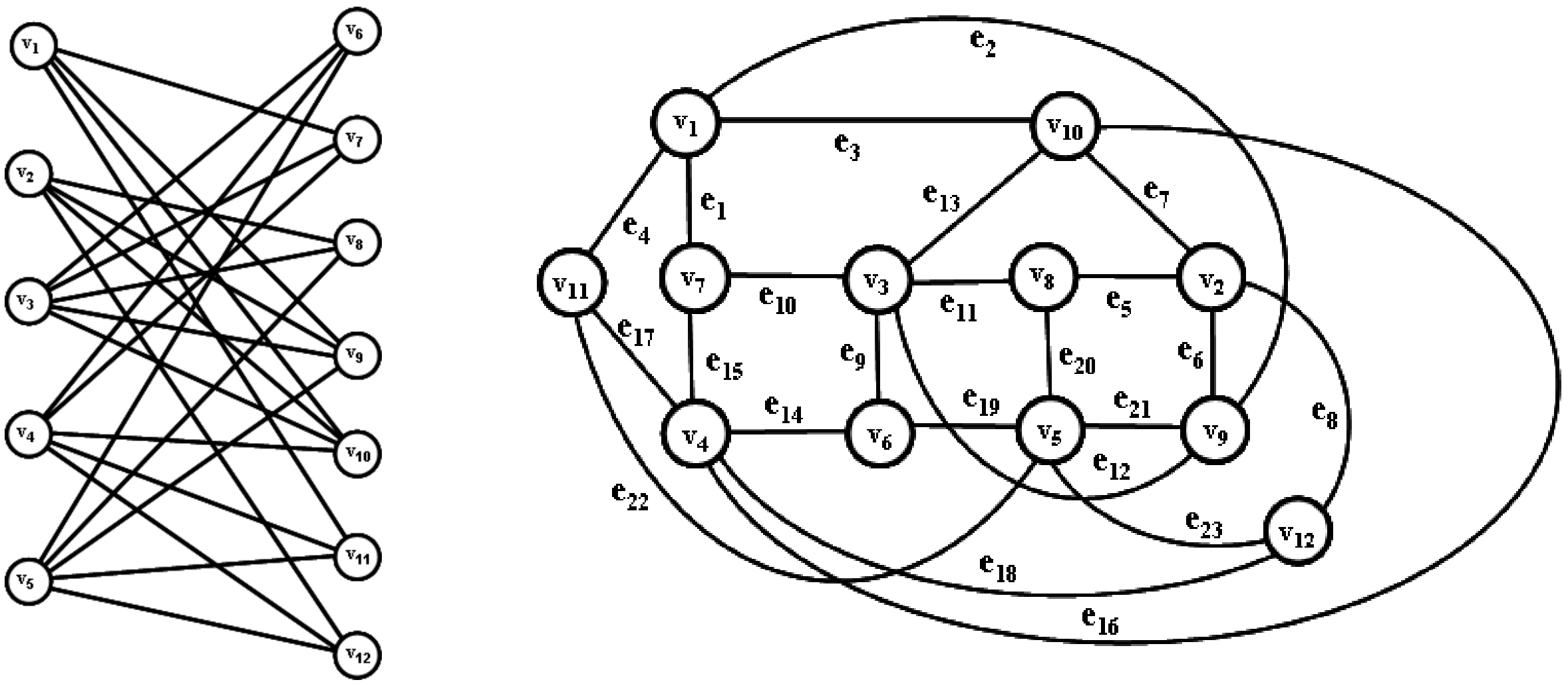

Рис. 5.7. Двудольный граф $G_{20}$.

Смежность графа $G_{20}$:

$v_1 = \{v_7,v_9,v_{10},v_{11}\}$;
$v_2 = \{v_8,v_9,v_{10},v_{12}\}$;
$v_3 = \{v_6,v_7,v_8,v_9,v_{10}\}$;
$v_4 = \{v_6,v_7,v_{10},v_{11},v_{12}\}$;
$v_5 = \{v_6,v_8,v_9,v_{11},v_{12}\}$;
$v_6 = \{v_3,v_4,v_5\}$;
$v_7 = \{v_1,v_3,v_4\}$;
$v_8 = \{v_2,v_3,v_5\}$;
$v_9 = \{v_1,v_2,v_3,v_5\}$;
$v_{10} = \{v_1,v_2,v_3,v_4\}$;
$v_{11} = \{v_1,v_4,v_5\}$;
$v_{12} = \{v_2,v_4,v_5\}$.

Инцидентность графа $G_{20}$:

$v_1 = \{e_1.e_2.e_3.e_4\}$;
$v_2 = \{e_5.e_6.e_7.e_8\}$;
$v_3 = \{e_9.e_{10}.e_{11}.e_{12}.e_{13}\}$;
$v_4 = \{e_{14}.e_{15}.e_{16}.e_{17}.e_{18}\}$;
$v_5 = \{e_{19}.e_{20}.e_{21}.e_{22}.e_{23}\}$;
$v_6 = \{e_9.e_{14}.e_{19}\}$;
$v_7 = \{e_1.e_{10}.e_{15}\}$;
$v_8 = \{e_5.e_{11}.e_{20}\}$;
$v_9 = \{e_2.e_6.e_{12}.e_{21}\}$;
$v_{10} = \{e_3.e_7.e_{13}.e_{16}\}$;
$v_{11} = \{e_4.e_{17}.e_{22}\}$;
$v_{12} = \{e_8.e_{18}.e_{23}\}$.

Базовые реберные разрезы графа $G_{20}$:

$w_0(e_1) = \{e_2,e_3,e_4,e_{10},e_{15}\}$;
$w_0(e_2) = \{e_1,e_3,e_4,e_6,e_{12},e_{21}\}$;
$w_0(e_3) = \{e_1,e_2,e_4,e_7,e_{13},e_{16}\}$;
$w_0(e_4) = \{e_1,e_2,e_3,e_{17},e_{22}\}$;
$w_0(e_5) = \{e_6,e_7,e_8,e_{11},e_{20}\}$;
$w_0(e_6) = \{e_2,e_5,e_7,e_8,e_{12},e_{21}\}$;
$w_0(e_7) = \{e_3,e_5,e_6,e_8,e_{13},e_{16}\}$;
$w_0(e_8) = \{e_5,e_6,e_7,e_{18},e_{23}\}$;
$w_0(e_9) = \{e_{10},e_{11},e_{12},e_{13},e_{14},e_{19}\}$;
$w_0(e_{10}) = \{e_1,e_9,e_{11},e_{12},e_{13},e_{15}\}$;
$w_0(e_{11}) = \{e_5,e_9,e_{10},e_{12},e_{13},e_{20}\}$;
$w_0(e_{12}) = \{e_2,e_6,e_9,e_{10},e_{11},e_{13},e_{21}\}$;
$w_0(e_{13}) = \{e_3,e_7,e_9,e_{10},e_{11},e_{12},e_{16}\}$;
$w_0(e_{14}) = \{e_9,e_{15},e_{16},e_{17},e_{18},e_{19}\}$;
$w_0(e_{15}) = \{e_1,e_{10},e_{14},e_{16},e_{17},e_{18}\}$;
$w_0(e_{16}) = \{e_3,e_7,e_{13},e_{14},e_{15},e_{17},e_{18}\}$;
$w_0(e_{17}) = \{e_4,e_{14},e_{15},e_{16},e_{18},e_{22}\}$;
$w_0(e_{18}) = \{e_8,e_{14},e_{15},e_{16},e_{17},e_{23}\}$;
$w_0(e_{19}) = \{e_9,e_{14},e_{20},e_{21},e_{22},e_{23}\}$;
$w_0(e_{20}) = \{e_5,e_{11},e_{19},e_{21},e_{22},e_{23}\}$;
$w_0(e_{21}) = \{e_2,e_6,e_{12},e_{19},e_{20},e_{22},e_{23}\}$;
$w_0(e_{22}) = \{e_4,e_{17},e_{19},e_{20},e_{21},e_{23}\}$;

$w_0(e_{23}) = \{e_8,e_{18},e_{19},e_{20},e_{21},e_{22}\}$.

Смежность реберного графа $L(G_{20})$, где вершины образы ребер $G_{20}$:

$e_1 = \{e_2,e_3,e_4,e_{10},e_{15}\}$;
$e_2 = \{e_1,e_3,e_4,e_6,e_{12},e_{21}\}$;
$e_3 = \{e_1,e_2,e_4,e_7,e_{13},e_{16}\}$;
$e_4 = \{e_1,e_2,e_3,e_{17},e_{22}\}$;
$e_5 = \{e_6,e_7,e_8,e_{11},e_{20}\}$;
$e_6 = \{e_2,e_5,e_7,e_8,e_{12},e_{21}\}$;
$e_7 = \{e_3,e_5,e_6,e_8,e_{13},e_{16}\}$;
$e_8 = \{e_5,e_6,e_7,e_{18},e_{23}\}$;
$e_9 = \{e_{10},e_{11},e_{12},e_{13},e_{14},e_{19}\}$;
$e_{10} = \{e_1,e_9,e_{11},e_{12},e_{13},e_{15}\}$;
$e_{11} = \{e_5,e_9,e_{10},e_{12},e_{13},e_{20}\}$;
$e_{12} = \{e_2,e_6,e_9,e_{10},e_{11},e_{13},e_{21}\}$;
$e_{13} = \{e_3,e_7,e_9,e_{10},e_{11},e_{12},e_{16}\}$;
$e_{14} = \{e_9,e_{15},e_{16},e_{17},e_{18},e_{19}\}$;
$e_{15} = \{e_1,e_{10},e_{14},e_{16},e_{17},e_{18}\}$;
$e_{16} = \{e_3,e_7,e_{13},e_{14},e_{15},e_{17},e_{18}\}$;
$e_{17} = \{e_4,e_{14},e_{15},e_{16},e_{18},e_{22}\}$;
$e_{18} = \{e_8,e_{14},e_{15},e_{16},e_{17},e_{23}\}$;
$e_{19} = \{e_9,e_{14},e_{20},e_{21},e_{22},e_{23}\}$;
$e_{20} = \{e_5,e_{11},e_{19},e_{21},e_{22},e_{23}\}$;
$e_{21} = \{e_2,e_6,e_{12},e_{19},e_{20},e_{22},e_{23}\}$;
$e_{22} = \{e_4,e_{17},e_{19},e_{20},e_{21},e_{23}\}$;
$e_{23} = \{e_8,e_{18},e_{19},e_{20},e_{21},e_{22}\}$.

Количество вершин в реберном графе = 23
Количество ребер в реберном графе = 69
Количество изометрических циклов в реберном графе = 75
Изометрические циклы в реберном графе в виде вершин :

$c_1^L = \{e_1,e_2,e_3\}$; тройка;
$c_2^L = \{e_1,e_2,e_4\}$; тройка;
$c_3^L = \{e_1,e_2,e_{10},e_{12}\}$; изометрический цикл;
$c_4^L = \{e_1,e_3,e_4\}$; тройка;
$c_5^L = \{e_1,e_3,e_{10},e_{13}\}$; изометрический цикл;
$c_6^L = \{e_1,e_3,e_{15},e_{16}\}$; изометрический цикл;
$c_7^L = \{e_1,e_4,e_{15},e_{17}\}$; изометрический цикл;
$c_8^L = \{e_1,e_{10},e_{15}\}$; тройка;
$c_9^L = \{e_2,e_3,e_4\}$; тройка;
$c_{10}^L = \{e_2,e_3,e_6,e_7\}$; изометрический цикл;
$c_{11}^L = \{e_2,e_3,e_{12},e_{13}\}$; изометрический цикл;
$c_{12}^L = \{e_2,e_4,e_{21},e_{22}\}$; изометрический цикл;
$c_{13}^L = \{e_2,e_6,e_{12}\}$; тройка;
$c_{14}^L = \{e_2,e_6,e_{21}\}$; тройка;
$c_{15}^L = \{e_2,e_{12},e_{21}\}$; тройка;
$c_{16}^L = \{e_3,e_4,e_{16},e_{17}\}$; изометрический цикл;

$c_{17}^{L} = \{e_3, e_7, e_{13}\}$; тройка;

$c_{18}^{L} = \{e_3, e_7, e_{16}\}$; тройка;

$c_{19}^{L} = \{e_3, e_{13}, e_{16}\}$; тройка;

$c_{20}^{L} = \{e_4, e_{17}, e_{22}\}$; тройка;

$c_{21}^{L} = \{e_5, e_6, e_7\}$; тройка;

$c_{22}^{L} = \{e_5, e_6, e_8\}$; тройка;

$c_{23}^{L} = \{e_5, e_6, e_{11}, e_{12}\}$; изометрический цикл;

$c_{24}^{L} = \{e_5, e_6, e_{20}, e_{21}\}$; изометрический цикл;

$c_{25}^{L} = \{e_5, e_7, e_8\}$; тройка;

$c_{26}^{L} = \{e_5, e_7, e_{11}, e_{13}\}$; изометрический цикл;

$c_{27}^{L} = \{e_5, e_8, e_{20}, e_{23}\}$; изометрический цикл;

$c_{28}^{L} = \{e_5, e_{11}, e_{20}\}$; тройка;

$c_{29}^{L} = \{e_6, e_7, e_8\}$; тройка;

$c_{30}^{L} = \{e_6, e_7, e_{12}, e_{13}\}$; изометрический цикл;

$c_{31}^{L} = \{e_6, e_8, e_{21}, e_{23}\}$; изометрический цикл;

$c_{32}^{L} = \{e_6, e_{12}, e_{21}\}$; тройка;

$c_{33}^{L} = \{e_7, e_8, e_{16}, e_{18}\}$; изометрический цикл

$c_{34}^{L} = \{e_7, e_{13}, e_{16}\}$; тройка;

$c_{35}^{L} = \{e_8, e_{18}, e_{23}\}$; тройка;

$c_{36}^{L} = \{e_9, e_{10}, e_{11}\}$; тройка;

$c_{37}^{L} = \{e_9, e_{10}, e_{12}\}$; тройка;

$c_{38}^{L} = \{e_9, e_{10}, e_{13}\}$; тройка;

$c_{39}^{L} = \{e_9, e_{10}, e_{14}, e_{15}\}$; изометрический цикл;

$c_{40}^{L} = \{e_9, e_{11}, e_{12}\}$; тройка;

$c_{41}^{L} = \{e_9, e_{11}, e_{13}\}$; тройка;

$c_{42}^{L} = \{e_9, e_{11}, e_{19}, e_{20}\}$; изометрический цикл;

$c_{43}^{L} = \{e_9, e_{12}, e_{13}\}$; тройка;

$c_{44}^{L} = \{e_9, e_{12}, e_{19}, e_{21}\}$; изометрический цикл;

$c_{45}^{L} = \{e_9, e_{13}, e_{14}, e_{16}\}$; изометрический цикл;

$c_{46}^{L} = \{e_9, e_{14}, e_{19}\}$; тройка;

$c_{47}^{L} = \{e_{10}, e_{11}, e_{12}\}$; тройка;

$c_{48}^{L} = \{e_{10}, e_{11}, e_{13}\}$; тройка;

$c_{49}^{L} = \{e_{10}, e_{12}, e_{13}\}$; тройка;

$c_{50}^{L} = \{e_{10}, e_{13}, e_{15}, e_{16}\}$; изометрический цикл;

$c_{51}^{L} = \{e_{11}, e_{12}, e_{13}\}$; тройка;

$c_{52}^{L} = \{e_{11}, e_{12}, e_{20}, e_{21}\}$; изометрический цикл;

$c_{53}^{L} = \{e_{14}, e_{15}, e_{16}\}$; тройка;

$c_{54}^{L} = \{e_{14}, e_{15}, e_{17}\}$; тройка;

$c_{55}^{L} = \{$ e14,e15,e18$\}$; тройка;
$c_{56}^{L} = \{$ e14,e16,e17$\}$; тройка;
$c_{57}^{L} = \{$ e14,e16,e18$\}$; тройка;
$c_{58}^{L} = \{$ e14,e17,e18$\}$; тройка;
$c_{59}^{L} = \{$ e14,e17,e19,e22$\}$; изометрический цикл;
$c_{60}^{L} = \{$ e14,e18,e19,e23$\}$; изометрический цикл;
$c_{61}^{L} = \{$ e15,e16,e17$\}$; тройка;
$c_{62}^{L} = \{$ e15,e16,e18$\}$; тройка;
$c_{63}^{L} = \{$ e15,e17,e18$\}$; тройка;
$c_{64}^{L} = \{$ e16,e17,e18$\}$; тройка;
$c_{65}^{L} = \{$ e17,e18,e22,e23$\}$; изометрический цикл;
$c_{66}^{L} = \{$ e19,e20,e21$\}$; тройка;
$c_{67}^{L} = \{$ e19,e20,e22$\}$; тройка;
$c_{68}^{L} = \{$ e19,e20,e23$\}$; тройка;
$c_{69}^{L} = \{$ e19,e21,e22$\}$; тройка;
$c_{70}^{L} = \{$ e19,e21,e23$\}$; тройка;
$c_{71}^{L} = \{$ e19,e22,e23$\}$; тройка;
$c_{72}^{L} = \{$ e20,e21,e22$\}$; тройка;
$c_{73}^{L} = \{$ e20,e21,e23$\}$; тройка;
$c_{74}^{L} = \{$ e20,e22,e23$\}$; тройка;
$c_{75}^{L} = \{$ e21,e22,e23$\}$. тройка.

Кортеж весов ребер: $\xi_L(\mathrm{G}_{20}) =$
= <8,10,11,7,8,11,10,7,11,11,11,15,15,11,11,14,11,10,11,11,14,10,11>;
Кортеж весов вершин: $\zeta_L(\mathrm{G}_{20}) =$
= <6,36,63,57,57,33,30,30,50,50,28,28>
Вектор весов ребер: $\mathrm{F}(\xi_L(\mathrm{G}_{20})) =$
= (7,7,8,8,10,10,10,10,11,11,11,11,11,11,11,11,11,11,11,14,14,15,15);
Вектор весов вершин: $\mathrm{F}(\xi_L(\mathrm{G}_{20})) =$
= (8,28,30,30,33,36,36,50,50,57,57,63).

Из рассмотренного следует,

$$card(C_\tau(L(G))) = card(C_\tau) + k_3 + k_4 > \nu(L(G)). \quad (5.3)$$

То есть, количество изометрических циклов реберного графа L(G) равно количеству изометрических циклов графа G плюс количество частей центральных разрезов мощность три (троек), без учета дубль-циклов графа G.

Качественный состав реберного графа определяется характерными подмножествами графа, такими как:

- центральные разрезы графа G;
- базовые реберные разрезы графа G;

- множество изометрических циклов графа G;
- добавка, состоящая из подмножества дубль - циклов графа G.

Реберный граф связывает в единое целое элементы подпространства разрезов и подпростравства циклов графа G. Кроме того реберный граф обладает свойством отображать часть центральных разрезов графа G (тройка) в треугольные циклы графа L(G), а ребра изометрических циклов графа G в вершины графа L(G) (рис. 5.8 – 5.9).

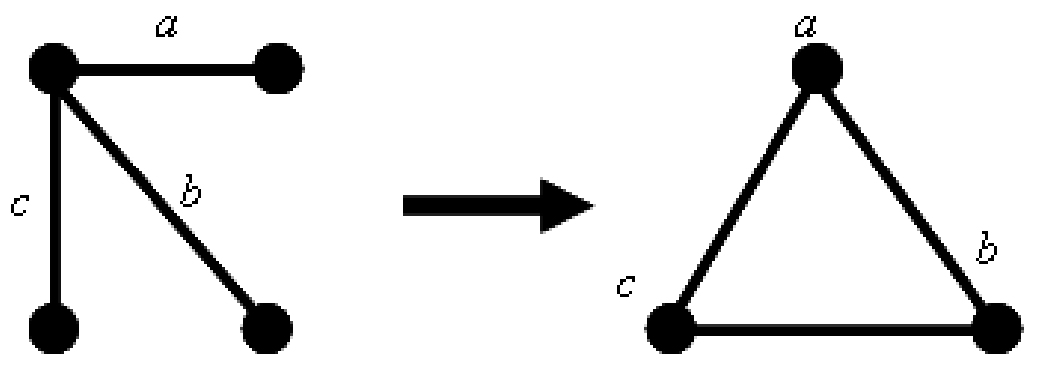

Рис. 5.8. Отображение трехреберного суграфа в цикл.

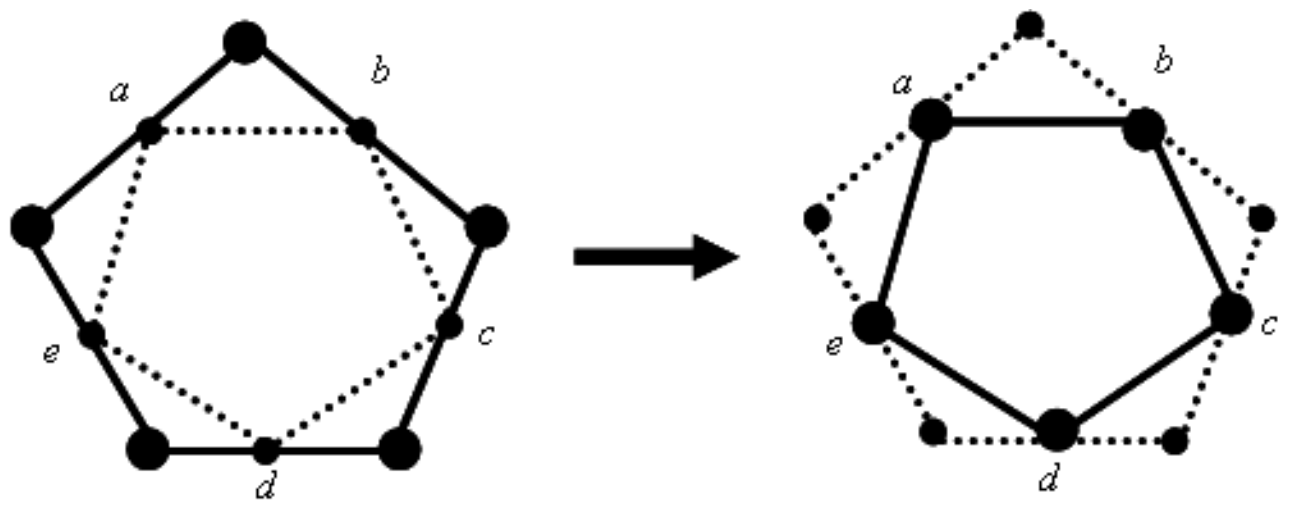

Рис. 5.9. Отображение реберного цикла в цикл с вершинами.

Теорема Уитни утверждает, что выделенное подмножество суграфов графа G является необходимым и достаточным условием для распознавания изоморфизма графов. С вычислительной точки зрения, такое выделение подмножества суграфов более оправдано, чем перебор всего множества суграфов графа G.

Следует заметить, что в состав изометрических циклов реберных графов $L(G_3)$ и $L(G_4)$ не входят образы дубль-циклов.

## Комментарии

Структуру реберного графа L(G) можно характеризовать элементами графа G. Образами центральных разрезов реберного графа L(G) могут служить базовые реберные разрезы графа G. Образами изометрических циклов реберного графа L(G) могут служить множество троек центральных разрезов графа G, множество изометрических циклов графа G и подмножество дубль-циклов графа G. Для применения условий теоремы Уитни об установлении изоморфизма реберных графов L(G) и L(H), построен цифровой инвариант, связывающий разрезы и циклы графа G и графа H в единое целое. Если цифровые инварианты реберных графов G и H равны, то по теореме Уитни графы изоморфны.

Следует заметить, что в структуру инварианта, построенного на множестве изометрических циклов реберного графа L(G), входят как центральные разрезы, так и изометрические циклы исходного графа G. А вот дубль–циклы могут присутствовать, но могут и не присутствовать в структуре реберного графа. Кольцевое сложение дубль-циклов, порождает топологический рисунок фрагмента плоской части графа G.

# Глава 6. ПОДПРОСТРАНСТВА ЦИКЛОВ И РАЗРЕЗОВ ГРАФА

## 6.1. Изоморфизм n-мерных пространств

**Определение 6.1.** Линейные пространства R и $R^*$ называются изоморфными, если между векторами $x \in$ R и векторами $x^* \in R^*$ можно установить взаимно однозначное соответствие $x \leftrightarrow x^*$ так, что если вектору $x$ соответствует вектор $x^*$, а вектору $y$ соответствует вектор $y^*$, то

1. вектору $x + y$ соответствует вектор $x^* + y^*$;
2. вектору $\lambda x$ соответствует вектор $\lambda x^*$;
3. из определения изоморфизма следует, что если $x,y,...$ – векторы из R, а $x^*,y^*,...$ – соответствующие им вектора из $R^*$, то равенство $\lambda x + \mu y + ... + ... = 0$ равносильно равенству $\lambda x^* + \mu y^* + ... + ... = 0$. Следовательно, линейно независимым векторам из R соответствуют линейно независимые векторы из $R^*$ и наоборот.

Возникает вопрос, какие пространства изоморфны между собой и какие нет.

Два пространства различной размерности заведомо не изоморфны друг другу.

В самом деле, пусть R и R* изоморфны. Из сделанного выше замечания следует, что максимальное число линейно независимых векторов в R и $R^*$ одно и то же, т.е. размерности пространства R и $R^*$ равны. Следовательно, пространства различной размерности не могут быть изоморфны между собой.

**Теорема 6.1.** Все пространства, имеющие одну и ту же размерность $n$, изоморфны друг другу [4].

Пусть задан вектор

$$x = \lambda y + \mu z + ... + \xi v \tag{6.1}$$

Если вектор $x$ выражается через векторы $y,z,...,v$, то $x$ есть линейная комбинация векторов $y,z,...,v$.

Доказательство. Пусть R и $R^*$ – два $n$-мерных пространства. Выберем в R базис $e_1,e_2,...,e_n$ и в $R^*$ $e_1^*,e_2^*,...,e_n^*$. Поставим в соответствие вектору

$$x = \xi_1 e_1 + \xi_2 e_2 + ... + \xi_n e_n \tag{6.2}$$

вектор $x^* = \xi_1 e_1^* + \xi_2 e_2^* + ... + \xi_n e_n^*$, т.е. линейную комбинацию векторов $e_i^*$ с теми же коэффициентами. Это соответствие взаимно однозначно. В самом деле, каждый вектор $x$ может быть однозначно представлен в виде (6.1). Поэтому числа $\xi_i$, а значит, и вектор $x^*$ определяются по вектору $x$ однозначно. Ввиду равноправности пространств R и $R^*$ в нашем построении каждому $x^*$ отвечает элемент из R и притом только один.

Из установленного закона соответствия сразу следует, что если $x \leftrightarrow x^*$ и $y \leftrightarrow y^*$, то $x + y = x^* + y^*$ и $\lambda x \leftrightarrow \lambda x^*$. Изоморфизм пространств R и R$^*$, таким образом, доказан.

Итак, единственной существенной характеристикой конечномерного линейного пространства является его размерность [4].

Перейдем к рассмотрению примеров. Рассматривая примеры, мы иногда для удобства будем обозначать ребра числами, а номера ребер в векторах будут читаться справа налево.

***Пример 6.1.*** Рассмотрим граф $G_{21}$ = (V,E) представленный на рис. 6.1.

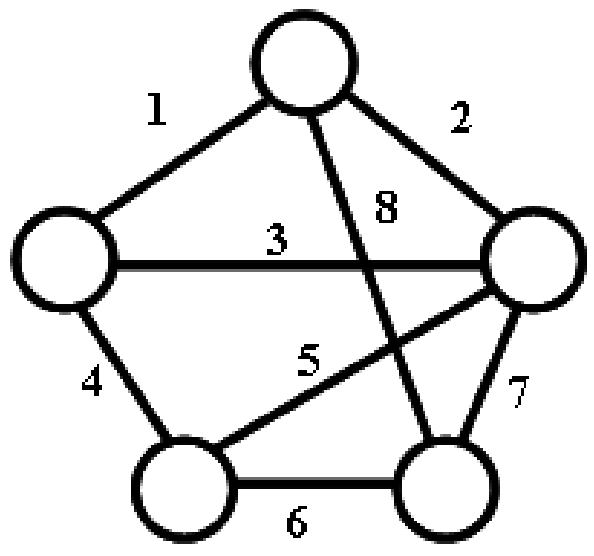

Рис. 6.1. Граф $G_{21}$.

Запишем базис пространства £$_G$ в виде однореберных суграфов:

$e_1$ = (0,0,0,0,0,0,0,1);
$e_2$ = (0,0,0,0,0,0,1,0);
$e_3$ = (0,0,0.0,0,1,0,0);
$e_4$ = (0,0,0,0,1,0,0,0);
$e_5$ = (0,0,0,1,0,0,0,0);
$e_6$ = (0,0,1,0,0,0,0,0);
$e_7$ = (0,1,0,0,0,0,0,0);
$e_8$ = (1,0,0,0,0,0,0,0).

Любой вектор пространства £$_G$, может быть представлен в виде:

$f_0$ = (0,0,0,0,0,0,0,0); $f_1$ = (0,0,0,0,0,0,0,1);
$f_2$ = (0,0,0,0,0,0,1,0); $f_3$ = (0.0,0,0,0,0.1,1);
…………………………………..……………..
$f_{254}$ = (1,1,1,1,1,1,1,0); $f_{255}$ = (1,1,1,1,1,1,1,1).

Для данного графа можно записать следующий базис:

$q_1$ = (0,0,0,0,0,0,0,1);
$q_2$ = (0,0,0,0,0,0,1,1);
$q_3$ = (0,0,0,0,0,1,1,1);
$q_4$ = (0,0,0,0,1,1,1,1);
$q_5$ = (0,0,0,1,1,1,1,1);
$q_6$ = (0,0,1,1,1,1,1,1);
$q_7$ = (0,1,1,1,1,1,1,1);
$q_8$ = (1,1,1,1,1,1,1,1).

Размер данного базиса равен восьми, что не противоречит лемме о векторном пространстве.

Следующая система независимых векторов:

$c_1 = 1\times e_1 + 1\times e_2 + 1\times e_3 + 0\times e_4 + 0\times e_5 + 0\times e_6 + 0\times e_7 + 0\times e_8;$
$c_2 = 0\times e_1 + 0\times e_2 + 1\times e_3 + 1\times e_4 + 1\times e_5 + 0\times e_6 + 0\times e_7 + 0\times e_8;$
$c_3 = 0\times e_1 + 0\times e_2 + 0\times e_3 + 0\times e_4 + 1\times e_5 + 1\times e_6 + 1\times e_7 + 0\times e_8;$
$c_4 = 0\times e_1 + 1\times e_2 + 0\times e_3 + 0\times e_4 + 0\times e_5 + 0\times e_6 + 1\times e_7 + 1\times e_8.$

характеризующая циклы графа и их количество, точно также не противоречит лемме о векторном пространстве.

***Пример 6.2.*** Рассмотрим два графа G и $G^*$, представленных на рис. 6.2.

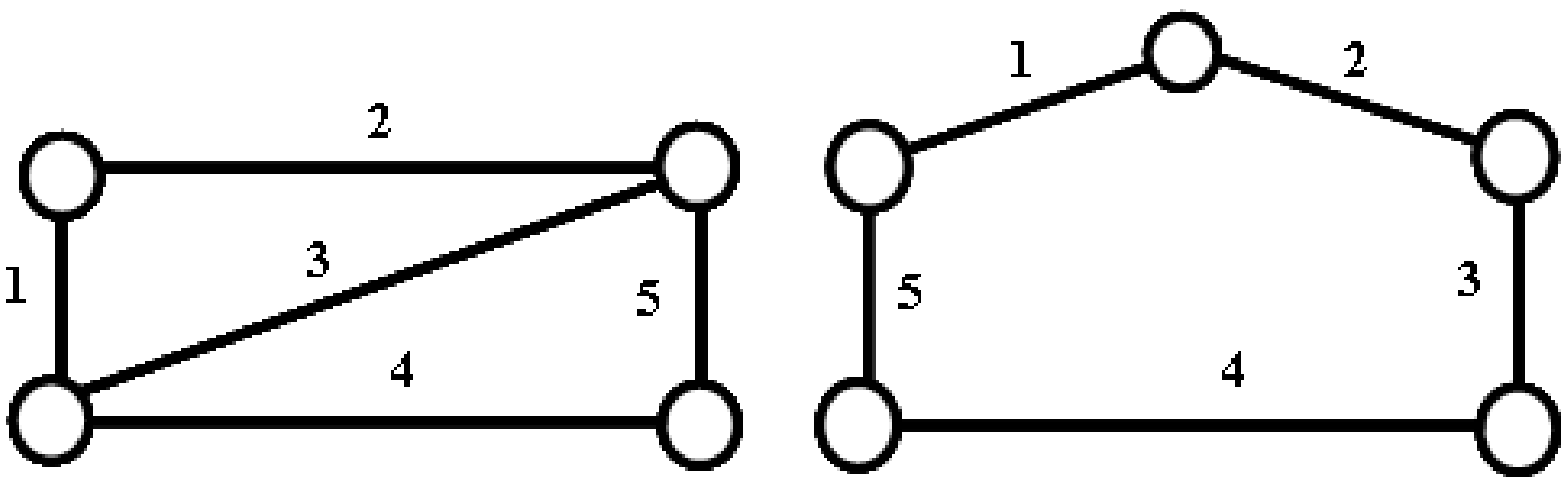

Рис.6.2. Графы G и $G^*$.

Рассмотрим пространство ребер для графа G и пространство ребер графа $G^*$.

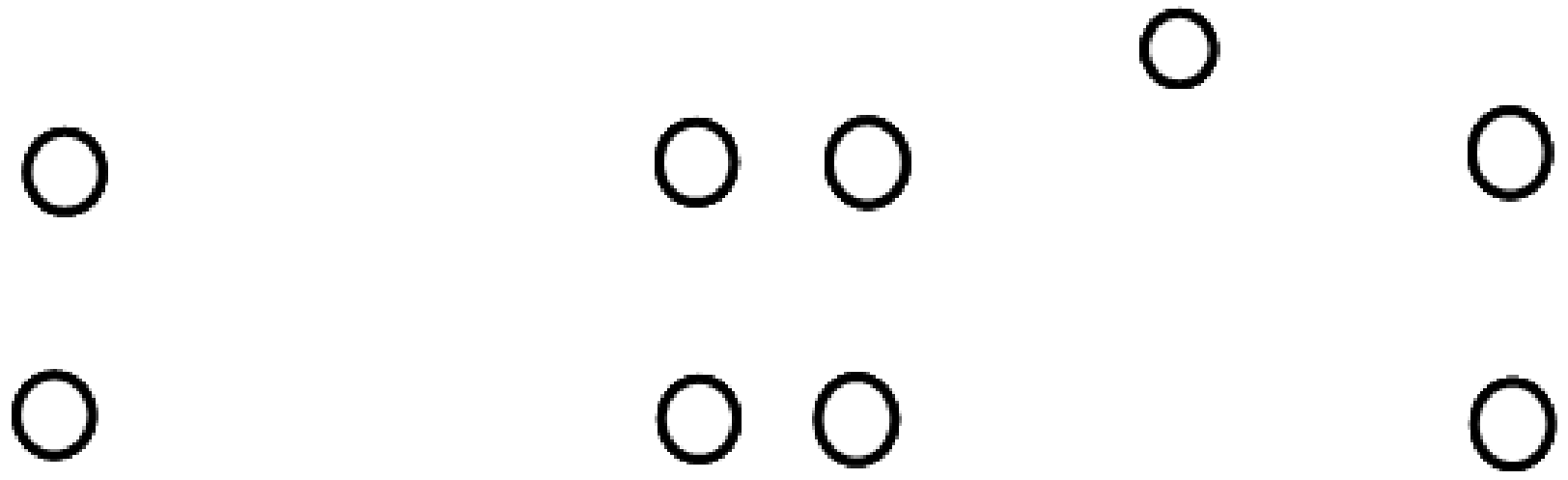

$f_0^*$ =

$f_0$ = (0,0,0,0,0) (0,0,0,0,0)

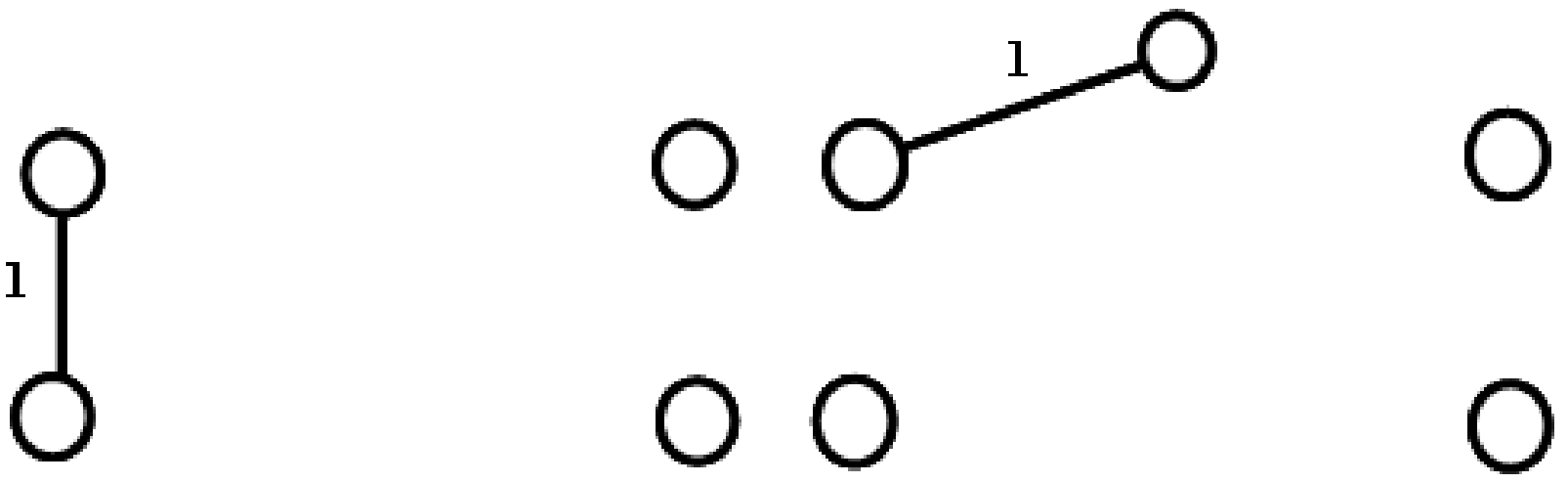

$f_1$ = (0,0,0,0,1) $f_1^*$ = (0,0,0,0,1)

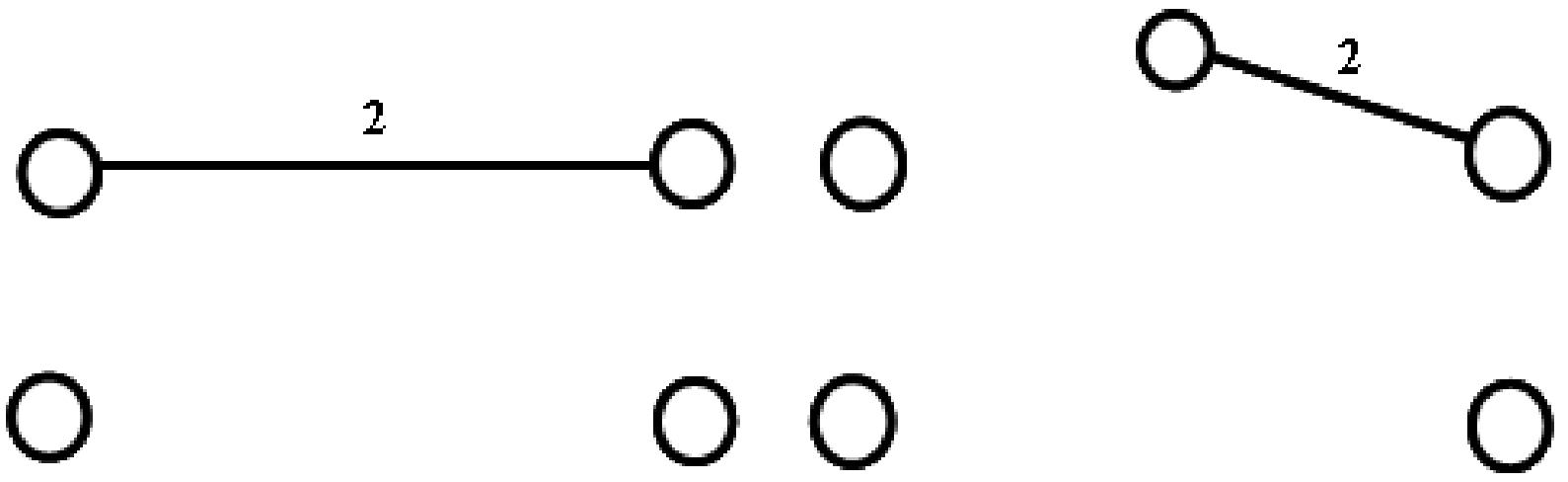

$f_2$ = (0,0,0,1,0) $f_2^*$ = (0,0,0,1,0)

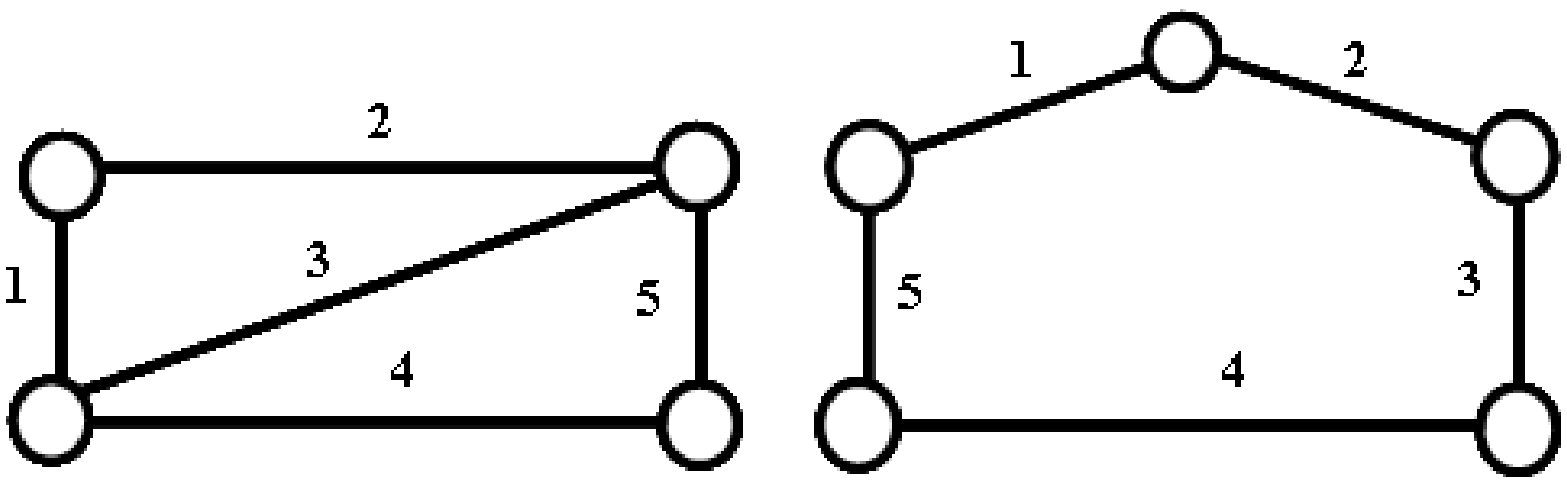

$f_{31}$ = (1,1,1,1,1) $f_{31}^*$ = (1,1,1,1,1)

Таким образом, пространства ребер графов G и $G^*$ – изоморфны, т.к. $m = m^*$.

## 6.2. Разложение пространства в прямую сумму подпространств.

Пусть задано два подпространства n-мерного пространства R. Обозначим их как $R_1$ и $R_2$.

**Определение 6.2.** Если каждый вектор *x* пространства R можно единственным образом представить как сумму двух векторов $x = x_1+x_2$, где $x_1 \in R_1$, а $x_2 \in R_2$, то пространство R можно представить в виде суммы подпространств: $R = R_1 + R_2$.

**Теорема 6.2** [4]. Для того чтобы пространство R разлагалось в прямую сумму подпространств $R_1$ и $R_2$ достаточно, чтобы:

1. Подпространства $R_1$ и $R_2$ имели только один общий вектор $x=0$ (нулевой вектор).
2. Сумма размерностей этих подпространств была равна размерности пространства R.

*Доказательство.* Выберем некоторый базис $e_1, e_2,...,e_k$ в подпространстве $R_1$ и базис $f_1, f_2,\ldots, f_l$ в подпространстве $R_2$. Поскольку сумма размерностей $R_1$ и $R_2$ есть *n,* то общее число этих векторов $k + l = n$.

Покажем, что векторы $e_1,e_2,...,e_k, f_1, f_2,\ldots, f_l$ линейно независимы, т.е. образуют базис пространства R. Действительно, пусть $\lambda_1 e_1 + ... + \lambda_k e_k + \mu_1 f_1 + ... + \mu_l f_l = 0$, следовательно $\lambda_1 e_1 + ... + \lambda_k e_k = - \mu_1 f_1 - ... - \mu_l f_l$ . Левая часть этого равенства есть вектор из $R_1$, а правая – из $R_2$.

Так как по условию единственный общий вектор $R_1$ и $R_2$ есть нулевой вектор, то

$$\lambda_1 e_1 + ... + \lambda_k e_k = 0;$$

$$\mu_1 f_1 + ... + \mu_l f_l = 0; \qquad (6.3)$$

Но каждый из наборов $e_1,e_2,...,e_k$ и $f_1, f_2,\ldots, f_l$ состоит из линейно независимых векторов, так как это базисы в $R_1$ и $R_2$. Поэтому из первого равенства (6.3) следует, что $\lambda_1 = \lambda_2 = \ldots = \lambda_k = 0$, а из второго – что $\mu_1 = \mu_2 = \ldots = \mu_l = 0$.

Следовательно, система $e_1, e_2,...,e_k, f_1, f_2,\ldots, f_l$ состоит из *n* линейно независимых векторов, т.е. это есть базис в пространстве R.

Мы доказали, что при выполнении условий теоремы существует базис, первые $k$ векторов которого образуют базис в $R_1$, а последние $l$ – базис в $R_2$.

Произвольный вектор $x$ из R можно разложить по векторам этого базиса:

$$x = \lambda_1 e_1 + ... + \lambda_k e_k + \mu_1 f_1 + ... + \mu_l f_l ,$$

при этом

$$x_1 = \lambda_1 e_1 + ... + \lambda_k e_k \in R_1$$

$$x_2 = \mu_1 f_1 + ... + \mu_l f_l \in R_2.$$

Таким образом, $x = x_1 + x_2$, где $x_1 \in R_1$, а $x_2 \in R_2$. Покажем, что это разложение единственно. Предположим, что существует два разложения:

$$x = x_1 + x_2, \text{ где } x_1 \in R_1, \text{ а } x_2 \in R_2$$

и

$$x^* = x_1^* + x_2^*, \text{ где } x_1^* \in R_1,\ x_2^* \in R_2$$

Вычитая второе равенство из первого, получаем $x_1 - x_1^* + x_2 - x_2^* = 0$ или $x_1 - x_1^* = x_2^* - x_2$.

Так как вектор, стоящий в левой части равенства, принадлежит $R_1$, а вектор, стоящий в правой части, принадлежит $R_2$, то каждый из этих векторов равен нулю. Единственность разложения доказана.

***Пример 6.3.*** Рассмотрим граф $G_{22}$ представленный на рис. 6.3.

Фундаментальная система циклов:

$c_1^*$ = (0,1,0,0,0,1,1) = {1,2,6};
$c_2^*$ = (0,1,1,0,1,0,0) = {3,5,6};
$c_3^*$ = (1,0,0,1,1,0,0) = {3,4,7}.

Фундаментальная система разрезов:

$s_1^*$ = (0,0,0,0,0,1,1) = {1,2};
$s_2^*$ = (1,0,1,0,1,0,0) = {3,5,7};
$s_3^*$ = (1,0,0,1,0,0,0) = {4,7};
$s_4^*$ = (0,1,1,0,0,0,1) = {1,5,6}.

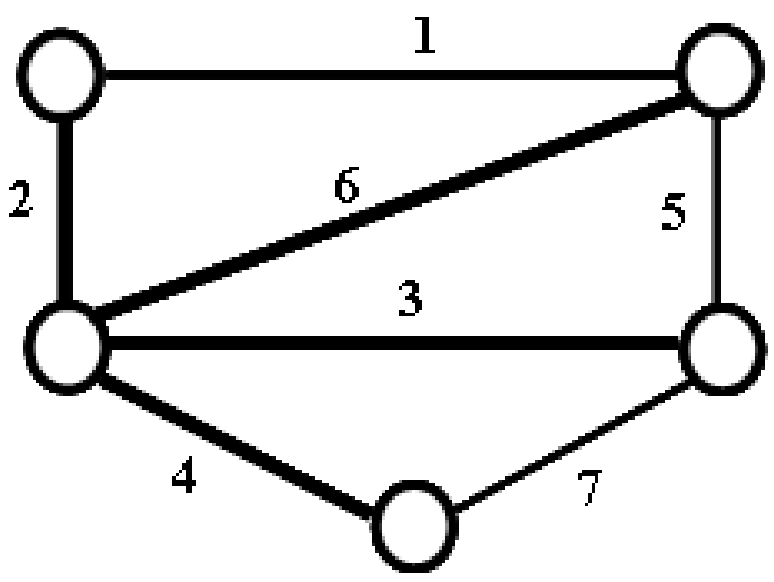

Рис. 6.3. Граф $G_{22}$.

Образуем подпространство квазициклов $C(G_{22})$:

$c_1 = c_1^* = (0,1,0,0,0,1,1) = \{1,2,6\};$

$c_2 = c_2^* = (0,1,1,0,1,0,0) = \{3,5,6\};$

$c_3 = c_3^* = (1,0,0,1,1,0,0) = \{3,4,7\};$

$c_4 = c_1^* \oplus c_2^* = (0,0,1,0,1,1,1) = \{1,2,3,5\};$

$c_5 = c_1^* \oplus c_3^* = (1,1,0,1,1,1,1) = \{1,2,3,4,6,7\};$

$c_6 = c_2^* \oplus c_3^* = (1,1,1,1,0,0,0) = \{4,5,6,7\};$

$c_7 = c_1^* \oplus c_2^* \oplus c_3^* = (1,0,1,1,0,1,1) = \{1,2,4,5,7\}.$

Образуем подпространство разрезов $S(G_{22})$:

$s_1 = s_1^* = (0,0,0,0,0,1,1) = \{1,2\};$

$s_2 = s_2^* = (1,0,1,0,1,0,0) = \{3,5,7\};$

$s_3 = s_3^* = (1,0,0,1,0,0,0) = \{4,7\};$

$s_4 = s_4^* = (0,1,1,0,0,0,1) = \{1,5,6\};$

$s_5 = s_1^* \oplus s_2^* = (1,0,1,0,1,1,1) = \{1,2,3,5,7\};$

$s_6 = s_1^* \oplus s_3^* = (1,0,0,1,0,1,1) = \{1,2,4,7\};$

$s_7 = s_1^* \oplus s_4 = (0,1,1,0,0,1,0) = \{2,5,6\};$

$s_8 = s_2^* \oplus s_3^* = (0,0,1,1,1,0,0) = \{3,4,5\};$

$s_9 = s_2^* \oplus s_4^* = (1,1,0,0,1,0,1) = \{1,3,6,7\};$

$s_{10} = s_3 \oplus s_4^* = (1,1,1,1,0,0,1) = \{1,4,5,6,7\};$

$s_{11} = s_1^* \oplus s_2^* \oplus s_3^* = (0,0,1,1,1,1,1) = \{1,2,3,4,5\};$

$s_{12} = s_1^* \oplus s_2^* \oplus s_4^* = (1,1,0,0,1,1,0) = \{2,3,6,7\};$

$s_{13} = s_1^* \oplus s_3^* \oplus s_4^* = (1,1,1,1,0,1,0) = \{2,4,5,6,7\};$

$s_{14} = s_2^* \oplus s_3^* \oplus s_4^* = (0,1,0,1,1,0,1) = \{1,3,4,6\};$

$s_{15} = s_1^* \oplus s_2^* \oplus s_3^* \oplus s_4^* = (0,1,0,1,1,1,0) = \{2,3,4,6\}.$

Составим сумму этих двух подпространств:

$\alpha_0 = (0,0,0,0,0,0,0) = c_0 \oplus s_0 = \varnothing;$

$\alpha_1 = (0,0,0,0,0,0,1) = s_{10} \oplus c_6 = \{1,4,5,6,7\} \oplus \{4,5,6,7\} = \{1\};$

$\alpha_2 = (0,0,0,0,0,1,0) = s_{13} \oplus c_6 = \{2,4,5,6,7\} \oplus \{4,5,6,7\} = \{2\};$

$\alpha_3 = (0,0,0,0,0,1,1) = s_1 = \{1,2\};$

$\alpha_4 = (0,0,0,0,1,0,0) = s_3 \oplus c_3 = \{4,7\} \oplus \{3,4,7\} = \{3\};$

$\alpha_5 = (0,0,0,0,1,0,1) = s_4 \oplus c_2 = \{1,5,6\} \oplus \{3,5,6\} = \{1,3\};$

$\alpha_6 = (0,0,0,0,1,1,0) = s_7 \oplus c_2 = \{2,5,6\} \oplus \{3,5,6\} = = \{2,3\};$

$\alpha_7 = (0,0,0,0,1,1,1) = s_6 \oplus c_3 = \{1,2,4,7\} \oplus \{3,4,7\} = \{1,2,3\};$

$\alpha_8 = (0,0,0,1,0.0,0) = s_{11} \oplus c_4 = \{1,2,3,4,5\} \oplus \{1,2.3,5\} = \{4\};$

$\alpha_9 = (0,0,0,1,0,0,1) = s_{12} \oplus c_5 = \{2,3,6,7\} \oplus \{1,2,3,4,6,7\} = \{1,4\};$

$\alpha_{10} = (0,0,0,1,0,1,0) = s_9 \oplus c_5 = \{1,3,6,7\} \oplus \{1,2,3,4,6,7\} = \{2,4\};$

$\alpha_{11} = (0,0,0,1,0,1,1) = s_8 \oplus c_4 = \{3,4,5\} \oplus \{1,2,3,5\} = \{1,2,4\};$

$\alpha_{12} = (0,0,0,1,1,0,0) = s_5 \oplus c_7 = \{1,2,3,5,7\} \oplus \{1,2,4,5,7\} = \{3,4\};$

$\alpha_{13} = (0,0,0,1,1,0,1) = s_{15} \oplus c_1 = \{2,3,4,6\} \oplus \{1,2,6\} = \{1,3,4\};$

$\alpha_{14} = (0,0,0,1,1,1,0) = s_{14} \oplus c_1 = \{1,3,4,6\} \oplus \{1,2,6\} = \{2,3,4\};$

$\alpha_{15}$ = (0,0,0,1,1,1,1) = $s_2 \oplus c_7$ = {3,5,7} ⊕ {1,2,4,5,7} = {1,2,3,4};
$\alpha_{16}$ = (0,0,1,0,0,0,0) = $s_6 \oplus c_7$ = {1,2,4,7} ⊕ {1,2,4,5,7} = {5};
$\alpha_{17}$ = (0,0,1,0,0,0,1) = $s_7 \oplus c_1$ = {2,5,6} ⊕ {1,2,6} = {1,5};
$\alpha_{18}$ = (0,0,1,0,0,1,0) = $s_4 \oplus c_1$ = {1,5,6} ⊕ {1,2,6} = {2,5};
$\alpha_{19}$ = (0,0,1,0,0,1,1) = $s_3 \oplus c_7$ = {4,7} ⊕ {1,2,4,5,7} = {1,2,5};
$\alpha_{20}$ = (0,0,1,0,1,0,0) = $s_1 \oplus c_4$ = {1,2} ⊕ {1,2,3,5} = {3,5};
$\alpha_{21}$ = (0,0,1,0,1,0,1) = $s_{13} \oplus c_5$ = {2,4,5,6,7} ⊕ {1,2,3,4,6,7} = {1,3,5};
$\alpha_{22}$ = (0,0,1,0,1,1,0) = $s_{10} \oplus c_5$ = {1,4,5,6,7} ⊕ {1,2,3,4,6,7} = {2,3,5};
$\alpha_{23}$ = (0,0,1,0,1,1,1) = $c_4$ = {1,2,3,5};
$\alpha_{24}$ = (0,0,1,1,0,0,0) = $s_2 \oplus c_3$ = {3,5,7} ⊕ {3,4,7} = {4,5};
$\alpha_{25}$ = (0,0,1,1,0,0,1) = $s_{14} \oplus c_2$ = {1,3,4,6) ⊕ {3,5,6) = {1,4,5};
$\alpha_{26}$ = (0,0,1,1,0,1,0) = $s_{15} \oplus c_2$ = {2,3,4,6} ⊕ {3,5,6} = {2,4,5};
$\alpha_{27}$ = (0,0,1,1,0,1,1) = $s_5 \oplus c_3$ = {1,2,3,5,7} ⊕ {3,4,7} = {1,2,4,5};
$\alpha_{28}$ = (0,0,1,1,1,0,0) = $s_8$ = {3,4,5};
$\alpha_{29}$ = (0,0,1,1,1,0,1) = $s_9 \oplus c_6$ = {1,3,6,7} ⊕ {4,5,6,7} = {1,3,4,5};
$\alpha_{30}$ = (0,0,1,1,1,1,0) = $s_{12} \oplus c_6$ = {2,3,6,7} ⊕ {4,5,6,7} = {2,3,4,5};
$\alpha_{31}$ = (0,0,1,1,1,1,1) = $s_{11}$ = {1,2,3,4,6};
$\alpha_{32}$ = (0,1,0,0,0,0,0) = $s_1 \oplus c_1$ = {1,2} ⊕ {1,2,6} = {6};
$\alpha_{33}$ = (0,1,0,0,0,0,1) = $s_{13} \oplus c_7$ = {2,4,5,6,7} ⊕ {1,2,4,5,7} = {1,6};
$\alpha_{34}$ = (0,1,0,0,0,1,0) = $s_{10} \oplus c_7$ = {1,4,5,6,7} ⊕ {1,2,4,5,7} = {2,6};
$\alpha_{35}$ = (0,1,0,0,0,1,1) = $c_1$ = {1,2,6};
$\alpha_{36}$ = (0,1,0,0,1,0,0) = $s_6 \oplus c_5$ = {1,2,4,7} ⊕ {1,2,3,4,6,7} = {3,6};
$\alpha_{37}$ = (0,1,0,0,1,0,1) = $s_7 \oplus c_4$ = {2,5,6} ⊕ {1,2,3,5} = {1,3,6};
$\alpha_{38}$ = (0,1,0,0,1,1,0) = $s_4 \oplus c_4$ = {1,5,6} ⊕ {1,2,3,5} = {2,3,6};
$\alpha_{39}$ = (0,1,0,0,1,1,1) = $s_3 \oplus c_5$ = {4,7} ⊕ {1,2,3,4,6,7} = {2,4,5};
$\alpha_{40}$ = (0,1,0,1,0,0,0) = $s_8 \oplus c_2$ = {3,4,5} ⊕ {3,5,6} = {4,6};
$\alpha_{41}$ = (0,1,0,1,0,0,1) = $s_9 \oplus c_3$ = {1,3,6,7} ⊕ {3,4,7} = {1,4,6};
$\alpha_{42}$ = (0,1,0,1,0,1,0) = $s_{12} \oplus c_3$ = {2,3,6,7} ⊕ {3,4,7} = {2,4,6};
$\alpha_{43}$ = (0,1,0,1,0,1,1) = $s_{11} \oplus c_2$ = {1,2,3,4,5} ⊕ {3,5,6} = {1,2,4,6};
$\alpha_{44}$ = (0,1,0,1,1,0,0) = $s_2 \oplus c_6$ = {3,5,7} ⊕ {4,5,6,7} = {3,4,6};
$\alpha_{45}$ = (0,1,0,1,1,0,1) = $s_{14}$ = {1,3,4,6};
$\alpha_{46}$ = (0,1,0,1,1,1,0) = $s_{15}$ = {2,3,4,6};
$\alpha_{47}$ = (0,1,0,1,1,1,1) = $s_5 \oplus c_6$ = {1,2,3,5,7} ⊕ {4,5,6,7} = {1,2,3,4,6};
$\alpha_{48}$ = (0,1,1,0,0,0,0) = $s_3 \oplus c_6$ = {4,7} ⊕ {4,5,6,7} = {5,6};
$\alpha_{49}$ = (0,1,1,0,0,0,1) = $s_4$ = {1,5,6};
$\alpha_{50}$ = (0,1,1,0,0,1,0) = $s_7$ = {2,5,6};
$\alpha_{51}$ = (0,1,1,0,0,1,1) = $s_6 \oplus c_6$ = {1,2,4,7} ⊕ {4,5,6,7} = {1,2,5,6};
$\alpha_{52}$ = (0,1,1,0,1,0,0) = $c_2$ = {3,5,6};
$\alpha_{53}$ = (0,1,1,0,1,0,1) = $s_{10} \oplus c_3$ = {1,4,5,6,7} ⊕ {3,4,7} = {1,3,5,6};
$\alpha_{54}$ = (0,1,1,0,1,1,0) = $s_{13} \oplus c_3$ = {2,4,5,6,7} ⊕ {3,4,7} = {2,3,5,6};

$\alpha_{55}$ = (0,1,1,0,1,1,1) = $s_1 \oplus c_2$ = {1,2} ⊕ {3,5,6} = {1,2,3,5,6};

$\alpha_{56}$ = (0,1,1,1,0,0,0) = $s_5 \oplus c_5$ = {1,2,3,5,7} ⊕ {1,2,3,4,6,7} = {4,5,6};

$\alpha_{57}$ = (0,1,1,1,0,0,1) = $s_{15} \oplus c_4$ = {2,3,4,6} ⊕ {1,2,3,5} = {1,4,5,6};

$\alpha_{58}$ = (0,1,1,1,0,1,0) = $s_{14} \oplus c_4$ = {1,3,4,6} ⊕ {1,2,3,5} = {2,4,5,6};

$\alpha_{59}$ = (0,1,1,1,0,1,1) = $s_2 \oplus c_5$ = {3,5,7} ⊕ {1,2,3,4,6,7} = {1,2,4,5,6};

$\alpha_{60}$ = (0,1,1,1,1,0,0) = $s_{11} \oplus c_1$ = {1,2,3,4,5} ⊕ {1,2,6} = {3,4,5,6};

$\alpha_{61}$ = (0,1,1,1,1,0,1) = $s_{12} \oplus c_7$ = {2,3,6,7} ⊕ {1,2,4,5,7} = {1,3,4,5,6};

$\alpha_{62}$ = (0,1,1,1,1,1,0) = $s_9 \oplus c_7$ = {1,3,6,7} ⊕ {1,2,4,5,7} = {2,3,4,5,6};

$\alpha_{63}$ = (0,1,1,1,1,1,1) = $s_8 \oplus c_1$ = {3,4,5} ⊕ {1,2,6} = {1,2,3,4,5,6};

$\alpha_{64}$ = (1,0,0,0,0,0,0) = $s_5 \oplus c_4$ = {1,2,3,5,7} ⊕ {1,2,3,5} = {7};

$\alpha_{65}$ = (1,0,0,0,0,0,1) = $s_{15} \oplus c_5$ = {2,3,4,6} ⊕ {1,2,3,4,6,7} = {1,7};

$\alpha_{66}$ = (1,0,0,0,0,1,0) = $s_{14} \oplus c_5$ = {1,3,4,6} ⊕ {1,2,3,4,6,7} = {2,7};

$\alpha_{67}$ = (1,0,0,0,0,1,1) = $s_2 \oplus c_4$ = {3,5,7} ⊕ {1,2,3,5} = {1,2,7};

$\alpha_{68}$ = (1,0,0,0,1,0,0) = $s_{11} \oplus c_7$ = {1,2,3,4,5} ⊕ {1,2,4,5,7} = {3,7};

$\alpha_{69}$ = (1,0,0,0,1,0,1) = $s_{12} \oplus c_1$ = {2,3,6,7} ⊕ {1,2,6} = {1,3,7};

$\alpha_{70}$ = (1,0,0,0,1,1,0) = $s_9 \oplus c_1$ = {1,3,6,7} ⊕ {1,2,6} = {2,3,7};

$\alpha_{71}$ = (1,0,0,0,1,1,1) = $s_8 \oplus c_7$ = {3,4,5} ⊕ {1,2,4,5,7} = {1,2,3,7};

$\alpha_{72}$ = (1,0,0,1,0,0,0) = $s_3$ = {4,7};

$\alpha_{73}$ = (1,0,0,1,0,0,1) = $s_4 \oplus c_6$ = {1,5,6} ⊕ {4,5,6,7} = {1,4,7};

$\alpha_{74}$ = (1,0,0,1,0,1,0) = $s_7 \oplus c_6$ = {2,5,6} ⊕ {4,5,6,7} = {2,4,7};

$\alpha_{75}$ = (1,0,0,1,0,1,1) = $s_6$ = {1,2,4,7};

$\alpha_{76}$ = (1,0,0,1,1,0,0) = $c_3$ = {3,4,7};

$\alpha_{77}$ = (1,0,0,1,1,0,1) = $s_{10} \oplus c_2$ = {1,4,5,6,7} ⊕ {3,5,6} = {1,3,4,7};

$\alpha_{78}$ = (1,0,0,1,1,1,0) = $s_{13} \oplus c_2$ = {2,4,5,6,7}⊕ {3,5,6} = {2,3,4,7};

$\alpha_{79}$ = (1,0,0,1,1,1,1) = $s_1 \oplus c_3$ = {1,2} ⊕ {3,4,7} = {1,2,3,4,7};

$\alpha_{80}$ = (1,0,1,0,0,0,0) = $s_8 \oplus c_3$ = {3,4,5} ⊕ {3,4,7} = {5,7};

$\alpha_{81}$ = (1,0,1,0,0,0,1) = $s_9 \oplus c_2$ = {1,3,6,7} ⊕ {3,5,6} = {1,5,7};

$\alpha_{82}$ = (1,0,1,0,0,1,0) = $s_{12} \oplus c_2$ = {2,3,6,7} ⊕ {3,5,6} = {2,5,7};

$\alpha_{83}$ = (1,0,1,0,0,1,1) = $s_{12} \oplus c_3$ = {1,2,3,4,5} ⊕ {3,4,7} = {1,2,5,7};

$\alpha_{84}$ = (1,0,1,0,1,0,0) = $s_2$ = {3,5,7};

$\alpha_{85}$ = (1,0,1,0,1,0,1) = $s_{14} \oplus c_6$ = {1,3,4,6} ⊕ {4,5,6,7} = {1,3,5,7};

$\alpha_{86}$ = (1,0,1,0,1,1,0) = $s_{15} \oplus c_6$ = {2,3,4,6} ⊕ {4,5,6,7} = {2,3,5,7};

$\alpha_{87}$ = (1,0,1,0,1,1,1) = $s_5$ = {1,2,3,5,7};

$\alpha_{88}$ = (1,0,1,1,0,0,0) = $s_1 \oplus c_7$ = {1,2} ⊕ {1,2,4,5,7} = {4,5,7};

$\alpha_{89}$ = (1,0,1,1,0,0,1) = $s_{13} \oplus c_1$ = {2,4,5,6,7} ⊕ {1,2,6} = {1,4,5,7};

$\alpha_{90}$ = (1,0,1,1,0,1,0) = $s_{10} \oplus c_1$ = {1,4,5,6,7} ⊕ (1,2,6} = {2,4,5,7};

$\alpha_{91}$ = (1,0,1,1,0,1,1) = $c_7$ = {1,2,4,5,7};

$\alpha_{92}$ = (1,0,1,1,1,0,0) = $s_6 \oplus c_4$ = {1,2,4,7} ⊕ {1,2,3,5} = {3,4,5,7};

$\alpha_{93}$ = (1,0,1,1,1,0,1) = $s_7 \oplus c_5$ = {2,5,6} ⊕ {1,2,3,4,6,7} = {1,3,4,5,7};

$\alpha_{94}$ = (1,0,1,1,1,1,0) = $s_4 \oplus c_5$ = {1,5,6} ⊕ {1,2,3,4,6,7} = {2,3,4,5,7};

$\alpha_{95}$ = (1,0,1,1,1,1,1) = $s_3 \oplus c_4$ = {4,7} ⊕ {1,2,3,5} = {1,2,3,4,5,7};

$\alpha_{96}$ = (1,1,0,0,0,0,0) = $s_2 \oplus c_2$ = {3,5,7} ⊕ {3,5,6} = {6,7};

$\alpha_{97}$ = (1,1,0,0,0,0,1) = $s_{14} \oplus c_3$ = {1,3,4,6} ⊕ {3,4,7} = {1,6,7};

$\alpha_{98}$ = (1,1,0,0,0,1,0) = $s_{15} \oplus c_3$ = {2,3,4,6} ⊕ {3,4,7} = {2,6,7};

$\alpha_{99}$ = (1,1,0,0,0,1,1) = $s_5 \oplus c_2$ = {1,2,3,5,7} ⊕ {3,5,6} = {1,2,6,7};

$\alpha_{100}$ = (1,1,0,0,1,0,0) = $s_8 \oplus c_6$ = {3,4,5} ⊕ {4,5,6,7} = {3,6,7};

$\alpha_{101}$ = (1,1,0,0,1,0,1) = $s_9$ = {1,3,6,7};

$\alpha_{102}$ = (1,1,0,0,1,1,0) = $s_{12}$ = {2,3,6,7};

$\alpha_{103}$ = (1,1,0,0,1,1,1) = $s_{11} \oplus c_6$ = {1,2,3,4,5} ⊕ {4,5,6,7} = {1,2,3,6,7};

$\alpha_{104}$ = (1,1,0,1,0,0,0) = $s_6 \oplus c_1$ = {1,2,4,7} ⊕ {1,2,6} = {4,6,7};

$\alpha_{105}$ = (1,1,0,1,0,0,1) = $s_7 \oplus c_7$ = {2,5,6} ⊕ {1,2,4,5,7} = {1,4,6,7};

$\alpha_{106}$ = (1,1,0,1,0,1,0) = $s_4 \oplus c_7$ = {1,5,6} ⊕ {1,2,4,5,7} = {2,4,6,7};

$\alpha_{107}$ = (1,1,0,1,0,1,1) = $s_3 \oplus c_1$ = {4,7}⊕ {1,2,6} = {1,2,4,6,7};

$\alpha_{108}$ = (1,1,0,1,1,0,0) = $s_1 \oplus c_5$ = {1,2} ⊕ {1,2,3,4,6,7} = {3,4,6,7};

$\alpha_{109}$ = (1,1,0,1,1,0,1) = $s_{13} \oplus c_4$ = {2,4,5,6,7} ⊕ {1,2,3,5} = {1,3,4,6,7};

$\alpha_{110}$ = (1,1,0,1,1,1,0) = $s_{10} \oplus c_4$ = {1,4,5,6,7} ⊕ {1,2,3,5} = {2,3,4,6,7};

$\alpha_{111}$ = (1,1,0,1,1,1,1) = $c_5$ = {1,2,3,4,6,7};

$\alpha_{112}$ = (1,1,1,0,0,0,0) = $s_{11} \oplus c_5$ = {1,2,3,4,5} ⊕ {1,2,3,4,6,7} = {5,6,7};

$\alpha_{113}$ = (1,1,1,0,0,0,1) = $s_{12} \oplus c_4$ = {2,3,6,7} ⊕ {1,2,3,5} = {1,5,6,7};

$\alpha_{114}$ = (1,1,1,0,0,1,0) = $s_9 \oplus c_4$ = {1,3,6,7} ⊕ {1,2,3,5} = {2,5,6,7};

$\alpha_{115}$ = (1,1,1,0,0,1,1) = $s_8 \oplus c_5$ = {3,4,5} ⊕ {1,2,3,4,6,7} = {1,2,5,6,7};

$\alpha_{116}$ = (1,1,1,0,1,0,0) = $s_5 \oplus c_1$ = {1,2,3,5,7} ⊕ {1,2,6} = {3,5,6,7};

$\alpha_{117}$ = (1,1,1,0,1,0,1) = $s_{15} \oplus c_7$ = {2,3,4,6} ⊕ {1,2,4,5,7} = {1,3,5,6,7};

$\alpha_{118}$ = (1,1,1,0,1,1,0) = $s_{14} \oplus c_7$ = {1,3,4,6} ⊕ {1,2,4,5,7} = {2,3,5,6,7};

$\alpha_{119}$ = (1,1,1,0,1,1,1) = $s_2 \oplus c_1$ = {3,5,7} ⊕ {1,2,6} = {1,2,3,5,6,7};

$\alpha_{120}$ = (1,1,1,1,0,0,0) = $c_6$ = {4,5,6,7};

$\alpha_{121}$ = (1,1,1,1,0,0,1) = $s_{10}$ = {1,4,5,6,7};

$\alpha_{122}$ = (1,1,1,1,0,1,0) = $s_{13}$ = {2,4,5,6,7};

$\alpha_{123}$ = (1,1,1,1,0,1,1) = $s_1 \oplus c_6$ = {1,2} ⊕ (4,5,6,7} = {1,2,4,5,6,7};

$\alpha_{124}$ = (1,1,1,1,1,0,0) = $s_3 \oplus c_2$ = {4,7} ⊕ {3,5,6} = {3,4,5,6,7};

$\alpha_{125}$ = (1,1,1,1,1,0,1) = $s_4 \oplus c_3$ = {1,5,6} ⊕ {3,4,7} = {1,3,4,5,6,7};

$\alpha_{126}$ = (1,1,1,1,1,1,0) = $s_7 \oplus c_3$ = {2,5,6} ⊕ {3,4,7} = {2,3,4,5,6,7};

$\alpha_{127}$ = (1,1,1,1,1,1,1) = $s_6 \oplus c_2$ = {1,2,4,7} ⊕ {3,5,6} = {1,2,3,4,5,6,7}.

Сумма этих подпространств есть пространство суграфов:

$$\dim £_G = \dim C + \dim S = 3+4=7.$$

Координатное представление вектора характеризующего суграф графа $G_{22}$ будем записывать в круглых скобках. Например, запись (1,1,0,0,1,0,1) адекватна записи множества ребер суграфа $\{e_1,e_3,e_6,e_7\}$ и читается справа налево. То же множество ребер суграфа в целях экономии будем записывать в виде цифр {1,3,6,7}, но каждый раз оговаривая вид записи.

## 6.3. Сумма и пересечение подпространств

Допустим, задано два произвольных подпространства $R_1$ и $R_2$ линейного пространства R.

Легко проверить, что совокупность векторов, принадлежащих обоим этим пространствам, также есть подпространство $R_0$ пространства R.

Это подпространство называется пересечением $R_1$ и $R_2$ и обозначается как $R_0 = R_1 \cap R_2$

Например, если $R_1$ и $R_2$ – два двухмерных подпространства трехмерного пространства, т. е. две плоскости, проходящие через начало координат, то $R_1 \cap R_2$ есть одномерное пространство, т. е. прямая, по которой пересекаются эти плоскости [8].

По двум подпространствам $R_1$ и $R_2$ можно построить еще одно подпространство, которое называется их суммой. Оно определяется следующим образом. Векторами этого подпространства является всевозможные суммы вида

$$x = x_1 + x_2, \tag{6.4}$$

где $x_1 \in R_1$, а $x_2 \in R_2$.

Легко проверить, что элементы вида (1.10) образуют подпространство. Это подпространство $R^*$ называется суммой пространств $R_1$ и $R_2$ и обозначается как $R^* = R_1 + R_2$.

Заметим, что в отличие от прямой суммы двух подпространств запись элемента из R в виде (1.9) может быть неоднозначной. Имеет место следующая теорема.

**Теорема 6.3.** Пусть заданы два подпространства $R_1$ и $R_2$ пространства R. Тогда сумма размерностей $R_1$ и $R_2$ равна размерности их суммы плюс размерность пересечения.

*Доказательство.* Выберем в пересечении $R_0 = R_1 \cap R_2$ базис

$$e_1, e_2, ..., e_k \tag{6.5}$$

Дополним этот базис с одной стороны до базиса в $R_1$

$$e_1, e_2, ..., e_k, f_1, f_2, \ldots, f_l \tag{6.6}$$

и с другой стороны – до базиса в $R_2$

$$e_1, e_2, ..., e_k, g_1, g_2, \ldots, g_m. \tag{6.7}$$

Покажем, что векторы $f_1, f_2, \ldots, f_l, e_1, e_2, ..., e_k, g_1, g_2, \ldots, g_m$ образуют базис в сумме $R^* = R_1 + R_2$.

Сначала покажем, что эти векторы линейно независимы. Действительно, пусть

$$\lambda_1 f_1 + \ldots + \lambda_l f_l + \mu_1 e_1 + \ldots + \mu_k e_k + \vartheta_1 g_1 + \ldots + \vartheta_m g_m = 0. \tag{6.8}$$

Тогда $\lambda_1 f_1 + \ldots + \lambda_l f_l + \mu_1 e_1 + \ldots + \mu_k e_k = -\vartheta_1 g_1 - \ldots - \vartheta_m g_m$

Левая часть этого равенства есть вектор из $\mathbf{R}_1$, правая – из $\mathbf{R}_2$. Таким образом, эта правая часть есть одновременно вектор из $\mathbf{R}_1$ и из $\mathbf{R}_2$. т.е. принадлежит $\mathbf{R}_0$ и, значит, выражается как линейная комбинация базиса $e_1, e_2, \ldots, e_k$ подпространства $\mathbf{R}_0$:

$$-\vartheta_1 g_1 - \ldots - \vartheta_m g_m = c_1 e_1 + c_2 e_2 + \ldots + c_k e_k$$

В силу линейной независимости векторов (6.8) это возможно только тогда, когда все коэффициенты – нули. В частности, $\vartheta_1 = \vartheta_2 = \ldots = \vartheta_m = 0$, т.е. $\lambda_1 f_1 + \ldots + \lambda_l f_l + \mu_1 e_1 + \ldots + \mu_k e_k = 0$.

Из линейной независимости векторов (6.7) получаем, что все коэффициенты $\lambda_1, \lambda_2, \ldots, \lambda_l$, $\mu_1, \ldots, \mu_k$ равны нулю. Таким образом, линейная независимость системы (6.8) доказана.

Покажем теперь, что всякий вектор $x \in R^*$ выражается как линейная комбинация векторов этой системы.

По определению $R^*$ вектор $x$ можно представить в виде $x = x_1 + x_2$, где $x_1 \in R_1$, $x_2 \in R_2$.

Так как $x_1 \in R_1$, то его можно представить как линейную комбинацию векторов (6.6). Аналогично $x_2 \in R_2$, можно представить как линейную комбинацию векторов (6.7). Складывая, получим, что вектор $x$ может быть представлен как линейная комбинация системы (6.8).

Итак, мы получили, что, с одной стороны, векторы $f_1, f_2, \ldots, f_l$ и $e_1, e_2, \ldots, e_k$, а также $g_1, g_2, \ldots, g_m$ линейно независимы и, с другой стороны, всякий вектор из $R^*$ есть их линейная комбинация. Отсюда следует, что эти векторы образуют базис в $R^*$. Итак, мы имеем $k$ векторов (6.5), образующих базис в $R_0$, $k + l$ векторов (6.6), образующих базис в $R^*$, $k + m$ векторов (6.7), образующих базис в $R_2$, $k + l + m$ векторов (6.8), образующих базис в $R^* = R_1 + R_2$. Утверждение теоремы обращается, таким образом, в тождество:

$$(k + l) + (k + m) = (k + l + m) + k.$$

Теорема доказана.

***Пример 6.4.*** В качестве примера, рассмотрим граф, представленный на рис. 6.4.

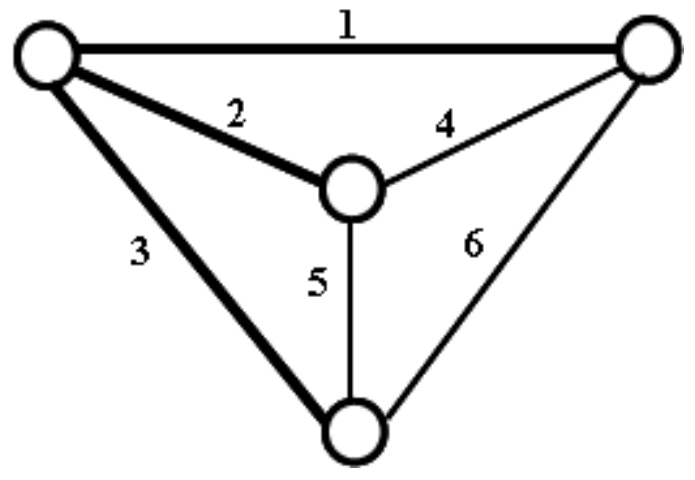

Рис. 6.4. Граф $K_4$.

Система фундаментальных циклов имеет вид:

$c_1^* = (0,0,1,0,1,1) = \{1,2,4\}$;
$c_2^* = (0,1,0,1,1,0) = \{2,3,5\}$;
$c_3^* = (1,0,0,1,0,1) = \{1,3,6\}$.

Система фундаментальных разрезов имеет вид:

$s_1^* = (1,0,1,0,0,1) = \{1,4,6\}$;
$s_2^* = (0,1,1,0,1,0) = \{2,4,5\}$;
$s_3^* = (1,1,0,1,0,0) = \{3,5,6\}$.

Образуем подпространство квазициков $C(K_4)$:

$c_1 = c_1^* = (0,0,1,0,1,1) = (1,2,4\}$;
$c_2 = c_2^* = (0,1,0,1,1,0) = \{2,3,5\}$;
$c_3 = c_3^* = (1,0,0,1,0,1) = \{1,3,6\}$;
$c_4 = c_1^* \oplus c_2^* = (0,1,1,1,0,1) = \{1,3,4,5\}$;
$c_5 = c_1^* \oplus c_3^* = (1,0,1,1,1,0) = \{2,3,4,6\}$;
$c_6 = c_2^* \oplus c_3^* = (1,1,0,0,1,1) = \{1,2,5,6\}$;
$c_7 = c_1^* \oplus c_2^* \oplus c_3^* = (1,1,1,0,0,0) = \{4,5,6\}$.

Образуем подпространство квалиразрезов $S(K_4)$:

$s_1 = s_1^* = (1,0,1,0,0,1) = \{1,4,6\}$;
$s_2 = s_2^* = (0,1,1,0,1,0) = \{2,4,5\}$;
$s_3 = s_3^* = (1,1,0,1,0,0) = \{3,5,6\}$;
$s_4 = s_1^* \oplus s_2^* = (1,1,0,0,1,1) = \{1,2,5,6\}$;
$s_5 = s_1^* \oplus s_3^* = (0,1,1,1,0,1) = \{1,3,4,5\}$;
$s_6 = s_2^* \oplus s_3^* = (1,0,1,1,1,0) = \{2,3,4,6\}$;
$s_7 = s_1^* \oplus s_2^* \oplus s_3^* = (0,0,0,1,1,1) = \{1,2,3\}$.

В данном примере присутствуют общие пересекающиеся суграфы. Образуем подпространство $R^0 = C \cap S$

$R_1^0 = (1,1,0,0,1,1) = \{1,2,5,6\}$;
$R_2^0 = (0,1,1,1,0,1) = \{1,3,4,5\}$;
$R_3^0 = (1,0,1,1,1,0) = \{2,3,4,6\}$.

Данное подпространство имеет базис, состоящий из двух векторов $R_1^0$ и $R_2^0$. Размерность этого подпространства $\dim R^0 = 2$.

Подпространство $R^* = C \oplus S$:

$\alpha_0 = (0,0,0,0,0,0) = c_5 \oplus s_6 = c_6 \oplus s_4 = c_4 \oplus s_5 = \{1,2,5,6\} \oplus \{1,2,5,6\} =$

$= \{1,3,4,5\} \oplus \{1,3,4,5\} = \{2,3,4,6\} \oplus \{2,3,4,6\} = \varnothing$;

$\alpha_7 = (0,0,0,1,1,1) = c_5 \oplus s_1 = c_6 \oplus s_3 = c_4 \oplus s_2 = s_7 = \{1,2,5,6\} \oplus \{3,5,6\} =$

$= \{1,3,4,5\} \oplus \{2,4,5\} = \{2,3,4,6\} \oplus \{1,4,6\} = \{1,2,3\};$

$\alpha_{11} = (0,0,1,0,1,1) = c_2 \oplus s_5 = c_3 \oplus s_6 = c_7 \oplus s_4 = c_1 = \{2,3,5\} \oplus \{1,3,4,5\} =$

$= \{1,3,6\} \oplus \{2,3,4,6\} = \{4,5,6\} \oplus \{1,2,5,6\} = \{1,2,4\};$

$\alpha_{12} = (0,0,1,1,0,0) = c_1 \oplus s_7 = c_2 \oplus s_2 = c_3 \oplus s_1 = c_2 \oplus s_3 = \{1,2,4\} \oplus$

$\oplus \{1,2,3\} = \{2,3,5\} \oplus \{2,4,5\} = \{1,3,6\} \oplus \{1,4,6\} = \{4,5,6\} \oplus \{3,5,6\} = \{3,4\};$

$\alpha_{17} = (0,1,0,0,0,1) = c_1 \oplus s_2 = c_2 \oplus s_7 = c_3 \oplus s_3 = c_7 \oplus s_1 = \{1,2,4\} \oplus \{2,4,5\} =$

$= \{2,3,5\} \oplus \{1,2,3\} = \{1,3,6\} \oplus \{3,5,6\} = \{4,5,6\} \oplus \{1,4,6\} = \{1,5\};$

$\alpha_{22} = (0,1,0,1,1,0) = c_1 \oplus s_5 = c_3 \oplus s_4 = c_7 \oplus s_6 = c_2 = \{1,2,4\} \oplus \{1,3,4,5\} =$

$= \{1,3,6\} \oplus \{1,2,5,6\} = \{4,5,6\} \oplus \{2,3,4,6\} = \{2,3,5\};$

$\alpha_{26} =(0,1,1,0,1,0) = c_6 \oplus s_1 = c_4 \oplus s_7 = c_5 \oplus s_3 = s_2 = \{1,2,5,6\} \oplus \{1,4,6\} =$

$= \{1,3,4,5\} \oplus \{1,2,3\} = \{2,3,4,6\} \oplus \{3,5,6\} = \{2,4,5\};$

$\alpha_{29} = (0,1,1,1,0,1) = c_6 \oplus s_6 = c_4 = s_5 = \{1,2,5,6\} \oplus \{2,3,4,6\} = \{1,3,4,5\};$

$\alpha_{34} = (1,0,0,0,0.1) = c_1 \oplus s_2 = c_2 \oplus s_7 = c_3 \oplus s_3 = c_7 \oplus s_1 = \{1,2,4\} \oplus \{2,4,5\} =$

$= \{2,3,5\} \oplus \{1,2,3\} = \{1,3,6\} \oplus \{3,5,6\} = \{4,5,6\} \oplus \{1,4,6\} = \{1,5\};$

$\alpha_{37} = (1,0,0,1,0,1) = c_1 \oplus s_6 = c_2 \oplus s_4 = c_7 \oplus s_5 = c_3 = \{1,2,4\} \oplus \{2,3,4,6\} =$

$= \{2,3,5\} \oplus \{1,2,5,6\} = \{4,5,6\} \oplus \{1,3,4,5\} = \{1,3,6\};$

$\alpha_{41} = (1,0,1,0,0,1) = c_6 \oplus s_2 = c_4 \oplus s_3 = c_5 \oplus s_7 = s_1 = \{1,2,5,6\}$ □ $\{2,4,5\} =$

$= \{1,3,4,5\} \oplus \{3,5,6\} = \{2,3,4,6\} \oplus \{1,2,3\} = \{1,4,6\};$

$\alpha_{43} = (1,0,1,1,1,0) = c_6 \oplus s_5 = c_5 = s_6 = \{1,2,5,6\} \oplus \{1,3,4,5\} = \{2,3,4,6\};$

$\alpha_{46} = (1,1,0,0,1,1) = c_4 \oplus s_6 = c_5 \oplus s_5 = c_6 = s_4 = \{1,3,4,5\} \oplus \{2,3,4,6\} = \{1,2,5,6\};$

$\alpha_{51} = (1,1,0,1,0,0) = c_6 \oplus s_7 = c_4 \oplus s_1 = c_5 \oplus s_2 = s_3 = \{1,2,5,6\} \oplus \{1,2,3\} =$

$= \{1,3,4,5\} \oplus \{1,4,6\} = \{2,3,4,6\} \oplus \{2,4,5\} = \{3,5,6\};$

$\alpha_{56} = (1,1,1,0,0,0) = c_1 \oplus s_4 = c_2 \oplus s_6 = c_3 \oplus s_5 = c_7 = \{1,2,4\} \oplus \{1,2,5,6\} =$

$= \{2,3,5\} \oplus \{2,3,4,6\} = \{1,3,6\} \oplus \{1,3,4,5\} =\{4,5,6\};$

$\alpha_{63} = (1,1,1,1,1,1) = c_1 \oplus s_3 = c_2 \oplus s_1 = c_3 \oplus s_2 = c_7 \oplus s_7 = \{2,3,5\} \oplus \{1,4,6\} =$

= {1,3,6} ⊕ {2,4,5} = {4,5,6} ⊕ {1,2,3} = {1,2,4} ⊕ {3,5,6} ={1,2,3,4,5,6}.

В качестве базиса подпространства R*, можно выбрать множество $\{\alpha 11, \alpha 12, \alpha 22, \alpha 37\}$.

Размерность подпространства R*:

dim R* = dim C(G) + dim S(G) - dim (C(G) ∩ S(G)) = 3+3-2=4.

## 6.4. Связь между подпространствами циклов и разрезов

Приведем некоторые теоремы, устанавливающие связь между подпространством циклов и подпространством разрезов.

**Теорема 6.4** [34]. Любой суграф в подпространстве циклов графа **G** имеет четное число общих ребер с любым суграфом в подпространстве разрезов того же графа.

**Теорема 6.5** [34]. Суграф графа G принадлежит подпространству циклов графа G, если он имеет четное число общих ребер с любым суграфом в подпространстве разрезов того же графа.

*Доказательство.* Предположим, что G – связный граф. Пусть T – остов графа G. Обозначим ветви остова T как $t_1,t_2,\ldots$,а хорды – как $b_1,b_2,\ldots$ Рассмотрим любой суграф Q графа G, имеющий четное число общих ребер с любым суграфом в подпространстве разрезов графа G. Предположим, не нарушая общности, что суграф *c* содержит хорды $b_1,b_2,\ldots,b_r$. Обозначим как Q′ кольцевую сумму базисных циклов $c_1,c_2,\ldots,c_r$ по отношению к хордам $b_1,b_2,\ldots,b_r$. Очевидно, что Q′ содержит только хорды $b_1,b_2,\ldots,b_r$. Следовательно, Q' ⊕ Q не содержит ни одной хорды. Поскольку Q′ – кольцевая сумма нескольких циклов графа **G**, то это множество имеет четное число общих ребер с каждым суграфом в подпространстве разрезов графа G. Из того, что Q также обладает этим свойством следует, что им обладает и Q' ⊕ Q. Теперь убедимся, что Q' ⊕ Q – пустое множество. Если это не так, то Q' ⊕ Q состоит только из ветвей. Пусть $a_i$ – любая ветвь в Q' ⊕ Q. Тогда $a_i$ – единственное общее ребро между Q' ⊕ Q и базисным разрезающим множеством по отношению к $a_i$. Однако это невозможно, так как Q' ⊕ Q должно иметь четное число общих ребер с любым разрезающим множеством. Таким образом, Q' ⊕ Q должно быть пустым. Другими словами, $Q = Q' = c_1 \oplus c_2 \oplus \ldots \oplus c_r$ и, следовательно, Q принадлежит подпространству циклов графа G.

***Пример 6.5.*** Согласно теореме 6.4, цикл и разрезающее множество связного графа имеют четное число общих ребер. В качестве примера рассмотрим пересечение элементов подпространства квазициклов и подпространства квалиразрезов для графа $G_{22}$ на рис. 6.3.

$c_1 \cap s_1 = \{1,2,6,\} \cap \{1,2\} = \{1,2\};$ $|\{1,2\}| = 2;$
$c_1 \cap s_2 = \{1,2,6\} \cap \{3,5,7\} = \varnothing;$
$c_1 \cap s_3 = \{1,2,6\} \cap \{4,7\} = \varnothing;$
$c_1 \cap s_4 = \{1,2,6\} \cap \{1,5,6\} = \{1,6\};$ $|\{1,6\}| = 2;$
$c_1 \cap s_5 = \{1,2,6\} \cap \{1,2,3,5,7\} = \{1,2\};$ $|\{1,2\}| = 2;$
$c_1 \cap s_6 = \{1,2,6\} \cap \{1,2,4,7\} = \{1,2\};$ $|\{1,2\}| = 2;$
$c_1 \cap s_7 = \{1,2,6\} \cap \{2,5,6\} = \{2,6\};$ $|\{2,6\}| = 2;$
$c_1 \cap s_8 = \{1,2,6\} \cap \{3,4,5\} = \varnothing;$
$c_1 \cap s_9 = \{1,2,6\} \cap \{1,3,6,7\} = \{1,6\};$ $|\{1,6\}| = 2;$
$c_1 \cap s_{10} = \{1,2,6\} \cap \{1,4,5,6,7\} = \{1,6\};$ $|\{1,6\}| = 2;$
$c_1 \cap s_{11} = \{1,2,6\} \cap \{1,2,3,4,5\} = \{1,2\};$ $|\{1,2\}| = 2;$
$c_1 \cap s_{12} = \{1,2,6\} \cap \{2,3,6,7\} = \{2,6\};$ $|\{2,6\}| = 2;$
$c_1 \cap s_{13} = \{1,2,6\} \cap \{2.4,5,6,7\} = \{2,6\};$ $|\{2,6\}| = 2;$
$c_1 \cap s_{14} = \{1,2,6\} \cap \{1,3,4,6\} = \{1,6\};$ $|\{1,6\}| = 2;$
$c_1 \cap s_{15} = \{1,2,6\} \cap \{2,3,4,6\} = \{2,6\};$ $|\{2,6\}| = 2;$
$c_2 \cap s_1 = \{3,5,6,\} \cap \{1,2\} = \varnothing;$
$c_2 \cap s_2 = \{3,5,6\} \cap \{3,5,7\} = \{3,5\};$ $|\{3,5\}| = 2;$
$c_2 \cap s_3 = \{3,5,6\} \cap \{4,7\} = \varnothing;$
$c_2 \cap s_4 = \{3,5,6\} \cap \{1,5,6\} = \{5,6\};$ $|\{5,6\}| = 2;$
$c_2 \cap s_5 = \{3,5,6\} \cap \{1,2,3,5,7\} = \{3,5\};$ $|\{3,5\}| = 2;$
$c_2 \cap s_6 = \{3,5,6\} \cap \{1,2,4,7\} = \varnothing;$
$c_2 \cap s_7 = \{3,5,6\} \cap \{2,5,6\} = \{5,6\};$ $|\{5,6\}| = 2;$
$c_2 \cap s_8 = \{3,5,6\} \cap \{3,4,5\} = \{3,5\};$ $|\{3,5\}| = 2;$
$c_2 \cap s_9 = \{3,5,6\} \cap \{1,3,6,7\} = \{3,6\};$ $|\{3,6\}| = 2;$
$c_2 \cap s_{10} = \{3,5,6\} \cap \{1,4,5,6,7\} = \{5,6\};$ $|\{5,6\}| = 2;$
$c_2 \cap s_{11} = \{3,5,6\} \cap \{1,2,3,4,5\} = \{3,5\};$ $|\{3,5\}| = 2;$
$c_2 \cap s_{12} = \{3,5,6\} \cap \{2,3,6,7\} = \{3,6\};$ $|\{3,6\}| = 2;$
$c_2 \cap s_{13} = \{3,5,6\} \cap \{2,4,5,6,7\} = \{5,6\};$ $|\{5,6\}| = 2;$
$c_2 \cap s_{14} = \{3,5,6\} \cap \{1,3,4,6\} = \{3,6\};$ $|\{3,6\}| = 2;$
$c_2 \cap s_{15} = \{3,5,6\} \cap \{2,3,4,6\} = \{3,6\};$ $|\{3,6\}| = 2;$
$c_3 \cap s_1 = \{3,4,7\} \cap \{1,2\} = \varnothing;$
$c_3 \cap s_2 = \{3,4,7\} \cap \{3,5,7\} = \{3,4\};$ $|\{3,4\}| = 2;$
$c_3 \cap s_3 = \{3,4,7\} \cap \{4,7\} = \{3,4\};$ $|\{3,4\}| = 2;$
$c_3 \cap s_4 = \{3,4.7\} \cap \{1,5,6\} = \varnothing;$
$c_3 \cap s_5 = \{3,4,7\} \cap \{1,2,3,5,7\} = \{3,7\};$ $|\{3,7\}| = 2;$
$c_3 \cap s_6 = \{3,4,7\} \cap \{1,2,4,7\} = \{4,7\};$ $|\{4,7\}| = 2;$
$c_3 \cap s_7 = \{3,4,7\} \cap \{2,5,6\} = \varnothing;$
$c_3 \cap s_8 = \{3,4,7) \cap (3,4,5\} = \{3,4\};$ $|\{3,4\}| = 2;$
$c_3 \cap s_9 = \{3,4,7\} \cap \{1,3,6,7\} = \{3,7\};$ $|\{3,7\}| = 2;$
$c_3 \cap s_{10} = \{3,4,7\} \cap \{1,4,5,6,7\} = \{4,7\};$ $|\{4,7\} = 2;$
$c_3 \cap s_{11} = \{3,4,7\} \cap \{1,2,3,4,5\} = \{3,4\};$ $|\{3,4\}| = 2;$
$c_3 \cap s_{12} = \{3,4,7\} \cap \{2,3,6,7\} = \{3,7\};$ $|\{3,7\}| = 2;$
$c_3 \cap s_{13} = \{3,4,7\} \cap \{2,4,5,6,7\} = \{4,7\};$ $|\{4,7\}| = 2;$
$c_3 \cap s_{14} = \{3,4,7\} \cap \{1,3,4,6\} = \{3,4\};$ $|\{3,4\}| = 2;$
$c_3 \cap s_{15} = \{3,4,7\} \cap \{2,3,4,6\} = \{3,4\};$ $|\{3,4\}| = 2;$
$c_4 \cap s_1 = \{1,2,3,5\} \cap \{1,2\} = \{1,2\};$ $|\{1,2\}| = 2;$
$c_4 \cap s_2 = \{1,2,3,5\} \cap \{3,5,7\} = \{3,5\};$ $|\{3,5\}| = 2;$
$c_4 \cap s_3 = \{1,2,3,5\} \cap \{4,7\} = \varnothing;$
$c_4 \cap s_4 = \{1,2,3,5\} \cap \{1,5,6\} = \{1,5\};$ $|\{1,5\}| = 2;$
$c_4 \cap s_5 = \{1,2,3,5\} \cap \{1,2,3,5,7\} = \{1,2,3,5\};$ $|\{1,2,3,5\}| = 4;$
$c_4 \cap s_6 = \{1,2,3,5\} \cap \{1,2,4,7\} = \{1,2\};$ $|\{1,2\}| = 2;$
$c_4 \cap s_7 = \{1,2,3,5\} \cap \{2,5,6\} = \{2,5\};$ $|\{2,5\}| = 2;$

$c_4 \cap s_8 = \{1,2,3,5\} \cap \{3,4,5\}= \{3,5\};$ $|\{3,5\}| = 2;$
$c_4 \cap s_9 = \{1,2,3,5\} \cap \{1,3,6,7\}=\{1,3\};$ $|\{1,3\}| = 2;$
$c_4 \cap s_{10} = \{1,2,3,5\} \cap \{1,4,5,6,7\} = \{1,5\};$ $|\{1,5\}| =2 ;$
$c_4 \cap s_{11} = \{1,2,3,5\} \cap \{1,2,3,4,5\} = \{1,2,3,5\};$ $|\{1,2,3,5\}| = 4;$
$c_4 \cap s_{12} = \{1,2,3,5\} \cap \{2,3,6,7\}= \{2,3\};$ $|\{2,3\}| = 2;$
$c_4 \cap s_{13} = \{1,2,3,5\} \cap \{2,4,5,6,7\} = \{2,5\};$ $|\{2,5\}| = 2;$
$c_4 \cap s_{14} = \{1,2,3,5\} \cap \{1,3,4,6\}= \{1,3\};$ $|\{1,3\}| = 2;$
$c_4 \cap s_{15} = \{1,2,3,5\} \cap \{2,3,4,6\}= \{2,3\};$ $|\{2,3\}| = 2;$
$c_5 \cap s_1 = \{1,2,3,4,6,7\} \cap \{1,2\} = \{1,2\};$ $|\{1,2\}| = 2;$
$c_5 \cap s_2 = \{1,2,3,4,6,7\} \cap \{3,5,7\} = \{3,7\};$ $|\{3,7\}| = 2;$
$c_5 \cap s_3 = \{1,2,3,4,6,7\} \cap \{4,7\} = \{4,7\};$ $|\{4,7\}| = 2;$
$c_5 \cap s_4 = \{1,2,3,4,6,7\} \cap \{1,5,6\} = \{1,6\};$ $|\{1,6\}| = 2;$
$c_5 \cap s_5 = \{1,2,3,4,6,7\} \cap \{1,2,3,5,7\} = \{1,2,3,7\};$ $|\{1,2,3,7\}| = 4;$
$c_5 \cap s_6 = \{1,2,3,4,6,7\} \cap \{1,2,4,7\} = \{1,2,4,7\};$ $|\{1,2,4,7\}| = 4;$
$c_5 \cap s_7 = \{1,2,3,4,6,7\} \cap \{2,5,6\} = \{2,6\};$ $|\{2,6\}| = 2;$
$c_5 \cap s_8 = \{1,2,3,4,6,7\} \cap \{3,4,5\} = \{3,4\};$ $|\{3,4\}| = 2;$
$c_5 \cap s_9 = \{1,2,3,4,6,7\} \cap \{1,3,6,7\} = \{1,3,6,7\};$ $|\{1,3,6,7\}| = 4;$
$c_5 \cap s_{10} = \{1,2,3,4,6,7\} \cap \{1,4,5,6,7\} = \{1,4,6,7\};$ $|\{1,4,6,7\}| = 4;$
$c_5 \cap s_{11} = \{1,2,3,4,6,7\} \cap \{1,2,3,4,5\} = \{1,2,3,4\};$ $|\{1,2,3,4\}| = 4;$
$c_5 \cap s_{12} = \{1,2.3,4,6.7\} \cap \{2,3,6,7\} = \{2,3,6,7\};$ $|\{2,3,6,7\}| = 4;$
$c_5 \cap s_{13} = \{1,2,3,4.6,7\} \cap \{2,4,5,6,7\} = \{2,4,6,7\};$ $|\{2,4,6,7\}| = 4;$
$c_5 \cap s_{14} = \{1,2,3,4,6,7\} \cap \{1,3,4,6\} = \{1,3,4,6\};$ $|\{1,3,4,6\}| = 4;$
$c_5 \cap s_{15} = \{1,2,3,4,6,7\} \cap \{2,3,4,6\} = \{2.3,4,6\};$ $|\{2,3,4,6\}| = 4;$
$c_6 \cap s_1 = \{4,5,6,7\} \cap \{1,2\} = \varnothing;$
$c_6 \cap s_2 = \{4,5,6,7\} \cap \{3,5,7\} = \{5,7\};$ $|\{3,5\}| = 2;$
$c_6 \cap s_3 = \{4,5,6,7\} \cap \{4,7\} = \{4,7\};$ $|\{4,7\}| = 2;$
$c_6 \cap s_4 = \{4,5,6,7\} \cap \{1,5,6\} = \{5,6\};$ $|\{5,6\}| = 2;$
$c_6 \cap s_5 = \{4,5,6,7\} \cap \{1,2,3,5,7\} = \{5,7\};$ $|\{5,7\}| = 2;$
$c_6 \cap s_6 = \{4,5,6,7\} \cap \{1,2,4,7\} = \{4,7\};$ $|\{4,7\}| = 2;$
$c_6 \cap s_7 = \{4,5,6,7\} \cap \{2,5,6\} = \{5,6\};$ $|\{5,6\}| = 2;$
$c_6 \cap s_8 = \{4,5,6,7\} \cap \{3,4,5\} = \{4,5\};$ $|\{4,5\}| = 2;$
$c_6 \cap s_9 = \{4,5,6,7\} \cap \{1,3,6,7\} = \{6,7\};$ $|\{6,7\}| = 2;$
$c_6 \cap s_{10} = \{4,5,6,7\} \cap \{1,4,5,6,7\} = \{4,5,6,7\};$ $|\{4,5,6,7\}| = 4;$
$c_6 \cap s_{11} = \{4,5,6,7\} \cap \{1,2,3,4,5\} = \{4,5\};$ $|\{4,5\}| = 2;$
$c_6 \cap s_{12} = \{4,5,6,7\} \cap \{2,3,6,7\} = \{6,7\};$ $|\{6,7\}| = 2;$
$c_6 \cap s_{13} = \{4,5,6,7\} \cap \{2,4,5,6,7\} = \{4,5,6,7\};$ $|\{4,5,6,7\}| = 4;$
$c_6 \cap s_{14} = \{4,5,6,7\} \cap \{1,3,4,6\}= \{4,6\};$ $|\{4,6\}| = 2;$
$c_6 \cap s_{15} = \{4,5,6,7\} \cap \{2,3,4,6\} = \{4,6\};$ $|\{4,6\}| = 2;$
$c_7 \cap s_1 = \{1,2,4,5,7\} \cap \{1,2\}=\{1,2\};$ $|\{1,2\}|=2;$
$c_7 \cap s_2 = \{1,2,4,5.7\} \cap \{3,5,7\}= \{5,7\};$ $|\{5,7\}| = 2;$
$c_7 \cap s_3 = \{1,2,4.5,7\} \cap \{4,7\} = \{4,7\};$ $|\{4,7\}| = 2;$
$c_7 \cap s_4 = \{1,2,4,5,7\} \cap \{1,5,6\}= \{1,5\};$ $|\{1,5\}| = 2;$
$c_7 \cap s_5 = \{1,2,4,5,7\} \cap \{1,2,3,5,7\} =\{1,2,5,7\};$ $|\{1,2,5,7\}| = 4;$
$c_7 \cap s_6 = \{1,2,4,5,7\} \cap \{1,2,4,7\} = \{1,2,4,7\};$ $|\{1,2,4,7\}| = 4;$
$c_7 \cap s_7 = \{1,2,4,5,7\} \cap \{2,5,6\} = \{2,5\};$ $|\{2,5\}| = 2;$
$c_7 \cap s_8 = \{1,2,4,5,7\} \cap \{3,4,5\} = \{4,5\};$ $|\{4,5\}| = 2;$
$c_7 \cap s_9 = \{1,2,4,5,7\} \cap \{1,3,6,7\} = \{1,7\};$ $|\{1,7\}| = 2;$
$c_7 \cap s_{10} = \{1,2,4,5,7\} \cap \{1,4,5,6,7\} = \{1,4,5,7\};$ $|\{1,4,5,7\}| = 4;$
$c_7 \cap s_{11} = \{1,2,4,5,7\} \cap \{1,2,3,4,5\}= \{1,2,4,5\};$ $|\{1,2,4,5\}| = 4;$
$c_7 \cap s_{12} = \{1,2,4,5,7\} \cap \{2,3,6,7\} = \{2,7\};$ $|\{2,7\}| = 2;$
$c_7 \cap s_{13} = \{1,2,4,5,7\} \cap \{2,4,5,6,7\} = \{2,4,5,7\};$ $|\{2,4,5,7\}| = 4;$
$c_7 \cap s_{14} = \{1,2,4,5,7\} \cap \{1,3,4,6\} = \{1,4\};$ $|\{1,4\}| = 2;$

$c_7 \cap s_{15} = \{1,2,4,5,7\} \cap \{2,3,4,6\} = \{2,4\}$;   $|\{2,4\}| = 2$.

Как видим из примера, пересечение цикла и разрезающих множеств связного графа имеют четное число общих ребер. Общие пересекающиеся суграфы отсутствуют.

## 6.5. Ортогональность подпространств циклов и разрезов

Каждое $n$-мерное векторное пространство над полем F изоморфно векторному пространству всех $n$-векторов над тем же полем. Следовательно, векторное пространство $£_G$ графа G изоморфно векторному пространству всех $m$-векторов над полем GF(2), где $m$ – число ребер графа G.

Пусть $u_1, u_2, ..., u_m$ – ребра графа G. Предположим, что мы сопоставили каждому реберно-порожденному суграфу $G_i$, графа G такой $m$-вектор $w_i$, что j-й элемент $w_i$ равен 1 тогда и только тогда, когда ребро $u_j$ принадлежит суграфу $G_i$. Тогда кольцевая сумма $G_i \oplus \xi\ G_j$ двух суграфов $G_i$ и $G_j$ будет соответствовать $m$-вектору $w_i + w_j$, являющемуся суммой по mod 2 векторов $w_i$, и $w_j$. Легко заметить, что описанное соответствие действительно определяет изоморфизм между $£_G$ и векторным пространством всех $m$-векторов над полем GF(2). В самом деле, если мы выберем $\{u_1\},\{u_2\},...,\{u_m\}$ в качестве базисных векторов для пространства $£_G$, то элементами $w_i$ будут координаты $G_i$, связанные с этим базисом.

При определении этого изоморфизма мы опять использовали символ $£_G$ для обозначения векторного пространства всех $m$-векторов, сопоставленных суграфам графа G. Пусть C(G) обозначает подпространство $m$-векторов, представляющих суграфы в подпространстве циклов графа G, а S(G) – подпространство, представляющее суграфы в подпространстве разрезов графа G.

**Определение 6.3.** *Скалярным произведением* $w_1 = (\alpha_1, \alpha_2, ..., \alpha_m)$ и $w_2 = (\beta_1, \beta_2, ..., \beta_m)$, обозначаемым $<w_1, w_2>$, называется скаляр, определяемый выражением $<w_1, w_2> = \alpha_1 \cdot \beta_1 + \alpha_2 \cdot \beta_2 + ... + \alpha_m \cdot \beta_m$.

Например, если $w_1 = (0,1,0,0,1)$ и $w_2 = (1,0,1,1,1)$, то $<w_1, w_2> = (0 \cdot 1 + 1 \cdot 0 + 0 \cdot 1 + 0 \cdot 1 + 1 \cdot 1) = 0+0+0+0+1 = 1$. Векторы $w_1$ и $w_2$ ортогональны друг другу, если $<w_1, w_2> = 0$, где 0 – аддитивный нулевой элемент поля GF(2). Например, векторы $w_1 = (1,1,0,1,1)$ и $w_2 = (1,1,1,0,0)$ ортогональны над полем GF(2), поскольку $<w_1, w_2> = (1 \cdot 1 + 1 \cdot 1 + 0 \cdot 1 + 1 \cdot 0 + 1 \cdot 0) = 1+1+0+0+0 = 0 \pmod 2$.

Два подпространства $W_1$ и $W_2$ пространства W являются ортогональными подпространствами этого же пространства, если каждый вектор одного пространства ортогонален любому вектору другого подпространства. Два подпространства $W_1$ и $W_2$ пространства W называются ортогональными дополнениями пространства W, если они ортогональны друг другу и их прямая сумма равна векторному пространству W.

Рассмотрим два таких вектора $w_i$ и $w_j$, что вектор $w_i$ находится в пространстве C(G), а вектор $w_j$ – в пространстве S(G). Из того факта, что любой подграф в пространстве C(G) имеет четное число общих ребер с произвольным суграфом в пространстве S(G) следует, что скалярное произведение $<w_i, w_j>$ векторов $w_i$ и $w_j$ равно сумме по mod 2 четного числа единиц. Это означает, что $<w_i, w_j> = 0$. Иначе говоря, $m$-векторы в пространстве C(G) ортогональны подобным векторам в пространстве S(G). Таким образом, имеет место следующая теорема:

**Теорема 6.6** [34]. Подпространства циклов и разрезов графа ортогональны.

Рассмотрим прямую сумму C(G) + S(G). Мы знаем, что dim(C(G) + S(G)) = dim(C(G)) + + dim(S(G)) - dim(C(G) $\cap$ S(G)). Поскольку dim(C(G)) + dim(S(G)) = $m$, получаем, что dim(C(G) + S(G)) = $m$ - dim(C(G) $\cap$ S(G)). Теперь ортогональные подпространства C(G) и S(G) будут также и ортогональными дополнениями $£_G$ тогда и только тогда, когда dim(C(G) + S(G)) = $m$. Иными словами, C(G) и S(G) будут ортогональными дополнениями в том и только в том случае, если dim(C(G) $\cap$ S(G)) = 0, т.е. C(G) $\cap$ S(G) – нулевой вектор (все элементы которого равны нулю). Поэтому мы получаем следующую теорему:

**Теорема 6.7** [34]. Подпространства C(G) и S(G) циклов и разрезов графа являются ортогональными дополнениями тогда и только тогда, когда C(G) $\cap$ S(G) – нулевой вектор.

Пусть C(G) и S(G) – ортогональные дополнения. Это означает, что каждый вектор в пространстве $£_G$ можно представить кольцевой суммой $w_i + w_j$, где вектор $w_i$, принадлежит пространству C(G), а вектор $w_j$ – пространству S(G). Другими словами, каждый суграф графа G можно представить кольцевой суммой двух суграфов, один из которых принадлежит подпространству циклов, а другой – подпространству разрезов. В частности, сам граф G можно представить таким же образом.

Предположим, что C(G) и S(G) не являются ортогональными дополнениями. Тогда, очевидно, существует такой суграф, который нельзя представить как кольцевую сумму суграфов в пространствах C(G) и S(G). Возникает вопрос: можно ли в этом случае представить граф **G** кольцевой суммой суграфов, принадлежащих подпространствам C(G) и S(G)? Ответом является следующая теорема:

**Теорема 6.8.** Любой граф G можно представить в виде кольцевой суммы двух суграфов, один из которых принадлежит подпространству циклов, а другой – подпространству разрезов графа G. Доказательство этой теоремы можно найти в работах

[30,43,50].

***Пример 6.5.*** В графе G условие теоремы 6.8 выполняется

$s_6 \oplus c_2 = \{1,2,4,7\} \oplus \{3,5,6\} = \{1,2,3,4,5,6,7\}$.

Для графа $K_4$ условме теоремы 6.8 также выполнимо

$c_3^* \oplus s_2^* = \{1,3,6\} \oplus \{2,4,5\} = \{1,2,3,4,5,6\}$.

## 6.6. Изоморфизм суграфов графа

Пусть имеется два графа G и H с изоморфными подпространства разрезов S(G) и S(H) и изоморфными подпространствами циклов C(G) и C(H). Согласно алгебраической теории графов, пространства изоморфны, если одинаковы размеры этих пространств [4].

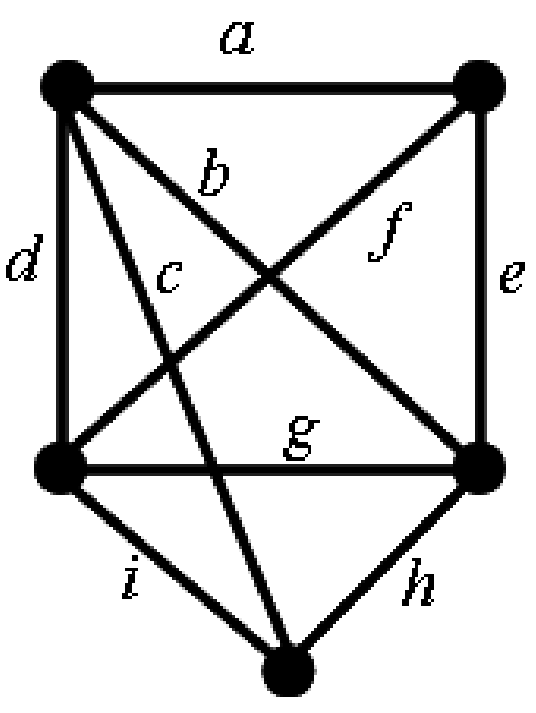

Рис. 6.5. Граф.

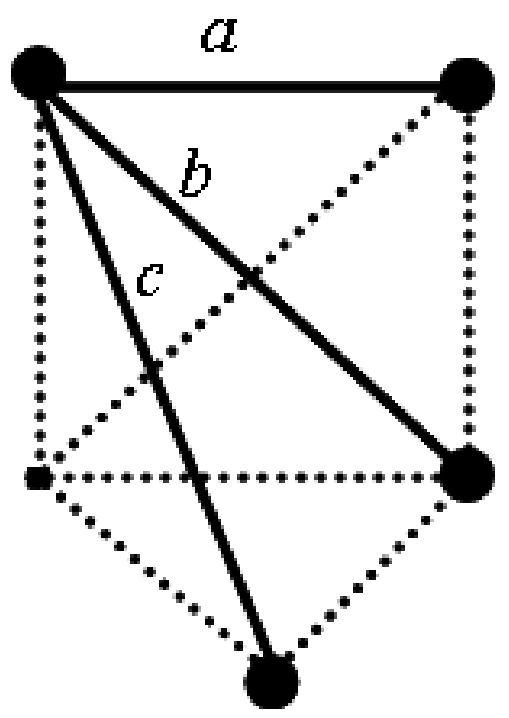

Рис. 6.6. Суграф {a,b,c}.

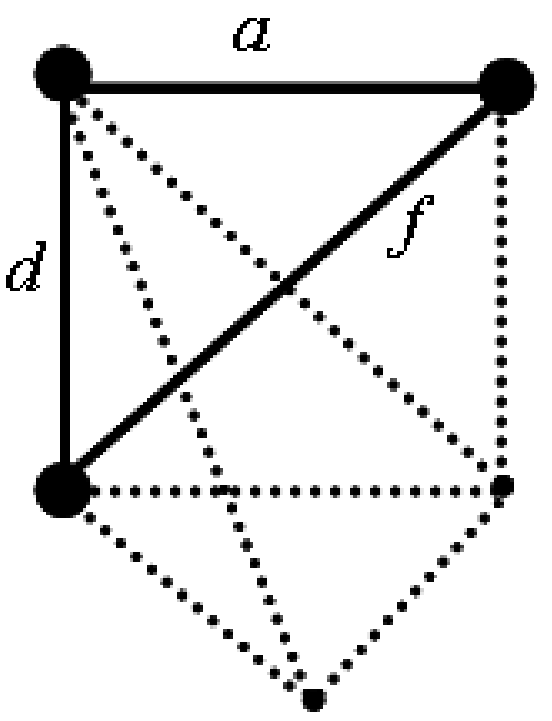

Рис. 6.7. Суграф {a,d,f}.

Изоморфизм суграфов определяется не только составом входящих в него ребер, но и всем подмножеством вершин суграфа графа.

Ребро графа может быть представлено в множественном виде как базовый реберный разрез графа. В этом случае реберный граф L(G) следует рассматривать как образ графа G.

Заметим, что подпространства S (или подпространства C) графов G и H, могут быть изоморфны, а множество суграфов пространства суграфов графа – неизоморфно.

Рассмотрим графы $G_{15}$ и $G_{16}$ (рис. 6.8 и 6.9) или (рис. 4.9 и 4.10).

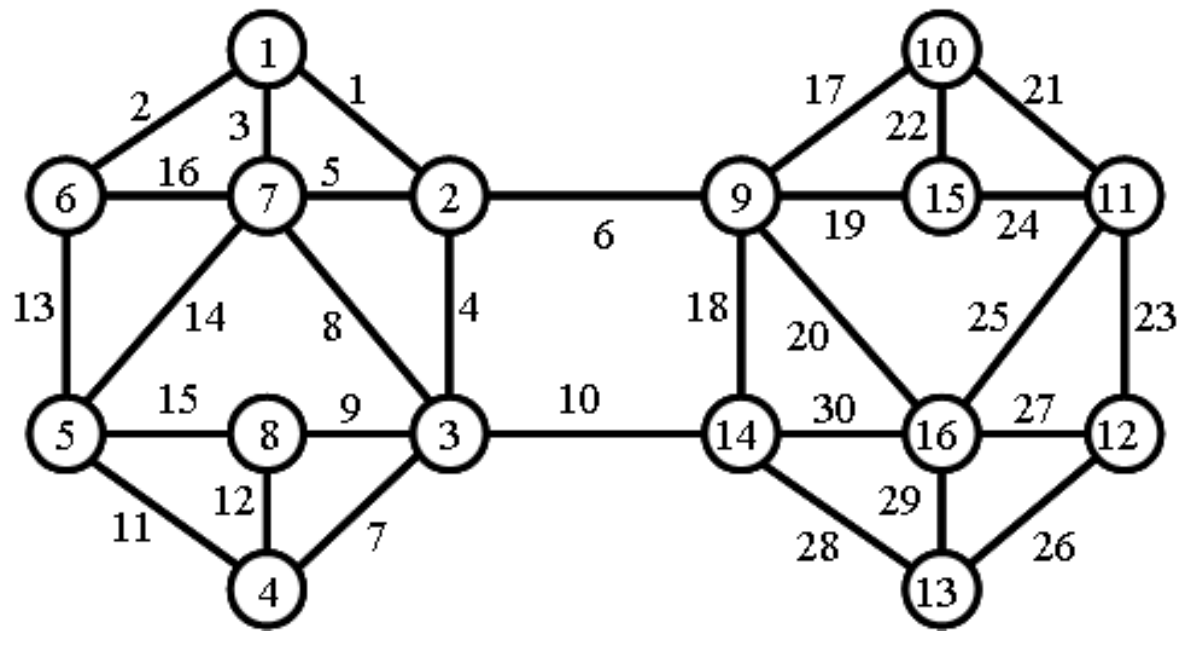

Рис. 6.8. Граф $G_{15}$.

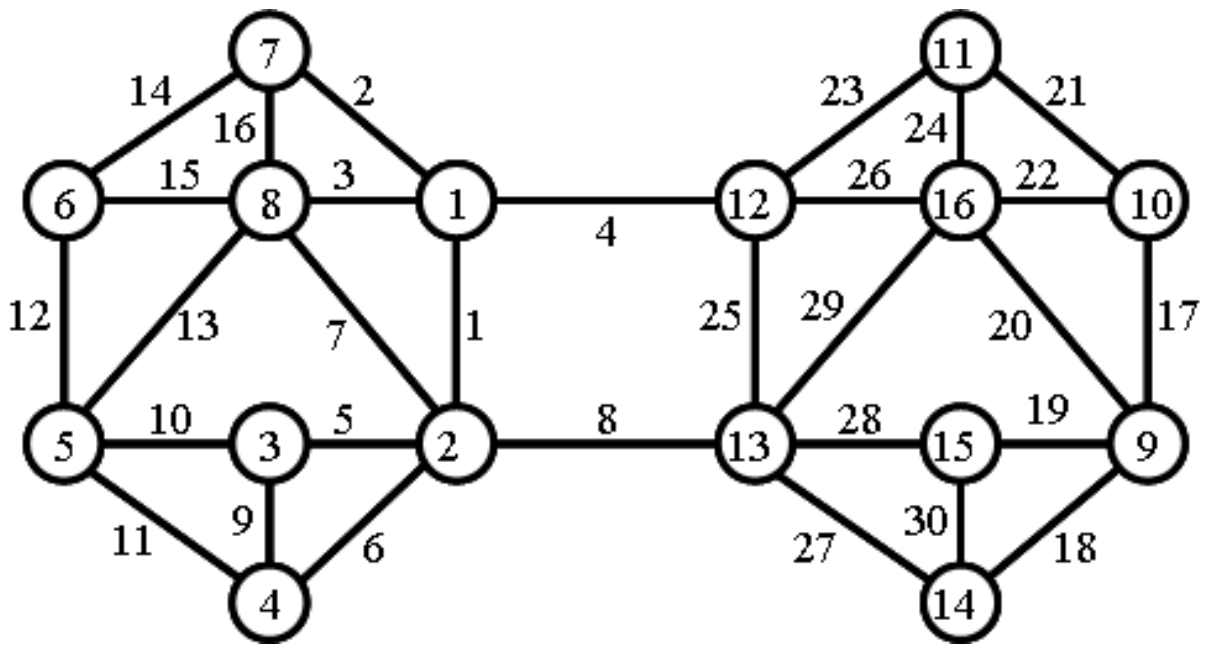

Рис. 6.9. Граф $G_{16}$.

Покажем, что подпространства циклов графов $G_{15}$ и $G_{16}$ изоморфны, но, графы $G_{15}$ и

$G_{16}$ не изоморфны, так как отличаются значениями весов элементов графов в векторных инвариантах спектра реберных разрезов.

Вектор весов ребер графа $G_{15}$: $F_w(\xi_0(G_{15})) =$

$= (4,4,12,12,12,12,13,13,14,14,14,14,15,15,15,15,18,18,18,18,19,19,19,19,19,19,21,21,23,23) =$
$= (2\times 4,4\times 12,2\times 13,4\times 14,4\times 15,4\times 18,6\times 19,2\times 21,2\times 23)$;

Вектор весов вершин графа $G_{15}$: $F_w(\zeta_0(G_{15})) =$

$= (37,37,37,37,37,37,45,45,55,55,80,80,81,81,100,100) =$
$= (6\times 37,2\times 45,2\times 55,2\times 80,2\times 81,2\times 100)$.

Вектор весов ребер графа $G_{16}$: $F_w(\xi_0(G_{16})) =$

$= (4,4,12,12,12,12,13,13,14,14,14,14,15,15,15,15,16,17,17,18,18,19,19,20,20,20,20,21,21,28) =$
$= (2\times 4,4\times 12,2\times 13,4\times 14,4\times 15,1\times 16,2\times 17,2\times 18,2\times 19,4\times 20,2\times 21,1\times 28)$;

Вектор весов вершин графа $G_{16}$: $F_w(\zeta_0(G_{16})) =$

$= (37,37,38,38,38,38,44,44,55,55,72,72,80,80,108,108) =$
$= (2\times 37,4\times 38,2\times 44,2\times 55,2\times 72,2\times 80,2\times 108)$.

Так как веса элементов суграфов базовых реберных разрезов $W_0(G_{15})$ и $W_0(G_{16})$ графов не равны, следовательно, не существует бинарного соответствия между суграфами базовых реберных разрезов графов $G_{15}$ и $G_{16}$.

Количество изометрических циклов в графах $G_{15}$ и $G_{16}$ = 17

$IC(G_{15}) = F_\tau(\xi_0(G_{15})\,\&\,F_\tau(\zeta_0(G_{15}) =$
$= (4,4,4,4,4,4,4,4,5,5,5,5,5,5,5,5,5,5,11,11,12,12,13,13,13,13,15,15,17,17)\,\&$
$\&\,(14,14,14,14,28,28,30,30,34,34,36,36,44,44,44,44)$.

$IC(G_{16}) = F_\tau(\xi_0(G_{16})\,\&\,F_\tau(\zeta_0(G_{16}) =$
$= (4,4,4,4,4,4,4,4,5,5,5,5,5,5,5,5,5,5,11,11,12,12,13,13,13,13,15,15,17,17)\,\&$
$\&\,(14,14,14,14,28,28,30,30,34,34,36,36,44,44,44,44)$.

Так как веса элементов (ребер и вершин) множества суграфов базовых реберных циклов $T_0(G_{15})$ и $T_0(G_{16})$ графов равны и определены, и проверка устанавливает равенство весов изометрических циклов графа, следовательно существует бинарное соответствие между суграфами базовых реберных циклов. В свою очередь, во множестве изометрических циклов существует базис, состоящий из изометрических циклов. Тогда суграфы подпространства циклов будут определены как линейная комбинация базисов. Следовательно, суграфы подпространства циклов изоморфны. Если нельзя установить равенство весов изометрических циклов, то суграфы подпространства циклов графов $G_{15}$ и $G_{16}$ не изоморфны.

Перейдем к рассмотрению реберных графов $L(G_{15})$ и $L(G_{16})$. Будем рассматривать составные части реберных графов $L(G_{15})$ и $L(G_{16})$.

Базовые реберные разрезы графа $G_{15}$

$w_0(e_1) = \{e_2,e_3,e_4,e_5,e_6\}$;
$w_0(e_2) = \{e_1,e_3,e_{13},e_{16}\}$;
$w_0(e_3) = \{e_1,e_2,e_5,e_8,e_{14},e_{16}\}$;
$w_0(e_4) = \{e_1,e_5,e_6,e_7,e_8,e_9,e_{10}\}$;
$w_0(e_5) = \{e_1,e_3,e_4,e_6,e_8,e_{14},e_{16}\}$;
$w_0(e_6) = \{e_1,e_4,e_5,e_{17},e_{18},e_{19},e_{20}\}$;
$w_0(e_7) = \{e_4,e_8,e_9,e_{10},e_{11},e_{12}\}$;
$w_0(e_8) = \{e_3,e_4,e_5,e_7,e_9,e_{10},e_{14},e_{16}\}$;
$w_0(e_9) = \{e_4,e_7,e_8,e_{10},e_{12},e_{15}\}$;
$w_0(e_{10}) = \{e_4,e_7,e_8,e_9,e_{18},e_{28},e_{30}\}$;
$w_0(e_{11}) = \{e_7,e_{12},e_{13},e_{14},e_{15}\}$;
$w_0(e_{12}) = \{e_7,e_9,e_{11},e_{15}\}$;
$w_0(e_{13}) = \{e_2,e_{11},e_{14},e_{15},e_{16}\}$;
$w_0(e_{14}) = \{e_3,e_5,e_8,e_{11},e_{13},e_{15},e_{16}\}$;
$w_0(e_{15}) = \{e_9,e_{11},e_{12},e_{13},e_{14}\}$;
$w_0(e_{16}) = \{e_2,e_3,e_5,e_8,e_{13},e_{14}\}$;
$w_0(e_{17}) = \{e_6,e_{18},e_{19},e_{20},e_{21},e_{22}\}$;
$w_0(e_{18}) = \{e_6,e_{10},e_{17},e_{19},e_{20},e_{28},e_{30}\}$;
$w_0(e_{19}) = \{e_6,e_{17},e_{18},e_{20},e_{22},e_{24}\}$;
$w_0(e_{20}) = \{e_6,e_{17},e_{18},e_{19},e_{25},e_{27},e_{29},e_{30}\}$;
$w_0(e_{21}) = \{e_{17},e_{22},e_{23},e_{24},e_{25}\}$;
$w_0(e_{22}) = \{e_{17},e_{19},e_{21},e_{24}\}$;
$w_0(e_{23}) = \{e_{21},e_{24},e_{25},e_{26},e_{27}\}$;
$w_0(e_{24}) = \{e_{19},e_{21},e_{22},e_{23},e_{25}\}$;
$w_0(e_{25}) = \{e_{20},e_{21},e_{23},e_{24},e_{27},e_{29},e_{30}\}$;
$w_0(e_{26}) = \{e_{23},e_{27},e_{28},e_{29}\}$;
$w_0(e_{27}) = \{e_{20},e_{23},e_{25},e_{26},e_{29},e_{30}\}$;
$w_0(e_{28}) = \{e_{10},e_{18},e_{26},e_{29},e_{30}\}$;
$w_0(e_{29}) = \{e_{20},e_{25},e_{26},e_{27},e_{28},e_{30}\}$;
$w_0(e_{30}) = \{e_{10},e_{18},e_{20},e_{25},e_{27},e_{28},e_{29}\}$.

Базовые реберные разрезы графа $G_{16}$

$w_0(e_1) = \{e_2,e_3,e_4,e_5,e_6,e_7,e_8\}$;
$w_0(e_2) = \{e_1,e_3,e_4,e_{14},e_{16}\}$;
$w_0(e_3) = \{e_1,e_2,e_4,e_7,e_{13},e_{15},e_{16}\}$;
$w_0(e_4) = \{e_1,e_2,e_3,e_{23},e_{25},e_{26}\}$;
$w_0(e_5) = \{e_1,e_6,e_7,e_8,e_9,e_{10}\}$;
$w_0(e_6) = \{e_1,e_5,e_7,e_8,e_9,e_{11}\}$;
$w_0(e_7) = \{e_1,e_3,e_5,e_6,e_8,e_{13},e_{15},e_{16}\}$;
$w_0(e_8) = \{e_1,e_5,e_6,e_7,e_{25},e_{27},e_{28},e_{29}\}$;
$w_0(e_9) = \{e_5,e_6,e_{10},e_{11}\}$;
$w_0(e_{10}) = \{e_5,e_9,e_{11},e_{12},e_{13}\}$;
$w_0(e_{11}) = \{e_6,e_9,e_{10},e_{12},e_{13}\}$;
$w_0(e_{12}) = \{e_{10},e_{11},e_{13},e_{14},e_{15}\}$;
$w_0(e_{13}) = \{e_3,e_7,e_{10},e_{11},e_{12},e_{15},e_{16}\}$;
$w_0(e_{14}) = \{e_2,e_{12},e_{15},e_{16}\}$;
$w_0(e_{15}) = \{e_3,e_7,e_{12},e_{13},e_{14},e_{16}\}$;
$w_0(e_{16}) = \{e_2,e_3,e_7,e_{13},e_{14},e_{15}\}$;
$w_0(e_{17}) = \{e_{18},e_{19},e_{20},e_{21},e_{22}\}$;
$w_0(e_{18}) = \{e_{17},e_{19},e_{20},e_{27},e_{30}\}$;
$w_0(e_{19}) = \{e_{17},e_{18},e_{20},e_{28},e_{30}\}$;
$w_0(e_{20}) = \{e_{17},e_{18},e_{19},e_{22},e_{24},e_{26},e_{29}\}$;
$w_0(e_{21}) = \{e_{17},e_{22},e_{23},e_{24}\}$;
$w_0(e_{22}) = \{e_{17},e_{20},e_{21},e_{24},e_{26},e_{29}\}$;
$w_0(e_{23}) = \{e_4,e_{21},e_{24},e_{25},e_{26}\}$;
$w_0(e_{24}) = \{e_{20},e_{21},e_{22},e_{23},e_{26},e_{29}\}$;
$w_0(e_{25}) = \{e_4,e_8,e_{23},e_{26},e_{27},e_{28},e_{29}\}$;
$w_0(e_{26}) = \{e_4,e_{20},e_{22},e_{23},e_{24},e_{25},e_{29}\}$;
$w_0(e_{27}) = \{e_8,e_{18},e_{25},e_{28},e_{29},e_{30}\}$;
$w_0(e_{28}) = \{e_8,e_{19},e_{25},e_{27},e_{29},e_{30}\}$;
$w_0(e_{29}) = \{e_8,e_{20},e_{22},e_{24},e_{25},e_{26},e_{27},e_{28}\}$;
$w_0(e_{30}) = \{e_{18},e_{19},e_{27},e_{28}\}$.

Изометрические циклы графа $G_{15}$

$c_1 = \{e_1,e_3,e_5\}$;
$c_2 = \{e_2,e_3,e_{16}\}$;
$c_3 = \{e_4,e_5,e_8\}$;
$c_4 = \{e_4,e_6,e_{10},e_{18}\}$;
$c_5 = \{e_7,e_8,e_{11},e_{14}\}$;
$c_6 = \{e_7,e_9,e_{12}\}$;
$c_7 = \{e_8,e_9,e_{14},e_{15}\}$;
$c_8 = \{e_{11},e_{12},e_{15}\}$;
$c_9 = \{e_{13},e_{14},e_{16}\}$;
$c_{10} = \{e_{17},e_{19},e_{22}\}$;
$c_{11} = \{e_{17},e_{20},e_{21},e_{25}\}$;
$c_{12} = \{e_{18},e_{20},e_{30}\}$;
$c_{13} = \{e_{19},e_{20},e_{24},e_{25}\}$;
$c_{14} = \{e_{21},e_{22},e_{24}\}$;
$c_{15} = \{e_{23},e_{25},e_{27}\}$;
$c_{16} = \{e_{26},e_{27},e_{29}\}$;
$c_{17} = \{e_{28},e_{29},e_{30}\}$.

Изометрические циклы графа $G_{16}$

$c_1 = \{e_1,e_3,e_7\}$;
$c_2 = \{e_1,e_4,e_8,e_{25}\}$;
$c_3 = \{e_2,e_3,e_{16}\}$;
$c_4 = \{e_5,e_6,e_9\}$;
$c_5 = \{e_5,e_7,e_{10},e_{13}\}$;
$c_6 = \{e_6,e_7,e_{11},e_{13}\}$;
$c_7 = \{e_9,e_{10},e_{11}\}$;
$c_8 = \{e_{12},e_{13},e_{15}\}$;
$c_9 = \{e_{14},e_{15},e_{16}\}$;
$c_{10} = \{e_{17},e_{20},e_{22}\}$;
$c_{11} = \{e_{18},e_{19},e_{30}\}$;
$c_{12} = \{e_{18},e_{20},e_{27},e_{29}\}$;
$c_{13} = \{e_{19},e_{20},e_{28},e_{29}\}$;
$c_{14} = \{e_{21},e_{22},e_{24}\}$;
$c_{15} = \{e_{23},e_{24},e_{26}\}$;
$c_{16} = \{e_{25},e_{26},e_{29}\}$;
$c_{17} = \{e_{27},e_{28},e_{30}\}$.

Количество базовых реберных разрезов графов $G_{15}$ и $G_{16}$ = 30.

Количество вершин в реберных графах $L(G_{15})$ и $L(G_{16})$ = 30.

Количество ребер в реберных графах $L(G_{15})$ и $L(G_{16})$= 88.

Количество изометрических циклов в реберных графах $L(G_{15})$ и $L(G_{16})$= 89.

Рассмотрим дубль-циклы графов $G_{15}$ и $G_{16}$.

| Дубль-циклы в графе $G_{15}$ | Дубль-циклы в графе $G_{16}$ |
|---|---|
| $d_1^L = \{e_1,e_2,e_5,e_{16}\}$ | $d_1^L = \{e_1,e_2,e_7,e_{16}\}$ |
| $d_2^L = \{e_1,e_3,e_4,e_8\}$ | $d_2^L = \{e_2,e_3,e_{14},e_{15}\}$ |
| $d_3^L = \{e_2,e_3,e_{13},e_{14}\}$ | $d_3^L = \{e_5,e_6,e_{10},e_{11}\}$ |
| $d_4^L = \{e_7,e_9,e_{11},e_{15}\}$ | $d_4^L = \{e_{12},e_{13},e_{14},e_{16}\}$ |
| $d_5^L = \{e_{17},e_{19},e_{21},e_{24}\}$ | $d_5^L = \{e_{17},e_{20},e_{21},e_{24}\}$ |
| $d_6^L = \{e_{18},e_{20},e_{28},e_{29}\}$ | $d_6^L = \{e_{18},e_{19},e_{27},e_{28}\}$ |
| $d_7^L = \{e_{23},e_{25},e_{26},e_{29}\}$ | $d_7^L = \{e_{21},e_{22},e_{23},e_{26}\}$ |
| $d_8^L = \{e_{26},e_{27},e_{28},e_{30}\}$ | $d_8^L = \{e_{23},e_{24},e_{25},e_{29}\}$ |

Кортеж весов ребер: $\xi_L(G_{15})$ =

= <7,5,11,12,12,10,10,16,10,10,7,4,6,13,7,10,10,12,10,16,7,4,6,7,13,5,10,7,11,12>.

Кортеж весов вершин : $\zeta_L(G_{15})$ = <23,41,58,21,33,21,62,21,58,21,33,21,23,41,21,62>.

Вектор весов ребер: $F(\xi_L(G_{15}))$ =

= (4,4,5,5,6,6,7,7,7,7,7,7,10,10,10,10,10,10,10,10,11,11,12,12,12,12,13,13,16,16).

Вектор весов вершин: $F(\zeta_L(G_{15}))$ = (21,21,21,21,21,21,23,23,33,33,41,41,58,58,62,62).

Кортеж весов ребер: $\xi_L(G_{16})$ =

= <12,7,12,7,10,10,16,13,4,7,7,6,13,5,10,11,6,7,7,13,5,10,7,11,12,12,10,10,16,4>.

Кортеж весов вершин : $\zeta_L(G_{16})$ = <38,61,21,21,33,21,23,62,33,21,23,38,61,21,21,62>.

Вектор весов ребер : $F(\xi_L(G_{16}))$ =

= (4,4,5,5,6,6,7,7,7,7,7,7,7,10,10,10,10,10,10,11,11,12,12,12,12,13,13,13,16,16).

Вектор весов вершин : $F(\zeta_L(G_{16})$ = (21,21,21,21,21,21,23,23,33,33,38,38,61,61,62,62),

Структуры реберного графа L(G) устанавливают различие для неизоморфных графов, даже в случае установления изоморфизма только суграфов подпространства разрезов или в случае установления изоморфизма только суграфов подпространства циклов графов.

## 6.7. Базисы подпространств и реберные графы

Если пространство $£_G$ разлагается в прямую сумму подпространств циклов и разрезов, то подпространства разрезов и циклов тесно связаны между собой и представляют жесткую систему. То есть, невозможно построить пространства суграфов £G используя подпространство разрезов одного графа и подпространство циклов другого графа. Если графы G и H изоморфны, то должны существовать изоморфные базисы и в подпространстве разрезов, и в подпространстве циклов этих графов, которые порождают элементы пространства суграфов $£_G$ и $£_H$. Таким образом, если найдены изоморфные базисы подпространства разрезов графов G и H, то для построения изоморфных пространств $£_G$ и $£_H$ должны быть найдены соответствующие изоморфные базисы подпространств циклов графов G и H.

Таким образом, для графов, которые не имеют пересекающихся элементов в подпространствах циклов и разрезов, можно установить следующее утверждение.

**Лемма 6.1.** Несепарабельные неориентированные графы G и H, имеющие одинаковое количество вершин и ребер и не имеющие пересекающихся элементов в подпространствах циклов C(G) и разрезов S(G), соответственно C(H) и S(H), изоморфны тогда и только тогда, когда изоморфны их подпространства циклов $C(G) \cong C(H)$ и разрезов $S(G) \cong S(H)$.

Доказательство.

Необходимость. Если изоморфны графы G и H, то существуют изоморфные базисы подпространств циклов C(G) и C(H), изоморфные базисы подпространств разрезов S(G) и S(H) графов, которые, в свою очередь, порождают изоморфные подпространства разрезов и циклов графов G и H.

Достаточность. Если изоморфны суграфы подпространств циклов C(G) и C(H), и изоморфны суграфы подпространств разрезов S(G) и S(H) графов G и H, то изоморфны суграфы пространств суграфов $£_G$ и $£_H$, так как каждый элемент пространств суграфов $£_G$ и $£_H$ может быть образован как результат кольцевого сложения одного элемента из подпространств циклов C(G) соответственно C(H) и другого элемента из подпространств разрезов S(G) соответственно S(H), согласно теореме 6.7. Другими словами, если элементы подпространств циклов $C(G) \cong C(H)$ и $S(G) \cong S(H)$, то изоморфна их кольцевая сумма $C(G) \oplus S(G) \equiv C(H) \oplus S(H)$. Следовательно, изоморфны пространства суграфов $£_G$ и $£_H$ графов G и H. Но тогда изоморфны ребра графов G и H, что является изоморфизмом графов G и H.

Лемма доказана.

## 6.7. Интегральный инвариант графа

Так как выделение множества изометрических циклов в реберном графе L(G) довольно громозко по объему, рассмотрим возможность сокращения информации.

Известно, что инвариант, основанный на векторе локальных степеней графа, не является полным инвариантом для определения изоморфизма графов [20].

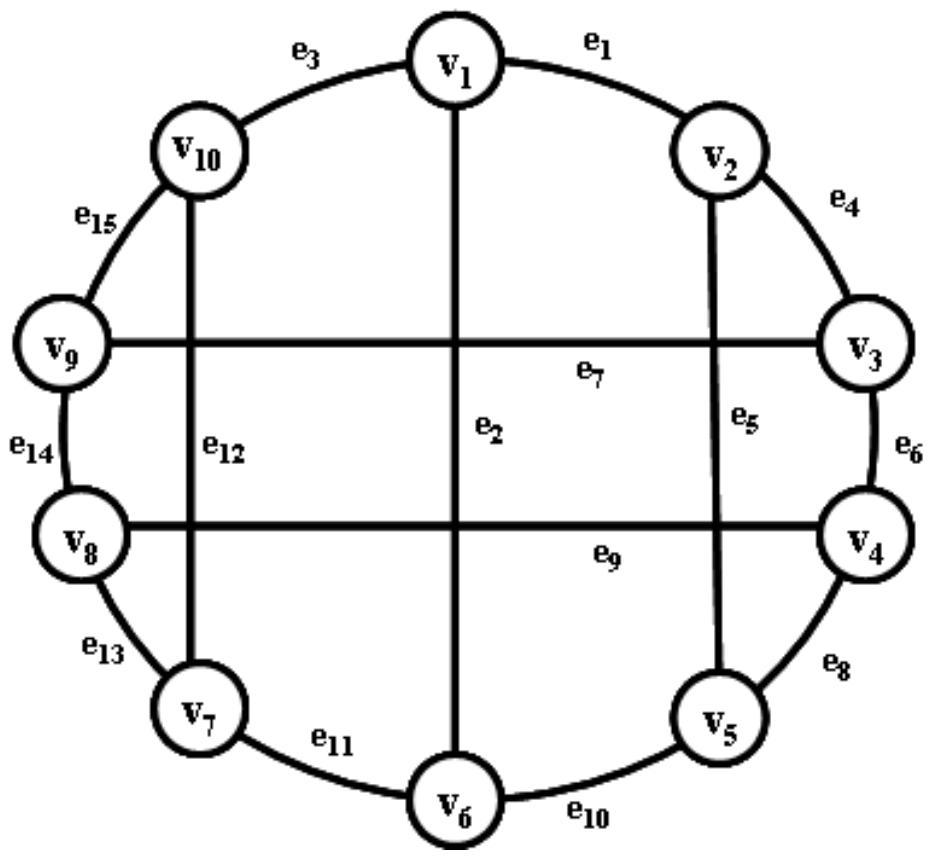

Рис. 6.10. Граф $G_{23}$.

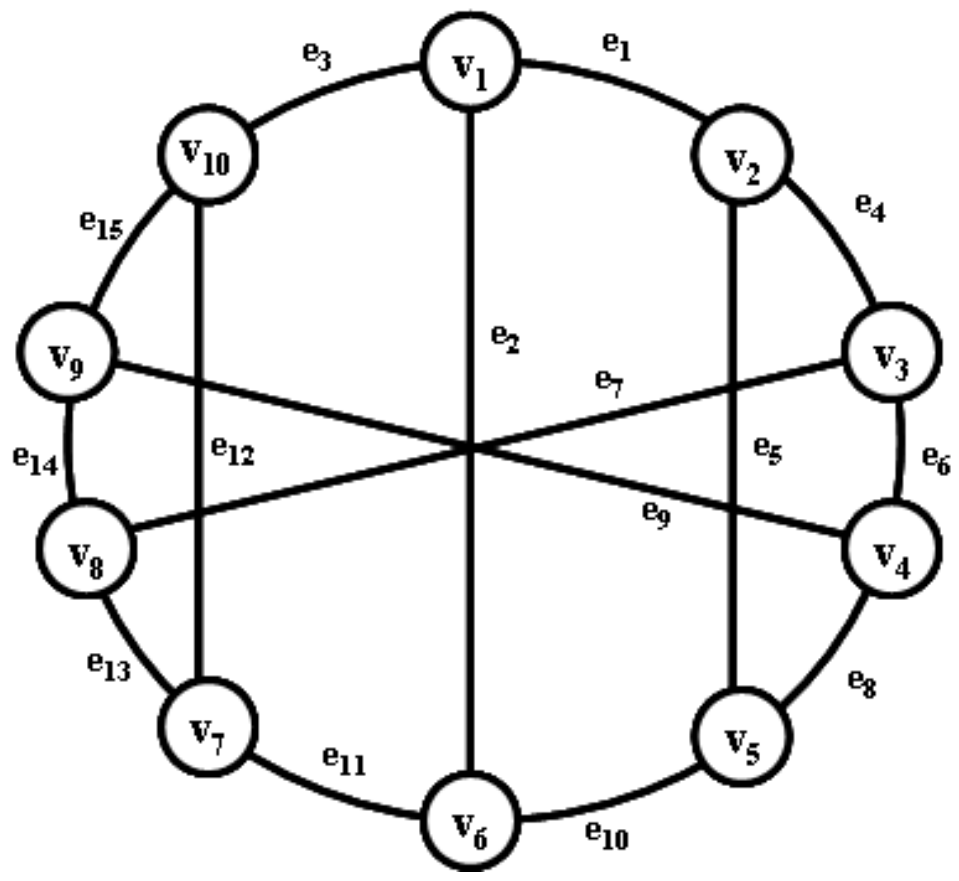

Рис. 6.11. Граф $G_{24}$.

Инвариант реберных разрезов графа, также не является полным инвариантом графа с точки зрения распознавания изоморфизма [20]. В качестве примера приведем пару неизоморфных графов имеющих одинаковое значение весов в векторных инвариантах спектра реберных разрезов (рис. 6.10 и 6.11).

Вектор весов ребер графа: $F_w(\xi(G_{23}))$ = (12,12,12,12,12,14,14,14,14,14,14,14,14,14,14);

Вектор весов вершин графа: $F_w(\zeta(G_{23}))$ = (40,40,40,40,40,40,40,40,40,40).

Вектор весов ребер графа: $F_w(\xi(G_{24}))$ = (12,12,12,12,12,14,14,14,14,14,14,14,14,14,14);

Вектор весов вершин графа: : $F_w(\zeta(G_{24}))$ = (40,40,40,40,40,40,40,40,40,40).

Однако, количество изометрических циклов в графах $G_{23}$ и $G_{24}$ различно. Соответственно, различны и спектры реберных циклов графов.

Количество изометрических циклов в графе $G_{23}$ = 7.

Изометрические циклы в графе $G_{23}$ в виде ребер:

$c_1 = \{e_1,e_2,e_5,e_{10}\}$;
$c_2 = \{e_1,e_3,e_4,e_7,e_{15}\}$;
$c_3 = \{e_2,e_3,e_{11},e_{13}\}$;
$c_4 = \{e_4,e_5,e_6,e_8\}$;
$c_5 = \{e_6,e_7,e_9,e_{14}\}$;
$c_6 = \{e_8,e_9,e_{10},e_{11},e_{12}\}$;
$c_7 = \{e_{12},e_{13},e_{14},e_{15}\}$.

Вектор весов ребер для спектра реберных циклов: $F_\tau(\xi(G_{23}))$ =

= (6,6,6,6,6,7,7,7,7,7,7,7,7,7,7);

Вектор весов вершин: $F_\tau(\zeta\,(G_{23}))$ = (20,20,20,20,20,20,20,20,20,20).

Количество изометрических циклов в графе $G_{24}$ = 15.

Изометрические циклы в графе $G_{24}$ в виде ребер:

$c_1 = \{e_1,e_2,e_5,e_{10}\}$;
$c_2 = \{e_1,e_3,e_4,e_7,e_{12},e_{13}\}$;
$c_3 = \{e_1,e_3,e_4,e_6,e_9,e_{15}\}$;
$c_4 = \{e_1,e_3,e_4,e_7,e_{14},e_{15}\}$;
$c_5 = \{e_2,e_3,e_{11},e_{13}\}$;
$c_6 = \{e_1,e_3,e_5,e_8,e_9,e_{15}\}$;
$c_7 = \{e_1,e_2,e_4,e_7,e_{11},e_{12}\}$;
$c_8 = \{e_4,e_5,e_6,e_8\}$;
$c_9 = \{e_6,e_7,e_9,e_{14}\}$;
$c_{10} = \{e_4,e_5,e_7,e_{10},e_{11},e_{12}\}$;
$c_{11} = \{e_8,e_9,e_{10},e_{11},e_{12},e_{14}\}$;
$c_{12} = \{e_2,e_3,e_8,e_9,e_{10},e_{15}\}$;
$c_{13} = \{e_8,e_9,e_{10},e_{11},e_{13},e_{15}\}$;
$c_{14} = \{e_6,e_7,e_8,e_{10},e_{11},e_{12}\}$;
$c_{15} = \{e_{12},e_{13},e_{14},e_{15}\}$.

Вектор весов ребер для спектра реберных циклов: $F_\tau(\xi\,(G_{24}))$ =

= (8,8,8,8,8,8,8,8,8,8,8,8,8,8,8);

Вектор весов вершин: $F_\tau(\zeta\,(G_{24}))$ = (24,24,24,24,24,24,24,24,24,24).

Векторный инвариант спектра реберных циклов также не может служить полным инвариантом, для распознавания изоморфизма графов. В качестве примера можно привести два неизоморфных графа $G_{15}$ и $G_{16}$, имеющих равное значение параметров в векторном инварианте реберных циклов (рис. 6.8 и рмс. 6.9)

Таким образом, можно сказать следующее:

Если значения весов элементов графа в векторном инварианте спектра реберных разрезов неизоморфных графов совпадают, то возможно не совпадают значения весов элементов в векторных инвариантах спектра реберных циклов графов. Если значения весов элементов в векторном инварианте спектра реберных циклов неизоморфных графов совпадают, то возможно не совпадают значения весов элементов в векторных инвариантох спектра реберных разрезов графов.

Следуя такой логике рассуждений, следует объединить инварианты реберных разрезов графа и инварианты реберных циклов графа, что косвенно сделано в цифровом векторном инварианте.

**Определение 6.4.** Будем называть *интегральным инвариантом графа* результат объединения векторных инвариантов спектра реберных разрезов графа и векторных инвариантов спектра реберных циклов графа.

Например, интегральный инвариант для графа $G_{16}$, будем иметь следующий вид:

$$F_w(\xi(G_{16})) \& F_w(\zeta(G_{16})) \& F_\tau(\xi(G_{16})) \& F_\tau(\zeta(G_{16})) =$$

$$= (2\times 4,4\times 12,2\times 13,4\times 14,4\times 15,1\times 16,2\times 17,2\times 18,2\times 19,4\times 20,2\times 21,1\times 28)\ \&$$

$$\&\ (2\times 37,4\times 38,2\times 44,2\times 55,2\times 72,2\times 80,2\times 108)\ \&.$$

$$\&\ (8\times 4,10\times 5,2\times 11,2\times 12,4\times 13,2\times 15,2\times 17)\ \&$$

$$\&\ (4\times 14,2\times 28,2\times 30,2\times 34,2\times 36,4\times 44).$$

Постараемся ответить на следующий вопрос: существует ли одинаковые значения весов элементов в интегральном инварианте для двух неизоморфных неориентированных несепарабельных графов?

Для графов, не имеющих пересекающихся элементов в подпространствах разрезов и циклов ответ очевиден. Рассмотрим более подробно структуру цифрового инварианта реберного графа L(G). Цифровой инвариант реберного графа L(G) является полным инвариантом графа G, так как построен на основании теоремы Уитни о изоморфизме графов. Определен состав векторного цифрового инварианта, состоящий из множества изометрических циклов графа L(G), образы которого соответствуют прообразам центральных разрезов графа G и прообразам изометрических циклов графа G и дубль-циклам графа G, являющихся кольцевым сложение .изометрических циклов.

Поэтому можно предположить, что в общем случае интегральный инвариант графа есть другая форма записи цифрового инварианта реберного графа без учета дубль-циклов.

Действительно рассматривая векторные цифровые инварианты графов Петерсена и авудольного графа (см. главу 5), можно заметить отсутствик в структурах дубль-циклов графа.

В качестве проверки предположения, что можно обойтись без учета дубль циклов графа, рассмотрим графы $G_{15}$ и $G_{16}$. Различие весов ребер в спектре реберных разрезов чувствуется уже на этапе формирования базовых реберных разрезов при одинаковом значении весов ребер в спектре реберных циклов:

Вектор весов ребер уровня 1 графа $G_{15}$: $F_\tau(\xi(G_{15}) = (4\times 4,8\times 5,8\times 6,8\times 7,2\times 8)$;

Вектор весов ребер уровня 1 графа $G_{16}$: $F_\tau(\xi(G_{16}) = (4\times 4,8\times 5,9\times 6,6\times 7,3\times 8)$.

Это различие и определяет разницу весов ребер и вершин в цифровых инвариантах реберных графов.

Вектор весов ребер в векторном цифровом инварианте: F( $\xi_L$ (G15)) =

= $(2\times4,2\times5,2\times6,6\times7,8\times10,2\times11,4\times12,2\times13,2\times16)$;

Вектор весов вершин в векторном цифровом инварианте: F( $\zeta_L$ (G15)) =

= $(6\times21,2\times23,2\times33,2\times41,2\times58.2\times62)$.

Вектор весов ребер в цифровом инварианте: F( $\xi_L$( G16)) =:

= $(2\times4,2\times5,2\times6,7\times7,6\times10,2\times11,4\times12,3\times13,2\times16)$;

Вектор весов вершин в цифровом инварианте: F( $\zeta_L$ (G16)) =

= $(6\times21,2\times23,2\times33,2\times38,2\times61.2\times62)$.

Таким образом, существует некая структурная аналогия между цифровым инвариантов реберного графа IL(G) и интегральным инвариантом графа IS(G) & IC(G). На основании вышесказанного интегральный инвариант графа является полным инвариантом только для несепарабельных графов, не имеющих пересекающихся элементов в подпространствах разрезов и циклов. В других случаях интегральный инвариант применяется по аналогии с цифровым векторным инвариантом, исключая дубль-циклы из состава структур. В свою очередь, дубль-циклы определяются кольцевым суммированием изометрических циклов графа.

Определим состав интегрального инварианта графа G:

- матрица инциденций графа B(G);
- $n$ – количество вершин графа G;
- $m$ – количество ребер графа G;
- множество изометрических циклов $C_\tau$ графа G;
- $\xi$ (w(G)) – кортеж весов ребер для спектра реберных разрезов графа;
- $\zeta$ ( w (G)) – кортеж весов вершин для спектра реберных разрезов графа.
- $\xi$ ( $\tau$ (G)) – кортеж весов ребер для спектра реберных циклов графа;
- $\zeta$ ( $\tau$ (G)) – кортеж весов вершин для спектра реберных циклов графа.
- IS(G) = $F_w(\xi(G))\,\&\,F_w(\zeta(G))$ – векторный инвариант спектра реберных разрезов графа G;
- IC(G) = $F_\tau(\xi(G))\,\&\,F_\tau(\zeta(G))$ – векторный инвариант спектра реберных циклов графа G;
- $F_w(\xi(G))\,\&\,F_w(\zeta(G))\,\&\,F_\tau(\xi(G))\,\&\,F_\tau(\zeta(G))$ – векторный интегральный инвариант графа G.

**Определение 6.5.** Инварианты графа, построенные на основе понятия веса элемента, будем называть *векторными инвариантами графа*.

К векторным инвариантам относится цифровой инвариант реберного графа IL(G), инварианты уровней спектра реберных разрезов графа IS(G), инварианты уровней спектра реберных циклов графа IC(G), интегральный инвариант графа. Векторный инвариант спектра будем считать суммарным инвариантом уровней спектра графа.

## Комментарии

Так как существуют одинаковые инварианты реберных разрезов для неизоморфных графов, и существуют одинаковые инварианты реберных циклов для неизоморфных графов, то для описания структуры графов представлен интегральный инвариант графа. Интегральный инвариант одновременно отражает свойства элементов и подпространства разрезов и подпространства циклов графа и является полным инвариантом. Интегральный инвариант характеризуется меньшим объемом информации, по сравнению с цифровым инвариантом графа.

## Глава 7. Симметрическая устойчивость весов вершин графа

### 7.1. Векторные инварианты вершин графа

Основная цель задачи распознавания графов - это не только определение соответствия между сравниваемыми элементами двух графов, но и определение соответствия элементов в самом графе [1,2].

Как правило, до вычисления спектра реберных разрезов или реберных циклов, количество уровней неизвестно. Естественно, что вычисление элементов спектра реберных разрезов требует определенных затрат машинного времени и памяти. Поэтому при создании программного обеспечения, требуется определить оптимальную величину количества уровней в спектре реберных разрезов графа.

Наличие оператора $W_\lambda$ в пространстве суграфов £(G) графа G позволяет осуществлять преобразование векторного пространства суграфов само в себя. То есть имеет смысл рассматривать операторы $W_\lambda^2 = W_\lambda \times W_\lambda, W_\lambda^3 = W_\lambda \times W_\lambda \times W_\lambda, \ldots, W_\lambda^q = W_\lambda \times W_\lambda \times \ldots \times W_\lambda$.

Такое преобразование представляет собой циклический процесс, так как, начиная с какого-то момента $W^q = W$, содержимое матриц начинает повторяться. Другими словами, такой процесс характеризуется ограниченным количеством уровней в процессе порождения реберных разрезов и влияет на вычислительную сложность алгоритма вычисления.

Нас интересует вопрос: что будет, если значение парамера $q$ будет превосходить количество ребер при ограничении $W^q = W_0 \to W^q = \varnothing$?

Очевидно, что поиск решения задачи распознования стркуктуры графа – это установление различения весов вершин в графе. И здесь следует установить минимальное количество уровней для уменьшения вычислительной сложности алгоритма распознавания. Рассмотрим процесс образования уровней с точки зрения различения весов вершин для каждого уровня.

В качестве примера рассмотрим следующий граф $G_{25}$ (рис. 7.1).

Центральные разрезы графа:

$s_1 = \{e_1,e_2,e_3\}$;
$s_2 = \{e_1,e_4,e_5\}$;
$s_3 = \{e_2,e_6,e_7\}$;
$s_4 = \{e_6,e_8,e_9\}$;
$s_5 = \{e_8,e_{10},e_{11}\}$;
$s_6 = \{e_4,e_{10},e_{12}\}$;
$s_7 = \{e_3,e_9,e_{13}\}$;
$s_8 = \{e_7,e_{11},e_{14}\}$;
$s_9 = \{e_5,e_{13},e_{15}\}$;

$s_{10} = \{e_{12},e_{14},e_{15}\}$.

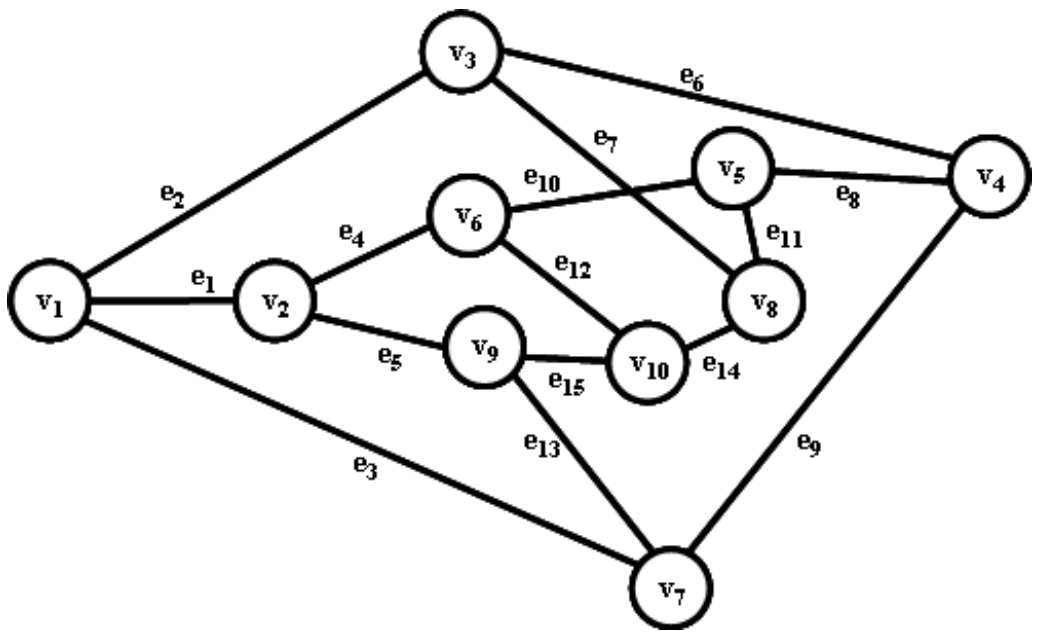

Рис. 7.1. Граф $G_{25}$.

Базовые реберные разрезы графа $G_{25}$:

$w_0(e_1) = \{e_2,e_3,e_4,e_5\}$;
$w_0(e_2) = \{e_1,e_3,e_6,e_7\}$;
$w_0(e_3) = \{e_1,e_2,e_9,e_{13}\}$;
$w_0(e_4) = \{e_1,e_5,e_{10},e_{12}\}$;
$w_0(e_5) = \{e_1,e_4,e_{13},e_{15}\}$;
$w_0(e_6) = \{e_2,e_7,e_8,e_9\}$;
$w_0(e_7) = \{e_2,e_6,e_{11},e_{14}\}$;
$w_0(e_8) = \{e_6,e_9,e_{10},e_{11}\}$;
$w_0(e_9) = \{e_3,e_6,e_8,e_{13}\}$;
$w_0(e_{10}) = \{e_4,e_8,e_{11},e_{12}\}$;
$w_0(e_{11}) = \{e_7,e_8,e_{10},e_{14}\}$;
$w_0(e_{12}) = \{e_4,e_{10},e_{14},e_{15}\}$;
$w_0(e_{13}) = \{e_3,e_5,e_9,e_{15}\}$;
$w_0(e_{14}) = \{e_7,e_{11},e_{12},e_{15}\}$;
$w_0(e_{15}) = \{e_5,e_{12},e_{13},e_{14}\}$.

Будем рассматривать только первые три уровня в спектре реберных разрезов.

Кортежи весов ребер по уровням:

$\xi(w(l_0))=$ <4,4,4,4,4,4,4,4,4,4,4,4,4,4,4>;
$\xi(w(l_1))=$ <10,10,8,10,8,8,10,10,10,10,8,8,10,10,10>;
$\xi(w(l_2))=$ <6,6,4,6,4,4,6,6,6,6,4,4,6,6,6>.

Кортежи весов вершин по уровням:

$\zeta(w(l_0))=$ <12,12,12,12,12,12,12,12,12,12>;
$\zeta(w(l_1))=$ <28,28,28,28,28,28,28,28,28,28>;
$\zeta(w(l_2))=$ <16,16,16,16,16,16,16,16,16,16>.

Значения весов вершин в кортежах могут меняться от уровня к уровню, но значения весов в кортеже не меняетсясостав блока не меняется.

Рассмотрим граф $G_2$, представленный на рис. 7.2.

Кортеж весов вершин для графа $G_8$ по уровням:

$\varsigma(w(l_0))=$ <15,23,23,22,15,22>;
$\varsigma(w(l_1))=$ <16,28,28,24,16,24>;
$\varsigma(w(l_2))=$ <20,20,20,24,20,24>.

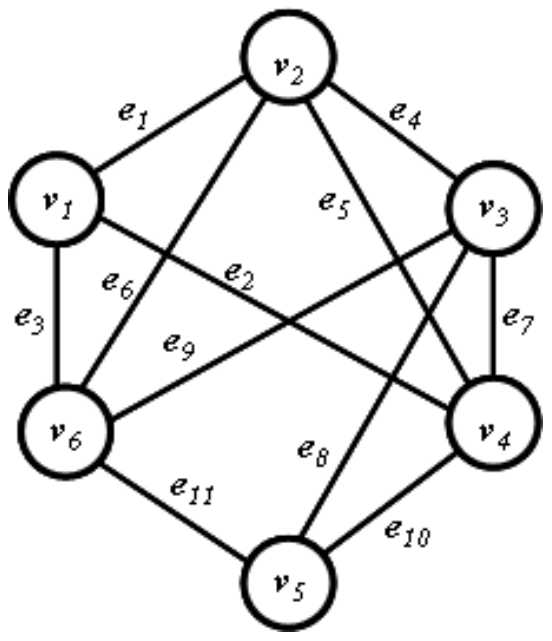

Рис. 7.2. Граф $G_2$.

Как мы видим, на уровне 2 происходит расхождение значений весов. На первых лвух уровнях пара вершин ($v_1$,$v_5$) имеют равный вес, пара вершин ($v_2$,$v_3$) также имеют равный вес, пара ($v_4$,$v_6$) имеет равный вес. Однако, на уровне 2 пары ($v_1$,$v_5$) и ($v_1$,$v_5$) имеют равные веса. .

Рассмотрим граф $G_{26}$ (рис. 7.3). Данный граф образован из полного графа $K_6$ путем удаления шести ребер.

Центральные разрезы графа $G_9$:

$s_1 = \{e_1,e_2,e_3\}$;
$s_2 = \{e_1,e_4,e_5\}$;
$s_3 = \{e_4,e_6,e_7\}$;
$s_4 = \{e_2,e_6,e_8\}$;
$s_5 = \{e_7,e_8,e_9\}$;
$s_6 = \{e_3,e_5,e_9\}$.

Множество базовых реберных разрезов графа $G_9$:

$w_0(e_1) = s_1 \oplus s_2 = \{e_2,e_3,e_4,e_5\}$;
$w_0(e_2) = s_1 \oplus s_4 = \{e_1,e_3,e_6,e_8\}$;
$w_0(e_3) = s_1 \oplus s_6 = \{e_1,e_2,e_5,e_9\}$;
$w_0(e_4) = s_2 \oplus s_3 = \{e_1,e_5,e_6,e_7\}$;
$w_0(e_5) = s_2 \oplus s_6 = \{e_1,e_3,e_4,e_9\}$;
$w_0(e_6) = s_3 \oplus s_4 = \{e_2,e_4,e_7,e_8\}$;
$w_0(e_7) = s_3 \oplus s_5 = \{e_4,e_6,e_8,e_9\}$;
$w_0(e_8) = s_4 \oplus s_5 = \{e_2,e_6,e_7,e_9\}$;
$w_0 (e_9) = s_5 \oplus s_6 = \{e_3,e_5,e_7,e_8\}$.

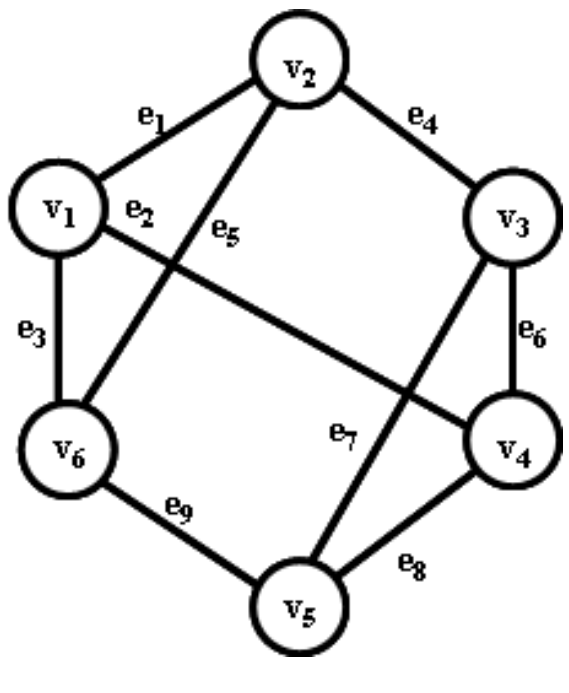

Рис. 7.3. Граф $G_{26}$.

Множество реберных разрезов 1-го уровня графа $G_{26}$:

$w_1(e_1) = s_3 \oplus s_4 = \{e_2,e_4,e_7,e_8\}$;
$w_1(e_2) = \varnothing$;
$w_1(e_3) = s_4 \oplus s_5 = \{e_2,e_6,e_7,e_9\}$;
$w_1(e_4) = \varnothing$;
$w_1(e_5) = s_3 \oplus s_5 = \{e_4,e_6,e_8,e_9\}$;
$w_1(e_6) = s_1 \oplus s_2 = \{e_2,e_3,e_4,e_5\}$;
$w_1(e_7) = s_2 \oplus s_6 = \{e_1,e_3,e_4,e_9\}$;
$w_1(e_8) = s_1 \oplus s_6 = \{e_1,e_2,e_5,e_9\}$;
$w_1(e_9) = \varnothing$.

Количество уровней в данном графе равно двум.

Составим кортежи весов вершин для уровней:

$\zeta(w(l_0)) =$ <12,12,12,12,12,12>;

$\zeta(w(l_1)) =$ <20,20,20,20,20,20>.

В этом случае, проявляется свойство подмножества вершин сохранять пропорциональное значение весов для любого уровня спектра реберных разрезов. В данном случае в каждом уровне множества веса всех вершин равны.

Рассмотрим следующий граф $G_{27}$ (рис. 7.4).

Центральные разрезы графа:

$s_1 = \{e_1,e_2,e_3,e_4,e_5\}$;
$s_2 = \{e_6,e_7,e_8,e_9\}$;
$s_3 = \{e_6,e_{10},e_{11},e_{12},e_{13}\}$;
$s_4 = \{e_1,e_{14},e_{15},e_{16}\}$;
$s_5 = \{e_{10},e_{17},e_{18},e_{19},e_{20}\}$;
$s_6 = \{e_2,e_{11},e_{17},e_{21},e_{22}\}$;
$s_7 = \{e_3,e_7,e_{12},e_{14},e_{18},e_{21}\}$;
$s_8 = \{e_4,e_8,e_{15},e_{19},e_{23}\}$;
$s_9 = \{e_5,e_9,e_{13},e_{16},e_{20},e_{22},e_{23}\}$.

Количество ребер в графе = 23;

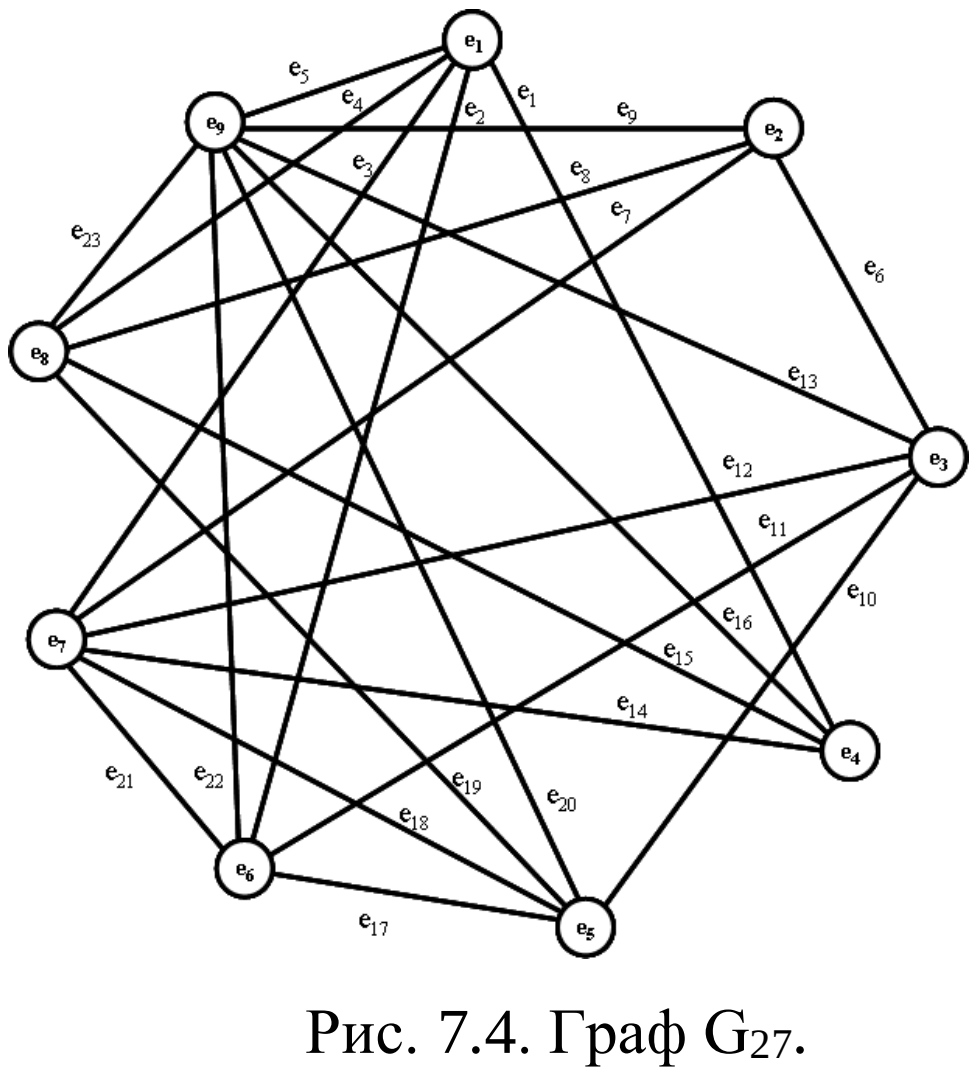

Рис. 7.4. Граф $G_{27}$.

Базовые реберные разрезы графа:

$w_0(e_1) = \{e_2,e_3,e_4,e_5,e_{14},e_{15},e_{16}\}$;
$w_0(e_2) = \{e_1,e_3,e_4,e_5,e_{11},e_{17},e_{21},e_{22}\}$;
$w_0(e_3) = \{e_1,e_2,e_4,e_5,e_7,e_{12},e_{14},e_{18},e_{21}\}$;
$w_0(e_4) = \{e_1,e_2,e_3,e_5,e_8,e_{15},e_{19},e_{23}\}$;
$w_0(e_5) = \{e_1,e_2,e_3,e_4,e_9,e_{13},e_{16},e_{20},e_{22},e_{23}\}$;
$w_0(e_6) = \{e_7,e_8,e_9,e_{10},e_{11},e_{12},e_{13}\}$;
$w_0(e_7) = \{e_3,e_6,e_8,e_9,e_{12},e_{14},e_{18},e_{21}\}$;
$w_0(e_8) = \{e_4,e_6,e_7,e_9,e_{15},e_{19},e_{23}\}$;
$w_0(e_9) = \{e_5,e_6,e_7,e_8,e_{13},e_{16},e_{20},e_{22},e_{23}\}$;
$w_0(e_{10}) = \{e_6,e_{11},e_{12},e_{13},e_{17},e_{18},e_{19},e_{20}\}$;
$w_0(e_{11}) = \{e_2,e_6,e_{10},e_{12},e_{13},e_{17},e_{21},e_{22}\}$;
$w_0(e_{12}) = \{e_3,e_6,e_7,e_{10},e_{11},e_{13},e_{14},e_{18},e_{21}\}$;
$w_0(e_{13}) = \{e_5,e_6,e_9,e_{10},e_{11},e_{12},e_{16},e_{20},e_{22},e_{23}\}$;
$w_0(e_{14}) = \{e_1,e_3,e_7,e_{12},e_{15},e_{16},e_{18},e_{21}\}$;
$w_0(e_{15}) = \{e_1,e_4,e_8,e_{14},e_{16},e_{19},e_{23}\}$;
$w_0(e_{16}) = \{e_1,e_5,e_9,e_{13},e_{14},e_{15},e_{20},e_{22},e_{23}\}$;
$w_0(e_{17}) = \{e_2,e_{10},e_{11},e_{18},e_{19},e_{20},e_{21},e_{22}\}$;
$w_0(e_{18}) = \{e_3,e_7,e_{10},e_{12},e_{14},e_{17},e_{19},e_{20},e_{21}\}$;
$w_0(e_{19}) = \{e_4,e_8,e_{10},e_{15},e_{17},e_{18},e_{20},e_{23}\}$;
$w_0(e_{20}) = \{e_5,e_9,e_{10},e_{13},e_{16},e_{17},e_{18},e_{19},e_{22},e_{23}\}$;
$w_0(e_{21}) = \{e_2,e_3,e_7,e_{11},e_{12},e_{14},e_{17},e_{18},e_{22}\}$;
$w_0(e_{22}) = \{e_2,e_5,e_9,e_{11},e_{13},e_{16},e_{17},e_{20},e_{21},e_{23}\}$;
$w_0(e_{23}) = \{e_4,e_5,e_8,e_9,e_{13},e_{15},e_{16},e_{19},e_{20},e_{22}\}$.

Кортежи весов спектра реберных разрезов графа $G_{27}$ по уровням:

$\xi(w(l_0))$ = <7,8,9,8,10,7,8,7,9,8,8,9,10,8,7,9,8,9,8,10,9,10,10>;
$\xi(w(l_1))$ = <9,13,10,12,10,13,5,13,13,7,9,12,10,5,13,11,8,9,11,11,13,11,8>;
$\xi(w(l_2))$ = <12,12,12,13,10,11,8,15,14,13,13,13,13,10,11,12,12,4,15,12,14,12,9>;
$\xi(w(l_3))$ = <12,12,14,10,11,13,11,11,12,15,14,8,13,12,14,15,13,15,5,14,14,13,13>;
$\xi(w(l_4))$ = <4,13,11,10,14,12,7,12,14,11,11,11,8,13,10,14,10,12,11,5,10,11,14>’;
$\xi(w(l_5))$ = <17,11,12,13,14,10,13,14,11,12,12,15,5,13,8,13,10,11,14,9,11,11,11>;
$\xi(w(l_6))$ = <10,12,12,13,13,8,11,10,12,14,13,5,8,10,5,9,13,13,10,10,10,9,12>;
$\xi(w(l_7))$ = <10,11,13,8,11,13,14,11,12,10,11,11,11,11,10,13,13,15,12,13,10,8,13>;
$\xi(w(l_8))$ = <13,10,13,5,12,12,10,12,9,11,7,14,7,12,12,15,8,11,13,12,13,8,11>;
$\xi(w(l_9))$ = <9,13,12,8,5,12,8,14,9,12,9,10,13,13,13,14,5,8,10,9,11,12,7>;
$\xi(w(l_{10}))$ = <12,13,11,11,10,13,5,5,4,14,12,12,11,11,13,12,8,10,8,9,12,13,9>;
$\xi(w(l_{11}))$ = <13,13,8,7,13,11,10,9,15,5,12,13,14,11,10,12,11,14,9,13,13,10,12>;
$\xi(w(l_{12}))$ = <8,7,7,13,11,12,12,10,12,9,4,10,10,9,11,13,7,11,13,9,10,10,12>;
$\xi(w(l_{13}))$ = <11,6,13,11,15,11,13,13,11,10,17,8,8,8,12,12,13,12,14,10,17,13,4>;
$\xi(w(l_{14}))$ = <15,16,15,14,15,10,13,12,13,13,10,9,5,8,13,10,11,10,11,12,11,13,17>;
$\xi(w(l_{15}))$ = <5,16,11,10,12,13,10,9,15,12,10,13,10,12,11,9,14,13,12,14,11,8,10>;
$\xi(w(l_{16}))$ = <11,9,12,8,12,13,11,9,11,9,13,14,12,13,11,15,10,13,13,11,11,5,10>;
$\xi(w(l_{17}))$ = <14,12,14,5,14,13,12,13,8,9,9,11,13,10,9,14,8,14,11,12,12,8,13>;
$\xi(w(l_{18}))$ = <10,14,10,10,11,7,13,9,10,13,12,12,13,10,8,9,5,13,7,12,12,11,9>;
$\xi(w(l_{19}))$ = <12,13,11,12,14,6,11,10,14,9,13,13,10,13,8,14,10,12,9,13,14,7,12>;
$\xi(w(l_{20}))$ = <13,11,12,13,12,16,11,12,11,10,8,11,11,13,12,11,12,13,12,13,15,13,13>;
$\xi(w(l_{21}))$ = <10,14,14,13,13,16,9,14,12,12,11,7,12,8,13,13,13,13,12,12,10,11,8>;
$\xi(w(l_{22}))$ = <8,13,5,10,11,9,8,11,10,14,15,9,13,5,10,11,13,9,4,11,12,14,11>;

$\xi$ (w($l_{23}$)) = <9,14,9,11,14,12,8,12,13,11,5,12,11,8,10,13,10,13,17,8,13,10,15>;
$\xi$ (w($l_{24}$)) = <13,7,10,12,11,14,12,12,13,12,11,12,11,11,13,12,11,14,10,7,11,8,5>;
$\xi$ (w($l_{25}$)) = <14,13,13,13,14,13,13,13,14,12,14,4,9,7,13,10,12,12,10,13,10,5,11>;
$\xi$ (w($l_{26}$)) = <11,11,12,11,11,11,10,13,13,13,10,17,8,13,8,12,13,14,13,15,15,10,14>;
$\xi$ (w($l_{27}$)) = <12,13,9,11,11,14,10,12,12,13,12,10,8,11,5,9,11,11,9,11,14,12,10>;
$\xi$ (w($l_{28}$)) = <13,12,9,9,13,13,13,11,13,12,13,10,12,14,8,12,11,12,12,12,11,13,12>;
$\xi$ (w($l_{29}$)) = <11,10,13,8,13,14,13,8,13,11,10,13,13,10,11,12,9,12,13,14,11,13,13>.

Спектр реберных разрезов графа $G_{27}$ имеет 30 уровней. Для определения различения мы рассмотрим только три из них.

Кортежи весов вершин для первых трех уровней:

$\zeta(w(l_0))$ = <42,31,42,34,43,44,52,40,69>;
$\zeta(w(l_1))$ = <52,44,51,38,46,53,54,55,74>;
$\zeta(w(l_2))$ = <59,63,48,45,56,63,64,63,82>.

Как видим с увеличение номера яруса различение весов вершин увеличивается.

Рассмотрим граф $G_{11}$ (рис. 7.5).

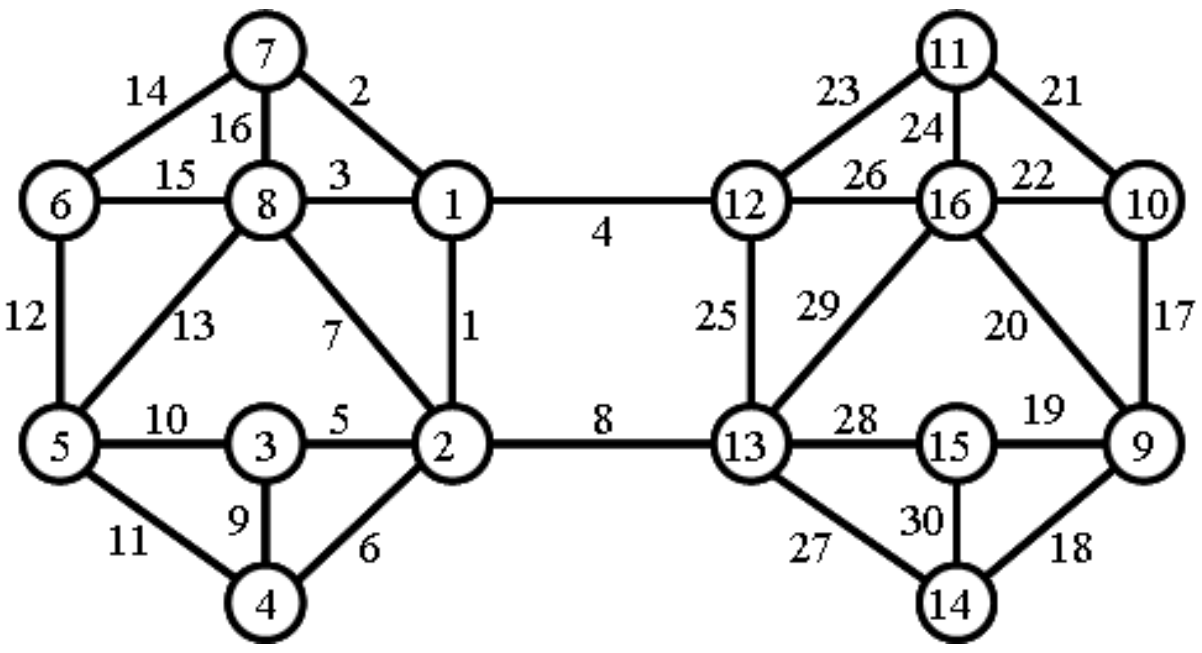

Рис. 7.5. Граф $G_{16}$.

Центральные разрезы графа $G_{16}$:

$s_1 = \{e_1,e_2,e_3,e_4\}$; $s_2 = \{e_1,e_5,e_6,e_7,e_8\}$; $s_3 = \{e_5,e_9,e_{10}\}$; $s_4 = \{e_6,e_9,e_{11}\}$; $s_5 = \{e_{10},e_{11},e_{12},e_{13}\}$;
$s_6 = \{e_{12},e_{14},e_{15}\}$; $s_7 = \{e_2,e_{14},e_{16}\}$; $s_8 = \{e_3,e_7,e_{13},e_{15},e_{16}\}$; $s_9 = \{e_{17},e_{18},e_{19},e_{20}\}$;
$s_{10} = \{e_{17},e_{21},e_{22}\}$; $s_{11} = \{e_{21},e_{23},e_{24}\}$; $s_{12} = \{e_4,e_{23},e_{25},e_{26}\}$; $s_{13} = \{e_8,e_{25},e_{27},e_{28},e_{29}\}$;
$s_{14} = \{e_{18},e_{27},e_{30}\}$; $s_{15} = \{e_{19},e_{28},e_{30}\}$; $s_{16} = \{e_{20},e_{22},e_{24},e_{26},e_{29}\}$.

Кортежи весов ребер по уровням:

$\xi$(w($l_0$)) = <7,5,7,6,6,6,8,8,4,5,5,5,7,4,6,6,5,5,5,7,4,6,5,6,7,7,6,6,8,4>;
$\xi$(w($l_1$)) = <14,12,11,10,14,14,11,20,0,9,9,7,8,8,7,9,7,9,9,8,8,7,12,9,14,11,14,14,11,0>;
$\xi$(w($l_2$)) = <10,12,9,14,12,12,9,10,0,8,8,9,7,13,14,13,9,8,8,7,13,14,12,13,10,9,12,12,9,0>;
$\xi$(w($l_3$)) = <10,17,15,16,11,11,11,14,0,14,14,12,14,12,10,14,12,14,14,14,12,10,17,14,10,15, 11,11,11,0>;
$\xi$(w($l_4$)) = <11,15,16,14,14,14,15,16,0,15,15,12,10,17,10,13,12,15,15,10,17,10,15,13,11,16, 14,14,15,0>;
$\xi$(w($l_5$)) = <16,11,10,14,17,17,14,14,0,19,19,17,10,16,11,17,17,19,19,10,16,11,11,17,16,10, 17,17,14,0>;
$\xi$(w($l_6$)) = <7,9,19,12,14,14,16,14,0,13,13,15,11,14,16,14,15,13,13,11,14,16,9,14,7,19, 14,14,16,0>;

$\xi(\mathrm{w}(l_7))$ = <14,12,15,10,16,16,11,12,0,15,15,11,16,12,7,9,11,15,15,16,12,7,12,9,14,15,16,16,11,0>;

$\xi(\mathrm{w}(l_8))$ = <10,12,9,14,12,12,9,10,0,8,8,9,7,13,14,13,9,8,8,7,13,14,12,13,10,9,12,12,9,0>;

$\xi(\mathrm{w}(l_9))$ = <10,17,15,16,11,11,11,14,0,14,14,12,14,12,10,14,12,14,14,14,12,10,17,14,10,15,11,11,11,0>;

$\xi(\mathrm{w}(l_{10}))$ = <11,15,16,14,14,14,15,16,0,15,15,12,10,17,10,13,12,15,15,10,17,10,15,13,11,16,14,14,15,0>;

$\xi(\mathrm{w}(l_{11}))$ = <16,11,10,14,17,17,14,14,0,19,19,17,10,16,11,17,17,19,19,10,16,11,11,17,16,10,17,17,14,0>;

$\xi(\mathrm{w}(l_{13}))$ = <7,8,18,10,13,13,15,12,0,13,13,15,11,14,16,14,15,13,13,11,14,16,8,14,7,18,13,13,15,0>;

$\xi(\mathrm{w}(l_{14}))$ = <13,11,13,10,14,14,10,10,0,14,14,11,16,11,7,8,11,14,14,16,11,7,11,8,13,13,14,14,10,0>.

Кортежи весов вершин для первых трех уровней:

$\zeta(\mathrm{w}(l_0))$ = <25,35,15,15,22,15,15,34,22,15,15,25,35,15,15>;

$\zeta(\mathrm{w}(l_1))$ = <47,73,23,23,33,22,29,46,33,22,29,47,73,23,23>;

$\zeta(\mathrm{w}(l_2))$ = <45,53,20,20, 32,36,38,52,32,36,38,45,53,20,20>.

На уровне 1 произошло разделение на блоки с разными весами, что увеличивает степень различения верщин графа. Уровень 2 остается без изменений.

В общем можно отметить, что разделение блоков вершин с равными весами с увеличением номеров уровней увеличивает степень различения, а объединение уменьшает степень различения. Таким образом, различение весов вершин графа влияет на количество уровней в спектре реберных разрезов и определяет его структуру.

Из рассмотренного материала следует, что можно ограничиться первыми уровнями в спектре реберных разрезов. Поэтому корректно ввести второе ограничение на процесс порождения реберных разрезов: количество уровней в спектре реберных разрезов не должно превышать двух.

Таким образом, процесс порождения реберных разрезов можно ограничить двумя факторами. Первое ограничение – оператор $\mathrm{W}_\lambda^q$ состоит из пустых строк или $\mathrm{W}_\lambda^q = \varnothing$ если $\mathrm{W}_\lambda^q = \mathrm{W}_\lambda$, а второе ограничение – величина уровней $l$ в спектре реберных разрезов не должна превышать двойки.

С учетом последнего, вычислительную сложность алгоритма построения реберных разрезов можно определить следующим образом: количество элементарных операций для построения суграфа как элемента спектра реберных разрезов определяется как произведение количества ребер суграфа предыдущего уровня на количество характеристических векторов размера $m$. То есть, вычислительную сложность определения суграфа как элемента спектра реберных разрезов можно определить как $m^2$. Тогда количество суграфов в столбце определяется как $m \times m^2 = m^3$. В свою очередь, размер матрицы $W_S$ можно определить в предположении, что количество уровней не превосходит двойки. Итого, для построения матрицы $W_S(G)$ нужно применить $4 \times m^3 = 4m^3$ элементарных операций с одновременным вычислением веса ребер. Вычислительную сложность алгоритма построения инварианта реберных разрезов графа, можно определить относительно ребер графа как $O(m^3)$.

Задача определения инварианта реберных разрезов графа относится к классу P – полиномиальных алгоритмов.

### 7.2. Симметрическая устойчивость векторного инварианта

В качестве еще одного примера рассмотрим граф $G_{28}$ представленный на рис. 7.6,a. Спектр реберных разрезов графа:

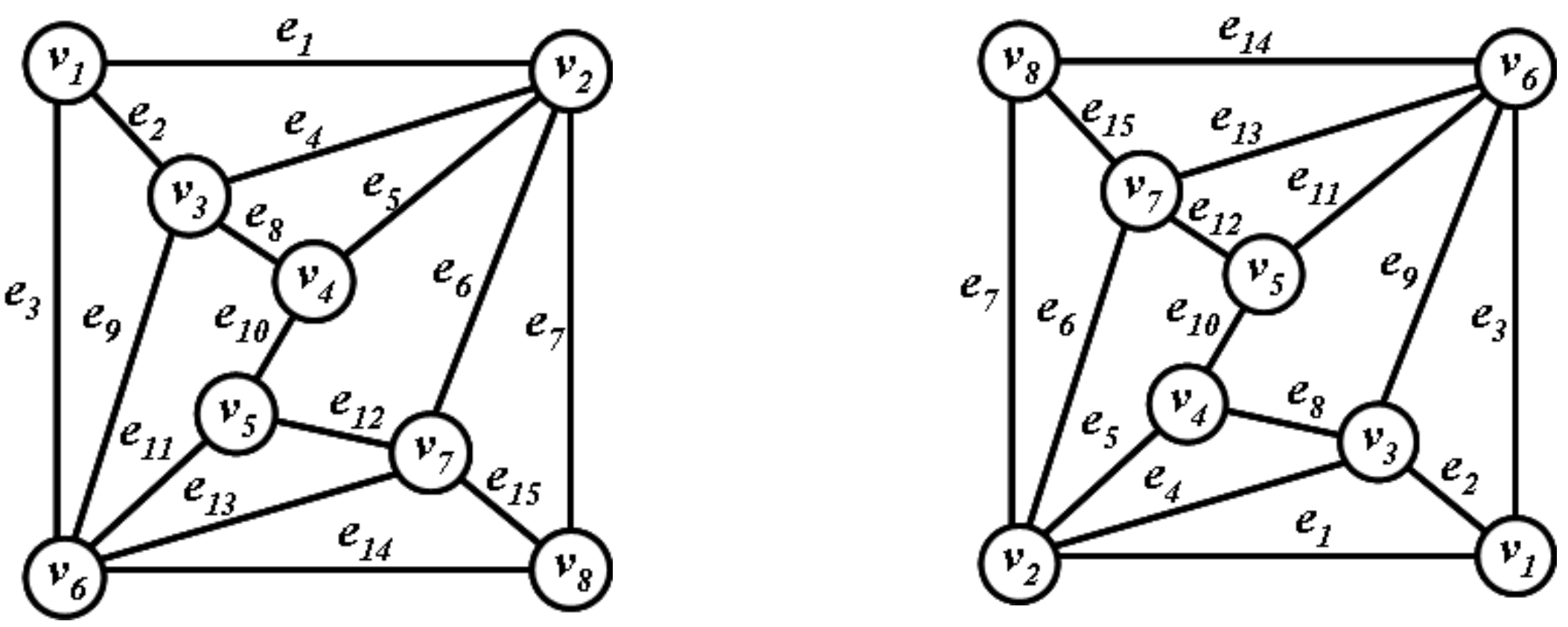

a) b)

Рис. 7.6. Автоморфизм графа $G_{28}$.

Смежность графа $G_{28}$:

вершина $v_1$: $v_2$ $v_3$ $v_6$
вершина $v_2$: $v_1$ $v_3$ $v_4$ $v_7$ $v_8$
вершина $v_3$: $v_1$ $v_2$ $v_4$ $v_6$
вершина $v_4$: $v_2$ $v_3$ $v_5$
вершина $v_5$: $v_4$ $v_6$ $v_7$
вершина $v_6$: $v_1$ $v_3$ $v_5$ $v_7$ $v_8$
вершина $v_7$: $v_2$ $v_5$ $v_6$ $v_8$
вершина $v_8$: $v_2$ $v_6$ $v_7$

Инцидентность графа $G_{28}$:

ребро $e_1$: $(v_1,v_2)$ или $(v_2,v_1)$
ребро $e_2$: $(v_1,v_3)$ или $(v_3,v_1)$
ребро $e_3$: $(v_1,v_6)$ или $(v_6,v_1)$
ребро $e_4$: $(v_2,v_3)$ или $(v_3,v_2)$
ребро $e_5$: $(v_2,v_4)$ или $(v_4,v_2)$
ребро $e_6$: $(v_2,v_7)$ или $(v_7,v_2)$
ребро $e_7$: $(v_2,v_8)$ или $(v_8,v_2)$
ребро $e_8$: $(v_3,v_4)$ или $(v_4,v_3)$
ребро $e_9$: $(v_3,v_6)$ или $(v_6,v_3)$
ребро $e_{10}$: $(v_4,v_5)$ или $(v_5,v_4)$
ребро $e_{11}$: $(v_5,v_6)$ или $(v_6,v_5)$
ребро $e_{12}$: $(v_5,v_7)$ или $(v_7,v_5)$
ребро $e_{13}$: $(v_6,v_7)$ или $(v_7,v_6)$
ребро $e_{14}$: $(v_6,v_8)$ или $(v_8,v_6)$
ребро $e_{15}$: $(v_7,v_8)$ или $(v_8,v_7)$

Базовые реберные разрезы графа $G_{28}$:

$w_0(e_1) = \{e_2,e_3,e_4,e_5,e_6,e_7\}$;
$w_0(e_2) = \{e_1,e_3,e_4,e_8,e_9\}$;
$w_0(e_3) = \{e_1,e_2,e_9,e_{11},e_{13},e_{14}\}$;
$w_0(e_4) = \{e_1,e_2,e_5,e_6,e_7,e_8,e_9\}$;
$w_0(e_5) = \{e_1,e_4,e_6,e_7,e_8,e_{10}\}$;
$w_0(e_6) = \{e_1,e_4,e_5,e_7,e_{12},e_{13},e_{15}\}$;
$w_0(e_7) = \{e_1,e_4,e_5,e_6,e_{14},e_{15}\}$;
$w_0(e_8) = \{e_2,e_4,e_5,e_9,e_{10}\}$;
$w_0(e_9) = \{e_2,e_3,e_4,e_8,e_{11},e_{13},e_{14}\}$;
$w_0(e_{10}) = \{e_5,e_8,e_{11},e_{12}\}$;
$w_0(e_{11}) = \{e_3,e_9,e_{10},e_{12},e_{13},e_{14}\}$;
$w_0(e_{12}) = \{e_6,e_{10},e_{11},e_{13},e_{15}\}$;
$w_0(e_{13}) = \{e_3,e_6,e_9,e_{11},e_{12},e_{14},e_{15}\}$;
$w_0(e_{14}) = \{e_3,e_7,e_9,e_{11},e_{13},e_{15}\}$;
$w_0(e_{15}) = \{e_6,e_7,e_{12},e_{13},e_{14}\}$.

Кортеж весов ребер уровня 0: $\xi(w(l_0))$ = <6,5,7,7,6,7,6,5,7,4,6,5,7,6,5>;
Кортеж весов ребер уровня 1: $\xi(w(l_1))$ = <9,5,7,8,9,6,7,7,6,10,9,7,8,9,5>;
Кортеж весов ребер уровня 2: $\xi(w(l_2))$ = <8,7,4,7,5,7,4,6,7,8,5,6,7,8,7>;
Кортеж весов ребер уровня 3: $\xi(w(l_3))$ = <10,6,10,4,10,4,10,10,4,0,10,10,4,10,6>
Кортеж весов ребер уровня 4: $\xi(w(l_4))$ = <8,10,8,10,10,10,8,8,10,0,10,8,10,8,10>;
Кортеж весов ребер уровня 5: $\xi(w(l_5))$ = <0,8,0,8,8,8,0,0,8,0,8,0,8,0,8>.
Общий кортеж весов ребер : $\xi(w)$ = <41,41,35,44,48,42,35,36,42,22,48,36,44,41,41>;
Кортеж весов вершин : $\varsigma(G_{28})$ = <117,210,163,106,106,210,163,117>.

Из кортежа весов вершин следует, что в графе $G_{28}$ имеются равные веса вершин $o_1 = \{v_1,v_8\}$, $o_2 = \{v_2,v_6\}$, $o_3 = \{v_3,v_7\}$, $o_4 = \{v_4,v_5\}$.

Рассмотрим кортеж весов верщин только для 0-го и 1-го уровней:

Кортеж весов ребер уровня 0: <6,5,7,7,6,7,6,5,7,4,6,5,7,6,5>;
Кортеж весов ребер уровня 1: <9,5,7,8,9,6,7,7,6,10,9,7,8,9,5>.

Кортеж суммарных весов ребер: <15,10,13,15,15,13,13,12,13,14,15,12,15,15,10>.

Кортеж весов вершин только для 0-го и 1-го уровней:

<38,71,50,41,41,71,50,38>.

Рассмотрим кортеж весов верщин только для 2-го и 3-го уровней:

Кортеж весов ребер уровня 2: <8,7,4,7,5,7,4,6,7,8,5,6,7,8,7>;
Кортеж весов ребер уровня 3: <10,6,10,4,10,4,10,10,4,0,10,10,4,10,6>

Кортеж суммарных весов ребер: <18,13,14,11,15,11,14,16,11,8,15,16,11,18,13>.

Кортеж весов вершин только для 2-го и 3-го уровней:

<45,69,51,39,39,69,51,45>.

Рассмотрим кортеж весов верщин только для 4-го и 5-го уровней:

Кортеж весов ребер уровня 4: <8,10,8,10,10,10,8,8,10,0,10,8,10,8,10>;
Кортеж весов ребер уровня 5: <0,8,0,8,8,8,0,0,8,0,8,0,8,0,8>.

Кортеж суммарных весов ребер: <8,18,8,18,18,18,8,8,18,0,18,8,18,8,18>.

Кортеж весов вершин только для 4-го и 5-го уровней:

<34,70,62,26,26,70,62,34>.

Рассмотрим кортеж весов ребер для 0-го уровня:

Кортеж весов ребер уровня 0: <6,5,7,7,6,7,6,5,7,4,6,5,7,6,5>.

Кортеж весов вершин только для 0-го уровня: <17,32,24,15,15,32,24,17>.

Рассмотрим кортеж весов ребер для 1-го уровня:

Кортеж весов ребер уровня 1: <9,5,7,8,9,6,7,7,6,10,9,7,8,9,5>.

Кортеж весов вершин только для 1-го уровня: <21,39,26,26,26,39,26,21>.

Как видно из примера, равенство весов для пар вершин $o_1 = \{v_1,v_8\}$, $o_2 = \{v_2,v_6\}$, $o_3 = \{v_3,v_7\}$, $o_4 = \{v_4,v_5\}$ сохраняется для всех случаев. Перестановка вершин во множествах вершин с равным весом, характеризует изоморфный граф с другой нумераций вершин (см. рис. 7.6.b).

Результат сравнения величин, в кортеже весов вершин для одиночного уровня, может быть не совсем корректно определен (см. кортеж вершин 1-го уровня). Поэтому, для вычисления весов вершин в спектре реберных разрезов, лучше всего пользоваться парами уровней. А с точки зрения вычислительной сложности спектра реберных разрезов, можно обойтись только первым и вторым уровнем. Данное правило справедливо для всего множества несепарабельных графов.

Построим множество изометрических циклов графа $G_{28}$:

$c_1 = \{e_1,e_2,e_4\}$; $c_2 = \{e_1,e_3,e_6,e_{13}\}$;
$c_3 = \{e_1,e_3,e_7,e_{14}\}$; $c_4 = \{e_2,e_3,e_9\}$;
$c_5 = \{e_4,e_5,e_8\}$; $c_6 = \{e_4,e_6,e_9,e_{13}\}$;
$c_7 = \{e_4,e_7,e_9,e_{14}\}$; $c_8 = \{e_5,e_6,e_{10},e_{12}\}$;
$c_9 = \{e_6,e_7,e_{15}\}$; $c_{10} = \{e_8,e_9,e_{10},e_{11}\}$;
$c_{11} = \{e_{11},e_{12},e_{13}\}$; $c_{12} = \{e_{13},e_{14},e_{15}\}$.

Вершинами реберного графа $L(G_{28})$ являются ребра графа $G_{28}$. Запишем смежность реберного графа $L(G_{12})$ в виде ребер графа $G_{28}$ [23]:

$v_1(L(G_{28})) = e_1$: $e_2$ $e_3$ $e_4$ $e_5$ $e_6$ $e_7$
$v_2(L(G_{28})) = e_2$: $e_1$ $e_3$ $e_4$ $e_8$ $e_9$
$v_3(L(G_{28})) = e_3$: $e_1$ $e_2$ $e_9$ $e_{11}$ $e_{13}$ $e_{14}$
$v_4(L(G_{28})) = e_4$: $e_1$ $e_2$ $e_5$ $e_6$ $e_7$ $e_8$ $e_9$
$v_5(L(G_{28})) = e_5$: $e_1$ $e_4$ $e_6$ $e_7$ $e_8$ $e_{10}$
$v_6(L(G_{28})) = e_6$: $e_1$ $e_4$ $e_5$ $e_7$ $e_{12}$ $e_{13}$ $e_{15}$
$v_7(L(G_{28})) = e_7$: $e_1$ $e_4$ $e_5$ $e_6$ $e_{14}$ $e_{15}$
$v_8(L(G_{28})) = e_8$: $e_2$ $e_4$ $e_5$ $e_9$ $e_{10}$
$v_9(L(G_{28})) = e_9$: $e_2$ $e_3$ $e_4$ $e_8$ $e_{11}$ $e_{13}$ $e_{14}$
$v_{10}(L(G_{28})) = e_{10}$: $e_5$ $e_8$ $e_{11}$ $e_{12}$
$v_{11}(L(G_{28})) = e_{11}$: $e_3$ $e_9$ $e_{10}$ $e_{12}$ $e_{13}$ $e_{14}$
$v_{12}(L(G_{28})) = e_{12}$: $e_6$ $e_{10}$ $e_{11}$ $e_{13}$ $e_{15}$
$v_{13}(L(G_{28})) = e_{13}$: $e_3$ $e_6$ $e_9$ $e_{11}$ $e_{12}$ $e_{14}$ $e_{15}$
$v_{14}(L(G_{28})) = e_{14}$: $e_3$ $e_7$ $e_9$ $e_{11}$ $e_{13}$ $e_{15}$
$v_{15}(L(G_{28})) = e_{15}$: $e_6$ $e_7$ $e_{12}$ $e_{13}$ $e_{14}$

Количество ребер в реберном графе $L(G_{28})) = 44$.

Количество изометрических циклов в реберном графе $L(G_{28})) = 50$.

Изометрические циклы графа $L(G_{28})$), относительно его вершин, можно записать в виде ребер графа $G_{28}$:

$\varphi: c_1^L \rightarrow \{e_1,e_2,e_3\}$; $\varphi: c_2^L \rightarrow \{e_1,e_2,e_4\}$;

$\varphi: c_3^L \rightarrow \{e_1,e_2,e_5,e_8\}$; $\varphi: c_4^L \rightarrow \{e_1,e_3,e_4,e_9\}$;

$\varphi: c_5^L \rightarrow \{e_1,e_3,e_5,e_{10},e_{11}\}$; $\varphi: c_6^L \rightarrow \{e_1,e_3,e_6,e_{13}\}$;

$\varphi: c_7^L \rightarrow \{e_1,e_3,e_7,e_{14}\}$; $\varphi: c_8^L \rightarrow \{e_1,e_4,e_5\}$;

$\varphi: c_9^L \rightarrow \{e_1,e_4,e_6\}$; $\varphi: c_{10}^L \rightarrow \{e_1,e_4,e_7\}$;

$\varphi: c_{11}^L \rightarrow \{e_1,e_5,e_6\}$; $\varphi: c_{12}^L \rightarrow \{e_1,e_5,e_7\}$;

$\varphi: c_{13}^L \rightarrow \{e_1,e_6,e_7\}$; $\varphi: c_{14}^L \rightarrow \{e_2,e_3,e_9\}$;

$\varphi: c_{15}^L \rightarrow \{e_2,e_4,e_8\}$; $\varphi: c_{16}^L \rightarrow \{e_2,e_4,e_9\}$;

$\varphi: c_{17}^L \rightarrow \{e_2,e_8,e_9\}$; $\varphi: c_{18}^L \rightarrow \{e_3,e_9,e_{11}\}$;

$\varphi: c_{19}^L \rightarrow \{e_3,e_9,e_{13}\}$; $\varphi: c_{20}^L \rightarrow \{e_3,e_9,e_{14}\}$;

$\varphi: c_{21}^L \rightarrow \{e_3,e_{11},e_{13}\}$; $\varphi: c_{22}^L \rightarrow \{e_3,e_{11},e_{14}\}$;

$\varphi: c_{23}^L \rightarrow \{e_3,e_{13},e_{14}\}$; $\varphi: c_{24}^L \rightarrow \{e_4,e_5,e_6\}$;

$\varphi: c_{25}^L \rightarrow \{e_4,e_5,e_7\}$; $\varphi: c_{26}^L \rightarrow \{e_4,e_5,e_8\}$;

$\varphi: c_{27}^L \rightarrow \{e_4,e_6,e_7\}$; $\varphi: c_{28}^L \rightarrow \{e_4,e_6,e_9,e_{13}\}$;

$\varphi: c_{29}^L \rightarrow \{e_4,e_7,e_9,e_{14}\}$; $\varphi: c_{30}^L \rightarrow \{e_4,e_8,e_9\}$;

$\varphi: c_{31}^L \rightarrow \{e_5,e_6,e_7\}$; $\varphi: c_{32}^L \rightarrow \{e_5,e_6,e_{10},e_{12}\}$;

$\varphi: c_{33}^L \rightarrow \{e_5,e_7,e_{10},e_{11},e_{14}\}$; $\varphi: c_{34}^L \rightarrow \{e_5,e_8,e_{10}\}$;

$\varphi: c_{35}^L \rightarrow \{e_6,e_7,e_{13},e_{14}\}$; $\varphi: c_{36}^L \rightarrow \{e_6,e_7,e_{15}\}$;

$\varphi: c_{37}^L \rightarrow \{e_6,e_{12},e_{13}\}$; $\varphi: c_{38}^L \rightarrow \{e_6,e_{12},e_{15}\}$;

$\varphi: c_{39}^L \rightarrow \{e_6,e_{13},e_{15}\}$; $\varphi: c_{40}^L \rightarrow \{e_7,e_{14},e_{15}\}$;

$\varphi: c_{41}^L \rightarrow \{e_8,e_9,e_{10},e_{11}\}$; $\varphi: c_{42}^L \rightarrow \{e_9,e_{11},e_{13}\}$;

$\varphi: c_{43}^L \rightarrow \{e_9,e_{11},e_{14}\}$; $\varphi: c_{44}^L \rightarrow \{e_9,e_{13},e_{14}\}$;

$\varphi: c_{45}^L \rightarrow \{e_{10},e_{11},e_{12}\}$; $\varphi: c_{46}^L \rightarrow \{e_{11},e_{12},e_{13}\}$;

$\varphi: c_{47}^L \rightarrow \{e_{11},e_{12},e_{14},e_{15}\}$; $\varphi: c_{48}^L \rightarrow \{e_{11},e_{13},e_{14}\}$;

$\varphi: c_{49}^L \rightarrow \{e_{12},e_{13},e_{15}\}$; $\varphi: c_{50}^L \rightarrow \{e_{13},e_{14},e_{15}\}$.

Множество изометрических циклов $C_\tau^L$ графа $L(G_{28})$ можно построить как объединение следующих трех соответствующих подмножеств суграфов графа $G_{28}$:

$$S_R^L \cup C_\tau^L \cup C_d^L$$

Подмножество изометрических циклов соответствующих множеству центральных разрезов графа $G_{28}$, будем обозначать как $S_R^L$:

$\varphi: c_1^L \rightarrow s_1 = \{e_1,e_2,e_3\}$;

$\varphi: (c_8^L \cup c_9^L \cup c_{10}^L \cup c_{11}^L \cup c_{12}^L \cup c_{13}^L \cup c_{24}^L \cup c_{25}^L \cup c_{27}^L \cup c_{31}^L) \rightarrow s_2 = \{e_1,e_4,e_5,e_6,e_7\}$;

$\varphi: (c_{15}^L \cup c_{16}^L \cup c_{17}^L \cup c_{30}^L) \rightarrow s_3 = \{e_2,e_4,e_8,e_9\}$;

$\varphi: c_{34}^L \rightarrow s_4 = \{e_5,e_8,e_{10}\}$;

$\varphi: c_{45}^L \rightarrow s_5 = \{e_{10},e_{11},e_{12}\}$;

$\varphi: (c_{18}^L \cup c_{19}^L \cup c_{20}^L \cup c_{21}^L \cup c_{22}^L \cup c_{23}^L \cup c_{42}^L \cup c_{43}^L \cup c_{44}^L \cup c_{48}^L) \rightarrow s_6 = \{e_3,e_9,e_{11},e_{13},e_{14}\}$;

$\varphi:(c^L_{37} \cup c^L_{38} \cup c^L_{39} \cup c^L_{49}) \to s_7 = \{e_6,e_{12},e_{13},e_{15}\}$;

$\varphi: c^L_{40} \to s_8 = \{e_7,e_{14},e_{15}\}$;

Подмножество изометрических циклов графа $L(G_{28})$ соответствующих множеству изометрических циклов графа $G_{28}$, будем обозначать как $C^L_\tau$:

$\varphi: c^L_2 \to c_1 = \{e_1,e_2,e_4\}$;

$\varphi: c^L_6 \to c_2 = \{e_1,e_3,e_6,e_{13}\}$;

$\varphi: c^L_7 \to c_3 = \{e_1,e_3,e_7,e_{14}\}$;

$\varphi: c^L_{14} \to c_4 = \{e_2,e_3,e_9\}$;

$\varphi: c^L_{26} \to c_5 = \{e_4,e_5,e_8\}$;

$\varphi: c^L_{28} \to c_6 = \{e_4,e_6,e_9,e_{13}\}$;

$\varphi: c^L_{29} \to c_7 = \{e_4,e_7,e_9,e_{14}\}$;

$\varphi: c^L_{32} \to c_8 = \{e_5,e_6,e_{10},e_{12}\}$;

$\varphi: c^L_{36} \to c_9 = \{e_6,e_7,e_{15}\}$;

$\varphi: c^L_{41} \to c_{10} = \{e_8,e_9,e_{10},e_{11}\}$.

$\varphi: c^L_{46} \to c_{11} = \{e_{11},e_{12},e_{13}\}$;

$\varphi: c^L_{50} \to c_{12} = \{e_{13},e_{14},e_{15}\}$.

Подмножество изометрических циклов соответствующих множеству суграфов графа $G_{12}$, будем обозначать как $C^L_d$:

$\varphi: c^L_3 \to c_1 \oplus c_5 = c_4 \oplus c_8 \oplus c_9 \oplus c_{10} \oplus c_{11} \oplus c_{12} = \{e_1,e_2,e_5,e_8\}$;

$\varphi: c^L_4 \to c_1 \oplus c_4 = c_5 \oplus c_8 \oplus c_9 \oplus c_{10} \oplus c_{11} \oplus c_{12} = \{e_1,e_3,e_4,e_9\}$;

$\varphi: c^L_5 \to c_1 \oplus c_4 \oplus c_5 \oplus c_{10} = c_3 \oplus c_8 \oplus c_9 \oplus c_{11} \oplus c_{12} = \{e_1,e_3,e_5,e_{10},e_{11}\}$.

$\varphi: c^L_{33} \to c_1 \oplus c_3 \oplus c_4 \oplus c_5 \oplus c_{10} = c_8 \oplus c_9 \oplus c_{11} \oplus c_{12} = \{e_5,e_7,e_{10},e_{11},e_{14}\}$;

$\varphi: c^L_{35} \to c_9 \oplus c_{12} = c_1 \oplus c_4 \oplus c_5 \oplus c_8 \oplus c_{10} \oplus c_{11} = \{e_6,e_7,e_{13},e_{14}\}$;

$\varphi: c^L_{47} \to c_{11} \oplus c_{12} = c_1 \oplus c_3 \oplus c_4 \oplus c_5 \oplus c_8 \oplus c_9 \oplus c_{10} = \{e_{11},e_{12},e_{14},e_{15}\}$.

Определим кортеж весов ребер и вершин исходя из множества изометрических циклов реберного графа $L(G_{28})$.

Кортеж весов ребер: <13,7,12,14,12,14,12,7,14,6,12,7,14,13,7>.

Кортеж весов вершин: <32,65,42,25,25,65,42,32>.

Заметим, что равные веса вершин определяются следующими парами вершин $o_1 = \{v_1,v_8\}$, $o_2 = \{v_2,v_6\}$, $o_3 = \{v_3,v_7\}$, $o_4 = \{v_4,v_5\}$.

Определим кортеж весов ребер для множества изометрических циклов графа $G_{28}$:

$\xi_\tau(G_{28}) =$ <3,2,3,4,2,4,3,2,4,2,2,2,4,3,2>.

Кортеж весов ребер порождает кортеж весов вершин графа: <8,16,12,6,6,16,12,8>.

И снова, равные веса вершин определяются следующими парами вершин $o_1 = \{v_1,v_8\}$, $o_2 = \{v_2,v_6\}$, $o_3 = \{v_3,v_7\}$, $o_4 = \{v_4,v_5\}$.

Желательно, для множества изометрических циклов найти представление в зависимости от ребер графа. Такое представление существует и определяется как кольцевая сумма изометрических циклов проходящих по ребру.

Множество базовых реберных циклов графа $G_{28}$:

$\tau_0(e_1) = \{e_2,e_3,e_4,e_6,e_{13}\}$;

$\tau_0(e_2) = \{e_1,e_3,e_4,e_9\}$;

$\tau_0(e_3) = \{e_1,e_2,e_6,e_9,e_{13}\}$;

$\tau_0(e_4) = \{e_1,e_2,e_5,e_6,e_7,e_8,e_{13},e_{14}\}$;

$\tau_0(e_5) = \{e_4,e_6,e_8,e_{10},e_{12}\}$;

$\tau_0(e_6) = \{e_1,e_3,e_4,e_5,e_7,e_9,e_{10},e_{12},e_{15}\}$;

$\tau_0(e_7) = \{e_4,e_6,e_9,e_{14},e_{15}\}$;

$\tau_0(e_8) = \{e_4,e_5,e_9,e_{10},e_{11}\}$;

$\tau_0(e_9) = \{e_2,e_3,e_6,e_7,e_8,e_{10},e_{11},e_{13},e_{14}\}$;

$\tau_0(e_{10}) = \{e_5,e_6,e_8,e_9,e_{11},e_{12}\}$;

$\tau_0(e_{11}) = \{e_8,e_9,e_{10},e_{12},e_{13}\}$;

$\tau_0(e_{12}) = \{e_5,e_6,e_{10},e_{11},e_{13}\}$;

$\tau_0(e_{13}) = \{e_1,e_3,e_4,e_9,e_{11},e_{12},e_{14},e_{15}\}$;

$\tau_0(e_{14}) = \{e_4,e_7,e_9,e_{13},e_{15}\}$;

$\tau_0(e_{15}) = \{e_6,e_7,e_{13},e_{14}\}$.

Элементы множества базовых циклов графа можно представить в виде, симметричной относительно главной диагонали, матрицы базовых реберных циклов $T(G_{28})$ и рассматривать ее как нильпотентный оператор. Множества нильпотентных операторов различных степеней определяют спектр реберных циклов графа. Но так как цифровой инвариант графа включает только множество изометрических циклов, то по аналогии, спектр реберных циклов может состоять только из одного базового уровня.

Определим кортеж весов ребер для спектра реберных циклов:

$\xi(\tau(G_{28})) = <5,4,5,8,5,9,5,5,9,6,5,5,8,5,4>$.

Кортеж весов вершин: $\zeta(\tau(G_{28})) = <14,32,26,16,16,32,26,14>$.

И вновь, равные веса вершин определяются следующими парами вершин $o_1 = \{v_1,v_8\}$, $o_2 = \{v_2,v_6\}$, $o_3 = \{v_3,v_7\}$, $o_4 = \{v_4,v_5\}$.

Кортеж весов вершин для спектра реберных разрезов графа:

$\zeta(w(G_{28})) = <38,71,50,41,41,71,50,38>$.

Кортеж весов вершин цифрового инварианта графа:

$\zeta_L(G_{28}) = <32,65,42,25,25,65,42,32>$.

Кортеж весов вершин для множества изометрических циклов графа:

$\zeta(C_\tau(G_{28}))$ = <8,16,12,6,6,16,12,8>.

Кортеж весов вершин для спектра реберных циклов графа:

$\zeta(\tau(G_{28}))$ = <14,32,26,16,16,32,26,14>.

Из сравнения результатов величин кортежей весов для различного представления графа $G_{28}$ следует, что, несмотря на их различие, существует устойчивые пары вершин одинакового веса. Для графа $G_{28}$ - это следующие подмножества вершин $o_1 = \{v_1,v_8\}$, $o_2 = \{v_2,v_6\}$, $o_3 = \{v_3,v_7\}$, $o_4 = \{v_4,v_5\}$.

**Определение 7.1.** Способность сохранять симметрическую устойчивость подмножества элементов с равным весом в векторных инвариантах, будем называть *устойчивостью вершин инвариантов.*

**Определение 7.2.** Устойчивое подмножество вершин, имеющие равные веса соответствуют *орбитами* графа, и обозначаются символом «o».

Из рассмотренного также следует, что для формирования спектра реберных разрезов достаточно иметь первые два уровня $w_0(G)$ и $w_1(G)$, а для формирования спектра реберных циклов достаточно одного базового уровня $\tau_0$(G), что соответствует цифровому векторному инварианта графа.

На данном примере графа $G_{28}$ рассмотрены основные свойства структуры графа, общие для всего множества несепарабельных графов.

## Комментарии

В данной главе рассматривается природное свойство векторного инварианта графа, симметричная устойчивость весов элементов. Симметричная устойчивость векторных инвариантов позволяет выделить подмножества вершин с равным весом. Следует заметить, что это свойство распространяется на все виды векторных инвариантов.

# Глава 8. ОПРЕДЕЛЕНИЕ ИЗОМОРФИЗМА ГРАФОВ

## 8.1. Инварианты графа

Рассмотрим различные классы связных неориентированных графов без петель и кратных ребер. Для определения изоморфизма графов будем применять следующие инварианты:

1. векторный интегральный инвариант графа G;
2. векторный цифровой инвариант реберного графа L(G).

В целях детального рассмотрения структур графа, наровне с векторным интегральным инвариантом графа будем применять векторный инвариант реберных разрезов графа и векторный инвариант уровня спектра реберных разрезов.

## 8.2. Класс сильно регулярных графов

***Пример 8.1.*** Рассмотрим сильно регулярные графы Lattice graph с параметрами srg(16,6,2,2) и Shrikhande graph с параметрами srg(16,6,2,2).

Количество уровней в спектре реберных разрезов графа $G_1$ =2.

Кортеж весов ребер графа $G_1$: $\xi_w(G_{29})$ = <48×26>.

Кортеж весов вершин графа $G_1$: $\zeta_w(G_{29})$ = <16×156>.

Векторный инвариант реберных разрезов графа $F_w(\xi(G_{29}))$ & $F_w(\zeta(G_{29}))$ = (48×26) & (16×156).

Количество уровней в спектре реберных разрезов графа $G_{30}$ =2.

Кортеж весов ребер графа $G_2$: $\xi_w(G_{30})$ = <48×34>.

Кортеж весов вершин графа $G_2$: $\zeta_w(G_{30})$ = <16×204>.

Инвариант реберных разрезов графа $F_w(\xi(G_{30}))$ & $F_w(\zeta(G_{30}))$ = (48×26) & (16×156).

Сравнивая базовые инварианты графов, приходим к выводу: графы $G_{29}$ и $G_{30}$ не изоморфны. Для распознования различения графов, достаточно произвести сравнение только векторных инвариантов реберных спектров графов.

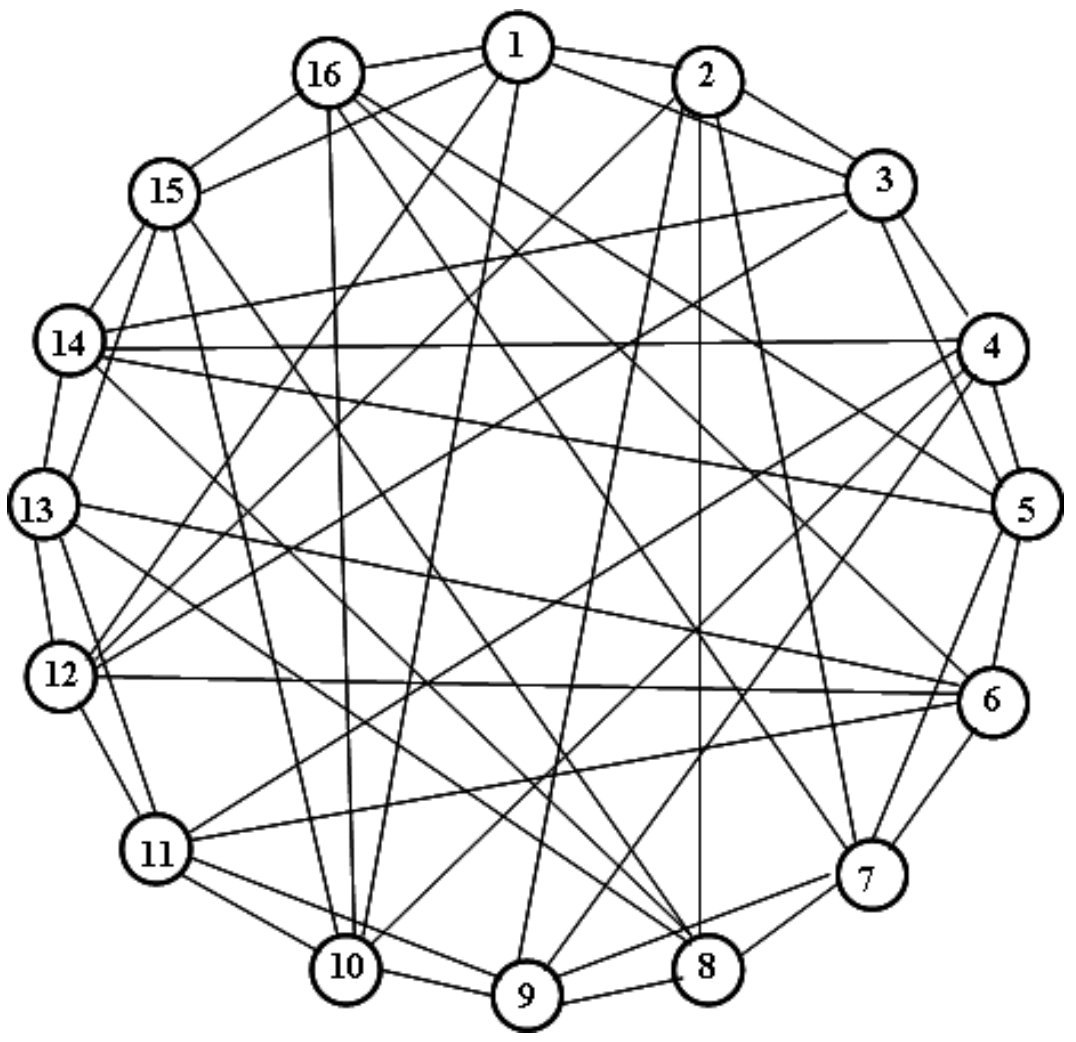

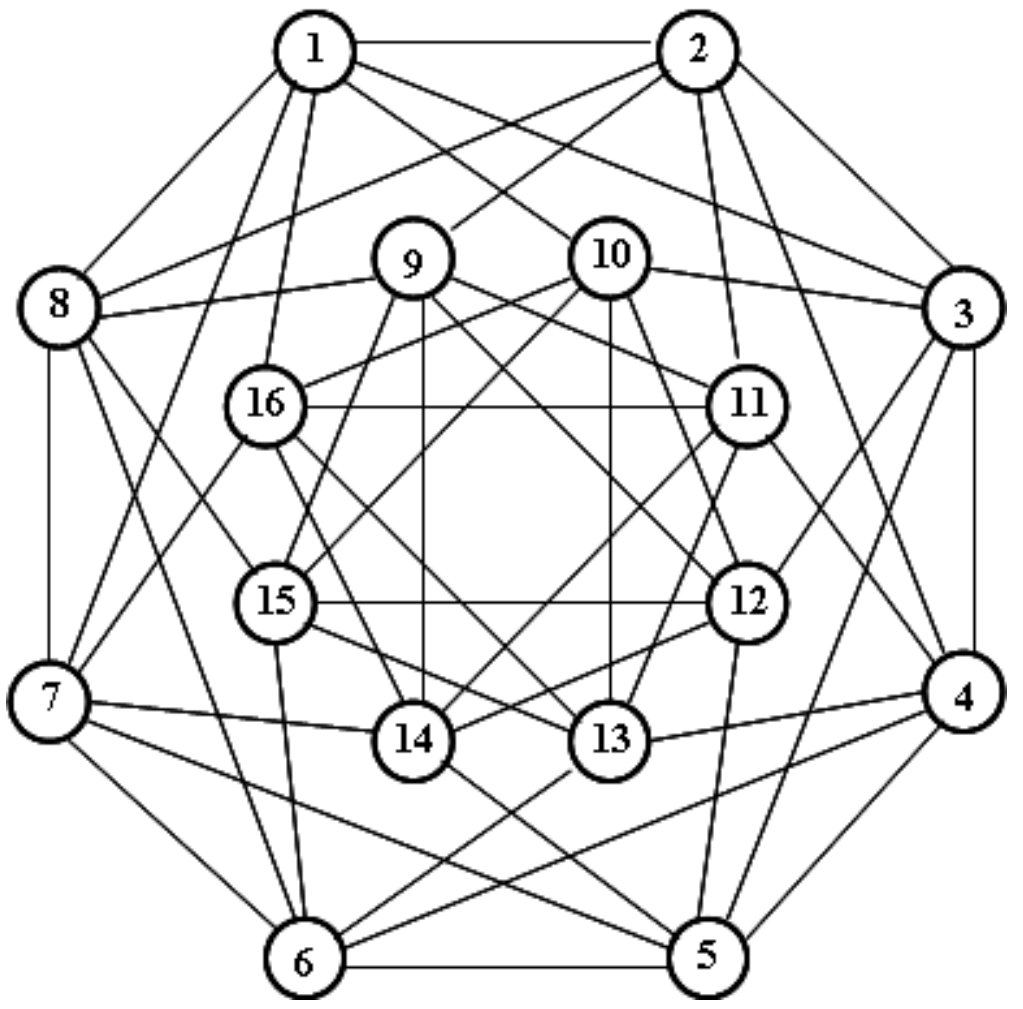

Рис. 8.1. Сильно регулярный граф $G_{29}$ (lattice graph) с параметрами srg(16,6,2,2).

Рис. 8.2. Сильно регулярный граф $G_{30}$ (shrikhande graph) с параметрами srg(16,6,2,2).

***Пример 8.2.*** Рассмотрим следующие сильно регулярные графы:

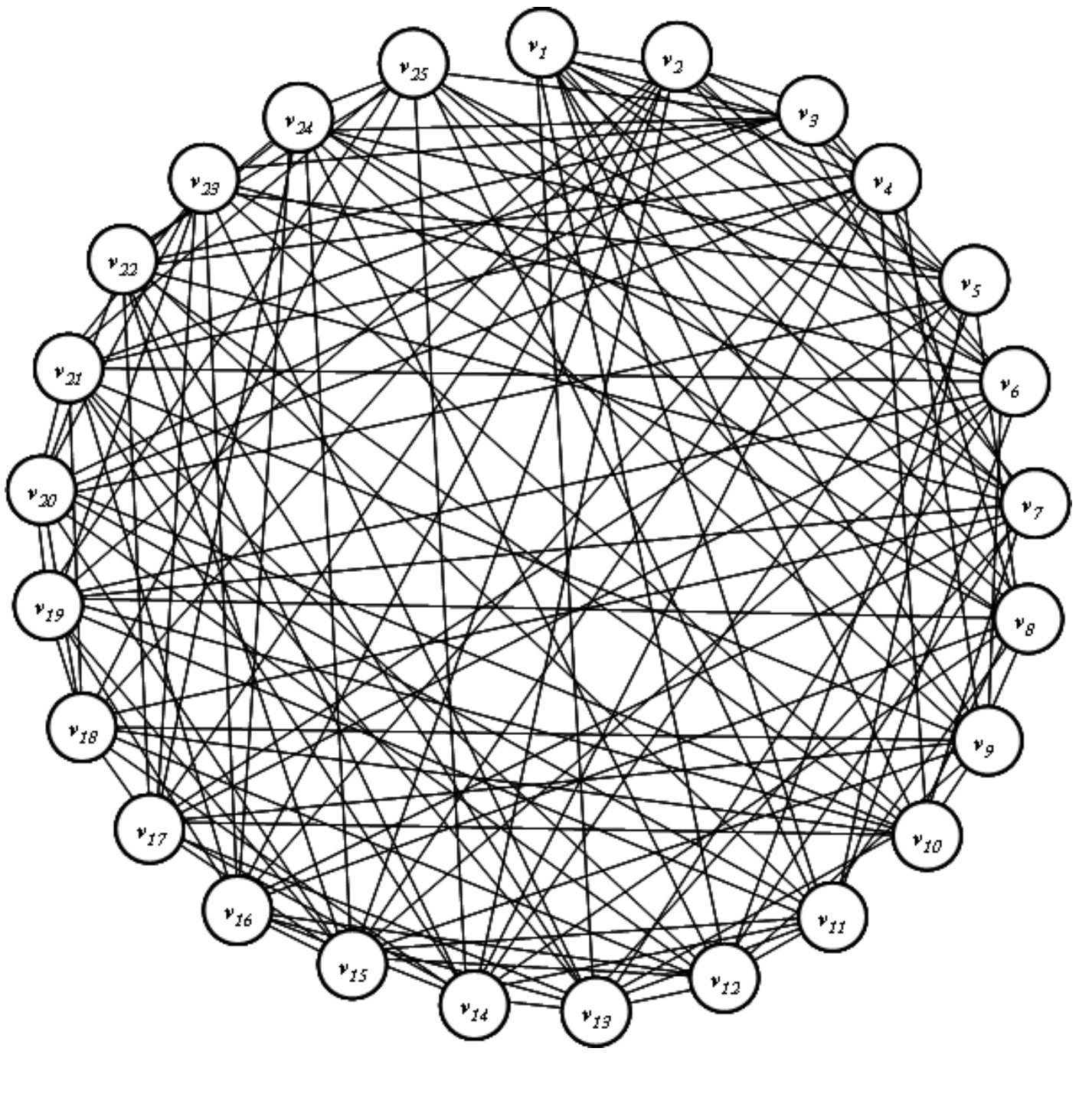

Рис. 8.3. Сильно регулярный граф $G_{31}$ с параметрами srg(25,12,5,6).

```
0111111111111000000000000
1011111000000111111000000
1101111000000000000111111
1110000111000111000111000
1110000100110100110100110
1110000010101010101010101
1110000001011001011001011
1001100011100101001000111
1001010101010010110010011
1001001110001000111101100
1000110100011110001101001
1000101010101011010110010
1000011001110101100011100
0101100100101001110011001
0101010010110001101101010
0101001100011110001010110
0100110011001110010001110
0100101011010100101110001
0100011101100011010100101
0011100001110010011011100
0011010010011101010100101
0011001001101110100100011
0010110101001001101110010
0010101110010011100001101
0010011110100100011011010
```

Рис. 8.4. Матрица смежностей графа $G_{31}$.

Количество уровней в спектре реберных разрезов графа $G_{31}$ =2.

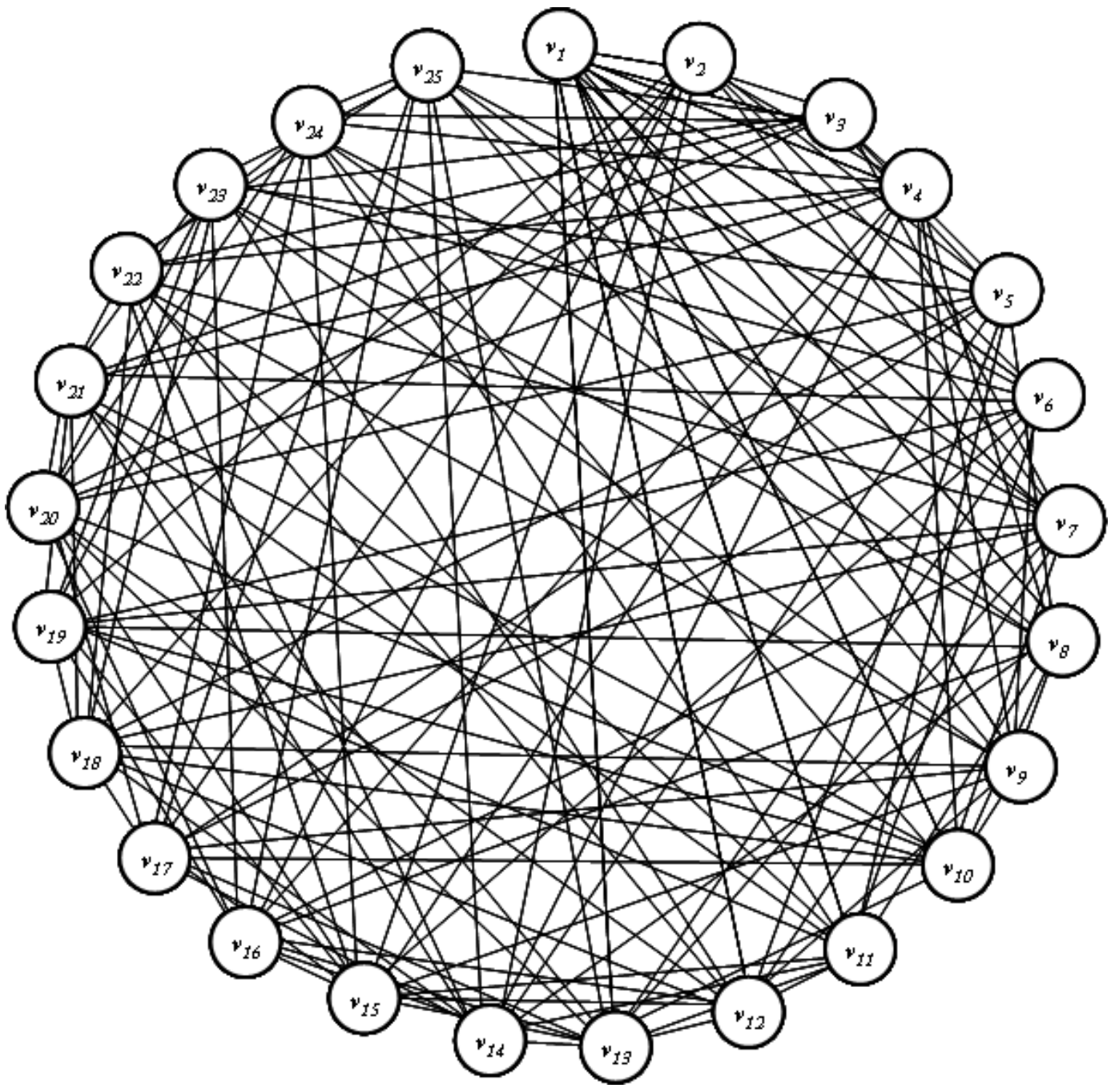

```
0111111111111000000000000
1011111000000111111000000
1101111000000000000111111
1110000111000111000111000
1110000100110100110100110
1110000010101010101010101
1110000001011001011001011
1001100011100101001000111
1001010101010010110001101
1001001110001000111110010
1000110100011110001011010
1000101010101011010101100
1000011001110101100110001
0101100100101001110011001
0101010010110001101101010
0101001100011110001100101
0100110011001110010100011
0100101011010100101011100
0100011101100011010010110
0011100001011011100010110
0011010001101100011101100
0011001010110110010010011
0010110110010001011110001
0010101101100010101101001
0010011110001101100001110
```

Рис. 8.5. Сильно регулярный граф $G_{32}$ с параметрами srg(25,12,5,6).

Рис. 8.6. Матрица смежностей графа $G_{32}$.

Инвариант реберных разрезов графа $F_w(\xi(G_{32}))$ & $F_w(\zeta(G_{32})) = (77\times92,68\times100,5\times108)$ &

& $(21\times1128,2\times1176,2\times1200)$.

Количество уровней в спектре реберных разрезов графа $G_{32}$ =2.

Инвариант реберных разрезов графа $F_w(\xi(G_{32}))$ & $F_w(\zeta(G_{32}) = (85\times92,64\times100,1\times108)$ &

& $(7\times1128,17\times1152,1\times1176)$.

Сравнивая инварианты реберных разрезов графов, приходим к выводу: графы $G_3$ и $G_4$ не изоморфны.

## 8.3. Класс изоспектральных графов

***Пример 8.3.*** Рассмотрим следующую пару изоспектральных графов (рис. 8.7 и 8.8).

Общим для указанных графов спектром собственных значений для матрицы смежностей является вектор (4,2,2,2,0,0,0,-2,-2,-2,-2,-2).

Инвариант реберных разрезов графа $G_{33}$: $F_w(\xi(G_{33}))$ & $F_w(\zeta(G_{33})) = (24\times30)$ &

& $(12\times120)$.

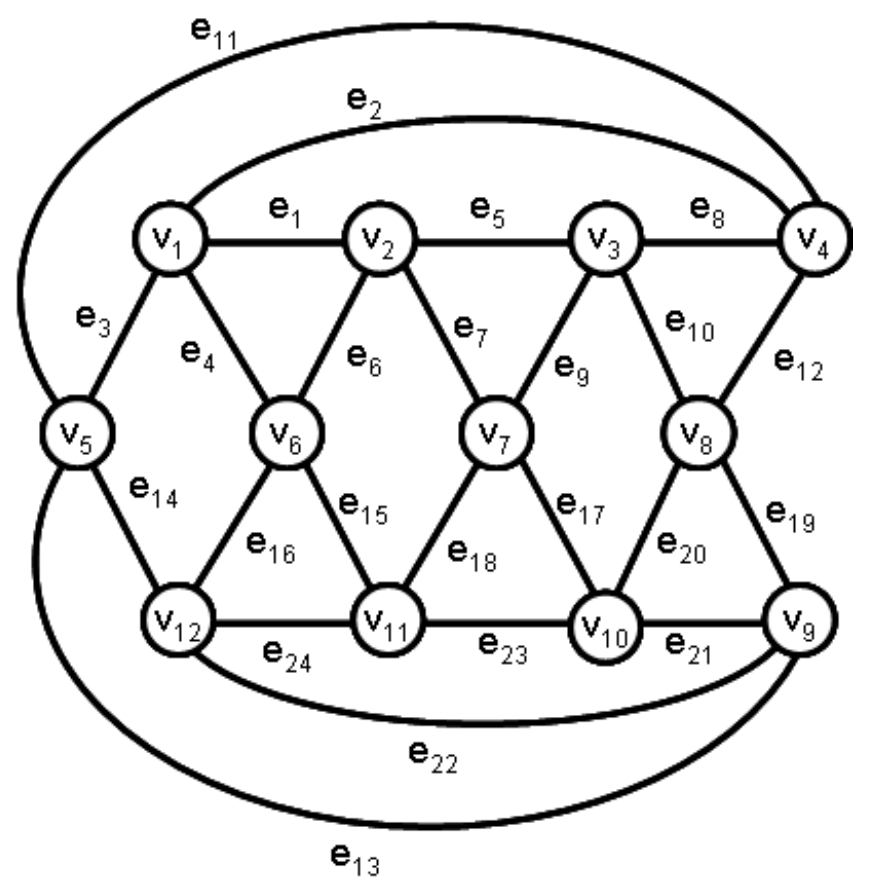

Рис. 8.7. Граф $G_{33}$.

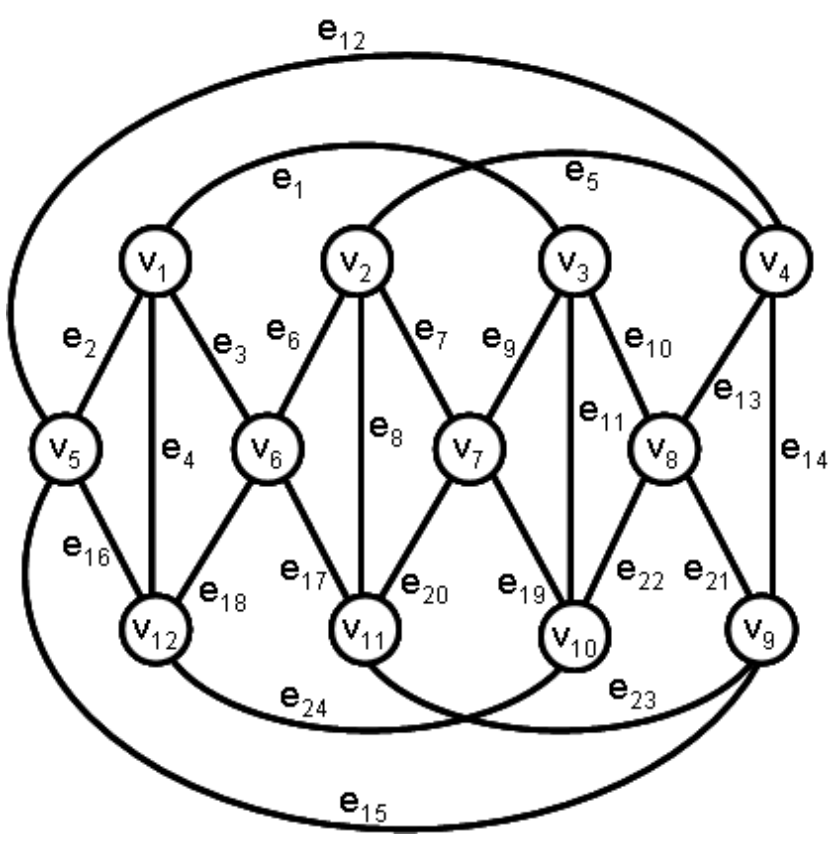

Рис. 8.8. Граф $G_{34}$.

Инвариант реберных разрезов графа $G_{34}$: $F_w(\xi(G_{34}))$ & $F_w(\zeta(G_{34}))$ = (8 × 14, 16 × 26) & & (8 × 80, 4 × 104).

Изоспектральные графы $G_{33}$ и $G_{34}$ не изоморфны, так как их инварианты спектра реберных разрезов не совпадают.

***Пример 8.4.*** Рассмотрим следующие изоспектральные графы (рис. 8.9 и 8.10).

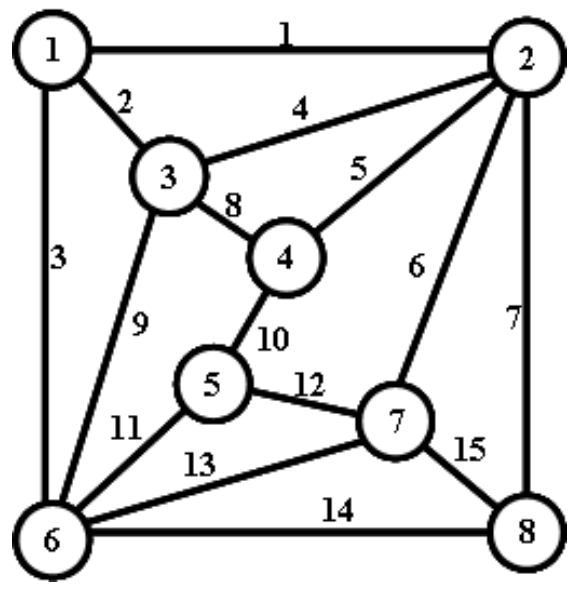

Рис. 8.9. Граф $G_{35}$.

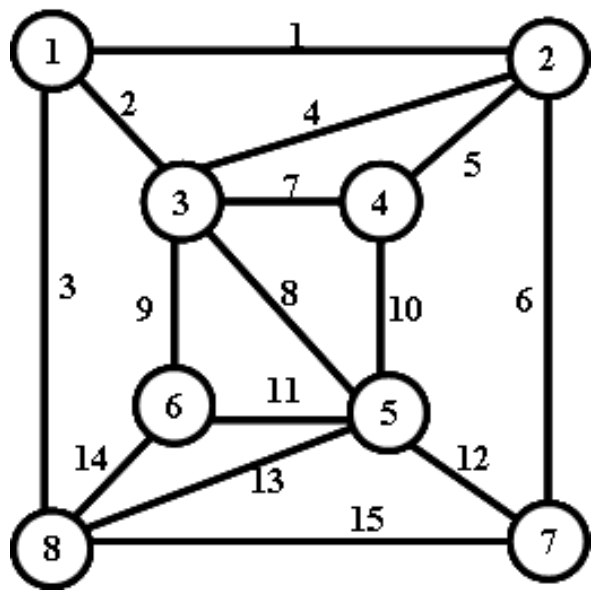

Рис. 8.10. Граф $G_{36}$.

Количество уровней в спектре реберных разрезов графа $G_{35}$ =6.

Кортеж весов ребер уровня 1: $\xi(w(l_0))$=<6,5,7,7,6,7,6,5,7,4,6,5,7,6,5>.

Кортеж весов ребер уровня 2: $\xi(w(l_1))$=<9,5,7,8,9,6,7,7,6,10,9,7,8,9,5>.

Кортеж весов ребер уровня 3: $\xi(w(l_2))$=<8,7,4,7,5,7,4,6,7,8,5,6,7,8,7>.

Кортеж весов ребер уровня 4: $\xi(w(l_3))$=<10,6,10,4,10,4,10,10,4,0,10,10,4,10,6>.

Кортеж весов ребер уровня 5: $\xi(w(l_4))$=<8,10,8,10,10,10,8,8,10,0,10,8,10,8,10>.

Кортеж весов ребер уровня 6: $\xi(w(l_5))$=<0,8,0,8,8,8,0,0,8,0,8,0,8,0,8>.

Суммарный кортеж весов ребер графа $G_{35}$:

$\xi_w(G_7)$ = = <41,41,35,44,48,42,35,36,42,22,48,36,44, 41,41>.

Кортеж весов вершин графа $G_{35}$: $\zeta_w(G_{35})$ = <117,210,163,106,106,210,163,117>.

Инвариант графа $F_w(\xi\, G_{35}))$ & $F_w(\zeta\ (G_{35}))$ =

= (22,35,35,36,36,41,41,41,41,42,42,44,44,48,48) & (106,106,117,117,163,163,210,210).

Количество уровней в спектре реберных разрезов графа $G_{36}$ =6.

Кортеж весов ребер уровня 1: $\xi(w(l_0))$=<6,6,5,7,5,5,6,8,6,6,6,6,7,5,5>.

Кортеж весов ребер уровня 2: $\xi(w(l_1))$=<8,7,6,9,9,6,6,6,8,8,6,7,9,9,8>.

Кортеж весов ребер уровня 3: $\xi(w(l_2))$=<5,7,9,6,9,9,11,12,7,7,11,7,6,9,5>.

Кортеж весов ребер уровня 4: $\xi(w(l_3))$=<8,6,8,10,6,8,8,0,8,8,8,6,10,6,8>.

Кортеж весов ребер уровня 5: $\xi(w(l_4))$=<6,10,6,12,10,6,6,0,6,6,6,10,12,10,6>.

Кортеж весов ребер уровня 6: $\xi(w(l_5))$=<12,12,12,0,12,12,12,0,12,12,12,12,0,12,12>.

Суммарный кортеж весов ребер графа $G_{36}$: $\xi_w(G_{36})$ =

= <44,48,46,44,51,46,49,26,47,47,49,48,44, 51,44>.

Кортеж весов вершин графа $G_{36}$: $\zeta_w(G_{36})$ = <138,185,214,147,214,147,138,185>.

Инвариант графа $F_w(\xi\,(G_{36}))$ & $F_w(\zeta\ (G_{36}))$ =

= (26,44,44,44,44,46,46,47,47,48,48,49,49,51,51) & (138,138,147,147,185,185,214,214).

Изоспектральные графы $G_{35}$ и $G_{36}$ не изоморфны, так как их инварианты не совпадают.

## 8.4. Класс плоских графов

***Пример 8.5.*** Сравним графы $G_{15}$ и $G_{16}$.

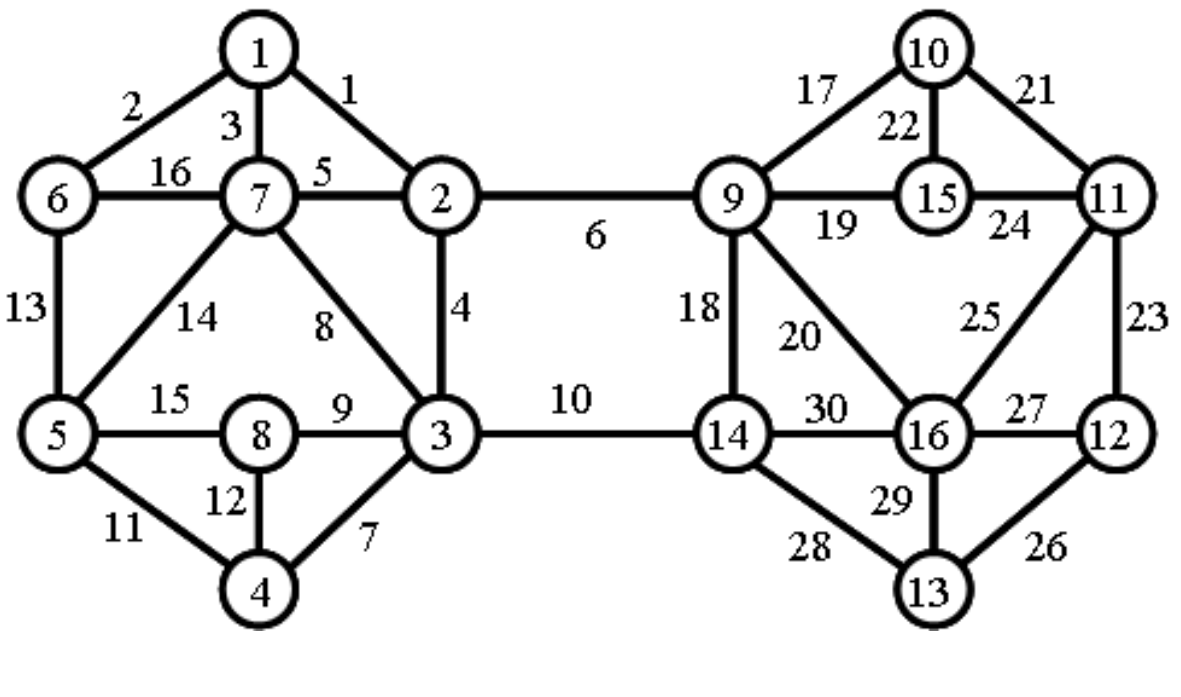

Рис. 8.11. Граф $G_{15}$.

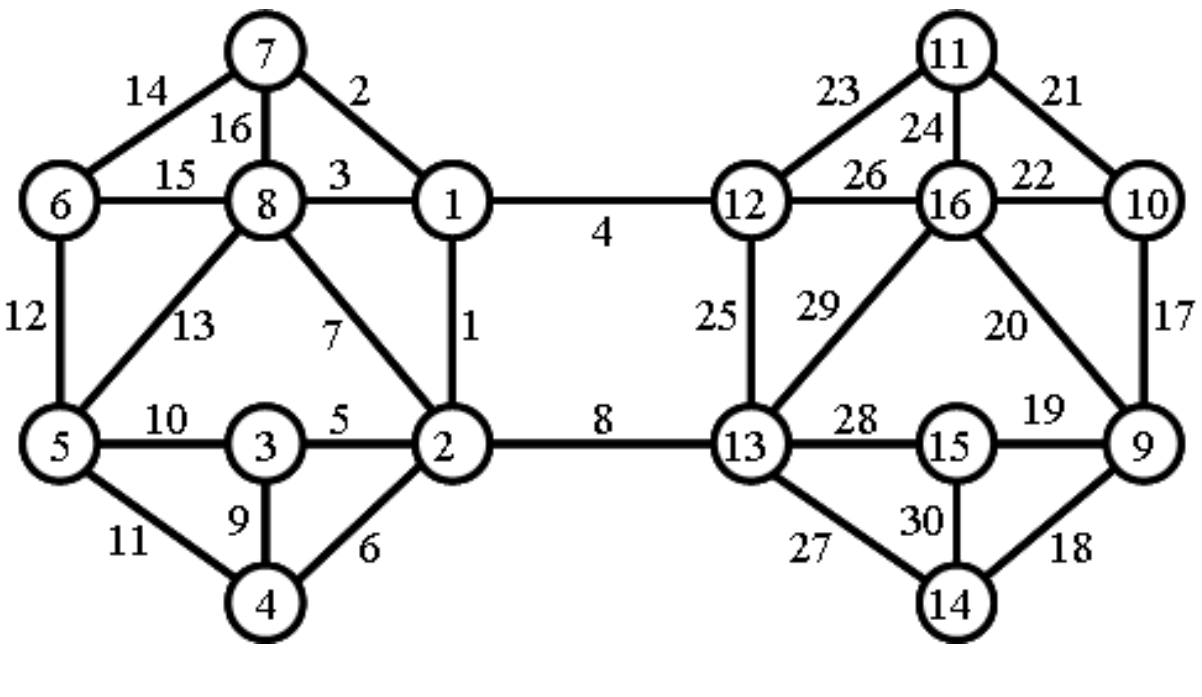

Рис. 8.12. Граф $G_{16}$.

Согласно алгебраической теории графов, пространства изоморфны, если одинаковы размеры этих пространств [4,34]. И тогда имеется два подпространства пространства суграфов: подпространство циклов C(G) и подпространство разрезов S(G). Рассмотрим случай, когда функции весов ребер для инварианта реберных циклов двух графов равны. Будут ли равны функции весов для инвариантов спектра реберных разрезов? Проверим чувствительность (способность реагировать на изменения) алгоритма распознавания изоморфизма, основанного на свойствах реберных разрезов и реберных циклов. С этой целью поменяем местами ребра $e_{20}$ и $e_{25}$ в графе $G_9$ на ребра $e_{20}$ и $e_{29}$ в графе $G_{10}$.

Изометрические циклы графа $G_{15}$

$c_1 = \{e_1,e_3,e_5\}$; $c_2 = \{e_2,e_3,e_{16}\}$; $c_3 = \{e_4,e_5,e_8\}$; $c_4 = \{e_4,e_6,e_{10},e_{18}\}$; $c_5 = \{e_7,e_8,e_{11},e_{14}\}$;
$c_6 = \{e_7,e_9,e_{12}\}$; $c_7 = \{e_8,e_9,e_{14},e_{15}\}$; $c_8 = \{e_{11},e_{12},e_{15}\}$; $c_9 = \{e_{13},e_{14},e_{16}\}$; $c_{10} = \{e_{17},e_{19},e_{22}\}$;
$c_{11} = \{e_{17},e_{20},e_{21},e_{25}\}$; $c_{12} = \{e_{18},e_{20},e_{30}\}$; $c_{13} = \{e_{19},e_{20},e_{24},e_{25}\}$; $c_{14} = \{e_{21},e_{22},e_{24}\}$;
$c_{15} = \{e_{23},e_{25},e_{27}\}$; $c_{16} = \{e_{26},e_{27},e_{29}\}$; $c_{17} = \{e_{28},e_{29},e_{30}\}$.

Количество циклов = $(12\times 3, 5\times 4)$.

Центральные разрезы графа $G_{15}$:

$s_1 = \{e_1,e_2,e_3\}$; $s_2 = \{e_1,e_4,e_5,e_6\}$; $s_3 = \{e_4,e_7,e_8,e_9,e_{10}\}$; $s_4 = \{e_7,e_{11},e_{12}\}$; $s_5 = \{e_{11},e_{13},e_{14},e_{15}\}$;
$s_6 = \{e_2,e_{13},e_{16}\}$; $s_7 = \{e_3,e_5,e_8,e_{14},e_{16}\}$; $s_8 = \{e_9,e_{12},e_{15}\}$; $s_9 = \{e_6,e_{17},e_{18},e_{19},e_{20}\}$;
$s_{10} = \{e_{17},e_{21},e_{22}\}$; $s_{11} = \{e_{21},e_{23},e_{24},e_{25}\}$; $s_{12} = \{e_{23},e_{26},e_{27}\}$; $s_{13} = \{e_{26},e_{28},e_{29}\}$;
$s_{14} = \{e_{10},e_{18},e_{28},e_{30}\}$; $s_{15} = \{e_{19},e_{22},e_{24}\}$; $s_{16} = \{e_{20},e_{25},e_{27},e_{29},e_{30}\}$.

Базовые реберные разрезы $G_{15}$:

$w_0(e_1) = \{e_2,e_3,e_4,e_5,e_6\}$;
$w_0(e_2) = \{e_1,e_3,e_{13},e_{16}\}$;
$w_0(e_3) = \{e_1,e_2,e_5,e_8,e_{14},e_{16}\}$;
$w_0(e_4) = \{e_1,e_5,e_6,e_7,e_8,e_9,e_{10}\}$;
$w_0(e_5) = \{e_1,e_3,e_4,e_6,e_8,e_{14},e_{16}\}$;
$w_0(e_6) = \{e_1,e_4,e_5,e_{17},e_{18},e_{19},e_{20}\}$;
$w_0 (e_7) = \{e_4,e_8,e_9,e_{10},e_{11},e_{12}\}$;
$w_0(e_8) = \{e_3,e_4,e_5,e_7,e_9,e_{10},e_{14},e_{16}\}$;
$w_0(e_9) = \{e_4,e_7,e_8,e_{10},e_{12},e_{15}\}$;
$w_0(e_{10}) = \{e_4,e_7,e_8,e_9,e_{18},e_{28},e_{30}\}$;
$w_0(e_{11}) = \{e_7,e_{12},e_{13},e_{14},e_{15}\}$;
$w_0(e_{12}) = \{e_7,e_9,e_{11},e_{15}\}$;
$w_0(e_{13}) = \{e_2,e_{11},e_{14},e_{15},e_{16}\}$;
$w_0(e_{14}) = \{e_3,e_5,e_8,e_{11},e_{13},e_{15},e_{16}\}$;
$w_0(e_{15}) = \{e_9,e_{11},e_{12},e_{13},e_{14}\}$;
$w_0(e_{16}) = \{e_2,e_3,e_5,e_8,e_{13},e_{14}\}$;
$w_0(e_{17}) = \{e_6,e_{18},e_{19},e_{20},e_{21},e_{22}\}$;
$w_0(e_{18}) = \{e_6,e_{10},e_{17},e_{19},e_{20},e_{28},e_{30}\}$;
$w_0(e_{19}) = \{e_6,e_{17},e_{18},e_{20},e_{22},e_{24}\}$;
$w_0(e_{20}) = \{e_6,e_{17},e_{18},e_{19},e_{25},e_{27},e_{29},e_{30}\}$;
$w_0(e_{21}) = \{e_{17},e_{22},e_{23},e_{24},e_{25}\}$;
$w_0(e_{22}) = \{e_{17},e_{19},e_{21},e_{24}\}$;
$w_0(e_{23}) = \{e_{21},e_{24},e_{25},e_{26},e_{27}\}$;
$w_0(e_{24}) = \{e_{19},e_{21},e_{22},e_{23},e_{25}\}$;
$w_0(e_{25}) = \{e_{20},e_{21},e_{23},e_{24},e_{27},e_{29},e_{30}\}$;
$w_0(e_{26}) = \{e_{23},e_{27},e_{28},e_{29}\}$;
$w_0(e_{27}) = \{e_{20},e_{23},e_{25},e_{26},e_{29},e_{30}\}$;
$w_0(e_{28}) = \{e_{10},e_{18},e_{26},e_{29},e_{30}\}$;
$w_0(e_{29}) = \{e_{20},e_{25},e_{26},e_{27},e_{28},e_{30}\}$;

$w_0(e_{30}) = \{e_{10},e_{18},e_{20},e_{25},e_{27},e_{28},e_{29}\}$.

Реберные разрезы 2-го уровня:

$w_1(e_1) = \{e_1,e_2,e_5,e_7,e_8,e_9,e_{10},e_{13},e_{16},e_{17},e_{18},e_{19},e_{20}\}$;
$w_1(e_2) = \{e_1,e_4,e_5,e_6,e_{11},e_{13},e_{14},e_{15}\}$;
$w_1(e_3) = \{e_4,e_7,e_8,e_9,e_{10},e_{11},e_{13},e_{14},e_{15}\}$;
$w_1(e_4) = \{e_2,e_3,e_4,e_5,e_7,e_9,e_{10},e_{11},e_{15},e_{17},e_{19},e_{20},e_{28},e_{30}\}$;
$w_1(e_5) = \{e_1,e_2,e_4,e_5,e_{11},e_{14},e_{15},e_{16},e_{17},e_{18},e_{19},e_{20}\}$;
$w_1(e_6) = \{e_2,e_6,e_7,e_9,e_{14},e_{16},e_{17},e_{18},e_{19},e_{20},e_{21},e_{24},e_{25},e_{27},e_{28},e_{29}\}$;
$w_1(e_7) = \{e_1,e_3,e_4,e_6,e_8,e_{10},e_{11},e_{13},e_{15},e_{16},e_{18},e_{28},e_{30}\}$;
$w_1(e_8) = \{e_1,e_3,e_7,e_9,e_{10},e_{14},e_{16},e_{18},e_{28},e_{30}\}$;
$w_1(e_9) = \{e_1,e_3,e_4,e_6,e_8,e_{10},e_{11},e_{13},e_{15},e_{16},e_{18},e_{28},e_{30}\}$;
$w_1(e_{10}) = \{e_1,e_3,e_4,e_7,e_8,e_9,e_{10},e_{11},e_{14},e_{15},e_{16},e_{17},e_{19},e_{25},e_{26},e_{27}\}$;
$w_1(e_{11}) = \{e_2,e_3,e_4,e_5,e_7,e_9,e_{10},e_{11},e_{15}\}$;
$w_1(e_{12}) = \varnothing$;
$w_1(e_{13}) = \{e_1,e_2,e_3,e_7,e_9,e_{13},e_{14}\}$;
$w_1(e_{14}) = \{e_2,e_3,e_5,e_6,e_8,e_{10},e_{13},e_{14}\}$;
$w_1(e_{15}) = \{e_2,e_3,e_4,e_5,e_7,e_9,e_{10},e_{11},e_{15}\}$;
$w_1(e_{16}) = \{e_1,e_5,e_6,e_7,e_8,e_9,e_{10}\}$;
$w_1(e_{17}) = \{e_1,e_4,e_5,e_6,e_{10},e_{18},e_{20},e_{21},e_{23},e_{24},e_{27},e_{28},e_{29}\}$;
$w_1(e_{18}) = \{e_1,e_5,e_6,e_7,e_8,e_9,e_{17},e_{18},e_{19},e_{21},e_{24},e_{26},e_{29},e_{30}\}$;
$w_1(e_{19}) = \{e_1,e_4,e_5,e_6,e_{10},e_{18},e_{20},e_{21},e_{23},e_{24},e_{27},e_{28},e_{29}\}$;
$w_1(e_{20}) = \{e_1,e_4,e_5,e_6,e_{17},e_{19},e_{25},e_{27},e_{28},e_{29}\}$;
$w_1(e_{21}) = \{e_6,e_{17},e_{18},e_{19},e_{21},e_{24},e_{26},e_{29},e_{30}\}$;
$w_1(e_{22}) = \varnothing$;
$w_1(e_{23}) = \{e_{17},e_{19},e_{23},e_{25},e_{26},e_{28},e_{29}\}$;
$w_1(e_{24}) = \{e_6,e_{17},e_{18},e_{19},e_{21},e_{24},e_{26},e_{29},e_{30}\}$;
$w_1(e_{25}) = \{e_6,e_{10},e_{20},e_{23},e_{25},e_{26},e_{29},e_{30}\}$;
$w_1(e_{26}) = \{e_{10},e_{18},e_{21},e_{23},e_{24},e_{25},e_{28},e_{30}\}$;
$w_1(e_{27}) = \{e_6,e_{10},e_{17},e_{19},e_{20},e_{28},e_{30}\}$;
$w_1(e_{28}) = \{e_4,e_6,e_7,e_8,e_9,e_{17},e_{19},e_{20},e_{23},e_{26},e_{27},e_{28},e_{30}\}$;
$w_1(e_{29}) = \{e_6,e_{17},e_{18},e_{19},e_{20},e_{21},e_{23},e_{24},e_{25}\}$;
$w_1(e_{30}) = \{e_4,e_7,e_8,e_9,e_{18},e_{21},e_{24},e_{25},e_{26},e_{27},e_{28},e_{30}\}$.

Реберные разрезы 3-го уровня:

$w_2(e_1) = \{e_2,e_3,e_5,e_7,e_8,e_9,e_{17},e_{19},e_{20},e_{21},e_{24},e_{25},e_{27},e_{29},e_{30}\}$;
$w_2(e_2) = \{e_1,e_3,e_4,e_6,e_8,e_{10},e_{11},e_{13},e_{15},e_{16},e_{17},e_{18},e_{19},e_{20}\}$;
$w_2(e_3) = \{e_1,e_2,e_5,e_6,e_7,e_8,e_9,e_{13},e_{16},e_{18},e_{28},e_{30}\}$;
$w_2(e_4) = \{e_2,e_6,e_{10},e_{13},e_{16},e_{20},e_{21},e_{24},e_{26},e_{29},e_{30}\}$;
$w_2(e_5) = \{e_1,e_3,e_6,e_{13},e_{16},e_{17},e_{18},e_{19},e_{20},e_{21},e_{24},e_{25},e_{27},e_{28},e_{29}\}$;
$w_2(e_6) = \{e_2,e_3,e_4,e_5,e_7,e_9,e_{13},e_{14},e_{17},e_{18},e_{19},e_{20},e_{21},e_{24},e_{25},e_{26},e_{27}\}$;
$w_2(e_7) = \{e_1,e_3,e_6,e_{10},e_{13},e_{16},e_{20},e_{25},e_{26},e_{27},e_{28},e_{30}\}$;
$w_2(e_8) = \{e_1,e_2,e_3,e_{10},e_{17},e_{18},e_{19},e_{25},e_{26},e_{27},e_{28},e_{30}\}$;
$w_2(e_9) = \{e_1,e_3,e_6,e_{10},e_{13},e_{16},e_{20},e_{25},e_{26},e_{27},e_{28},e_{30}\}$;
$w_2(e_{10}) = \{e_2,e_4,e_7,e_8,e_9,e_{11},e_{14},e_{15},e_{16},e_{17},e_{18},e_{19},e_{23},e_{25},e_{26},e_{29},e_{30}\}$;
$w_2(e_{11}) = \{e_2,e_{10},e_{13},e_{16},e_{18},e_{28},e_{30}\}$;
$w_2(e_{12}) = \varnothing$;
$w_2(e_{13}) = \{e_2,e_3,e_4,e_5,e_6,e_7,e_9,e_{11},e_{15}\}$;
$w_2(e_{14}) = \{e_6,e_{10},e_{17},e_{19},e_{20},e_{28},e_{30}\}$;
$w_2(e_{15}) = \{e_2,e_{10},e_{13},e_{16},e_{18},e_{28},e_{30}\}$;
$w_2(e_{16}) = \{e_2,e_3,e_4,e_5,e_7,e_9,e_{10},e_{11},e_{15},e_{17},e_{19},e_{20},e_{28},e_{30}\}$;
$w_2(e_{17}) = \{e_1,e_2,e_5,e_6,e_8,e_{10},e_{14},e_{16},e_{23},e_{27},e_{28},e_{29}\}$;

$w_2(e_{18}) = \{e_2,e_3,e_5,e_6,e_8,e_{10},e_{11},e_{15},e_{23},e_{26},e_{27}\}$;
$w_2(e_{19}) = \{e_1,e_2,e_5,e_6,e_8,e_{10},e_{14},e_{16},e_{23},e_{27},e_{28},e_{29}\}$;
$w_2(e_{20}) = \{e_1,e_2,e_4,e_5,e_6,e_7,e_9,e_{14},e_{16},e_{26},e_{28},e_{29}\}$;
$w_2(e_{21}) = \{e_1,e_4,e_5,e_6,e_{23},e_{26},e_{27}\}$;
$w_2(e_{22}) = \varnothing$;
$w_2(e_{23}) = \{e_{10},e_{17},e_{18},e_{19},e_{21},e_{24},e_{26},e_{29},e_{30}\}$;
$w_2(e_{24}) = \{e_1,e_4,e_5,e_6,e_{23},e_{26},e_{27}\}$;
$w_2(e_{25}) = \{e_1,e_5,e_6,e_7,e_8,e_9,e_{10}\}$;
$w_2(e_{26}) = \{e_4,e_6,e_7,e_8,e_9,e_{10},e_{18},e_{20},e_{21},e_{23},e_{24},e_{27},e_{28},e_{29}\}$;
$w_2(e_{27}) = \{e_1,e_5,e_6,e_7,e_8,e_9,e_{17},e_{18},e_{19},e_{21},e_{24},e_{26},e_{29},e_{30}\}$;
$w_2(e_{28}) = \{e_3,e_5,e_7,e_8,e_9,e_{11},e_{14},e_{15},e_{16},e_{17},e_{19},e_{20},e_{26},e_{29},e_{30}\}$;
$w_2(e_{29}) = \{e_1,e_4,e_5,e_{10},e_{17},e_{19},e_{20},e_{23},e_{26},e_{27},e_{28},e_{30}\}$;
$w_2(e_{30}) = \{e_1,e_3,e_4,e_7,e_8,e_9,e_{10},e_{11},e_{14},e_{15},e_{16},e_{23},e_{27},e_{28},e_{29}\}$.

Реберные разрезы 4-го уровня:

$w_3(e_1) = \{e_1,e_2,e_3,e_{11},e_{13},e_{14},e_{15},e_{17},e_{19},e_{21},e_{24}\}$;
$w_3(e_2) = \{e_1,e_3,e_{10},e_{13},e_{16},e_{18},e_{21},e_{24},e_{25},e_{27},e_{29},e_{30}\}$;
$w_3(e_3) = \{e_1,e_2,e_3,e_6,e_{18},e_{20},e_{25},e_{26},e_{27}\}$;
$w_3(e_4) = \{e_6,e_7,e_9,e_{10},e_{11},e_{15},e_{17},e_{19},e_{20},e_{21},e_{23},e_{24},e_{25},e_{26},e_{29},e_{30}\}$;
$w_3(e_5) = \{e_6,e_{11},e_{13},e_{14},e_{15},e_{17},e_{18},e_{19},e_{20},e_{21},e_{24},e_{25},e_{26},e_{27}\}$;
$w_3(e_6) = \{e_3,e_4,e_5,e_{11},e_{14},e_{15},e_{16},e_{18},e_{20},e_{21},e_{23},e_{24},e_{26},e_{27}\}$;
$w_3(e_7) = \{e_4,e_7,e_8,e_9,e_{11},e_{13},e_{14},e_{15},e_{18},e_{21},e_{23},e_{24},e_{25},e_{28},e_{30}\}$;
$w_3(e_8) = \{e_7,e_9,e_{10},e_{13},e_{14},e_{18},e_{23},e_{27},e_{29},e_{30}\}$;
$w_3(e_9) = \{e_4,e_7,e_8,e_9,e_{11},e_{13},e_{14},e_{15},e_{18},e_{21},e_{23},e_{24},e_{25},e_{28},e_{30}\}$;
$w_3(e_{10}) = \{e_2,e_4,e_8,e_{11},e_{13},e_{15},e_{16},e_{18},e_{21},e_{24},e_{25},e_{27},e_{29},e_{30}\}$;
$w_3(e_{11}) = \{e_1,e_4,e_5,e_6,e_7,e_9,e_{10},e_{11},e_{15},e_{17},e_{18},e_{19},e_{25},e_{26},e_{27},e_{28},e_{30}\}$;
$w_3(e_{12}) = \varnothing$;
$w_3(e_{13}) = \{e_1,e_2,e_5,e_7,e_8,e_9,e_{10},e_{13},e_{16},e_{17},e_{18},e_{19},e_{20}\}$;
$w_3(e_{14}) = \{e_1,e_5,e_6,e_7,e_8,e_9,e_{17},e_{18},e_{19},e_{21},e_{24},e_{26},e_{29},e_{30}\}$;
$w_3(e_{15}) = \{e_1,e_4,e_5,e_6,e_7,e_9,e_{10},e_{11},e_{15},e_{17},e_{18},e_{19},e_{25},e_{26},e_{27},e_{28},e_{30}\}$;
$w_3(e_{16}) = \{e_2,e_6,e_{10},e_{13},e_{16},e_{20},e_{21},e_{24},e_{26},e_{29},e_{30}\}$;
$w_3(e_{17}) = \{e_1,e_4,e_5,e_{11},e_{13},e_{14},e_{15},e_{17},e_{18},e_{19},e_{20},e_{21},e_{23},e_{24},e_{25}\}$;
$w_4(e_{18}) = \{e_2,e_3,e_5,e_6,e_7,e_8,e_9,e_{10},e_{11},e_{13},e_{14},e_{15},e_{17},e_{19},e_{21},e_{24}\}$;
$w_4(e_{19}) = \{e_1,e_4,e_5,e_{11},e_{13},e_{14},e_{15},e_{17},e_{18},e_{19},e_{20},e_{21},e_{23},e_{24},e_{25}\}$;
$w_4(e_{20}) = \{e_3,e_4,e_5,e_6,e_{13},e_{16},e_{17},e_{19},e_{23},e_{25}\}$;
$w_4(e_{21}) = \{e_1,e_2,e_4,e_5,e_6,e_7,e_9,e_{10},e_{14},e_{16},e_{17},e_{18},e_{19},e_{21},e_{24},e_{28},e_{30}\}$;
$w_4(e_{22}) = \varnothing$;
$w_4(e_{23}) = \{e_4,e_6,e_7,e_8,e_9,e_{17},e_{19},e_{20},e_{23},e_{26},e_{27},e_{28},e_{30}\}$;
$w_4(e_{24}) = \{e_1,e_2,e_4,e_5,e_6,e_7,e_9,e_{10},e_{14},e_{16},e_{17},e_{18},e_{19},e_{21},e_{24},e_{28},e_{30}\}$;
$w_4(e_{25}) = \{e_2,e_3,e_4,e_5,e_7,e_9,e_{10},e_{11},e_{15},e_{17},e_{19},e_{20},e_{28},e_{30}\}$;
$w_4(e_{26}) = \{e_3,e_4,e_5,e_6,e_{11},e_{14},e_{15},e_{16},e_{23},e_{27},e_{28},e_{29}\}$;
$w_4(e_{27}) = \{e_2,e_3,e_5,e_6,e_8,e_{10},e_{11},e_{15},e_{23},e_{26},e_{27}\}$;
$w_4(e_{28}) = \{e_7,e_9,e_{11},e_{15},e_{21},e_{23},e_{24},e_{25},e_{26},e_{28},e_{29}\}$;
$w_4(e_{29}) = \{e_2,e_4,e_8,e_{10},e_{14},e_{16},e_{26},e_{28},e_{29}\}$;
$w_4(e_{30}) = \{e_2,e_4,e_7,e_8,e_9,e_{10},e_{11},e_{14},e_{15},e_{16},e_{21},e_{23},e_{24},e_{25}\}$.

Реберные разрезы 5-го уровня:

$w_4(e_1) = \{e_2,e_3,e_4,e_5,e_6,e_7,e_9,e_{11},e_{15}\}$;
$w_4(e_2) = \{e_1,e_5,e_7,e_8,e_9,e_{11},e_{13},e_{14},e_{15},e_{20},e_{25},e_{27},e_{29},e_{30}\}$;
$w_4(e_3) = \{e_1,e_5,e_6,e_8,e_{10},e_{13},e_{14},e_{17},e_{19},e_{21},e_{23},e_{24},e_{26},e_{27}\}$;
$w_4(e_4) = \varnothing$;
$w_4(e_5) = \{e_1,e_2,e_3,e_4,e_7,e_8,e_9,e_{10},e_{11},e_{13},e_{14},e_{15},e_{17},e_{19},e_{21},e_{23},e_{24},e_{26},e_{27}\}$;

$w_4(e_6) = \{e_1,e_3,e_4,e_8,e_{14},e_{16},e_{17},e_{18},e_{19},e_{25},e_{27},e_{29},e_{30}\}$;
$w_4(e_7) = \{e_1,e_2,e_4,e_5,e_{10},e_{13},e_{16},e_{20},e_{21},e_{23},e_{24},e_{27},e_{29},e_{30}\}$;
$w_4(e_8) = \{e_2,e_3,e_4,e_5,e_6,e_{11},e_{13},e_{14},e_{15},e_{17},e_{19},e_{21},e_{23},e_{24},e_{26},e_{27}\}$;
$w_4(e_9) = \{e_1,e_2,e_4,e_5,e_{10},e_{13},e_{16},e_{20},e_{21},e_{23},e_{24},e_{27},e_{29},e_{30}\}$;
$w_4(e_{10}) = \{e_3,e_4,e_5,e_7,e_9,e_{14},e_{16},e_{18},e_{20},e_{25},e_{27},e_{28},e_{29}\}$;
$w_4(e_{11}) = \{e_1,e_2,e_5,e_8,e_{14},e_{16},e_{17},e_{19},e_{20},e_{23},e_{27},e_{29},e_{30}\}$; $w_5(e_{12}) = \varnothing$;
$w_4(e_{13}) = \{e_2,e_3,e_5,e_7,e_8,e_9,e_{17},e_{19},e_{20},e_{21},e_{24},e_{25},e_{27},e_{29},e_{30}\}$;
$w_4(e_{14}) = \{e_2,e_3,e_5,e_6,e_8,e_{10},e_{11},e_{15},e_{23},e_{26},e_{27}\}$;
$w_4(e_{15}) = \{e_1,e_2,e_5,e_8,e_{14},e_{16},e_{17},e_{19},e_{20},e_{23},e_{27},e_{29},e_{30}\}$;
$w_4(e_{16}) = \{e_6,e_7,e_9,e_{10},e_{11},e_{15},e_{17},e_{19},e_{20},e_{21},e_{23},e_{24},e_{25},e_{26},e_{29},e_{30}\}$;
$w_4(e_{17}) = \{e_3,e_5,e_6,e_8,e_{11},e_{13},e_{15},e_{16},e_{18},e_{23},e_{26},e_{27},e_{28},e_{30}\}$;
$w_4(e_{18}) = \varnothing$;
$w_4(e_{19}) = \{e_3,e_5,e_6,e_8,e_{11},e_{13},e_{15},e_{16},e_{18},e_{23},e_{26},e_{27},e_{28},e_{30}\}$;
$w_4(e_{20}) = \{e_2,e_7,e_9,e_{10},e_{11},e_{13},e_{15},e_{16},e_{18},e_{21},e_{23},e_{24},e_{25},e_{26},e_{29},e_{30}\}$;
$w_4(e_{21}) = \{e_3,e_5,e_7,e_8,e_9,e_{13},e_{16},e_{20},e_{25},e_{26},e_{27},e_{28},e_{30}\}$;
$w_4(e_{22}) = \varnothing$;
$w_4(e_{23}) = \{e_3,e_5,e_7,e_8,e_9,e_{11},e_{14},e_{15},e_{16},e_{17},e_{19},e_{20},e_{26},e_{29},e_{30}\}$;
$w_4(e_{24}) = \{e_3,e_5,e_7,e_8,e_9,e_{13},e_{16},e_{20},e_{25},e_{26},e_{27},e_{28},e_{30}\}$;
$w_4(e_{25}) = \{e_2,e_6,e_{10},e_{13},e_{16},e_{20},e_{21},e_{24},e_{26},e_{29},e_{30}\}$;
$w_4(e_{26}) = \{e_3,e_5,e_8,e_{14},e_{16},e_{17},e_{19},e_{20},e_{21},e_{23},e_{24},e_{25},e_{28},e_{30}\}$;
$w_4(e_{27}) = \{e_2,e_3,e_5,e_6,e_7,e_8,e_9,e_{10},e_{11},e_{13},e_{14},e_{15},e_{17},e_{19},e_{21},e_{24}\}$;
$w_4(e_{28}) = \{e_{10},e_{17},e_{18},e_{19},e_{21},e_{24},e_{26},e_{29},e_{30}\}$;
$w_4(e_{29}) = \{e_2,e_6,e_7,e_9,e_{10},e_{11},e_{13},e_{15},e_{16},e_{20},e_{23},e_{25},e_{28},e_{30}\}$;
$w_4(e_{30}) = \{e_2,e_6,e_7,e_9,e_{11},e_{13},e_{15},e_{16},e_{17},e_{18},e_{19},e_{20},e_{21},e_{23},e_{24},e_{25},e_{26},e_{28},e_{29}\}$.

Реберные разрезы 6-го уровня:

$w_5(e_1) = \varnothing$;
$w_5(e_2) = \{e_1,e_4,e_5,e_6,e_7,e_9,e_{13},e_{14},e_{17},e_{19},e_{21},e_{24}\}$;
$w_5(e_3) = \{e_4,e_8,e_{10},e_{13},e_{14},e_{17},e_{19},e_{21},e_{24}\}$;
$w_5(e_4) = \varnothing$;
$w_5(e_5) = \{e_1,e_2,e_4,e_5,e_7,e_9,e_{14},e_{16},e_{18},e_{20},e_{21},e_{24}\}$;
$w_5(e_6) = \{e_2,e_6,e_{11},e_{14},e_{15},e_{16},e_{18},e_{20},e_{25},e_{27},e_{28},e_{29}\}$;
$w_5(e_7) = \{e_1,e_2,e_4,e_5,e_{11},e_{14},e_{15},e_{16},e_{18},e_{20},e_{23},e_{25},e_{26},e_{28},e_{29}\}$;
$w_5(e_8) = \{e_1,e_3,e_{10},e_{11},e_{14},e_{15},e_{16},e_{17},e_{18},e_{19},e_{21},e_{24},e_{28},e_{30}\}$;
$w_5(e_9) = \{e_1,e_2,e_4,e_5,e_{11},e_{14},e_{15},e_{16},e_{18},e_{20},e_{23},e_{25},e_{26},e_{28},e_{29}\}$;
$w_5(e_{10}) = \{e_1,e_3,e_4,e_8,e_{10},e_{14},e_{16},e_{21},e_{24},e_{25},e_{26},e_{27}\}$;
$w_5(e_{11}) = \{e_4,e_6,e_7,e_8,e_9,e_{11},e_{13},e_{14},e_{15},e_{20},e_{23},e_{25},e_{26},e_{29},e_{30}\}$;
$w_5(e_{12}) = \varnothing$;
$w_5(e_{13}) = \{e_1,e_2,e_3,e_{11},e_{13},e_{14},e_{15},e_{17},e_{19},e_{21},e_{24}\}$;
$w_5(e_{14}) = \{e_2,e_3,e_5,e_6,e_7,e_8,e_9,e_{10},e_{11},e_{13},e_{14},e_{15},e_{17},e_{19},e_{21},e_{24}\}$;
$w_5(e_{15}) = \{e_4,e_6,e_7,e_8,e_9,e_{11},e_{13},e_{14},e_{15},e_{20},e_{23},e_{25},e_{26},e_{29},e_{30}\}$;
$w_5(e_{16}) = \varnothing$;
$w_5(e_{17}) = \{e_1,e_2,e_3,e_4,e_8,e_{13},e_{14},e_{18},e_{21},e_{24},e_{25},e_{26},e_{27},e_{28},e_{30}\}$;
$w_5(e_{18}) = \varnothing$;
$w_5(e_{19}) = \{e_1,e_2,e_3,e_4,e_8,e_{13},e_{14},e_{18},e_{21},e_{24},e_{25},e_{26},e_{27},e_{28},e_{30}\}$;
$w_5(e_{20}) = \{e_1,e_4,e_5,e_6,e_7,e_9,e_{11},e_{15},e_{21},e_{24},e_{25},e_{27},e_{28},e_{29}\}$;
$w_5(e_{21}) = \{e_2,e_3,e_5,e_8,e_{10},e_{13},e_{14},e_{17},e_{18},e_{19},e_{20},e_{21},e_{23},e_{24},e_{25}\}$;
$w_5(e_{22}) = \varnothing$;
$w_5(e_{23}) = \{e_7,e_9,e_{11},e_{15},e_{21},e_{23},e_{24},e_{25},e_{26},e_{28},e_{29}\}$;
$w_5(e_{24}) = \{e_2,e_3,e_5,e_8,e_{10},e_{13},e_{14},e_{17},e_{18},e_{19},e_{20},e_{21},e_{23},e_{24},e_{25}\}$;
$w_5(e_{25}) = \{e_6,e_7,e_9,e_{10},e_{11},e_{15},e_{17},e_{19},e_{20},e_{21},e_{23},e_{24},e_{25},e_{26},e_{29},e_{30}\}$;
$w_5(e_{26}) = \{e_7,e_9,e_{10},e_{11},e_{15},e_{17},e_{18},e_{19},e_{23},e_{25},e_{28},e_{30}\}$;

$w_5(e_{27}) = \varnothing$;
$w_5(e_{28}) = \varnothing$;
$w_5(e_{29}) = \{e_6,e_7,e_9,e_{11},e_{15},e_{18},e_{20},e_{23},e_{25}\}$;
$w_5(e_{30}) = \{e_4,e_8,e_{11},e_{15},e_{17},e_{18},e_{19},e_{25},e_{26},e_{27},e_{28},e_{30}\}$.

Изометрические циклы графа $G_{16}$:

$c_1 = \{e_1,e_3,e_7\}$; $c_2 = \{e_1,e_4,e_8,e_{25}\}$; $c_3 = \{e_2,e_3,e_{16}\}$; $c_4 = \{e_5,e_6,e_9\}$; $c_5 = \{e_5,e_7,e_{10},e_{13}\}$;
$c_6 = \{e_6,e_7,e_{11},e_{13}\}$; $c_7 = \{e_9,e_{10},e_{11}\}$; $c_8 = \{e_{12},e_{13},e_{15}\}$; $c_9 = \{e_{14},e_{15},e_{16}\}$;
$c_{10} = \{e_{17},e_{20},e_{22}\}$; $c_{11} = \{e_{18},e_{19},e_{30}\}$; $c_{12} = \{e_{18},e_{20},e_{27},e_{29}\}$; $c_{13} = \{e_{19},e_{20},e_{28},e_{29}\}$;
$c_{14} = \{e_{21},e_{22},e_{24}\}$; $c_{15} = \{e_{23},e_{24},e_{26}\}$; $c_{16} = \{e_{25},e_{26},e_{29}\}$; $c_{17} = \{e_{27},e_{28},e_{30}\}$.

Количество циклов = $(12\times 3, 5\times 4)$.

Центральные разрезы графа $G_{16}$:

$s_1 = \{e_1,e_2,e_3,e_4\}$; $s_2 = \{e_1,e_5,e_6,e_7,e_8\}$; $s_3 = \{e_5,e_9,e_{10}\}$; $s_4 = \{e_6,e_9,e_{11}\}$; $s_5 = \{e_{10},e_{11},e_{12},e_{13}\}$;
$s_6 = \{e_{12},e_{14},e_{15}\}$; $s_7 = \{e_2,e_{14},e_{16}\}$; $s_8 = \{e_3,e_7,e_{13},e_{15},e_{16}\}$; $s_9 = \{e_{17},e_{18},e_{19},e_{20}\}$;
$s_{10} = \{e_{17},e_{21},e_{22}\}$; $s_{11} = \{e_{21},e_{23},e_{24}\}$; $s_{12} = \{e_4,e_{23},e_{25},e_{26}\}$; $s_{13} = \{e_8,e_{25},e_{27},e_{28},e_{29}\}$;
$s_{14} = \{e_{18},e_{27},e_{30}\}$; $s_{15} = \{e_{19},e_{28},e_{30}\}$; $s_{16} = \{e_{20},e_{22},e_{24},e_{26},e_{29}\}$.

Базовые реберные разрезы графа $G_{16}$:

$w_0(e_1) = \{e_2,e_3,e_4,e_5,e_6,e_7,e_8\}$;
$w_0(e_2) = \{e_1,e_3,e_4,e_{14},e_{16}\}$;
$w_0(e_3) = \{e_1,e_2,e_4,e_7,e_{13},e_{15},e_{16}\}$;
$w_0(e_4) = \{e_1,e_2,e_3,e_{23},e_{25},e_{26}\}$;
$w_0(e_5) = \{e_1,e_6,e_7,e_8,e_9,e_{10}\}$;
$w_0(e_6) = \{e_1,e_5,e_7,e_8,e_9,e_{11}\}$;
$w_0(e_7) = \{e_1,e_3,e_5,e_6,e_8,e_{13},e_{15},e_{16}\}$;
$w_0(e_8) = \{e_1,e_5,e_6,e_7,e_{25},e_{27},e_{28},e_{29}\}$;
$w_0(e_9) = \{e_5,e_6,e_{10},e_{11}\}$;
$w_0(e_{10}) = \{e_5,e_9,e_{11},e_{12},e_{13}\}$;
$w_0(e_{11}) = \{e_6,e_9,e_{10},e_{12},e_{13}\}$;
$w_0(e_{12}) = \{e_{10},e_{11},e_{13},e_{14},e_{15}\}$;
$w_0(e_{13}) = \{e_3,e_7,e_{10},e_{11},e_{12},e_{15},e_{16}\}$;
$w_0(e_{14}) = \{e_2,e_{12},e_{15},e_{16}\}$;
$w_0(e_{15}) = \{e_3,e_7,e_{12},e_{13},e_{14},e_{16}\}$;
$w_0(e_{16}) = \{e_2,e_3,e_7,e_{13},e_{14},e_{15}\}$;
$w_0(e_{17}) = \{e_{18},e_{19},e_{20},e_{21},e_{22}\}$;
$w_0(e_{18}) = \{e_{17},e_{19},e_{20},e_{27},e_{30}\}$;
$w_0(e_{19}) = \{e_{17},e_{18},e_{20},e_{28},e_{30}\}$;
$w_0(e_{20}) = \{e_{17},e_{18},e_{19},e_{22},e_{24},e_{26},e_{29}\}$;
$w_0(e_{21}) = \{e_{17},e_{22},e_{23},e_{24}\}$;
$w_0(e_{22}) = \{e_{17},e_{20},e_{21},e_{24},e_{26},e_{29}\}$;
$w_0(e_{23}) = \{e_4,e_{21},e_{24},e_{25},e_{26}\}$;
$w_0(e_{24}) = \{e_{20},e_{21},e_{22},e_{23},e_{26},e_{29}\}$;
$w_0(e_{25}) = \{e_4,e_8,e_{23},e_{26},e_{27},e_{28},e_{29}\}$;
$w_0(e_{26}) = \{e_4,e_{20},e_{22},e_{23},e_{24},e_{25},e_{29}\}$;
$w_0(e_{27}) = \{e_8,e_{18},e_{25},e_{28},e_{29},e_{30}\}$;
$w_0(e_{28}) = \{e_8,e_{19},e_{25},e_{27},e_{29},e_{30}\}$;
$w_0(e_{29}) = \{e_8,e_{20},e_{22},e_{24},e_{25},e_{26},e_{27},e_{28}\}$;
$w_0(e_{30}) = \{e_{18},e_{19},e_{27},e_{28}\}$.

Реберные разрезы 2-го уровня:

$w_1(e_1) = \{e_1,e_3,e_5,e_6,e_8,e_{10},e_{11},e_{14},e_{16},e_{23},e_{26},e_{27},e_{28},e_{29}\}$;

$w_1(e_2) = \{e_2,e_3,e_5,e_6,e_7,e_8,e_{12},e_{14},e_{15},e_{23},e_{25},e_{26}\}$;
$w_1(e_3) = \{e_1,e_2,e_3,e_{10},e_{11},e_{13},e_{14},e_{15},e_{23},e_{25},e_{26}\}$;
$w_1(e_4) = \{e_5,e_6,e_{13},e_{14},e_{15},e_{20},e_{21},e_{22},e_{27},e_{28}\}$;
$w_1(e_5) = \{e_1,e_2,e_4,e_7,e_8,e_{10},e_{11},e_{12},e_{15},e_{16},e_{25},e_{27},e_{28},e_{29}\}$;
$w_1(e_6) = \{e_1,e_2,e_4,e_7,e_8,e_{10},e_{11},e_{12},e_{15},e_{16},e_{25},e_{27},e_{28},e_{29}\}$;
$w_1(e_7) = \{e_2,e_5,e_6,e_8,e_{13},e_{15},e_{16},e_{25},e_{27},e_{28},e_{29}\}$;
$w_1(e_8) = \{e_1,e_2,e_5,e_6,e_7,e_{10},e_{11},e_{13},e_{15},e_{16},e_{18},e_{19},e_{20},e_{22},e_{23},e_{24},e_{25},e_{27},e_{28},e_{29}\}$;
$w_1(e_9) = \varnothing$;
$w_1(e_{10}) = \{e_1,e_3,e_5,e_6,e_8,e_{10},e_{11},e_{14},e_{16}\}$;
$w_1(e_{11}) = \{e_1,e_3,e_5,e_6,e_8,e_{10},e_{11},e_{14},e_{16}\}$;
$w_1(e_{12}) = \{e_2,e_5,e_6,e_{12},e_{13},e_{14},e_{16}\}$;
$w_1(e_{13}) = \{e_3,e_4,e_7,e_8,e_{12},e_{13},e_{14},e_{16}\}$;
$w_1(e_{14}) = \{e_1,e_2,e_3,e_4,e_{10},e_{11},e_{12},e_{13}\}$;
$w_1(e_{15}) = \{e_2,e_3,e_4,e_5,e_6,e_7,e_8\}$;
$w_1(e_{16}) = \{e_1,e_5,e_6,e_7,e_8,e_{10},e_{11},e_{12},e_{13}\}$;
$w_1(e_{17}) = \{e_{17},e_{20},e_{21},e_{23},e_{24},e_{27},e_{28}\}$;
$w_1(e_{18}) = \{e_8,e_{18},e_{19},e_{21},e_{24},e_{25},e_{26},e_{27},e_{28}\}$;
$w_1(e_{19}) = \{e_8,e_{18},e_{19},e_{21},e_{24},e_{25},e_{26},e_{27},e_{28}\}$;
$w_1(e_{20}) = \{e_4,e_8,e_{17},e_{20},e_{21},e_{24},e_{26},e_{29}\}$;
$w_1(e_{21}) = \{e_4,e_{17},e_{18},e_{19},e_{20},e_{23},e_{25},e_{26}\}$;
$w_1(e_{22}) = \{e_4,e_8,e_{23},e_{26},e_{27},e_{28},e_{29}\}$;
$w_1(e_{23}) = \{e_1,e_2,e_3,e_8,e_{17},e_{21},e_{22},e_{23},e_{26},e_{27},e_{28},e_{29}\}$;
$w_1(e_{24}) = \{e_8,e_{17},e_{18},e_{19},e_{20},e_{25},e_{27},e_{28},e_{29}\}$;
$w_1(e_{25}) = \{e_2,e_3,e_5,e_6,e_7,e_8,e_{18},e_{19},e_{21},e_{24},e_{25},e_{26},e_{27},e_{28}\}$;
$w_1(e_{26}) = \{e_1,e_2,e_3,e_{18},e_{19},e_{20},e_{21},e_{22},e_{23},e_{25},e_{26}\}$;
$w_1(e_{27}) = \{e_1,e_4,e_5,e_6,e_7,e_8,e_{17},e_{18},e_{19},e_{22},e_{23},e_{24},e_{25},e_{29}\}$;
$w_1(e_{28}) = \{e_1,e_4,e_5,e_6,e_7,e_8,e_{17},e_{18},e_{19},e_{22},e_{23},e_{24},e_{25},e_{29}\}$;
$w_1(e_{29}) = \{e_1,e_5,e_6,e_7,e_8,e_{20},e_{22},e_{23},e_{24},e_{27},e_{28}\}$;
$w_1(e_{30}) = \varnothing$.

Реберные разрезы 3-го уровня:

$w_2(e_1) = \{e_{12},e_{14},e_{15},e_{18},e_{19},e_{21},e_{23},e_{24},e_{27},e_{28}\}$;
$w_2(e_2) = \{e_3,e_4,e_5,e_6,e_7,e_{14},e_{16},e_{20},e_{21},e_{22},e_{25},e_{29}\}$;
$w_2(e_3) = \{e_2,e_{12},e_{15},e_{16},e_{20},e_{21},e_{22},e_{27},e_{28}\}$;
$w_2(e_4) = \{e_2,e_5,e_6,e_{12},e_{13},e_{14},e_{16},e_{17},e_{20},e_{21},e_{23},e_{24},e_{27},e_{28}\}$;
$w_2(e_5) = \{e_2,e_4,e_{12},e_{15},e_{16},e_{18},e_{19},e_{20},e_{22},e_{24},e_{25},e_{26}\}$;
$w_2(e_6) = \{e_2,e_4,e_{12},e_{15},e_{16},e_{18},e_{19},e_{20},e_{22},e_{24},e_{25},e_{26}\}$;
$w_2(e_7) = \{e_2,e_{14},e_{16},e_{18},e_{19},e_{20},e_{22},e_{23},e_{24}\}$;
$w_2(e_8) = \{e_{10},e_{11},e_{13},e_{14},e_{15},e_{18},e_{19},e_{20},e_{21},e_{22}\}$;
$w_2(e_9) = \varnothing$;
$w_2(e_{10}) = \{e_8,e_{12},e_{14},e_{15},e_{25},e_{27},e_{28},e_{29}\}$;
$w_2(e_{11}) = \{e_8,e_{12},e_{14},e_{15},e_{25},e_{27},e_{28},e_{29}\}$;
$w_2(e_{12}) = \{e_1,e_3,e_4,e_5,e_6,e_{10},e_{11},e_{14},e_{16}\}$;
$w_2(e_{13}) = \{e_4,e_8,e_{23},e_{26},e_{27},e_{28},e_{29}\}$;
$w_2(e_{14}) = \{e_1,e_2,e_4,e_7,e_8,e_{10},e_{11},e_{12},e_{15},e_{16},e_{23},e_{25},e_{26}\}$;
$w_2(e_{15}) = \{e_1,e_3,e_5,e_6,e_8,e_{10},e_{11},e_{14},e_{16},e_{23},e_{26},e_{27},e_{28},e_{29}\}$;
$w_2(e_{16}) = \{e_2,e_3,e_4,e_5,e_6,e_7,e_{12},e_{14},e_{15},e_{25},e_{27},e_{28},e_{29}\}$;
$w_2(e_{17}) = \{e_4,e_{18},e_{19},e_{21},e_{24},e_{25},e_{26},e_{27},e_{28}\}$;
$w_2(e_{18}) = \{e_1,e_5,e_6,e_7,e_8,e_{17},e_{21},e_{22}\}$;
$w_2(e_{19}) = \{e_1,e_5,e_6,e_7,e_8,e_{17},e_{21},e_{22}\}$;
$w_2(e_{20}) = \{e_2,e_3,e_4,e_5,e_6,e_7,e_8\}$;
$w_2(e_{21}) = \{e_1,e_2,e_3,e_4,e_8,e_{17},e_{18},e_{19},e_{22},e_{23},e_{24},e_{25},e_{29}\}$;

$w_2(e_{22}) = \{e_2,e_3,e_5,e_6,e_7,e_8,e_{18},e_{19},e_{21},e_{24},e_{25},e_{26},e_{27},e_{28}\}$;
$w_2(e_{23}) = \{e_1,e_4,e_7,e_{13},e_{14},e_{15},e_{21},e_{24},e_{26},e_{27},e_{28},e_{29}\}$;
$w_2(e_{24}) = \{e_1,e_4,e_5,e_6,e_7,e_{17},e_{21},e_{22},e_{23},e_{26},e_{27},e_{28},e_{29}\}$;
$w_2(e_{25}) = \{e_2,e_5,e_6,e_{10},e_{11},e_{14},e_{16},e_{17},e_{21},e_{22}\}$;
$w_2(e_{26}) = \{e_5,e_6,e_{13},e_{14},e_{15},e_{17},e_{22},e_{23},e_{24}\}$;
$w_2(e_{27}) = \{e_1,e_3,e_4,e_{10},e_{11},e_{13},e_{15},e_{16},e_{17},e_{22},e_{23},e_{24}\}$;
$w_2(e_{28}) = \{e_1,e_3,e_4,e_{10},e_{11},e_{13},e_{15},e_{16},e_{17},e_{22},e_{23},e_{24}\}$;
$w_2(e_{29}) = \{e_2,e_{10},e_{11},e_{13},e_{15},e_{16},e_{21},e_{23},e_{24}\}$;
$w_2(e_{30}) = \varnothing$.

Реберные разрезы 4-го уровня:

$w_3(e_1) = \{e_2,e_3,e_4,e_7,e_{10},e_{11},e_{17},e_{20},e_{25},e_{29}\}$;
$w_3(e_2) = \{e_1,e_3,e_8,e_{10},e_{11},e_{12},e_{13},e_{14},e_{16},e_{17},e_{18},e_{19},e_{21},e_{22},e_{23},e_{26},e_{29}\}$;
$w_3(e_3) = \{e_1,e_2,e_3,e_4,e_{10},e_{11},e_{12},e_{13},e_{17},e_{20},e_{21},e_{23},e_{24},e_{27},e_{28}\}$;
$w_3(e_4) = \{e_1,e_3,e_5,e_6,e_{10},e_{11},e_{14},e_{16},e_{18},e_{19},e_{21},e_{24},e_{25},e_{26},e_{27},e_{28}\}$;
$w_3(e_5) = \{e_4,e_8,e_{10},e_{11},e_{12},e_{13},e_{20},e_{22},e_{24},e_{26},e_{29}\}$;
$w_3(e_6) = \{e_4,e_8,e_{10},e_{11},e_{12},e_{13},e_{20},e_{22},e_{24},e_{26},e_{29}\}$;
$w_3(e_7) = \{e_1,e_7,e_{12},e_{13},e_{21},e_{23},e_{24},e_{25},e_{27},e_{28},e_{29}\}$;
$w_3(e_8) = \{e_2,e_5,e_6,e_{12},e_{13},e_{14},e_{16},e_{17},e_{20},e_{21},e_{23},e_{24},e_{27},e_{28}\}$;
$w_3(e_9) = \varnothing$;
$w_3(e_{10}) = \{e_1,e_2,e_3,e_4,e_5,e_6,e_{10},e_{11},e_{18},e_{19},e_{20},e_{22},e_{23},e_{24}\}$;
$w_3(e_{11}) = \{e_1,e_2,e_3,e_4,e_5,e_6,e_{10},e_{11},e_{18},e_{19},e_{20},e_{22},e_{23},e_{24}\}$;
$w_3(e_{12}) = \{e_2,e_3,e_5,e_6,e_7,e_8,e_{12},e_{14},e_{15},e_{23},e_{25},e_{26}\}$;
$w_3(e_{13}) = \{e_2,e_3,e_5,e_6,e_7,e_8,e_{18},e_{19},e_{21},e_{24},e_{25},e_{26},e_{27},e_{28}\}$;
$w_3(e_{14}) = \{e_2,e_4,e_8,e_{12},e_{15},e_{16},e_{20},e_{21},e_{22},e_{23},e_{26},e_{29}\}$;
$w_3(e_{15}) = \{e_{12},e_{14},e_{15},e_{18},e_{19},e_{21},e_{23},e_{24},e_{27},e_{28}\}$;
$w_3(e_{16}) = \{e_2,e_4,e_8,e_{14},e_{16},e_{18},e_{19},e_{20},e_{22},e_{24},e_{26},e_{27},e_{28},e_{29}\}$;
$w_3(e_{17}) = \{e_1,e_2,e_3,e_8,e_{17},e_{21},e_{22},e_{23},e_{26},e_{27},e_{28},e_{29}\}$;
$w_3(e_{18}) = \{e_2,e_4,e_{10},e_{11},e_{13},e_{15},e_{16},e_{18},e_{19},e_{23},e_{25},e_{26},e_{27},e_{28}\}$;
$w_3(e_{19}) = \{e_2,e_4,e_{10},e_{11},e_{13},e_{15},e_{16},e_{18},e_{19},e_{23},e_{25},e_{26},e_{27},e_{28}\}$;
$w_3(e_{20}) = \{e_1,e_3,e_5,e_6,e_8,e_{10},e_{11},e_{14},e_{16},e_{23},e_{26},e_{27},e_{28},e_{29}\}$;
$w_3(e_{21}) = \{e_2,e_3,e_4,e_7,e_8,e_{13},e_{14},e_{15},e_{17},e_{22},e_{23},e_{24}\}$;
$w_3(e_{22}) = \{e_2,e_5,e_6,e_{10},e_{11},e_{14},e_{16},e_{17},e_{21},e_{22}\}$;
$w_3(e_{23}) = \{e_2,e_3,e_7,e_8,e_{10},e_{11},e_{12},e_{14},e_{15},e_{17},e_{18},e_{19},e_{20},e_{21},e_{24},e_{25},e_{26}\}$;
$w_3(e_{24}) = \{e_3,e_4,e_5,e_6,e_7,e_8,e_{10},e_{11},e_{13},e_{15},e_{16},e_{21},e_{23},e_{24}\}$;
$w_3(e_{25}) = \{e_1,e_4,e_7,e_{12},e_{13},e_{18},e_{19},e_{23},e_{26},e_{29}\}$;
$w_3(e_{26}) = \{e_2,e_4,e_5,e_6,e_{12},e_{13},e_{14},e_{16},e_{17},e_{18},e_{19},e_{20},e_{23},e_{25},e_{26}\}$;
$w_3(e_{27}) = \{e_3,e_4,e_7,e_8,e_{13},e_{15},e_{16},e_{17},e_{18},e_{19},e_{20}\}$;
$w_3(e_{28}) = \{e_3,e_4,e_7,e_8,e_{13},e_{15},e_{16},e_{17},e_{18},e_{19},e_{20}\}$;
$w_3(e_{29}) = \{e_1,e_2,e_5,e_6,e_7,e_{14},e_{16},e_{17},e_{20},e_{25},e_{29}\}$;
$w_3(e_{30}) = \varnothing$.

Реберные разрезы 5-го уровня:

$w_4(e_1) = \{e_3,e_4,e_7,e_8,e_{10},e_{11},e_{14},e_{16},e_{17},e_{21},e_{22}\}$;
$w_4(e_2) = \{e_3,e_5,e_6,e_7,e_{10},e_{11},e_{13},e_{15},e_{16},e_{21},e_{24},e_{26},e_{27},e_{28},e_{29}\}$;
$w_4(e_3) = \{e_1,e_2,e_7,e_8,e_{10},e_{11},e_{12},e_{15},e_{16},e_{18},e_{19},e_{21},e_{23},e_{24},e_{27},e_{28}\}$;
$w_4(e_4) = \{e_1,e_5,e_6,e_7,e_{12},e_{14},e_{15},e_{17},e_{21},e_{22},e_{25},e_{27},e_{28},e_{29}\}$;
$w_4(e_5) = \{e_2,e_4,e_8,e_{10},e_{11},e_{12},e_{13},e_{14},e_{16},e_{18},e_{19},e_{23},e_{26},e_{29}\}$;
$w_4(e_6) = \{e_2,e_4,e_8,e_{10},e_{11},e_{12},e_{13},e_{14},e_{16},e_{18},e_{19},e_{23},e_{26},e_{29}\}$;
$w_4(e_7) = \{e_1,e_2,e_3,e_4,e_{12},e_{14},e_{15},e_{17},e_{18},e_{19},e_{22},e_{23},e_{24},e_{27},e_{28}\}$;
$w_4(e_8) = \{e_1,e_3,e_5,e_6,e_{10},e_{11},e_{14},e_{16},e_{18},e_{19},e_{21},e_{24},e_{25},e_{26},e_{27},e_{28}\}$;
$w_4(e_9) = \varnothing$;

$w_4(e_{10}) = \{e_1,e_2,e_3,e_5,e_6,e_8,e_{13},e_{14},e_{15},e_{21},e_{24},e_{26},e_{27},e_{28},e_{29}\}$;
$w_4(e_{11}) = \{e_1,e_2,e_3,e_5,e_6,e_8,e_{13},e_{14},e_{15},e_{21},e_{24},e_{26},e_{27},e_{28},e_{29}\}$;
$w_4(e_{12}) = \{e_3,e_4,e_5,e_6,e_7,e_{14},e_{16},e_{20},e_{21},e_{22},e_{25},e_{29}\}$;
$w_4(e_{13}) = \{e_2,e_5,e_6,e_{10},e_{11},e_{14},e_{16},e_{17},e_{21},e_{22}\}$;
$w_4(e_{14}) = \{e_1,e_4,e_5,e_6,e_7,e_8,e_{10},e_{11},e_{12},e_{13},e_{17},e_{18},e_{19},e_{20},e_{23},e_{25},e_{26}\}$;
$w_4(e_{15}) = \{e_2,e_3,e_4,e_7,e_{10},e_{11},e_{17},e_{20},e_{25},e_{29}\}$;
$w_4(e_{16}) = \{e_1,e_2,e_3,e_5,e_6,e_8,e_{12},e_{13},e_{18},e_{19},e_{23},e_{26},e_{29}\}$;
$w_4(e_{17}) = \{e_1,e_4,e_7,e_{13},e_{14},e_{15},e_{21},e_{24},e_{26},e_{27},e_{28},e_{29}\}$;
$w_4(e_{18}) = \{e_3,e_5,e_6,e_7,e_8,e_{14},e_{16},e_{20},e_{21},e_{22},e_{23},e_{25},e_{26},e_{27},e_{28}\}$;
$w_4(e_{19}) = \{e_3,e_5,e_6,e_7,e_8,e_{14},e_{16},e_{20},e_{21},e_{22},e_{23},e_{25},e_{26},e_{27},e_{28}\}$;
$w_4(e_{20}) = \{e_{12},e_{14},e_{15},e_{18},e_{19},e_{21},e_{23},e_{24},e_{27},e_{28}\}$;
$w_4(e_{21}) = \{e_1,e_2,e_3,e_4,e_8,e_{10},e_{11},e_{12},e_{13},e_{17},e_{18},e_{19},e_{20},e_{25},e_{27},e_{28},e_{29}\}$;
$w_4(e_{22}) = \{e_1,e_4,e_7,e_{12},e_{13},e_{18},e_{19},e_{23},e_{26},e_{29}\}$;
$w_4(e_{23}) = \{e_3,e_5,e_6,e_7,e_{14},e_{16},e_{18},e_{19},e_{20},e_{22},e_{24},e_{26},e_{27},e_{28},e_{29}\}$;
$w_4(e_{24}) = \{e_2,e_3,e_7,e_8,e_{10},e_{11},e_{17},e_{20},e_{23},e_{25},e_{26},e_{27},e_{28}\}$;
$w_4(e_{25}) = \{e_4,e_8,e_{12},e_{14},e_{15},e_{18},e_{19},e_{21},e_{24},e_{26},e_{29}\}$;
$w_4(e_{26}) = \{e_2,e_5,e_6,e_8,e_{10},e_{11},e_{14},e_{16},e_{17},e_{18},e_{19},e_{22},e_{23},e_{24},e_{25},e_{29}\}$;
$w_4(e_{27}) = \{e_2,e_3,e_4,e_7,e_8,e_{10},e_{11},e_{17},e_{18},e_{19},e_{20},e_{21},e_{23},e_{24}\}$;
$w_4(e_{28}) = \{e_2,e_3,e_4,e_7,e_8,e_{10},e_{11},e_{17},e_{18},e_{19},e_{20},e_{21},e_{23},e_{24}\}$;
$w_4(e_{29}) = \{e_2,e_4,e_5,e_6,e_{10},e_{11},e_{12},e_{15},e_{16},e_{17},e_{21},e_{22},e_{23},e_{25},e_{26}\}$;
$w_4(e_{30}) = \varnothing$.

Реберные разрезы 6-го уровня:

$w_5(e_1) = \{e_3,e_4,e_5,e_6,e_7,e_8,e_{10},e_{11},e_{12},e_{13},e_{14},e_{16},e_{18},e_{19},e_{27},e_{28}\}$;
$w_5(e_2) = \{e_5,e_6,e_{10},e_{11},e_{17},e_{18},e_{19},e_{20},e_{21},e_{23},e_{24}\}$;
$w_5(e_3) = \{e_1,e_3,e_4,e_{12},e_{15},e_{16},e_{17},e_{20},e_{27},e_{28}\}$;
$w_5(e_4) = \{e_1,e_3,e_5,e_6,e_{13},e_{15},e_{16},e_{20},e_{22},e_{24},e_{25},e_{26},e_{27},e_{28}\}$;
$w_5(e_5) = \{e_1,e_2,e_4,e_7,e_8,e_{10},e_{11},e_{14},e_{16},e_{18},e_{19},e_{21},e_{24},e_{25},e_{26},e_{27},e_{28}\}$;
$w_5(e_6) = \{e_1,e_2,e_4,e_7,e_8,e_{10},e_{11},e_{14},e_{16},e_{18},e_{19},e_{21},e_{24},e_{25},e_{26},e_{27},e_{28}\}$;
$w_5(e_7) = \{e_1,e_5,e_6,e_7,e_8,e_{10},e_{11},e_{13},e_{14},e_{15},e_{17},e_{18},e_{19},e_{20}\}$;
$w_5(e_8) = \{e_1,e_5,e_6,e_7,e_{12},e_{14},e_{15},e_{17},e_{21},e_{22},e_{25},e_{27},e_{28},e_{29}\}$;
$w_5(e_9) = \varnothing$;
$w_5(e_{10}) = \{e_1,e_2,e_5,e_6,e_7,e_{12},e_{15},e_{16},e_{17},e_{18},e_{19},e_{20},e_{21},e_{23},e_{24},e_{25},e_{27},e_{28},e_{29}\}$;
$w_5(e_{11}) = \{e_1,e_2,e_5,e_6,e_7,e_{12},e_{15},e_{16},e_{17},e_{18},e_{19},e_{20},e_{21},e_{23},e_{24},e_{25},e_{27},e_{28},e_{29}\}$;
$w_5(e_{12}) = \{e_1,e_3,e_8,e_{10},e_{11},e_{12},e_{13},e_{14},e_{16},e_{17},e_{18},e_{19},e_{21},e_{22},e_{23},e_{26},e_{29}\}$;
$w_5(e_{13}) = \{e_1,e_4,e_7,e_{12},e_{13},e_{18},e_{19},e_{23},e_{26},e_{29}\}$;
$w_5(e_{14}) = \{e_1,e_5,e_6,e_7,e_8,e_{12},e_{14},e_{15},e_{17},e_{18},e_{19},e_{22},e_{23},e_{24},e_{27},e_{28}\}$;
$w_5(e_{15}) = \{e_3,e_4,e_7,e_8,e_{10},e_{11},e_{14},e_{16},e_{17},e_{21},e_{22}\}$;
$w_5(e_{16}) = \{e_1,e_3,e_4,e_5,e_6,e_{10},e_{11},e_{12},e_{15},e_{16},e_{18},e_{19},e_{21},e_{23},e_{24},e_{27},e_{28}\}$;
$w_5(e_{17}) = \{e_2,e_3,e_7,e_8,e_{10},e_{11},e_{12},e_{14},e_{15},e_{17},e_{18},e_{19},e_{20},e_{21},e_{24},e_{25},e_{26}\}$;
$w_5(e_{18}) = \{e_1,e_2,e_5,e_6,e_7,e_{10},e_{11},e_{12},e_{13},e_{14},e_{16},e_{17},e_{22},e_{23},e_{24},e_{25},e_{27},e_{28},e_{29}\}$;
$w_5(e_{19}) = \{e_1,e_2,e_5,e_6,e_7,e_{10},e_{11},e_{12},e_{13},e_{14},e_{16},e_{17},e_{22},e_{23},e_{24},e_{25},e_{27},e_{28},e_{29}\}$;
$w_5(e_{20}) = \{e_2,e_3,e_4,e_7,e_{10},e_{11},e_{17},e_{20},e_{25},e_{29}\}$;
$w_5(e_{21}) = \{e_2,e_5,e_6,e_8,e_{10},e_{11},e_{12},e_{15},e_{16},e_{17},e_{21},e_{22},e_{25},e_{27},e_{28},e_{29}\}$;
$w_5(e_{22}) = \{e_4,e_8,e_{12},e_{14},e_{15},e_{18},e_{19},e_{21},e_{24},e_{26},e_{29}\}$;
$w_5(e_{23}) = \{e_2,e_{10},e_{11},e_{12},e_{13},e_{14},e_{16},e_{18},e_{19},e_{27},e_{28}\}$;
$w_5(e_{24}) = \{e_2,e_4,e_5,e_6,e_{10},e_{11},e_{14},e_{16},e_{17},e_{18},e_{19},e_{22},e_{24},e_{25},e_{26},e_{27},e_{28}\}$;
$w_5(e_{25}) = \{e_4,e_5,e_6,e_8,e_{10},e_{11},e_{17},e_{18},e_{19},e_{20},e_{21},e_{24},e_{26},e_{27},e_{28},e_{29}\}$;
$w_5(e_{26}) = \{e_4,e_5,e_6,e_{12},e_{13},e_{17},e_{22},e_{24},e_{25},e_{26}\}$;
$w_5(e_{27}) = \{e_1,e_3,e_4,e_5,e_6,e_8,e_{10},e_{11},e_{14},e_{16},e_{18},e_{19},e_{21},e_{23},e_{24},e_{25},e_{29}\}$;
$w_5(e_{28}) = \{e_1,e_3,e_4,e_5,e_6,e_8,e_{10},e_{11},e_{14},e_{16},e_{18},e_{19},e_{21},e_{23},e_{24},e_{25},e_{29}\}$;
$w_5(e_{29}) = \{e_8,e_{10},e_{11},e_{12},e_{13},e_{18},e_{19},e_{20},e_{21},e_{22},e_{25},e_{27},e_{28},e_{29}\}$;

$w_5(e_{30}) = \varnothing$.

Реберные разрезы 7-го уровня:

$w_6(e_1) = \{e_4,e_8,e_{23},e_{26},e_{27},e_{28},e_{29}\}$;
$w_6(e_2) = \{e_4,e_{18},e_{19},e_{21},e_{24},e_{25},e_{26},e_{27},e_{28}\}$;
$w_6(e_3) = \{e_5,e_6,e_8,e_{10},e_{11},e_{12},e_{14},e_{15},e_{17},e_{18},e_{19},e_{20},e_{21},e_{23},e_{24},e_{25},e_{27},e_{28},e_{29}\}$;
$w_6(e_4) = \{e_1,e_2,e_7,e_{13},e_{15},e_{16},e_{20},e_{22},e_{23},e_{24},e_{25},e_{29}\}$;
$w_6(e_5) = \{e_3,e_7,e_8,e_{10},e_{11},e_{12},e_{15},e_{16},e_{17},e_{21},e_{22},e_{23},e_{25},e_{26}\}$;
$w_6(e_6) = \{e_3,e_7,e_8,e_{10},e_{11},e_{12},e_{15},e_{16},e_{17},e_{21},e_{22},e_{23},e_{25},e_{26}\}$;
$w_6(e_7) = \{e_4,e_5,e_6,e_{10},e_{11},e_{12},e_{14},e_{15},e_{17},e_{18},e_{19},e_{20},e_{21},e_{24},e_{25},e_{26}\}$;
$w_6(e_8) = \{e_1,e_3,e_5,e_6,e_{13},e_{15},e_{16},e_{20},e_{22},e_{24},e_{25},e_{26},e_{27},e_{28}\}$;
$w_6(e_9) = \varnothing$;
$w_6(e_{10}) = \{e_3,e_5,e_6,e_7,e_{12},e_{15},e_{16},e_{20},e_{21},e_{22},e_{23},e_{26},e_{29}\}$;
$w_6(e_{11}) = \{e_3,e_5,e_6,e_7,e_{12},e_{15},e_{16},e_{20},e_{21},e_{22},e_{23},e_{26},e_{29}\}$;
$w_6(e_{12}) = \{e_3,e_5,e_6,e_7,e_{10},e_{11},e_{13},e_{15},e_{16},e_{21},e_{24},e_{26},e_{27},e_{28},e_{29}\}$;
$w_6(e_{13}) = \{e_4,e_8,e_{12},e_{14},e_{15},e_{18},e_{19},e_{21},e_{24},e_{26},e_{29}\}$;
$w_6(e_{14}) = \{e_3,e_7,e_{13},e_{15},e_{16},e_{17},e_{18},e_{19},e_{20},e_{23},e_{26},e_{27},e_{28},e_{29}\}$;
$w_6(e_{15}) = \{e_3,e_4,e_5,e_6,e_7,e_8,e_{10},e_{11},e_{12},e_{13},e_{14},e_{16},e_{18},e_{19},e_{27},e_{28}\}$;
$w_6(e_{16}) = \{e_4,e_5,e_6,e_8,e_{10},e_{11},e_{12},e_{14},e_{15},e_{17},e_{20},e_{23},e_{26},e_{29}\}$;
$w_6(e_{17}) = \{e_3,e_5,e_6,e_7,e_{14},e_{16},e_{18},e_{19},e_{20},e_{22},e_{24},e_{26},e_{27},e_{28},e_{29}\}$;
$w_6(e_{18}) = \{e_2,e_3,e_7,e_{13},e_{14},e_{15},e_{17},e_{22},e_{24},e_{26},e_{27},e_{28},e_{29}\}$;
$w_6(e_{19}) = \{e_2,e_3,e_7,e_{13},e_{14},e_{15},e_{17},e_{22},e_{24},e_{26},e_{27},e_{28},e_{29}\}$;
$w_6(e_{20}) = \{e_3,e_4,e_7,e_8,e_{10},e_{11},e_{14},e_{16},e_{17},e_{21},e_{22}\}$;
$w_6(e_{21}) = \{e_2,e_3,e_5,e_6,e_7,e_{10},e_{11},e_{12},e_{13},e_{20},e_{22},e_{24},e_{26},e_{29}\}$;
$w_6(e_{22}) = \{e_4,e_5,e_6,e_8,e_{10},e_{11},e_{17},e_{18},e_{19},e_{20},e_{21},e_{24},e_{26},e_{27},e_{28},e_{29}\}$;
$w_6(e_{23}) = \{e_1,e_3,e_4,e_5,e_6,e_{10},e_{11},e_{14},e_{16}\}$;
$w_6(e_{24}) = \{e_2,e_3,e_4,e_7,e_8,e_{12},e_{13},e_{17},e_{18},e_{19},e_{21},e_{22},e_{27},e_{28}\}$;
$w_6(e_{25}) = \{e_2,e_3,e_4,e_5,e_6,e_7,e_8\}$;
$w_6(e_{26}) = \{e_1,e_2,e_5,e_6,e_7,e_8,e_{10},e_{11},e_{12},e_{13},e_{14},e_{16},e_{17},e_{18},e_{19},e_{21},e_{22},e_{27},e_{28}\}$;
$w_6(e_{27}) = \{e_1,e_2,e_3,e_8,e_{12},e_{14},e_{15},e_{17},e_{18},e_{19},e_{22},e_{24},e_{26},e_{29}\}$;
$w_6(e_{28}) = \{e_1,e_2,e_3,e_8,e_{12},e_{14},e_{15},e_{17},e_{18},e_{19},e_{22},e_{24},e_{26},e_{29}\}$;
$w_6(e_{29}) = \{e_1,e_3,e_4,e_{10},e_{11},e_{12},e_{13},e_{14},e_{16},e_{17},e_{18},e_{19},e_{21},e_{22},e_{27},e_{28}\}$;
$w_6(e_{30}) = \varnothing$.

Реберные разрезы 8-го уровня:

$w_7(e_1) = \{e_2,e_3,e_5,e_6,e_7,e_8,e_{18},e_{19},e_{21},e_{24},e_{25},e_{26},e_{27},e_{28}\}$;
$w_7(e_2) = \{e_1,e_2,e_3,e_8,e_{17},e_{21},e_{22},e_{23},e_{26},e_{27},e_{28},e_{29}\}$;
$w_7(e_3) = \{e_1,e_2,e_3,e_5,e_6,e_{10},e_{11},e_{20},e_{21},e_{22},e_{23},e_{25},e_{26},e_{27},e_{28}\}$;
$w_7(e_4) = \{e_{10},e_{11},e_{13},e_{14},e_{15},e_{18},e_{19},e_{20},e_{21},e_{22}\}$;
$w_7(e_5) = \{e_1,e_3,e_5,e_6,e_{12},e_{13},e_{14},e_{16},e_{18},e_{19},e_{20},e_{21},e_{22},e_{23},e_{25},e_{26}\}$;
$w_7(e_6) = \{e_1,e_3,e_5,e_6,e_{12},e_{13},e_{14},e_{16},e_{18},e_{19},e_{20},e_{21},e_{22},e_{23},e_{25},e_{26}\}$;
$w_7(e_7) = \{e_1,e_7,e_8,e_{10},e_{11},e_{18},e_{19},e_{20},e_{22},e_{23},e_{24}\}$;
$w_7(e_8) = \{e_1,e_2,e_7,e_{13},e_{15},e_{16},e_{20},e_{22},e_{23},e_{24},e_{25},e_{29}\}$;
$w_7(e_9) = \varnothing$;
$w_7(e_{10}) = \{e_3,e_4,e_7,e_{12},e_{13},e_{14},e_{16},e_{17},e_{18},e_{19},e_{20},e_{25},e_{27},e_{28},e_{29}\}$;
$w_7(e_{11}) = \{e_3,e_4,e_7,e_{12},e_{13},e_{14},e_{16},e_{17},e_{18},e_{19},e_{20},e_{25},e_{27},e_{28},e_{29}\}$;
$w_7(e_{12}) = \{e_5,e_6,e_{10},e_{11},e_{17},e_{18},e_{19},e_{20},e_{21},e_{23},e_{24}\}$;
$w_7(e_{13}) = \{e_4,e_5,e_6,e_8,e_{10},e_{11},e_{17},e_{18},e_{19},e_{20},e_{21},e_{24},e_{26},e_{27},e_{28},e_{29}\}$;
$w_7(e_{14}) = \{e_4,e_5,e_6,e_{10},e_{11},e_{17},e_{20},e_{23},e_{25},e_{26},e_{27},e_{28}\}$;
$w_7(e_{15}) = \{e_4,e_8,e_{23},e_{26},e_{27},e_{28},e_{29}\}$;
$w_7(e_{16}) = \{e_5,e_6,e_8,e_{10},e_{11},e_{17},e_{20},e_{25},e_{29}\}$;
$w_7(e_{17}) = \{e_2,e_{10},e_{11},e_{12},e_{13},e_{14},e_{16},e_{18},e_{19},e_{27},e_{28}\}$;

$w_7(e_{18}) = \{e_1,e_4,e_5,e_6,e_7,e_{10},e_{11},e_{12},e_{13},e_{17},e_{20},e_{21},e_{24},e_{26},e_{29}\}$;
$w_7(e_{19}) = \{e_1,e_4,e_5,e_6,e_7,e_{10},e_{11},e_{12},e_{13},e_{17},e_{20},e_{21},e_{24},e_{26},e_{29}\}$;
$w_7(e_{20}) = \{e_3,e_4,e_5,e_6,e_7,e_8,e_{10},e_{11},e_{12},e_{13},e_{14},e_{16},e_{18},e_{19},e_{27},e_{28}\}$;
$w_7(e_{21}) = \{e_1,e_2,e_3,e_4,e_5,e_6,e_{12},e_{13},e_{18},e_{19},e_{27},e_{28}\}$;
$w_7(e_{22}) = \{e_2,e_3,e_4,e_5,e_6,e_7,e_8\}$;
$w_7(e_{23}) = \{e_2,e_3,e_5,e_6,e_7,e_8,e_{12},e_{14},e_{15},e_{23},e_{25},e_{26}\}$;
$w_7(e_{24}) = \{e_1,e_7,e_8,e_{12},e_{13},e_{18},e_{19},e_{27},e_{28}\}$;
$w_7(e_{25}) = \{e_1,e_3,e_5,e_6,e_8,e_{10},e_{11},e_{14},e_{16},e_{23},e_{26},e_{27},e_{28},e_{29}\}$;
$w_7(e_{26}) = \{e_1,e_2,e_3,e_5,e_6,e_{13},e_{14},e_{15},e_{18},e_{19},e_{23},e_{25},e_{26},e_{27},e_{28}\}$;
$w_7(e_{27}) = \{e_1,e_2,e_3,e_{10},e_{11},e_{13},e_{14},e_{15},e_{17},e_{20},e_{21},e_{24},e_{25},e_{26},e_{27},e_{28}\}$;
$w_7(e_{28}) = \{e_1,e_2,e_3,e_{10},e_{11},e_{13},e_{14},e_{15},e_{17},e_{20},e_{21},e_{24},e_{25},e_{26},e_{27},e_{28}\}$;
$w_7(e_{29}) = \{e_2,e_8,e_{10},e_{11},e_{13},e_{15},e_{16},e_{18},e_{19},e_{25},e_{29}\}$;
$w_7(e_{30}) = \varnothing$.

Реберные разрезы 9-го уровня:

$w_8(e_1) = \{e_2,e_5,e_6,e_{10},e_{11},e_{14},e_{16},e_{17},e_{21},e_{22}\}$;
$w_8(e_2) = \{e_1,e_4,e_7,e_{13},e_{14},e_{15},e_{21},e_{24},e_{26},e_{27},e_{28},e_{29}\}$;
$w_8(e_3) = \{e_5,e_6,e_{13},e_{14},e_{15},e_{17},e_{22},e_{23},e_{24}\}$;
$w_8(e_4) = \varnothing$;
$w_8(e_5) = \{e_1,e_3,e_4,e_{10},e_{11},e_{13},e_{15},e_{16},e_{17},e_{22},e_{23},e_{24}\}$;
$w_8(e_6) = \{e_1,e_3,e_4,e_{10},e_{11},e_{13},e_{15},e_{16},e_{17},e_{22},e_{23},e_{24}\}$;
$w_8(e_7) = \{e_2,e_{10},e_{11},e_{13},e_{15},e_{16},e_{21},e_{23},e_{24}\}$;
$w_8(e_8) = \varnothing$;
$w_8(e_9) = \varnothing$;
$w_8(e_{10}) = \{e_1,e_5,e_6,e_7,e_8,e_{17},e_{21},e_{22}\}$;
$w_8(e_{11}) = \{e_1,e_5,e_6,e_7,e_8,e_{17},e_{21},e_{22}\}$;
$w_8(e_{12}) = \{e_4,e_{18},e_{19},e_{21},e_{24},e_{25},e_{26},e_{27},e_{28}\}$;
$w_8(e_{13}) = \{e_2,e_3,e_4,e_5,e_6,e_7,e_8\}$;
$w_8(e_{14}) = \{e_1,e_2,e_3,e_4,e_8,e_{17},e_{18},e_{19},e_{22},e_{23},e_{24},e_{25},e_{29}\}$;
$w_8(e_{15}) = \{e_2,e_3,e_5,e_6,e_7,e_8,e_{18},e_{19},e_{21},e_{24},e_{25},e_{26},e_{27},e_{28}\}$;
$w_8(e_{16}) = \{e_1,e_4,e_5,e_6,e_7,e_{17},e_{21},e_{22},e_{23},e_{26},e_{27},e_{28},e_{29}\}$;
$w_8(e_{17}) = \{e_1,e_3,e_4,e_5,e_6,e_{10},e_{11},e_{14},e_{16}\}$;
$w_8(e_{18}) = \{e_8,e_{12},e_{14},e_{15},e_{25},e_{27},e_{28},e_{29}\}$;
$w_8(e_{19}) = \{e_8,e_{12},e_{14},e_{15},e_{25},e_{27},e_{28},e_{29}\}$;
$w_8(e_{20}) = \{e_4,e_8,e_{23},e_{26},e_{27},e_{28},e_{29}\}$;
$w_8(e_{21}) = \{e_1,e_2,e_4,e_7,e_8,e_{10},e_{11},e_{12},e_{15},e_{16},e_{23},e_{25},e_{26}\}$;
$w_8(e_{22}) = \{e_1,e_3,e_5,e_6,e_8,e_{10},e_{11},e_{14},e_{16},e_{23},e_{26},e_{27},e_{28},e_{29}\}$;
$w_8(e_{23}) = \{e_3,e_4,e_5,e_6,e_7,e_{14},e_{16},e_{20},e_{21},e_{22},e_{25},e_{29}\}$;
$w_8(e_{24}) = \{e_2,e_3,e_4,e_5,e_6,e_7,e_{12},e_{14},e_{15},e_{25},e_{27},e_{28},e_{29}\}$;
$w_8(e_{25}) = \{e_{12},e_{14},e_{15},e_{18},e_{19},e_{21},e_{23},e_{24},e_{27},e_{28}\}$;
$w_8(e_{26}) = \{e_2,e_{12},e_{15},e_{16},e_{20},e_{21},e_{22},e_{27},e_{28}\}$;
$w_8(e_{27}) = \{e_2,e_4,e_{12},e_{15},e_{16},e_{18},e_{19},e_{20},e_{22},e_{24},e_{25},e_{26}\}$;
$w_8(e_{28}) = \{e_2,e_4,e_{12},e_{15},e_{16},e_{18},e_{19},e_{20},e_{22},e_{24},e_{25},e_{26}\}$;
$w_8(e_{29}) = \{e_2,e_{14},e_{16},e_{18},e_{19},e_{20},e_{22},e_{23},e_{24}\}$;
$w_8(e_{30}) = \varnothing$.

Реберные разрезы 10-го уровня:

$w_9(e_1) = \{e_1,e_4,e_7,e_{12},e_{13},e_{18},e_{19},e_{23},e_{26},e_{29}\}$;
$w_9(e_2) = \{e_2,e_3,e_7,e_8,e_{10},e_{11},e_{12},e_{14},e_{15},e_{17},e_{18},e_{19},e_{20},e_{21},e_{24},e_{25},e_{26}\}$;
$w_9(e_3) = \{e_2,e_4,e_5,e_6,e_{12},e_{13},e_{14},e_{16},e_{17},e_{18},e_{19},e_{20},e_{23},e_{25},e_{26}\}$;
$w_9(e_4) = \varnothing$;
$w_9(e_5) = \{e_3,e_4,e_7,e_8,e_{13},e_{15},e_{16},e_{17},e_{18},e_{19},e_{20}\}$;

$w_9(e_6) = \{e_3,e_4,e_7,e_8,e_{13},e_{15},e_{16},e_{17},e_{18},e_{19},e_{20}\}$;
$w_9(e_7) = \{e_1,e_2,e_5,e_6,e_7,e_{14},e_{16},e_{17},e_{20},e_{25},e_{29}\}$;
$w_9(e_8) = \varnothing$;
$w_9(e_9) = \varnothing$;
$w_9(e_{10}) = \{e_2,e_4,e_{10},e_{11},e_{13},e_{15},e_{16},e_{18},e_{19},e_{23},e_{25},e_{26},e_{27},e_{28}\}$;
$w_9(e_{11}) = \{e_2,e_4,e_{10},e_{11},e_{13},e_{15},e_{16},e_{18},e_{19},e_{23},e_{25},e_{26},e_{27},e_{28}\}$;
$w_9(e_{12}) = \{e_1,e_2,e_3,e_8,e_{17},e_{21},e_{22},e_{23},e_{26},e_{27},e_{28},e_{29}\}$;
$w_9(e_{13}) = \{e_1,e_3,e_5,e_6,e_8,e_{10},e_{11},e_{14},e_{16},e_{23},e_{26},e_{27},e_{28},e_{29}\}$;
$w_9(e_{14}) = \{e_2,e_3,e_4,e_7,e_8,e_{13},e_{14},e_{15},e_{17},e_{22},e_{23},e_{24}\}$;
$w_9(e_{15}) = \{e_2,e_5,e_6,e_{10},e_{11},e_{14},e_{16},e_{17},e_{21},e_{22}\}$;
$w_9(e_{16}) = \{e_3,e_4,e_5,e_6,e_7,e_8,e_{10},e_{11},e_{13},e_{15},e_{16},e_{21},e_{23},e_{24}\}$;
$w_9(e_{17}) = \{e_2,e_3,e_5,e_6,e_7,e_8,e_{12},e_{14},e_{15},e_{23},e_{25},e_{26}\}$;
$w_9(e_{18}) = \{e_1,e_2,e_3,e_4,e_5,e_6,e_{10},e_{11},e_{18},e_{19},e_{20},e_{22},e_{23},e_{24}\}$;
$w_9(e_{19}) = \{e_1,e_2,e_3,e_4,e_5,e_6,e_{10},e_{11},e_{18},e_{19},e_{20},e_{22},e_{23},e_{24}\}$;
$w_9(e_{20}) = \{e_2,e_3,e_5,e_6,e_7,e_8,e_{18},e_{19},e_{21},e_{24},e_{25},e_{26},e_{27},e_{28}\}$;
$w_9(e_{21}) = \{e_2,e_4,e_8,e_{12},e_{15},e_{16},e_{20},e_{21},e_{22},e_{23},e_{26},e_{29}\}$;
$w_9(e_{22}) = \{e_{12},e_{14},e_{15},e_{18},e_{19},e_{21},e_{23},e_{24},e_{27},e_{28}\}$;
$w_9(e_{23}) = \{e_1,e_3,e_8,e_{10},e_{11},e_{12},e_{13},e_{14},e_{16},e_{17},e_{18},e_{19},e_{21},e_{22},e_{23},e_{26},e_{29}\}$;
$w_9(e_{24}) = \{e_2,e_4,e_8,e_{14},e_{16},e_{18},e_{19},e_{20},e_{22},e_{24},e_{26},e_{27},e_{28},e_{29}\}$;
$w_9(e_{25}) = \{e_2,e_3,e_4,e_7,e_{10},e_{11},e_{17},e_{20},e_{25},e_{29}\}$;
$w_9(e_{26}) = \{e_1,e_2,e_3,e_4,e_{10},e_{11},e_{12},e_{13},e_{17},e_{20},e_{21},e_{23},e_{24},e_{27},e_{28}\}$;
$w_9(e_{27}) = \{e_4,e_8,e_{10},e_{11},e_{12},e_{13},e_{20},e_{22},e_{24},e_{26},e_{29}\}$;
$w_9(e_{28}) = \{e_4,e_8,e_{10},e_{11},e_{12},e_{13},e_{20},e_{22},e_{24},e_{26},e_{29}\}$;
$w_9(e_{29}) = \{e_1,e_7,e_{12},e_{13},e_{21},e_{23},e_{24},e_{25},e_{27},e_{28},e_{29}\}$;
$w_9(e_{30}) = \varnothing$.

Реберные разрезы 11-го уровня:

$w_{10}(e_1) = \{e_4,e_8,e_{12},e_{14},e_{15},e_{18},e_{19},e_{21},e_{24},e_{26},e_{29}\}$;
$w_{10}(e_2) = \{e_3,e_5,e_6,e_7,e_{14},e_{16},e_{18},e_{19},e_{20},e_{22},e_{24},e_{26},e_{27},e_{28},e_{29}\}$;
$w_{10}(e_3) = \{e_2,e_5,e_6,e_8,e_{10},e_{11},e_{14},e_{16},e_{17},e_{18},e_{19},e_{22},e_{23},e_{24},e_{25},e_{29}\}$;
$w_{10}(e_4) = \varnothing$;
$w_{10}(e_5) = \{e_2,e_3,e_4,e_7,e_8,e_{10},e_{11},e_{17},e_{18},e_{19},e_{20},e_{21},e_{23},e_{24}\}$;
$w_{10}(e_6) = \{e_2,e_3,e_4,e_7,e_8,e_{10},e_{11},e_{17},e_{18},e_{19},e_{20},e_{21},e_{23},e_{24}\}$;
$w_{10}(e_7) = \{e_2,e_4,e_5,e_6,e_{10},e_{11},e_{12},e_{15},e_{16},e_{17},e_{21},e_{22},e_{23},e_{25},e_{26}\}$;
$w_{10}(e_8) = \varnothing$;
$w_{10}(e_9) = \varnothing$;
$w_{10}(e_{10}) = \{e_3,e_5,e_6,e_7,e_8,e_{14},e_{16},e_{20},e_{21},e_{22},e_{23},e_{25},e_{26},e_{27},e_{28}\}$;
$w_{10}(e_{11}) = \{e_3,e_5,e_6,e_7,e_8,e_{14},e_{16},e_{20},e_{21},e_{22},e_{23},e_{25},e_{26},e_{27},e_{28}\}$;
$w_{10}(e_{12}) = \{e_1,e_4,e_7,e_{13},e_{14},e_{15},e_{21},e_{24},e_{26},e_{27},e_{28},e_{29}\}$;
$w_{10}(e_{13}) = \{e_{12},e_{14},e_{15},e_{18},e_{19},e_{21},e_{23},e_{24},e_{27},e_{28}\}$;
$w_{10}(e_{14}) = \{e_1,e_2,e_3,e_4,e_8,e_{10},e_{11},e_{12},e_{13},e_{17},e_{18},e_{19},e_{20},e_{25},e_{27},e_{28},e_{29}\}$;
$w_{10}(e_{15}) = \{e_1,e_4,e_7,e_{12},e_{13},e_{18},e_{19},e_{23},e_{26},e_{29}\}$;
$w_{10}(e_{16}) = \{e_2,e_3,e_7,e_8,e_{10},e_{11},e_{17},e_{20},e_{23},e_{25},e_{26},e_{27},e_{28}\}$;
$w_{10}(e_{17}) = \{e_3,e_4,e_5,e_6,e_7,e_{14},e_{16},e_{20},e_{21},e_{22},e_{25},e_{29}\}$;
$w_{10}(e_{18}) = \{e_1,e_2,e_3,e_5,e_6,e_8,e_{13},e_{14},e_{15},e_{21},e_{24},e_{26},e_{27},e_{28},e_{29}\}$;
$w_{10}(e_{19}) = \{e_1,e_2,e_3,e_5,e_6,e_8,e_{13},e_{14},e_{15},e_{21},e_{24},e_{26},e_{27},e_{28},e_{29}\}$;
$w_{10}(e_{20}) = \{e_2,e_5,e_6,e_{10},e_{11},e_{14},e_{16},e_{17},e_{21},e_{22}\}$;
$w_{10}(e_{21}) = \{e_1,e_4,e_5,e_6,e_7,e_8,e_{10},e_{11},e_{12},e_{13},e_{17},e_{18},e_{19},e_{20},e_{23},e_{25},e_{26}\}$;
$w_{10}(e_{22}) = \{e_2,e_3,e_4,e_7,e_{10},e_{11},e_{17},e_{20},e_{25},e_{29}\}$;
$w_{10}(e_{23}) = \{e_3,e_5,e_6,e_7,e_{10},e_{11},e_{13},e_{15},e_{16},e_{21},e_{24},e_{26},e_{27},e_{28},e_{29}\}$;
$w_{10}(e_{24}) = \{e_1,e_2,e_3,e_5,e_6,e_8,e_{12},e_{13},e_{18},e_{19},e_{23},e_{26},e_{29}\}$;
$w_{10}(e_{25}) = \{e_3,e_4,e_7,e_8,e_{10},e_{11},e_{14},e_{16},e_{17},e_{21},e_{22}\}$;

$w_{10}(e_{26}) = \{e_1,e_2,e_7,e_8,e_{10},e_{11},e_{12},e_{15},e_{16},e_{18},e_{19},e_{21},e_{23},e_{24},e_{27},e_{28}\}$;
$w_{10}(e_{27}) = \{e_2,e_4,e_8,e_{10},e_{11},e_{12},e_{13},e_{14},e_{16},e_{18},e_{19},e_{23},e_{26},e_{29}\}$;
$w_{10}(e_{28}) = \{e_2,e_4,e_8,e_{10},e_{11},e_{12},e_{13},e_{14},e_{16},e_{18},e_{19},e_{23},e_{26},e_{29}\}$;
$w_{10}(e_{29}) = \{e_1,e_2,e_3,e_4,e_{12},e_{14},e_{15},e_{17},e_{18},e_{19},e_{22},e_{23},e_{24},e_{27},e_{28}\}$;
$w_{10}(e_{30}) = \varnothing$.

Реберные разрезы 12-го уровня:

$w_{11}(e_1) = \{e_4,e_5,e_6,e_8,e_{10},e_{11},e_{17},e_{18},e_{19},e_{20},e_{21},e_{24},e_{26},e_{27},e_{28},e_{29}\}$;
$w_{11}(e_2) = \{e_2,e_{10},e_{11},e_{12},e_{13},e_{14},e_{16},e_{18},e_{19},e_{27},e_{28}\}$;
$w_{11}(e_3) = \{e_4,e_5,e_6,e_{12},e_{13},e_{17},e_{22},e_{24},e_{25},e_{26}\}$;
$w_{11}(e_4) = \varnothing$;
$w_{11}(e_5) = \{e_1,e_3,e_4,e_5,e_6,e_8,e_{10},e_{11},e_{14},e_{16},e_{18},e_{19},e_{21},e_{23},e_{24},e_{25},e_{29}\}$;
$w_{11}(e_6) = \{e_1,e_3,e_4,e_5,e_6,e_8,e_{10},e_{11},e_{14},e_{16},e_{18},e_{19},e_{21},e_{23},e_{24},e_{25},e_{29}\}$;
$w_{11}(e_7) = \{e_8,e_{10},e_{11},e_{12},e_{13},e_{18},e_{19},e_{20},e_{21},e_{22},e_{25},e_{27},e_{28},e_{29}\}$;
$w_{11}(e_8) = \varnothing$;
$w_{11}(e_9) = \varnothing$;
$w_{11}(e_{10}) = \{e_1,e_2,e_5,e_6,e_7,e_{10},e_{11},e_{12},e_{13},e_{14},e_{16},e_{17},e_{22},e_{23},e_{24},e_{25},e_{27},e_{28},e_{29}\}$;
$w_{11}(e_{11}) = \{e_1,e_2,e_5,e_6,e_7,e_{10},e_{11},e_{12},e_{13},e_{14},e_{16},e_{17},e_{22},e_{23},e_{24},e_{25},e_{27},e_{28},e_{29}\}$;
$w_{11}(e_{12}) = \{e_2,e_3,e_7,e_8,e_{10},e_{11},e_{12},e_{14},e_{15},e_{17},e_{18},e_{19},e_{20},e_{21},e_{24},e_{25},e_{26}\}$;
$w_{11}(e_{13}) = \{e_2,e_3,e_4,e_7,e_{10},e_{11},e_{17},e_{20},e_{25},e_{29}\}$;
$w_{11}(e_{14}) = \{e_2,e_5,e_6,e_8,e_{10},e_{11},e_{12},e_{15},e_{16},e_{17},e_{21},e_{22},e_{25},e_{27},e_{28},e_{29}\}$;
$w_{11}(e_{15}) = \{e_4,e_8,e_{12},e_{14},e_{15},e_{18},e_{19},e_{21},e_{24},e_{26},e_{29}\}$;
$w_{11}(e_{16}) = \{e_2,e_4,e_5,e_6,e_{10},e_{11},e_{14},e_{16},e_{17},e_{18},e_{19},e_{22},e_{24},e_{25},e_{26},e_{27},e_{28}\}$;
$w_{11}(e_{17}) = \{e_1,e_3,e_8,e_{10},e_{11},e_{12},e_{13},e_{14},e_{16},e_{17},e_{18},e_{19},e_{21},e_{22},e_{23},e_{26},e_{29}\}$;
$w_{11}(e_{18}) = \{e_1,e_2,e_5,e_6,e_7,e_{12},e_{15},e_{16},e_{17},e_{18},e_{19},e_{20},e_{21},e_{23},e_{24},e_{25},e_{27},e_{28},e_{29}\}$;
$w_{11}(e_{19}) = \{e_1,e_2,e_5,e_6,e_7,e_{12},e_{15},e_{16},e_{17},e_{18},e_{19},e_{20},e_{21},e_{23},e_{24},e_{25},e_{27},e_{28},e_{29}\}$;
$w_{11}(e_{20}) = \{e_1,e_4,e_7,e_{12},e_{13},e_{18},e_{19},e_{23},e_{26},e_{29}\}$;
$w_{11}(e_{21}) = \{e_1,e_5,e_6,e_7,e_8,e_{12},e_{14},e_{15},e_{17},e_{18},e_{19},e_{22},e_{23},e_{24},e_{27},e_{28}\}$;
$w_{11}(e_{22}) = \{e_3,e_4,e_7,e_8,e_{10},e_{11},e_{14},e_{16},e_{17},e_{21},e_{22}\}$;
$w_{11}(e_{23}) = \{e_5,e_6,e_{10},e_{11},e_{17},e_{18},e_{19},e_{20},e_{21},e_{23},e_{24}\}$;
$w_{11}(e_{24}) = \{e_1,e_3,e_4,e_5,e_6,e_{10},e_{11},e_{12},e_{15},e_{16},e_{18},e_{19},e_{21},e_{23},e_{24},e_{27},e_{28}\}$;
$w_{11}(e_{25}) = \{e_3,e_4,e_5,e_6,e_7,e_8,e_{10},e_{11},e_{12},e_{13},e_{14},e_{16},e_{18},e_{19},e_{27},e_{28}\}$;
$w_{11}(e_{26}) = \{e_1,e_3,e_4,e_{12},e_{15},e_{16},e_{17},e_{20},e_{27},e_{28}\}$;
$w_{11}(e_{27}) = \{e_1,e_2,e_4,e_7,e_8,e_{10},e_{11},e_{14},e_{16},e_{18},e_{19},e_{21},e_{24},e_{25},e_{26},e_{27},e_{28}\}$;
$w_{11}(e_{28}) = \{e_1,e_2,e_4,e_7,e_8,e_{10},e_{11},e_{14},e_{16},e_{18},e_{19},e_{21},e_{24},e_{25},e_{26},e_{27},e_{28}\}$;
$w_{11}(e_{29}) = \{e_1,e_5,e_6,e_7,e_8,e_{10},e_{11},e_{13},e_{14},e_{15},e_{17},e_{18},e_{19},e_{20}\}$;
$w_{11}(e_{30}) = \varnothing$.

Реберные разрезы 13-го уровня:

$w_{12}(e_1) = \varnothing$;
$w_{12}(e_2) = \{e_1,e_3,e_4,e_5,e_6,e_{10},e_{11},e_{14},e_{16}\}$;
$w_{12}(e_3) = \{e_1,e_2,e_5,e_6,e_7,e_8,e_{10},e_{11},e_{12},e_{13},e_{14},e_{16},e_{17},e_{18},e_{19},e_{21},e_{22},e_{27},e_{28}\}$;
$w_{12}(e_4) = \varnothing$;
$w_{12}(e_5) = \{e_1,e_2,e_3,e_8,e_{12},e_{14},e_{15},e_{17},e_{18},e_{19},e_{22},e_{24},e_{26},e_{29}\}$;
$w_{12}(e_6) = \{e_1,e_2,e_3,e_8,e_{12},e_{14},e_{15},e_{17},e_{18},e_{19},e_{22},e_{24},e_{26},e_{29}\}$;
$w_{12}(e_7) = \{e_1,e_3,e_4,e_{10},e_{11},e_{12},e_{13},e_{14},e_{16},e_{17},e_{18},e_{19},e_{21},e_{22},e_{27},e_{28}\}$;
$w_{12}(e_8) = \varnothing$;
$w_{12}(e_9) = \varnothing$;
$w_{12}(e_{10}) = \{e_2,e_3,e_7,e_{13},e_{14},e_{15},e_{17},e_{22},e_{24},e_{26},e_{27},e_{28},e_{29}\}$;
$w_{12}(e_{11}) = \{e_2,e_3,e_7,e_{13},e_{14},e_{15},e_{17},e_{22},e_{24},e_{26},e_{27},e_{28},e_{29}\}$;
$w_{12}(e_{12}) = \{e_3,e_5,e_6,e_7,e_{14},e_{16},e_{18},e_{19},e_{20},e_{22},e_{24},e_{26},e_{27},e_{28},e_{29}\}$;
$w_{12}(e_{13}) = \{e_3,e_4,e_7,e_8,e_{10},e_{11},e_{14},e_{16},e_{17},e_{21},e_{22}\}$;

$w_{12}(e_{14}) = \{e_2,e_3,e_5,e_6,e_7,e_{10},e_{11},e_{12},e_{13},e_{20},e_{22},e_{24},e_{26},e_{29}\}$;
$w_{12}(e_{15}) = \{e_4,e_5,e_6,e_8,e_{10},e_{11},e_{17},e_{18},e_{19},e_{20},e_{21},e_{24},e_{26},e_{27},e_{28},e_{29}\}$;
$w_{12}(e_{16}) = \{e_2,e_3,e_4,e_7,e_8,e_{12},e_{13},e_{17},e_{18},e_{19},e_{21},e_{22},e_{27},e_{28}\}$;
$w_{12}(e_{17}) = \{e_3,e_5,e_6,e_7,e_{10},e_{11},e_{13},e_{15},e_{16},e_{21},e_{24},e_{26},e_{27},e_{28},e_{29}\}$;
$w_{12}(e_{18}) = \{e_3,e_5,e_6,e_7,e_{12},e_{15},e_{16},e_{20},e_{21},e_{22},e_{23},e_{26},e_{29}\}$;
$w_{12}(e_{19}) = \{e_3,e_5,e_6,e_7,e_{12},e_{15},e_{16},e_{20},e_{21},e_{22},e_{23},e_{26},e_{29}\}$;
$w_{12}(e_{20}) = \{e_4,e_8,e_{12},e_{14},e_{15},e_{18},e_{19},e_{21},e_{24},e_{26},e_{29}\}$;
$w_{12}(e_{21}) = \{e_3,e_7,e_{13},e_{15},e_{16},e_{17},e_{18},e_{19},e_{20},e_{23},e_{26},e_{27},e_{28},e_{29}\}$;
$w_{12}(e_{22}) = \{e_3,e_4,e_5,e_6,e_7,e_8,e_{10},e_{11},e_{12},e_{13},e_{14},e_{16},e_{18},e_{19},e_{27},e_{28}\}$;
$w_{12}(e_{23}) = \{e_4,e_{18},e_{19},e_{21},e_{24},e_{25},e_{26},e_{27},e_{28}\}$;
$w_{12}(e_{24}) = \{e_4,e_5,e_6,e_8,e_{10},e_{11},e_{12},e_{14},e_{15},e_{17},e_{20},e_{23},e_{26},e_{29}\}$;
$w_{12}(e_{25}) = \varnothing$;
$w_{12}(e_{26}) = \{e_5,e_6,e_8,e_{10},e_{11},e_{12},e_{14},e_{15},e_{17},e_{18},e_{19},e_{20},e_{21},e_{23},e_{24},e_{25},e_{27},e_{28},e_{29}\}$;
$w_{12}(e_{27}) = \{e_3,e_7,e_8,e_{10},e_{11},e_{12},e_{15},e_{16},e_{17},e_{21},e_{22},e_{23},e_{25},e_{26}\}$;
$w_{12}(e_{28}) = \{e_3,e_7,e_8,e_{10},e_{11},e_{12},e_{15},e_{16},e_{17},e_{21},e_{22},e_{23},e_{25},e_{26}\}$;
$w_{12}(e_{29}) = \{e_4,e_5,e_6,e_{10},e_{11},e_{12},e_{14},e_{15},e_{17},e_{18},e_{19},e_{20},e_{21},e_{24},e_{25},e_{26}\}$;
$w_{12}(e_{30}) = \varnothing$.

Реберные разрезы 14-го уровня:

$w_{13}(e_1) = \varnothing$;
$w_{13}(e_2) = \varnothing$;
$w_{13}(e_3) = \{e_1,e_2,e_3,e_5,e_6,e_{13},e_{14},e_{15},e_{18},e_{19},e_{23},e_{25},e_{26},e_{27},e_{28}\}$;
$w_{13}(e_4) = \varnothing$;
$w_{13}(e_5) = \{e_1,e_2,e_3,e_{10},e_{11},e_{13},e_{14},e_{15},e_{17},e_{20},e_{21},e_{24},e_{25},e_{26},e_{27},e_{28}\}$;
$w_{13}(e_6) = \{e_1,e_2,e_3,e_{10},e_{11},e_{13},e_{14},e_{15},e_{17},e_{20},e_{21},e_{24},e_{25},e_{26},e_{27},e_{28}\}$;
$w_{13}(e_7) = \{e_2,e_8,e_{10},e_{11},e_{13},e_{15},e_{16},e_{18},e_{19},e_{25},e_{29}\}$;
$w_{13}(e_8) = \varnothing$;
$w_{13}(e_9) = \varnothing$;
$w_{13}(e_{10}) = \{e_1,e_4,e_5,e_6,e_7,e_{10},e_{11},e_{12},e_{13},e_{17},e_{20},e_{21},e_{24},e_{26},e_{29}\}$;
$w_{13}(e_{11}) = \{e_1,e_4,e_5,e_6,e_7,e_{10},e_{11},e_{12},e_{13},e_{17},e_{20},e_{21},e_{24},e_{26},e_{29}\}$;
$w_{13}(e_{12}) = \{e_2,e_{10},e_{11},e_{12},e_{13},e_{14},e_{16},e_{18},e_{19},e_{27},e_{28}\}$;
$w_{13}(e_{13}) = \{e_3,e_4,e_5,e_6,e_7,e_8,e_{10},e_{11},e_{12},e_{13},e_{14},e_{16},e_{18},e_{19},e_{27},e_{28}\}$;
$w_{13}(e_{14}) = \{e_1,e_2,e_3,e_4,e_5,e_6,e_{12},e_{13},e_{18},e_{19},e_{27},e_{28}\}$;
$w_{13}(e_{15}) = \varnothing$;
$w_{13}(e_{16}) = \{e_1,e_7,e_8,e_{12},e_{13},e_{18},e_{19},e_{27},e_{28}\}$;
$w_{13}(e_{17}) = \{e_5,e_6,e_{10},e_{11},e_{17},e_{18},e_{19},e_{20},e_{21},e_{23},e_{24}\}$;
$w_{13}(e_{18}) = \{e_3,e_4,e_7,e_{12},e_{13},e_{14},e_{16},e_{17},e_{18},e_{19},e_{20},e_{25},e_{27},e_{28},e_{29}\}$;
$w_{13}(e_{19}) = \{e_3,e_4,e_7,e_{12},e_{13},e_{14},e_{16},e_{17},e_{18},e_{19},e_{20},e_{25},e_{27},e_{28},e_{29}\}$;
$w_{13}(e_{20}) = \{e_4,e_5,e_6,e_8,e_{10},e_{11},e_{17},e_{18},e_{19},e_{20},e_{21},e_{24},e_{26},e_{27},e_{28},e_{29}\}$;
$w_{13}(e_{21}) = \{e_4,e_5,e_6,e_{10},e_{11},e_{17},e_{20},e_{23},e_{25},e_{26},e_{27},e_{28}\}$;
$w_{13}(e_{22}) = \varnothing$;
$w_{13}(e_{23}) = \varnothing$;
$w_{13}(e_{24}) = \{e_5,e_6,e_8,e_{10},e_{11},e_{17},e_{20},e_{25},e_{29}\}$;
$w_{13}(e_{25}) = \varnothing$;
$w_{13}(e_{26}) = \{e_1,e_2,e_3,e_5,e_6,e_{10},e_{11},e_{20},e_{21},e_{22},e_{23},e_{25},e_{26},e_{27},e_{28}\}$;
$w_{13}(e_{27}) = \{e_1,e_3,e_5,e_6,e_{12},e_{13},e_{14},e_{16},e_{18},e_{19},e_{20},e_{21},e_{22},e_{23},e_{25},e_{26}\}$;
$w_{13}(e_{28}) = \{e_1,e_3,e_5,e_6,e_{12},e_{13},e_{14},e_{16},e_{18},e_{19},e_{20},e_{21},e_{22},e_{23},e_{25},e_{26}\}$;
$w_{13}(e_{29}) = \{e_1,e_7,e_8,e_{10},e_{11},e_{18},e_{19},e_{20},e_{22},e_{23},e_{24}\}$;
$w_{13}(e_{30}) = \varnothing$.

Для графа $G_{15}$ количество уровней в спектре реберных разрезов = 6, а для графа $G_{16}$ это 14-ть уровней. Графы не изоморфны, так как отличаются количеством уровней в спектре реберных разрезов.

Таким образом, алгоритм, основанный на свойствах реберных разрезов, сумел определить неизоморфность графов и в случае, когда существует равенство весов инварианта реберных циклов. Кроме того, перестановка всего двух ребер приводит к существенному изменению количества уровней в спектре реберных разрезов. Интересным представляется вопрос расхождения реберных разрезов в зависимости от уровня для неизоморфных графов. Для сравнения выберем реберные разрезы для ребра $e_3$ графа $G_{15}$ и разрезы для ребра $e_{16}$ графа $G_{16}$ (рис. 8.13 – 8.16). Рассматривая рисунки можно сказать, что расхождение для ребра $e_3$ графа $G_{15}$ и ребра $e_{16}$ графа $G_{16}$ происходит на 3-уровне.

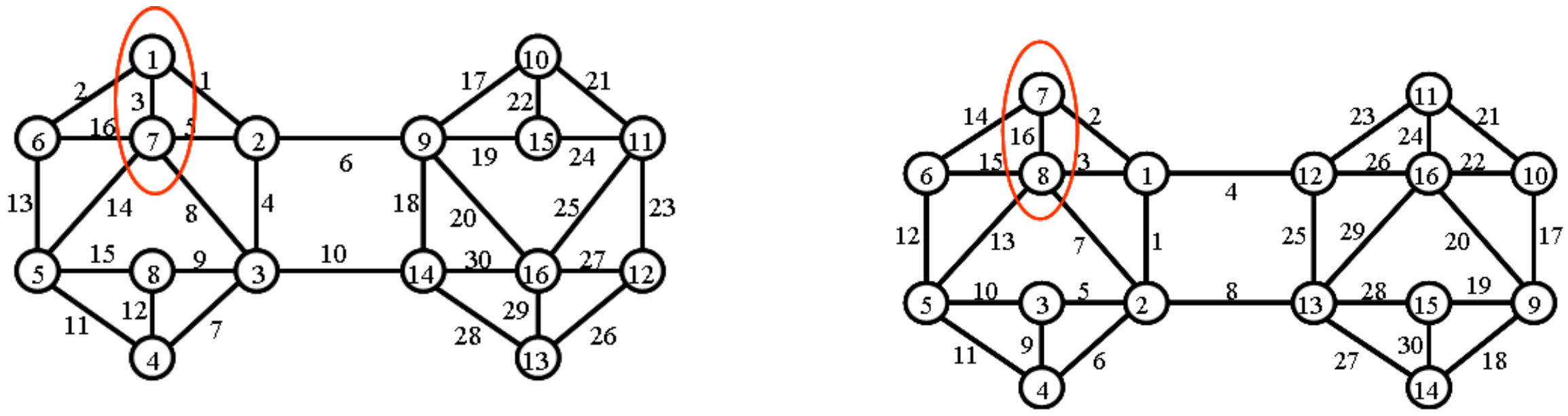

Рис. 8.13. Реберные разрезы для 1-го уровня.

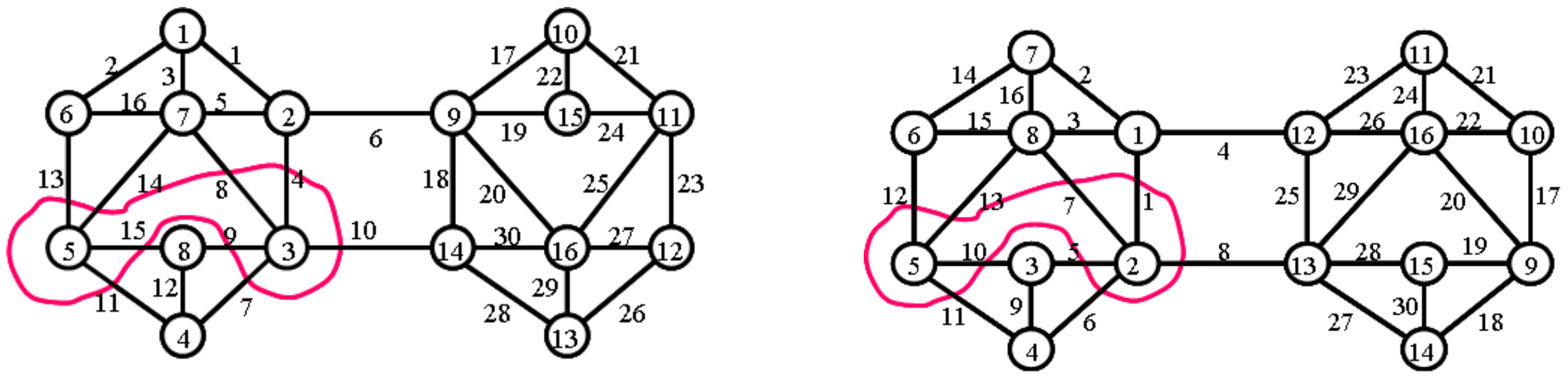

Рис. 8.14. Реберные разрезы для 2-го уровня.

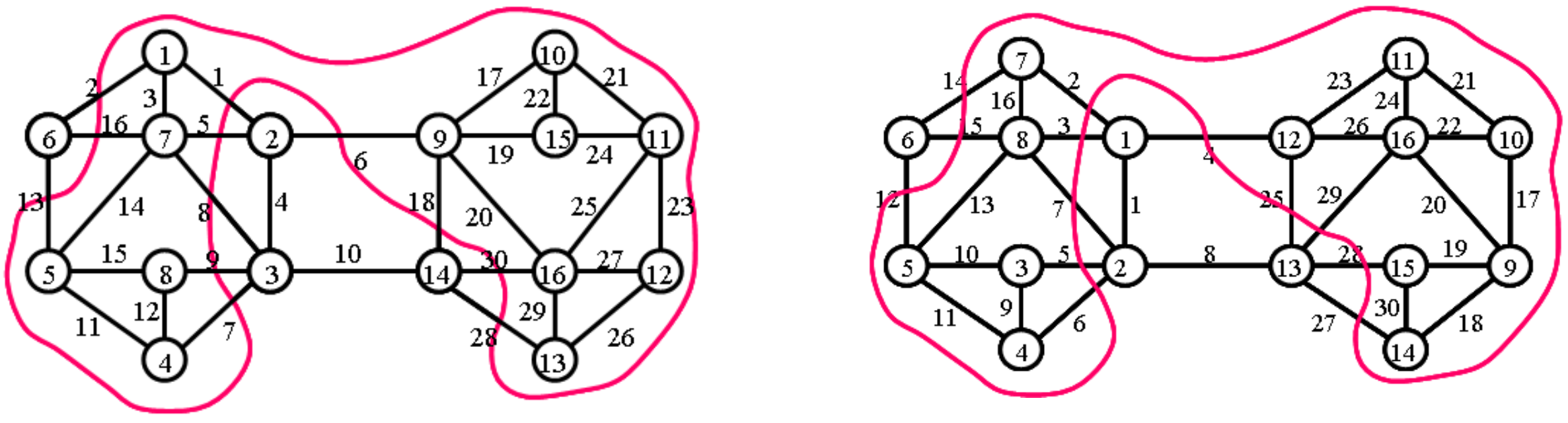

Рис. 8.15. Реберные разрезы для 3-го уровня.

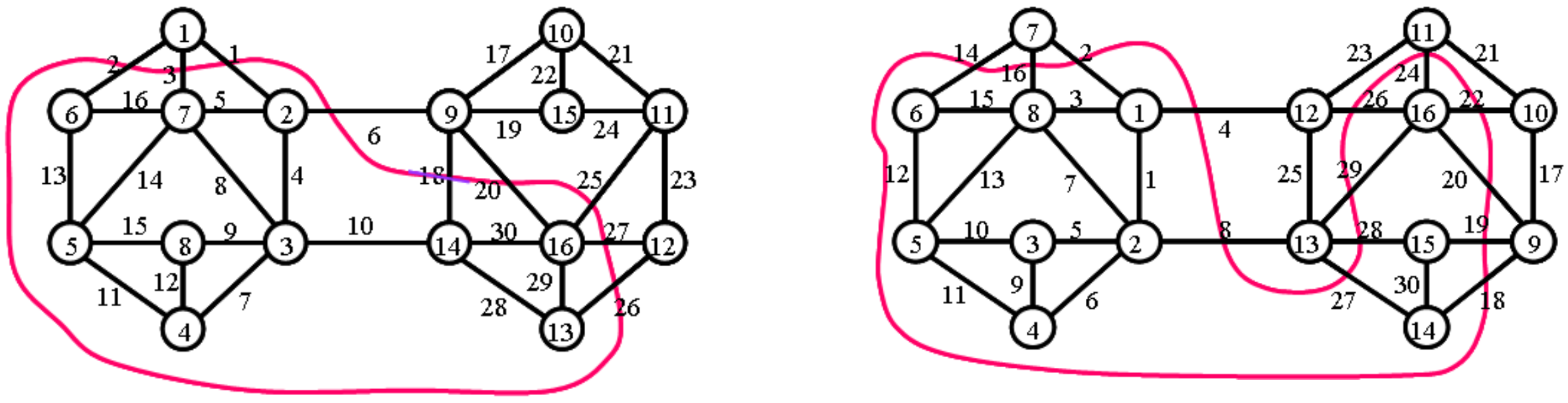

Рис. 8.16. Реберные разрезы для 4-го уровня.

***Пример 8.6.*** Рассмотрим формирование инварианта спектра реберных разрезов для плоских графов $G_{37}$ и $G_{38}$ (рис. 8.17 и 8.18).

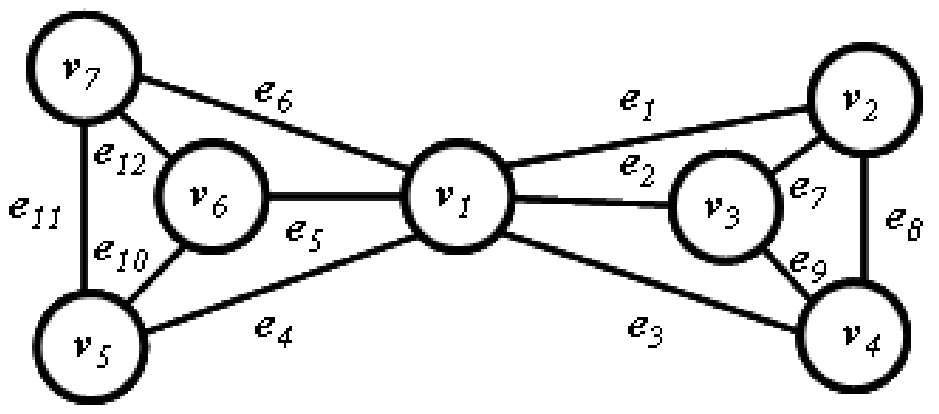

Рис. 8.17. Граф $G_{37}$.

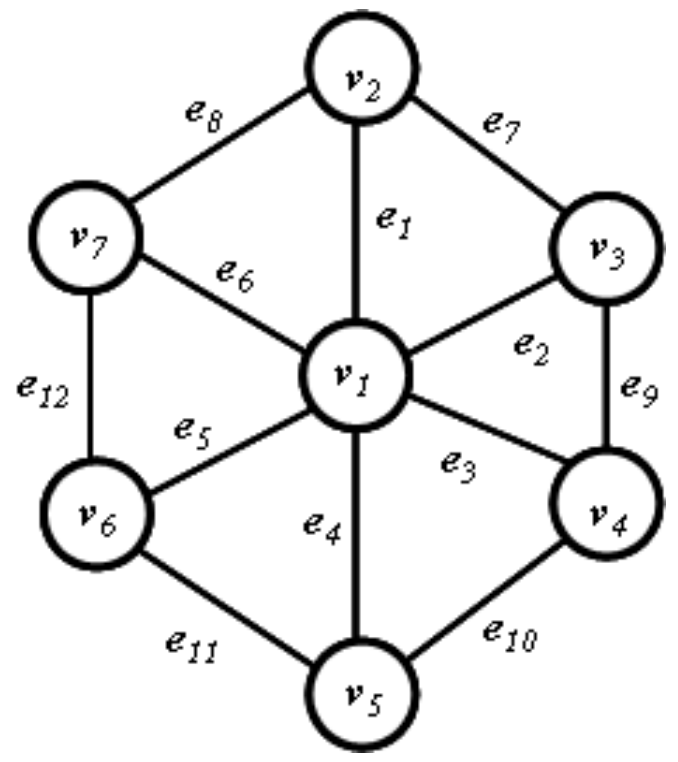

Рис. 8.18. Граф $G_{38}$.

Инвариант реберных разрезов графа $G_{37}$: $F_w(\xi(G_{37}))$ & $F_w(\zeta(G_{37})) = (6\times 4,6\times 13)$ &

& $(1\times 78,6\times 21)$ ,

Инвариант реберных разрезов графа $G_{38}$: $F_w(\xi(G_{38}))$ & $F_w(\zeta(G_{38})) = (6\times 22,6\times 24)$ &

& $(6\times 68,1\times 144\}$.

Графы не изоморфны.

## 8.5. Деревья графа

***Пример 8.7.*** Рассмотрим формирование инварианта спектра реберных разрезов для дерева графа.

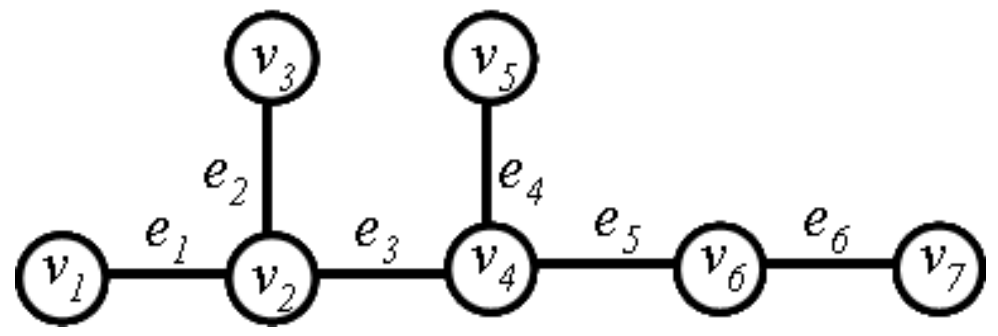

Рис. 8.19. Дерево графа $G_{39}$.

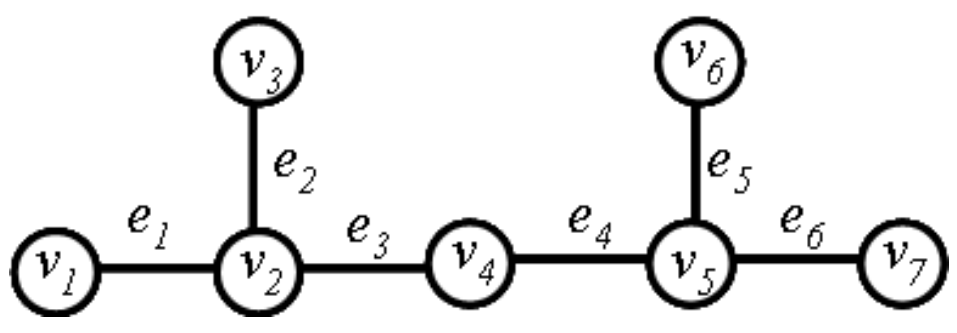

Рис. 8.20. Дерево графа $G_{40}$.

Кортеж весов ребер : $\xi_w(G_{39})$ = <10,10,12,12,10,10>;

Кортеж весов вершин : $\zeta_w(G_{39}) = <10,32,10,24,32,10,10>$.

Ивариант графа $G_{13}$: $F_w(\xi(G_{39}))$ & $F_w(\zeta(G_{39})) = <4\times 10,2\times 12>$ & $<4\times 10,1\times 24,2\times 32>$.

Кортеж весов ребер : $\xi_w(G_{40}) = <18,18,20,20,16,16>$;

Кортеж весов вершин : $\zeta_w(G_{40}) = <18,56,18,56,20,32,16>$.

Инвариант графа $G_{40}$: $F_w(\xi(G_{40}))$ & $F_w(\zeta(G_{40})) = <2\times 16,2\times 18,2\times 20>$ &

& $<1\times 16,2\times 18,1\times 20,1\times 32,2\times 56>$.

Графы не изоморфны.

## 8.6. Сравнительный анализ графов

***Пример 8.8.*** Рассмотрим следующую совокупность графов на 6 вершин.

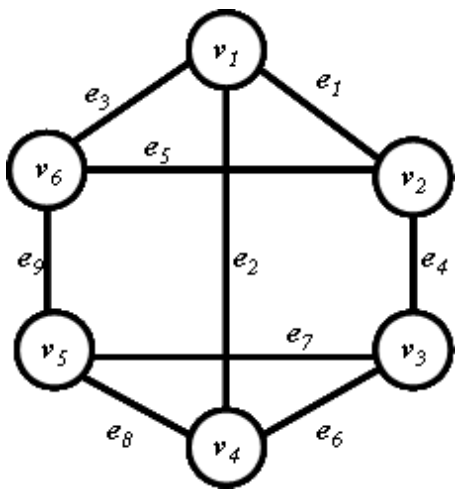

Рис. 8.21. Граф $G_{41}$.

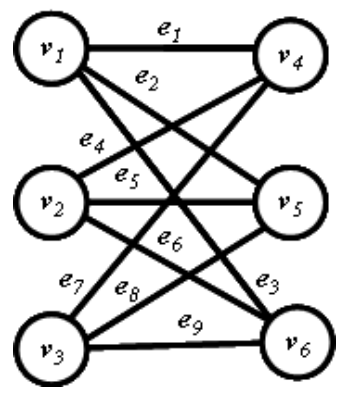

Рис. 8.22. Граф $G_{42}$.

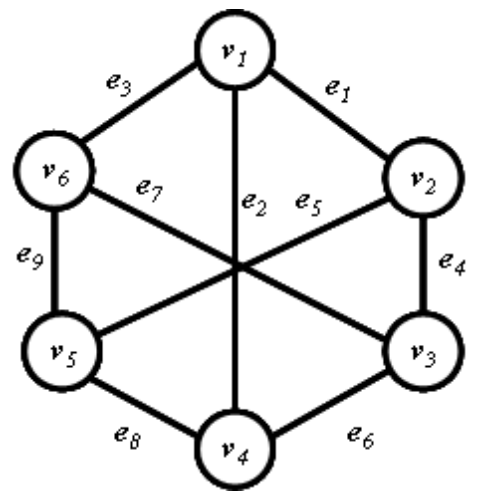

Рис. 8.23. Граф $G_{43}$.

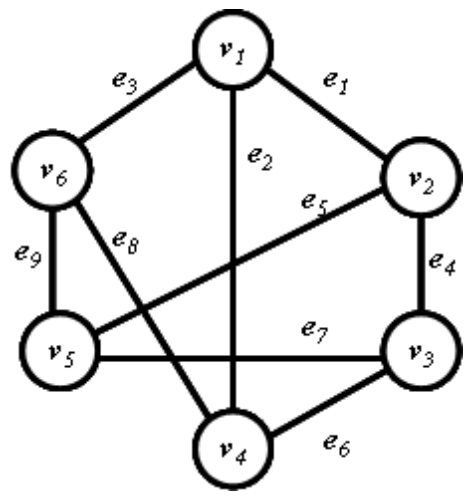

Рис. 8.24. Граф $G_{44}$.

a). Инварианты графа $G_{41}$

Количество вершин в графе = 6;

Количество ребер в графе = 9;

Спектр реберных разрезов:

Кортеж весов ребер уровня 1: <4,4,4,4,4,4,4,4,4>;

Кортеж весов ребер уровня 2: <2,4,2,4,2,2,2,2,4>;

Суммарный кортеж весов ребер : <6,8,6,8,6,6,6,6,8>;

Кортеж весов вершин : <20,20,20,20,20,20>;

Спектр реберных циклов:

Кортеж весов ребер: <5,6,5,6,5,5,5,5,6>.

Кортеж весов вершин : <16,16,16,16,16,16>;

Интегральный инвариант графа

$F_w(\xi(G_{41}))$ & $F_w(\zeta(G_{41}))$ & $F_\tau(\xi(G_{41}))$ & $F\tau(\zeta(G_{41})) =$

$= (6\times 6,3\times 8)$ & $(9\times 20)$ & $(6\times 5,3\times 6)$ & $(9\times 16)$.

b). Инварианты графа $G_{42}$:

Количество вершин в графе = 6;

Количество ребер в графе = 9;

Спектр реберных разрезов:

Кортеж весов ребер уровня 1: <4,4,4,4,4,4,4,4,4>;

Кортеж весов ребер : <4,4,4,4,4,4,4,4,4>;

Кортеж весов вершин : <12,12,12,12,12,12>;

Спектр реберных циклов:

Кортеж весов ребер : <4,4,4,4,4,4,4,4,4>;

Кортеж весов вершин : <12,12,12,12,12,12>;

Интегральный инвариант графа $F_w(G_{42})$ & $F_w(G_{42})$ & $F_\tau(\xi(G_{42}))$ & $F_\tau(\zeta(G_{42}))$ =

= $(9\times4)$ & $(9\times12)$ & $(9\times4)$ & $(9\times12)$.

c). Инварианты графа $G_{43}$:

Количество вершин в графе = 6;

Количество ребер в графе = 9;

Спектр реберных разрезов:

Кортеж весов ребер уровня 1: <4,4,4,4,4,4,4,4,4>;

Кортеж весов ребер : <4,4,4,4,4,4,4,4,4>;

Кортеж весов вершин : <12,12,12,12,12,12>;

Спектр реберных циклов:

Кортеж весов ребер : <4,4,4,4,4,4,4,4,4>;

Кортеж весов вершин : <12,12,12,12,12,12>;

Интегральный инвариант графа $F_w(G_{43})$ & $F_w(G_{43})$ & $F_\tau(\xi(G_{43}))$ & $F_\tau(\zeta(G_{43}))$ =

= $(9\times4)$ & $(9\times12)$ & $(9\times4)$ & $(9\times12)$.

d). Инварианты графа $G_{44}$:

Количество вершин в графе = 6;

Количество ребер в графе = 9;

Спектр реберных разрезов:

Кортеж весов ребер уровня 1: <4,4,4,4,4,4,4,4,4>;

Кортеж весов ребер уровня 2: <4,2,2,2,2,4,2,2,4>;

Суммарный кортеж весов: <8,6,6,6,6,8,6,6,8>

Кортеж весов вершин : <20,20,20,20,20,20>;

Спектр реберных циклов:

Кортеж весов ребер: <6,5,5,5,5,6,5,5,6>.

Кортеж весов вершин : <16,16,16,16,16,16>;

Интегральный инвариант графа $F_w(G_{44})$ & $F_w(G_{44})$ & $F_\tau(\xi(G_{44}))$ & $F_\tau(\zeta(G_{44}))$ =

$= (6\times 6{,}3\times 8)$ & $(9\times 20)$ & $(6\times 5{,}3\times 6)$ & $(9\times 16)$

***Вывод:*** Графы $G_{42}$,$G_{43}$ и графы $G_{41}$,$G_{44}$ –изоморфны между собой.

***Пример 8.9.*** Рассмотрим следующую совокупность графов на 8 вершин (рис. 8.25-8.35).

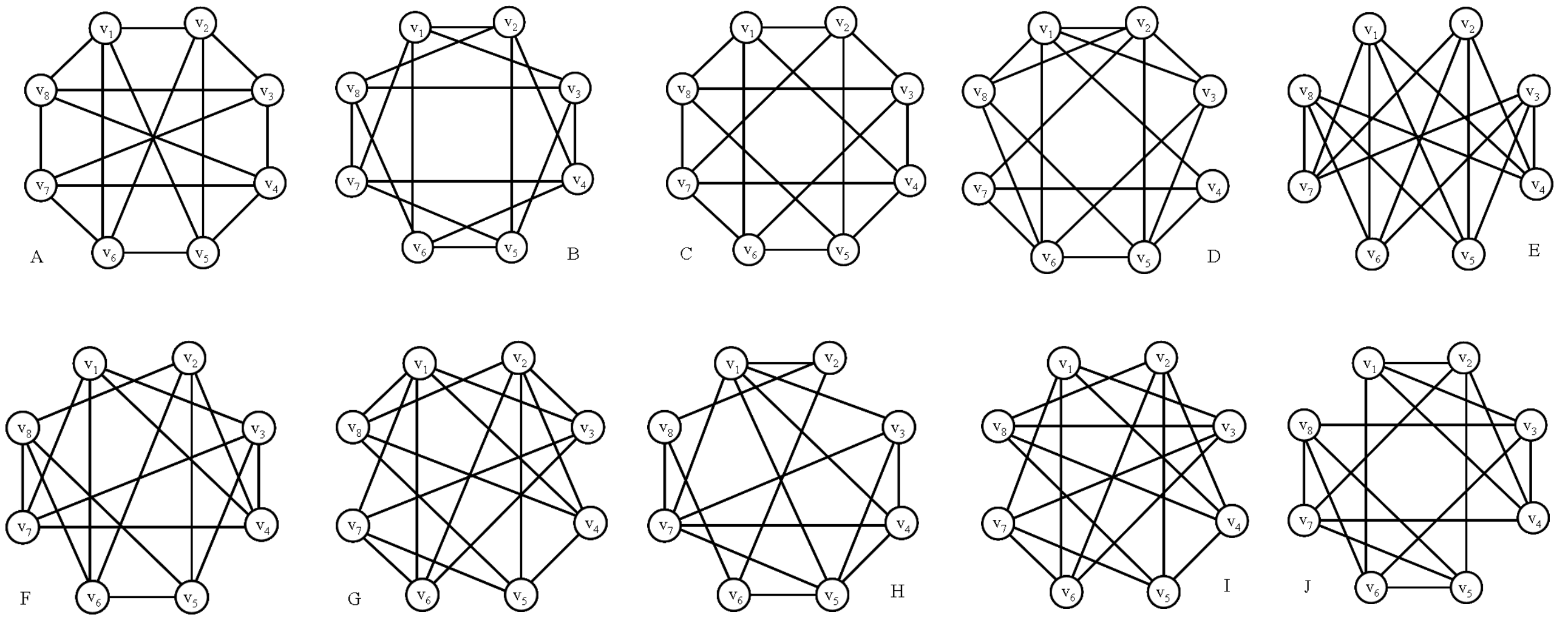

Рис. 8.25. Совокупность 8-ми вершинных графов.

а). Граф $G_{45}$.

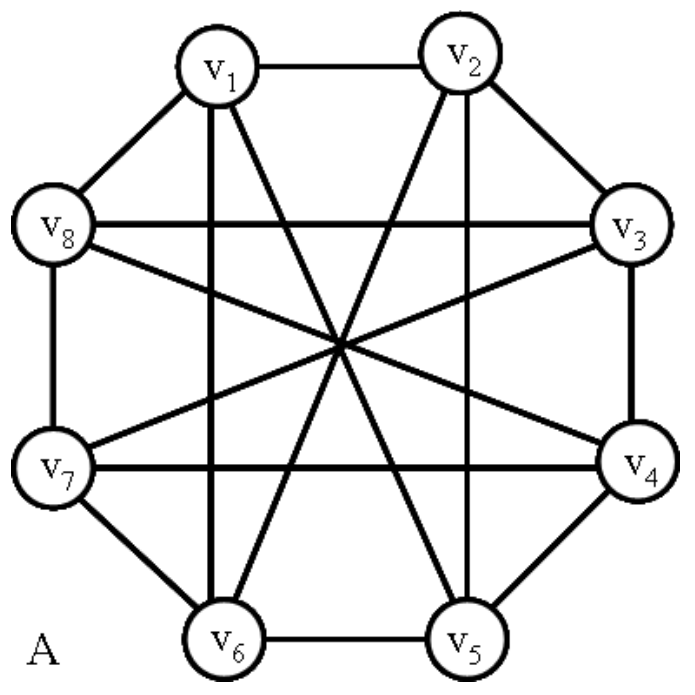

Рис. 8.26. Граф $G_{45}$.

Количество вершин в графе = 8;

Количество ребер в графе = 16;

Спектр реберных разрезов:

Кортеж весов ребер уровня 1: <6,6,6,6,6,6,6,6,6,6,6,6,6,6,6,6>;

Кортеж весов ребер уровня 2: <8,8,8,0,0,8,8,8,8,8,0,8,8,8,0,8>;

Суммарный кортеж весов ребер : <14,14,14,6,6,14,14,14,14,14,6,14,14,14,6,14>;

Кортеж весов вершин : <48,48,48,48,48,48,48,48>;

Спектр реберных циклов:

Кортеж весов ребер : <8,8,8,6,6,8,8,8,8,8,6,8,8,8,6,8>

Кортеж весов вершин : <30,30,30,30,30,30,30,30>;

Интегральный инвариант графа $F_w(\xi(G_{45}))$ & $F_w(\zeta(G_{45}))$ & $F_\tau(\xi(G_{45}))$ & $F_\tau(\zeta(G_{45}))$ =

= $(4\times 6, 12\times 14)$ & $(8\times 48)$ & $(4\times 6, 12\times 8)$ & $(8\times 30)$.

b). Граф $G_{46}$.

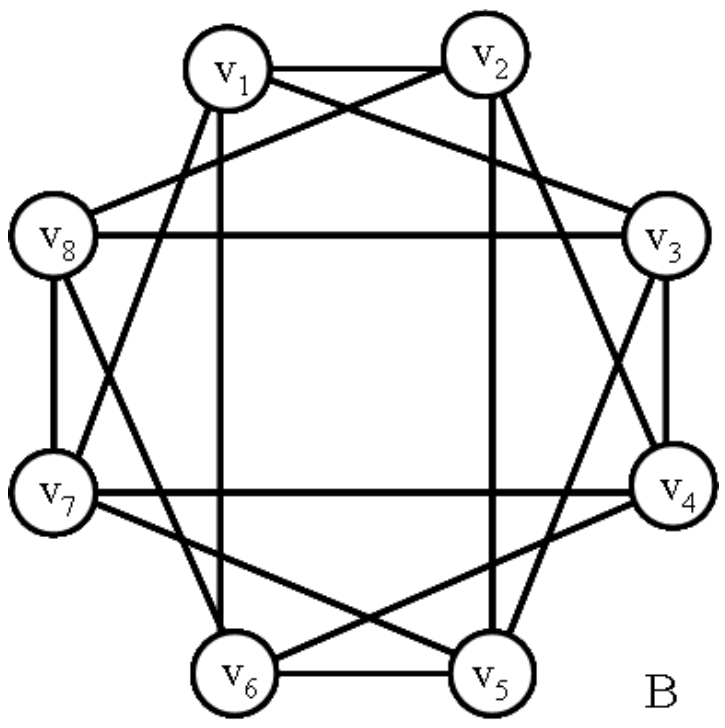

Рис. 8.27. Граф $G_{46}$.

Количество вершин в графе = 8;

Количество ребер в графе = 16;

Спектр реберных разрезов:

Кортеж весов ребер уровня 1: <6,6,6,6,6,6,6,6,6,6,6,6,6,6,6,6>;

Кортеж весов ребер : <6,6,6,6,6,6,6,6,6,6,6,6,6,6,6,6>;

Кортеж весов вершин : <24,24,24,24,24,24,24,24>;

Спектр реберных циклов:

Кортеж весов ребер : <16,16,16,16,16,16,16,16,16,16,16,16,16,16,16,16>

Кортеж весов вершин : <64,64,64,64,64,64,64,64>;

Основной базовый инвариант графа

$F_w(\xi(G_{46}))$ & $F_w(\zeta(G_{46}))$ & $F_\tau(\xi(G_{46}))$ & $F_\tau(\zeta(G_{46}))$ =

= $(16\times 6)$ & $(8\times 24)$ & $(16\times 16)$ & $(8\times 64)$.

c). Граф $G_{47}$.

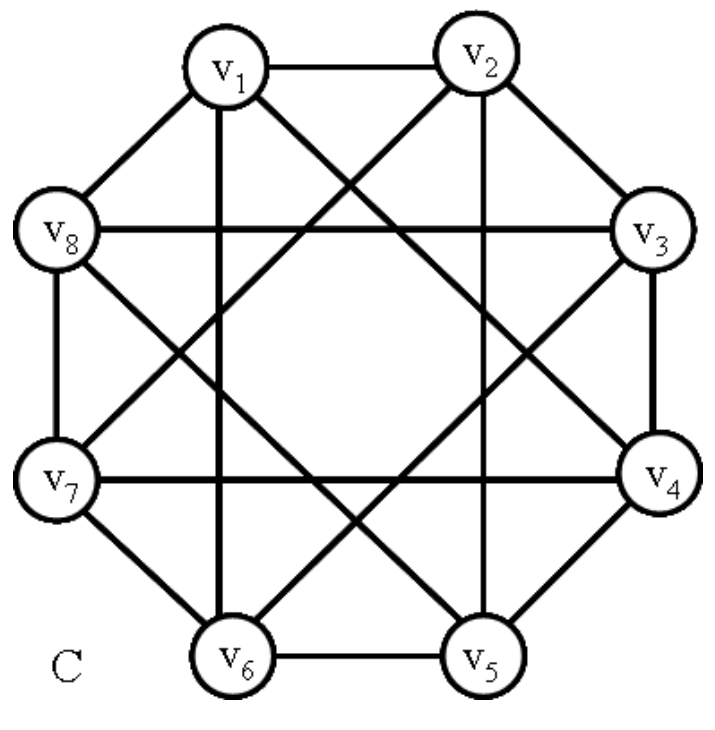

Рис. 8.28. Граф $G_{47}$.

Количество вершин в графе = 8;

Количество ребер в графе = 16;

Спектр реберных разрезов:

Кортеж весов ребер уровня 1: <6,6,6,6,6,6,6,6,6,6,6,6,6,6,6,6>;

Кортеж весов ребер : <6,6,6,6,6,6,6,6,6,6,6,6,6,6,6,6>;

Кортеж весов вершин : <24,24,24,24,24,24,24,24>;

Спектр реберных циклов:

Кортеж весов ребер : <16,16,16,16,16,16,16,16,16,16,16,4,16,16,16,16>

Кортеж весов вершин : <64,64,64,64,64,64,64,64>;

Интегральный инвариант графа

$F_w(\xi(G_{47}))$ & $F_w(\zeta(G_{47}))$ & $F_\tau(\xi(G_{47}))$ & $F_\tau(\zeta(G_{47}))$ = $(16\times 6)$ & $(8\times 24)$ &

& $(16\times 16)$ & $(8\times 64)$.

d). Граф $G_{48}$.

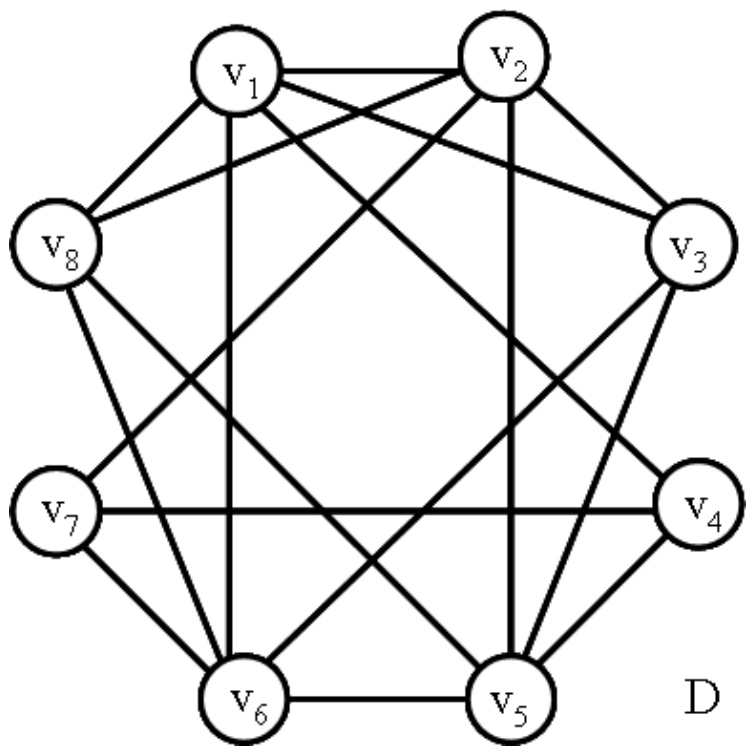

Рис. 8.29. Граф $G_{48}$.

Количество вершин в графе = 8;

Количество ребер в графе = 17;

Спектр реберных разрезов:

Кортеж весов ребер уровня 1: <8,7,6,8,7,7,8,6,7,7,7,6,4,8,7,6,7>;

Кортеж весов ребер уровня 2: <6,9,5,6,9,9,6,5,9,9,9,5,8,6,9,5,9>;

Кортеж весов ребер уровня 3: <6,5,5,6,5,5,6,5,5,5,5,5,8,6,5,5,5>;

Суммарный кортеж весов ребер : <20,21,16,20,21,21,20,16,21,21,21,16,20,20,21,16,21>;

Кортеж весов вершин : <98,98,84,52,98,98,52,84>;

Спектр реберных циклов:

Кортеж весов ребер : <8,8,8,8,8,8,8,8,8,8,8,8,4,8,8,8,8>

Кортеж весов вершин : <40,40,32,20,40,40,20,32>;

Интегральный инвариант графа $F_w(\xi(G_{48}))$ & $F_w(\zeta(G_{48}))$ & $F_\tau(\xi(G_{48}))$ & $F_\tau(\zeta(G_{48}))$ =

= $(4\times 16, 4\times 20, 8\times 21)$ & $(2\times 52, 2\times 84, 4\times 98)$ & $(1\times 4, 16\times 8)$ & $(2\times 20, 2\times 32, 4\times 40)$.

e). Граф $G_{49}$.

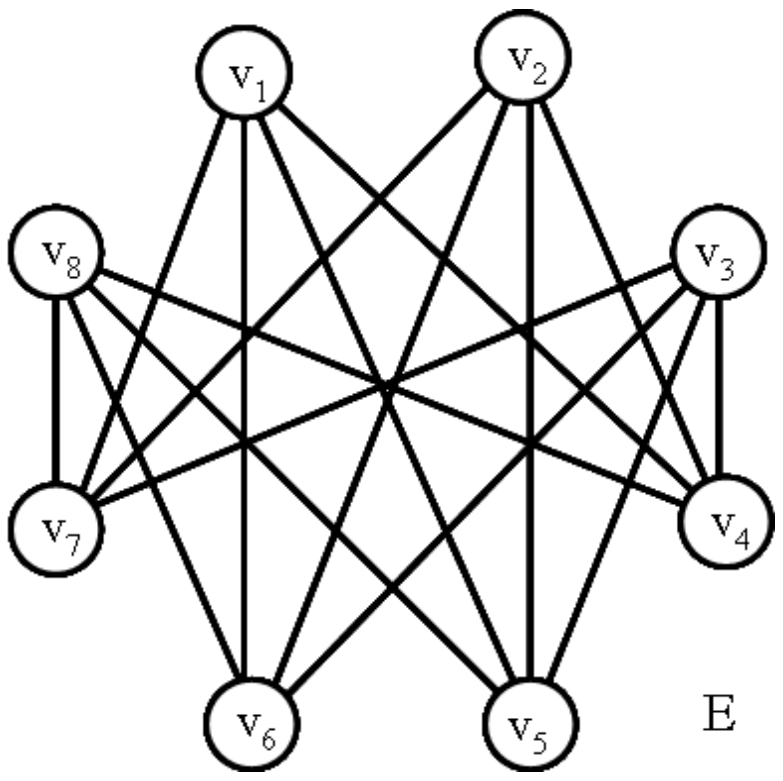

Рис. 8.30. Граф $G_{49}$.

Количество вершин в графе = 8;

Количество ребер в графе = 16;

Спектр реберных разрезов:

Кортеж весов ребер уровня 1: <6,6,6,6,6,6,6,6,6,6,6,6,6,6,6,6>;

Кортеж весов ребер : <6,6,6,6,6,6,6,6,6,6,6,6,6,6,6,6>;

Кортеж весов вершин : <24,24,24,24,24,24,24,24>;

Спектр реберных циклов:

Кортеж весов ребер : <16,16,16,16,16,16,16,16,16,16,16,4,16,16,16,16>

Кортеж весов вершин : <64,64,64,64,64,64,64,64>;

Интегральный инвариант графа $F_w(\xi(G_{49}))$ & $F_w(\zeta(G_{49}))$ & $F_\tau(\xi(G_{49}))$ & $F_\tau(\zeta(G_{49}))$ =

= $(16\times 6)$ & $(8\times 24)$ & $(16\times 16)$ & $(8\times 64)$.

f). Граф $G_{50}$.

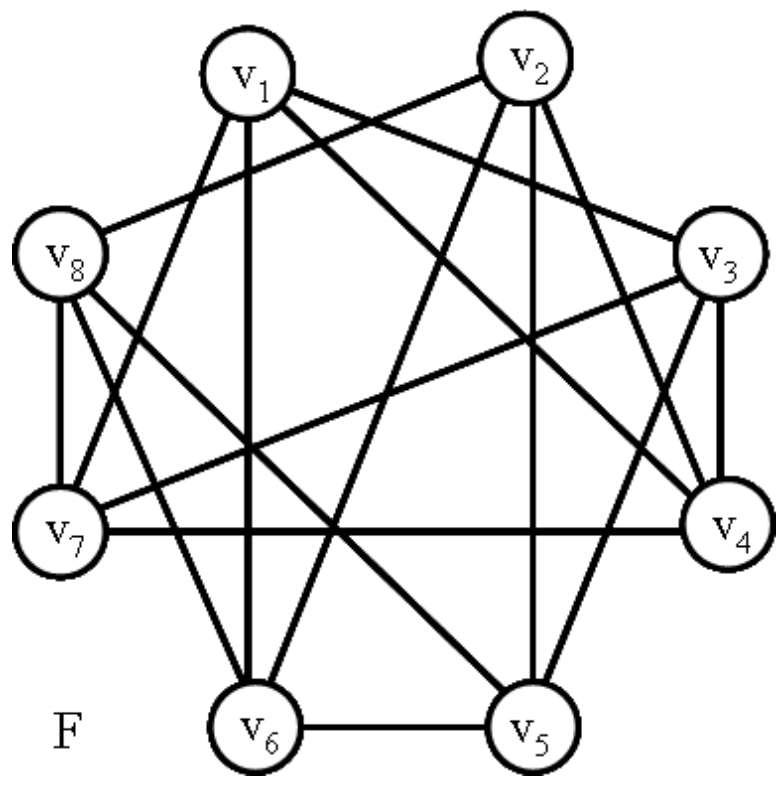

Рис. 8.31. Граф $G_{50}$.

Количество вершин в графе = 8;

Количество ребер в графе = 16;

Спектр реберных разрезов:

Кортеж весов ребер уровня 1: <6,6,6,6,6,6,6,6,6,6,6,6,6,6,6,6>;

Кортеж весов ребер уровня 2: <8,8,0,8,0,8,8,8,8,0,8,8,8,8,8,0>;

Суммарный кортеж весов ребер : <14,14,6,14,6,14,14,14,14,6,14,14,14,14,14,6>;

Кортеж весов вершин : <48,48,48,48,48,48,48,48>;

Спектр реберных циклов:

Кортеж весов ребер : <8,8,6,8,6,8,8,8,8,6,8,8,8,8,8,6>

Кортеж весов вершин : <30,30,30,30,30,30,30,30>;

Интегральный инвариант графа $F_w(\xi(G_{50}))$ & $F_w(\zeta(G_{50}))$ & $F_\tau(\xi(G_{50}))$ & $F_\tau(\zeta(G_{50}))$ =

$= (4\times 6, 12\times 14)$ & $(8\times 48)$ & $(4\times 6, 12\times 8)$ & $(8\times 30)$.

g). Граф $G_{51}$.

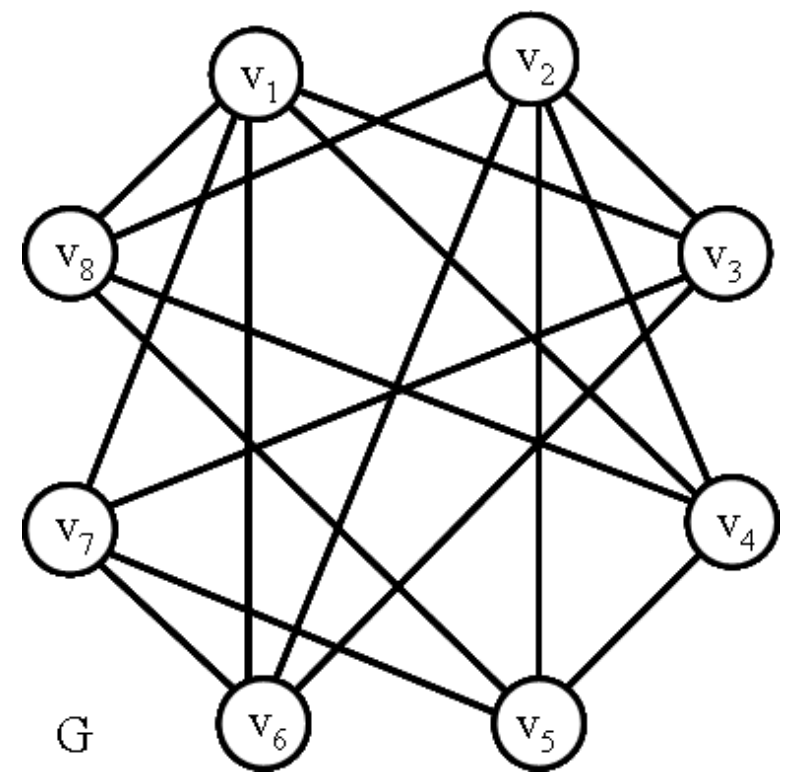

Рис. 8.32. Граф $G_{51}$.

Количество вершин в графе = 8;

Количество ребер в графе = 17;

Спектр реберных разрезов:

Кортеж весов ребер уровня 1: <7,7,7,7,7,7,7,7,7,7,6,6,6,6,6,6,6>;

Кортеж весов ребер уровня 2: <8,6,8,8,6,6,8,8,6,8,6,8,8,6,0,8,8>;

Кортеж весов ребер уровня 3: <5,5,5,4,5,5,5,4,5,5,4,5,5,4,0,5,5>;

Кортеж весов ребер уровня 4: <5,5,5,4,5,5,5,4,5,5,4,5,5,4,0,5,5>;

Суммарный кортеж весов ребер : <25,23,25,23,23,23,25,23,23,25,20,24,24,20,6,24,24>;

Кортеж весов вершин : <119,119,92,92,77,92,77,92>;

Спектр реберных циклов:

Кортеж весов ребер : <8,9,8,8,9,9,8,8,9,8,7,8,8,7,8,8,8>

Кортеж весов вершин : <42,42,32,32,32,32,32,32>;

Интегральный инвариант графа $F_w(\xi(G_{51}))$ & $F_w(\zeta(G_{51}))$ & $F_\tau(\xi(G_{51}))$ & $F_\tau(\zeta(G_{51}))$ –

$= (1\times 6, 2\times 20, 6\times 23, 4\times 24, 4\times 25)$ & $(2\times 77, 4\times 92, 2\times 119)$ & $(2\times 7, 11\times 8, 4\times 9)$ &

& $(6\times 32, 2\times 42)$.

h). Граф $G_{52}$.

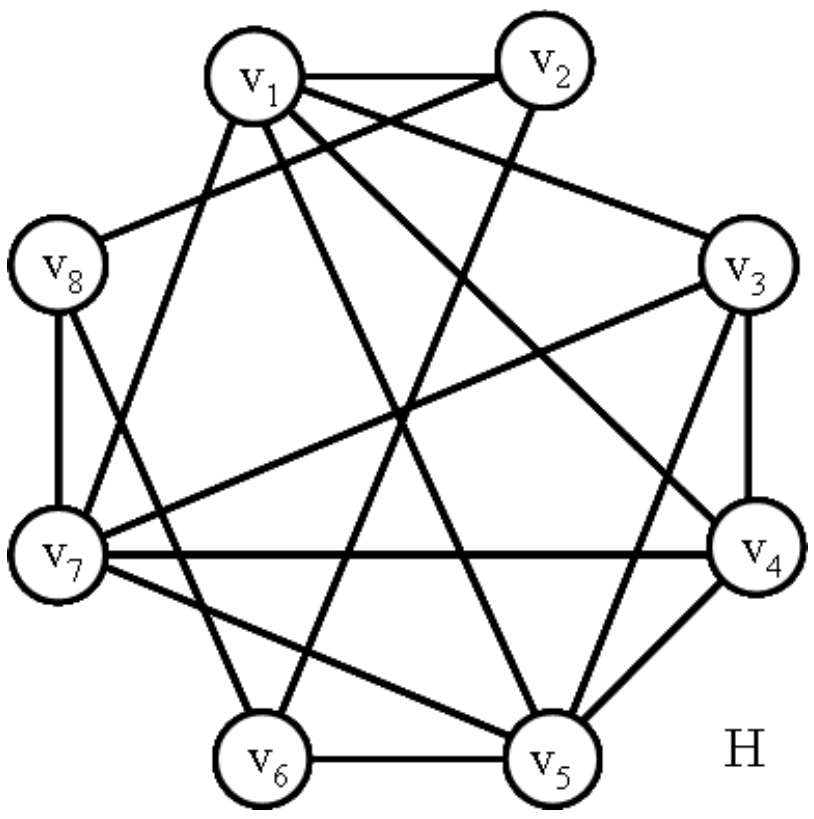

Рис. 8.33. Граф $G_{52}$.

Количество вершин в графе = 8;

Количество ребер в графе = 17;

Спектр реберных разрезов:

Кортеж весов ребер уровня 1: <6,7,7,8,8,4,4,6,7,7,7,7,6,8,4,6>;

Кортеж весов ребер уровня 2: <6,5,5,2,2,6,6,6,5,5,5,5,6,2,6,6>;

Суммарный кортеж весов ребер : <12,12,12,10,10,10,10,12,12,12,12,12,12,10,10,12>;

Кортеж весов вершин : <56,32,48,48,56,32,56,32>;

Спектр реберных циклов:

Кортеж весов ребер : <6,4,4,9,9,5,5,0,4,4,4,4,6,9,5,6>;

Кортеж весов вершин : <32,16,12,12,32,16,32,16>;

Интегральный инвариант графа $F_w(\xi(G_{52}))$ & $F_w(\zeta(G_{52}))$ & $F_\tau(\xi(G_{52}))$ & $F_\tau(\zeta(G_{52}))$ =

$= (6\times 10, 10\times 12)$ & $(3\times 32, 2\times 48, 3\times 56)$ & $(1\times 0, 6\times 4, 3\times 5, 3\times 6, 4\times 9)$ & $(2\times 12, 3\times 16, 3\times 32)$.

i). Граф $G_{53}$.

Количество вершин в графе = 8;

Количество ребер в графе = 16;

Спектр реберных разрезов:

Кортеж весов ребер уровня 1: <6,6,6,6,6,6,6,6,6,6,6,6,6,6,6,6>;

Кортеж весов ребер уровня 2: <8,0,8,8,8,8,0,8,8,8,0,8,8,0,8,8>;

Суммарный кортеж весов ребер : <14,6,14,14,14,14,6,14,14,14,6,14,14,6,14,14>

Кортеж весов вершин : <48,48,48,48,48,48,48,48>;

Спектр реберных циклов:

Кортеж весов ребер : <8,6,8,8,8,8,6,8,8,8,6,8,8,6,8,8>;

Кортеж весов вершин : <30,30,30,30,30,30,30,30>;

Интегральный инвариант графа $F_w(\xi(G_{53}))$ & $F_w(\zeta(G_{53}))$ & $F_\tau(\xi(G_{53}))$ & $F_\tau(\zeta(G_{53}))$ =

= $(4\times 6, 12\times 14)$ & $(8\times 48)$ & $(4\times 6, 12\times 8)$ & $(8\times 30)$.

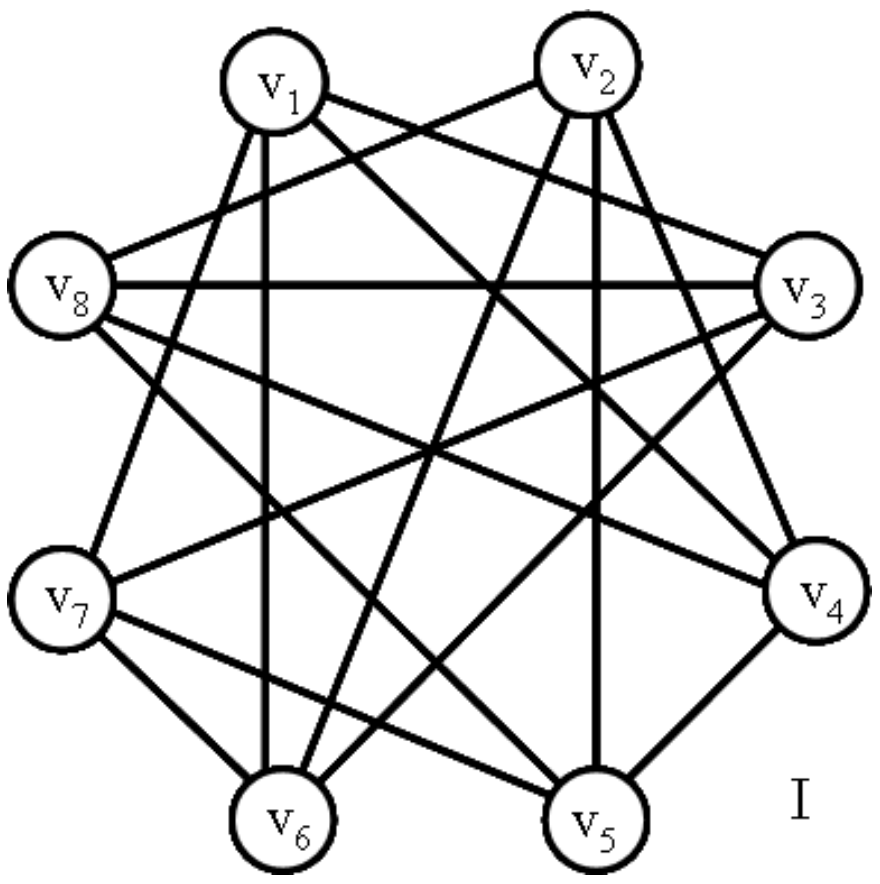

Рис. 8.34. Граф $G_{53}$.

j). Граф $G_{54}$.

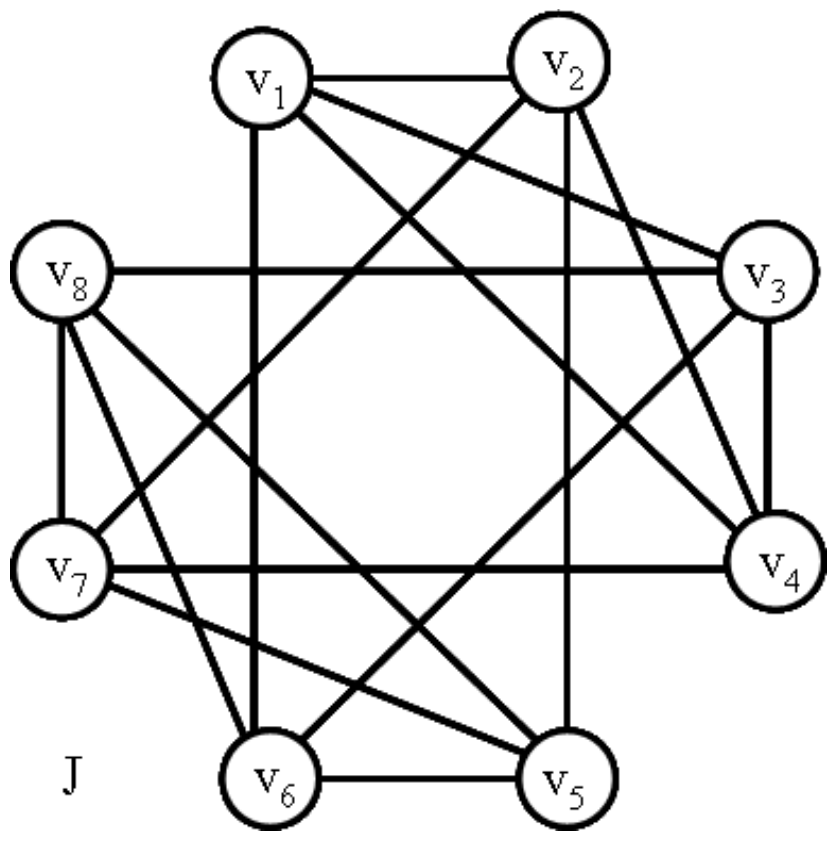

Рис. 8.35. Граф $G_{54}$.

Количество вершин в графе = 8;

Количество ребер в графе = 16;

Спектр реберных разрезов:

Кортеж весов ребер уровня 1: <6,6,6,6,6,6,6,6,6,6,6,6,6,6,6,6>

Кортеж весов ребер уровня 2: <8,10,10,8,10,8,10,8,10,8,8,8,10,10,10,8>;

Кортеж весов ребер уровня 3: <8,12,12,8,12,8,12,8,12,8,8,8,12,12,12,8>;

Суммарный кортеж весов ребер : <22,28,28,22,28,22,28,22,28,22,22,22,28,28,28,22>;

Кортеж весов вершин : <100,100,100,100,100,100,100,100>;

Спектр реберных циклов:

Кортеж весов ребер : <5,4,4,5,4,5,4,5,4,5,5,5,4,4,4,5>;

Кортеж весов вершин : <18,18,18,18,18,18,18,18>;

Интегральный инвариант графа $F_w(\xi(G_{54}))$ & $F_w(\zeta(G_{54}))$ & $F_\tau(\xi(G_{54}))$ & $F_\tau(\zeta(G_{54}))$ =

= $(8\times 22, 8\times 28)$ & $(8\times 100)$ & $(8\times 4, 8\times 5)$ & $(8\times 18)$.

**Вывод:**

Графы $G_{45}, G_{50}, G_{53}$ – изоморфны между собой,

графы $G_{46}, G_{47}, G_{49}$ – также изоморфны между собой.

***Пример 8.10***. Определим цифровой инвариант реберных графов $L(G_{29})$ и $L(G_{30})$.

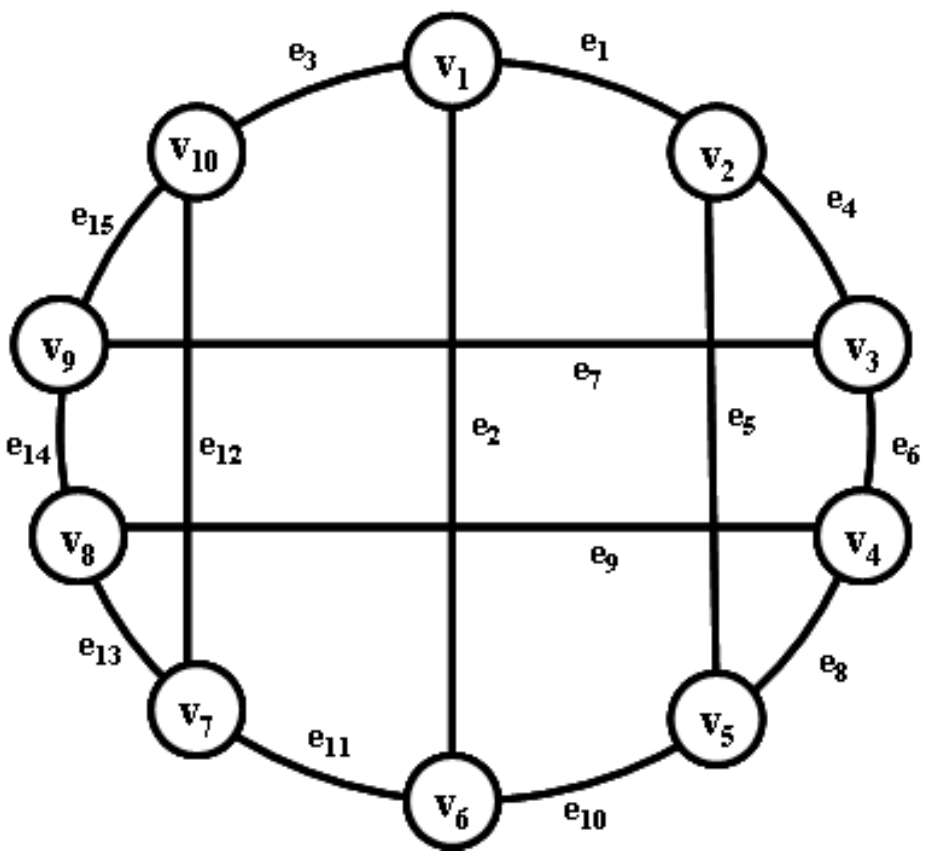

Рис. 8.36. Граф $G_{23}$.

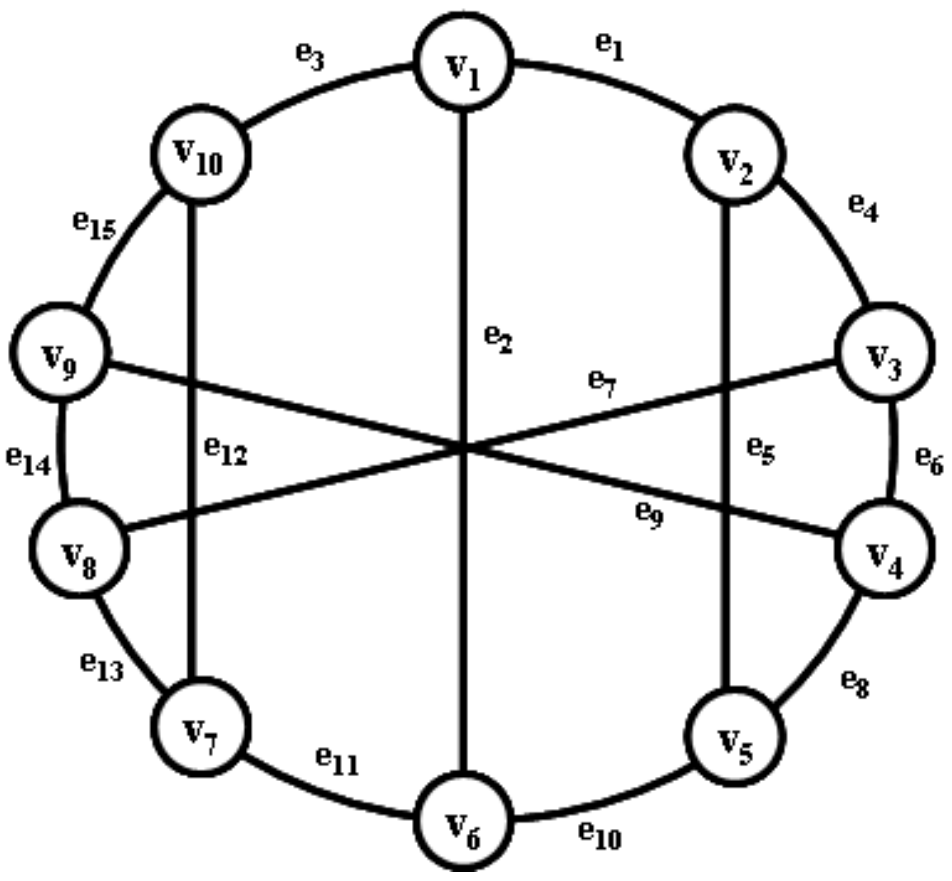

Рис. 8.37. Граф $G_{24}$.

Количество изометрических циклов в реберном графе $L(G_{29})$ = 27

Кортеж весов ребер $\xi_L(G_{23})$= <5,5,4,4,7,7,7,7,3,7,7,7,4,4,5,5>.

Кортеж весов вершин $\zeta_L(G_{23})$= <14,16,16,28,14,28,14,16,16,14>.

Вектор весов ребер : $F(\xi_L(G_{23}))$ = (3,4,4,4,4,5,5,5,5,7,7,7,7,7,7,7).

Вектор весов вершин : $F(\zeta_L(G_{23}))$ = (14,14,14,14,16,16,16,16,28,28).

Цифровой инвариант реберного графа $IL(G_{23}) = F(\xi_L(G_{23}))$ & $F(\zeta_L(G_{23}))$ =

= $(1\times 3, 4\times 4, 4\times 5, 7\times 7)$ & $(4\times 14, 4\times 16, 2\times 28)$.

Количество изометрических циклов в реберном графе $L(G_{24})$= 27

Кортеж весов ребер $\xi_L(G_{24})$= <5,5,4,4,7,7,7,5,5,7,4,7,7,4,5,5>.

Кортеж весов вершин $\zeta_L(G_{24})$= <14,16,16,26,16,16,26,16,16,14>.

Вектор весов ребер : $F(\xi_L(G_{24}))$ = (4,4,4,4,5,5,5,5,5,5,7,7,7,7,7,7)

Вектор весов вершин : $F(\zeta_L(G_{24}))$ = (14,14,16,16,16,16,16,16,26,26).

Цифровой инвариант реберного графа $IL(G_{24}) = F(\xi_L(G_{24}))$ & $F(\zeta_L(G_{24}))$ =

$= (4\times4,\ 6\times5,\ 6\times7)\ \&\ (2\times14,\ 6\times16,\ 2\times26).$

## Комментарии

В данной главе рассмотрены примеры решения задачи изоморфизма графов с применением инварианта реберных разрезов и инварианта реберных циклов. Основой для различения внутренних структур графов является спектр реберных разрезов и спектр реберных циклов графа. Показано, что спектр реберных разрезов, как и спектр реберных циклов, не зависят от перенумерации вершин графа.

Применение операции кольцевого сложения для суграфов реберных разрезов позволяет построить спектр реберных разрезов графа в виде элементов подпространства разрезов S(G). Это позволяет построить цепочку реберных разрезов для каждого отдельного ребра. Полученные цепочки реберных разрезов для ребра позволяют проводить сравнительный анализ для графов имеющих одинаковое количество вершин и ребер. Количество подмножеств в спектре реберных разрезов определяет меру участия данного ребра в структуре графа и называется весом ребра. На основе понятия веса ребра в спектре реберных разрезов строится векторный инвариант графа относительно спектра реберных разрезов и определяются кортежи весов ребер и вершин.

Множество изометрических циклов графа позволяет построить спектр реберных циклов графа в виде элементов подпространства циклов C(G). Это позволяет построить базовые реберные циклы относительно ребер графа. Базовые реберные циклы в сочетании со спектром базовых разрезов позволяют проводить сравнительный анализ структур графа.

Множество изометрических циклов реберного графа L(G) позволяет построить цифровой инвариант реберного графа и определить структуру графа G.

Рассмотрены примеры определения изоморфизма для некоторых классов графов.

## Выводы

В работе показано, что наряду с вектором локальных степеней графа, для различения внутренней структуры графа применимы следующие векторные инварианты, для построения которых применимы алгоритмы с полиномиальной вычислительной сложностью:

- IS(G) – векторный инвариант спектра реберных разрезов графа G;
- $IS_l$(G) – векторный инвариант уровня *l* спектра реберных разрезов графа G;
- IC(G) – векторный инвариант реберных циклов графа G;
- $IC_\tau$(G) – векторный инвариант множества изометрических циклов несепарабельного графа G;
- IS(G) & IC(G) – векторный интегральный инвариант графа G;
- IL(G) – векторный цифровой инвариант реберного графа L(G).

Данная система векторных инвариантов может быть использована только для неориентированных несепарабельных графов, так как в ее основе применяются методы цикломатики графов.

Основой для построения математических моделей векторных инвариантов служит структура реберного графа L(G) (согласно теореме Уитни), состоящая из прообразов центральных разрезов и изометрических циклов графа G. Очевидно, что векторный цифровой инвариант реберного графа возможно представить как полный инвариант графа.

Показано, что данные векторные инварианты не зависят от нумерации вершин графа.

Применение операции кольцевого сложения для суграфов реберных разрезов позволяет построить спектр реберных разрезов графа в виде элементов подпространства разрезов S(G). Это позволяет построить цепочку реберных разрезов для каждого отдельного ребра.

Полученные цепочки реберных разрезов для ребра позволяют проводить сравнительный анализ для графов имеющих одинаковое количество вершин и ребер. Множество суграфов спектра реберных разрезов определяет меру участия данного ребра в структуре графа и называется весом ребра. На основе построения кортежей весов ребер и вершин строится векторный инвариант спектра реберных разрезов графа.

Множество изометрических циклов графа позволяет построить спектр реберных циклов графа в виде элементов подпространства циклов C(G). Это позволяет получить реберные циклы для определения векторного инварианта реберных циклов графа

Существуют неизоморфные графы, имеющие одинаковые значения весов реберных разрезов и существуют неизоморфные графы, имеющие одинаковые значения весов реберных циклов. Поэтому для различения структур таких графов следует применять интегральные инварианты графа. Однако применение интегрального инварианта, согласно лемме о изоморфизме пространства суграфов, теоретически возможно только для графов имеющих непересекающиеся элементы подпростронства разрезов S(G) и подпространства циклов C(G). Применение интегрального инварианта для графов имеющих пересекающиеся элементы подпространства разрезов и подпространства циклов проблематично, так как теоретически не обосновано.

Векторный инвариант реберного графа строится на основе структуры реберного графа L(G) и может быть описан триадами составляющих центральных разрезов графа G, изометрическими циклами графа G и суграфами дубль-циклов графа G. Суграфы дубль-циклов характеризуют плоские фрагменты состоящие из системы изометрических циклов графа G.

Векторный интегральный инвариант, хотя это теоретически не обоснованно, можно применять и для оценки графов, имеющих пересекающиеся суграфы, принадлежащие подпространствам разрезов и циклов. Так как интегральный инвариант строится теми же структурами графа G, которые составляют реберный граф L(G), но отличается способом построения и строится без учета влияния дубль-циклов.

Вычислительную сложность алгоритма определения инварианта спектра реберных разрезов графа можно определить относительно ребер графа как $O(qm^3)$. В случае, когда веса вершин определяются только первыми двумя уровнями спектра реберных разрезов, вычислительную сложность алгоритма можно определить как $O(m^3)$. Задача определения инварианта реберных разрезов графа относится к классу P – полиномиальных алгоритмов.

Вычислительную сложность алгоритма определения инварианта реберных циклов графа, можно определить как количество операций для построения множества изометрических циклов $O(n^4)$. Задача определения инварианта реберных циклов графа также относится к классу P – полиномиальных алгоритмов.

Интегральный инвариант объединяет в единое целое и инвариант реберных разрезов и инвариант реберных циклов – необходимое условие для решения задачи различения структуры графов. Вычислительную сложность алгоритма построения интегрального инварианта можно рассматривать как сумму вычислительных сложностей построения инвариантов реберных разрезов и реберных циклов $O(n^6) + O(n^4)$.

Вычислительная сложность алгоритма построения векторного инварианта реберного графа в основном определяется сложностью формирования изометрических циклов реберного графа и составляет величину $O(m^4)$ или $O(n^8)$.

Для распознавания изоморфизма деревьев следует применять только инвариант – спектр реберных разрезов, без ограничения на количество уровней.

В свою очередь, задача распознавания изоморфизма графов может быть распространена на решение следующих задач теории графов:

- определение изоморфного суграфа;
- поиск бинарного соответствия между двумя графами G и H;
- определение группы автоморфизмов графа *Aut*(G).

Данные задачи тесно связаны между собой основой построения векторного интегрального инварианта и вектора цифрового инварианта реберного графа, что позволяет проводить решение данных задач.

Рассмотрены примеры определения изоморфизма для некоторых классов графов.

По сути, рассмотренные математические модели и теории, говорят о том, что векторный цифровой инвариант можно рассматривать как полный инвариант графа. По осторожному мнению авторов, утверждение о том, что интегральный инвариант является полным инвариантом графа несколько преждевременно, так как к настоящему времени не накоплено достаточного фактического материала для такого утверждения. Однако, векторный интегральный инвариант можно рассматривать как аналог векторного цифрового инварианта без учета структур соответствующих дубль-циклам графа G.

## СПИСОК ИСПОЛЬЗУЕМОЙ ЛИТЕРАТУРЫ